\documentclass{memo-l}
\usepackage{etex} 
\usepackage{amssymb}
\usepackage{amscd}
\usepackage[noadjust]{cite}
\usepackage{booktabs}
\usepackage{url}
\usepackage{hyphenat}
\usepackage{mathtools}
\usepackage{cancel} 
\usepackage{ifpdf}
\usepackage[T1]{fontenc}
\usepackage[utf8]{inputenc}
\usepackage[curve]{xypic}
\entrymodifiers={+!!<0pt,\the\fontdimen22\textfont2>}
\usepackage{enumitem}
\usepackage{calc}
\ifpdf
  \usepackage[pdftex,final]{graphicx}
  \usepackage[bookmarks=true, bookmarksopen=true,%
    bookmarksdepth=3,bookmarksopenlevel=1,%
    colorlinks=true,%
    linkcolor=blue,%
    citecolor=blue,%
    filecolor=blue,%
    menucolor=blue,%
    urlcolor=blue]{hyperref}
  \hypersetup{pdftitle={Bordered Heegaard Floer homology}}
  \hypersetup{pdfauthor={Robert Lipshitz, Peter S. Ozsváth, and Dylan
      P. Thurston}}
\else
  \usepackage[dvips,final]{graphicx}
  \usepackage[dvips,lmargin=1in,rmargin=1in,tmargin=1in,bmargin=1in]{geometry}
  \usepackage[draft]{hyperref}
\fi

\usepackage{tikz}
\usetikzlibrary{arrows}
\tikzset{aar/.style={->, thick}}
\tikzset{taar/.style={double, double equal sign distance, -implies}}
\tikzset{amar/.style={->, dashed}}
\tikzset{amarp/.style={->, dotted}}
\tikzset{dmar/.style={->, dashed}}

\usepackage[sanitize=none, hyperfirst=false]{glossaries}

\makeatletter
\newcommand{\@dotsep}{4.5}
\renewcommand*{\dotfill}{
  \leavevmode\leaders
  \hbox{$\m@th \mkern \@dotsep mu\hbox{.}\mkern \@dotsep mu$}
  \hfill\kern\z@
}
\def\@tocline#1#2#3#4#5#6#7{\relax
  \ifnum #1>\c@tocdepth 
  \else
    \par \addpenalty\@secpenalty\addvspace{#2}%
    \begingroup \hyphenpenalty\@M
    \@ifempty{#4}{%
      \@tempdima\csname r@tocindent\number#1\endcsname\relax
    }{%
      \@tempdima#4\relax
    }%
    \parindent\z@ \leftskip#3\relax \advance\leftskip\@tempdima\relax
    \rightskip\@pnumwidth plus4em \parfillskip-\@pnumwidth
    #5\leavevmode\hskip-\@tempdima #6\nobreak\relax
    \dotfill\hbox to\@pnumwidth{\@tocpagenum{#7}}\par
    \nobreak
    \endgroup
  \fi}
\def\l@figure{\@tocline{0}{3pt plus2pt}{0pt}{2.8pc}{}}
\makeatother

\newread\testin

\graphicspath{{draws/}{mpdraws/}{}}
\makeatletter
\def\input@path{{}{draws/}}
\makeatother

\def\mathcenter#1{%
  \vcenter{\hbox{$#1$}}%
}

\def\mfig#1{
        \mathcenter{\includegraphics[scale=.8333]{#1}}
}

\def\mfigb#1{
        \mathcenter{\includegraphics[scale=.8333,trim=-1 -1 -1 -1]{#1}}
}

\DeclareRobustCommand{\widebar}[1]{\overline{#1}{}}
\DeclareRobustCommand{\widedblbar}[1]{\overline{\overline{#1}}{}}

\makeatletter
\newcommand\mi@kern[1]{%
  \settowidth\@tempdima{$\mi@obj^{#1}$}
  \kern-\@tempdima
  #1
  \settowidth\@tempdima{$\mi@obj$}
  \kern\@tempdima
}

\newtoks\mi@toksp
\newtoks\mi@toksb
\DeclareRobustCommand{\manyindices}[5]{
  \def\mi@obj{#5}
  \mi@toksp\expandafter{\mi@kern{#2}}
  \mi@toksb\expandafter{\mi@kern{#1}}
  \@mathmeasure4\textstyle{#5_{#1}^{#2}}
  \@mathmeasure6\textstyle{#5_{#3}^{#4}}
  \dimen0-\wd6 \advance\dimen0\wd4
  \@mathmeasure8\textstyle{\hphantom{{}_{#1}^{#2}}#5^{\the\mi@toksp#4}_{\the\mi@toksb#3}}
  \hbox to \dimen0{}{\kern-\dimen0\box8}
}
\makeatother 

\usepackage{ifpdf}
\ifpdf
  \let\textalt\texorpdfstring
\else
  \newcommand{\textalt}[2]{#1}
\fi

\newcommand{\RR}{\mathbb R}
\newcommand{\CC}{\mathbb C}
\newcommand{\ZZ}{\mathbb Z}
\newcommand{\QQ}{\mathbb Q}

\newcommand{\DD}{\mathbb D}
\newcommand{\FF}{\mathbb F}
\newcommand{\NN}{\mathbb N}
\newcommand{\Ker}{\mathrm{Ker}}
\newcommand{\Image}{\mathrm{Im}}

\newcommand{\connectsum}{\mathbin \#}

\newcommand{\co}{\nobreak\mskip2mu\mathpunct{}\nonscript
  \mkern-\thinmuskip{:}\penalty300\mskip6muplus1mu\relax}
\newcommand{\semico}{;\penalty 300 }
\newcommand{\eps}{\epsilon}
\newcommand{\abs}[1]{{\lvert #1 \rvert}}

\newcommand{\OneHalf}{{\textstyle\frac{1}{2}}}
\newcommand{\OneQuart}{{\textstyle\frac{1}{4}}}

\newcommand{\lbracket}{[}
\newcommand{\rbracket}{]}

\newcommand{\isom}{\simeq}

\newcommand{\bdy}{\partial}
\newcommand\piBig{{\widetilde\pi}_2}

\newcommand{\llbracket}{[\![}
\newcommand{\rrbracket}{]\!]}

\newcommand\Unit{\mathbf I}
\newcommand\SetS{\mathbf{s}}
\newcommand\SetT{\mathbf{t}}
\newcommand\SetU{\mathbf{u}}
\newcommand{\spinc}{\mathfrak s}
\DeclareMathOperator{\divis}{div}

\DeclareMathOperator{\Sym}{Sym}

\DeclareMathOperator{\Hom}{Hom}

\DeclareMathOperator{\rank}{rank}
\DeclareMathOperator{\im}{im}
\newcommand{\spin}{\text{spin}}
\newcommand{\SpinC}{\spin^c}
\DeclareMathOperator{\BSpinC}{Bspin^c}

\DeclareMathOperator{\ind}{ind}
\DeclareMathOperator{\eDim}{ind}
\DeclareMathOperator{\sing}{sing}

\newcommand\AGrading{{\widetilde \smallGroup}}
\newcommand\AGradingD{\widetilde{S}_D}
\newcommand\Agr{\widetilde\gr}

\DeclareMathOperator{\inv}{inv}
\DeclareMathOperator{\Inv}{Inv}
\DeclareMathOperator{\ev}{ev}
\DeclareMathOperator{\gr}{gr}
\DeclareMathOperator{\tgr}{\underline{\gr}}
\DeclareMathOperator{\grT}{gr^{\DT}}
\DeclareMathOperator{\Area}{Area}
\DeclareMathOperator{\br}{br} 

\DeclareMathOperator{\Int}{Int}
\newcommand{\rel}{{\mathrm{rel}}} 
\newcommand{\emb}{{\mathrm{emb}}} 
\newcommand{\Winding}{\mathcal W}

\newcommand{\op}{\mathrm{op}}

\DeclareMathOperator{\BSO}{\mathit{BSO}}

\theoremstyle{plain}
\numberwithin{equation}{chapter}
\newtheorem{theorem}[equation]{Theorem}
\newtheorem{citethm}[equation]{Theorem}
\newtheorem{proposition}[equation]{Proposition}
\newtheorem{lemma}[equation]{Lemma}

\newtheorem{corollary}[equation]{Corollary}

\newtheorem{observation}[equation]{Observation}

\theoremstyle{definition}
\newtheorem{definition}[equation]{Definition}

\newtheorem{construction}[equation]{Construction}

\theoremstyle{remark}
\newtheorem{example}[equation]{Example}
\newtheorem{remark}[equation]{Remark}

\renewcommand{\thefigure}{\thechapter.\arabic{figure}}

\newcounter{step}
\newcounter{Step}

\newenvironment{steps}{%
  \begingroup%
  \setcounter{step}{0}%
  \newcommand{\step}{\stepcounter{Step}\refstepcounter{step}\textbf{Step \thestep.} }%
}{\endgroup}

\newcounter{claim}
\newcounter{Claim}

\newenvironment{claims}{%
  \begingroup%
  \setcounter{claim}{0}%
  \newcommand{\claim}{\stepcounter{Claim}\refstepcounter{claim}\textbf{Claim \theclaim.} }%
}{\endgroup}

\newcommand{\HF}{\mathit{HF}}
\newcommand{\HFa}{\widehat {\HF}}

\newcommand{\CFa}{\widehat {\mathit{CF}}}
\newcommand{\tCFa}{\underline{\CFa}}
\newcommand{\HFKa}{\widehat {\mathit{HFK}}}
\newcommand{\x}{\mathbf x}
\newcommand{\y}{\mathbf y}
\newcommand{\z}{\mathbf z}
\newcommand{\w}{\mathbf w}

\newcommand\HH{\mathit{HH}}
\newcommand\Hochschild\HH

\newcommand{\MCG}{\mathit{MCG}}
\newcommand\Cvert{C^{\mathrm{v}}}
\newcommand\Chor{C^{\mathrm{h}}}
\newcommand\Hvert{H^{\mathrm{v}}}
\newcommand\Hhor{H^{\mathrm{h}}}
\newcommand{\partialvert}{\partial^{\mathrm{v}}}
\newcommand{\partialhor}{\partial^{\mathrm{h}}}
\newcommand{\Ainf}{\mathcal A_\infty}
\newcommand{\Alg}{\mathcal{A}}

\newcommand{\bigAlg}{\overline{\Alg}}
\newcommand{\cA}{{\mathcal{A}}}
\newcommand{\cB}{{\mathcal{B}}}
\newcommand{\Idem}{\mathcal{I}}
\newcommand{\alphas}{{\boldsymbol{\alpha}}}
\newcommand{\betas}{{\boldsymbol{\beta}}}
\newcommand{\gammas}{{\boldsymbol{\gamma}}}
\newcommand{\rhos}{{\boldsymbol{\rho}}}
\newcommand{\etas}{\boldsymbol{\eta}}
\newcommand{\sigmas}{{\boldsymbol{\sigma}}}
\newcommand{\bSigma}{\widebar{\Sigma}}
\newcommand{\balpha}{\widebar{\alpha}}
\newcommand{\balphas}{\widebar{\alphas}}
\newcommand{\cM}{\mathcal{M}}
\newcommand{\Mod}{\cM}

\newcommand{\cN}{\mathcal{N}}
\newcommand{\ocN}{\widebar{\mathcal{N}}}
\newcommand{\tcN}{\widetilde{\mathcal{N}}}
\newcommand{\ocM}{\widebar{\cM}{}} 
\newcommand{\oocM}{\widedblbar{\cM}{}}
\newcommand{\ocMM}{\widebar{\cMM}{}} 
\newcommand{\tcMM}{\widetilde{\mathcal{MM}}}
\newcommand{\cMM}{\mathcal{MM}}
\newcommand{\tcM}{\widetilde{\mathcal{M}}}
\newcommand{\CFD}{\mathit{CFD}}
\newcommand{\CFDD}{\mathit{CFDD}}
\newcommand{\CFA}{\mathit{CFA}}

\newcommand{\CFDA}{\mathit{CFDA}}
\newcommand{\CFDAa}{\widehat{\CFDA}}
\newcommand{\CFAA}{\mathit{CFAA}}
\newcommand{\CFAAa}{\widehat{\CFAA}}
\newcommand{\CFDa}{\widehat{\CFD}}

\newcommand{\CFDm}{\CFD^-}
\newcommand{\CFAm}{\CFA^-}
\newcommand{\CFK}{\mathit{CFK}}
\newcommand{\CFKa}{\widehat{\CFK}}
\newcommand{\CFKm}{\CFK^-}
\newcommand{\HFK}{\mathit{HFK}}
\newcommand{\HFKm}{\HFK^-}
\newcommand{\gCFKm}{\mathit{gCFK}^-}
\newcommand{\gCFKa}{\widehat{\mathit{gCFK}}}
\newcommand{\CFDDa}{\widehat{\CFDD}}
\newcommand{\tCFDa}{\underline{\CFDa}}
\newcommand{\CFAa}{\widehat{\CFA}}
\newcommand{\tCFAa}{\underline{\CFAa}}
\newcommand{\s}{\mathfrak{s}}
\newcommand{\Source}{{S^{\mspace{1mu}\triangleright}}}
\newcommand{\SourceSub}[1]{{S_{#1}^{\mspace{1mu}\triangleright}}}
\newcommand{\biSource}{T^{\mspace{1mu}\rotatebox{90}{$\scriptstyle\lozenge$}}}
\newcommand{\TSource}{{T^\Delta}}
\newcommand{\closedSource}{S^{\mspace{1mu}\lozenge}}
\newcommand{\cZ}{\mathcal{Z}}
\newcommand{\PtdMatchCirc}{\cZ}
\newcommand{\PMC}{\PtdMatchCirc}

\newcommand{\CircPts}{{\mathbf{a}}}
\newcommand{\glue}{\mathbin{\natural}}

\newcommand\DGA{A}
\newcommand{\dg}{\textit{dg} }
\newcommand{\cDGA}{\mathcal{A}}
\newcommand{\cModule}{\mathcal{M}}
\newcommand{\cNodule}{\mathcal{N}}
\newcommand\Id{\mathbb{I}}
\newcommand\Ground{\mathbf k}

\newcommand\Mas{\ind}
\newcommand\DTP{\mathbin{\widetilde\otimes}}
\newcommand\DT{\boxtimes}
\newcommand\Gen{\mathfrak{S}}
\newcommand\tGen{\underline{\mathfrak{S}}}
\renewcommand{\S}{\Gen}
\newcommand\ModFlow{\tcM}
\newcommand\Tensor{\mathcal T}
\newcommand\Zmod[1]{\mathbb{Z}/{#1}\mathbb{Z}}
\newcommand{\Field}{\FF_2}
\DeclareMathOperator{\nbd}{nbd}
\DeclareMathOperator{\sgn}{sgn}
\newcommand\UnparModFlow{\Mod}
\newcommand{\geqo}{\mathord{\geq}}
\newcommand{\leqo}{\mathord{\leq}}
\newcommand{\dbar}{\bar{\partial}}
\newcommand{\Heegaard}{\mathcal{H}}
\newcommand{\HD}{\Heegaard}

\newcommand{\HDst}{\HD_{\mathit{st}}}
\newcommand{\cF}{\mathcal{F}}
\renewcommand{\th}{^\text{th}}
\newcommand{\st}{^\text{st}}
\newcommand{\bigGroup}{G'}
\newcommand{\smallGroup}{G}
\newcommand{\ABigGrSet}{S'_A}
\newcommand{\tABigGrSet}{\underline{S}'_A}
\newcommand{\DBigGrSet}{S'_D}
\newcommand{\tDBigGrSet}{\underline{S}'_D}
\newcommand{\DSmallGrSet}{S_D}
\newcommand{\tDSmallGrSet}{\underline{S}_D}
\newcommand{\ASmallGrSet}{S_A}
\newcommand{\tASmallGrSet}{\underline{S}_A}
\newcommand{\grb}{\gr'}
\newcommand{\tgrb}{\underline{\gr}'}
\newcommand{\gb}{g'}
\newcommand{\gs}{g}
\newcommand{\Ags}{\widetilde{\gs}}
\newcommand{\APeriodics}{\widetilde{P}}
\newcommand{\Torus}{\mathbb{T}}

\DeclareMathOperator{\Mor}{Mor}
\DeclareMathOperator{\ob}{Ob}

\newcommand{\tsic}{\mathit{tsic}}
\newcommand{\cU}{\mathcal{U}}
\newcommand{\cV}{\mathcal{V}}
\newcommand{\arcz}{\mathbf{z}}
\newcommand{\pd}{\mathrm{pd}}
\newcommand{\HalfPlane}{\mathbb{H}}
\newcommand{\tB}{\widetilde{B}}

\newcommand{\symmdiff}{\mathbin{\Delta}}

\newcommand{\DDm}{\textit{DD}}
\newcommand{\DAm}{\textit{DA}}

\newcommand{\AAm}{\textit{AA}} 

\newcommand\strands[2]{\bigl\langle\begin{smallmatrix}
#1 \\ #2
\end{smallmatrix}\bigr\rangle}
\newcommand\reebchords[2]{\bigl\{\begin{smallmatrix}
#1 \\ #2
\end{smallmatrix}\bigr\}}
\makeatletter
\newcommand\honestalg[3]{\bigl\lbracket
\begin{smallmatrix} #1\@ifempty{#3}{}{&#3} \\ #2 \end{smallmatrix}
\bigr\rbracket}
\newcommand{\sos}[3]{\mathbin{{}_{#1}\mathord#2_{#3}}}
\newcommand\honestalga[3]{\bigl\lbracket
\begin{smallmatrix} #1&#3 \\ #2 \end{smallmatrix}
\bigr\rbracket}
\makeatother
\newcommand{\lab}[1]{$\scriptstyle #1$}

\newcommand\gxx{x_2x_3}
\newcommand\gxy{x_2y_3}
\newcommand\gyx{y_2x_3}
\newcommand\gyy{y_2y_3}

\newcommand{\colorused}{\relax}

\newcommand{\lsub}[2]{{}_{#1}#2}

\makeglossaries

\makeindex

\usepackage{xr-hyper}
\input{flip}

\begin{document}
\frontmatter
\title{Bordered Heegaard Floer homology}

\author[Lipshitz]{Robert Lipshitz}
\thanks{RL was supported by an NSF Mathematical Sciences Postdoctoral
  Research Fellowship, NSF grants DMS-0905796 and DMS-1149800 and a Sloan Research Fellowship.}
\address{Department of Mathematics, University of Oregon\\
  Eugene, OR 97403}
\email{lipshitz@uoregon.edu}

\author[Ozsv\'ath]{Peter S.~Ozsv\'ath}
\thanks{PSO was supported by NSF grants DMS-0505811, DMS-0804121, and DMS-1105810.}
\address {Department of Mathematics, Princeton University\\ 
Princeton, NJ 08544}
\email {petero@math.princeton.edu}

\author[Thurston]{Dylan P.~Thurston}
\thanks {DPT was supported by NSF
  grants DMS-1008049 and DMS-1358638 and a Sloan Research Fellowship.}
\address{Department of Mathematics, Indiana University\\
  Bloomington, IN 47405}
\email{dpthurst@indiana.edu}

\date{August 10, 2014}


\maketitle

\tableofcontents

\begin{abstract}
  We construct Heegaard Floer theory for 3-manifolds with connected
  boundary.  The theory associates to an oriented, parametrized two-manifold a
  differential graded algebra.  For a three-manifold with parametrized
  boundary, the invariant comes in two different versions, one of
  which (type $D$) is a module over the algebra and the other of which
  (type $A$) is an $\Ainf$ module. Both are well-defined up to chain
  homotopy equivalence. For a decomposition of a 3-manifold into two
  pieces, the $\Ainf$ tensor product of the type $D$ module of one
  piece and the type $A$ module from the other piece is $\HFa$ of the
  glued manifold.

  As a special case of the construction, we specialize to the case of
  three-manifolds with torus boundary. This case can be used to give
  another proof of the surgery exact triangle for $\HFa$. We relate
  the bordered Floer homology of a three-manifold with torus boundary
  with the knot Floer homology of a filling.
\end{abstract}

\let\contentsname\listfigurename
\listoffigures

\input{glossary}

\mainmatter
\activateflip

\chapter{Introduction}
\label{chap:intro}

\section{Background}

\colorused

Since the pioneering work of Simon Donaldson, techniques from gauge
theory have taken a central role in the study of smooth four-manifold
topology~\cite{DonKron}. His numerical invariants, associated to
closed, smooth four-manifolds with $b_2^+ > 1$, have shed much light on our
understanding of differential topology in dimension four. Moreover,
these invariants, and the subsequent closely-related Seiberg-Witten
invariants~\cite{Witten} and Heegaard Floer invariants~\cite{OS04:HolomorphicDisks}
all fit into a formal framework reminiscent of the ``topological
quantum field theories'' proposed by Witten~\cite{TQFT}.  Crudely speaking, these
theories have the following form.  To a closed three-manifold $Y$ one
associates a (suitably graded) abelian group, the {\em Floer homology
  of $Y$}, and, to a four-manifold~$W$ with boundary identified with
$Y$, a homology class in the Floer homology of the boundary. If
a closed, smooth four-manifold $X$ decomposes along $Y$ into a union
of two four-manifolds with boundary $X_1$ and $X_2$, then the
numerical (Donaldson, Seiberg-Witten, or Heegaard Floer) invariant of
the closed four-manifold is obtained as a suitable pairing between the relative
invariants coming from $X_1$ and~$X_2$ in a corresponding version of
Floer homology of~$Y$.

As the name suggests, the first such construction was proposed by Andreas
Floer (for a restricted class of three-manifolds) as a tool for
studying Donaldson's theory. A complete construction of the
corresponding three-dimensional Floer theory for Seiberg-Witten
invariants was given by Kronheimer and
Mrowka~\cite{KronheimerMrowka}. Heegaard Floer homology was defined
by Zolt{\'a}n Szab{\'o} and the second author.

The aim of the present work is to perform a corresponding
construction one dimension lower. Specifically, we produce an
invariant which, loosely speaking, associates to a parametrized,
closed, oriented
surface~$F$ a differential graded algebra $\Alg(F)$ and associates to a
three-manifold whose boundary is identified with $F$ a differential
graded module over
$\Alg(F)$.
When a closed, oriented three-manifold can be decomposed
along $F$ into two pieces $Y_1$ and $Y_2$, a suitable variant of Floer
homology is gotten as a pairing of the differential graded modules
associated to $Y_1$ and $Y_2$.

We give now a slightly more detailed version of
this picture, starting with some more remarks about Heegaard Floer
homology, and then a more precise sketch of the invariants constructed
in the present book.

Recall that there are several variants of Heegaard Floer homology
stemming from the fact that, in its most basic form, the Heegaard
Floer homology of a three-manifold is the homology of a chain complex
defined over a polynomial ring in an indeterminant~$U$.  The full theory can
be promoted to construct invariants for closed,
smooth four-manifolds~\cite{OS06:HolDiskFour}, similar in character (and
conjecturally equal) to Seiberg-Witten invariants.  In this
book, we focus on the specialization of Heegaard Floer homology for
three-manifolds in the case where $U=0$, giving the three-manifold
invariant denoted $\HFa(Y)$. The corresponding simplified theory is
not rich enough to construct interesting closed four-manifold
invariants, but it does already contain interesting geometric
information about the underlying three-manifold including, for example,
full information about the minimal genera
of embedded surfaces in~$Y$~\cite{OS04:GenusBounds,Ni:spheres}.

With this background in place, we proceed as follows. Let $F$ be a
closed, oriented two-manifold. Equip $F$ with some additional
data, in the form of a minimal handle decomposition and a basepoint,
data which we refer to as a {\em pointed matched circle} (see
Definition~\ref{def:PointedMatchedCircle}) and denote by $\PMC$. 
This handle decomposition can be
thought of as giving a parameterization of $F$. We associate to $\PMC$
a differential graded algebra~$\Alg(\PMC)$.

Next, fix a
three-manifold with
boundary~$Y$ together with an orientation\hyp preserving diffeomorphism
$\phi\co F\to\partial Y$. To this data $(Y,\phi)$ we associate a left
differential module over $\Alg(-\PMC)$, the \emph{type~$D$ module of $Y$},
denoted $\CFDa(Y, \phi\co F\to Y)$ or, less precisely, $\CFDa(Y)$. 
To $(Y,\phi)$ we also associate the \emph{type~$A$ module of $Y$},
denoted $\CFAa(Y,\phi\co F\to Y)$ or
$\CFAa(Y)$. The type $A$ module of $Y$ invariant is a more general type
of object: it is a right $\Ainf$ module over $\Alg(\PMC)$. Both
invariants are
defined by counting holomorphic curves.
The definitions of both of these modules depend on a number
of auxiliary choices, including compatible Heegaard diagrams and
associated choices of almost complex structures. However, as is
typical of Floer homology theories, the underlying ($\Ainf$) homotopy equivalence
types of the modules are independent of these choices,
giving rise to topological invariants of our
three-manifold with boundary.  When speaking informally about these
invariants, we refer to them collectively as the {\em bordered
  Heegaard Floer invariants} of~$Y$.

The relationship with closed invariants is encapsulated as follows.
Suppose that $Y_1$ and $Y_2$ are two three-manifolds with boundary,
where $\partial Y_1$ is identified with $F$ and $\partial Y_2$ is
identified with $-F$. Then we can form a closed three-manifold $Y$ by
identifying $Y_1$ and $Y_2$ along their boundaries.  The pairing
theorem for bordered Heegaard Floer homology states that the
Heegaard Floer complex $\HFa$ of $Y$ is homotopy equivalence to the
$\Ainf$ tensor product of $\CFAa(Y_1)$ with $\CFDa(Y_2)$.

Recall that $\HFa$ can be calculated algorithmically for any closed
three-manifold, thanks to the important work of Sarkar and
Wang~\cite{SarkarWang07:ComputingHFhat}.  Nevertheless, the present
work (and forthcoming extensions) makes it possible to compute it in
infinite families: by cutting a 3-manifold into simpler pieces along
surfaces and computing an invariant for each piece, we can reduce the
computation to a computation for each piece and an algebraic
computation.  Indeed, we expect this to lead to a more efficient
algorithm for
computation in general.

Moreover, bordered
Heegaard Floer homology provides a conceptual framework for
organizing the structure of $\HFa$. We give a few of these properties in
this book, and return to further applications elsewhere
(see Section~\ref{sec:further-developments}).

\section{The bordered Floer homology package}

The formal algebraic setting for the constructions herein lie
outside the working toolkit for the typical low-dimensional
topologist: $\Ainf$ modules over differential graded algebras, and
their $\Ainf$ tensor products (though, we should point out, these
objects are now familiar in symplectic topology,
cf.~\cite{Fukaya,Kontsevich,Seidel02:FukayaDef}). We will recall the basics before
giving precise statements of our results. But before doing this, we
sketch the outlines of the geometry which underpins the constructions,
giving informal statements of the basic
results of this package.

Central to this geometric picture is the cylindrical reformulation of
Heegaard Floer homology~\cite{Lipshitz06:CylindricalHF}.  Recall that Heegaard Floer homology in its
original incarnation~\cite{OS04:HolomorphicDisks} is an invariant associated to a Heegaard diagram
for a three-dimensional manifold $Y$. Start from a Heegaard diagram
for $Y$, specified by an oriented surface $\Sigma$, equipped with two
$g$-tuples of attaching circles $\{\alpha_1,\dots,\alpha_g\}$ and
$\{\beta_1,\dots,\beta_g\}$, as well as an additional basepoint $z\in \Sigma$
in the complement of these circles. One then considers a version of
Lagrangian Floer homology in the $g$-fold symmetric product
of~$\Sigma$. In the cylindrical reformulation, disks in the $g$-fold
symmetric product are reinterpreted as holomorphic curves in the
four-manifold $\Sigma\times [0,1]\times \RR$, with coordinates
$(x,s,t)$, satisfying certain constraints (dictated by the attaching
circles) at $s\in \{0,1\}$ and asymptotic constraints as $t$ approaches
$\pm\infty$.

Starting from this cylindrical reformulation, suppose that our
Heegaard diagram for $Y$ is equipped with a closed, separating curve
$Z$ in $\Sigma$. We suppose that $Z$ meets the following combinatorial
requirements: it is disjoint from the $\beta$-curves---so cutting
$\Sigma$ along $Z$ and the $\beta$-curves results in a surface with
two components---and although we do not require that $Z$ be disjoint
from the $\alpha$-curves, we still require that the result of cutting
$\Sigma$ along $Z$ and the $\alpha$-curves results in a surface with
two components. The constructions in this book emerge when one considers limits of
holomorphic curves in $\Sigma\times[0,1]\times \RR$ as the complex
structure on $\Sigma$ is pinched along $Z$.

\begin{figure}
  \centerline{\includegraphics[scale=.83333]{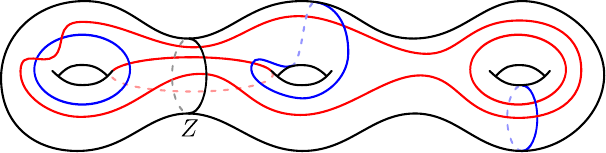}}
  \caption[Example of cutting a Heegaard diagram]{\textbf{An
      example of cutting a Heegaard diagram.}  The Heegaard diagram
    shown may be degenerated by cutting along (or pinching) the curve marked~$Z$.
    Each half of the result will be a bordered Heegaard diagram.}
  \label{fig:CutHeegaard}
\end{figure}

In the limit as the circle $Z$ is pinched to a node, the
holomorphic curves limit to holomorphic curves in $((\Sigma_1\amalg
\Sigma_2) \times [0,1]\times\RR)$, where here $\Sigma_1$ and
$\Sigma_2$ are the two components of $\Sigma-Z$ equipped with a cylindrical ends. The
holomorphic curves from these two halves have constrained limiting
behavior as they enter the ends of $\Sigma_1$ and~$\Sigma_2$: 
there are finitely many values of $t$ where our holomorphic curves 
are asymptotic to arcs in the ideal boundaries $\partial
\Sigma_i$. The data of
the asymptotics (the values of $t$ and the arcs in $\partial
\Sigma_i$) must match in the limit.  Thus,
we may reconstruct holomorphic curve counts for
$\Sigma\times[0,1]\times \RR$ from holomorphic curve counts in
$\Sigma_i\times[0,1]\times \RR$.

The data on the Heegaard surface specified by $Z$ has a more
intrisic, three\hyp dimensional interpretation.  The $\alpha$-circles meet
$Z$ generically in a collection of points, which come in pairs
(those points which belong to the same $\alpha_i$). This is the
structure of a \emph{matched circle} mentioned above.  It specifies a closed
surface~$F$.  This can be seen in terms of Morse functions, as follows.
Think about the Heegaard diagram  as
the intersections of the ascending and descending disks of a
self-indexing Morse function on~$Y$ with the middle level in the usual
way.  We then construct
$F$ from $Z$ by flowing $Z$ backwards and forwards under the
Morse flow.  Concretely, $F$ is obtained from a disk with boundary
$Z$ by attaching one-handles along the matched
pairs of points, and then filling in the remaining disk with a
two-handle. (Our assumptions guarantee there is only one two-handle to
be added in the end.)  The two halves of a
Heegaard diagram specify three-manifolds $Y_1$ and $Y_2$ meeting along
their common boundary~$F$.

In general terms, the differential graded algebra $\Alg(\PMC)$
is constituted from the Reeb chords of a corresponding matched
circle~$\PMC$.
(The reader should not be thrown by the use of the term ``Reeb chord'' here: the ambient contact manifold here is the circle, and hence the Reeb chords
are simply arcs in a circle, with boundaries on the feet of the
$1$-handles. That said, the terminology is not entirely pointless; see
Chapter~\ref{chap:structure-moduli}.)
The product on the algebra and its differential are defined
combinatorially, but rigged to agree with
(relatively simple) holomorphic curve counts in the cylinder 
which interpolates between Reeb chords. The precise algebra is defined in
Chapter~\ref{chap:algebra}.

The type $D$ module for the component $Y_2$ is defined as a chain
complex over $\Alg(\PMC)$, generated by $g_2$-tuples of intersection
points of various of the $\alpha_i$ and $\beta_j$ contained
in~$\Sigma_2$, with a differential given as a weighted count of rigid
holomorphic curves, where the coefficient in $\Alg(\PMC)$ of a rigid
holomorphic curve measures the asymptotics of its Reeb chords at its
boundary.  The precise definition is given in
Chapter~\ref{chap:type-d-mod}.  In that chapter, we also prove a more
precise version (Theorem~\ref{thm:D-invariance}) of the following:

\begin{theorem}
  \label{intro:D-invariance}
  Let $Y_2$ be a three-manifold equipped with an orientation\hyp preserving
  diffeomorphism $\phi\co F \to -\partial Y_2$. The
  homotopy equivalence class of the associated differential module
  $\CFDa(Y_2,\phi\co F \to -\partial Y_2)$ is a topological invariant of the three-manifold $Y_2$ with boundary parameterized by $-F$.
\end{theorem}

For the type~$A$ module for the component~$Y_1$, on the other hand, the
algebraic operations are defined by counting holomorphic
curves that are rigid when subject to certain height constraints on
their
Reeb chords. Specifically, we require that certain clumps of these
Reeb chords occur at the same height.  The proper algebraic set-up for
these operations is that of an $\Ainf$ module over $\Alg(\PMC)$, and the
precise definition of the type $A$ module is given in
Chapter~\ref{chap:type-a-mod}, where we also prove a more precise
version (Theorem~\ref{thm:A-invariance}) of the following:

\begin{theorem}
  \label{intro:A-invariance}
  Let $Y_1$ be a three-manifold equipped with an
  orientation\hyp preserving diffeomorphism $\phi\co F \to \partial
  Y_1$.  The $\Ainf$ homotopy equivalence class of the $\Ainf$ module
  $\CFAa(Y_1,\phi\co F\to \partial Y_1)$ is a topological invariant
  of the three-manifold $Y_1$ with boundary parameterized by~$F$.
\end{theorem}

Of course, differential graded modules over a differential graded
algebra are special cases of $\Ainf$ modules. Moreover, there is a
pairing, the {\em $\Ainf$ tensor product}, between two $\Ainf$ modules
$M_1$ and $M_2$, giving rise to a chain complex $M_1\DTP M_2$,
whose chain homotopy type
depends only on the chain homotopy types of the two factors; this
generalizes the derived tensor product of ordinary modules.

The Heegaard Floer homology of $Y$ can be reconstituted from the
bordered Heegaard Floer homology of its two components according to the
following pairing theorem:

\begin{theorem}
  \label{thm:TensorPairing}
  Let $Y_1$ and $Y_2$ be two three-manifolds with parameterized
  boundary $\partial Y_1=F=-\partial Y_2$, where $F$ is specified
  by the pointed matched circle $\PMC$.  Fix corresponding
  bordered Heegaard diagrams for $Y_1$ and $Y_2$. Let $Y$ be the
  closed three-manifold obtained by gluing $Y_1$ and $Y_2$ along $F$.
  Then
  $\CFa(Y)$ is homotopy equivalent to the $\Ainf$ tensor product of
  $\CFAa(Y_1)$ and $\CFDa(Y_2)$.
  In particular,
  \[\HFa(Y)\cong
H_*\left(\CFAa(Y_1)\DTP_{\Alg(\PMC)}\CFDa(Y_2)\right).\]
\end{theorem}

We give two proofs of the above theorem. One proof
(in Chapter~\ref{chap:nice-diagrams}) makes use of the
powerful technique of Sarkar and Wang~\cite{SarkarWang07:ComputingHFhat}: we construct
convenient Heegaard diagrams for $Y_1$ and $Y_2$, where the
holomorphic curve counts can be calculated combinatorially. For such
diagrams, the $\Ainf$ structure of the type $A$ module
simplifies immensely (higher multiplications all vanish), $\Ainf$
tensor products coincide with traditional tensor products, and the
proof of the pairing theorem becomes quite simple. The other proof (in
Chapter~\ref{chap:tensor-prod})
involves a rescaling argument to identify $\CFa(Y)$ with another model
for the $\Ainf$ tensor product. This second proof gives geometric
insight into the connection between the analysis and the algebra. It also leads to
pairing theorems for triangle and polygon maps \cite{LOT:DCov2, LOTCobordisms}.

\section{On gradings}

One further surprising aspect of the theory that deserves further comment
is the structure of gradings.

Gradings in the Heegaard Floer theory of closed $3$-manifolds have an
unconventional form: for each $\SpinC$ structure there is a relative
grading taking values in a cyclic abelian group. Alternatively, this
can be viewed as an absolute grading valued in some $\ZZ$-set. Indeed,
for the case of Floer homology of Seiberg-Witten
monopoles~\cite{KronheimerMrowka}, this $\ZZ$-set has an elegant
formulation as the set of isotopy classes of non-vanishing vector
fields on the $3$-manifold.

Gradings in bordered Floer homology have a correspondingly even less
conventional form. The algebra of a surface is graded by the
Heisenberg group associated to the surface's intersection form.  (This
group also has a geometric interpretation; see
Remark~\ref{rmk:geom-gradings}.)  The gradings for bordered Heegaard
Floer modules take values in sets with actions of this Heisenberg
group.  This is explained in Chapter~\ref{chap:gradings}.

\section{The case of three-manifolds with torus boundary}

As explained in Chapter~\ref{chap:TorusBoundary}, for three-manifolds
with boundary the torus~$T$, the (appropriate summand of the) algebra~$\Alg(T)$
is particularly simple: it is finite-dimensional (an eight-dimensional
subalgebra of $4 \times 4$ upper-triangular matrices) and has
vanishing differential.

The bordered Heegaard Floer homology of a solid torus bounded by $T$ can
be easily calculated; see
also \cite[Section 5.3]{Lipshitz06:BorderedHF}.
A fundamental result in Heegaard
Floer homology, the {\em surgery exact triangle}, gives a long exact
sequence relating the Heegaard Floer homology groups of three
three-manifolds which are obtained as different fillings of the same
three-manifold with torus boundary. (The result was first proved
in~\cite{OS04:HolDiskProperties}, inspired by a corresponding result
of Floer for the case of his instanton homology, see~\cite{FloerTriangles}, compare
also~\cite{KMOS} for a Seiberg-Witten analogue.) This result
can now be seen as a consequence of the pairing theorem, together
with a short exact sequence relating three type $D$ modules
associated to three different solid tori with boundary~$T$.

In a related vein, recall~\cite{OS04:Knots,Rasmussen03:Knots} that there is a
construction of Heegaard Floer homology theory for knots $K$ in a
three-manifold $Y$.  Information about this knot Floer homology turns
out to determine the bordered Floer homology of the knot complement.
This is stated precisely and worked out in Theorem~\ref{thm:HFKtoHFD}
in Section~\ref{sec:CFK-to-CFD} for large framings and in
Theorem~\ref{thm:HFKtoHFDframed} in Appendix~\ref{app:Bimodules} for
general framings, using some results from \cite{LOT2}.
This result provides a number of explicit, non-trivial examples
of bordered Floer homology.

Indeed, using this connection between knot Floer
homology and bordered Heegaard Floer homology, and combining it with
an adaptation of Theorem~\ref{thm:TensorPairing} (stated in
Theorem~\ref{thm:PairingKnot} below), we are able to calculate knot
Floer homology groups of satellite knots in terms of the filtered chain
homotopy type of the knot filtration of the pattern knot, together
with some information (a type $A$ module) associated to the pattern.
We illustrate this in some concrete examples, where the type $A$ module
can be calculated explicitly. These ideas are further pursued
in~\cite{Petkova:cableofthin,Levine12:slicingBing,Levine12:doubling}.

\section{Previous work}

Some of the material in this book, especially from
Chapters~\ref{chap:heegaard-diagrams-boundary}
and~\ref{chap:structure-moduli}, first appeared in the first
author's Ph.D. thesis~\cite{Lipshitz06:BorderedHF}.  As in that work, we construct here an
invariant of bordered manifolds that lives in a suitable category of
chain complexes. By modifying the algebra and modules introduced in
that work appropriately we
are able to reconstruct the invariant of a closed 3-manifold. Related
constructions in the case of manifolds with torus boundary have been
worked out by Eftekhary~\cite{Eftekhary08:Splicing}. For the case of
satellite knots, there is extensive work by Hedden,
including~\cite{HeddenThesis}, which we used as a check of some of our
work.  In a different direction, an invariant
for manifolds-with-boundary equipped with sutures on the boundary
has been given by Andr{\'a}s Juh{\'a}sz~\cite{Juhasz06:Sutured}.

An introduction to some of the structures used in this book, in the
form of a toy model illustrating key features, is given in a separate
paper~\cite{LOT0}.

\section{Further developments}\label{sec:further-developments}
The theory presented here is in some ways rather limited as a theory of
3-manifolds with boundary: we deal only with a single, connected
boundary component, for the somewhat limited $\HFa$ theory.
We return to these points in other papers:
\begin{itemize}
\item For a $3$-manifold $M$ with two boundary components, together
  with par\-a\-me\-tri\-za\-tions of these boundary components and a framed arc
  connecting them we associate several types of
  bimodules~\cite{LOT2}. These bimodules are well-defined up to
  homotopy equivalence, and satisfy appropriate pairing theorems. Some
  parts of this theory are sketched in Appendix~\ref{app:Bimodules}.
\item As a special case, these bimodules allow us to change the
  parametrization of the boundary. (For instance, if $F$ is a torus
  this corresponds to changing the framing on a knot.) This gives a
  representation of the mapping class group of $F$ on the derived
  category of modules over~$\Alg(\PMC)$.  This mapping class group
  action is faithful~\cite{LOT13:faith}. (Bimodules for generators of
  the mapping class group of the torus are given explicitly in
  Appendix~\ref{app:Bimodules}.)
\item The type $A$ and $D$ versions of the bordered invariant of a
  $3$-manifold with boundary are also related by certain dualizing
  bimodules~\cite{LOT2}. (Again, this is sketched in a little more
  detail in Appendix~\ref{app:Bimodules}.) They are also related by
  another kind of duality, using the
  $\Hom$-functor~\cite{LOTHomPair}. In particular, $\CFAa(Y)$ can be
  calculated explicitly from $\CFDa(Y)$, and vice-versa.
\item Decomposing a three-manifold into simple pieces, and calculating
  the (bi\discretionary{-)}{}{)}mod\-ules associated to the pieces,
  one can calculate the
  Heegaard Floer invariant $\HFa$ for any 3-manifold, as well as the
  bordered invariants of any bordered
  $3$-manifold~\cite{LOT4}. Similar techniques, together with a
  pairing theorem for polygon maps, allow one to compute the spectral
  sequence from Khovanov homology to $\HFa$ of the branched double
  cover and, more generally, the link surgery spectral sequence
  \cite{LOT:DCov1, LOT:DCov2}. These methods also allow one to compute
  the maps on $\HFa$ associated to $4$-dimensional
  cobordisms~\cite{LOTCobordisms}.
\item Andr\'as Juh\'asz has constructed an invariant of sutured
  manifolds, called sutured Floer homology \cite{Juhasz06:Sutured}, which has been remarkably
  effective at studying geometric questions
  \cite{Juhasz08:SuturedDecomp}. Rumen Zarev has developed
  a \emph{bordered-sutured Floer theory} generalizing both bordered
  Floer homology and sutured Floer homology. This theory allows one to
  recover the sutured invariants by tensoring the bordered Floer
  invariants with a module recording the sutures \cite{Zarev09:BorSut}. In turn, it also
  shows that much of bordered Floer theory can be captured in terms of
  sutured Floer homology~\cite{Zarev:JoinGlue}.
\end{itemize}

\section{Organization}

In Chapter~\ref{chap:ainfinity} we briefly recall the language of
$\Ainf$ modules. Most of the material in that chapter is standard,
with the possible exception of a particular model, denoted $\DT$, for
the $\Ainf$ tensor product. This model is applicable when one of the
two factors has a particular simple algebraic structure, which we call
a \emph{type~$D$ structure}; the modules $\CFDa$ have this form.

In Chapter~\ref{chap:algebra} we construct the differential graded
algebra associated to a closed, oriented surface.
Chapter~\ref{chap:heegaard-diagrams-boundary} collects a number of
basic properties of Heegaard diagrams for three-manifolds with
(parameterized) boundary, including when two such diagrams represent
the same $3$-manifold. The generators of the bordered Floer modules
are also introduced there, and we see how to associate
$\SpinC$-structures to these and discuss domains connecting them (when
they represent the same $\SpinC$-structure).  In
Chapter~\ref{chap:structure-moduli}, we collect the technical tools
for moduli spaces of holomorphic curves which are counted in the
algebraic structure underlying bordered Floer homology. 

With this background in place, we proceed in
Chapter~\ref{chap:type-d-mod} to the definition of the type $D$
module, and establish its invariance properties sketched in
Theorem~\ref{intro:D-invariance}.  In Chapter~\ref{chap:type-a-mod} we
treat the case of type $A$ modules, establishing a precise version of
Theorem~\ref{intro:A-invariance}.  We turn next to the pairing
theorem. Chapter~\ref{chap:nice-diagrams} gives one proof of
Theorem~\ref{thm:TensorPairing}, using so-called nice diagrams (in the
sense of \cite{SarkarWang07:ComputingHFhat}).  As a bi-product, this
gives an algorithm to calculate the bordered Floer homology of a
three-manifold.  We give a second, more analytic proof of
Theorem~\ref{thm:TensorPairing} in Chapter~\ref{chap:tensor-prod}.

While gradings on the algebra are discussed in
Chapter~\ref{chap:algebra}, discussion of gradings on the invariants
of bordered $3$-manifolds and a graded version of the pairing theorem
are deferred to Chapter~\ref{chap:gradings},

We conclude the book proper with a chapter on the torus boundary case,
Chapter~\ref{chap:TorusBoundary}. After reviewing the algebra
associated to the torus, we use a computation of the bordered Floer
homology for solid tori to give a quick proof of the surgery exact
triangle for $\HFa$. In this chapter, we also relate this invariant
with knot Floer homology. We start by showing how to recapture much of
knot Floer homology from the bordered Floer invariants of a knot
complement. We then go the other direction, showing that the bordered
Floer homology of knot complements in $S^3$ (with suitable framings)
can be computed from knot Floer homology; this is
Theorem~\ref{thm:HFKtoHFD}.

In Appendix~\ref{app:Bimodules}, we give a quick introduction
to~\cite{LOT2}, which gives a generalization of the present work to
the case of three-manifolds with two boundary components. Using an
important special case, one can describe the dependence of the
bordered invariants on the parameterization of the boundary. We give
explicit answers for this dependence in the case of torus
boundary. Using these computations, we are able to extend
Theorem~\ref{thm:HFKtoHFD} to arbitrary framings.

\section*{Acknowledgements}
We thank Dror Bar-Natan, Eaman Eftekhary, Yakov Eliashberg, Matthew
Hedden, Mikhail
Khovanov, Aaron Lauda, Dusa McDuff, Tim Perutz and Zolt{\'a}n
Szab{\'o} for helpful conversations. We are also grateful to Atanas
Atanasov, Tova
Brown, Christopher Douglas, Tom Hockenhull, Matthew
Hedden, Jennifer Hom, Yank{\i} Lekili, Adam Levine, Andy Manion, Ciprian Manolescu, Yi
Ni, Ina Petkova, Vera Vert\'esi, Chuen-Ming Mike Wong, and Chris Xiu Yang
for useful remarks on
versions of this manuscript. We also thank the valiant first referee for a careful
reading and many helpful suggestions, and the other referees for further corrections and suggestions.


\chapter{\textalt{$\Ainf$}{A-infty} structures}
\label{chap:ainfinity}

In this work, we extensively use of the notion of an
$\Ainf$ module. Although $\Ainf$ notions, first introduced by Stasheff
in the study of $H$-spaces~\cite{Stasheff}, have become commonplace
now in symplectic geometry (see for
example~\cite{Fukaya,Kontsevich,Seidel02:FukayaDef}), they might not be so familiar
to low-dimensional topologists. In Sections~\ref{sec:AinfModules}
and~\ref{sec:AinfTensorProducts}, we review the notions of primary
importance to us now, in particular sketching the $\Ainf$ tensor
product.  Keller has another pleasant
exposition of this material in \cite{AinftyAlg}.

The type~$D$ module $\CFDa(\Heegaard)$ defined in
Chapter~\ref{chap:type-d-mod} is an ordinary differential graded
module, not an $\Ainf$ module, so the reader may wish to skip this
chapter and concentrate on the type~$D$ structure at first.
Similarly, the type~$A$ module $\CFAa(\Heegaard)$ of a nice
diagram~$\Heegaard$ (in the technical sense, see Definition~\ref{def:NiceDiagram},
which in turn is motivated by~\cite{SarkarWang07:ComputingHFhat}) 
is an ordinary differential graded module.
However, the notion of grading by a non-commutative group
(Section~\ref{sec:grad-non-comm}) is used in all cases.

In Section~\ref{sec:TypeDModules},
we introduce a further algebraic structure which naturally gives rise
to a module, which we call a \emph{type $D$ structure}.  
In Section~\ref{sec:DT}, we
also give an explicit and smaller construction of the $\Ainf$
tensor product when one of the two factors is induced from a type $D$
structure. Although type $D$ structures have appeared in various guises
elsewhere (see Example~\ref{ex:FreeModule} and Remark~\ref{rmk:OtherTypeDs}), 
this tensor product is apparently new.
As the name suggests, the type~$D$ module of a bordered
three-manifold comes from a type~$D$ structure.

Finally, in Section~\ref{sec:grad-non-comm} we introduce gradings of
differential graded algebras or $\Ainf$ algebras with values in a
non-commutative group.

\section{\textalt{$\Ainf$}{A-infty} algebras and modules}
\label{sec:AinfModules}

Although we will be working over differential graded (\textit{dg}) algebras
(which are less general than $\Ainf$ algebras) we recall the
definition of $\Ainf$ algebras here as it makes a useful warm-up for
defining for
$\Ainf$ modules, which we will need in the sequel.

To avoid complicating matters at first, we will work in the category
of $\ZZ$-graded complexes over a fixed, commutative ground ring
$\gls*{Ground}$,
which we
assume to have characteristic two. (For our purposes, this ground ring $\Ground$
will typically be a direct sum of copies of $\Field=\Zmod{2}$.)
In particular, we consider modules
$M$ over $\Ground$ which are graded by the integers, so that
$$M=\bigoplus_{d\in\ZZ} M_d.$$
If $M$ is a graded module and $m\in \ZZ$, we define
$M[n]$\glsadd{shiftV}
to be the graded module defined by
$M[n]_d=M_{d-n}$.

Note that since we are working over characteristic two and don't need
to worry about signs, the following
discussion also carries over readily to the ungraded setting. The
modules in our present application are graded, but not by $\ZZ$. We return
to this point in Section~\ref{sec:grad-non-comm}.

\begin{definition}
  Fix a ground ring $\Ground$ with characteristic two.
  An \emph{$\Ainf$ algebra}
  \index{$\Ainf$!algebra|see{algebra, $\Ainf$}}\index{algebra!$\Ainf$}%
  $\gls*{cDGA}$
  over~$\Ground$
  is a graded $\Ground$-module~$\gls*{A}$,
  equipped with
  $\Ground$-linear multiplication maps
  $$\gls*{musubi}\co \DGA^{\otimes i}
  \to \DGA[2-i]$$
  defined for all $i\geq 1$, satisfying the
  compatibility conditions
  \begin{equation}
    \label{eq:AinfAlgebraAinfRelation}
    \sum_{i+j=n+1}\sum_{\ell=1}^{n-j+1}\mu_i(a_1\otimes \dots\otimes a_{\ell-1}\otimes \mu_j(a_\ell\otimes \dots\otimes a_{\ell+j-1})\otimes a_{\ell+j}\otimes \dots\otimes a_{n}) 
        =0
  \end{equation}
  for each~$n \ge 1$.
  Here, $\gls*{Atensori}$ denotes the $\Ground$-module
  $\overbrace{\DGA\otimes_{\Ground}\dots\otimes_{\Ground}\DGA}^{i}$.
  Note that throughout this chapter, all 
  tensor products are over $\Ground$ unless otherwise specified.
  We use $\cDGA$ for the $\Ainf$ algebra and $\DGA$ for its
  underlying $\Ground$-module.
  An $\Ainf$ algebra $\Alg$ is \emph{strictly unital}
  \index{strictly unital!$\Ainf$ algebra}\index{algebra!$\Ainf$!strictly unital}%
  if there is an element
  $1\in\Alg$ with the property that $\mu_2(a,1)=\mu_2(1,a)=a$ and
  $\mu_i(a_1,\dots,a_i)=0$ if $i\ne 2$ and $a_j=1$ for some $j$.
\end{definition}

In particular, an $\Ainf$ algebra is a chain complex over $\Ground$,
with differential $\mu_1$. In the case where all $\mu_i=0$ for $i>2$,
an $\Ainf$ algebra is just a differential graded algebra over
$\Ground$, with differential $\mu_1$ and (associative) multiplication
$\mu_2$. We have assumed that $\Ground$ has characteristic two; in the
more general case, the compatibility equation must be taken with signs
(see, e.g.,~\cite{AinftyAlg}).

We can think of the algebraic operations graphically as follows.  The
module $A^{\otimes i}$ is denoted by drawing $i$ parallel,
downward-oriented strands. The multiplication operation $\mu_i$ is
represented by an oriented, planar tree with one vertex, $i$ incoming strands
and one outgoing strand. With this convention, the compatibility
relation with $n$ inputs can be visualized as follows. Fix a tree $T$
with one vertex, $n$ incoming edges, and one outgoing one. Consider
next all the planar trees $S$ (with two vertices) with the property that if
we contract one edge in $S$, we obtain $T$.  Each such tree $S$ represents
a composition of multiplication maps. The compatibility condition
states that the sum of all the maps gotten by these resolutions $S$
vanishes. For example,
\[
\mathcenter{
\begin{tikzpicture}[x=.3in]
  \node at (0,0) (tc) {};
  \node at (-1,0) (tl) {};
  \node at (1,0) (tr) {};
  \node at (-.5,-1) (mu1) {$\mu_1$};
  \node at (0,-2) (mu3) {$\mu_3$};
  \node at (0,-3) (bc) {};
  \draw[->] (tc) to (mu3);
  \draw[->] (tr) to (mu3);
  \draw[->] (tl) to (mu1);
  \draw[->] (mu1) to (mu3);
  \draw[->] (mu3) to (bc);
\end{tikzpicture}}
\!\!+\!\!
\mathcenter{
\begin{tikzpicture}[x=.3in]
  \node at (0,0) (tc) {};
  \node at (-1,0) (tl) {};
  \node at (1,0) (tr) {};
  \node at (0,-1) (mu1) {$\mu_1$};
  \node at (0,-2) (mu3) {$\mu_3$};
  \node at (0,-3) (bc) {};
  \draw[->] (tc) to (mu1);
  \draw[->] (tr) to (mu3);
  \draw[->] (tl) to (mu3);
  \draw[->] (mu1) to (mu3);
  \draw[->] (mu3) to (bc);
\end{tikzpicture}}
\!\!+\!\!
\mathcenter{
\begin{tikzpicture}[x=.3in]
  \node at (0,0) (tc) {};
  \node at (-1,0) (tl) {};
  \node at (1,0) (tr) {};
  \node at (.5,-1) (mu1) {$\mu_1$};
  \node at (0,-2) (mu3) {$\mu_3$};
  \node at (0,-3) (bc) {};
  \draw[->] (tr) to (mu1);
  \draw[->] (tc) to (mu3);
  \draw[->] (tl) to (mu3);
  \draw[->] (mu1) to (mu3);
  \draw[->] (mu3) to (bc);
\end{tikzpicture}}
\!\!+\!\!
\mathcenter{
\begin{tikzpicture}[x=.3in]
  \node at (0,0) (tc) {};
  \node at (-1,0) (tl) {};
  \node at (1,0) (tr) {};
  \node at (0,-2) (mu1) {$\mu_1$};
  \node at (0,-1) (mu3) {$\mu_3$};
  \node at (0,-3) (bc) {};
  \draw[->] (tr) to (mu3);
  \draw[->] (tc) to (mu3);
  \draw[->] (tl) to (mu3);
  \draw[->] (mu3) to (mu1);
  \draw[->] (mu1) to (bc);
\end{tikzpicture}}
\!\!+\!\!
\mathcenter{
\begin{tikzpicture}[x=.3in]
  \node at (0,0) (tc) {};
  \node at (-1,0) (tl) {};
  \node at (1,0) (tr) {};
  \node at (0,-1) (mu2a) {$\mu_2$};
  \node at (0,-2) (mu2b) {$\mu_2$};
  \node at (0,-3) (bc) {};
  \draw[->] (tl) to (mu2a);
  \draw[->] (tc) to (mu2a);
  \draw[->] (tr) to (mu2b);
  \draw[->] (mu2a) to (mu2b);
  \draw[->] (mu2b) to (bc);
\end{tikzpicture}}
\!\!+\!\!
\mathcenter{
\begin{tikzpicture}[x=.3in]
  \node at (0,0) (tc) {};
  \node at (-1,0) (tl) {};
  \node at (1,0) (tr) {};
  \node at (0,-1) (mu2a) {$\mu_2$};
  \node at (0,-2) (mu2b) {$\mu_2$};
  \node at (0,-3) (bc) {};
  \draw[->] (tr) to (mu2a);
  \draw[->] (tc) to (mu2a);
  \draw[->] (tl) to (mu2b);
  \draw[->] (mu2a) to (mu2b);
  \draw[->] (mu2b) to (bc);
\end{tikzpicture}}
\!\!=0.
\]
Note that the last two terms are the usual associativity relation, and
the remaining terms say that associativity holds only up to a homotopy.

Another way to think of the compatibility conditions (Equation~\eqref{eq:AinfAlgebraAinfRelation}) uses
the tensor algebra
\[\Tensor^*(\DGA[1])\coloneqq\bigoplus_{n=0}^{\infty}A^{\otimes
  n}[n].\]
\glsadd{TensorA}%
The maps~$\mu_i$ can be combined into a single map
\[\mu\co\Tensor^*(\DGA[1]) \to \DGA[2]\]
with the convention that $\mu_0 = 0$.  We can also construct
a natural degree~1 endomorphism $\gls*{barD}\co \Tensor^*(\DGA[1])\to \Tensor^*(\DGA[1])$ by
\begin{equation}
\label{eq:DefDBar}
{\widebar D}(a_1\otimes\dots\otimes a_n)
= \sum_{j=1}^n \sum_{\ell=1}^{n-j+1} a_1\otimes\dots
\otimes \mu_j(a_\ell\otimes\dots\otimes a_{\ell+j-1})\otimes\dots\otimes a_n.
\end{equation}
Then the compatibility condition is that
\begin{align*}
  \mu \circ \widebar{D} &= 0\\
  \shortintertext{or equivalently}
  \widebar{D} \circ \widebar{D} &= 0.
\end{align*}
These conditions can also be expressed graphically as
\[
\mathcenter{
\begin{tikzpicture}
  \node at (0,0) (tc) {};
  \node at (0,-1) (D) {$\overline{D}$};
  \node at (0,-2) (mu) {$\mu$};
  \node at (0,-3) (bc) {};
  \draw[taar] (tc) to (D);
  \draw[taar] (D) to (mu);
  \draw[aar] (mu) to (bc);
\end{tikzpicture}}
\qquad\mathcenter{\text{or}}\qquad
\mathcenter{
\begin{tikzpicture}
  \node at (0,0) (tc) {};
  \node at (0,-1) (D) {$\overline{D}$};
  \node at (0,-2) (D2) {$\overline{D}$};
  \node at (0,-3) (bc) {};
  \draw[taar] (tc) to (D);
  \draw[taar] (D) to (D2);
  \draw[taar] (D2) to (bc);
\end{tikzpicture}}.
\]
Here, the doubled arrows indicate elements of $\Tensor^*(\DGA[1])$,
while the single arrows indicate elements of $\DGA$.

At some points later we will want to allow only finitely many nonzero
$\mu_i$. We give this condition a name:
\begin{definition}
  An $\Ainf$ algebra $(A,\{\mu_i\}_{i=1}^\infty)$ is \emph{operationally bounded}
  \index{operationally bounded $\Ainf$ algebra}\index{algebra!$\Ainf$!operationally bounded}%
  \index{bounded!operationally b.~$\Ainf$ algebra}%
  if $\mu_i=0$ for $i$ sufficiently large.
\end{definition}

\begin{definition}\label{def:ainf-module}
  A \emph{(right) $\Ainf$ module}
  \index{$\Ainf$!module|see{module, $\Ainf$}}%
  \index{module!$\Ainf$}%
  $\gls*{cModule}$
  over $\cDGA$ is a graded
  $\Ground$-module $\gls*{M}$,
  equipped with operations
  $$\gls*{msubi}\co M\otimes \DGA^{\otimes (i-1)}
  \to M[2-i],$$
  defined for all $i \geq 1$, satisfying the compatibility conditions
  \begin{equation}\label{eq:Ainf-mod-compat}
  \begin{split}
    0 &= \sum_{i+j=n+1}\!m_i(m_j(\x\otimes a_1\otimes \dots\otimes a_{j-1})
        \otimes \dots\otimes a_{n-1})\\
    &\quad+\!\sum_{i+j=n+1}\sum_{\ell=1}^{n-j}m_i(\x\otimes a_1\otimes \dots\otimes a_{\ell-1}\otimes \mu_j(a_\ell\otimes \dots\otimes a_{\ell+j-1})\otimes \dots\otimes a_{n-1}).
  \end{split}
  \end{equation}
  An $\Ainf$ module $\cModule$ over a strictly unital $\Ainf$ algebra
  is said to be {\em strictly unital}
  \index{strictly unital!$\Ainf$ module}\index{module!$\Ainf$!strictly unital}%
  if for any $\x\in M$,
  $m_2(\x\otimes 1)=\x$ and
  $m_i(\x\otimes a_1\otimes \dots\otimes a_{i-1}) = 0$
  if $i> 2$ and some $a_j=1$.  The module $\cModule$ is said to be
  {\em bounded}
  \index{bounded!$\Ainf$ module}\index{module!$\Ainf$!bounded}%
  if $m_i=0$ for all sufficiently large $i$.
\end{definition}

Graphically, the module $M\otimes \DGA^{\otimes i-1}$ is represented by
$i$ parallel
strands, where the leftmost strand is colored with $M$ (while the others
are colored by $\DGA$). The multiplication map $m_{i}$ is represented
by the one-vertex tree with $i$ incoming strands, the leftmost of which
is $M$-colored, and whose output is also $M$-colored. The compatibility
condition can be thought of the same as before, except
now we have distinguished the leftmost strands with~$M$.
Thus, the compatibility conditions can be thought of as stating the
vanishing of maps induced by sums of planar trees, with a fixed number of
inputs (and one output).  

As before, we can also write this condition in terms of the tensor
algebra.
The maps $m_i$ (and the $\mu_i$) can be combined
to form a degree $1$ map 
$\gls*{barm}\co M\otimes \Tensor^*(\DGA[1])\to M\otimes
\Tensor^*(\DGA[1])$,
using the following analogue of Equation~\eqref{eq:DefDBar}:
\begin{align}
\label{def:DefOfMBar}
{\widebar m}(\x_1\otimes a_2 \otimes\dots\otimes a_n)
&= \sum_{\ell=1}^n
m_\ell(\x\otimes a_2\otimes\dots\otimes a_{\ell})\otimes \dots\otimes a_n
\nonumber \\
&+
\sum_{j=1}^n \sum_{\ell=1}^{n-j+1} \x\otimes\dots
\otimes \mu_j(a_\ell\otimes\dots\otimes a_{\ell+j-1})\otimes\dots\otimes a_n.
\end{align}
The $\Ainf$ relation is equivalent to the condition that
${\overline m}\circ{\overline m} = 0$.

The $\Ainf$ relation in tensor form can also be  expressed
graphically. To do this
conveniently, let $\Delta\co
\Tensor^*\DGA\to\Tensor^*\DGA\otimes\Tensor^*\DGA$ be the canonical
comultiplication 
\[
\Delta(a_1\otimes\cdots\otimes
a_n)=\sum_{m=0}^{n}(a_1\otimes\cdots\otimes
a_m)\otimes(a_{m+1}\otimes\cdots\otimes a_n).
\]
If we use a dashed line to denote an element of $\cModule$
then the compatibility condition for an $\Ainf$ module can be drawn as
\[
\mathcenter{
\begin{tikzpicture}
  \node at (0,0) (tc) {};
  \node at (2,0) (tr) {};
  \node at (1,-1) (Delta) {$\Delta$};
  \node at (0,-2) (ma) {$m$};
  \node at (0,-3) (mb) {$m$};
  \node at (0,-4) (bc) {};
  \draw[taar] (tr) to (Delta);
  \draw[taar] (Delta) to (ma);
  \draw[taar] (Delta) to (mb);
  \draw[amar] (tc) to (ma);
  \draw[amar] (ma) to (mb);
  \draw[amar] (mb) to (bc);
\end{tikzpicture}}
\mathcenter{+}
\mathcenter{
\begin{tikzpicture}
  \node at (0,0) (tc) {};
  \node at (2,0) (tr) {};
  \node at (1,-1) (D) {$\overline{D}$};
  \node at (0,-2) (m) {$m$};
  \node at (0,-3) (bc) {};
  \draw[taar] (tr) to (D);
  \draw[taar] (D) to (m);
  \draw[amar] (tc) to (m);
  \draw[amar] (m) to (bc);
\end{tikzpicture}
}
\mathcenter{\quad=\quad 0}.
\]

\begin{example}
  Let $\cModule$ be a right $\Ainf$ module over $\cDGA$. Suppose
  moreover that, for all $i>2$, $m_i=0$ and $\mu_i=0$. Then not only
  is $\cDGA$ a differential graded algebra, but also $\cModule$ is a
  module over $\cDGA$ (in the traditional sense), equipped with a
  differential satisfying the Leibniz rule with respect to the algebra
  action on the module, i.e., a \emph{differential graded module}.
  \index{differential graded module}%
  \index{module!differential graded}%
  \index{dg module}%
\end{example}

\begin{definition}\label{def:AinfMorphism}
  Let $\cModule$ and $\cModule'$ be strictly unital
  right $\Ainf$ modules over a strictly unital $\Ainf$ algebra. A {\em
    strictly unital homomorphism $f$ of $\Ainf$ modules},
  \index{strictly unital!homomorphism of $\Ainf$ modules}%
  \index{homomorphism!of $\Ainf$ modules}%
  \index{$\Ainf$!homomorphism|see{homomorphism of $\Ainf$ modules}}%
  or simply an
  $\Ainf$ homomorphism, is a collection of maps
  $$f_i\co M\otimes \DGA^{\otimes (i-1)}
  \to M'[1-i]$$
  indexed by $i\geq 1$,
  satisfying the  compatibility conditions
  \begin{align*}
    0 &=\sum_{i+j=n+1}\!m_i'(f_j(\x\otimes a_1\otimes \dots\otimes a_{j-1})\otimes \dots\otimes a_{n-1})\\
    &\quad+\!\sum_{i+j=n+1}\!f_i(m_j(\x\otimes a_1\otimes \dots\otimes a_{j-1})\otimes \dots\otimes a_{n-1}) \\
    &\quad+\!\sum_{i+j=n+1}\sum_{\ell=1}^{n-j}f_i(\x\otimes a_1\otimes \dots\otimes a_{\ell-1}\otimes \mu_j(a_\ell\otimes \dots\otimes a_{\ell+j-1})\otimes \dots\otimes a_{n-1})
  \end{align*}
  and the unital condition
  $$f_i(\x\otimes a_1\otimes\dots\otimes a_{i-1})=0$$ if
  $i>1$ and some $a_j=1$.

  We call a strictly unital homomorphism $\gls*{fsubi}$
  of $\Ainf$ modules
  \emph{bounded} if $f_i=0$ for $i$ sufficiently large.
  \index{bounded!homomorphism of $\Ainf$ modules}%
\end{definition}
  
The compatibility condition can be formulated in terms of 
the tensor algebra as follows.
Promote the maps $\{f_i\}_{i=1}^{\infty}$ to a degree $0$ map 
$${\overline f}\co M\otimes\Tensor^*(\DGA[1])\to M'\otimes \Tensor^*(\DGA[1])$$
by the formula:
\begin{equation}
{\widebar f}(\x_1\otimes a_1 \otimes\dots\otimes a_n)
= \sum_{\ell=0}^n
f_{\ell+1}(\x\otimes a_1\otimes\dots\otimes a_{\ell})\otimes \dots\otimes a_n.
\label{eq:DefOfFBar} 
\end{equation}
The $\Ainf$ relation now is equivalent to the condition that the map
${\overline f}$ induces a chain map from the chain complex
$M\otimes\Tensor^*(\DGA[1])$ with its induced differential ${\overline
  m}$ (Equation~\eqref{def:DefOfMBar}) to the chain complex $M'\otimes
\Tensor^*(\DGA[1])$ with its induced differential ${\overline m}'$.

The compatibility condition for $\Ainf$ homomorphisms can be
drawn as
\[
\mathcenter{
\begin{tikzpicture}
  \node at (0,0) (tc) {};
  \node at (2,0) (tr) {};
  \node at (1,-1) (Delta) {$\Delta$};
  \node at (0,-2) (ma) {$m$};
  \node at (0,-3) (mb) {$f$};
  \node at (0,-4) (bc) {};
  \draw[taar] (tr) to (Delta);
  \draw[taar] (Delta) to (ma);
  \draw[taar] (Delta) to (mb);
  \draw[amar] (tc) to (ma);
  \draw[amar] (ma) to (mb);
  \draw[amarp] (mb) to (bc);
\end{tikzpicture}}
\mathcenter{+}
\mathcenter{
\begin{tikzpicture}
  \node at (0,0) (tc) {};
  \node at (2,0) (tr) {};
  \node at (1,-1) (Delta) {$\Delta$};
  \node at (0,-2) (ma) {$f$};
  \node at (0,-3) (mb) {$m'$};
  \node at (0,-4) (bc) {};
  \draw[taar] (tr) to (Delta);
  \draw[taar] (Delta) to (ma);
  \draw[taar] (Delta) to (mb);
  \draw[amar] (tc) to (ma);
  \draw[amarp] (ma) to (mb);
  \draw[amarp] (mb) to (bc);
\end{tikzpicture}}
\mathcenter{+}
\mathcenter{
\begin{tikzpicture}
  \node at (0,0) (tc) {};
  \node at (2,0) (tr) {};
  \node at (1,-1) (D) {$\overline{D}$};
  \node at (0,-2) (m) {$f$};
  \node at (0,-3) (bc) {};
  \draw[taar] (tr) to (D);
  \draw[taar] (D) to (m);
  \draw[amar] (tc) to (m);
  \draw[amarp] (m) to (bc);
\end{tikzpicture}
}
\mathcenter{=0}.
\]
Here, dashed lines represent $\cM$ and dotted lines represent~$\cModule'$.

For example, for any $\Ainf$ module $\cModule$, the \emph{identity
  homomorphism}~$\gls*{Id}$
\index{homomorphism!of $\Ainf$ modules!identity}%
is the map with
\begin{align*}
  \Id_1(\x)&\coloneqq\x \\
  \Id_i(\x\otimes A^{\otimes i-1})&\coloneqq 0\quad (i > 0).
\end{align*}

If $f$ is an $\Ainf$ homomorphism from $\cModule$ to $\cModule'$,
and $g$ is an $\Ainf$ homomorphism from $\cModule'$ to
$\cModule''$, we can form their {\em composite}
\index{composite of $\Ainf$ homomorphisms}%
\index{homomorphism!of $\Ainf$ modules!composition}%
$g\circ f$,
defined by
  $$(g\circ f)_n(\x\otimes (a_1\otimes\dots\otimes a_{n-1}))
  \coloneqq\!\sum_{i+j=n+1}\! g_j (f_i(\x\otimes a_1\otimes\dots\otimes
  a_{i-1})\otimes\dots\otimes a_{n-1}).$$

\begin{definition}
  Given any collection of maps 
  $$h_i \co M\otimes \DGA^{\otimes i-1} \to
  M'[-i]$$ 
  with $h_i(\x\otimes a_1\otimes\dots \otimes a_{i-1})=0$ if $i > 1$ and some $a_j=1$,
  we can construct a strictly unital $\Ainf$ homomorphism~$f$ by the expression
  \begin{multline*}
    f_n(\x\otimes a_1 \otimes\dots\otimes a_{n-1})=\\
  \begin{aligned}
    &\sum_{i+j=n+1}\!h_i(m_j(\x\otimes a_1\otimes \dots\otimes a_{j-1})\otimes \dots\otimes a_{n-1})
\\
    &+\!\sum_{i+j=n+1}\!m_i'(h_j(\x\otimes a_1\otimes \dots\otimes a_{j-1})\otimes \dots\otimes a_{n-1})\\
    &+\!\sum_{i+j=n+1}\sum_{\ell=1}^{n-j}h_i(\x\otimes a_1\otimes \dots\otimes a_{\ell-1}\otimes \mu_j(a_\ell\otimes \dots\otimes a_{\ell+j-1})\otimes \dots\otimes a_{n-1}).
   \end{aligned}
  \end{multline*}
   
   If an $\Ainf$ homomorphism~$f$ can be obtained from a
   map~$h$ in this way, we call $f$ \emph{null homotopic}.
   \index{null homotopic|see{homomorphism, null homotopic}}%
   \index{homomorphism!of $\Ainf$ modules!null homotopic}%
   Two
   $\Ainf$ homomorphisms 
   $f,g\co \cModule \to \cModule'$ are {\em
     homotopic} if their difference is null homotopic.
   \index{homotopic|see{homomorphism, homotopic}}
   \index{homomorphism!of $\Ainf$ modules!homotopic}%

   Two $\Ainf$ modules $\cModule$ and $\cModule'$ are ($\Ainf$)
   \emph{homotopy
     equivalent}
   \index{homotopy equivalent|see{module, homotopy equivalent}}
   \index{$\Ainf$!module!homotopy equivalent}%
   if there are $\Ainf$ homomorphisms $f$ from
   $\cModule$ to $\cModule'$ and $g$ from $\cModule'$ to $\cModule$
   such that $f\circ g$ and $g\circ f$ are homotopic to the identity.
\end{definition}

Suppose that $f,g\co \cModule\to\cModule'$ are homotopic via a homotopy
$h$. Then, we can promote the components
of $h$, $\{h_i\}_{i=1}^{\infty}$ to a map 
$${\overline h}\co M\otimes\Tensor^*(\DGA[1])\to M'\otimes \Tensor^*(\DGA[1])$$
following Equation~\eqref{eq:DefOfFBar}.
The condition that $h$ gives the homotopy between $f$ and $g$ is equivalent
to the more familiar condition
$${\overline h}\circ {\overline m}+{\overline m}'\circ {\overline h}
= {\overline f}+{\overline g}.$$
Graphically, this equation is
\[
\mathcenter{
\begin{tikzpicture}
  \node at (0,0) (tc) {};
  \node at (2,0) (tr) {};
  \node at (1,-1) (Delta) {$\Delta$};
  \node at (0,-2) (ma) {$m$};
  \node at (0,-3) (mb) {$h$};
  \node at (0,-4) (bc) {};
  \draw[taar] (tr) to (Delta);
  \draw[taar] (Delta) to (ma);
  \draw[taar] (Delta) to (mb);
  \draw[amar] (tc) to (ma);
  \draw[amar] (ma) to (mb);
  \draw[amarp] (mb) to (bc);
\end{tikzpicture}}
\mathcenter{+}
\mathcenter{
\begin{tikzpicture}
  \node at (0,0) (tc) {};
  \node at (2,0) (tr) {};
  \node at (1,-1) (Delta) {$\Delta$};
  \node at (0,-2) (ma) {$h$};
  \node at (0,-3) (mb) {$m'$};
  \node at (0,-4) (bc) {};
  \draw[taar] (tr) to (Delta);
  \draw[taar] (Delta) to (ma);
  \draw[taar] (Delta) to (mb);
  \draw[amar] (tc) to (ma);
  \draw[amarp] (ma) to (mb);
  \draw[amarp] (mb) to (bc);
\end{tikzpicture}}
\mathcenter{+}
\mathcenter{
\begin{tikzpicture}
  \node at (0,0) (tc) {};
  \node at (2,0) (tr) {};
  \node at (1,-1) (D) {$\overline{D}$};
  \node at (0,-2) (m) {$h$};
  \node at (0,-3) (bc) {};
  \draw[taar] (tr) to (D);
  \draw[taar] (D) to (m);
  \draw[amar] (tc) to (m);
  \draw[amarp] (m) to (bc);
\end{tikzpicture}
}
\mathcenter{\quad=\quad}
\mathcenter{
  \begin{tikzpicture}
    \node at (0,0) (tc) {};
    \node at (1,0) (tr) {};
    \node at (0,-1) (f) {$f$};
    \node at (0,-2) (bc) {};
    \draw[taar] (tr) to (f);
    \draw[amar] (tc) to (f);
    \draw[amarp] (f) to (bc);
  \end{tikzpicture}
}
\mathcenter{+}
\mathcenter{
  \begin{tikzpicture}
    \node at (0,0) (tc) {};
    \node at (1,0) (tr) {};
    \node at (0,-1) (f) {$g$};
    \node at (0,-2) (bc) {};
    \draw[taar] (tr) to (f);
    \draw[amar] (tc) to (f);
    \draw[amarp] (f) to (bc);
  \end{tikzpicture}
}.
\]

There are analogous definitions for left $\Ainf$ modules, maps, and
homotopies.

\section{\textalt{$\Ainf$}{A-infty} tensor products}
\label{sec:AinfTensorProducts}

Consider the category of strictly unital $\Ainf$ modules over
a \dg algebra $\DGA$, whose morphism sets are the strictly unital
$\Ainf$ homomorphisms. In this
category, the usual na\"{\i}ve tensor product $M \otimes_A N$ of a
right ($\Ainf$) module~$M$ and a left module~$N$ is not available: the
relation
$(m\cdot a)\otimes n\sim m\otimes (a\cdot n)$
used in the definition of the tensor product is not
transitive in the $\Ainf$ setting. In fact, even in the original
setting of differential graded modules, this tensor product is often not the right one: if $M$ and
$M'$ are quasi-isomorphic, then $M \otimes_A N$ and $M' \otimes_A N$
are not necessarily quasi-isomorphic.

To produce a
functor which works for $\Ainf$ modules,
we pass instead to the
$\Ainf$ tensor product, a generalization of the derived tensor product.

\begin{definition}
  Let $\cDGA$ be an $\Ainf$ algebra over~$\Ground$,
  $\cModule$ be a right $\Ainf$ module over~$\cDGA$
  and $\cNodule$ be a left $\Ainf$ module over~$\cDGA$.  Then their
  \emph{$\Ainf$ tensor product}
  \index{$\Ainf$!tensor product|see{tensor product, $\Ainf$}}%
  \index{tensor product!$\Ainf$}%
  is the
  chain complex
  \[
  \cModule\gls*{DTP}_\cDGA \cNodule \coloneqq
    M \otimes \Tensor^*(\DGA[1]) \otimes N
  \]
  equipped with the boundary operator
  \begin{multline*}
  \partial (\x \otimes a_1 \otimes \cdots \otimes a_n \otimes \y) \coloneqq\\
  \begin{aligned}
    &\sum_{i=1}^{n+1} m_i(\x\otimes a_1\otimes\cdots\otimes a_{i-1})
       \otimes\cdots\otimes a_n\otimes \y\\
    &+\sum_{i=1}^n \sum_{\ell=1}^{n-i+1}
       \x\otimes a_1 \otimes\cdots\otimes\mu_i(a_\ell\otimes\cdots\otimes a_{\ell+i-1})\otimes\cdots\otimes a_n\otimes\y\\
    &+\sum_{i=1}^{n+1} \x \otimes a_1 \otimes \cdots \otimes m_i(a_{n-i+2}\otimes\cdots\otimes a_n\otimes\y).
  \end{aligned}
  \end{multline*}
\end{definition}

\begin{lemma}
  The operator~$\bdy$ on $\cModule \DTP \cNodule$ defined above
  satisfies $\bdy^2 = 0$.
\end{lemma}

\begin{proof}
  This is a straightforward consequence of the relations defining
  $\Ainf$ algebras and modules.
\end{proof}

An important case is when $\cNodule = \cDGA$.  The resulting complex
is called the \emph{bar resolution}~$\gls*{barModule}$
\index{bar resolution}%
of~$\cModule$.  As a $\Ground$-module it is
\[
\gls*{barM} \cong M \otimes \Tensor^*(\DGA[1]) \otimes \DGA
\cong (M \otimes\Tensor^+(\DGA[1]))[-1]
\]
where
\glsadd{TensorPlusA}
$\Tensor^+(\DGA[1]) \coloneqq \bigoplus_{i=1}^\infty
(\DGA[1])^{\otimes i}$
is $\Tensor^*(\DGA[1])$ without the first summand of~$\Ground$.  $\widebar{\cM}$
has the further structure of a right $\Ainf$ module
over~$\cDGA$, with multiplications $m_1 = \partial$ as defined above
and, for $i \ge 2$,
\begin{multline*}m_i((\x\otimes a_1\otimes\dots\otimes a_n)\otimes b_1\otimes\dots \otimes b_{i-1})\coloneqq \\
\sum_{\ell=1}^{n}
\x\otimes a_1\otimes \dots\otimes a_{n-\ell}
\otimes \mu_{i+\ell-1}(a_{n-\ell+1}\otimes\dots \otimes a_n\otimes
b_1\otimes\dots\otimes b_{i-1}).
\end{multline*}
Note the range of summation: the multiplication $\mu_{i+\ell-1}$ is
applied to all of the $b_j$'s and at least one~$a_j$.

The following lemma is again straightforward and standard.

\begin{lemma}\label{lem:bar-resolution-module}
  ${\widebar \cModule}$ is an $\Ainf$ module over~$\Alg$.  
\end{lemma}

Notice that if $\cDGA$ is an honest differential graded algebra then
${\widebar \cModule}$ is an honest differential graded module, i.e.,
all of the higher $m_i$ vanish.  Furthermore, if $\cM$ is also an honest
differential graded module, this is the usual bar resolution.  If
in addition $\cN$ is an honest differential graded module, then we could
also define $\cM \DTP \cN$ as $\widebar{\cM} \otimes_\DGA \cN$,
the na\"{\i}ve tensor product with the bar resolution,
but this is not possible in general.

\begin{remark}
  One way to understand the products on $\widebar{\cM}$ and
  Lemma~\ref{lem:bar-resolution-module} is via $\Ainf$ bimodules.  If
  $\cA$ and $\cB$ are two $\Ainf$ algebras, an $\cA$--$\cB$
  bimodule~$\cN$ is a $\Ground$-module~$N$ with maps
  \[
  m_{i,j} \co A^{\otimes i} \otimes N \otimes B^{\otimes j}
    \to N
  \]
  satisfying a natural version of the compatibility conditions.  In
  particular, $\cA$ is an $\cA$--$\cA$ bimodule in a natural way, and
  $\widebar{\cModule}$ is just $\cModule\DTP_\cA\cA$. A
  more general version of Lemma~\ref{lem:bar-resolution-module} is
  that if $\cM$ is a $\cA$ module and $\cN$ is an $\cA$--$\cB$
  bimodule, then $\cM \otimes_\cA \cN$
  is a right $\cB$ module.
\end{remark}

\begin{proposition}
  \label{prop:BarResolution}
  If $\cDGA$ and $\cModule$ are strictly
  unital, then ${\widebar \cModule}$ is
  homotopy equivalent to $\cModule$.
\end{proposition}

\begin{proof}
  We define an $\Ainf$ module map~$\phi$ from ${\widebar \cModule}$ to
  $\cModule$ by
  \[
  \phi_i((\x\otimes a_1\otimes\cdots\otimes a_n) \otimes b_1\otimes\cdots\otimes b_{i-1})
  \coloneqq
      m_{i+n}(\x\otimes a_1\otimes\cdots\otimes a_n \otimes b_1 \otimes\cdots\otimes b_{i-1}).
  \]
  Similarly, we define an $\Ainf$ module map~$\psi$ from $\cM$ to
  ${\widebar \cM}$ by
  \[
    \psi_i(\x\otimes a_1\otimes \dots \otimes a_{i-1})
    \coloneqq \x\otimes a_1\otimes\dots\otimes a_{i-1}\otimes 1.
  \]
  Also define maps for a homotopy~$h$ from $\widebar\cModule$ to $\widebar\cModule$ by
  \[
  h_i((\x\otimes a_1\otimes \dots\otimes a_n) \otimes b_1\otimes\cdots\otimes b_{i-1})
  \coloneqq
      \x\otimes a_1\otimes\dots\otimes a_n\otimes
      b_1\otimes\dots\otimes b_{i-1} \otimes 1.
  \]
  It is straightforward to verify that $\phi\circ \psi$ is the identity map,
  while $\psi\circ \phi$ is homotopic, via $h$, to the identity map.
\end{proof}

\begin{proposition}
  An $\Ainf$ map $f\co\cM \to \cM'$ induces a chain map $\gls*{DTPfId}
  \co \cM\DTP\cN \to \cM'\DTP\cN$.
  If $f$ is null homotopic, so is $f\DTP\Id_{\cN}$.
  In particular, if $\cModule_1$ and $\cModule_2$ are homotopy equivalent right
  $\Ainf$ modules and $\cNodule_1$ and $\cNodule_2$ are homotopy
  equivalent left $\Ainf$ modules, then $\cModule_1\DTP \cNodule_1$ is
  homotopy equivalent to $\cModule_2\DTP \cNodule_2$.
\end{proposition}

\begin{proof}
  For $f$ a collection of maps
  $f_i\co M\otimes A^{\otimes (i-1)} \to M[c-i]$ for a fixed $c\in\{0,1\}$,
  define $f \DTP \Id_{\cN}$ by
  \[
  (f \DTP \Id_{\cN})(\x\otimes a_1\otimes\dots\otimes a_n\otimes\y) \coloneqq
    \sum_{i=1}^{n+1}f_i(\x\otimes a_1\dots\otimes a_{i-1})
      \otimes a_i\otimes\dots\otimes a_n\otimes\y.
  \]
  If $f$ is an $\Ainf$ map, so is $f
  \DTP \Id_{\cN}$, and if $h$ is a null homotopy of~$f$,
  then $h \DTP \Id_{\cN}$ is a null homotopy of $f \DTP \Id_{\cN}$.
  The same construction works, \emph{mutatis mutandis}, for chain maps
  and homotopies on the right.
\end{proof}

\section{Type \textalt{$D$}{D} structures}
\label{sec:TypeDModules}

There is a smaller model for the derived tensor product in the case
where one of the two factors has a special type: where it is the
$\Ainf$ module associated to a type $D$ structure
(Definitions~\ref{def:TypeD}). In this section, we introduce these
type $D$ structures and develop their basic properties.

We first specialize to the case where $\cDGA$ is a differential graded
algebra, which is all that we need for our present purposes.

\begin{definition}
  \label{def:TypeD}
  Fix a \dg algebra $\cDGA$.
  Let $N$ be a graded $\Ground$-module, equipped with a map
  $$\gls*{ddiff}_N\co N \to (\DGA\otimes N)[1],$$
  satisfying the compatibility condition that
   $$(\mu_2\otimes
   \Id_N)\circ(\Id_A\otimes\delta^1)\circ\delta^1+(\mu_1\otimes
   \Id_N)\circ\delta^1\co N\to\DGA\otimes N$$
   vanishes. (When the module $N$ is fixed, we drop the subscript $N$ from
   the notation for $\delta^1$.) 
  We call the pair $(N,\delta^1_N)$ a {\em type $D$
    structure over $A$ with base ring $\Ground$}.
  \index{module!type $D$}%
  \index{type $D$!module|see{module, type $D$}}%
  \glsadd{typedstr}%
  Let 
  $(N_1,\delta^1_{N_1})$ and $(N_2,\delta^1_{N_1})$ be two type $D$ structures.
  A $\Ground$-module map $\psi^1\co N_1\to \DGA\otimes N_2$
  is a {\em $D$-structure homomorphism}
  \index{homomorphism!of type $D$ modules}%
  \index{type $D$!homomorphism|see{homomorphism, type $D$}}%
  if 
  $$(\mu_2\otimes \Id_{N_2})\circ (\Id_{\DGA}\otimes \psi^1)\circ \delta^1_{N_1}
  +(\mu_2\otimes \Id_{N_2})\circ (\Id_{\DGA}\otimes \delta^1_{N_2})\circ \psi^1
  +(\mu_1\otimes \Id_{N_2})\circ\psi^1=0.$$
  Given type $D$ homomorphisms
  $\phi^1\co N_1\to \DGA\otimes N_2$ and $\psi^1\co N_2\to\DGA\otimes N_3$,
  their {\em composite} $(\psi^1\circ \phi^1)$ is defined by
  $(\mu_2\otimes \Id_{N_3})\circ (\Id_\DGA\otimes \psi^1)\circ \phi^1$.
  A {\em homotopy}
  \index{homomorphism!of type $D$ modules!homotopy of}%
  $h$ between
  two $D$-structure homomorphisms
  $\psi^1$ and $\phi^1$ from 
  $(N_1,\delta^1_{N_1})$ to $(N_2,\delta^1_{N_2})$ is a $\Ground$-module
  map $h\co N_1\to (A \otimes N_2)[-1]$
  with 
  $$(\mu_2\otimes \Id_{N_2})\circ (\Id_{\DGA}\otimes h)\circ \delta^1_{N_1}
+(\mu_2\otimes \Id_{N_2})\circ (\Id_{\DGA}\otimes \delta^1_{N_2})\circ h
+(\mu_1\otimes \Id_{N_2})\circ h=\psi^1-\phi^1.$$
\end{definition}

\begin{example}
  \label{ex:FreeModule}
  Assume $\DGA$ is a \dg algebra over
  $\Ground=\Field$. Let $M$ be a \dg module which is free as an $\DGA$-module.
  For a fixed basis of $M$ over $\DGA$, let $X$ denote the
  $\Ground$-span of that basis.
  The restriction of the boundary operator to $X$ gives a map
  $$ \delta^1\co X \to \Alg\otimes X=M $$
  which satisfies the relations of a type $D$ structure.  Maps between
  such modules are $D$-structure homomorphisms.
\end{example}

\begin{lemma}
  \label{lemma:AssociatedTypeD}
  If $(N,\delta^1)$ is a type $D$ structure, then 
  $\DGA\otimes N$ can be given the structure of a left $\DGA$ module,
  with 
  \begin{align*}
  m_1(a \otimes \y)&\coloneqq[(\mu_2\otimes\Id_{N})\circ(\Id_{\DGA}\otimes\delta^1)
  + \mu_1\otimes\Id_{N}](a \otimes \y) \\
  m_2(a_1 \otimes (a \otimes \y))&\coloneqq \mu_2(a_1\otimes a) \otimes \y.
  \end{align*}
  Moreover, if $\psi^1\co N_1 \to A\otimes N_2$ is a $D$-structure homomorphism,
  then there is an induced map of differential graded modules from $A\otimes N_1$ to $A\otimes N_2$,
  defined by
  \[(a \otimes \y)\mapsto
    (m_2 \otimes \Id_{N_2})\circ(\Id_A \otimes \psi^1) \]
  Similarly, a homotopy between two type $D$-structure homomorphisms induces
  a chain homotopy between the associated chain maps.
\end{lemma}

\begin{proof}
  The proof is straightforward. 
\end{proof}

Given a type $D$ structure $(N,\delta^1)$ over a \dg algebra $\DGA$, we denote the $\DGA$-module
from Lemma~\ref{lemma:AssociatedTypeD} by~$\gls*{cNodule}$.

Note that the structure equation for a type $D$
structure $(N,\delta^1)$ is equivalent to the condition that $m_1$ on 
$\cNodule$ is, indeed, a differential.

As a partial converse to Lemma~\ref{lemma:AssociatedTypeD}, we have the
following lemma. Recall that if $N_1$ and $N_2$ are two \dg modules
over a \dg algebra, then a \dg homomorphism is a chain map $\psi\co
N_1\to N_2$ which commutes with the algebra action. Equivalently, it
is an $\Ainf$ homomorphism $\{\psi_i\}$ whose components $\phi_i$
vanish for all $i>1$.

\begin{lemma}
  \label{lemma:ConverseAssociatedTypeD}
  Let $N_1$ and $N_2$ be two type $D$ structures over a \dg algebra
  $\DGA$, and $\cNodule_1$ and
  $\cNodule_2$ be their associated \dg modules.  The correspondence
  from Lemma~\ref{lemma:AssociatedTypeD} gives an isomorphism between
  the space of type $D$ homomorphisms from $N_1$ to $N_2$ with the
  space of \dg homomorphisms from $\cNodule_1$ to
  $\cNodule_2$. Moreover, two type $D$ homomorphisms are 
  homotopic if and only if the corresponding \dg homomorphisms
  are homotopic by an $\DGA$-equivariant homotopy.
\end{lemma}

\begin{proof}
  Given $\psi\co \cNodule_1\to\cNodule_2$, we can realize
  $\psi$ as the homomorphism associated to the map
  $\psi^1\co N_1\to \Alg\otimes N_2$ determined by
  $\psi^1(\x)=\phi(1\otimes \x)$. Similarly, if 
  $h\co \cNodule_1\to\cNodule_2$ is an $\DGA$-equivariant
  homotopy (i.e., $h(a\cdot x)=a\cdot h(x)$ for all $a\in\DGA$ and $x\in \cNodule_1$,
  then we can define
  $h^1(\x)=h(1\otimes \x)$ for any $\x\in N_1$. It is straightforward to verify that
  these maps satisfy the required properties.
\end{proof}

\begin{definition}
  \label{def:BoundedTypeD}
Let $A$ be an $\Ainf$ algebra and $N$ a graded left $\Ground$-module.
A map $\delta^1 \co N \to (A \otimes
N)[1]$ can
be iterated to construct maps
$$\gls*{ddiffk}\co N \to (A^{\otimes k} \otimes N)[k]$$
inductively by
\begin{align*}
  \delta^0 &= \Id_N \\
  \delta^i &= \left(\Id_{A^{\otimes(i-1)}} \otimes
    \delta^1\right)\circ \delta^{i-1}.
\end{align*}

By construction, these maps satisfy the basic equation
\[
  (\Id_{A^{\otimes j}}\otimes \delta^{i})\circ \delta^j = 
        \delta^{i+j}
\]
for all $i, j\geq 0$,
or graphically
\[
\mathcenter{
\begin{tikzpicture}
  \node at (0,0) (tc) {};
  \node at (0,-1) (delta) {$\delta$};
  \node at (-1,-2) (Delta) {$\Delta$};
  \node at (-2,-3) (bll) {};
  \node at (-1,-3) (bl) {};
  \node at (0,-3) (bc) {};
  \draw[dmar] (tc) to (delta);
  \draw[dmar] (delta) to (bc);
  \draw[taar] (delta) to (Delta);
  \draw[taar] (Delta) to (bll);
  \draw[taar] (Delta) to (bl);
\end{tikzpicture}}
\mathcenter{\quad=\quad}
\mathcenter{
\begin{tikzpicture}
  \node at (0,0) (tc) {};
  \node at (0,-1) (deltaa) {$\delta$};
  \node at (0,-2) (deltab) {$\delta$};
  \node at (-2,-3) (bll) {};
  \node at (-1,-3) (bl) {};
  \node at (0,-3) (bc) {};
  \draw[dmar] (tc) to (deltaa);
  \draw[dmar] (deltaa) to (deltab);
  \draw[dmar] (deltab) to (bc);
  \draw[taar] (deltaa) to (bll);
  \draw[taar] (deltab) to (bl);
\end{tikzpicture}}.
\]

  We say that $\delta^1$ is \emph{bounded} if for all $\x\in N$,
  \index{bounded!type $D$ structure}%
  \index{module!type $D$!bounded}%
  there is an
  integer $n$ so that for all $i\geq n$, $\delta^i(\x)=0$.

  We call a type $D$ structure \emph{bounded} if its structure map $\delta^1$ is.
\end{definition}

For a bounded map~$\delta^1$, we can think of the $\delta^k$ as fitting 
together to form a map
\begin{align*}
\gls*{ddiffall}&\co N \to \Tensor^*(\DGA[1])\otimes N,\\
\delta(\x)&\coloneqq\sum_{k=0}^{\infty} \delta^k(\x).
\end{align*}
If $\delta^1$ is not bounded then, instead, the $\delta^k$ fit
together to form a map to the completed tensor algebra:
\begin{align*}
\delta&\co N \to \overline{\Tensor}^*(\DGA[1])\otimes
N\coloneqq \left(\prod_{i=0}^\infty (\DGA[1])^{\otimes i}\right) \otimes N,\\
\delta(\x)&\coloneqq\prod_{k=0}^{\infty} \delta^k(\x).
\end{align*}

We can also give a tensor algebra version of
Definition~\ref{def:TypeD}.  This also allows us to generalize it to
the $\Ainf$ setting.
\begin{definition}\label{def:TypeDGeneral}
Fix an $\Ainf$ algebra $\cDGA$ and a pair $(N,\delta^1)$ as
above. Assume that either $(N,\delta^1)$ is bounded or $\cDGA$ is
operationally bounded. Then we say $(N,\delta^1)$ defines a
type~$D$ structure if
\begin{align*}
  (\mu \otimes \Id_N) \circ\delta &= 0,\\
\shortintertext{or equivalently}
  ({\widebar D}\otimes\Id_N)\circ \delta&=0,
\end{align*}
or graphically
\[
\mathcenter{
\begin{tikzpicture}
  \node at (0,0) (tc) {};
  \node at (0,-1) (delta) {$\delta$};
  \node at (-1,-2) (D) {$\overline{D}$};
  \node at (-2,-3) (bl) {};
  \node at (0,-3) (bc) {};
  \draw[dmar] (tc) to (delta);
  \draw[dmar] (delta) to (bc);
  \draw[taar] (delta) to (D);
  \draw[taar] (D) to (bl);
\end{tikzpicture}}
\;=\; 0.
\]

For two type $D$
structures $(N_1, \delta^1_{N_1})$ and $(N_2, \delta^1_{N_2})$, a map $\psi^1
\co N_1 \to A \otimes N_2$ 
\index{homomorphism!of type $D$ modules}%
can be used to construct maps \begin{align*}
\gls*{typedmapk} &\co N_1 \to (\DGA^{\otimes k} \otimes N_2)[k-1]\\
  \psi^k(\x) &\coloneqq\! \sum_{i+j=k-1}\!
    (\Id_{A^{\otimes(i+1)}}\otimes\delta^j_{N_2}) \circ
    (\Id_{A^{\otimes i}}\otimes\psi^1) \circ \delta^i_{N_1}.
\end{align*}
The map $\psi^1$ is \emph{bounded} if $\psi^k=0$ for $k$ sufficiently
\index{bounded!type $D$ homomorphism}%
\index{homomorphism!of type $D$ modules!bounded}%
large. (This is automatic if $N_1$ and $N_2$ are bounded.)  We can define
\[
  \gls*{typedmapall}(\x) \coloneqq \sum_{k=1}^\infty \psi^k(\x), 
\]
which
maps to $\Tensor^*(\DGA[1])\otimes N$ if $\psi^1$ is bounded and
$\overline{\Tensor}^*(\DGA[1])\otimes N$ in general.  Assuming either
$\psi^1$ is bounded or $\cDGA$ is operationally bounded, we say $\psi$ is a
\emph{$D$-structure homomorphism} if
\[
(\mu \otimes \Id_{N_2}) \circ \psi = 0.
\]
\index{homomorphism! of type $D$ modules}%
Similarly, \emph{$\psi$ is homotopic to $\phi$} if
\index{homomorphism!of type $D$ modules!homotopic}%
\[
(\mu \otimes \Id_{N_2}) \circ h = \psi - \phi
\]
for $h$ constructed from some map $h^1 \co N \to (\DGA \otimes
N_2)[-1]$ in the analogous way (and with $h^1$ bounded if $\cDGA$ is
not operationally bounded).

The main complication with type $D$ structures over an $\Ainf$ algebra is that
the category of type $D$ structures over an $\Ainf$ algebra is an
$\Ainf$ category, with composition given graphically as follows:
\[
\circ_k(f_1,\dots,f_k)\;=\;
  \mathcenter{\begin{tikzpicture}[x=1cm,y=32pt]
      \node at (0,1) (blank0) {};
      \node at (0,0) (d0) {$\delta^{M_0}$} ;
      \node at (0,-1) (h1) {$(f_1)^1$};
      \node at (0,-2) (d1) {$\delta^{M_1}$};
      \node at (0,-3) (h2) {$(f_2)^1$};
      \node at (0,-4) (dots) {$\vdots$};
      \node at (0,-5) (hk) {$(f_k)^1$};
      \node at (0,-6) (dk) {$\delta^{M_k}$};
      \node at (-1,-2) (place2) {};
      \node at (-1,-3) (place3) {};
      \node at (-1,-4) (place4) {};
      \node at (-1,-5) (place5) {$\vdots$};
      \node at (-1,-6) (place6) {};
      \node at (-1,-7) (mu) {$\mu$};
      \node at (-1,-8) (blank) {};
      \node at (0,-8) (blank2) {};
      \draw[dmar] (d0) to (h1);
      \draw[dmar] (blank0) to (d0);
      \draw[dmar] (h1) to (d1);
      \draw[dmar] (d1) to (h2);
      \draw[dmar] (h2) to (dots);
      \draw[dmar] (dots) to (hk);
      \draw[dmar] (hk) to (dk);
      \draw[dmar] (dk) to (blank2);
      \draw[taar, bend right=15] (d0) to (place2);
      \draw[aar] (h1) to (place2);
      \draw[taar] (d1) to (place3);
      \draw[aar] (h2) to (place4);
      \draw[aar] (hk) to (place6);
      \draw[taar] (dk) to (mu);
      \draw[taar] (place2) to (place3);
      \draw[taar] (place3) to (place4);
      \draw[taar] (place4) to (place5);
      \draw[taar] (place5) to (place6);
      \draw[taar] (place6) to (mu);
      \draw[aar] (mu) to (blank);
    \end{tikzpicture}}
\]
Note that the condition $\circ_1(f)=0$ is exactly the condition of being a type $D$ homomorphism.
\end{definition}

Specializing Definition~\ref{def:TypeDGeneral} to the case of a \dg
algebra, we obtain Definition~\ref{def:TypeD}.

If $(N_1,\delta^1_{N_1})$ and $(N_2,\delta^1_{N_2})$ are both bounded then
any homomorphism $\psi^1\co N_1\to A\otimes N_2$ between them is
automatically bounded.

There is a generalization of Lemma~\ref{lemma:AssociatedTypeD}.

\begin{lemma}
  \label{lemma:AssociatedTypeDInf}
  Let $\Alg$ be an operationally bounded $\Ainf$ algebra, and 
  $(N,\delta^1)$ a type $D$ structure over it. Then 
  $\DGA \otimes N$ can be given the structure of a left $\Ainf$
  module~$\gls*{cNodule}$,   
  with
  \[
  m_i \coloneqq \sum_{k=0}^\infty (\mu_{i+k}\otimes \Id_N)\circ
  (\Id_{A^{\otimes i}} \otimes \delta^k).
  \]
  Moreover, if $\psi^1\co N_1 \to \DGA\otimes N_2$ is a $D$-structure homomorphism,
  then there is an induced $\Ainf$ homomorphism from $\cN_1$ to
  $\cN_2$, whose $i\th$ component is given by
  \[\sum_{k=1}^\infty (\mu_{i+k+1}\otimes \Id_{N_1})\circ \psi^k.
  \]
  Similarly, a homotopy between two $D$-structure homomorphisms induces
  a homotopy between the associated $\Ainf$ maps.
\end{lemma}

As in the case for differential graded algebras,
we will denote the $\Ainf$ module associated to a type $D$ structure $(N,\delta^1)$
constructed in Lemma~\ref{lemma:AssociatedTypeD} by ~$\cNodule$.

\begin{remark}
  \label{rmk:OtherTypeDs} 
  \index{idempotents!minimal}%
  Type $D$ structures have appeared
  elsewhere in various guises.  For example, let $\mathcal{C}$ be a
  small \dg category over $\Field$. Define an algebra $\Alg(\mathcal{C})$ by
  \[
  \Alg(\mathcal{C})=\bigoplus_{a,b,\in\ob(\mathcal{C})}\Mor(a,b),
  \]
  with multiplication given by composition if defined and zero
  otherwise. View $\Alg$ as a \dg algebra over
  $\Ground=\bigoplus_{a\in\ob(\mathcal{C})}\Field$. Then twisted
  complexes
  (see~\cite{BondalKapranov}) over
  $\mathcal{C}$ are the same as suitably graded type $D$ structures over
  $\Alg(\mathcal{C})$. (Note that this differs slightly from the
  definition in~\cite{SeidelBook} in that we have no condition about
  being lower-triangular.)

  Alternatively, type $D$ structures can be viewed as differential
  comodules over the bar resolution, thought of as a coalgebra.
  See~\cite{LefevreAInfinity,KellerLefevre}.
\end{remark}

\section{Another model for the \textalt{$\Ainf$}{A-infty} tensor product}
\label{sec:DT}

Having introduced type $D$ structures, we now introduce a pairing
$\DT$ between $\Ainf$ modules and type $D$ structures. When defined,
this pairing is homotopy equivalent to the derived tensor product with
the $\Ainf$
module associated to the type $D$ structure, but is smaller; for
instance, it is finite-dimensional if the modules are.  As we will see in
Chapter~\ref{chap:tensor-prod}, this new tensor product also matches better
with the geometry and helps give a proof of the pairing theorem.

Throughout this section, we restrict attention to the case where
$\Ainf=\allowbreak(A,\allowbreak\{\mu_i\}_{i=1}^\infty)$ is
operationally bounded. Note that
this includes the case of \dg algebras.

\begin{definition}\label{def:DT}\index{differential!on $\DT$ product}
  Let $\cModule$ be a right $\Ainf$ module over $\cDGA$ and $(N,\delta^1)$ a
  type $D$ structure. Suppose moreover that either $\cModule$ is a bounded
  $\Ainf$ module or $(N,\delta^1)$ is a bounded type $D$ structure. Then, we can
  form the $\Ground$-module $\cModule\gls*{DT} N= M\otimes_{\Ground} N$,
  equipped with the endomorphism
  $$\gls*{DTdiff}(\x\otimes \y) \coloneqq
  \sum_{k=0}^{\infty} (m_{k+1}\otimes\Id_N)
  (\x\otimes\delta^{k}(\y)).$$ 
  We call $(M \otimes N,\partial^\DT)$ the \emph{box tensor product} of $\cModule$
  and $(N,\delta^1)$.
  \index{tensor product!box}%
  \index{box product|see{tensor product, box}}%
\end{definition}
Graphically, the differential on the box tensor product is given by:
\[
\bdy^\DT\;=\;
\mathcenter{
\begin{tikzpicture}
  \node at (0,0) (tl) {};
  \node at (1,0) (tr) {};
  \node at (1,-1) (delta) {$\delta$};
  \node at (0,-2) (m) {$m$};
  \node at (0,-3) (bl) {};
  \node at (1,-3) (br) {};
  \draw[dmar] (tr) to (delta);
  \draw[dmar] (delta) to (br);
  \draw[amar] (tl) to (m);
  \draw[amar] (m) to (bl);
  \draw[taar] (delta) to (m);
\end{tikzpicture}.
}
\]

\begin{example}
  If $\Alg$ is a \dg algebra and $\cModule$ is a differential module,
  then $\cModule \DT N$ is the ordinary tensor product $\cModule \otimes_\Alg
  \cNodule$.
\end{example}

\begin{example}\label{ex:explicitDT}
  Let $\Alg$ be a finite-dimensional \dg algebra over a ring
  $\Ground=\Field e_1\oplus\dots\oplus\Field e_m$,
  a direct sum of copies of $\Field$.
  Choose a $\Field$-basis $\{a_1,\dots,a_n\}$ for $\Alg$ which is
  homogeneous with respect to the action by $\Ground$,
  in the sense that for each $a_i$, there are elements $e_{\ell_i}$ and
  $e_{r_i}$ with $a_i=e_{\ell_i} a_i
  e_{r_i}$. 
  Given a type $D$ structure $(N,\delta^1)$ over $\Alg$, we have maps
  $D_i\co e_{\ell_i} N\to e_{r_i} N$ for $i=1,\dots n$ uniquely determined by the formula
  \[
  \delta^1=\sum_i a_i\otimes D_i.
  \]
  If $M$ is a right $\Ainf$ module over $\DGA$
  then we can write the differential $\partial^{\DT}$ more explicitly as
\begin{equation}
  \label{eq:ExplicitDT}
  \partial^{\DT}(\x\otimes \y) = \sum m_{k+1}(\x, a_{i_1},\dots, a_{i_{k}})\otimes (D_{i_k}\circ\dots\circ D_{i_1})(\y),
\end{equation}
where the sum is taken over all sequences
$i_1,\dots,i_k$ of elements in $\{1,\dots,n\}$ (including the empty
sequence with $k=0$).
\end{example}

\begin{lemma}
  \label{lem:ThetaComplex}
  Let $\cModule$ be an $\Ainf$ module and $N$ a type $D$ structure. Assume
  that either $M$ or $N$ is bounded. Then the $\Ground$-module
  $M\otimes N$ equipped with the endomorphism $\partial^{\DT}$ is a
  chain complex.
\end{lemma}
\begin{proof}
  The proof is easiest to understand in the graphical notation:
\[
\mathcenter{
\begin{tikzpicture}
  \node at (0,0) (tl) {};
  \node at (2,0) (tr) {};
  \node at (2,-1) (deltaa) {$\delta$};
  \node at (0,-2) (ma) {$m$};
  \node at (2,-3) (deltab) {$\delta$};
  \node at (0,-4) (mb) {$m$};
  \node at (0,-5) (bl) {};
  \node at (2,-5) (br) {};
  \draw[dmar] (tr) to (deltaa);
  \draw[dmar] (deltaa) to (deltab);
  \draw[dmar] (deltab) to (br);
  \draw[amar] (tl) to (ma);
  \draw[amar] (ma) to (mb);
  \draw[amar] (mb) to (bl);
  \draw[taar] (deltaa) to (ma);
  \draw[taar] (deltab) to (mb);
\end{tikzpicture}}
\;=\;
\mathcenter{
\begin{tikzpicture}
  \node at (0,0) (tl) {};
  \node at (2,0) (tr) {};
  \node at (2,-1) (delta) {$\delta$};
  \node at (1,-2) (Delta) {$\Delta$};
  \node at (0,-3) (ma) {$m$};
  \node at (0,-4) (mb) {$m$};
  \node at (0,-5) (bl) {};
  \node at (2,-5) (br) {};
  \draw[dmar] (tr) to (delta);
  \draw[dmar] (delta) to (br);
  \draw[amar] (tl) to (ma);
  \draw[amar] (ma) to (mb);
  \draw[amar] (mb) to (bl);
  \draw[taar] (delta) to (Delta);
  \draw[taar] (Delta) to (ma);
  \draw[taar] (Delta) to (mb);
\end{tikzpicture}}
\;=\;
\mathcenter{
\begin{tikzpicture}
  \node at (0,0) (tl) {};
  \node at (2,0) (tr) {};
  \node at (2,-1) (delta) {$\delta$};
  \node at (1,-2) (D) {$\overline{D}$};
  \node at (0,-3) (m) {$m$};
  \node at (0,-4) (bl) {};
  \node at (2,-4) (br) {};
  \draw[dmar] (tr) to (delta);
  \draw[dmar] (delta) to (br);
  \draw[amar] (tl) to (m);
  \draw[amar] (m) to (bl);
  \draw[taar] (delta) to (D);
  \draw[taar] (D) to (m);
\end{tikzpicture}}
\;=\;0.
\]
The first equality uses the definition of $\delta$ (or that $\delta$
is a comodule map), the second the structure equation for an
$\Ainf$ module, and the third the structure equation for a type $D$ structure.
\end{proof}

\begin{example}
  \label{ex:Modulification}
  Fix an $\Ainf$ algebra $\Alg$ and a left type $D$ structure $N$ over $\Alg$.
  We can view $\Alg$ also as a right $\Ainf$ module over $\Alg$.
  In this case, the chain complex $\Alg\DT N$ constructed above
  coincides with the chain complex underlying the $\Ainf$ module associated
  to $(N,\delta^1)$ from Lemma~\ref{lemma:AssociatedTypeD} (in the case where
  $\Alg$ is a \dg algebra) or~\ref{lemma:AssociatedTypeDInf} (in the general case).
\end{example}

\begin{lemma}
  \label{lem:ThetaComplex2}
  Let 
  $f\co \cModule_1\to \cModule_2$ be a homomorphism of $\Ainf$
  modules and $\phi^1\co N_1\to \DGA\otimes N_2$ 
  a homomorphism of type $D$ structures. Then:
  \begin{enumerate}
  \item If either $\cModule_1$, $\cModule_2$ and $f$ are bounded or
    $(N_1,\delta^1_{N_1})$ is
    bounded then $f$ induces a chain map $f\DT\Id_{N_1}\co \cModule_1\DT N_1\to
    \cModule_2\DT N_1$.
    \glsadd{DTphiId}%
    \index{tensor product!box!of homomorphisms}%
  \item If either $(N_1,\delta^1_{N_1})$ and $(N_2,\delta^1_{N_2})$ are bounded or $\cModule_1$
    is bounded then $\phi^1$ induces a chain map $\Id_{\cModule_1}\DT{\phi^1}\co \cModule_1\DT N_1\to \cModule_1\DT N_2$.
  \item Suppose $f'\co \cModule_1\to \cModule_2$ is
    another homomorphism of $\Ainf$ modules which is homotopic to
    $f$. Suppose further that either $(N_1,\delta^1_{N_1})$ is bounded or
    $\cModule$, $f$, $f'$ and the homotopy between them are all
    bounded. Then $f\DT\Id_{N_1}$ is homotopic to $f'\DT\Id_{N_1}$.
  \item Suppose ${\widetilde{\psi}}^1\co N_1\to \DGA\otimes N_2$ is another homomorphism
    of type $D$ structures which is homotopic to $\psi^1$. Suppose
    further that either $\cModule_1$ is bounded or $\psi^1$, ${\widetilde{\psi}}^1$ and
    the homotopy between them are all bounded. Then $\Id_{\cModule_1}\DT{\psi^1}$ and
    $\Id_{\cModule_1}\DT{\widetilde{\psi}}^1$ are homotopic.
  \item Under boundedness conditions so that all of the $\DT$~products
    are defined, given a third map $g\co \cModule_2\to\cModule_3$, then
    $(g\circ f)\DT\Id_{N_1}$ is homotopic to $(g\DT\Id_{N_1})\circ
    (f\DT\Id_{N_1})$. 
  \item Under boundedness conditions so that all of the $\DT$~products
    are defined, given a third map $\psi^1 \co (N_2,\delta^1_{N_2})\to
    (N_3,\delta^1_{N_3})$, then
    $\Id_{\cModule_1}\DT(\psi^1\circ\phi^1)$ is homotopic to $(\Id_{\cModule_1}\DT\psi^1)\circ
    (\Id_{\cModule_1}\DT\phi^1)$.
  \item Under boundedness conditions so that all of the $\DT$~products
    are defined, $(f\DT\Id_{N_1})\circ(\Id_{\cModule_1}\DT\phi^1)$ is homotopic
    to $(\Id_{\cModule_1}\DT\phi^1)\circ (f\DT\Id_{N_1})$.
  \end{enumerate}
\end{lemma}
\begin{proof}
  Given $f$, we define the induced map $f\otimes \Id_N$
  by 
  $$(f\otimes \Id_N)(\x\otimes \y) \coloneqq
        \sum_{k=0}^{\infty} (f_{k+1}\otimes\Id_N) \circ
          (\x\otimes\delta^{k}(\y)),$$
  where $f_k$ are the components of $f$. 
  Similarly, given a homotopy $h$ with components~$h_k$
  from $\phi$ to $\phi'$, we define
  a homotopy $h\DT\Id_N$ from $f\DT\Id_N$ to $f'\DT\Id_N$ by
  $$(h\DT \Id_N)(\x\otimes \y) = 
  \sum_{k=0}^{\infty} (h_{k+1}\otimes\Id_N) \circ
  (\x\otimes\delta^{k}(\y)).$$
  Graphically, the maps $f\DT \Id_N$ and $h\DT \Id_N$
  are given by
\[
\mathcenter{
\begin{tikzpicture}
  \node at (0,0) (tl) {};
  \node at (1,0) (tr) {};
  \node at (1,-1) (delta) {$\delta$};
  \node at (0,-2) (m) {$f$};
  \node at (0,-3) (bl) {};
  \node at (1,-3) (br) {};
  \draw[dmar] (tr) to (delta);
  \draw[dmar] (delta) to (br);
  \draw[amar] (tl) to (m);
  \draw[amarp] (m) to (bl);
  \draw[taar] (delta) to (m);
\end{tikzpicture}}
\qquad\mathcenter{\text{and}}\qquad
\mathcenter{
\begin{tikzpicture}
  \node at (0,0) (tl) {};
  \node at (1,0) (tr) {};
  \node at (1,-1) (delta) {$\delta$};
  \node at (0,-2) (m) {$h$};
  \node at (0,-3) (bl) {};
  \node at (1,-3) (br) {};
  \draw[dmar] (tr) to (delta);
  \draw[dmar] (delta) to (br);
  \draw[amar] (tl) to (m);
  \draw[amarp] (m) to (bl);
  \draw[taar] (delta) to (m);
\end{tikzpicture}}.
\]

  Given $\phi^1$, we define $\Id_{\cModule}\DT\phi^1$ by
  \[
  (\Id_{\cModule}\DT \phi^1)(\x\otimes \y) = 
        \sum_{k=1}^{\infty} (m_{k+1}\otimes\Id_{N'}) \circ
          (\x\otimes\phi^{k}(\y)),
  \]
  with $\phi^k$ constructed from the map $\phi^1 \co N \to
  A \otimes N'$ as in Definition~\ref{def:TypeDGeneral}. 
  A homotopy $h$ from $\phi^1$ to ${\widetilde\phi}^1$ is promoted
  to a homotopy from $\Id_{\cModule}\DT \phi^1$ to 
  $\Id_{\cModule}\DT {\widetilde\phi}^1$ by a similar formula.

  Verification of the properties of these maps is a straightforward
  modification of the argument that $\partial^{\DT}\circ \partial^{\DT}=0$.
\end{proof}

\begin{example}
  \label{ex:Modulification-cont}
  Continuing from Example~\ref{ex:Modulification},
  if $\psi^1\co N_1\to \DGA\otimes N_2$  is a homomorphism of type
  $D$ structures, then $\Id_{\Alg}\DT \psi^1\co \Alg\DT N_1\to \Alg\DT N_2$
  is the chain map underlying the $\Ainf$ homomorphism from $\cN_1$ to $\cN_2$
  from  Lemmas~\ref{lemma:AssociatedTypeD} or~\ref{lemma:AssociatedTypeDInf}.
\end{example}

\begin{proposition}\label{prop:IdentifyDT}
  Let $\cModule$ be an $\Ainf$ module over an operationally bounded
  $\Ainf$ algebra $\cDGA$ and $(N,\delta^1)$ be a type $D$
  structure. Suppose that either $\cM$ is a bounded $\Ainf$ module or
  $(N,\delta^1)$ is bounded. In the former case, suppose
  moreover that $(N,\delta^1)$ is homotopy equivalent to a bounded
  type $D$ structure.  Then the product $\cModule\DT N$ is
  homotopy equivalent to the $\Ainf$ tensor product $\cModule\DTP
  \cNodule$.
\end{proposition}

\begin{proof}
  The boundedness assumption implies that $\cModule \DT N$ is
  defined. Moreover, since $\DT$ respects homotopy equivalences (by
  Lemma~\ref{lem:ThetaComplex2}), we may assume that $N$ is bounded.
  By Proposition~\ref{prop:BarResolution}, $\cModule$ is $\Ainf$ homotopy
  equivalent to ${\widebar \cModule}$, and hence by
  Lemma~\ref{lem:ThetaComplex2}, $\cM\DT N\simeq {\widebar \cM}\DT N$.
  But this is the same complex with the same boundary operator as $\cM
  \DTP \cN$.  For instance, as $\Ground$-modules we have
  \[
  \widebar{\cM}\DT\cN \cong (M \otimes \Tensor^+(A[1])[-1] \otimes N
  \cong M \otimes \Tensor^*(A[1]) \otimes (A \otimes N)
  \cong \cM \DTP \cN.
  \]
  (Recall that the underlying space of $\cN$ is $A \otimes N$.  Also, the
  $\Ainf$ homotopy equivalence between $\cM$ and $\widebar{\cM}$ is
  not bounded, so we need the boundedness assumption on~$N$.)
\end{proof}

The following consequence of Lemma~\ref{lem:ThetaComplex2} will be
used in Section~\ref{sec:surg-exact-triangle}:

\begin{definition}
  \label{def:ShortExactSequence}
  Let $(N_1,\delta_1)$, $(N_2,\delta_2)$, and $(N_3,\delta_3)$ be
  three type $D$ structures over some fixed \dg algebra $\DGA$, and let 
  $\phi^1\co N_1\to \DGA\otimes N_2$ and
  $\psi^1\co N_2\to \DGA\otimes N_3$
  be two homomorphisms of
  type $D$ structures. We say these form a {\em short exact sequence}
  if the induced modules $\cN_1$, $\cN_2$, $\cN_3$ equipped with
  homomorphisms associated to between them $\Id_{\Alg}\otimes \phi^1$
  and $\Id_{\Alg}\otimes \psi^1$ fit into a short exact sequence
  of differential graded $\DGA$-modules.
  We abbreviate this, writing
  $
  0\longrightarrow N_1\stackrel{\phi^1}{\longrightarrow} N_2\stackrel{\psi^1}{\longrightarrow}
  N_3\longrightarrow 0.$
\end{definition}

\begin{proposition}
  \label{prop:InducedLongExactSequence}
  If $0\longrightarrow N_1\stackrel{\phi^1}{\longrightarrow} N_2\stackrel{\psi^1}{\longrightarrow}
  N_3\longrightarrow 0$
  is a short exact sequence of type $D$ structures over a \dg algebra
  and $\cModule$ is a bounded $\Ainf$ module, then
  there is an exact sequence in homology
  \[
    H_*(\cModule\DT N_1)\longrightarrow H_*(\cModule\DT N_2)\longrightarrow
    H_*(\cModule\DT N_3)\longrightarrow H_*(\cModule\DT N_1)[-1].
  \]
\end{proposition}
\index{tensor product!box!is triangulated functor}%
\begin{proof}
  The short exact sequence can be viewed as a $3$-step filtered type $D$ structure $N$. 
  Since $\DGA \DT N_3$ is a projective $\DGA$-module, this short exact sequence
  is split. The splitting can be viewed as a null-homotopy of the $3$-step filtered complex.

  Next, form $\cModule\DT N$ to get a three-step filtered complex,
  whose associated graded complex is $\bigoplus_{i=1}^3 \cModule\DT
  N_i$.  (It is perhaps worth noting that the differential of an element
  in $\cModule \DT N_1\subset \cModule\DT N$ can have a component in
  $\cModule\DT N_3$. In fact, the $\cModule\DT N_3$ component of the
  differential from $\cModule\DT N_1$ is precisely the null-homotopy
  of $(\Id_{\cModule}\DT \psi^1)\circ(\Id_{\cModule}\DT \phi^1)$.)  By
  Lemma~\ref{lem:ThetaComplex2}, this complex $\cModule\DT N$ is
  homotopy equivalent to the trivial complex; in particular,
  $\cModule\DT N$ is acyclic.  It follows that there is a long exact
  sequence on the homology groups of the terms in its associated
  graded complex. (This latter point is a small generalization of the
  fact that a short exact sequence of chain complexes induces a long
  exact sequence in homology: it is a degenerate case of the Leray
  spectral sequence.)
\end{proof}

\begin{remark}
  \index{tensor product!$\Ainf$!is triangulated functor}%
  Proposition~\ref{prop:InducedLongExactSequence} is an explicit form
  of the statement that $\cModule\DT\cdot$ is a triangulated
  functor. Since $\cModule\DT N\simeq \cModule\DTP\cNodule$
  in a natural way, this follows from the fact that the $\Ainf$ tensor
  product is a triangulated functor; over a \dg
  algebra, this, in turn, follows from the fact that the derived tensor product is a
  triangulated functor.
\end{remark}

\section{Gradings by non-commutative groups}
\label{sec:grad-non-comm}

In earlier sections, we considered the familiar case of $\ZZ$-graded
modules.   Here we record basic definitions of $\Ainf$ algebras
(including differential algebras) that are
graded over a possibly non-commutative group~$G$.  The definitions
are straightforward generalizations of standard definitions.  As
before, we assume that our ground ring~$\Ground$ has
characteristic~$2$ so that we need not worry about signs.

We start with ordinary differential graded algebras.

\begin{definition} \label{def:dga-nc-grading}
  Let $\gls*{Glambda}$ 
  be a pair of a group~$\gls*{Group}$ (written
  multiplicatively) with a distinguished
  element~$\gls*{lambda}$ in the center of~$G$.
  A \emph{differential algebra graded by~$(G,\lambda)$}
  \index{grading!differential algebra graded by a group}%
  is, 
  first of
  all, a differential
  algebra~$A$ with a grading $\gr$ of $\DGA$ by $G$ (as a set),
  i.e., a decomposition $\DGA=\bigoplus_{g\in G} \gls*{DGAsubg}$. 
  We say an element $a$ in $\DGA_g$ is {\em homogeneous of degree~$g$}
  \index{homogeneous element of group-graded algebra|see{grading, group-valued}}%
  \index{grading!group-valued, on algebra}%
  \index{grading!group-valued, on algebra!homogeneous element of}%
  \index{G-graded algebra@$G$-graded algebra|see{grading}}%
  and 
  write $\gr(a)=g$.  For homogeneous elements $a$ and~$b$, we
  further require the grading to be
  \begin{itemize}
  \item compatible with the product, i.e., $\gr(a\cdot b) =
    \gr(a)\gr(b)$, and
  \item compatible with the differential, i.e., $\gr(\partial
    a) = \lambda^{-1}\gr(a)$.
  \end{itemize}
\end{definition}

For an ordinary $\ZZ$-grading, we would take $(G,\lambda) = (\ZZ,1)$.
(We choose the differential to lower the grading to end up with
homology rather than cohomology.)
We require $\lambda$ to be central because otherwise the identity
\[\partial(ab) = \partial a \cdot b + a \cdot \partial b\]
would not be homogeneous in general.
The map $k \mapsto \lambda^k$  is a homomorphism from $\ZZ$ to~$G$, but
even if $\lambda$ is of infinite order in~$G$, there may not be a
map from~$G$ to~$\ZZ$ taking $\lambda$ to~$1$, so we do not in general get an
ordinary $\ZZ$-grading.

\begin{remark}
  When working over a field whose characteristic is not~$2$, we need a
  $\ZZ/2$ grading (i.e., a homomorphism from $G$ to $\ZZ/2$ taking
  $\lambda$ to~$1$) in order to, for instance, specify the signs in
  the Leibniz relation.  We will not consider this further here.
\end{remark}

There is a natural extension of this definition to $\Ainf$ algebras.
In general, if
$V$ and~$W$ are $G$-graded $\Ground$-modules, we can turn $V \otimes
W$ into a $G$-graded $\Ground$-module by setting
\begin{equation}\label{eq:grading-tensor}
\gr(v \otimes w) \coloneqq \gr(v) \gr(w)
\end{equation}
\index{grading!on tensor product of group-graded vector spaces}%
for homogeneous elements $v$ and~$w$ of $V$ and~$W$, respectively.  If
$V$ is a $G$-graded module and $n \in \ZZ$, define $\gls*{shiftV}$ 
to be the graded module with
\begin{equation}\label{eq:grading-shift-noncomm}
V[n]_g \coloneqq V_{\lambda^{-n}g}.
\end{equation}
Since $\lambda$ is central, we have $(V[n]) \otimes W \cong (V \otimes
W)[n]\cong V\otimes (W[n])$.
Now the previous definitions of
the multiplications maps~$m_i$ and their compatibility conditions carry
over as written.

We now turn to differential modules.
\begin{definition}\label{def:module-nc-grading}
  Let $A$ be a differential algebra graded by $G$, as in
  Definition~\ref{def:dga-nc-grading}.  Let $S$ be a set with a right
  $G$ action.
  A \emph{right differential $A$-module graded by $S$},
  \index{grading!differential module graded by $G$-set}%
  \index{set-graded|see{grading}}%
  \index{S-graded module@$S$-graded module}%
  or, simply, a right
  $S$-graded module,
  is a right differential $A$-module~$M$ with a grading $\gr$ of $M$
  by~$S$ (as a set)
  so that, for homogeneous elements $a \in A$ and $x \in M$,
  \begin{itemize}
  \item $\gr(xa) = \gr(x)\gr(a)$ and
  \item $\gr(\partial x) = \gr(x)\lambda^{-1}$.
  \end{itemize}
\end{definition}

If the action of $G$ on $S$ is transitive,
$S$ is determined by the stabilizer $\gls*{Gsubs}$ 
of any point $s \in
S$, and we can identify $S$ with $\gls*{GsubsCosets}$, 
the space of right cosets of~$G_s$.

More generally, if $V$ is graded by $S$ and $W$ is graded by $G$, then
Equation~\eqref{eq:grading-tensor} gives a natural grading of $V
\otimes W$ by~$S$ and Equation~\eqref{eq:grading-shift-noncomm}
defines a shift functor; then the previous definitions for $\Ainf$
modules, homomorphisms, and homotopies carry through.  (For
homomorphisms and homotopies, the two modules should be graded by the
same set.)  For
example, for $\Ainf$ modules $M$ over a graded $\Ainf$ algebra $\Alg$,
if $x$ is a homogeneous element of $M$ and $a_1,\dots,a_\ell$ are
homogeneous elements of $\Alg$, then
$y=m_{\ell+1}(x,a_1,\dots,a_\ell)$ is homogeneous of degree
\begin{equation}
  \label{eq:DegreeAinfAction}
  \gr(y)=\lambda^{\ell-1}\gr(x)\cdot \gr(a_1)\cdots \gr(a_\ell).
\end{equation}

Left modules are defined similarly, but with gradings by a set~$T$
with a left action by~$G$.

Suppose the pair $(G,\lambda)$ is such that $\lambda$ has infinite
order. Then there is a partial order on $G$ given by
$g\gls*{lessthan}h$ if $h=g\lambda^n$ for some $n>0$. So, it sometimes
makes sense to say that $\gr(a)<\gr(b)$ for homogeneous elements $a,b$
in a group-graded algebra. Similarly, for a $G$-set $S$, if the
subgroup $\langle\lambda\rangle$ acts freely on $S$ then there is a
partial order on $S$ given by $s<t$ if $t=\lambda^n s$ for some
$n>0$. Thus, for a set-graded module $M$, it sometimes makes sense to
write $\gr(x)\gls*{lessthan}\gr(y)$ for homogeneous elements $x,y\in
M$.

Finally we turn to tensor products.

\begin{definition}
  \index{grading!on tensor product of set-graded modules}%
  Let $A$ be a differential algebra graded by $(G,\lambda)$, let $M$ be
  a right $A$-module graded by a right $G$-set~$S$, and let $N$ be a
  left $A$-module graded by a left $G$-set~$T$.  Define $S \times_G
  T$ by
  \[
  \gls*{ProductGSets}
  \coloneqq (S \times T)/
    \{\,(s,gt) \sim (sg,t)\mid s \in S, t \in T, g \in G\,\}.
  \]
  Define a $(S \times_G T)$ grading on $M \otimes_\Ground
  N$ by
  \begin{equation}\label{eq:grading-tensor-left-right}
    \gls*{grTensorProd}
    \coloneqq [\gr(m) \times \gr(n)].
  \end{equation}
  By inspection, this descends to
  the na\"{\i}ve tensor product
  \[
  M \otimes_A N \coloneqq M \otimes_{\Ground} N/
      \langle\, m\otimes an - ma\otimes n\mid m \in M, n \in N, a \in A\,\rangle.
  \]
  Furthermore, there is a natural action of $\lambda$ (or any element
  in the center of~$G$) on $S \times_G T$ via
  \[
  \lambda \cdot [s \times t] \coloneqq [s\lambda \times t] =
    [s \times \lambda t].
  \]
  (This definition fails to descend to the quotient for non-central
  elements in place of~$\lambda$.) Again by inspection, the boundary
  operator
  \[
  \bdy(m \otimes_A n) \coloneqq (\bdy m) \otimes_A n +
    m \otimes_A (\bdy n)
  \]
  acts by $\lambda^{-1}$ on the gradings.  Thus $M \otimes_A N$ is a
  chain complex with a grading by $S \times_G T$, which is a set with a
  $\ZZ$-action.
\end{definition}

Using Equation~\eqref{eq:grading-tensor-left-right} and the same shift
functor as before, we can extend the definitions of both versions of
the tensor product, chain maps and homotopies to this setting as
well.

We now see a genuinely new feature of \dg algebras where the grading is
non-commutative: even if $\lambda$ has infinite order in $G$, $S$, and
$T$, it may not have infinite order on $S \times_G T$.
We illustrate this with a simple, and pertinent, example.  

\begin{example}\label{ex:order-tensor}
Let $H$ be a variant of the
Heisenberg group:
\index{Heisenberg group}%
\[
H = \langle\,\alpha, \beta, \lambda \mid
   \alpha\beta = \lambda^2\beta\alpha, \text{$\lambda$ is central}\,\rangle.
\]
(The usual Heisenberg group over the integers has $\lambda$ instead of
$\lambda^2$ in the first relation; with the exponent of 2,
the group supports a $\ZZ/2$ grading.)
Let $L$ be the subgroup of~$H$ generated by~$\alpha$.  Then $\lambda$
has infinite order on~$H$ and on $H/L$ and $L \backslash H$, the
spaces of left and right cosets of $L$.  However, the reduced product is
\begin{align*}
(L \backslash H) \times_H (H / L) &\isom L \backslash H / L\\
  &= \{\,L\beta^k\lambda^l L\mid k,l \in \ZZ\,\}.
\end{align*}
These double cosets are not all distinct.  In particular,
\begin{align*}
  L\beta^k\lambda^l L &= L\alpha\beta^k\lambda^l L\\
    &= L\beta^k\lambda^{l+2k}\alpha L\\
    &= L\beta^k\lambda^{l+2k}L.
\end{align*}
Equalities like this generate all equalities between the double cosets.
Thus the order of $\lambda$ on $L\beta^k\lambda^l L$ is $2k$, which is
not free when $k \ne 0$.
\end{example}

This example is exactly what we will see when we do $0$-surgery on a knot
in $S^3$. (See Section~\ref{sec:GradedPairingThm}.)


\chapter{The algebra associated to a pointed matched circle}
\label{chap:algebra}

The aim of the present chapter is to define the algebra $\Alg(\PMC)$
associated to the
boundary of a three-manifold, parametrized by the pointed matched
circle~$\PMC$. In
Section~\ref{sec:strands-algebra}, we start by defining the strands
algebra~$\Alg(n)$, which is a simplified analogue of this algebra.
The desired algebra~$\Alg(\PMC)$, introduced in
Section~\ref{sec:matched-circles},
is described as a subalgebra of $\Alg(n)$.  In
Section~\ref{sec:gradings-algebra}, we explain the grading
on~$\Alg(\PMC)$.

Although the algebras are defined combinatorially, the motivation for
their definitions comes from the behavior of holomorphic curves. The
reader might find it helpful to consult Examples
\ref{eg:dd1}--\ref{eg:dd3} and \ref{eg:a1}--\ref{eg:a4}, as well
as~\cite{LOT0}, in conjunction with this chapter.

Our algebras, by construction, have characteristic $2$. With a little more
work, the algebras could be lifted to give a construction over $\ZZ$;
for the algebra $\Alg(n)$ this is done in~\cite{Khovanov10:gl12}.
Indeed, with considerably more work, one should be able to define
bordered Floer homology over $\ZZ$.

\section{The strands algebra \textalt{\protect{$\Alg(n,k)$}}{A(n,k)}}
\label{sec:strands-algebra}

We will describe an algebra $\Alg(n,k)$ in three different ways,
which will be useful for different purposes.

\subsection{Algebraic definition}
\label{sec:algebra-algebraic}

For non-negative integers $n$ and $k$ with $n \ge k$, let $\gls*{Ank}$
be the $\Field$--vector space generated by partial permutations $(S, T, \phi)$,
where $S$ and $T$ are $k$-element subsets of the set $\{1,\dots,n\}$
(hereafter denoted $\gls*{nset}$) 
and $\phi$ is a bijection from $S$ to $T$, with
$i \le \phi(i)$ for all $i \in S$.
For brevity, we will adopt a 2-line notation for these
generators, so that, for instance, 
$\gls*{twoline}$
denotes the algebra element 
$
(\{1,2,4\}, \{3,4,5\}, \{1 \mapsto 5, 2 \mapsto 3, 4 \mapsto 4\})
$
inside $\Alg(3,5)$.

For $a = (S, T, \phi)$ a generator
of $\Alg(n,k)$, let $\gls*{invnum}$ or $\gls*{invnumphi}$ be the number of
\emph{inversions}
\index{inversions}%
of~$\phi$: the number of pairs $i, j \in S$ with $i
< j$ and $\phi(j) < \phi(i)$.  We make $\Alg(n,k)$ into an algebra using
composition.  For $a, b \in \Alg(n,k)$, with $a = (S,T,\phi)$, $b
= (T,U,\psi)$, and $\inv(\psi\circ \phi) = \inv(\phi) +
\inv(\psi)$, define $a\cdot b$ to be $(S, U, \psi\circ \phi)$.  In
all other cases, (i.e., if the range of $\phi$ is not the domain of
$\psi$, or if $\inv(\psi \circ \phi) \neq  \inv(\phi) +
\inv(\psi)$), the product of two generators is zero.

There is one idempotent $\gls*{IdemAnkS}$ in $\Alg(n,k)$ for each $k$-element
\index{idempotents!of $\Alg(n,k)$}\index{idempotents!minimal}\index{minimal idempotent|see{idempotents, minimal}}%
subset~$S$ of $[n]$; it is $(S, S, \Id_S)$.  It is easy to see that
the $I(S)$ are all the minimal idempotents in $\Alg(n,k)$. (\emph{Minimal
idempotents} are
the idempotents that can not be further decomposed as sums of
orthogonal idempotents).  Let
$\gls*{Idemnk} \subset \Alg(n,k)$ 
be the subalgebra generated by the idempotents.

We also define a (co)differential on $\Alg(n,k)$.
For $a =
(S,T,\phi)$ in $\Alg(n,k)$, let $\gls*{invset}$ 
be the set counted in
$\inv(\phi)$: the set $\{\,(i, j) \mid i < j, \phi(j) <
\phi(i)\,\}$. For $\sigma \in \Inv(\phi)$, let $\gls*{resolvecross}$
be defined by
\[
\phi_\sigma(k) =
\begin{cases}
  \phi(j) & k = i\\
  \phi(i) & k = j\\
  \phi(k) & \text{otherwise}.
\end{cases}
\]
Note that $\inv(\phi_\sigma) < \inv(\phi)$, since the inversion
$(i,j)$ was removed from $\Inv(\phi)$ and none are added.  It
is also possible that $\inv(\phi_\sigma)<\inv(\phi)-1$, since
more inversions might be removed.
Define
\index{differential!on $\Alg(n,k)$}%
\[
\partial a \coloneqq
  \!\!\!\sum_{\substack{\sigma\in\Inv(\phi)\\
                                \inv(\phi_\sigma) = \inv(\phi)-1}}
    \!\!\!(S, T, \phi_\sigma).
\]

\begin{lemma}\label{lem:Ank-is-dga}
  The set $\Alg(n,k)$, equipped with a product induced by composition
  of partial permutations and with the above-defined
  endomorphism~$\partial$, is a differential algebra.
\end{lemma}
We defer the proof until Section~\ref{sec:algebra-geometric}, in which we develop a convenient geometrical language for $\Alg(n,k)$.

\begin{definition}\label{def:strands-algebra}
  The differential algebra $\Alg(n,k)$ is called the \emph{strands
    algebra with $k$ strands and $n$ places.}
  \index{strands algebra}%
  \index{algebra!strands|see{strands algebra}}%
  The direct sum
  $\gls*{Algn}=\bigoplus_{k}\Alg(n,k)$ is called the \emph{strands
    algebra with $n$ places}.
\end{definition}

\begin{remark}
  The function $\inv$ provides a $\ZZ$-grading on $\Alg(n)$.
\end{remark}
\begin{remark}
  Setting the product in $\Alg(n,k)$ to be zero if the number of inversions in the
  composition $\psi \circ \phi$ decreases is reminiscent of the
  nilCoxeter algebra \cite{FominStanley94:nilCoxeter, Khovanov01:Nilcoxeter}.
\end{remark}
\subsection{Geometric definition}
\label{sec:algebra-geometric}

There is an alternate description of the algebra $\Alg(n)$, in terms
of \emph{strands diagrams}.
\index{strands diagram}%
 We represent a triple $(S, T, \phi)$ by a diagram with
$n$~dots on both the left and right (numbered from the bottom up), with a
set of strands representing $\phi$ connecting the subsets $S$ on the
left and $T$ on the right.  The condition that $i \le \phi(i)$ means
that strands stay horizontal or move up when read from left to
right.  For instance, the diagram
\[
\gls*{graphicalnotation}
\]
represents the algebra element $\strands{1&2&4}{5&3&4}$.

It is easy to see that $\inv(a)$ is the
minimal number of crossings in any strand diagram.
Accordingly, we may view basis elements of $\Alg(n)$ as strand diagrams as above, with
no pair of strands crossing more than once, considered up to
homotopy. The summand $\Alg(n,k)\subset\Alg(n)$ corresponds to strand
diagrams with exactly $k$ strands.

\index{strands algebra}%
We can define the product on $\Alg(n,k)$ in these pictures using
horizontal juxtaposition.  The product $a\cdot b$ of two algebra
elements is $0$
if the right side of~$a$ does not match the left side of~$b$;
otherwise, we place the two diagrams next to each other and join them
together. 
If in this juxtaposition two strands cross each other twice,
we set the product to be~$0$:
\begin{equation}\label{eq:double-cross-prod}
\mfigb{strands-10} = 0.
\end{equation}
This corresponds to the condition that,
with $a=(T,U,\phi)$ and $b=(S,T,\psi)$, we set
$a\cdot b=0$
when $\inv(\phi\circ\psi) \neq \inv(\phi)+\inv(\psi)$.
If there are no double crossings then the juxtaposition is the product.
 For instance, the algebra element above can be factored as
\[
\mfigb{strands-1} \cdot \mfigb{strands-2} =
  \mfigb{strands-0}.
\]

The boundary operator $\partial a$ also has a natural graphical
description.
\index{differential!on $\Alg(n,k)$}%
\index{smoothing a crossing}%
For any crossing in a strand diagram, we can consider \emph{smoothing}
it:
\[
\mfigb{strands-20} \mapsto \mfigb{strands-21}.
\]
The boundary of a strand diagram is the sum over all
ways of smoothing
one crossing of the diagram, where we remove those diagrams from the sum 
for which there are two strands which cross each other twice (compare Equation~\eqref{eq:double-cross-prod}). For instance, we have
\begin{equation}\label{eq:double-cross-boundary}
\partial\left(\mfigb{strands-25}\right) =
   \mfigb{strands-26} +
   \cancelto{0}{\mfigb{strands-27}} +
   \mfigb{strands-28}.
\end{equation}
It is straightforward to check that this agrees with the earlier
description.  

\begin{proof}[Proof of Lemma~\ref{lem:Ank-is-dga}]
  Consider the larger algebra $\bigAlg(n,k)$ generated over $\Field$
  by strands diagrams but where we do not set double crossings to
  $0$. (We also do not allow isotoping away double crossings.) It is
  obvious that $\bigAlg(n,k)$ is a differential algebra, i.e., that
  the product is associative and the differential satisfies the
  Leibniz rule and $\bdy^2=0$.  Let
  $\bigAlg_d(n,k)\subset\bigAlg(n,k)$ be the sub-vector space of
  $\Alg(n,k)$ generated by all strands diagrams with at least one
  double crossing. We claim that $\bigAlg_d(n,k)$ is a differential
  ideal in $\bigAlg(n,k)$; since $\Alg(n,k)\cong
  \bigAlg(n,k)/\bigAlg_d(n,k)$, the result follows. Clearly,
  $\bigAlg(n,k)$ is an ideal; and if $s\in\bigAlg_d(n,k)$ is a strand
  diagram with at least two double crossings then
  $\bdy(s)\in\bigAlg_d(n,k)$. If $s$ is a strand diagram with a single
  double crossing, there are two terms in $\bdy(s)$ which do not lie
  in $\bigAlg_d(n,k)$, but these terms cancel.
\end{proof}

\begin{remark}
  In combinatorial terms, the terms in the boundary of an algebra
  element~$a$ correspond to the edges in the Hasse diagram of the
  Bruhat order on the symmetric group that start from the
  permutation corresponding to~$a$.  More
  precisely, the \emph{Bruhat order}\index{Bruhat order} on the
  symmetric group may be defined by setting $w > w'$ if, when we write
  $w = s_1\dots s_r$ as a reduced product of adjacent transpositions,
  we can delete some of the $s_i$ to obtain a product
  representing~$w'$.  This turns out
  not to depend on the reduced expression for $w$ that we picked.
  (See, e.g., \cite[Section 2.2]{BB05:CombinatricsCoxeter}.)  If we
  represent a reduced expression for $w$ by drawing strands crossing
  each other at each adjacent transposition, then deleting an $s_i$
  corresponds to smoothing a crossing.  The \emph{Hasse
    diagram}\index{Hasse diagram} of a partially ordered set is the
  graph connecting
  two elements $w, w'$ if one covers the other, i.e., if $w > w'$ and
  there is no other element $w''$ with $w > w'' > w'$.  In the case of
  Bruhat order on the symmetric group, this is equivalent to $w > w'$
  and $\inv(w) = \inv(w')+1$ (see \cite[Theorem
  2.2.6]{BB05:CombinatricsCoxeter}, \textit{inter alia}), just as in
  $\Alg(n,k)$.
\end{remark}

\subsection{Reeb chords}
\label{sec:reeb-chords-def}

There is one more description of $\Alg(n,k)$, relating it to the Reeb
chords we will see from the moduli spaces.  Fix an oriented circle~$\gls*{Circle}$
with a set of $n$ distinct points $\CircPts = \{a_1, \dots, a_n\}
\subset Z$,
together with a basepoint $z \in Z\setminus \CircPts$.  We view $Z$ as a
contact $1$-manifold, and $\CircPts$ as a Legendrian
submanifold of~$Z$. A \emph{Reeb
  chord}~$\gls*{ReebChord}$
\index{Reeb chords}%
\index{chord|see{Reeb chords}}%
in $(Z\setminus z,\CircPts)$ is an immersed (and consequently
embedded) arc in $Z\setminus z$ with endpoints in~$\CircPts$, with
orientation induced by the orientation on~$Z$. The initial endpoint of
$\rho$ is denoted 
$\gls*{InitialEndpt}$ and the final endpoint $\gls*{TerminalEndpt}$.
A
set~$\gls*{ReebChords} = \{\rho_1,\dots,\rho_j\}$ of $j$~Reeb chords is said to be
\emph{consistent} 
\index{consistent set of Reeb chords}%
\index{Reeb chords!set of!consistent}%
if 
$\gls*{InitialEndpts}\coloneqq\{\rho_1^-,\dots,\rho_j^-\}$
and
$\gls*{TerminalEndpts}\coloneqq\{\rho_1^+,\dots,\rho_j^+\}$ are both sets with
$j$~elements (i.e., no two $\rho_i$ share initial or final
endpoints).  We adopt a two-line notation for consistent sets of Reeb
chords so that, for instance, 
\gls*{chordnotation}
is the set of
chords $\{[1,5],[2,3]\}$.

\begin{definition}
\label{def:alg0assoc-to-rhos}
Let $\rhos$ be a consistent set of Reeb chords. The \emph{strands algebra element associated
to $\rhos$},
\index{strands algebra!element associated to $a_0(\rhos)$}%
denoted $a_0(\rhos)\in\Alg(n)$, is defined by
\[
\gls*{asubzero}\coloneqq\sum_{\{S \mid S\cap (\rhos^-\cup\rhos^+)=\emptyset\}} (S\cup\rhos^-,S\cup\rhos^+,\phi_{S})\] 
where $\phi_{S}|_{S}=\Id$ and $\phi_{S}(\rho_i^-)=\rho_i^+$. Geometrically, this corresponds to taking the
sum over all ways of consistently adding horizontal strands to the set of strands
connecting each $\rho_i^-$ to $\rho_i^+$.
For convenience, define
$a_0(\rhos)$ to be~$0$ if $\rhos$ is not consistent.
\end{definition}

The strands algebra
$\Alg(n)$ is generated as an algebra by the elements $a_0(\rhos)$ and
idempotents~$I(S)$, with $S \subset \CircPts$. In
fact, the basis for $\Alg(n)$ over~$\Field$ that we constructed in
Section~\ref{sec:algebra-algebraic} are the non-zero elements of the form
$I(S)a_0(\rhos)$, where $S\subset \CircPts$ is a $k$-element subset,
and $\rhos$ is a consistent set of Reeb chords.

We give a description of the product on $\Alg(n,k)$ in terms of
this basis, after introducing two definitions.

\begin{definition}\label{def:interleaved-nested}
A pair of Reeb chords $(\rho,\sigma)$ is said to be \emph{interleaved}
\index{interleaved Reeb chords}%
\index{Reeb chords!interleaved}%
if $\rho^- < \sigma^- < \rho^+ < \sigma^+$.  They are said to be
\emph{nested} if $\rho^- < \sigma^- < \sigma^+ < \rho^+$.  An
\index{Reeb chords!nested}%
unordered pair of Reeb chords is said to be interleaved (respectively
nested) if either possible ordering is interleaved (respectively
nested).
\end{definition}

\begin{definition}
\label{def:Composable-one}
If $\rho$ and $\sigma$ are Reeb chords such that $\rho^+ = \sigma^-$,
we say that $\rho$ and~$\sigma$
\emph{abut} and define their {\em join} 
\index{Reeb chords!abutting}%
\index{Reeb chords!join of}%
\index{join!of Reeb chords}%
\index{abutting!Reeb chords}%
$\rho \gls*{uplusone} \sigma$ to be the Reeb chord
$[\rho^-,\sigma^+]$ obtained by concatenating the two chords.
\end{definition}
Note that abutting is not a symmetric relation.

\begin{definition}
\label{def:Composable-sets}
For two consistent sets of Reeb chords $\rhos$, $\sigmas$,
define $\rhos\gls*{uplus}\sigmas$ to be obtained
from the union $\rhos \cup
\sigmas$ by replacing every abutting pair $\rho_i \in \rhos$ and
$\sigma_j \in \sigmas$ by their join $\rho_i \uplus\sigma_j$.  If
$\rho_i \in \rhos$, let
\begin{equation*}
  \rho_i^{++} \coloneqq
  \begin{cases}
    \sigma_j^+&\textrm{$\rho_i$, $\sigma_j$ abut for some
      $\sigma_j\in\sigmas$}\\
    \rho_i^+ & \textrm{otherwise.}
  \end{cases}
\end{equation*}
\glsadd{TerminalTerminalEndpt}\glsadd{InitialInitialEndpt}%
Similarly define $\sigma_j^{--}$ for $\sigma_j\in\sigmas$.

For $\rhos$ and $\sigmas$ consistent sets of Reeb chords,
we say that $\rhos$ and $\sigmas$ are \emph{composable}
\index{Reeb chords!composable!sets of}%
\index{composable!sets of Reeb chords}%
 if $\rhos\uplus\sigmas$ is consistent and has no double crossings;
 precisely, if $\rhos\uplus\sigmas$ is consistent and
there is no pair $\rho_i\in\rhos$ and $\sigma_j\in\sigmas$ so that
$\rho_i^- < \sigma_j^{--}$,
$\sigma_j^- < \rho_i^+$,  and
$\rho_i^{++} < \sigma_j^+$.
Finally, we iterate this notion to sequences of consistent sets of Reeb chords
$\vec{\rhos}=(\rhos_1,\dots,\rhos_n)$: we say such a sequence is
{\em composable}
if for all $i=1,\dots,n-1$, the sets $\biguplus_{j=1}^{i}\rhos_j$ and $\rhos_{i+1}$
are composable.
\index{Reeb chords!composable!sequence of sets of}%
\index{composable!sequence of sets of Reeb chords}
\end{definition}
\begin{lemma}\label{lem:reeb-product}
  Given consistent sets of Reeb chords $\rhos$ and $\sigmas$ and a subset
  $S\subset[4k]$ with $I(S)a_0(\rhos)\neq 0$,
  we have
  \[
  I(S)a_0(\rhos)a_0(\sigmas) =
  \begin{cases}
    I(S)a_0(\rhos\uplus\sigmas)&\textrm{$\rhos$ and $\sigmas$
      composable}\\
    0 & \textrm{otherwise.}
  \end{cases}
  \]
\end{lemma}

\begin{proof}
  This follows directly from the definitions. 
\end{proof}

We now describe the differential on $\Alg(n)$ in terms of Reeb chords.
Observe that  there are two ways a term in the boundary operator
can appear, depending on whether we smooth a crossing between two
moving strands or between a moving strand and a horizontal strand.

\begin{figure}
  \centering
  \includegraphics{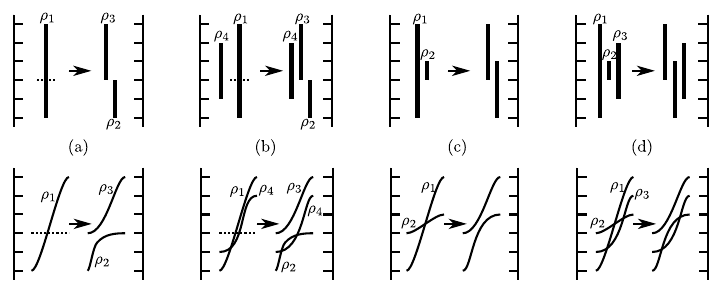}
  \caption[Splittings and shuffles of Reeb chords]{\textbf{Splittings and shuffles of Reeb chords.} (a) A splitting. (b) A weak splitting which is not a splitting. In both cases, the position at which the splitting occurs is marked with a dotted line. (c) A shuffle. (d) A weak shuffle which is not a shuffle. In all cases, the chords are shown on the top line, and the associated strands diagrams on the bottom line.}
  \label{fig:splitshuffle}
\end{figure}

\begin{definition}\label{def:split-shuffle}
  For a Reeb chord $\rho_1$, a \emph{splitting of $\rho_1$}
  \index{Reeb chords!splitting of}%
  is a pair of
  abutting chords $(\rho_2,\rho_3)$ so that $\rho_1 =
  \rho_2\uplus\rho_3$.  For a set $\rhos$ of Reeb chords, a \emph{weak
    splitting of $\rhos$} 
  \index{Reeb chords!set of!weak splitting of}%
is a set obtained from $\rhos$ by replacing
  a chord $\rho_1\in\rhos$ with two chords $\{\rho_2,\rho_3\}$, where
  $(\rho_2,\rho_3)$ is a splitting of $\rho_1$.
  A \emph{splitting of $\rhos$} is a weak splitting of
  \index{Reeb chords!set of!splitting of}%
  \index{Reeb chords!splitting of!set of}%
  $\rhos$, replacing $\rho_1$ with $\{\rho_2,\rho_3\}$,
  which is consistent
  and does not introduce double
  crossings, i.e., there is no chord~$\rho_4\in\rhos$ nested in
  $\rho_1$ and with $\rho_4^- < \rho_2^+ = \rho_3^- < \rho_4^+$. See Figure~\ref{fig:splitshuffle}.

  For a set of Reeb chords~$\rhos$, a \emph{weak shuffle of $\rhos$}
  \index{Reeb chords!set of!weak shuffle of}%
  is a set obtained by replacing a nested pair $(\rho_1,\rho_2)$ of
  chords in $\rhos$ with the corresponding interleaved pair
  $([\rho_1^-,\rho_2^+], [\rho_2^-,\rho_1^+])$.
  A \emph{shuffle of $\rhos$} is a weak shuffle with
  \index{Reeb chords!set of!shuffle of}%
  nested pair $(\rho_1,\rho_2)$ that does not introduce double
  crossings, i.e., there is no chord $\rho_3 \in \rhos$ with $\rho_2$
  nested in~$\rho_3$ and $\rho_3$ nested in~$\rho_1$. Again, see Figure~\ref{fig:splitshuffle}.
\end{definition}

\begin{lemma}\label{lem:reeb-differential}
  We have
  \[
  \partial a_0(\rhos) =
  \sum_{\substack{\textrm{$\rhos'\!$ a splitting}\\
                   \textrm{of $\rhos$}}} a_0(\rhos')
  +\sum_{\substack{\textrm{$\rhos'\!$ a shuffle}\\
                   \textrm{of $\rhos$}}} a_0(\rhos').
  \]
\end{lemma}

\begin{proof}
  This follows directly from the definitions.
\end{proof}

See also Lemma~\ref{lem:composable-gr} for alternate
characterizations of composability, splittings, and shuffles in terms
of a grading on the algebra.

\begin{remark}
  \label{rmk:alg-differential-reps-shuffles-and-joins}
  In the type~$A$ module, these two types of
  terms in the boundary, splittings and shuffles, appear respectively as
  join curve ends and odd shuffle curve ends,
  as defined in Definition~\ref{def:ends-moduli};
  see the proof of Proposition~\ref{prop:A-module-defined}.
\end{remark}

\section{Matched circles and their algebras}
\label{sec:matched-circles}

We now turn to the algebra associated to a parametrized surface.
We represent our boundary surfaces by \emph{matched circles}.  

\begin{definition}
	\label{def:PointedMatchedCircle}
A \emph{matched circle} is a triple $(Z, \CircPts, M)$ of an
\index{matched circle}%
\index{circle, matched|see{matched circle}}%
\index{matched circle!pointed|see{pointed matched circle}}%
oriented circle $Z$,
$4k$~points $\gls*{CircPts}=\{a_1,\dots,a_{4k}\}$ in $Z$ (as in
Section~\ref{sec:reeb-chords-def}), and a \emph{matching}: a $2$-to-$1$
\index{matching}%
function $\gls*{Matching}\co \CircPts\to[2k]$.
We require that performing
surgery along the $2k$ pairs
of points in~$Z$ (viewed as $0$-spheres) yields a single circle, rather than
several disjoint circles.  

A \emph{pointed matched circle}
\index{pointed matched circle}%
$\gls*{PMC}$ is a matched circle together with a basepoint
$\gls*{pmcz}\in Z\setminus \CircPts$. 
\end{definition}

\begin{definition}\label{def:pmc-to-from-surf}
  To each matched circle $\PtdMatchCirc$ we associate an oriented
  surface $\gls*{SurfPMC}$,  
  \index{pointed matched circle!surface associated to}%
  \index{surface associated to pointed matched circle}%
  or just $F$, of genus~$k$, by taking a disk with boundary~$Z$
  (respecting the orientation), attaching oriented $2$-dimensional
  $1$-handles along the pairs specified by $M$, and filling the
  resulting boundary circle with another disk.  Alternately, $F$ may
  be obtained from a polygon with $4k$ sides, with the
  points~$\CircPts$ in the middle of the sides, by gluing pairs of
  sides specified by~$M$ in an orientation\hyp preserving way.  The
  last disk glued in is a neighborhood of the unique resulting
  vertex. We call $F(\PMC)$ the \emph{surface associated to~$\PMC$}.
\end{definition}

If $-\PMC$ denotes the orientation reverse of $\PMC$ (i.e., the same
data except with the orientation on $Z$ reversed) then $F(-\PMC)$ is
the orientation reverse $-F(\PMC)$ of $F(\PMC)$,
canonically.\index{orientation reversal}%
\index{pointed matched circle!orientation reverse}

In the other direction from Definition~\ref{def:pmc-to-from-surf}, we
have the following:

\begin{definition}
  \label{def:PMCCompatibleMorse}
  Given a closed, orientable surface $F$ of genus $k$, let
  $\gls*{MorseFn}\co F\to \RR$ be a Morse function with a unique index
  $0$ critical point and a unique index $2$ critical point, and such
  that $f(p_i)=1$ for all index $1$ critical points
  $p_1,\dots,p_{2k}$. Suppose that $f^{-1}(3/2)$ is non-empty.  Fix
  also a Riemannian metric $\gls*{Metric}$ on $F$. Let $Z=f^{-1}(3/2)$,
  and let $\CircPts\subset Z$ denote the ascending spheres of the
  index $1$ critical points. Define $M(a_j)=i$ if $a_j$ is in the
  ascending sphere of $p_i$. Choose also a point $z\in
  Z\setminus\CircPts$. Then $\PMC=(Z,\CircPts,M,z)$ is a pointed
  matched circle, and $F(\PMC)\cong F$. We say that $(f,g)$ is
  \emph{compatible with the pointed matched circle $\PMC$}.
  \index{Morse function!compatible with pointed matched circle}%
\end{definition}

The points $\CircPts$ in a matched circle $(Z,\CircPts,M)$ inherit a cyclic
ordering from the
orientation of $Z$; in a pointed matched circle, this cyclic ordering is lifted to
an honest ordering, which we denote by~$\gls*{lessdot}$.

\begin{example}
  \label{ex:TorusMatchedCircle}
  Let $Z$ be a circle with four points, which we number $0,1,2,3$ in their cyclic ordering, and let
  $M\co \{1,\dots,4\}\to \{1,2\}$ be characterized by $M(i)\equiv i\pmod{2}$.
  Then $Z$ is a matched circle which represents the torus. Indeed, it
  is the unique matched circle 
  which represents the torus.
\end{example}

\begin{figure}
  \centering
  \includegraphics[scale=.8333]{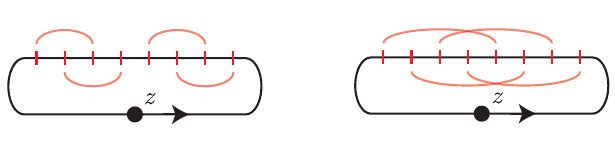}
  \caption[Pointed matched circles for the genus $2$ surface]{\textbf{Two pointed matched circles for the genus $2$ surface.} Left: the split pointed matched circle. Right: the antipodal pointed matched circle.}
  \label{fig:split-antipodal-pointed}
\end{figure}

\begin{example}
  Let $Z$ be a circle with eight marked points, which we number $\{0,\dots,7\}$ in their cyclic ordering.
  Let $M\co \{1,\dots,8\}\to\{1,\dots,4\}$ be characterized by $M(i)\equiv i\pmod{4}$.
  Then $Z$ is a matched circle which represents a surface of
  genus~$2$, called the \emph{antipodal} pointed matched circle. A
  different matched circle
  representing the same surface can be obtained by forming the connected sum (in the obvious way) of two 
  pointed matched circles for the torus; this is called the
  \emph{split} pointed matched circle.  See
  Figure~\ref{fig:split-antipodal-pointed}.
\end{example}

\begin{figure}
\includegraphics[scale=.8333]{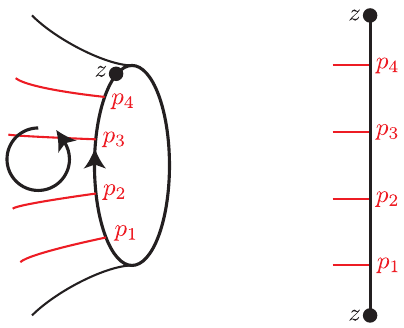}
\caption[Orientation of a pointed matched circle]{\textbf{Orientation of a pointed matched circle.} Left: the orientation of a Heegaard surface, and
  the induced orientation of its boundary. Right: the induced ordering of $\CircPts$: in
  this figure, $p_1\lessdot p_2\lessdot p_3\lessdot p_4$.}\label{fig:OrientationOrdering}
\end{figure}

\begin{definition}
  \label{def:AlgPMC}
  Fix a pointed matched circle $\PMC$.  For each subset $\SetS$ of
  $[2k]$ we have one idempotent $\gls*{IdemPMCS}$ in the strands algebra $\Alg(4k)$,
  given by 
  \begin{equation}\label{eq:idempotent-def-2}
    I(\SetS) \coloneqq \!\!\!\sum_{\substack{S\subset [2k]\\
        \text{$S$ a section of $M$ over $\SetS$}}}\!\!\!
    I(S),
  \end{equation}
  where by a \emph{section} of $M$ over $\SetS$ we mean a subset of $S$ that
  maps bijectively to $\SetS$ via~$M$.
  \index{section of matching}%
  That is, each element of~$\SetS$ defines a pair of points in~$\CircPts$; $I(\SetS)$
  is the sum over the $2^{|\SetS|}$ minimal idempotents in $\Alg(4k,|\SetS|)$
  that pick one point out of each of these pairs.  
  The {\em ring of idempotents associated to the pointed matched circle}
  is the algebra $\gls*{IdemPMC}$ generated by the elements $I(\SetS)$
  associated to all the subsets $\SetS\subset [2k]$.
  \index{idempotents!of $\Alg(\PMC)$}%
  Define
  \[
  \gls*{Unit}
  \coloneqq \sum_{\SetS \subset [2k]} I(\SetS)
  \]
  to be the unit in $\Idem(\PtdMatchCirc)$.
\end{definition}

We will now define $\Alg(\PtdMatchCirc)$, our
main algebra of interest.

\begin{definition}
  \label{def:arhos}
  For a given pointed matched circle~$\PMC$, the
  \emph{algebra associated to~$\PMC$}, denoted 
  $\gls*{AlgPMC}$,
  \index{algebra!associated to pointed matched circle}%
  is the subalgebra of $\bigoplus_{i=0}^{2k}\Alg(4k,i)$ generated (as
  an algebra) by $\Idem(\PtdMatchCirc)$ and by $\Unit a_0(\rhos)\Unit$ for every
  set of Reeb chords~$\rhos$.
  The {\em algebra element associated to $\rhos$}, denoted
  $\gls*{aofrhos}$, is the projection $\Unit \cdot a_0(\rhos)\cdot \Unit$ of $a_0(\rhos)$ to
  $\Alg(\PtdMatchCirc)$.
  
  Set
  \begin{align*}
    \gls*{AlgPMCi}&=\Alg(\PtdMatchCirc)\cap\Alg(4k,k+i)\\
    \gls*{IdemPMCi}&=\Idem(\PtdMatchCirc)\cap\Alg(4k,k+i),
  \end{align*}
  and call $\Alg(\PtdMatchCirc,i)$ the part of $\Alg(\PMC)$ with \emph{weight}
  $i$.
  \index{weight of algebra element}%
\end{definition}

It is easy to verify that $\Alg(\PtdMatchCirc)$ is closed under
multiplication and $\partial$.

In Chapter~\ref{chap:ainfinity}, we considered algebras
(and modules) defined over a ground ring~$\Ground$. From now on, we
will view $\Alg(\PMC)$ as defined over the ground
ring~$\Ground=\Idem(\PMC)$.  

The algebra
$\Alg(\PtdMatchCirc)$ decomposes (as a direct sum of differential
graded algebras) according to the weights of the algebra elements,
\[
\Alg(\PtdMatchCirc)=\bigoplus_{i=-k}^k\Alg(\PtdMatchCirc,i).
\]

There is a convenient basis for $\Alg(\PMC)$ over~$\Field$,
consisting of the non-zero
elements of
the form $I(\SetS)a(\rhos)$, where $\SetS\subset [2k]$
and $\rhos$ is a
consistent set of Reeb chords. 
In two-line notation, we express these \emph{basic generators}
\index{basic generator}%
$I(\SetS)a(\rhos)$ as
\[\honestalga{x_1&\cdots&x_k}{y_1&\cdots&y_k}{z_1&\cdots&z_l}
\coloneqq I(\{M(x_1),\dots,M(x_k),M(z_1),\dots,M(z_l)\}
a\bigl(\reebchords{x_1&\cdots&x_k}{y_1&\cdots&y_k}\bigr).
\]
Graphically, we denote the matching~$M$ by drawing dotted arcs on
the left and right of a strand diagram.  We draw the strands in
an element $\Alg(\PtdMatchCirc)$ as we did for elements of $\Alg(n)$ in
Section~\ref{sec:algebra-geometric},
except that we draw horizontal strands
(which, if we expand out the sum for $I(\SetS)$, appear in half the
terms) with dashed lines.  For instance, if $\PMC$ is the
split genus-two pointed matched circle, one basic generator is
\begin{equation}
\label{eq:genus-two-strands}
\mfigb{strands-30} \;=\;
  \mfigb{strands-33}\; + \;\mfigb{strands-34}
\end{equation}
which, in two-line notation, can be denoted
$\honestalg{1}{6}{2}$ or \gls*{honestnotation}. 
In particular, note
that in a non-zero basic element of $\Alg(\PMC)$, no initial
(respectively terminal) ends of strands can be matched with each other.

\begin{example}
  Let $\PMC$ be the unique pointed matched circle representing the
  torus. The algebra $\Alg(\PMC,0)$ is described in detail in
  Section~\ref{sec:torus-algebra}. The algebra $\Alg(\PMC,-1)$ is
  isomorphic to $\Field$, with generator $(\emptyset, \emptyset,\phi)$
  where $\phi$ is the unique bijection from $\emptyset$ to $\emptyset$.
  The algebra $\Alg(\PMC,1)$ has $\Field$-basis
  \[
  \bigl\lbracket
\begin{smallmatrix} 1 & 2\\ & \end{smallmatrix}\bigr\rbracket,
  \ \honestalg{1}{3}{2},\ \honestalg{1 & 2}{2 &
    3}{},\ \honestalg{2}{4}{1},\ \honestalg{2 & 3}{3 & 4}{},\ \honestalg{1
    & 2}{4 & 3}{},\ \honestalg{1 &2}{3 & 4}{}.
  \]
  The non-trivial differentials are given by 
  \[
  \bdy \honestalg{1}{3}{2}=\honestalg{1 & 2}{2 &
    3}{}, \qquad
  \bdy\honestalg{2}{4}{1}=\honestalg{2 & 3}{3 & 4}{},\qquad
  \bdy\honestalg{1
    & 2}{4 & 3}{}=\honestalg{1 &2}{3 & 4}{}.
  \]
  In particular, the homology $H_*(\Alg(\PMC,1))$ is $\Field$,
  generated by the idempotent. Therefore, the algebra $\Alg(\PMC,1)$ is formal
  (quasi-isomorphic to its homology).
  \index{formal \dg algebra}%

  It is not a coincidence that $H_*(\Alg(\PMC,1))\cong
  H_*(\Alg(\PMC,-1))$. In general, for a pointed matched circle $\PMC$
  for a surface of genus $k$, $\Alg(\PMC,-k)$ is $1$-dimensional while
  $\Alg(\PMC,k)$ is formal and has $1$-dimensional homology. (The
  first fact is trivial; the second takes a little work.) More
  generally, there is a duality relating $\Alg(\PMC,i)$ and
  $\Alg(\PMC',-i)$ for appropriate pairs of pointed matched circles
  $\PMC$ and $\PMC'$; see~\cite{LOTHomPair}. The fact that the algebra
  is one-dimensional in the extremal weight relates to the fact that
  $\HFKa$ is monic for fibered knots; see~\cite[Theorem 7]{LOT2} for further
  discussion.
\end{example}

Finally, we give interpretations of $\Alg(\PMC)$ directly in terms of
Reeb chords.
For any pointed matched circle~$\PMC$ and set of Reeb
chords~$\rhos$, if $a_0(\rhos)\neq 0$, then
$a(\rhos)\neq 0$ if and only if
$\abs{M(\rhos^-)} = \abs{M(\rhos^+)} = \abs{\rhos}$, i.e., if and only
if no two $\rho_i$ have matched initial or final endpoints.


\begin{lemma}
  \label{lem:VanishingProduct}
  Let $\rhos$ and $\sigmas$ be a pair of sets of Reeb chords so that
  $a_0(\rhos)\cdot a_0(\sigmas)\ne 0$. Then
  $a(\rhos)\cdot a(\sigmas)=0$ if and only if there is a matched pair
  $\{p,q\}$ in~$\PMC$ satisfying one of the following conditions:
  \begin{enumerate}
    \item \label{item:StartMatched}
      $\rhos\uplus\sigmas$ has a chord starting at $p$ and another one
      starting at $q$;
    \item \label{item:TerminateMatched}
      $\rhos\uplus\sigmas$ has a chord terminating at $p$ and another
      one terminating at $q$; or
    \item \label{item:BadProduct}
      $\rhos$ has a chord which terminates in~$p$
      while $\sigmas$ has a chord which starts in~$q$.
  \end{enumerate}
  Further, if none of
  Conditions~(\ref{item:StartMatched})--(\ref{item:BadProduct}) are
  satisfied then for any $\SetS$ such that $I(\SetS)\cdot a(\rhos)\neq 0$
  and $I(\SetS)\cdot a(\rhos\uplus\sigmas)\neq 0$ it follows that
  $I(\SetS)\cdot a(\rhos)\cdot a(\sigmas)\neq 0$.
\end{lemma}

\begin{proof}
  Any of the three conditions implies the vanishing
  of $a(\rhos)\cdot a(\sigmas)$:
  in view of Lemma~\ref{lem:reeb-product},
  the first is the case where $\Unit \cdot a_0(\rhos)\cdot a_0(\sigmas)=0$,
  the second is where $a_0(\rhos)\cdot a_0(\sigmas)\cdot\Unit =0$,
  and the third is where $a_0(\rhos)\cdot\Unit\cdot a_0(\sigmas) =0$.

  Conversely, suppose none of
  Conditions~(\ref{item:StartMatched})--(\ref{item:BadProduct})
  hold. Let $S\coloneqq (\rhos\uplus\sigmas)^-$ and $T\coloneqq
  (\rhos\uplus\sigmas)^+$.  The failure of
  Condition~(\ref{item:StartMatched})
  (respectively~(\ref{item:TerminateMatched})) implies that the
  restriction of $M$ to~$S$ (respectively~$T$) is injective, so $I(S)$
  (respectively $I(T)$) is a summand in $\Unit$.  Let $U\coloneqq
  (S\setminus \rhos^-)\cup\rhos^+$; observe that
  $I(S)a_0(\rhos)=I(S)a_0(\rhos)I(U)$. It follows from
  Lemma~\ref{lem:reeb-product} (and the definitions) that
  $U=(T\setminus \sigmas^+)\cup \sigmas^-$.  By the failure of
  Condition~(\ref{item:BadProduct}), the restriction of $M$ to~$U$ is
  injective, so $U$ is a summand of~$\Unit$.  Since $I(S)\cdot
  a_0(\rhos)\cdot I(U)\cdot a_0(\sigmas)\cdot I(T)\neq 0$,
  it follows that $a(\rhos)\cdot a(\sigmas)\neq 0$.

  For the last part of the statement, since $I(\SetS)\cdot
  a(\rhos\uplus\sigmas)\neq0$ there is a subset $S\subset[4k]$ such
  that $M(S)=\SetS$, $M|_S$ is injective and
  $(\rhos\uplus\sigmas)^-\subset S$. Moreover, by the failure of
  Condition~(\ref{item:BadProduct}), $U\coloneqq (S\setminus \rhos^-)\cup
  \rhos^+$ contains $\sigmas^-$, and since $I(\SetS)\cdot a(\rhos)\neq
  0$, $M|_{U}$ is injective. Similarly,
  $T\coloneqq (U\setminus \sigmas^-)\cup\sigmas^+=\bigl(S\setminus
  (\rhos\uplus\sigmas)^-\bigr)\cup (\rhos\uplus\sigmas)^+$ has $M|_T$
  injective since $I(\SetS)\cdot a(\rhos\uplus\sigmas)\neq 0$.
  So, setting $\SetT\coloneqq M(T)$ and $\SetU\coloneqq
  M(U)$, we have
  \[
  I(\SetS)\cdot a_0(\rhos)\cdot I(\SetU)\cdot a_0(\sigmas)\cdot I(\SetT)\neq 0,
  \]
  as desired.
\end{proof}

\begin{remark}
  One might wonder to what extent $\Alg(\PMC)$ is a topological invariant
  of the surface specified by $\PMC$. 
  In fact, if $\PMC$ and $\PMC'$ are two different pointed matched
  circles representing diffeomorphic surfaces, then neither
  $\Alg(\PMC)$ and $\Alg(\PMC')$ nor their homologies $H_*(\Alg(\PMC)$
  and $H_*(\Alg(\PMC'))$ are necessarily isomorphic; indeed, they need
  not have the same rank.  (See
  \cite[Theorem~\ref*{LOT2:thm:halg-support}]{LOT2}.)
  However, it turns out that the module categories
  of $\Alg(\PMC)$ and $\Alg(\PMC')$ are derived
  equivalent \cite[Theorem~\ref*{LOT2:thm:AlgebraDependsOnSurface}]{LOT2}.
\end{remark}

\begin{remark}
The summands $\Alg(\PtdMatchCirc,i)\subset \Alg(\PtdMatchCirc)$ for
$i\neq 0$ act trivially on the modules $\CFDa(Y)$ and $\CFAa(Y)$
defined in Chapters~\ref{chap:type-d-mod} and~\ref{chap:type-a-mod}. So,
in the present work, only the subalgebra
$\Alg(\PtdMatchCirc,0)\subset\Alg(\PtdMatchCirc)$ will be of
interest. The other summands of $\Alg(\PtdMatchCirc)$ do play an
important role in the case of three-manifolds with disconnected boundary; see~\cite{LOT2},
and also Appendix~\ref{app:Bimodules}.
\end{remark}

\begin{remark}
  One can think of the index $i$ in $\Alg(\PtdMatchCirc,i)$ as
  corresponding to $\SpinC$-structures on the surface
  $F(\PtdMatchCirc)$, compare~\cite{LOT2}. For a $3$-manifold $Y$ with
  connected boundary $\bdy Y$, only the middle $\SpinC$-structure on
  $\bdy Y$ extends over $Y$.
\end{remark}

\begin{remark}
  We get pointed matched circles by splitting a pointed Heegaard
  diagram in two along a circle~$Z$, satisfying a suitable
  combinatorial hypothesis; see Lemma~\ref{lem:cut-hd-cut-Y}.
  The points~$\CircPts$ are the
  intersections of $Z$ with the $\alpha$-circles and the matching~$M$ records
  which intersection points are on the same $\alpha$-circle.
  We declare one of the two components of the complement of $Z$ in the Heegaard
  surface to be the ``left''
  (west, or type~$A$)
  side, and let $Z$ inherit the orientation as the boundary of this side.
  See Figure~\ref{fig:OrientationOrdering}. 
\end{remark}
\section{Gradings}
\label{sec:gradings-algebra}

Although the algebra $\Alg(n,k)$ has a $\ZZ$-grading, for instance by the
crossing number (i.e., number of inversions), it turns out that the
subalgebra
$\Alg(\PtdMatchCirc)$ has no $\ZZ$-grading.

\begin{example}
For $\PMC$ the pointed matched circle representing a torus, consider the two algebra
elements in $\Alg(\PMC,1)$
\[
x = \honestalg{1}{3}{2} = \mfigb{strands-40}
\qquad\text{and}\qquad
y = \honestalg{2}{4}{1} = \mfigb{strands-41}.
\]
On one hand, we have
\begin{align*}
y\cdot x &= \mfigb{strands-42}
\intertext{but on the other hand,}
  \bdy x &= \mfigb{strands-43}\\[2pt]
  (\bdy x) \cdot y &= \mfigb{strands-44}.
\end{align*}
Hence $\bdy((\bdy x) \cdot y) = y \cdot x,$ which is inconsistent
with any $\ZZ$-grading in which $x$ and~$y$ are homogeneous.
\end{example}

Instead, we grade $\Alg(\PtdMatchCirc)$ by a non-commutative group.
A 
grading of a differential algebra~$A$ by a
non-commutative group~$G$ is a straightforward, but not very standard, notion.
Briefly, it is a grading of $A$ as an algebra by~$G$ together with a
distinguished central element $\lambda \in G$ so that for any
homogeneous elements $a,b \in A$, we have $\gr(a\cdot b)=\gr(a)\gr(b)$ and $\gr(\partial a)
=\lambda^{-1}\gr(a)$.  See Section~\ref{sec:grad-non-comm} for the
details of the algebraic setting.

In fact, we give two gradings of $\Alg(\PMC)$ by non-commutative
groups. In Section~\ref{sec:pregrading} we construct a canonical grading of
$\Alg(\PMC)$ by a group $\bigGroup(n)$. In
Section~\ref{sec:refined-grading} we refine this grading to a
grading by a subgroup $\smallGroup(\PMC)\subset\bigGroup(n)$. The
subgroup $\smallGroup(\PMC)$ has a more obvious topological
interpretation and behaves better with respect to the operation of
gluing Heegaard diagrams (see Section~\ref{sec:GradedPairingThm}),
but the definition of the refined grading on $\Alg(\PMC)$ depends on
some additional choices.

\subsection{Unrefined grading on \textalt{\protect{$\Alg(\PMC)$}}{A(Z)}}
\label{sec:pregrading}
\index{grading!on $\Alg(\PMC)$|(}%

We start by giving a grading~$\grb$ of $\Alg(n,k)$ with
values in a non-commutative group $\bigGroup(n)$. We then show that
this grading has the property that $\grb$
is unchanged by adding horizontal strands, and so descends to a
grading on $\Alg(\PMC)$.  

Let $Z' = Z \setminus
\{z\}$.  The group $\gls*{bigGroup}$ can be abstractly
defined as a central extension by $\ZZ$ of the relative
homology group $H_1(Z',\CircPts)$ (which is isomorphic to
$\ZZ^{n-1}$).
\index{grading!group!big}%
In particular, there is an exact sequence of groups
\[
\ZZ \overset{\lambda}{\hookrightarrow} \bigGroup(n)
   \overset{[\cdot]}{\twoheadrightarrow} H_1(Z',\CircPts)
\]
where in the first map $1$ maps to $\lambda$.

To define the central extension concretely we make some
more definitions.
For $p \in \CircPts$ and $\alpha \in H_1(Z',\CircPts)$, define the
\emph{multiplicity} $\gls*{multalphap}$ of $p$ in~$\alpha$ 
\index{multiplicity of point in homology class}%
to be the average
multiplicity with which $\alpha$ covers the regions on either side
of~$p$.  
Extend $m$ to a map $H_1(Z',\CircPts) \times H_0(\CircPts) \to
\OneHalf\ZZ$ bilinearly.

For $\alpha_1, \alpha_2 \in H_1(Z',\CircPts)$, define
\begin{equation}\label{eq:L-def}
\gls*{Linkingalpha}\coloneqq m(\alpha_2, \bdy\alpha_1) 
\end{equation}
where $\gls*{connectinghom}$ is the connecting homomorphism $\bdy\co
H_1(Z',\CircPts) \to H_0(\CircPts)$ from the long exact sequence for a
pair.  The map $L$ can be thought of as
the linking number of $\bdy\alpha_1$ and $\bdy\alpha_2$.  Also, for
$\alpha \in H_1(Z',\CircPts)$, define
$\epsilon(\alpha) \in (\OneHalf \ZZ)/\ZZ$ by
\glsadd{epsilonalpha}%
\begin{align*}
  \epsilon(\alpha)
    &\coloneqq \OneQuart \#(\text{parity changes in $\alpha$}) \pmod{1}.
\end{align*}
Here, a parity change in $\alpha$ is a point $p \in \CircPts$ so
that the multiplicity of~$\alpha$ to the left of~$p$ has different
parity than the multiplicity of~$\alpha$ to the right of~$p$, or,
equivalently, so that $m(\alpha,p)$ is a half-integer.  (Note that
there are always an even number of such points.)

\begin{definition}\label{def:big-gr-group}
  Define $\bigGroup(n)$ to be the group generated by pairs
  $\gls*{jalpha}$, where $j \in \OneHalf\ZZ$, $\alpha\in
  H_1(Z',\CircPts)$ and $j \equiv \epsilon(\alpha) \bmod 1$, with
  multiplication given by
  \begin{equation}\label{eq:Gn-mult-def}
    (j_1,\alpha_1)\cdot(j_2,\alpha_2) \coloneqq (j_1 + j_2 +
    L(\alpha_1, \alpha_2), \alpha_1 + \alpha_2).
  \end{equation}
  (Thus the map $L$ is a 2-cocycle on $H_1(Z',\CircPts)$ which defines
  the central extension~$\bigGroup(n)$.)  In an element $g = (j,
  \alpha)$ of $\bigGroup(n)$, we refer to $j$ as the \emph{Maslov
    component} of~$g$ and $\alpha$ as the \emph{$\SpinC$ component}
  of~$g$.  
  \index{Maslov|see{grading, Maslov}}%
  \index{spinc component of grading@$\SpinC$ component of grading}%
  \index{grading!Maslov component}%
  \index{grading!$\SpinC$ component}%
  The distinguished central element of $\bigGroup(n)$ is $\gls*{lambda} =
  (1,0)$. 
\end{definition}
To see that $\bigGroup(n)$ is, in fact, a group we make $L$ more explicit and
prove some elementary properties.

\begin{lemma}\label{lem:L-determinant}
  If $\alpha = \sum_{i=1}^{n-1} u_i[p_i, p_{i+1}]$ and $\beta = \sum_{i=1}^{n-1}
  v_i[p_i,p_{i+1}]$, then 
  \[
  L(\alpha,\beta) = \sum_{i=1}^{n-2} \frac{u_iv_{i+1} - u_{i+1}v_i}{2}.
  \]
  In particular, $L(\alpha,\beta) = -L(\beta,\alpha)$.
\end{lemma}
Note that $u_iv_{i+1} - u_{i+1}v_i$ can be interpreted as the
determinant of a $2 \times 2$ submatrix of the matrix
$\bigl(\begin{smallmatrix} u_1&u_2&\cdots&u_{n-1}\\v_1&v_2&\cdots&v_{n-1}
\end{smallmatrix}\bigr)$.

\begin{proof}
  Recall that $\bdy[p_i,p_{i+1}] = p_{i+1}-p_i$, from which it follows
  that $\bdy\alpha = \sum_{i=1}^n (u_{i-1}-u_i)p_i$, where we set
  $u_0 = u_n = 0$.  Then
  \begin{align*}
    L(\alpha, \beta)
      &= \sum_{i=1}^n \frac{(u_{i-1}-u_i)(v_{i-1}+v_i)}{2}\\
      &= \sum_{i=1}^n \frac{u_{i-1}v_i - u_iv_{i-1}}{2} +
        \sum_{i=1}^n\biggl(\frac{u_{i-1}v_{i-1}}{2} - \frac{u_iv_i}{2}\biggr)\\
      &= \sum_{i=1}^n \frac{u_{i-1}v_i - u_iv_{i-1}}{2}
  \end{align*}
  since the second sum in the second line telescopes.  The first and last terms
  in the final sum vanish, proving the result.
\end{proof}

\begin{lemma}\label{lem:epsilon-add}
  If $\alpha, \beta \in H_1(Z',\CircPts)$, then
  \[
  \epsilon(\alpha+\beta) \equiv
    \epsilon(\alpha) + \epsilon(\beta) + L(\alpha,\beta) \pmod{1}.
  \]
\end{lemma}
\begin{proof}
  Consider the terms in $L(\alpha,\beta) = m(\beta,\bdy\alpha)$
  corresponding to the different points $p \in \CircPts$.  A
  particular $p$ contributes a half-integer to $L(\alpha,\beta)$ when
  $p$ appears an odd number of times in $\bdy\alpha$ and $m(\beta,p)$
  is a half-integer, i.e., $p$ appears an odd number of times in
  $\bdy\beta$.  This case (when $p$ is odd in both $\bdy\alpha$ and
  $\bdy\beta$) is exactly when the contribution of $p$ to
  $\epsilon(\alpha+\beta) - \epsilon(\alpha) - \epsilon(\beta)$
  is $1/2 \pmod 1$.
\end{proof}

\begin{proposition}\label{prop:G-n-group}
  $\bigGroup(n)$ is a group.
\end{proposition}
\begin{proof}
  Lemma~\ref{lem:epsilon-add} guarantees that if $(j_1,\alpha_1)$ and
  $(j_2,\alpha_2)$ are in $\bigGroup(n)$, then the product is as well:
  \begin{align*}
    j_1 + j_2 + L(\alpha_1,\alpha_2) &\equiv \epsilon(\alpha_1) +
      \epsilon(\alpha_2) + L(\alpha_1,\alpha_2) \pmod{1}\\
    &\equiv \epsilon(\alpha_1 + \alpha_2) \pmod{1}.
  \end{align*}
  The inverse of $(j,\alpha)$ is easily seen to be $(-j,-\alpha)$.  We
  must also check that the product is associative:
\begin{multline*}
  ((j_1,\alpha_1)\cdot(j_2,\alpha_2))\cdot (j_3,\alpha_3) \\
    \begin{aligned}
      &=(j_1 + j_2 + j_3 + L(\alpha_1,\alpha_2) +
      L(\alpha_1+\alpha_2, \alpha_3), \alpha_1 + \alpha_2 +
      \alpha_3)\\
      &=(j_1 + j_2 + j_3 + L(\alpha_1,\alpha_2) +
      L(\alpha_1,\alpha_3) + L(\alpha_2,\alpha_3), \alpha_1 + \alpha_2 +
      \alpha_3)\\
    \end{aligned}
\end{multline*}
  where we use bilinearity of $L$.  Associating in the opposite order
  gives the same result.
\end{proof}

\begin{definition}\label{def:grading-algebra}
The grading $\grb(a)$ of an element $a\in\Alg(n,k)$
with starting idempotent~$S$ and ending idempotent $T$ is defined as
follows. 
The $\SpinC$ ($H_1(Z',\CircPts)$) component of the grading of an element $a = (S,T,\phi)$,
denoted~$[a]$,
is defined to  be the
element of the relative homology group $H_1(Z',\CircPts)$ given by summing up
the intervals corresponding to the strands:
\begin{align*}
  \gls*{homola} &\coloneqq \sum_{s \in S} [s,\phi(s)]. 
  \intertext{The Maslov component of the grading is given by}
  \gls*{iota} &\coloneqq \inv(a) - m([a], S).
  \intertext{So, the grading of an element $a$ is given by}
  \gls*{grprime}(a) &\coloneqq (\iota(a), [a]).
\end{align*}
Note that $\bdy[a] = T - S \in H_0(\CircPts)$, so by
Lemma~\ref{lem:L-determinant}, $m([a],S) = m([a],T)$.
\end{definition}

\begin{proposition}\label{prop:grb-is-grading}
  The function $\grb$ defines a grading (in the sense of
  Definition~\ref{def:dga-nc-grading}) on $\Alg(n,k)$, with $\lambda = (1,0)$.
\end{proposition}
\begin{proof}
  Let $a$ be an algebra element with starting idempotent~$S$ and
  ending idempotent~$T$.
  We first check that $\grb(a) \in \bigGroup(n)$.  Modulo~$1$, we have
\begin{align*}
  \iota(a) &\equiv m([a],S) \\
    &\equiv \OneHalf \#(S \setminus T) \equiv \OneHalf \#(T \setminus S)\\
    &\equiv \OneQuart \#(S \symmdiff T) \equiv \epsilon([a]),
\end{align*}
\glsadd{SymmDiff}%
where $S \symmdiff T$ is the symmetric difference of $S$ and~$T$:
\[
S\symmdiff T\coloneqq (S \setminus T) \cup (T \setminus S).
\]

For the behavior of $\grb$ under the action of $\partial$, observe that
$\iota(\partial a) = \iota(a) - 1$: neither $S$ nor $[a]$ is
changed by~$\partial$, so $m([a], S)$ is unchanged, and in the
definition of~$\iota$ only the
$\inv(a)$ term changes.  Thus with $\lambda
= (1, 0)$, we have $\grb(\partial a) = \lambda^{-1}\grb(a)$,
as desired.

We also check that $\grb$ is compatible with
multiplication: if $a$ is given by $(S,T,\phi_1)$ and $b$ is given by
$(T, U, \phi_2)$ with $a \cdot b \ne 0$, then
\begin{align*}
  [a\cdot b] &= [a] + [b]\\
  \iota(a\cdot b) &= \inv(a\cdot b) - m([a] + [b], S)\\
    &= \inv(a) + \inv(b) - m([a], S) - m([b], T) + m([b], T - S)\\
    &= \iota(a) + \iota(b) + m([b], \bdy[a])
\end{align*}
and so $\grb(a\cdot b) = \grb(a)\cdot\grb(b)$.
\end{proof}

\begin{proposition}\label{prop:grading-descends}
  For a consistent set~$\rhos$ of Reeb chords, $a_0(\rhos)$ is
  homogeneous with respect to $\grb$.  In particular,
  the grading $\grb$ descends to a grading on the subalgebra
  $\Alg(\PtdMatchCirc)$.
\end{proposition}
\begin{proof}
  Observe that adding a horizontal strand at $p \in \CircPts$ to an algebra
  element $a = (S, T, \phi)$ does not change $[a]$ and changes
  $\inv(a)$ and $m([a], S)$ by the same amount, $m([a], p)$.  Thus
  $\grb(a)$ is unchanged by adding a horizontal strand (moving from
  $\Alg(n,k)$ to $\Alg(n,k+1)$).  In particular the terms in the
  definition of $a_0(\rhos)$ all have the same grading, proving the
  first statement.  Since
  $\Alg(\PtdMatchCirc)$ is generated by elements of the form
  $I(\SetS)a_0(\rhos)$, and $I(\SetS)$ is homogeneous with grading $(0,0)$, the second
  statement follows.
\end{proof}

As in Section~\ref{sec:grad-non-comm}, there is a partial order on $\bigGroup(n)$ by $(j,\alpha)\gls*{lessthan}(j',\alpha')$ if $\alpha=\alpha'$ and $j<j'$.

\index{grading!on $\Alg(\PMC)$|)}
\subsection{Refined grading on \textalt{\protect{$\Alg(\PtdMatchCirc,i)$}}{A(Z,i)}}
\label{sec:refined-grading}\index{grading!on $\Alg(\PMC)$!refined|(}
We can also construct a refined grading $\gr$ on $\Alg(\PtdMatchCirc,i)$
with values in a smaller group $\smallGroup(\PtdMatchCirc)$.
\index{grading!group!small}%

The group \gls*{smallGroup} 
can be abstractly defined as a central
extension of $H_1(F) \cong \ZZ^{2k}$.  Again, there is a
sequence of groups
\begin{equation}\label{eq:smallGroup-exact-seq}
\ZZ \overset{\lambda}{\hookrightarrow} \smallGroup(\PtdMatchCirc)
  \overset{[\cdot]}{\twoheadrightarrow} H_1(F)
\end{equation}
where $1 \in \ZZ$ maps to $\lambda \in
\smallGroup(\PtdMatchCirc)$.
This central extension is defined as an abstract group as follows.
Let $\lambda$ denote the generator of the center of $\smallGroup(\PMC)$ and
for $g\in
\smallGroup(\PtdMatchCirc)$, let $[g]$ be its image in $H_1(F)$, as indicated
above.  Then, for $g, h \in \smallGroup(\PtdMatchCirc)$, we have
\begin{equation}
  \label{eq:commutators}
  gh = hg\lambda^{2([g] \cap [h])}.
\end{equation}
Equation~\eqref{eq:commutators} determines
$\smallGroup(\PtdMatchCirc)$ uniquely up to isomorphism.

\index{algebra!associated to pointed matched circle!grading on|see{grading}}%
In order to describe~$\gr$, first suppose that we have an
element~$\alpha$ of
$H_1(Z',\CircPts)$ for which
$M_*(\bdy \alpha) = 0$; that is, for each pair $p,q$ of points
in~$\CircPts$ identified by~$M$, the sum of the multiplicities of
$\bdy\alpha$ at~$p$ and at~$q$ is~$0$.  Then we can construct an
element $\gls*{Hofalpha}$ of $H_1(F(\PtdMatchCirc))$
from~$\alpha$ by embedding $Z$
in 
$F(\PtdMatchCirc)$ to make $\alpha$ a chain in
$F(\PtdMatchCirc)$, and then adding multiples of the cores of
the 1-handles we attached in the construction of
$F(\PtdMatchCirc)$ to form a cycle~$\widebar\alpha$.  Indeed, we
can identify $H_1(F)$ as $\ker (M_* \circ \bdy)$ inside
$H_1(Z',\CircPts)$.

With this setup, we define $G(\PMC)$ to be the subgroup of
$\bigGroup(4k)$ defined by
\[\{(j,\alpha)\in \bigGroup(4k)\mid M_*(\bdy(\alpha))=0\}.\]

\begin{lemma}
  $\smallGroup(\PMC)$ is a group that fits in the exact
  sequence~\eqref{eq:smallGroup-exact-seq} and satisfies
  equation~\eqref{eq:commutators}.
\end{lemma}
\begin{proof}
  The fact that $[\cdot]\co \bigGroup(4k) \to H_1(Z',\CircPts)$ is a
  homomorphism guarantees that $\smallGroup(\PMC)$ is a group.  The
  image of $G(\PMC)$ under $[\cdot]$ is, by definition,
  $\ker(M_* \circ \bdy) \simeq H_1(F)$ and the kernel of $[\cdot]$ is
  $\{(j,0) \mid j \in \ZZ\}$, so we have the desired exact sequence.
  Finally, for $\alpha,\beta \in \ker(M_* \circ \bdy)$, the
  corresponding elements $\overline{\alpha}, \overline{\beta}$ in
  $H_1(F)$ satisfy
  $$\overline{\alpha}\cap \overline{\beta} = m(\alpha, \beta),$$
  which implies Equation~\eqref{eq:commutators}.
\end{proof}

To define the grading $\gr$ on $\Alg(\PMC,i)$, pick an arbitrary base idempotent $I(\SetS_0) \in
\Alg(\PtdMatchCirc,i)$ and, for each other idempotent $I(\SetS)$, pick a
group element $\gls*{refinementdata} \in \bigGroup(4k)$ so that $M_*\bdy([\psi(\SetS)]) =
\SetS-\SetS_0$.
Then define
\begin{equation}\label{eq:small-grading}
  \gls*{smallgr}\bigl(I(\SetS)a(\rhos)I(\SetT)\bigr) \coloneqq \psi(\SetS)\grb(a(\rhos))\psi(\SetT)^{-1}.
\end{equation}
We will call the choices $\SetS_0$, $\psi(\SetS)$ the choice of
\emph{grading refinement data}.
\index{grading!refinement data}%
\index{refinement data|see{grading refinement data}}%

\begin{lemma}
  The function $\gr$ is a grading with values in $\smallGroup(\PMC) \subset \bigGroup(4k)$.
\end{lemma}

\begin{proof}
  To check the values are in $\smallGroup(\PMC)$, we compute:
\begin{align*}
  M_*\bdy [\gr(I(\SetS)a(\rhos)I(\SetT))] &= M_*\bdy[\psi(\SetS)] + M_*\bdy[a(\rhos)] - M_*\bdy[\psi(\SetT)]\\
    &= (\SetS-\SetS_0) + (\SetT - \SetS) - (\SetT-\SetS_0)\\
    &= 0.
\end{align*}
  It is straightforward to check that $\gr$ is a grading.
\end{proof}

\begin{remark}\label{rem:refined-grading-change}
\index{grading!refinement data!dependence on}%
  Different choices of grading refinement data give
  gradings on $\Alg(\PMC,i)$ that are conjugate in an
  idempotent-dependent way, in the following sense. For any two
  choices $(\SetS_0, \psi)$ and $(\SetS'_0, \psi')$ giving gradings
  $\gr_1$ and $\gr_2$ respectively, there is a
  function $\xi \co \Idem(\PMC,i) \to \smallGroup(\PMC)$ so that, if $a
  \in \Alg(\PMC)$ is homogeneous and satisfies $a = I(\SetS) a I(\SetT)$,
  \[
  \gr_2(a) = \xi(\SetS)\cdot\gr_1(a)\cdot\xi(\SetT)^{-1}.
  \]
  (The map $\xi$ is given by $\xi(\SetS) = \psi'(\SetS)\psi(\SetS)^{-1}$.)  As a
  result, the categories of $\Alg(\PMC,i)$-modules that are graded with
  respect to these different gradings are equivalent; see
  \cite[Proposition~\ref*{LOT2:prop:ChangeOfRefinementData}]{LOT2}
  for more details.
\end{remark}

\begin{remark}
  An alternate way to construct the refined grading, which we employ
  in Chapter~\ref{chap:TorusBoundary} for the case where $F$ is a
  torus, is to use a homomorphism from $\bigGroup(4k)$ to
  $\smallGroup(\PtdMatchCirc)$ that fixes $\smallGroup(\PtdMatchCirc)$
  as a subset of $\bigGroup(4k)$ (after extending scalars from $\ZZ$
  to~$\QQ$). For $k>1$, however, there is no homomorphism
  $\bigGroup(4k)\to\smallGroup(\PtdMatchCirc)$ fixing $\smallGroup(\PtdMatchCirc)$.
\end{remark}

\begin{remark}\label{rmk:geom-gradings}
  $\smallGroup(\PtdMatchCirc)$ has a topological
  interpretation as the group of homotopy classes of non-vanishing
  vector fields on $F \times [0,1]$ with fixed behavior on the
  boundary.  Precisely, fix a vector field $V_0$ on $F \times \RR$
  that is invariant by translation in $\RR$.  Then the homotopy
  classes of vector fields on $F \times [0,1]$ that agree with $V_0$
  on neighborhoods of $F \times \{0\}$ and $F \times \{1\}$ form a
  group by concatenation; this group is isomorphic to
  $\smallGroup(\PtdMatchCirc)$.  This is unsurprising in light of the fact that
  Seiberg-Witten Floer homology on $Y$ is graded by non-vanishing
  vector fields on~$Y$ \cite[Section~28]{KronheimerMrowka}.
\end{remark}
\index{grading!on $\Alg(\PMC)$!refined|)}


\chapter{Bordered Heegaard diagrams}
\label{chap:heegaard-diagrams-boundary}

Just as closed $3$-manifolds can be represented by
Heegaard diagrams, bordered $3$\hyp manifolds can be described by a
suitable type of Heegaard diagrams, which we call \emph{bordered Heegaard
diagrams}. Just as Heegaard Floer homology of closed $3$-manifolds is
associated to closed Heegaard diagrams, bordered Floer invariants of
$3$-manifolds are associated to bordered Heegaard diagrams. This
chapter is devoted to defining bordered Heegaard diagrams and
discussing their combinatorial / topological properties.

In Section~\ref{sec:BorderedDiagrams}, we define bordered $3$-manifolds
and bordered Heegaard diagrams (Definitions~\ref{def:bordered-3-mfld}
and~\ref{def:BorderedDiagram}, respectively), and discuss
the relevant existence and uniqueness properties for bordered diagrams. In
Section~\ref{sec:Bestiary}, we give some examples of bordered Heegaard
diagrams, to give the reader a concrete picture. In
Section~\ref{sec:homology-classes-generators}, we describe the
generators of the bordered Heegaard Floer complex, and
discuss the $\SpinC$ structures associated to generators, and the
homology classes connecting them. As in the closed case, one cannot
use an arbitrary bordered Heegaard diagram to define bordered Heegaard
Floer homology: rather, one needs to use bordered diagrams satisfying
certain combinatorial conditions ensuring that the sums encountered in
the differentials are finite. In
Section~\ref{sec:Admissibility}, we describe these criteria.
Finally, in Section~\ref{sec:GluingDiagrams}, we discuss how to
glue bordered Heegaard diagrams and their generators to obtain
Heegaard diagrams and generators for closed manifolds.

\section{Bordered Heegaard diagrams: definition, existence, and uniqueness}
\label{sec:BorderedDiagrams}
We recall first the diagrams relevant for Heegaard Floer homology
of closed three-manifolds.

\begin{definition}
A \emph{pointed Heegaard diagram} $\gls*{Heegaard}$ 
\index{Heegaard diagram!pointed}%
\index{pointed Heegaard diagram|see{Heegaard diagram, pointed}}%
\index{diagram, Heegaard|see{Heegaard diagram}}%
is a quadruple
$(\Sigma,\alphas,\betas,\gls*{z})$ consisting of
\begin{itemize}
\item a compact, oriented surface $\Sigma$ with no boundary, of
  some genus~$g$;
\item two $g$-tuples $\alphas = \{\alpha_1,\dots,\alpha_g\}$ and $\gls*{betas} =
  \{\beta_1,\dots,\beta_g\}$ of circles in
  $\Sigma$, with the circles in each tuple disjoint; and
\item a point $\gls*{z}$ in $\Sigma\setminus(\alphas\cup\betas)$,
\end{itemize}
such that the $\alpha$-circles and $\beta$-circles intersect
transversally and $\Sigma \setminus \alphas$ and $\Sigma \setminus
\betas$ are each connected.
\end{definition}
We will also sometimes refer to pointed Heegaard diagrams as in
Definition~\ref{sec:BorderedDiagrams} as \emph{closed} Heegaard
diagrams, to distinguish them from the bordered Heegaard diagrams
introduced soon.
\index{Heegaard diagram!closed}%

When discussing pointed Heegaard diagrams (and their bordered
variants), we slightly abuse notation, letting $\alphas$ denote both the set or
tuple of $\alpha$-curves and the union of the $\alpha$-curves, and
similarly for $\betas$. Note that the
$\alpha$-circles and $\beta$-circles are each maximal collections of
embedded circles with connected complement.  An equivalent requirement
to connectivity of the complement is that the $\alpha_i$ be
homologically linearly independent in $H_1(\Sigma)$, spanning a
Lagrangian subspace, and likewise for the~$\beta_i$.

Thus fortified, we turn to the bordered case.
\begin{definition}\label{def:bordered-3-mfld} A \emph{bordered
    $3$-manifold} is a triple $\gls*{borderedmfld}$ where $Y$ is a
  compact, oriented $3$-manifold with connected boundary $\bdy Y$\!,\,
  $\PMC$ is a pointed matched circle, and $\phi\co F(\PMC)\to \bdy Y$
  is an orientation-preserving homeomorphism.
  \index{bordered!$3$-manifold}%

  Two bordered $3$-manifolds $(Y,\PMC,\phi)$ and $(Y',\PMC',\phi')$ are
  called \emph{equivalent} if there is an orientation-preserving homeomorphism $\psi\co Y\to
  Y'$ so that $\phi'=\psi\circ\phi$.
  \index{bordered!$3$-manifold!equivalent}%
  \index{equivalent bordered $3$-manifolds}%
\end{definition}

Note that if $\phi\co F(\PMC)\to \partial Y$ and $\phi'\co F(\PMC)\to\partial Y$
are isotopic homeomorphisms, then, by viewing the isotopy as a homeomorphism
on a collar of $\partial Y$, we see that the bordered $3$-manifolds
$(Y,\PMC,\phi)$ and $(Y,\PMC,\phi')$ are equivalent.

We will often abuse notation and refer to a bordered $3$-manifold $Y$,
suppressing $\PMC$ and $\phi$. Also, by a bordered $3$-manifold we
will implicitly mean an equivalence class of bordered $3$-manifolds.

\begin{definition}
\label{def:BorderedDiagram}
\index{bordered!Heegaard diagram|see{Heegaard diagram, bordered}}%
\index{Heegaard diagram!bordered}%
A \emph{bordered Heegaard diagram}
is a quadruple
$\gls*{Heegaard}=(\gls*{bSigma},\widebar{\alphas},\gls*{betas},\gls*{z})$ consisting of
\begin{itemize}
\item a compact, oriented surface $\bSigma$ with one boundary component, of
  some genus~$g$;
\item a $g$-tuple of pairwise-disjoint circles
  $\gls*{betas}=\{\beta_1,\dots,\beta_g\}$ in the interior of~$\Sigma$;
  \glsadd{betai}
\item a $(g+k)$-tuple of pairwise-disjoint curves
  $\gls*{oalphas}$ in 
  $\bSigma$, split into $g-k$ circles
  $\gls*{alphasc}=(\alpha_1^c,\dots,\alpha_{g-k}^c)$ in the 
  interior of $\bSigma$ and $2k$ arcs
  $\gls*{alphasa}=(\balpha_1^a,\dots,\balpha_{2k}^a)$ in 
  $\bSigma$ with boundary on $\partial\bSigma$ (and transverse to
  $\partial{\bSigma}$); and
  $\glsadd{alphaia}\glsadd{alphaic}$
\item a point $\gls*{z}$ in $(\partial\bSigma)\setminus(\overline{\alphas}\cap\partial\bSigma)$,
\end{itemize}
such that the intersections are transverse and
$\bSigma\setminus \balphas$ and $\bSigma
\setminus \betas$ are connected.
\end{definition}

The conditions on the numbers of
$\beta$-circles and $\alpha$-circles and arcs again come naturally
from requiring that the tuples be maximal collections of curves of the
appropriate type with connected complement.  Again, an equivalent
requirement to connectivity of the complement is that the $\beta_i$ be
homologically linearly independent in $H_1(\bSigma)$ and that the
$\alpha_i$ be homologically linearly independent in
$H_1(\bSigma,\partial\bSigma)$.

The boundary of a bordered Heegaard diagram naturally has the
structure of a
pointed matched circle:
\begin{lemma}\label{lem:bdy-hd-is-pmc} Let
  $(\bSigma,\widebar{\alphas},\betas,z)$ be a bordered
  Heegaard diagram. Let $Z=\bdy\bSigma$ and
  $\CircPts=\widebar{\alphas}\cap\bdy\bSigma$. Define a
  matching $M$ on $\CircPts$ by
  $M(\alpha_i^a\cap\bdy\bSigma)=i$. Then $(Z,\CircPts,M,z)$ is
  a pointed matched circle, in the sense of
  Definition~\ref{def:PointedMatchedCircle}.
\end{lemma}
\begin{proof}
  We need to verify that the result $Z'$ of performing surgery on the
  points in~$\CircPts$ according to the matching $M$ gives a (connected)
  circle. 
  Let $\bSigma'$ be the result of performing surgery on $\bSigma$
  along each $\alpha$-circle, i.e., deleting a tubular neighborhood of
  each $\alpha$-circle and gluing in two disks to the result.
  The circle $Z'$ is homeomorphic to the boundary
  $B$ of a regular neighborhood of
  $\bdy\bSigma\cup\bigcup_{i=1}^{2k}\alpha_i^a$. But $B$, in
  turn, is homeomorphic to $\bdy
  (\bSigma'\setminus\nbd(\balphas))$, and
  $\bSigma'\setminus\nbd(\balphas)$ is a disk.
\end{proof}

\begin{definition}\label{def:bdy-of-bord-is-pmc}
  \index{pointed matched circle!boundary of bordered Heegaard diagram}%
  If $\HD$ is a bordered Heegaard diagram, we call
  the pointed matched circle appearing on its boundary (in
  the sense of Lemma~\ref{lem:bdy-hd-is-pmc}) the \emph{boundary of
    $\HD$} and denote it $\gls*{boundaryHD}$.
\end{definition}

\begin{figure}
  \centering
  \includegraphics[scale=.82]{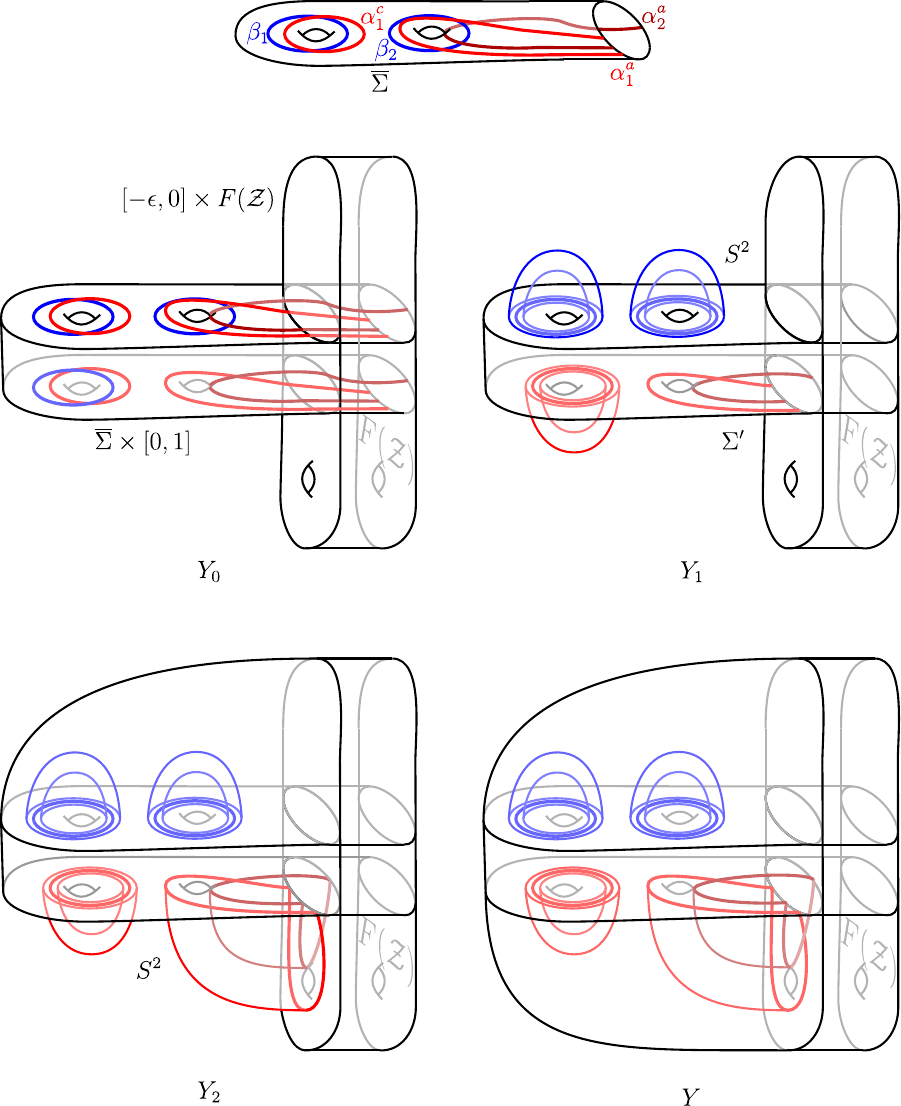}
  \caption[$3$-manifold represented by a bordered Heegaard
      diagram]{\textbf{The $3$-manifold represented by a bordered Heegaard
      diagram.} The bordered Heegaard diagram $\HD$ is shown at the top; the
  four steps in Construction~\ref{constr:bordered-hd-gives-mfld} are
  illustrated below.}
  \label{fig:bord-to-3mfld}
\end{figure}

\begin{construction}\label{constr:bordered-hd-gives-mfld}
  Let $\HD$ be a bordered
  Heegaard diagram. We will show how $\HD$ specifies a bordered
  $3$-manifold, with
  boundary parameterized by $F(\bdy\HD)$.

  The construction is illustrated in Figure~\ref{fig:bord-to-3mfld}.
  Write $\HD=(\bSigma,\widebar{\alphas},\betas,z)$, and write
  $\PMC=\bdy\HD$.
  Let $[-\epsilon,0]\times Z$ denote a closed collar neighborhood of
  $\bdy\bSigma$, so that $0 \times Z$ is identified with
  $\bdy\bSigma$. Choose also a closed tubular neighborhood
  $Z\times[0,1]$ of~$Z$ in $F(\PMC)$.

  Let $Y_0$ denote the three-manifold obtained by
  gluing $\bSigma\times [0,1]$ to $[-\epsilon,0]\times F(\PMC)$,
  by identifying
  $([-\epsilon,0]\times Z)\times [0,1] \subset \bSigma\times [0,1]$
  with
  $[-\epsilon,0]\times (Z\times [0,1])\subset [-\epsilon,0]\times F(\PMC)$.

  Next, attach a $3$-dimensional $2$-handle to each $\beta_i\times
  \{1\}\subset \bSigma\times[0,1]\subset Y_0$ and to each
  $\alpha_i^c\times\{0\}\subset \bSigma\times[0,1]\subset Y_0$. Call
  the result~$Y_1$. The manifold $Y_1$ has three boundary components:
  \begin{itemize}
  \item an $S^2$ which meets $\bSigma\times\{1\}$,
  \item a surface $\Sigma'$ of genus~$2k$ which meets
    $\bSigma\times\{0\}$, and
  \item a surface identified with $F(\PMC)$,
    given by $\{0\}\times F(\PMC)\subset [-\epsilon,0]\times F(\PMC)$.
  \end{itemize}
  Glue a $3$-ball to the $S^2$ boundary component of $Y_1$. 
  Join each
  $\alpha_i^a\times\{0\}\cap \Sigma'$ to the co-core of the corresponding handle
  in $\{-\epsilon\}\times F(\PMC)$ to form a closed curve, and attach
  a $3$-dimensional $2$-handle along each of these circles. Call the
  result of these attachings~$Y_2$. The manifold $Y_2$ has one $S^2$
  boundary component which meets $\bSigma\times\{0\}\subset Y_0$ and a
  boundary component identified with $F(\PMC)$. Glue a $3$-ball to the
  $S^2$ boundary component of $Y_2$ and call the result $Y$. Let
  $\phi$ be the obvious identification of $F(\PMC)$ with $\bdy
  Y$. Then $(Y,\PMC,\phi)$ is the desired bordered $3$-manifold.
\end{construction}

\begin{definition}\label{def:bord-mfld-repd-by-hd}
  Given a bordered Heegaard diagram we call the bordered $3$-manifold
  given by Construction~\ref{constr:bordered-hd-gives-mfld} the \emph{bordered
    $3$-manifold represented by $\HD$} and denote it $\gls*{YofHD}$.
  \index{bordered!$3$-manifold!represented by Heegaard diagram}%
  \index{Heegaard diagram!bordered!represents bordered $3$-manifold}%
\end{definition}

Construction~\ref{constr:bordered-hd-gives-mfld} introduces a few
extra pairs of canceling handles.  This is done to make it fit well
with the Morse theory below. A quicker (and equivalent) construction is the following:
$Y(\HD)$ is obtained by thickening the Heegaard surface $[0,1]\times
{\overline\Sigma}$, attaching three-dimensional two-handle to each
$\alpha_i^c\times \{0\}\times{\overline \Sigma}$, and a three-dimensional two-handle to each $\beta_i\times\{1\}\times {\overline \Sigma}$. The 
boundary of the result is naturally identified with $F(\PMC)$.

Conversely, any $3$-manifold with connected boundary is
specified by some bordered Heegaard diagram. One way to see this is in
terms of Morse theory, following~\cite[Section
2.2]{Lipshitz06:BorderedHF}:

\begin{definition}
\label{def:CompatibleMorseNew}
Fix a bordered $3$-manifold $(Y,\PMC,\phi)$.
We say that a pair consisting of a Riemannian metric $\gls*{Metric}$ on
$Y$ and a self-indexing Morse function
$\gls*{MorseFn}$ on $Y$ are \emph{compatible with $(Y,\PMC,\phi)$} if
\index{Morse function!compatible with bordered $3$-manifold}%
\begin{enumerate}
\item\label{item:compat-morse-1} the boundary of $Y$ is geodesic,
\item\label{item:compat-morse-2} the gradient vector field $\nabla f|_{\bdy Y}$ is tangent to $\bdy Y$,
\item\label{item:compat-morse-3} $f$ has a unique index $0$ and a unique index $3$ critical point, both of which
  lie on $\bdy Y$, and are the unique index $0$ and $2$ critical points
  of $f|_{\bdy Y}$, respectively,
\item\label{item:compat-morse-4} the index $1$ critical points of $f|_{\bdy Y}$ are also index
  $1$ critical points of $f$, and
\item\label{item:compat-morse-5} the pair $(f\circ\phi,\phi^*g)$ on $F(\PMC)$ are compatible with the
  pointed matched circle $\PMC$ in the sense of
  Definition~\ref{def:PMCCompatibleMorse}.
\end{enumerate}
\end{definition}

\begin{lemma}\label{lem:3mfld-heegaard}
  Any bordered $3$-manifold $(Y,\PMC,\phi)$ is represented by some
  bordered Heegaard diagram $\HD$.
\end{lemma}
\begin{proof}
  The first step is to find a Morse function $f$ and metric $g$ compatible
  with $(Y,\PMC,\phi)$ as in
  Definition~\ref{def:CompatibleMorseNew}. Choose a Morse
  function $f_F$ and metric $g_F$ on $F(\PMC)$ compatible with~$\PMC$ in the sense of
  Definition~\ref{def:PMCCompatibleMorse} which is self-indexing
  except that $f_F$ takes is $3$ on the index~$2$ critical
  point. Extend $f_F\circ \phi^{-1}$ and
  $\phi_*g_F$ to $f$ and~$g$ on a collar neighborhood of $\bdy Y$ satisfying
  Conditions (\ref{item:compat-morse-1}), (\ref{item:compat-morse-2})
  and~(\ref{item:compat-morse-4}), and so that the index $0$ and $2$
  critical points of $f_F\circ \phi^{-1}$ are index $0$ and $3$ critical
  points of~$f$. (Note that
  Condition~(\ref{item:compat-morse-5}) is automatically satisfied.)
  Extend $f$ and $g$ arbitrarily to
  the rest of~$Y$. Now consider the graph formed by flows between the
  index~$0$ and~$1$ critical points, which is connected.  Since the
  descending flows from boundary
  index~$1$ critical points remain on the boundary, we can and do
  cancel every interior index~$0$
  critical point with an interior index~$1$ critical
  point.  Similarly cancel interior index~$3$ critical points with
  index~$2$ critical points, so that Condition~(\ref{item:compat-morse-3})
  is satisfied.
  Modify the function in the interior so the result is
  self-indexing.

  Now, given the boundary compatible pair $(f,g)$, construct a
  bordered Heegaard diagram $(\bSigma,\balphas,\betas,z)$ as follows.
  Let the Heegaard surface $\overline{\Sigma}$ be $f^{-1}(3/2)$, the
  curves $\widebar{\alphas}$ be the intersection of the ascending
  disks of the index $1$ critical points of $f$ (including those on
  $\bdy Y$) with~$\overline{\Sigma}$, and the curves $\betas$ be the intersection of
  the descending disks of the index $2$ critical points of $f$ with~$\overline{\Sigma}$.
\end{proof}

The Morse theory description also implies that any two bordered
Heegaard diagrams for the same bordered $3$-manifold can be connected
by (bordered) Heegaard moves, as specified in the following:
\begin{proposition}\label{prop:heegaard-moves}
Any pair of bordered Heegaard diagrams for equivalent bordered
$3$\hyp manifolds can be made diffeomorphic by a sequence of
\index{Heegaard moves}%
\begin{itemize}
\item isotopies of the $\alpha$-curves and $\beta$-circles, not crossing
  $\partial\overline{\Sigma}$,
  \index{isotopy|see{Heegaard moves, isotopy}}%
  \index{Heegaard moves!isotopy}%
\item handleslides of $\alpha$-curves over $\alpha$-\emph{circles} and $\beta$-circles over
  $\beta$-circles, and
  \index{handleslide}%
  \index{Heegaard moves!handleslide}%
\item stabilizations and destabilizations in the interior of
  $\bSigma$.
  \index{stabilization}%
  \index{Heegaard moves!stabilization}%
  \index{destabilization}%
  \index{Heegaard moves!destabilization}%
\end{itemize}
\end{proposition}

Before giving the proof of
Proposition~\ref{prop:heegaard-moves}, for which we follow closely the
proof of \cite[Proposition~2.2]{OS04:HolomorphicDisks}, we include two
lemmas which will be useful for excluding stabilizing by index zero
(or three) critical points
(compare~\cite[Lemmas~2.3~and~2.4]{OS04:HolomorphicDisks}).

\begin{lemma}
	\label{lem:AddOneMoreCurve}
	Let $F$ be a compact surface of genus~$h$ with a single
        boundary component.  Let $\gammas=\{\gamma_1,\dots,\gamma_h\}$ be an
	$h$-tuple of embedded, pairwise disjoint, simple closed curves so that
	$F\setminus\gammas$ is a punctured disk.  Suppose moreover
        that $\delta$
	is a simple closed curve in $F$ which is disjoint from the
	$\gammas$. Then either $\delta$ is null-homologous or there is
	some $\gamma_i$ with the property that $\delta$ is isotopic
	to a curve obtained by handlesliding $\gamma_i$ across some
	collection of the $\gamma_j$ with $j\neq i$.
\end{lemma}

\begin{proof}
	If we surger out $\gamma_1,\dots,\gamma_{h}$, we replace $F$
	by the disk $D$ with $2h$ marked points
	$\{p_1,q_1,\dots,p_{h},q_{h}\}$ (where the pair of points
	$\{p_i,q_i\}$ corresponds to the zero-sphere which replaced
	the circle $\gamma_i$).  Now, $\delta$ can be viewed as a
	Jordan curve in the disk $D$. If $\delta$ separates some $p_i$ from
	its corresponding $q_i$, then it is easy to see that $\gamma$
	is isotopic to the curve gotten by handlesliding $\gamma_i$
	over some collection of the $\gamma_j$ with $j\neq
	i$. Otherwise, $\delta$ was null-homologous.
\end{proof}

\begin{lemma}
	\label{lem:HandleslideChoices}
	Let $F$ be a compact surface of genus $h$ with a single
        boundary component and $\{\gamma_1,\dots,\gamma_d\}$ a
        collection of embedded, pairwise disjoint, simple closed
        curves.  Then any two $h$-tuples of linearly-independent
        $\gamma_i$'s are related by a sequence of isotopies and
        handleslides.
\end{lemma}

\begin{proof}
	This is proved by induction on $h$. The case where $h=0$ is
	vacuously true. Next, suppose we have two subsets with an
	element in common, which we label $\gamma_1$. Surger
	out $\gamma_1$ to obtain a surface $F'$ with one lower
	genus. Isotopies in $F'$ which cross the two surgery points
	correspond to handleslides in $F$ across $\gamma_1$. Thus, in
	this case, the proof follows from the induction hypothesis.
	Finally, consider the case where we have two disjoint subsets
	$\{\gamma_1,\dots,\gamma_h\}$ and
	$\{\gamma_{h+1},\dots,\gamma_{2h}\}$. The curve
	$\gamma_{h+1}$ cannot be null-homologous, so, after renumbering,
	we can obtain $\gamma_{h+1}$ by handlesliding $\gamma_1$ over
	some of the $\{\gamma_2,\dots,\gamma_h\}$, according to Lemma~\ref{lem:AddOneMoreCurve}. Thus, we have reduced
	to the case where the two subsets are not disjoint.
\end{proof}

\begin{proof}[Proof of Proposition~\ref{prop:heegaard-moves}]
  The proof involves adapting standard handle calculus as in~\cite[Section 2.1]{OS04:HolomorphicDisks}.

  Any two bordered Heegaard diagrams $\HD_0$ and
  $\HD_1$ for equivalent bordered manifolds $(Y_i,\PMC_i,\phi_i)$ come from Morse functions and metrics $(f_0,g_0)$ and
  $(f_1,g_1)$ on $Y_i$ compatible with the bordering $(\PMC,\phi)$ of
  $Y_i$, in the sense of
  Definition~\ref{def:CompatibleMorseNew}. Using the equivalence
  between the two bordered manifolds, we can view $f_1$ and $g_1$ as
  being defined on $Y=Y_1$. Moreover, we may
  choose $f_0$, $g_0$, $f_1$ and $g_1$ so that for some neighborhood
  $N(\partial Y)$ of $\partial Y$, $f_1|_{N(\bdy Y)}=f_2|_{N(\bdy Y)}$
  and $g_1|_{N(\bdy Y)}=g_2|_{N(\bdy Y)}$.
  We may connect $f_0$ and $f_1$
  by a generic path $f_t$ of smooth functions which are Morse
  functions for generic~$t$ and so that the family is constant on $N(\bdy Y)$.

  Although the $f_t$ will not in general be self-indexing, there are
  still no flow lines from higher-index critical points to
  lower-index.  To $f_t$ at generic~$t$ we can still associate a
  generalized Heegaard diagram~$\HD_t$, possibly with extra $\alpha$-
  and $\beta$-circles.

  While $f_t$ remains a Morse function, the Heegaard diagram changes
  by isotopy of the $\alpha$- and $\beta$-curves.  When $f_t$ fails to
  be a Morse function, there is either a birth-death singularity
  involving a pair of critical points of adjacent index, or there is a
  flow line between two critical points of the same index.  The
  possibilities, and the way the generalized Heegaard diagram changes, are:
  \begin{description}
  \item[Index 0 and 1 birth--death] A disjoint
	$\alpha$-circle is added or deleted, without changing the
	span of the $\alpha$-circles.
  \item[Index 1 and 2 birth--death] The Heegaard diagram is stabilized or destabilized.
  \item[Index 2 and 3 birth--death] A disjoint
	 $\beta$-circle is added or deleted, without changing the span
	of the $\beta$-circles.
  \item[Flow line between index 1 critical points] A handleslide of an
    $\alpha$-circle or $\alpha$-arc over an $\alpha$-circle.  Since
    our conditions imply that the descending disk of the index~1
    critical points corresponding to $\alpha$-arcs lies entirely on the
    boundary, we cannot have a handleslide over an $\alpha$-arc.
  \item[Flow line between index 2 critical points] A handleslide of a
    $\beta$-circle over a $\beta$-circle.
  \end{description}
  Note also that the critical points on the boundary cannot be
  involved in any cancellations, since $f$ was fixed on all of
  $N(\partial Y)$.

  Now, given any two bordered
  Heegaard diagrams $\HD_0$ and $\HD_1$ for $Y$, in view
  of the above remarks, we can stabilize $\HD_0$ and $\HD_1$
  to obtain a pair of 
  bordered Heegaard diagrams $\HD_0'$ and $\HD_1'$ with of the same genus
  and with the following property. Each of $\HD_0'$ and $\HD_1'$
  can be extended by adding new $\alpha$- and $\beta$-circles if
  necessary to form a pair of generalized bordered Heegaard diagrams
  $\HD_0''$ and $\HD_1''$ respectively, which in turn can
  be connected by isotopies and handleslides (subject to the constraints from
  Proposition~\ref{prop:heegaard-moves},
  i.e., $\alpha$-arcs are allowed to slide over
  $\alpha$-circles, but not vice versa). The proposition is
  proved then once we show that for
  each generalized Heegaard diagram $\HD''$ of genus~$g$, 
  if $\HD_1$ and $\HD_2$ are any two
  ordinary bordered Heegaard diagrams obtained from $\HD''$
  by picking  out $g-k$
  of the $\alpha$-circles and $g$ of the $\beta$-circles, each
  homologically linearly independent,
  then $\HD_1$ and $\HD_2$ can be connected by handleslides and isotopies.
  
  To this end, consider any two $(g-k)$-tuples of the $\alpha$-circles
  of $\HD''$ chosen so that they, along with the $2k$ $\alpha$-arcs,
  span a $g+k$-dimensional subspace of $H_1(\bSigma,\partial\bSigma)$.
  The fact that
  these two $(g-k)$-tuples can be connected by isotopies and
  handleslides follows directly from
  Lemma~\ref{lem:HandleslideChoices}, applied to the surface
  $F={\widebar\Sigma}\setminus\{\alpha_1^a,\dots,\alpha_{2k}^a\}$.  It
  follows similarly that any two $g$-tuples of the $\beta$ circles
  which are linearly independent in $H_1(\overline\Sigma)$ can be connected by handleslides
  and isotopies.  This completes the proof.
\end{proof}

Let $(\gls*{Sigma},\gls*{noalphas},\betas,z)$ be the result of attaching
a \emph{cylindrical end}
\index{cylindrical end}%
\index{end!cylindrical|see{cylindrical end}}%
to $(\bSigma,\widebar{\alphas},\betas,z)$, in the sense of
symplectic field theory \cite{EGH00:IntroductionSFT}. Topologically
$\Sigma=\bSigma\setminus\partial\bSigma$, and
$\alpha_i^a=\balpha_i^a\setminus(\partial\balpha_i^a)$.
In due course, we will choose a conformal structure on $\Sigma$ making it a punctured
Riemann surface and so that each $\alpha_i^a$ is radial near the puncture. Abusing
terminology, we will often also refer to $(\Sigma,\alphas,\betas,z)$
as a bordered Heegaard diagram.
$\glsadd{noalphasa}$

\section{Examples of bordered Heegaard diagrams}
\label{sec:Bestiary}

Here are a few
examples of bordered Heegaard diagrams.

\begin{figure}
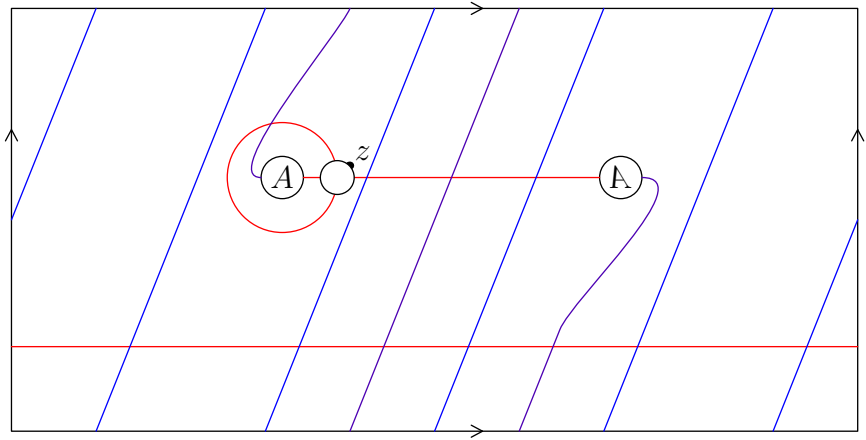

  \begin{center}
    $\mfig{trefoil-3}$
  \end{center}
  \caption[Bordered Heegaard diagram for the trefoil complement]{\textbf{Bordered Heegaard diagram for the trefoil complement.}  This
  is based on the observation that $+5$ surgery on the trefoil is the
  lens space $L(5,-1)$, with the knot forming the core of the solid
  torus a Berge knot.  The diagram is on a genus 2 surface,
  represented as a torus with a handle attached (at
  $A$).
  \index{trefoil}%
  \index{Heegaard diagram!bordered!for trefoil complement}%
}
  \label{fig:TrefoilComplement}%
\end{figure}

Fix an oriented surface $\Sigma$, equipped with a $g$-tuple of
pairwise disjoint, homologically independent curves, $\betas$, and a
$(g-1)$-tuple of pairwise disjoint, homologically independent curves
$\alphas^c=\{\alpha_1^c,\dots,\alpha_1^{g-1}\}$. Then
$(\Sigma,\alphas^c,\betas)$ is a Heegaard diagram for a
three-manifold with torus boundary, and indeed any such
three-manifold~$Y$ can be described by a Heegaard diagram of this
type. To turn such
a diagram into a bordered Heegaard diagram, we proceed as follows.
Fix an additional pair of circles $\gamma_1$ and $\gamma_2$ in
$\Sigma$ so that:
\begin{itemize}	
\item $\gamma_1$ and $\gamma_2$ intersect, transversally, in a
  single point~$p$ and
\item both of the homology classes
  $[\gamma_1]$ and $[\gamma_2]$ are homologically independent
  from $[\alpha_1^c],\dots,[\alpha_{g-1}^c]$.
\end{itemize}
Let $D$ be a disk around~$p$ which is disjoint from all the above
curves, except for $\gamma_1$ and~$\gamma_2$, each of which it meets
in a single arc. Then, the complement of $D$ specifies a bordered
Heegaard diagram for $Y$, for some parametrization of $\partial Y$.
A bordered Heegaard diagram for the trefoil complement is illustrated
in Figure~\ref{fig:TrefoilComplement}.
(See also Chapter~\ref{chap:TorusBoundary} for a further discussion of
bordered three-manifolds with torus boundary.)

\begin{figure}
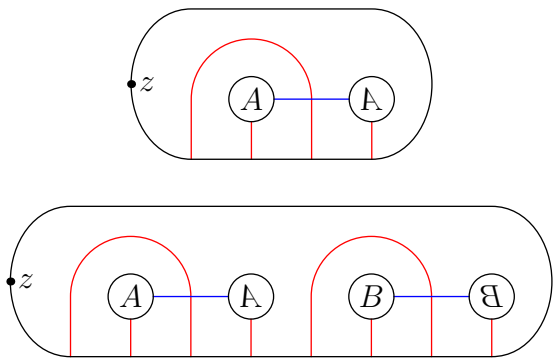

  \begin{gather*}
    \mfigb{torus-5}\\[10pt]\mfigb{torus-6}
  \end{gather*}
  \caption[Heegaard diagrams for handlebodies]{\textbf{Heegaard diagrams for handlebodies.} Top: a bordered genus one handlebody.  Bottom: the boundary
    connect sum of two copies of the handlebody on top.}
  \label{fig:BoundaryConnSum}
\end{figure}

As another example, consider boundary connect sums.
\index{boundary connect sum!of bordered Heegaard diagrams}%
\index{boundary connect sum!of bordered $3$-manifolds}%
Fix bordered
Heegaard diagrams $\HD_i=(\bSigma_i,\widebar{\alphas}_i,
\betas_i,z_i)$ for $Y_i$ with $i=1,2$. Take the boundary connect sum of
$\bSigma_1$ and $\bSigma_2$ along $z_1$ and $z_2$,
i.e., attach a rectangle to $\bSigma_1\cup \bSigma_2$
by gluing two opposite sides to $\partial{\Sigma}_1$ and
$\partial{\Sigma}_2$ along intervals containing $z_1$ and
$z_2$. Introduce a new basepoint $z$ on one of the two remaining boundary
components of the rectangle. (This involves a choice.)
This
procedure gives a bordered Heegaard diagram for the boundary connect
sum of $Y_1$ with $Y_2$. See
Figure~\ref{fig:BoundaryConnSum}, where we form the boundary connect
sum of two bordered Heegaard diagrams for genus-one handlebodies to
obtain a bordered Heegaard diagram for a genus-two handlebody.

\begin{figure}
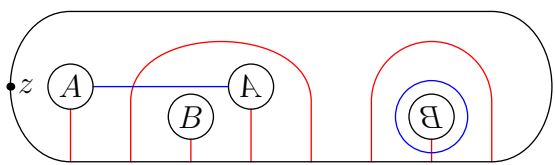

  \[
  \mfig{torus-7}
  \]
  \caption[Another Heegaard diagram for a genus $2$
      handlebody]{\textbf{Another bordered Heegaard diagram for a genus two
      handlebody,} using a different parameterization of the
    boundary.}
  \label{fig:GenusTwoInteresting}
\end{figure}

Any bordered Heegaard diagram in which
there are no $\alpha$ circles necessarily represents a genus $g$
handlebody. In Figure~\ref{fig:GenusTwoInteresting}, we have
illustrated a genus two handlebody with respect to a parameterization of the
genus two surface which is different from the one used in Figure~\ref{fig:BoundaryConnSum}.
\index{handlebody, bordered Heegaard diagrams for}%
\index{Heegaard diagram!bordered!for handlebody}%

\section{Generators, homology classes and
  \textalt{$\spin^c$}{spin-c} structures}
\label{sec:homology-classes-generators}
Fix a bordered Heegaard diagram~$\HD$, specifying some
three-manifold~$Y$.
\begin{definition}\label{def:generator}
  By a \emph{generator} of $\HD$
\index{generator}%
  we mean a $g$-element subset
  $\gls*{genx}=\{x_1,\dots,x_g\}$ so that 
  \begin{itemize}
    \item exactly one $\gls*{genxi}$ lies on each $\beta$-circle,
    \item exactly one $x_i$ lies on each $\alpha$-circle and
    \item at most one $x_i$ lies on each $\alpha$-arc.
  \end{itemize}
  Let $\gls*{generators}$ 
or $\gls*{generatorsTwo}$ denote the set of generators. 
Given a generator $\x$, let $\gls*{ox}$ 
denote the set of $\alpha$-arcs
which are occupied by $\x$, i.e., $o(\x) \coloneqq \{\,i \mid
\x\cap\alpha_i^a\neq\emptyset\,\}\subset[2k]$.
\index{occupied $\alpha$-arcs}%
\end{definition}

Soon (Chapter~\ref{chap:structure-moduli} and later) we will be
interested in holomorphic curves in $\Sigma\times I_s\times\RR_t$,
where $I_s = [0,1]$ is the unit interval with parameter~$s$ and
$\RR_t$ is $\RR$ with parameter~$t$.  These curves have boundary
on $\alphas\times\{1\}\times\RR_t$ and
$\betas\times\{0\}\times\RR_t$. They will be asymptotic to $g$-tuples
of chords $\x\times I_s$ and $\y\times I_s$ at $\pm\infty$, where
$\x$ and $\y$ are generators. Each such curve carries a relative homology
class, and we let $\pi_2(\x,\y)$ denote the set of these relative
homology classes. More precisely, we make the following definition.
\begin{definition}\label{def:pi2}
  Fix generators $\x$ and $\y$. Let $I_s = [0,1]$
  and~$I_t = [-\infty,\infty]$ be intervals. We work in the relative homology
  group
  \[
    H_2\bigl(\bSigma \times I_s \times I_t,
    \bigl((S_\alphas \cup S_\betas \cup S_\bdy) \times I_t\bigr) \cup
    (G_\x \times \{-\infty\}) \cup (G_\y \times \{\infty\})\bigr),
  \]
  where
  \begin{align*}
    S_\alphas &= \alphas \times \{1\} &
    S_\betas &= \betas \times \{0\} &
    S_\bdy &= (\bdy\bSigma\setminus z) \times I_s\\
    G_\x &= \x \times I_s &
    G_\y &= \y \times I_s.
  \end{align*}
  Let $\gls*{pitwo}$, the \emph{homology classes connecting $\x$ to
    $\y$},
\index{homology class!connecting generators}%
  denote those
  elements of this group which map to the relative fundamental class
  of $\x \times I_s \cup \y \times I_s$ under the composition of the
  boundary homomorphism and collapsing the remainder of the boundary.
\end{definition}

\begin{definition}
\label{def:Domain}
Given a homology class $\gls*{homolclass}\in\pi_2(\x,\y)$, projecting
to $\bSigma$ gives a well-defined element of
$H_2(\bSigma,\alphas\cup\betas\cup \partial{\overline{\Sigma}})$,
i.e., a linear combination of the components of $\bSigma \setminus
(\alphas \cup \betas)$, which we sometimes call \emph{regions}. We
call the linear combination of regions corresponding to $B$ the
\emph{domain} of~$B$.  \index{domain}\index{region}%
\end{definition}

It is not hard to see that a homology class $B$ is uniquely determined
by its domain.  

\begin{definition}
\label{def:LocMult}
Given any point $p\in
\Sigma\setminus(\alphas\cup\betas)$, the \emph{local multiplicity}
of~$B$ at~$p$ is the local degree of the two-chain at the point~$p$.
\index{local multiplicity}%
\end{definition}

By hypothesis, the coefficient in the domain of $B$ of the region
$\gls*{Dz}$ adjacent to $z\in\bdy\bSigma$ is~$0$.

\begin{remark}
  We use the notation $\pi_2$ for this object to agree with the original
  conventions~\cite{OS04:HolomorphicDisks}.
  In the cylindrical setting of Heegaard Floer homology this is a
  homology class rather than a homotopy class. Also, note that, unlike
  the closed case, we require the multiplicity of a domain to be $0$
  at $z$; this is because we are only concerned with $\HFa$.
\end{remark}

Notice that concatenation (in the $t$ factor) gives a product
$\gls*{concat}$ from $\pi_2(\x,\y)\times\pi_2(\y,\w)$ to
$\pi_2(\x,\w)$; this operation
corresponds to addition of domains.
\index{concatenation}%
If $\pi_2(\x,\y)$ is non-empty then
concatenation makes $\pi_2(\x,\y)$ into a free, transitive
$\pi_2(\x,\x)$-set. Elements of $\pi_2(\x,\x)$, or their domains, are called
\emph{periodic domains}.
\index{periodic domain|see{domain, periodic}}%
\index{domain!periodic}%

\begin{lemma}\label{lem:pi2-h2}
  There is a natural group isomorphism $\pi_2(\x,\x)\cong
  H_2(Y,\bdy Y)$.
\end{lemma}

\begin{proof}
  The proof is a simple adaptation from the closed case (\cite[Lem\-ma
  2.1]{Lipshitz06:CylindricalHF} in this cylindrical setting); the
  present case is also proved in~\cite[Lemma
  2.6.1]{Lipshitz06:BorderedHF}, which we repeat here.  Let $z'$ be a
  point in the interior
  of $\bSigma$ near $z$. Let
  $\Sigma'=(\bSigma/\bdy\bSigma)\setminus\{z'\}$.  There is an
  isomorphism
  \[
  \pi_2(\x,\x)\cong
  H_2(\Sigma'\times[0,1],(\balphas\times\{1\})\cup
  (\betas\times\{0\}))
  \]
  gotten by adding (or subtracting) a copy of
  $\x\times[0,1]\times\RR$. From the long exact sequence for the
  pair $(\Sigma'\times[0,1],(\balphas\times\{1\})\cup
  (\betas\times\{0\}))$, we have
  \[
    0 \to
    H_2(\Sigma'\times[0,1],(\balphas\times\{1\})\cup
      (\betas\times\{0\})) \to
    H_1((\balphas\times\{1\})\cup
      (\betas\times\{0\})) \to
    H_1(\Sigma').
  \]
  (Note that $H_2(\Sigma'\times[0,1])=0$.) Consequently, 
  \[
  \pi_2(\x,\x)\cong \ker\left(H_1(\balphas/\bdy\balphas)\oplus
  H_1(\betas)\to H_1(\bSigma/\bdy\bSigma)\right).
  \]
  But it follows by analyzing
  Construction~\ref{constr:bordered-hd-gives-mfld} that this kernel is exactly
  $H_2(Y,\bdy Y)$.
\end{proof}

We sometimes denote the image of
$B\in\pi_2(\x,\x)$ in $H_2(Y,\bdy Y)$ by~$\gls*{homolclasshomolclass}$.

We split up the boundary of the domain of $B$ into three pieces,
$\gls*{bdyalpha}$, contained in~$\alphas$, $\gls*{bdybeta}$, contained
in~$\betas$, and $\gls*{bdybdy}$, 
 contained in $\bdy \bSigma$.
We orient them so that the oriented boundary of the domain of $B$ is
$\bdy^\alpha(B) + \bdy^\beta(B) + \bdy^\bdy(B)$; in particular, this
implies that
$-\bdy(\bdy^\beta B)=\bdy(\bdy^\alpha B)+\bdy(\bdy^\bdy B)=
\y - \x$.  We will think of
$\bdy^\bdy B$ as an element of $H_1(\bdy \bSigma, \CircPts)$.

\begin{definition} A homology class $B$ connecting $\x$ to $\y$ is called
  \emph{provincial} if $\bdy^\bdy B=0$.  Let $\gls*{pitwobdy}$ denote
  the set of provincial homology classes connecting $\x$ to $\y$.
  \index{provincial|see{homology class, provincial}}%
  \index{homology class!provincial}%
  \index{domain!provincial}%
\end{definition}
In other words, a homology class $B\in\pi_2(\x,\y)$ belongs to $\pi_2^\bdy(\x,\y)$ if
and only if the domain of $B$ does not include any of the regions
adjacent to $\bdy\bSigma$.

The concatenation maps obviously restrict to maps
\index{concatenation}%
$\gls*{concat}\co \pi_2^\bdy(\x,\y)\times\pi_2^\bdy(\y,\w)\to\pi_2^\bdy(\x,\w)$, again
making $\pi_2^\bdy(\x,\y)$ into a free, transitive
$\pi_2^\bdy(\x,\x)$-set (when non-empty). Elements of
$\pi_2^\bdy(\x,\x)$, or their domains,
are called \emph{provincial periodic domains}.

\begin{lemma}\label{lem:pi2-h2-bdy}
  There is a natural group isomorphism $\pi_2^\bdy(\x,\x)\cong H_2(Y)$.
\end{lemma}

\begin{proof}
  Again, the proof is a simple adaptation of the closed case, and is
  given in \cite[Lemma 2.6.4]{Lipshitz06:BorderedHF}. There is an
  isomorphism
  \[
  \pi_2^\bdy(\x,\y)\cong H_2(\bSigma\times[0,1],(\alphas\times\{1\})\cup(\betas\times\{0\})).
  \]
  Using the long exact sequence for the pair
  $(\bSigma\times[0,1],(\alphas\times\{1\})\cup(\betas\times\{0\}))$, it
  follows that 
  \[
  H_2(\bSigma\times[0,1],(\alphas\times\{1\})\cup(\betas\times\{0\}))\cong
  \ker(H_1(\alphas)\oplus H_1(\betas)\to H_1(\bSigma)).
  \]
  The Mayer-Vietoris sequence for the decomposition of $Y$ along
  $\bSigma$ identifies this kernel with $H_2(Y)$.
\end{proof}

Abusing notation, we may also sometimes
denote the image of $B\in\pi_2^\bdy(\x,\x)$ in $H_2(Y)$ by~$[B]$.

Next we ask when $\pi_2(\x,\y)$ and $\pi_2^\bdy(\x,\y)$ are
non-empty. As in the closed case \cite[Section 2.6]{OS04:HolomorphicDisks}, the answer relates to
$\spin^c$ structures on $Y$. Also as in the closed case, our treatment
is close to Turaev's~\cite{Turaev97:spinc}.

Choose a self-indexing Morse function $f$ and metric $g$ on $Y$ inducing
the bordered Heegaard diagram $(\bSigma,\balphas,\betas)$. Let $\x$ be a
generator. Using $\nabla f$, $z$ and $\x$ we will construct a
non-vanishing vector field $\vec{v}_z(\x)$; this reduces the structure
group of $TY$ from $SO(3)$ to $SO(2)=U(1)$, and composing with the
inclusion $U(1)\hookrightarrow U(2)=\spin^c(3)$ given by
$(e^{i\theta})\mapsto
\bigl(\begin{smallmatrix}
  e^{i\theta} & 0\\
  0 & 1
\end{smallmatrix}\bigr)$
defines a $\spin^c$ structure on $Y$.
\index{spinc-structure@$\SpinC$-structure}%
\index{spinc-structure@$\SpinC$-structure!associated to vector field}%

Each $x_i\in\x$ lies on a flow line $\gamma_i$ of $\nabla f$ from an
index $1$ critical point of $f$ to an index $2$ critical point of
$f$. The point $z$ lies on a flow line $\gamma$ from the index $0$
critical point of $f$ to the index $3$ critical point of $f$. Let $B$
denote a regular neighborhood of
$\gamma_1\cup\dots\cup\gamma_g\cup\gamma$. Since each flow line
$\gamma_i$ or $\gamma$ runs between critical points of opposite
parity, it is possible to extend $\nabla f|_{Y\setminus B}$ to a
vector field $\vec{v}^{\,0}_z(\x)$ on $Y$ so that $\vec{v}^{\,0}_z(\x)$ is non-vanishing on
$B$. The vector field $\vec{v}^{\,0}_z(\x)$ still has exactly $k$ zeroes, one
for each $\balpha_i^a$ not containing an $x_i$. All of these zeroes lie on
$\bdy Y$, however, so one can modify $\vec{v}_z^{\,0}(\x)$ in a neighborhood of
these zeroes to produce a non-vanishing vector field
$\vec{v}_z(\x)$. This modification can be done in a standard way,
depending only on $o(\x)$, and not on the generator $\x$ itself.
The vector field $\vec{v}_z(\x)$ induces the $\spin^c$ structure $\gls*{szx}$. It
is routine to verify that $\s_z(\x)$ does not depend on the choices
made in its construction. 
\index{spinc-structure@$\SpinC$-structure!associated to generator}%

Notice that the restriction of $\s_z(\x)$ to a collar neighborhood of
$\bdy Y$ depends only on the set of occupied $\alpha$-arcs $o(\x)$, and not on $\x$
itself. That is, for $o$ a $k$-element subset of
$\{\balpha_1^a,\dots,\balpha_{2k}^a\}$, there is an induced
$\spin^c$ structure $\s_z(o)$ on a collar neighborhood of $\bdy
Y$.
\index{spinc-structure@$\SpinC$-structure!relative!associated to generator}%
Given $\x$, let $\gls*{szxrel}$ denote the relative $\spin^c$
structure on $(Y,\nbd(\bdy Y))$ induced by $(z,\x)$, relative to
$\s_z(o)$.\footnote{A $\SpinC$-structure on $Y$ is a lift of the
  classifying map $[TY]\co Y\to \BSO(3)$ to a map $\spinc\co Y\to
  \BSpinC(3)$. Given a $\SpinC$-structure $\s_V$ on a subset
  $V\subset Y$, a relative $\SpinC$-structure on $Y$ relative to
  $\s_V$ is a lift $\spinc\co Y\to\BSpinC(3)$ with
  $\spinc|_{V}=\s_V$. The relative $\SpinC$-structures on $(Y,V)$ form
  an affine copy of $H^2(Y,V)$.

  While we normally consider $\SpinC$-structures up to homotopy (and
  relative $\SpinC$-structures up to homotopy relative to $V$) to
  make sense of relative $\SpinC$-structures on $(Y,V)$ relative to
  $\spinc_V$, it is important that the
  lift $\spinc_V$ over $V$ be given on the nose, not just up to homotopy.
  \index{spinc-structure@$\SpinC$-structure!relative}%
  \index{nose}%
}

\begin{lemma}
  \label{lem:SpinCStructures}
  Given generators $\x$ and $\y$,
  $\pi_2(\x,\y)\neq\emptyset$ if and only if $\s_z(\x)=\s_z(\y)$. Further,
  $\pi_2^\bdy(\x,\y)\neq\emptyset$ if and only if $o(\x)=o(\y)$ and
  $\s_z^\rel(\x)=\s_z^\rel(\y)$.
\end{lemma}
\begin{proof}
  Again, this is a straightforward adaptation of the closed case, and
  is also proved in~\cite[Lemmas 2.6.2 and
  2.6.5]{Lipshitz06:BorderedHF}.

  On the one hand, the obstruction $\gls*{epsilonxy}$ to finding a
  homology class connecting $\x$ to $\y$ is gotten as follows: join
  $\x$ to $\y$ by a union of paths
  $\gamma_\alpha\subset\alphas\cup(\bdy\bSigma\setminus z)$ and by a
  union of paths $\gamma_\beta\subset\betas$. Then $\x$ and $\y$ can
  be connected by a homology class if and only
  $\gamma_\alpha-\gamma_\beta$ can be made null-homologous in
  $\bSigma$ by adding or subtracting some $\alpha$-curves and
  $\beta$-circles, i.e., if and only if the image $\epsilon(\x,\y)$ of
  $\gamma_\alpha-\gamma_\beta$ in
  $H_1(\bSigma\times[0,1],\balphas\times\{0\}\cup\betas\times\{1\}\cup(\bdy\bSigma\setminus z)\times[0,1])\cong
  H_1(Y,\bdy Y)$ vanishes.

  On the other hand, the difference class $\s_z(\y)-\s_z(\x)$ is
  gotten as follows. As in the definition of $\s_z(\x)$, after fixing
  a Morse function and metric inducing the Heegaard diagram, the
  generator $\x$ corresponds to $g$ flows
  $\gamma_{x_1},\dots,\gamma_{x_g}$ in $Y$, and the generator $\y$
  corresponds to $g$ flows $\gamma_{y_1},\dots,\gamma_{y_g}$ in
  $Y$. The difference class is given by
  \[
  \s_z(\y)-\s_z(\x)=[\gamma_{y_1}+\dots+\gamma_{y_g}-\gamma_{x_1}-\dots-\gamma_{x_g}]\in
  H_1(Y,\bdy Y)\cong H^2(Y).
  \]
  The chain
  $\gamma_{y_1}+\dots+\gamma_{y_g}-\gamma_{x_1}-\dots-\gamma_{x_g}$ is
  homologous in $(Y,\bdy Y)$ to the chain $\gamma_\alpha-\gamma_\beta$
  defining $\epsilon(\x,\y)$.
  
  The case for $\s_z^\rel$ is exactly analogous, except that now
  $\gamma_\alpha\subset \alphas$, the group
  $H_1(\bSigma\times[0,1],\balphas\times\{0\}\cup\betas\times\{1\}\cup(\bdy\bSigma\setminus z)\times[0,1])\cong
  H_1(Y,\bdy Y)$
  is
  replaced by
  $H_1(\bSigma\times[0,1],\balphas\times\{0\}\cup\betas\times\{1\})$,
  and the group $H_1(Y,\bdy Y)$ is
  replaced by $H_1(Y)\cong H^2(Y,\bdy Y)$.
\end{proof}

For $\s$ a $\spin^c$ structure on $Y$, define $\gls*{generatorsS} \subset
\S(\HD)$ to be those generators~$\x$ so that $\s_z(\x) = \s$.

\section{Admissibility criteria}
\label{sec:Admissibility}
As in~\cite{OS04:HolomorphicDisks}, we need certain additional
hypothesis to ensure the definitions of the differentials involve only
finite sums. The idea is as follows. Given a homology class
$B\in\pi_2(\x,\y)$, in Chapter~\ref{chap:structure-moduli} we will
construct a moduli space of holomorphic curves in the homology
class~$B$. In Chapters~\ref{chap:type-d-mod}
and~\ref{chap:type-a-mod}, we
will use counts these moduli spaces, when they are $0$-dimensional, to
define the invariants of bordered $3$-manifolds. For each $B$, such a
count will be finite. But we will sum over all $B$ (perhaps with given
$\bdy^\bdy$), so we must ensure
that only finitely many such $B$ can contribute. We will use the fact (Lemma~\ref{lem:holo-has-pos-domain}) that
if the moduli space of curves in class $B$ is non-empty then $B$ must
be positive, in the following sense:
\begin{definition}\label{def:positive}%
\index{positive homology class}\index{domain!positive}\index{homology class!positive}%
  A homology class $B\in\pi_2(\x,\y)$ is called \emph{positive} if the
  coefficient in the domain of $B$ of every component of
  $\bSigma\setminus(\balphas\cup\betas)$ is non-negative.
\end{definition}

We will have two different analogues of the
``weak admissibility'' criterion of \cite[Definition 4.10]{OS04:HolomorphicDisks},
one of which we call provincial
admissibility and the other of which we call simply
admissibility. We will see in
Lemmas~\ref{prop:provincial-admis-finiteness}
and~\ref{prop:admis-finiteness} that these
conditions ensure that the set of positive homology classes with
given $\bdy^\bdy$ (respectively all positive homology classes) is finite.

In Chapters~\ref{chap:type-d-mod} and \ref{chap:type-a-mod}, provincial
admissibility will be sufficient to define the invariants of bordered
$3$-manifolds. However, in order for the tensor product involved in
gluing two bordered $3$-manifolds to make sense, one of the two
manifolds has to be admissible and not just provincially
admissible.

\begin{definition}\label{def:provincial-admissibility}A bordered Heegaard diagram
  is called \emph{provincially admissible}
  if every provincial periodic domain has both positive and negative
  coefficients.
  \index{admissible!Heegaard diagram!provincially}%
\end{definition}

\begin{definition}\label{def:admissibility}A bordered Heegaard diagram
  is called \emph{admissible} if every
  periodic domain has both positive and negative coefficients.
  \index{admissible!Heegaard diagram}%
\end{definition}

\begin{proposition}\label{prop:admis-achieve-maintain}Every bordered Heegaard diagram is isotopic to an
  admissible bordered Heegaard diagram.  Further, any two admissible
  bordered Heegaard diagrams can be connected though a sequence of
  Heegaard moves in which every intermediate Heegaard diagram is
  admissible. The same statements hold if ``admissible'' is replaced
  by ``provincially admissible.''
\end{proposition}

\begin{proof}
	This follows by winding transverse to the $\beta$-circles, as
	in the case of closed manifolds~\cite[Section 5]{OS04:HolomorphicDisks}.
\end{proof}

As in the closed case, there are reformulations of the two
admissibility criteria:
\begin{lemma}\label{lemma:admiss-reform}A bordered Heegaard diagram
  $(\Sigma,\alphas,\betas)$ is admissible if and only if there is an
  area function $A$ on $\Sigma$ such that $A(P)=0$ for any periodic
  domain~$P$. A bordered Heegaard diagram $(\Sigma,\alphas,\betas)$
is provincially admissible if and only if there is an area
function~$A$ on $\Sigma$ such that $A(P)=0$ for any provincial
periodic domain~$P$.
  \index{admissible!Heegaard diagram!alternate characterization of}%
  \index{area function}%
\end{lemma}
\begin{proof}
The ``only if'' direction is immediate.  For the ``if'' direction,
the proof follows exactly as in the closed case, see
\cite[Lemma 4.12]{OS04:HolomorphicDisks}. For the
reader's benefit, we recast the result here in terms of Farkas's lemma 
from the theory of convex sets~\cite{Farkas}:
\index{Farkas's lemma}%
\begin{lemma}Let $V$ be a vector space and $\{p_i\}_{i=1}^N\subset V$ a finite
  set of vectors. Then either
  \begin{itemize}
  \item there is a nonzero linear functional $\ell\in V^*$ such that
    $\ell(p_i)\geq 0$ for all~$i$ or
  \item there are constants $c_i> 0$ such that $\sum_{i=1}^{N}c_{i}p_{i}=0$.
  \end{itemize}
\end{lemma}

Now, to prove the first statement, enumerate the regions
$R_1,\dots,R_N$ of $\Sigma$. Take $V^*$ to be the vector space
generated by the periodic domains in $\Sigma$, and $p_i\in V=V^{**}$
to be evaluation at the $i\th$ region $R_i$. By assumption, there is
no nonzero $P\in V^*$ such that $p_i(V)\geq 0$ for all
$i$. Consequently, there are constants $c_i>0$ such that $\sum_{i-1}^N
c_i p_i=0$. But then setting the area of $R_i$ to be $c_i$ is the
desired area function.

The second statement follows by just the same reasoning, using only
the provincial periodic domains.
\end{proof}

\begin{proposition}\label{prop:provincial-admis-finiteness}Suppose
  that $\HD$ is provincially admissible, as in
  Definition~\ref{def:provincial-admissibility}. Fix generators
  $\x,\y\in\S(\HD)$ and an element $h\in
  H_1(Z',\CircPts)$. Then there are finitely many
  $B\in\pi_2(\x,\y)$ such that $\bdy^\bdy B=h$ and all coefficients
  of~$B$ are non-negative.
\end{proposition}
\begin{proof}
  Let $A$ be an area form as guaranteed by
  Lemma~\ref{lemma:admiss-reform}. If $B,B'\in\pi_2(\x,\y)$ with
  $\bdy^\bdy B=\bdy^\bdy B'=h$ then $B-B'$ is a provincial periodic
  domain. Consequently the areas $A(B)$ and $A(B')$ are equal. But
  there are only finitely many positive domains of any given area.
\end{proof}

\begin{proposition}\label{prop:admis-finiteness}Suppose that $\HD$ is
  admissible, as in
  Definition~\ref{def:admissibility}. There are finitely many
  generators $\x, \y \in \S(\HD)$ and
  $B\in\pi_2(\x,\y)$ such that all coefficients of~$B$ are
  non-negative.
\end{proposition}
\begin{proof}
  There are only finitely many generators, so we may first fix $\x$ and~$\y$.
  Let $A$ be an area form as guaranteed by
  Lemma~\ref{lemma:admiss-reform}. If $B,B'\in\pi_2(\x,\y)$ then
  $B-B'$ is a periodic domain. Consequently the areas $A(B)$ and
  $A(B')$ are equal. But there are only finitely many positive domains
  of any given area.
\end{proof}

\section{Closed diagrams}
\label{sec:GluingDiagrams}
In this section we discuss the behavior of (bordered) Heegaard
diagrams, generators and homology classes under cutting and gluing.

We will use minus signs to denote orientation reversal, so, for
\index{orientation reversal}%
instance, $\gls*{minus}Y$ is $Y$ with its orientation reversed and $-\HD$ is
$\HD=(\bSigma,\balphas,\betas,z)$ with the orientation of $\bSigma$
reversed. Note that $Y(-\HD)\cong-Y(\HD)$ and $F(-\PMC)\cong-F(\PMC)$,
canonically.

Gluing bordered Heegaard diagrams corresponds to gluing bordered
$3$\hyp manifolds:
\index{gluing!bordered Heegaard diagrams}%
\index{gluing!bordered $3$-manifolds}%
\index{Heegaard diagram!bordered!gluing}%
\index{bordered!$3$-manifold!gluing}%
\begin{lemma}\label{lem:glue-hd-glue-Y}
  Let $\HD_1=(\bSigma_1,\balphas_1,\betas_1,z_1)$ and
  $\HD_2=(\bSigma_1,\balphas_1,\betas_1,z_2)$ be bordered Heegaard
  diagrams for bordered $3$-manifolds $(Y_1,\PMC,\phi_1)$ and
  $(Y_2,\PMC,\phi_2)$. Obtain a closed Heegaard diagram
  $-\HD_1\gls*{gluealongbdy} \HD_2$ by gluing $-\HD_1$ and $\HD_2$ along their
  boundaries, according to the markings from their pointed matched circles.
  \[
  -\HD_1\cup_\bdy \HD_2=(\bSigma_1\cup_\bdy\bSigma_2, \balphas_1\cup_\bdy\balphas_2,\betas_1\cup\betas_2,z_1=z_2).
  \]
  Then $-\HD_1\cup_\bdy\HD_2$ represents $-Y_1\cup_\bdy
  Y_2=-Y_1\cup_{\phi_2\circ\phi_1^{-1}}Y_2$.
\end{lemma}
\begin{proof}
  This is immediate from the Morse theory description.
\end{proof}

Cutting bordered Heegaard diagrams corresponds to cutting bordered
$3$\hyp manifolds:
\index{cutting!bordered Heegaard diagrams}%
\index{cutting!bordered $3$-manifolds}%
\index{Heegaard diagram!bordered!cutting}%
\index{bordered!$3$-manifold!cutting}%
\begin{lemma}\label{lem:cut-hd-cut-Y}
  Let $\HD=(\Sigma,\alphas,\betas,z)$ be a (closed) Heegaard diagram
  and $Z$ a  separating curve such that
  \begin{itemize}
  \item $Z$ is disjoint from the $\beta$-circles (so that
    $\Sigma\setminus(\betas\cup Z)$ has exactly two components),
  \item $\Sigma\setminus(\alphas\cup Z)$ has exactly two components, and
  \item $Z$ passes through the basepoint.
  \end{itemize}
  Write $\Sigma=\bSigma_L\cup_Z\bSigma_R$. Then
  $\HD_L=(\bSigma_L,\alphas\cap\bSigma_L,\betas\cap\bSigma_L,z)$ and
  $\HD_R=(\bSigma_R,\alphas\cap\bSigma_R,\betas\cap\bSigma_R,z)$ are
  bordered Heegaard diagrams. Further, there is a separating surface $F\subset
  Y(\HD)$ so that $Y(\HD_L)$ and $Y(\HD_R)$ are the two components of
  $Y(\HD)\setminus F$.
\end{lemma}
\begin{proof}
  It is immediate from the hypotheses that $\HD_L$ and $\HD_R$ are
  bordered Heegaard diagrams.  The fact that there is a surface
  $F\subset Y(\HD)$ as in the statement of the lemma follows from
  Lemma~\ref{lem:glue-hd-glue-Y}. Alternatively, $F$ is the union of
  the flow lines through $Z$ with respect to a Morse function and
  metric inducing $\HD$.
\end{proof}

Next we discuss how generators and homology classes behave under these
cutting operations. Let $\HD_1$ and $\HD_2$ be bordered Heegaard
diagrams with $\bdy\HD_1=-\bdy\HD_2$.
Fix generators $\x_1\in\S(\HD_1)$ and $\x_2\in\S(\HD_2)$.  If the
corresponding sets of occupied $\alpha$-arcs are
disjoint, i.e., $o(\x_1)\cap o(\x_2)=\emptyset$, then the union
$\x_1\cup\x_2$ can be viewed as a generator in $\S(\HD)$
for the Heegaard Floer complex of the closed manifold~$Y$. We will
call pairs of generators $(\x_1,\x_2)$ such that $o(\x_1\cap
o(\x_2)=\emptyset$ \emph{compatible pairs}.
\index{generator!compatible pair of}%
\index{compatible pair!of generators}%
\index{gluing!generators}%
\index{generator!gluing g.'s}%

\begin{lemma}
  \label{lem:HoClassFibProd}
  Let $({\mathcal H}_1,z)=(\Sigma_1,\alphas_1,\betas_1,z)$ and
  $({\mathcal H}_2,z)=(\Sigma_2,\alphas_2,\betas_2,z)$ be two bordered
  Heegaard diagrams for the bordered three-manifolds
  $(Y_1,-\PMC,\phi_1)$ and $(Y_2,\PMC,\phi_2)$, with compatibly marked
  boundaries.
  Given $\x_1,\y_1\in\S({\mathcal H}_1)$ and $\x_2,\y_2\in\S({\mathcal H}_2)$,
  we have a natural identification of
  $\pi_2(\x_1\cup\x_2,\y_1\cup\y_2)$ with the
  subset of
  $\pi_2(\x_1,\y_1)\times \pi_2(\x_2,\y_2)$ consisting of pairs $(B_1, B_2)$
  with $\bdy^{\bdy}(B_1)+\bdy^{\bdy}(B_2) = 0$. 
  Moreover, given $p\in\Sigma_2\setminus(\alphas\cup\betas)$
  and $B_i\in\pi_2(\x_i,\y_i)$ for $i=1,2$ which agree along
  the boundary, the local multiplicity of $B_2$ at $p$ coincides
  with the local multiplicity of the glued homology class 
  in $\pi_2(\x_1\cup\x_2,\y_1\cup\y_2)$ (under
  the above identification) at the point $p$, now thought
  of as a point in $\Sigma=\Sigma_1\cup\Sigma_2$.
\index{gluing!homology classes}%
\index{homology class!gluing}%
\end{lemma}

\begin{proof}
	The proof is straightforward.
\end{proof}

For $B_1$ and $B_2$ which agree along the boundary, let $B_1
\gls*{gluehomolclass} B_2$
be the homology class in $\pi_2(\x_1\cup\x_2,\y_1\cup\y_2)$ guaranteed
by Lemma~\ref{lem:HoClassFibProd}.

Finally, we give a criterion for when a closed Heegaard diagram is
admissible.
\begin{lemma}\label{lem:closed-admissible}
  Let $\HD_1$ and $\HD_2$ be two bordered Heegaard diagrams with
  $\bdy\HD_1 = -\bdy\HD_2$.  If $\HD_1$ is admissible and $\HD_2$ is
  provincially admissible, then $\HD = \HD_1 \cup_\bdy \HD_2$ is
  admissible.
\end{lemma}
\begin{proof}
  Let $B_1\glue B_2$ be a positive periodic domain for $\HD$.  Since
  $\HD_1$ is admissible, $B_1 = 0$ and so $\bdy^\bdy B_2 = 0$.  Then
  since $\HD_2$ is provincially admissible, $B_2=0$.
\end{proof}


\chapter{Moduli spaces}
\label{chap:structure-moduli}
Let $Y$ be a bordered $3$-manifold. In order to define invariants of
$Y$ we will choose a bordered Heegaard diagram
$\Heegaard=(\bSigma,\balphas,\betas,z)$ for $Y$ and count
$J$\hyp holomorphic curves in $\Sigma \times [0,1] \times \RR$ with
boundary conditions coming from $\balphas$ and $\betas$.  In this
chapter we give the technical results
concerning the moduli spaces of holomorphic curves that we use to
construct the invariants $\CFDa(\Heegaard)$ and $\CFAa(\Heegaard)$. We
begin in Section~\ref{sec:moduli-overview} with an outline of the main
results that we will need and the
ideas in the proof; the precise mathematics begins in
Section~\ref{sec:curves-in-sigma} below.

\section{Overview of the moduli spaces}\label{sec:moduli-overview}
In this chapter, we will
discuss $J$-holomorphic curves\footnote{With respect to an
  appropriate $J$; see Definition~\ref{def:admissible_J} below.}
\[
u\co (S,\bdy S)\to \left(\Sigma\times[0,1]\times\RR,
(\alphas\times\{1\}\times\RR)\cup(\betas\times\{0\}\times\RR)\right)
\]
where $S$ is a Riemann surface with boundary and punctures.  Here we
view $\Sigma$, the interior of~$\bSigma$, as a Riemann surface with
one puncture, denoted~$p$.  As
in Section~\ref{sec:homology-classes-generators}, the coordinate on
the $[0,1]$ factor is denoted~$s$ and the
coordinate on the $\RR$ factor is denoted~$\gls*{t}$.

We view the manifold $\Sigma\times[0,1]\times\RR$ as having three
different infinities: $\Sigma\times[0,1]\times\{+\infty\}$,
\index{infinity!of $\Sigma\times[0,1]\times\RR$}%
\index{infinity!east}%
\index{infinity!$\pm$}%
$\Sigma\times[0,1]\times\{-\infty\}$ and $p\times[0,1]\times\RR$,
which we refer to as $\gls*{plusinfty}$, $\gls*{minusinfty}$ and $\gls*{eastinfty}$ respectively.
(Here, $e$ stands for ``east''.)  The asymptotics of the holomorphic
curves~$u$ we consider are:
\begin{itemize}
\item At $\pm\infty$, $u$ is asymptotic to
  a $g$-tuple of chords of the form $x_i\times[0,1] \times \{\pm\infty\}$, where
  $\x=\{x_i\}_{i=1}^g$ is a generator in the sense of
  Definition~\ref{def:generator}.
\item At $e\infty$, $u$ is asymptotic to a finite collection of Reeb chords
  $\{\rho_i\times \{(1, t_i)\}\}$. Here $\rho_i$ is a Reeb chord
  in the ideal contact boundary $\bdy\bSigma$ (also called~$Z$) with
  endpoints on
  $\CircPts=\alphas\cap Z$
  (see Section~\ref{sec:reeb-chords-def}).  We refer to $\gls*{tsubi}$
  as the \emph{height} of the chord~$\gls*{rhosubi}$.
  \index{Reeb chord}%
  \index{height}%
\end{itemize}
\index{puncture!east}%
The set of east punctures of $u$ is
partially ordered by the heights of the corresponding Reeb chords.
\index{ordering of Reeb chords at $e\infty$}%
\index{ordered partition}%
\index{partition!ordered|see{ordered partition}}%
This induces an
\emph{ordered partition} $\gls*{OrderedPartition}=(P_1,P_2,\dots,P_\ell)$ of the east
punctures,
where each $\gls*{Part}$ consists of those punctures occurring at one
particular height.  (For $\CFDa(\Heegaard)$ we only need to consider
discrete partitions where each $P_i$ has only one puncture.)

With this background, we consider the moduli space
\[
\gls*{ModSpaceOrderedParam}
\]
\index{moduli space}%
(Definition~\ref{def:tcM}),
which consists of curves from the \emph{decorated source}~$\gls*{DecSource}$,
\index{source!decorated}%
\index{decorated source|see{source, decorated}}%
asymptotic to $\x$ at $-\infty$ and $\y$ and $+\infty$, in the
homology class~$\gls*{homolclass}$ (as in Section~\ref{sec:homology-classes-generators})
and respecting
the ordered partition~$\vec{P}$ of the Reeb chords.  Here the
source surface~$\gls*{SourceSurf}$ is decorated with labels describing its
asymptotics (Definition~\ref{def:decorated-source}).  We also
consider the reduced moduli space, its quotient
\[
\gls*{ModSpaceOrderedUnparam}
\coloneqq\tcM^B(\x,\y;\Source;\vec{P})/\RR
\]
\index{moduli space!reduced}\index{reduced!moduli space}%
by translation in the $t$ coordinate.

These moduli spaces are smooth manifolds
(Proposition~\ref{prop:transversality}) whose dimensions are easy to
compute (Equation~(\ref{eq:Index})).  In fact, for an
\emph{embedded} curve the dimension is independent of the source~$S$
and depends only on the asymptotics and homology class of the
curve (Section~\ref{sec:expected-dimensions}); this 
gives rise to the gradings on the algebra and on the modules.

As is typical in symplectic geometry, our invariants are defined by
looking only at \emph{rigid} moduli spaces, those
$\cM^B(\x,\y;\Source;\vec{P})$ which are
$0$-dimensional.
\index{moduli space!rigid}%
To show that $\bdy^2 = 0$ and for proofs of
invariance, we also need to consider $1$-dimensional moduli spaces,
and, in particular, we need appropriate compactifications.  These are
provided by the spaces $\gls*{ModSpaceOrderedCpct}$,
which include
\emph{holomorphic combs}. The idea is that, as in symplectic field
\index{holomorphic comb}%
\index{comb|see{holomorphic comb}}%
theory, some parts of a holomorphic curve may go to infinity relative
to other parts. In our case, this can result in a degeneration at
\index{degeneration of holomorphic curves}
$\pm\infty$ or at $e\infty$. When a curve degenerates at $\pm\infty$, the
result is two or more curves of the kind we have already
discussed. A degeneration at $e\infty$ creates a new kind of object: a
curve in
$\RR\times Z\times[0,1]\times\RR$. Section~\ref{sec:curves_at_east_infinity}
is devoted to discussing these curves at east infinity.  This
compactification is defined in Definition~\ref{def:ocM}, and a
schematic illustration of a holomorphic comb can be found in
Figure~\ref{fig:comb_schematic}.

\begin{figure}
\includegraphics[scale=.7143]{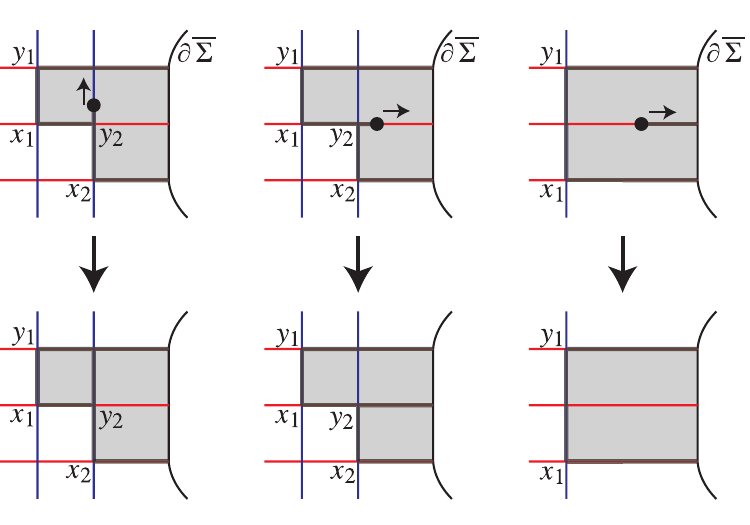}
\caption[Examples of three kinds of codimension $1$
  degenerations]{\textbf{Examples of three kinds of codimension $1$
  degenerations.} The large black dot represents a boundary branch
  point of $\pi_\Sigma\circ u$. Left: degenerating into a height $2$
  holomorphic comb. Center: degenerating a join curve. Right: degenerating a
  split curve. The diagrams show the projection of the curve in
  $\Sigma \times [0,1] \times \RR$ to $\Sigma$.}\label{fig:degen_examples}
\end{figure}

The main result we need for our purposes is a classification of the
codimension~1 degenerations
(Section~\ref{sec:degenerations-holomorphic-curves}).
Proposition~\ref{prop:restrict_degens_1} and
Lemma~\ref{lemma:NoBoundaryDoublePoints} say that for the curves we
consider, these degenerations come in four kinds:
\begin{itemize}
\item Degenerating into a height $2$ holomorphic comb at $\pm\infty$. These are
  \index{holomorphic comb!height two}%
  exactly the kinds of degenerations that occur in Heegaard
  Floer homology for closed three-manifolds.
\item A boundary branch point of $\pi_\Sigma\circ u$ can go off to $e\infty$,
  in the process splitting a Reeb chord in two. Below, this is
  referred to as degenerating a ``join curve'' at $e\infty$.
  \index{curve!join}%
  (From the point of view of the curve at $e\infty$, two Reeb chords
  come together and merge as we travel further away from the main surface.)
\item A collapse in the ordering of $\vec{P}$, i.e., the heights of
  two parts of $\vec{P}$ coming together. In the process, some
  ``split curves'' can degenerate at $e\infty$. These again
  \index{curve!split}%
  correspond to boundary branch points of $\pi_\Sigma\circ u$
  approaching~$e\infty$, but in a combinatorially different way from
  join curves.
  (In the curve at $e\infty$, one (or more) Reeb chord
  splits apart as we travel further away from the main surface.)
\item Degenerating a ``shuffle curve'' at east $\infty$. This
  \index{curve!shuffle}%
  degeneration involves either an interior branch point or two
  boundary branch points approaching east $\infty$, and is explained
  in more detail at the end of
  Section~\ref{sec:curves_at_east_infinity}; see in particular
  Figure~\ref{fig:shuffle2} for a domain in which a shuffle curve
  degenerates. This type of degeneration is not relevant for $\CFDa$,
  only for $\CFAa$.
  (In the curve at $e\infty$, two Reeb chords
  recombine with each other to create two different Reeb chords as we
  travel further away from the main surface.)
\end{itemize}
Figure~\ref{fig:degen_examples} shows examples in which the first three
kinds of degenerations occur. 

Our next goal is to obtain an identity relating the number of
degenerations of the four types above.  The usual way to do this is to
show that if the expected dimension $\ind(\x,\y;\Source;\vec{P})$ is~$2$,
then the space $\ocM^B(\x,\y\semico\Source\semico\vec{P})$ is a compact
one-manifold with boundary given by all holomorphic combs with the correct
combinatorics to occur as one of these degenerations.  The identity
would then follow from the fact that a compact one-manifold has an
even number of endpoints.

Unfortunately, available gluing technology seems ill equipped for this
approach, the difficulty being with split and shuffle curve
degenerations. (See Remark~\ref{remark:gluing-is-hard}
  and the beginning of Section~\ref{sec:codim-one-bound} for more
  discussion of this point.)
Instead, we use a somewhat indirect argument to show that there are an
even number of holomorphic combs
of the types that may appear as degenerations in an index~$2$ moduli
space (Theorem~\ref{thm:master_equation}).
This is all
that we will need later to define the invariants.

This chapter is organized as follows.
Section~\ref{sec:curves-in-sigma} defines the moduli spaces of
curves in $\Sigma \times [0,1] \times \RR$ of interest to us, and
Section~\ref{sec:curves_at_east_infinity} the moduli spaces for
curves at east infinity.  In Section~\ref{sec:combs-compact}, we put
these moduli spaces together to define the compactified moduli spaces
of so-called holomorphic combs and prove that they form
compactifications. Next, in Section~\ref{sec:combs-gluing} we prove
the gluing results for combs needed for the rest of the paper.  In
Section~\ref{sec:degenerations-holomorphic-curves} we restrict what
kind of degenerations can occur in codimension-one; in particular, we
show that only those degenerations for which we proved gluing results
are possible, leading to Theorem~\ref{thm:master_equation}.  Finally,
in Section~\ref{sec:expected-dimensions} we study the index of the
$\overline{\bdy}$-operator at an embedded curve, and use this to give
further restrictions on the codimension-one boundaries of moduli
spaces of embedded curves.

Fortified by this overview, we turn to the mathematics itself.

\section{Holomorphic curves in
  \textalt{$\Sigma\times[0,1]\times\RR$}{Sigma \texttimes\ [0,1] \texttimes\ R}}\label{sec:curves-in-sigma}
In this section, we will define various moduli spaces of
holomorphic curves in $\Sigma\times[0,1]\times\RR$ with names like
$\cM^B(\x,\y;\Source)$, $\cM^B(\x,\y;\Source;P)$ and
$\cM^B(\x,\y;\Source;\vec{P})$. In Chapters~\ref{chap:type-d-mod}
and~\ref{chap:type-a-mod} we will use the moduli spaces
$\cM^B(\x,\y;\Source;\vec{P})$ to define our invariants of
bordered $3$-manifolds. These are, consequently, the spaces in
which we will ultimately be interested. However, it is easier to
formulate several technical results by describing the spaces
$\cM^B(\x,\y;\Source;\vec{P})$ as subspaces of $\cM^B(\x,\y;\Source)$
and $\cM^B(\x,\y;\Source;P)$, which we introduce first.

As before, let $\HD$ be a bordered Heegaard diagram of
genus~$g$ representing a manifold~$Y$ with boundary of genus~$k$.
Choose a
symplectic form~$\gls*{SigmaSymp}$ on~$\gls*{Sigma}$, %
\index{symplectic form!on $\Sigma$}%
\index{cylindrical end}%
\index{complex structure!on $\Sigma$}%
\index{almost complex structure|see{complex structure}}%
with respect to which the
boundary of
$\Sigma$ is a cylindrical end.
Let $\gls*{SigmaCx}$ be a complex
structure on~$\Sigma$ compatible with~$\omega_\Sigma$.
We will further assume that the
$\alpha$-arcs $\gls*{noalphasa}$ are cylindrical near $p$ in the following sense. Fix a punctured neighborhood
$\gls*{Usubp}$ of $p$, and a symplectic identification $\gls*{pphi}:U_p\stackrel{\cong}{\to}
S^1\times(0,\infty)\subset T^*S^1$. Then we assume that both $j_\Sigma|_{U_p}$ and
$\phi(\alpha_i^a\cap U_p)$ are invariant with respect to the $\RR$-translation action
in $S^1\times(0,\infty)$.
With respect to $j_\Sigma$, the boundary is a puncture, 
\index{puncture!boundary of $\Sigma$ is a}%
which we denote~$\gls*{p}$.

Let $\gls*{Sigmae}$ be the result of filling in the puncture of $\Sigma$, so
$j_\Sigma$ induces a complex structure on $\Sigma_{\widebar e}$.
\index{complex structure!on $\Sigma_{\widebar e}$}%

We consider moduli spaces of
holomorphic curves in $\Sigma\times[0,1]\times\RR$. Let
\begin{align*}
  \gls*{piSigma}&\co\Sigma\times[0,1]\times\RR\to\Sigma,\\
  \gls*{piD}&\co\Sigma\times[0,1]\times\RR\to[0,1]\times\RR,\\
  \gls*{s}&\co\Sigma\times[0,1]\times\RR\to[0,1]\textrm{, and }\\
  \gls*{t}&\co\Sigma\times[0,1]\times\RR\to\RR  
\end{align*}
denote the projections. We equip $[0,1]\times\RR$ with the symplectic
form $\gls*{DSymp}\coloneqq ds\wedge dt$ 
\index{symplectic form!on $\DD$}%
and the complex structure
$\gls*{DCx}$ 
\index{complex structure!on $\DD$}%
with $j_\DD\bigl(\frac{\partial}{\partial s}\bigr) =
\frac{\partial}{\partial t}$, and equip $\Sigma\times[0,1]\times\RR$
with the split symplectic form $\pi_\Sigma^*(\omega_\Sigma) +
\pi_\DD^*(\omega_\DD)$.
\index{symplectic form!on $\Sigma\times[0,1]\times\RR$}%
Note also that there is an $\RR$-action on
$\Sigma\times[0,1]\times\RR$ by translation on the $t$~coordinate.
\index{$\RR$-action!on $\Sigma\times[0,1]\times\RR$}%

\begin{definition}\label{def:admissible_J}We say an
almost complex structure $J$ on $\Sigma\times[0,1]\times\RR$ is \emph{admissible}
if the following conditions are satisfied:
\index{complex structure!on $\Sigma\times[0,1]\times\RR$}%
\index{complex structure!admissible}\index{admissible!complex structure|see{complex structure, admissible}}%
\begin{enumerate}[label=(J-\arabic*),ref=J-\arabic*]
\item\index{(J1)--(J4)}\label{item:J1}The projection map $\pi_\DD$ is $J$-holomorphic.
\item\label{item:J2} For $\frac{\partial}{\partial s}$ and
  $\frac{\partial}{\partial t}$ the vector fields tangent to the
  fibers of $\pi_\Sigma$ induced by $s$ and $t$ respectively,
  $J\frac{\partial}{\partial s}=\frac{\partial}{\partial t}$ (and
  consequently the fibers of $\pi_\Sigma$ are $J$-holomorphic).
\item The $\RR$-action is $J$-holomorphic.
\item\index{(J1)--(J4)}\label{item:J4}The complex structure is split,
  i.e., $J=j_\Sigma\times j_\DD$, over some (fixed, $\RR$-invariant) neighborhood of
  $p\times[0,1]\times\RR$.
\end{enumerate}
\end{definition}
In point~\ref{item:J4} the neighborhood will be chosen small enough
that its closure is disjoint from the $\beta$-circles, and so that the
$\alpha$-arcs intersect it in lines (with respect to a holomorphic
identification with $\DD^2\setminus\{0\}$); any such neighborhood will
do. One example of an admissible almost complex structure is the
split complex structure $j_{\Sigma}\times j_{\DD}$.

Fix now an admissible $J$; later, we will assume $J$ is generic in certain senses.

\begin{definition}\label{def:decorated-source}
  \index{source!decorated}%
  By a \emph{decorated source $\gls*{DecSource}$} we mean:
\begin{enumerate}
\item \label{item:smooth_source}a
  smooth (not
  nodal) Riemann surface~$S$ with boundary and
  punctures on the boundary (i.e., a compact Riemann surface with boundary minus finitely many points on the boundary),
\item a labeling of each puncture of $S$ by one of $+$, $-$, or~$e$ and
  \index{puncture!east}%
\item a labeling of each $e$ puncture of $S$ by a Reeb chord in
  $(Z,\CircPts)$, as defined in
  Section~\ref{sec:reeb-chords-def}.
\end{enumerate}
We view two decorated sources as equivalent if there is an orientation-preserving diffeomorphism between them which respects the labelings of the punctures.
\end{definition}

Given a decorated source $\Source$, let $S_{\widebar e}$ denote the result of
filling in the $e$ punctures of $S$.
We will consider maps
\[
u\co (S,\partial S)\to (\Sigma\times[0,1]\times\RR,(\alphas\times\{1\}\times\RR)\cup(\betas\times\{0\}\times\RR))
\]
such that
\begin{enumerate}[label=(M-\arabic*),ref=M-\arabic*]
\item\index{(M1)--(M11)}\label{item:moduliFirst} The map $u$ is $(j,J)$-holomorphic with respect to some
  almost complex structure $j$ on $S$.
\item\label{item:moduliSecond} The map $u\co S\to\Sigma\times[0,1]\times\RR$ is proper.
\item The map $u$ extends to a proper map $u_{\widebar e}\co
  S_{\widebar e}\to\Sigma_{\widebar e}\times[0,1]\times\RR$.
\item  The map
  $u_{\widebar e}$ has finite energy
  \index{energy}%
  in the sense of Bourgeois,
  Eliashberg, Hofer, Wysocki and
  Zehnder~\cite{BEHWZ03:CompactnessInSFT} (see also Section~\ref{sec:combs-compact}).
\item \label{item:moduliNonconstant} $\pi_\DD\circ u_{\widebar e}$ is a $g$-fold branched
  cover. In particular, $\pi_\DD\circ u$ is non-constant on every component of $S$.
\item At each $-$-puncture $q$ of $S$, $\lim_{z\to q}(t\circ u)(z)=-\infty$.
\item At each $+$-puncture $q$ of $S$, $\lim_{z\to q}(t\circ u)(z)=+\infty$.
\item At each $e$ puncture $q$ of $S$, $\lim_{z\to q}(\pi_\Sigma\circ u)(z)$
  is the Reeb chord $\rho$ labeling $q$.
\item \label{item:moduliPenUlt} $\pi_\Sigma\circ u$ does not cover the region of
  $\Sigma$ adjacent to $z$.
\item \label{item:moduliLast}\label{item:moduliWeakBndryMntncty}
  For each $t\in\RR$ and each $i=1,\dots,g$,
  $u^{-1}(\beta_i\times\{0\}\times\{t\})$ consists of exactly one point. Similarly,
  for each $t\in\RR$ and each $i=1,\dots,g-k$,
  $u^{-1}(\alpha_i^c\times\{1\}\times\{t\})$ consists of exactly one point.
\end{enumerate}

We call Condition~(\ref{item:moduliWeakBndryMntncty}) \emph{weak boundary
  monotonicity.} 
\index{monotonicity|see{boundary monotonicity}}%
\index{boundary monotonicity!weak}%
Sometimes
later we will impose the following additional condition, which we call \emph{strong
  boundary monotonicity}:
\begin{enumerate}[resume*]
\item\index{(M1)--(M11)}\label{item:moduliStrongBndyMntncty} For each $t\in\RR$ and each $i=1,\ldots,2k$,
  $u^{-1}(\alpha_i^a\times\{1\}\times\{t\})$ consists of at most one point.
\end{enumerate}
\index{boundary monotonicity!strong}%
However, our moduli spaces are easiest to define if we do not
initially impose this condition.

It follows from the
Conditions~(\ref{item:moduliFirst})--(\ref{item:moduliLast}) that at
$-\infty$, $u$ is asymptotic to a $g$-tuple of chords of the form $x_i\times[0,1]$,
where $x_i\in\alphas\cap\betas$. By weak boundary monotonicity, the set $\x=\{x_i\}$
has the property that exactly one point of $\x$ lies on each $\alpha$-circle, and one
point on each $\beta$-circle. However, more than one point of $\x$ may lie on the same
$\alpha$-arc. We call such a $g$-tuple a \emph{generalized generator.}
\index{generator!generalized}%
\index{generalized!generator}%
If the strong boundary monotonicity condition is also satisfied then at
most one $x_i$ lies on each $\alpha$-arc, and hence $\x$ is a generator in the sense
of Section~\ref{sec:homology-classes-generators}. Exactly analogous statements hold
at $+\infty$.

Any curve $\gls*{u}$ satisfying Conditions~(\ref{item:moduliFirst})--(\ref{item:moduliLast})
belongs to some homology class
\index{homology class!of curve}%
 $\gls*{bracu}\in\pi_2(\x,\y)$. (Here, $\pi_2(\x,\y)$ is the obvious
generalization of the notion from Section~\ref{sec:homology-classes-generators} to
the case of generalized generators.) We collect these curves into moduli spaces:
\begin{definition}
Given generalized generators $\x$ and $\y$ and a
homology class $B\in\pi_2(\x,\y)$, let
\[
\gls*{ModSpaceOpenParam}
\]
denote the moduli space of curves satisfying
Conditions~(\ref{item:moduliFirst})--(\ref{item:moduliLast}), with decorated source
$\Source$, asymptotic to $\x$ at $-\infty$ and $\y$ at $+\infty$ and
in the homology class~$B$.
\end{definition}
\index{moduli space}%

The domains of holomorphic curves (in the sense of Definition~\ref{def:Domain})
are always positive.
\begin{lemma}\label{lem:holo-has-pos-domain}
  If the moduli space $\tcM^B(\x,\y;\Source)$ is
  non-empty then the homology class $B$ is positive (in the sense of
  Definition~\ref{def:positive}).
\end{lemma}
\begin{proof}
  Let $u$ be an element of $\tcM^B(\x,\y;\Source)$.  The multiplicity
  of $B$ at a region $R$ is given by the intersection number $u\cdot
  (\{p\}\times[0,1]\times\RR)$. By Condition~(\ref{item:J2}) of
  Definition~\ref{def:admissible_J}, the fiber
  $\{p\}\times[0,1]\times\RR$ is $J$-holomorphic, so $u\cdot
  (\{p\}\times[0,1]\times\RR)$ is the intersection number between two
  $J$-holomorphic curves, which is positive by~\cite{McDuff94:positivity}
  or~\cite{MicallefWhite95:intersection-positivity}.
\end{proof}

Given $u\in\tcM^B(\x,\y;\Source)$, for each puncture $q$ of $S$ labeled by~$e$,
let $\gls*{evq}(u)=t\circ u_{\widebar e}(q)$. This gives an evaluation map
$\ev_q\co \tcM^B(\x,\y;\Source)\to\RR.$
Set 
\[
\gls*{ev}=\,\,\prod_{\mathclap{q \in E(\Source)}}\,\,\ev_q\co \tcM^B(\x,\y;\Source)\to\RR^{E(\Source)},
\]
where $\gls*{EofS}$ (or just $\gls*{Eof}$, when unambiguous) is the set of
east punctures of~$\Source$.

We next use the evaluation maps to define certain subspaces of
$\tcM^B(\x,\y;\Source)$. Fix a partition $\gls*{Partition}=\{\gls*{Part}\}$ of~$E$. By the
\emph{partial
diagonal $\gls*{PartialDiagonal}$ in $\RR^E$} 
\index{partial diagonal}%
\index{diagonal, partial}%
we mean the subspace of $\RR^E$ defined by the set of 
equations $\{\,x_p=x_q \mid P_i\in P,\,\,p,q\in P_i\,\}$. That is,
$\Delta_P$ is obtained by
requiring all coordinates in each part $P_i$ of $P$ to be equal. Now, we use the
cycle $\Delta_P$ to cut down the moduli spaces:
\begin{definition}\label{def:open-moduli-spaces}
  Given generalized generators $\x$ and $\y$, $B\in\pi_2(\x,\y)$,
  a decorated source $\Source$, and a partition $P$ of the $e$ punctures of
  $\Source$, let
\[
\gls*{ModSpaceUnorderedParam}
\coloneqq\ev^{-1}(\Delta_P)\subset \tcM^B(\x,\y;\Source).
\]
\end{definition}
Notice that $\tcM^B(\x,\y;\Source)$ is the special case of $\tcM^B(\x,\y;\Source;P)$ when
$P$ is the discrete partition, with one element per part.
\index{moduli space}%

We turn now to the meaning of genericity for $J$.
\begin{proposition}\label{prop:transversality}
  In the space of $C^\infty$ admissible almost complex
  structures there is a dense set of $J$ with the property that the
  moduli spaces $\tcM^B(\x,\y;\Source)$ are transversally cut out by
  the \index{transversality}%
  $\widebar{\partial}$-equations. Indeed, for any countable set
  $\{M_i\}$ of submanifolds of $\RR^E$, there is a dense set of
  admissible $J$ which satisfy the further property that
  $\ev\co\tcM^B(\x,\y;\Source)\to\RR^E$ is transverse to all of
  the~$M_i$.
\end{proposition}
This result follows from standard arguments; see~\cite[Chapter
3]{MS04:HolomorphicCurvesSymplecticTopology} for a nice explanation in
a slightly different context, or~\cite[Proposition
3.8]{Lipshitz06:CylindricalHF} for the proof of the first half of the
statement in a closely related setting. Briefly, one considers the Banach
manifold $\mathcal{B}$ of triples $(j_S,J,u)$ where $j_S$ is a point
in the Teichm\"uller space of complex structures on $S$, $J$ is a
$C^\ell$ complex structure on $\Sigma\times[0,1]\times\RR$ satisfying
Conditions (\ref{item:J1})--(\ref{item:J4}), and $u\co S\to
\Sigma\times[0,1]\times\RR$ is a $W^{1,p}_\delta$ map satisfying
Conditions~(\ref{item:moduliSecond})--(\ref{item:moduliLast}). (Here,
$W^{1,p}_\delta$ denotes a weighted Sobolev space; see, for
example,~\cite[Section (8h)]{SeidelBook}. 
We refer the reader to McDuff-Salamon~\cite[Section
3.2]{MS04:HolomorphicCurvesSymplecticTopology} for a discussion of the
regularity of these Banach manifolds and bundles. In particular, $\ell$ must be sufficiently large, depending on the Fredholm index at the homology class under consideration; since we will eventually work with $C^\infty$ almost complex structures we elide this point.)

The $\overline{\bdy}$ map
defines a section of a vector bundle $\mathcal{E}$ over $\mathcal{B}$,
whose fiber is the space of $(0,1)$-forms on $S$ valued in
$u^*T(\Sigma\times[0,1]\times\RR)$. Applying elliptic
regularity one sees that this section $\overline{\partial}$ is
transverse to the $0$-section of $\mathcal{E}$, so
$\overline{\partial}^{-1}(0)$ is a Banach manifold. (A key point in
this proof is verifying that for any connected holomorphic curve, either the projection to $\Sigma$ is constant or the projection
to $\Sigma\times[0,1]$ is somewhere injective; see~\cite[Lemma
3.3]{Lipshitz06:CylindricalHF} for that argument.)
An infinite-dimensional version of Sard's theorem~\cite[Theorem
1.3]{Smale65:Sard} guarantees that a
generic complex structure $J$ on $\Sigma\times[0,1]\times\RR$
(satisfying Conditions (\ref{item:J1})--(\ref{item:J4})) is a regular
value of the projection $\pi\co \overline{\partial}^{-1}(0)\to
\mathcal{J}$ to the space $\mathcal{J}$  of almost complex structures
satisfying Conditions (\ref{item:J1})--(\ref{item:J4}). The moduli
space $\tcM^B(\x,\y;\Source)$ is the quotient of
$\pi^{-1}(J)$ by the mapping class group.
Condition~(\ref{item:moduliNonconstant}) implies that the mapping
class group acts freely on $\pi^{-1}(J)$, so
$\tcM^B(\x,\y;\Source)$ is a smooth manifold of the expected dimension
whenever $\pi^{-1}(J)$ is transversally cut out. This
condition is what we mean by saying that $\tcM^B(\x,\y;\Source)$ is
transversally cut out.

For the second half of the statement, one verifies further that as a
map from $\overline{\partial}^{-1}(0)$ to $\RR^E$, evaluation at the
marked points is a submersion; the result then
follows. See~\cite[Section
3.4]{MS04:HolomorphicCurvesSymplecticTopology} for a more detailed
discussion of this part in a closely related setting.%

By Proposition~\ref{prop:transversality}, we can choose $J$ satisfying the following definition:
\begin{definition}\label{def:sufficiently-generic}%
\index{sufficiently generic!complex structure}%
\index{complex structure!sufficiently generic}%
  An admissible almost complex structure is called \emph{sufficiently
    generic} if the moduli spaces $\tcM^B(\x,\y;\Source;P)$ are
  transversely cut out for all choices of $\x$, $\y$, $B$, $\Source$
  and $P$.
\end{definition}
From now on, we will always work with sufficiently generic almost
complex structures that are $C^\infty$, unless we explicitly say otherwise.

The expected dimension of $\tcM^B(\x,\y;\Source;P)$ is not hard to
calculate:
\index{dimension|see{expected dimension}}\index{expected dimension|seealso{index}}%
\index{expected dimension!of $\tcM^B(\x,\y;\Source;P)$}%
\begin{proposition}\label{Prop:Index}The expected dimension $\ind(B,\Source,P)$ of
  $\tcM^B(\x,\y;\Source;P)$ is
\begin{equation}\label{eq:Index}
\gls*{indSource}=g-\chi(S)+2e(B)+|P|.
\end{equation}
Here $\gls*{PartitionLength}$ denotes the number of parts in~$P$ and $\gls*{Euler}$ is the
Euler measure of the domain of~$B$.
\end{proposition}

Here the \emph{Euler measure}\index{Euler measure}
of a region in $\Sigma \setminus
(\alphas \cup \betas)$ is its Euler characteristic minus $1/4$ the
number of corners (intersections of $\alpha$-curves and $\beta$-curves
or $\alpha$-arcs and $\bdy \Sigma$), and is defined to be additive
under union.  The Euler measure can be viewed as the integral of the Gaussian
curvature of the region with
respect to any Riemannian metric for which the boundary,
the $\alpha$-curves and the $\beta$-circles are all geodesics,
meeting at right angles at the corners.

\begin{proof}[Proof of Proposition~\ref{Prop:Index}]
  This is proved in~\cite[Proposition
  4.5.1]{Lipshitz06:BorderedHF}; we repeat it here for the reader's
  convenience. In brief, the Proposition follows from the
  corresponding index formula for closed Heegaard Floer
  homology,~\cite[Formula (6)]{Lipshitz06:CylindricalHF}, by a
  standard doubling argument at the east punctures, cf.~\cite[Section
  5]{Bourgeois02:MorseBott}. (The formula~\cite[Formula
  (6)]{Lipshitz06:CylindricalHF} itself follows from the Riemann-Roch
  theorem by another doubling argument.)

  Let $\bSigma^{\cap}$ denote the result of gluing a collar $
  [0, \epsilon)\times (\bdy\bSigma) $ to $\bdy \bSigma$. Inside $[0, \epsilon)\times (\bdy
  \bSigma \setminus z)$, choose arcs
  connecting every pair of (ends of) $\alpha$-arcs, as shown in
  Figure~\ref{Figure:Capping}. Call the new arcs $\balphas^{\cap}_0$,
  and let $\alphas^\cap=\balphas\cup\balphas^\cap_0$.  (This
  ``caps off'' the Reeb chords at east infinity.) We call the result
  the \emph{capped diagram}.\index{Heegaard diagram!capped} By
  attaching cylindrical ends to them, we can think of the
  $\balphas^\cap_0$ as arcs $\alphas^\cap_0$ lying in
  $\RR\times \bdy\bSigma$.

  For each Reeb chord $\rho_i$ there is a holomorphic disk
  $D_{\rho_i}$ in $( \RR\times \bdy \bSigma  \times [0, 1] \times \RR,
  \alphas^\cap_0 \times \{1\} \times \RR)$ asymptotic to
  $\gamma_i$. Given a map $u\in\cM^B(\x,\y;\Source)$ where the $e$
  punctures of $\Source$ are decorated by Reeb chords
  $\rho_{i_1},\dots,\rho_{i_n}$ one can preglue $u$ to
  $\bigcup_{j=1}^n D_{\rho_{i_j}}$ to obtain an approximately holomorphic map
  $u^\cap \co S^\cap \to \Sigma^{\cap} \times [0, 1] \times \RR$.
  
  \begin{figure}
    \begin{centering}
      \includegraphics[scale=.83333]{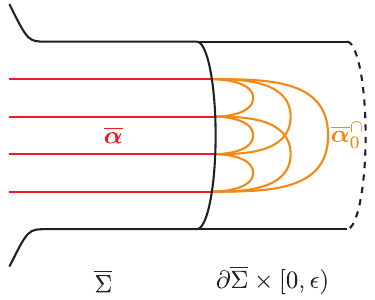}
      \caption[A capping operation]{\label{Figure:Capping}\textbf{A
          capping operation.} Part of the surface $\bSigma^\cap$ and the arcs
      $\balphas^\cap$ are shown; on the left side is $\bSigma$
      and $\balphas$, and on the right side is
      $[0,\epsilon)\times \bdy\bSigma$ and $\balphas^\cap_0$.}
    \end{centering}
  \end{figure}
  \colorused

  Each height constraint decreases both $|P|$ and the expected
  dimension by one. So, it suffices to prove the formula in the case
  that  $P$ is the discrete partition.
  Recall from~\cite[Formula (6), Section 4.1]{Lipshitz06:CylindricalHF}
  that in the closed case the index at a map $u \co S\to \Sigma \times
  [0, 1] \times \RR$ is given by $g - \chi (S) + 2 e(B)$.
  The closed formula applies to the linearization of the
  $\dbar$-operator at $u^\cap$, giving $\ind( u^\cap) = g - \chi
  (S^\cap) + 2 e (B^\cap) = g - \chi (S) + 2 e (B^\cap)$, where
  $B^\cap$ is the domain in $\bSigma^\cap$ corresponding to the domain
  $B$ in $\bSigma$. The index of $D \dbar$ at any of the
  disks $D_{\rho}$ is $1$. It follows that the index at $u$ is
  \[ 
  \ind (u) = g - \chi (S) + 2 ( e (B) + |P| / 2 ) - |P|+ |P|, 
   \]
  where the $|P| / 2$ comes from the effect of capping on the Euler measure,
  the $- |P|$ from the fact that $\ind (D_{\rho}) = 1$, and the $+ |P|$
  comes from the matching conditions at the punctures (because of the
  Morse-Bott nature of the gluing). This proves the proposition.
\end{proof}

Suppose $u\in\tcM^B(\x,\y;\Source;P)$. Then the $t$-coordinate of $u$
induces a partial ordering on the partition $P$. This combinatorial
data allows us to divide $\tcM^B(\x,\y;\Source;P)$ into different
strata.
\begin{definition}\label{def:tcM} Given $\x$, $\y$, and $\Source$ as
  above and an ordered partition
  \index{partition!ordered}%
  $\gls*{OrderedPartition}$ of the $e$ punctures of $\Source$, with $P$ the
  corresponding unordered partition, let 
  $\gls*{ModSpaceOrderedParam}$
  denote the
  (open) subset of $\tcM^B(\x,\y;\Source;P)$ consisting of those holomorphic curves
  for which the ordering of~$P$ induced by~$t$ agrees with the
  ordering given in~$\vec{P}$.
\end{definition}
\index{moduli space}%
That is, for $u\in\tcM^B(\x,\y\semico\Source\semico P)$ and $\vec{P}=(P_1,\dots,P_m)$,
we have $u\in \tcM^B(\x,\y\semico\Source\semico\vec{P})$ if and only if for all $i<i'$,
$q\in P_i$ and $q'\in P_{i'}$ implies $t\circ u(q)< t\circ u(q')$.

Given a partition~$P$ of $E(\Source)$,
there is an associated set~$\gls*{PartitionAssocChords}$ of
multi-sets of Reeb chords, given by replacing
each puncture of $\Source$ in $P$ by the associated Reeb
chord. Similarly, given an ordered partition $\vec{P}$ of $E(\Source)$
there is an
associated sequence~$\gls*{PartitionOrderedAssocChords}$ of multi-sets of Reeb chords.

\index{$\RR$-action!on $\Sigma\times[0,1]\times\RR$}\index{$\RR$-action!on $\tcM^B$}%
The $\RR$-action on $\Sigma\times[0,1]\times\RR$ by
translation on~$\RR$ induces an $\RR$-action on each
$\tcM^B(\x,\y\semico\Source\semico P)$ and $\tcM^B(\x,\y\semico\Source\semico\vec{P})$.  
We denote this action by $\gls*{tauR}$ for $R \in \RR$.  
The action is free except
in the trivial case that $\Source$ consists of $g$ disks with two boundary
punctures each, and $B=0$.  We say that a curve is \emph{stable}
\index{stable!holomorphic curve}%
if it
is not this trivial case.
We are interested primarily in the quotient
by this action:
\begin{definition}\label{def:cM}
Given $\x$, $\y$, $B$, $\Source$, and $P$ or
  $\vec{P}$ as above, let
\begin{align*}
\gls*{ModSpaceUnorderedUnparam}
&\coloneqq\tcM^B(\x,\y;\Source;P)/\RR\\
\gls*{ModSpaceOrderedUnparam}
&\coloneqq\tcM^B(\x,\y;\Source;\vec{P})/\RR.
\end{align*}
denote the reduced moduli spaces.
\index{moduli space!reduced}\index{reduced!moduli space}%
\end{definition}
The expected dimension of $\cM^B(\x,\y\semico\Source\semico P)$ (except in the 
trivial case $B=0$) is given by $\ind(B,\Source,P)-1$.
\index{expected dimension!of $\cM^B(\x,\y\semico\Source\semico P)$}%

Although the evaluation maps $\gls*{evq}$ (for $q \in E$) do not descend to
the quotient by translation, the difference between any two of them,
$\ev_p - \ev_q$, does descend.  We denote this difference by
\label{page:evpq}
\begin{align}
  \gls*{evpq}
  &\co \Mod^B(\x,\y;\Source;\vec{P}) \to \RR.
  \label{eq:DefEvPQ} \\
\intertext{We can also combine all the evaluation maps into a single map}
  \gls*{ev}&\co \cM^B(\x,\y;\Source;\vec{P})\to\gls*{REmodR}
  \nonumber
\end{align}
where $\RR$ acts diagonally by translation on $\RR^E$.

\section{Holomorphic curves in
  \textalt{$\RR\times Z\times[0,1]\times\RR$}
              {R \texttimes\ Z \texttimes\ [0,1]\ \texttimes\ R}}
\label{sec:curves_at_east_infinity}
In order to discuss the compactifications of the moduli spaces $\cM^B(\x,\y\semico\Source\semico P)$
and $\cM^B(\x,\y\semico\Source\semico\vec{P})$, we will need to consider certain holomorphic curves
``at east $\infty$,'' i.e., in $\RR\times Z \times[0,1]\times\RR$.
\index{east infinity!holomorphic curves at}%
We
endow $\RR\times Z\times[0,1]\times\RR$ with the obvious split
symplectic form. Inside this space, we have $4k$ Lagrangian planes
$\RR\times\mathbf{a}\times\{1\}\times\RR$.  The space
$\RR\times Z\times[0,1]\times\RR$ has four different ends: the
ends $\pm\infty$ in the first $\RR$-factor, which we call east and
west infinity
\index{east infinity!of $\RR\times Z\times[0,1]\times\RR$}%
\index{west infinity!of $\RR\times Z\times[0,1]\times\RR$}%
\index{infinity!of $\RR\times Z\times[0,1]\times\RR$}%
\index{infinity!east}%
\index{infinity!west}%
respectively, and the ends $\pm\infty$ in the second $\RR$-factor, which we call
$\pm\infty$. Note that there is an $(\RR\times\RR)$-action
\index{$\RR\times\RR$-action}%
on $\RR\times Z\times[0,1]\times\RR$, by translation in the two
$\RR$-factors.
There are projection maps $\pi_\Sigma$, $\pi_\DD$, $s$, and~$t$,
just as before.

Fix a split complex structure $J = j_\Sigma\times j_\DD$ on
$\RR\times Z\times[0,1]\times\RR$.
\index{complex structure!on $\RR\times Z\times[0,1]\times\RR$}%

\begin{definition}
By a \emph{bi-decorated source $\gls*{biDecSource}$} we mean:
\index{source!bi-decorated}%
\begin{enumerate}
\item A topological type of smooth (not nodal) Riemann surface $\gls*{biSourceSurf}$ with boundary and
  punctures on the boundary.  
\item A labeling of each puncture of $T$ by either $e$ or $w$ (east or
  west), and
  \index{east puncture|see{puncture, east}}%
  \index{west puncture|see{puncture, west}}%
  \index{puncture!east}%
  \index{puncture!west}%
\item A labeling of each puncture of $T$ by a Reeb chord~$\rho$ in
  $(Z,\mathbf{a})$.
\end{enumerate}
\end{definition}

Given a bi-decorated source $\biSource$, we will consider maps 
\[
\gls*{v}\co (T,\bdy T)\to (\RR\times (Z\setminus z)\times[0,1]\times\RR, \RR\times\mathbf{a}\times\{1\}\times\RR)
\]
satisfying the following conditions:
\begin{enumerate}
\item \label{item:e_map_first}The map $v$ is $(j,J)$-holomorphic with respect to some
  almost complex structure $j$ on $T$.
\item The map $v$ is proper.
\item At each west puncture $q$ of $T$ labeled by a Reeb chord~$\rho$,
  $\lim_{z\to q}\pi_\Sigma\circ u(z)$ is $\rho\subset \{-\infty\}\times Z$.
\item \label{item:e_map_last}At each east puncture $q$ of $T$ labeled by
  a Reeb chord~$\rho$,
  $\lim_{z\to q}\pi_\Sigma\circ u(z)$ is $\rho\subset \{+\infty\}\times Z$.
\end{enumerate}
Note that each component of such a holomorphic curve necessarily maps
to a single point under~$\pi_\DD$ by the maximum principle, since its
boundary maps entirely to $s=1$.
\index{east infinity!holomorphic curves at}%

Again, we collect these holomorphic curves into moduli spaces:
\begin{definition}For $\biSource$ a bi-decorated source
let 
$\gls*{ModSpaceOpenBi}$
denote the moduli space of holomorphic maps from $T$ satisfying
properties~(\ref{item:e_map_first})--(\ref{item:e_map_last}) above.
\end{definition}
\index{moduli space!of curves at east $\infty$}%

A holomorphic map $v$ is called \emph{stable}\index{stable holomorphic curve!at east $\infty$} if:
\begin{itemize}
\item The source of every component of $v$ on which $\pi_\Sigma\circ v$ is constant has no non-trivial automorphisms and
\item There is at least one component of $v$ which is not a twice-punctured disk.
\end{itemize}

We shall be interested in $\tcN(\biSource)$ only in the cases that the
corresponding maps are stable.

Suppose that $q$ is a puncture of $\biSource$. Then there is an evaluation map
$\gls*{evq}
\co \tcN(\biSource)\to\RR$ given by $\ev_{q}(v)=\lim_{z\to q}t\circ
v(z)$. (Here, writing $\lim$ is somewhat silly, since $t\circ v$ is constant on
each connected component of $\biSource$.) There are, consequently, evaluation maps
\begin{align*}
\gls*{evwest}
&\coloneqq\,\,\prod_{\mathclap{q \in W(\biSource)}}\,\,\ev_{q}\co \tcN(\biSource)\to\RR^{W(\biSource)}\quad\textrm{and}\\
\gls*{eveast}&\coloneqq\,\,\prod_{\mathclap{q \in E(\biSource)}}\,\,\ev_{q}\co \tcN(\biSource)\to\RR^{E(\biSource)}
\end{align*}
where $\gls*{WofT}$ or $\gls*{Wofbi}$ is the set of west punctures
of~$\biSource$ and $\gls*{EofT}$ or $\gls*{Eofbi}$ is the
set of east punctures.

We define sub-moduli spaces of $\tcN(\biSource)$ by cutting-down by partial
diagonals.
\begin{definition} Given a bi-decorated source $\biSource$ and partitions $\gls*{PartitionWest}$ and
  $\gls*{PartitionEast}$ of the west and east punctures of $\biSource$, respectively, let
\[
\gls*{ModSpaceBiUnorderedParam}
=(\ev_w\times\ev_e)^{-1}(\Delta_{P_w}\times\Delta_{P_e}).
\]
If $P_w$ is the discrete partition, we denote
$\tcN(\biSource;P_w,P_e)$ by $\gls*{ModSpaceBiUnorderedParamRDisc}$.
\end{definition}

The moduli spaces $\tcN(\biSource)$ can be understood
concretely. Indeed, they are determined by the $t$-coordinates of the
components of $T$, and a holomorphic map from $\biSource$ to $\RR\times Z$,
which gives $\pi_\Sigma\circ v$ (which in turn is determined by local branching
data over $\Sigma$). However, in general the
moduli spaces $\tcN(\biSource)$ are not transversally cut
out. One situation in which they are---and indeed this is the case which
is relevant for our applications---is the following:
\begin{proposition}\label{prop:east_transversality}Suppose $\biSource$ is a decorated
  source such that all components of $T$ are topological disks. Then $\tcN(\biSource)$
  is transversally cut out by the $\widebar{\partial}$ equation for any split complex
  structure on $\RR\times Z\times[0,1]\times\RR$.
\end{proposition}
\begin{proof}
  This is proved by explicit computation of the cokernel of the
  $\dbar$\hyp operator, which is identified with a certain sheaf
  cohomology group. See McDuff and Salamon~\cite[Section
  3.3]{MS04:HolomorphicCurvesSymplecticTopology} for a nice exposition
  of these ideas in the absolute case. (The relative case can be
  deduced from the absolute case via a doubling argument; see, for
  instance, Hofer, Lizan, and Sikorav~\cite[Section
  4]{HLS97:GenericityHoloCurves}.)  More details can also be found
  in~\cite[Lemma 4.1.2]{Lipshitz06:BorderedHF}.
\end{proof}

The $(\RR\times\RR)$-translation action on
\index{$\RR\times\RR$-action}%
$\RR\times Z\times[0,1]\times\RR$ induces an
$(\RR\times\RR)$-action
on $\tcN(\biSource)$ for any $\biSource$; if $\tcN(\biSource)$ is stable then (almost
by definition) this $(\RR\times\RR)$-action is free. As for the
earlier moduli spaces, we
are primarily interested in the
quotient space:
\begin{definition} Given a stable decorated source $\biSource$, let
  \[\gls*{ModSpaceBiOpenUnparam} =\tcN(\biSource)/(\RR\times\RR).\]
\end{definition}

Certain holomorphic curves in $\RR\times Z\times[0,1]\times\RR$ will play a
special role below. By a \emph{trivial component}
\index{trivial!component}%
we mean a component of $\biSource$
which is a topological disk with exactly two boundary punctures, one east and
one west, and both labeled by the same Reeb chord. Holomorphic maps from such
components are uninteresting. In particular, they are preserved by the
$\RR$-translation action on $\RR\times Z$.

A more interesting kind of component is a \emph{join component}.
\index{join!component|see{component, join}}%
\index{component!join}%
This
is a component of $\biSource$ which is a topological disk with two
west punctures and one east puncture.  In counterclockwise order,
suppose the punctures are labeled by $(e,\rho_e)$, $(w,\rho_1)$ and
$(w,\rho_2)$.  (We will sometimes refer to the $(w,\rho_1)$ puncture
as the \emph{top puncture}
\index{top puncture}%
\index{puncture!top}%
\index{bottom puncture}%
\index{puncture!bottom}%
of the join component, and the $(w,\rho_2)$
puncture as the \emph{bottom puncture}.) Then there is a holomorphic
map from this component if and only if
$\rho_e=\rho_2\uplus\rho_1$.  
If such a holomorphic map exists, it is unique up to translation.  A
curve consisting entirely of one join component and some number of
trivial components is called a \emph{join curve}.
\index{join!curve|see{curve, join}}%
\index{curve!join}%

Symmetrically, a \emph{split component} is a component of $\biSource$
which is a topological disk with two punctures labeled $e$ and one
labeled $w$. 
\index{component!split}%
\index{split!component|see{component, split}}%
In counterclockwise order,
suppose the punctures are labeled by $(w,\rho_w)$, $(e,\rho_1)$ and
$(e,\rho_2)$. (We will sometimes refer to the $(e,\rho_2)$ puncture as
the \emph{top puncture} of the split component, and the $(e,\rho_1)$
puncture as the \emph{bottom puncture}.) Then there is a holomorphic
\index{top puncture}%
\index{puncture!top}%
\index{bottom puncture}%
\index{puncture!bottom}%
map from such this component if and only if
$\rho_w=\rho_1\uplus\rho_2$.  
Again, if such a holomorphic map exists, it is unique up to
translation. A stable curve consisting entirely of some number of
split components and some number of trivial components is called a
\index{curve!split}%
\index{split!curve|see{curve, split}}%
\emph{split curve}.  (Note the definition is not symmetric with that
of join curves.)

\begin{figure}
\includegraphics[scale=.83333]{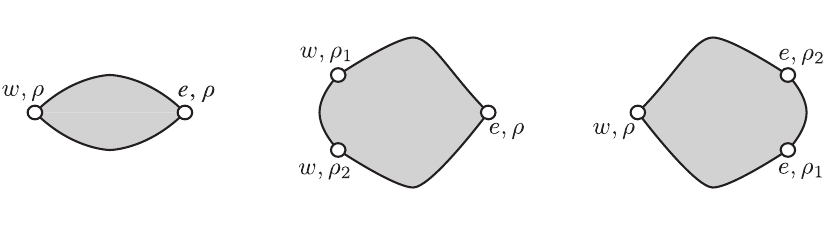}
\caption[Sources of curves at east $\infty$]{\textbf{Sources of curves at east $\infty$.} Left: a trivial component. Center: a join component. Right: a
  split component.}\label{fig:trivial_split_join_domain}
\end{figure}

The sources
of these components are illustrated in
Figure~\ref{fig:trivial_split_join_domain}.  Split and join components
can arise from the degenerations in Figure~\ref{fig:degen_examples}.

\begin{figure}
\includegraphics[scale=.83333]{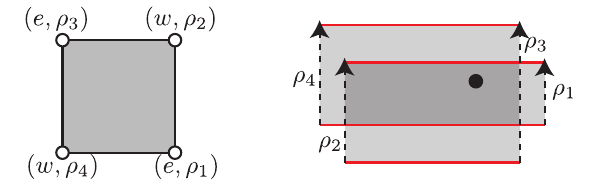}
\caption[Visualizing shuffle curves, I]{\textbf{Visualizing shuffle curves.} Left: a shuffle-component. Right: a schematic of the image of
  an odd shuffle-component projected to
  $\bdy\overline{\Sigma}\times\RR$.}\label{fig:shuffle1}
\end{figure}
\colorused
\begin{figure}
\includegraphics[scale=.83333]{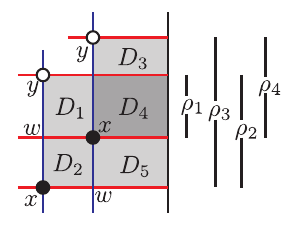}
\caption[Visualizing shuffle curves, II]{\textbf{A domain in a Heegaard diagram where a shuffle component
  degenerates off.} The domain is
  $B=D_1+D_2+D_3+2D_4+D_5\in\pi_2(\x,\y)$, as indicated by the
  shading. There is a $1$-parameter family of holomorphic curves in
  $\cM^B(\x,\y;\Source;P)$ where $P$ is the partition
  $\{\{\rho_1,\rho_3\}\}$. One end of this moduli space consists of a
  height $2$ holomorphic comb connecting $\x$ to $\w$ via $D_2$ and
  $\w$ to $\y$ via $D_1+D_3+2D_4+D_5$. At the other end, either an
  interior branch point of $\pi_\Sigma\circ u$ (in region $D_4$) or
  two boundary branch points of $\pi_\Sigma\circ u$ (between $D_3$ and
  $D_4$ or between $D_4$ and $D_5$) approach $\bdy \overline{\Sigma}$,
  resulting in the degeneration of a shuffle
  curve.}\label{fig:shuffle2}
\end{figure}

Finally, there is one more complicated kind of curve which we must
consider, a \emph{shuffle curve}. A \emph{shuffle component}
\index{shuffle!curve|see{curve, shuffle}}%
\index{shuffle!component|see{component, shuffle}}%
\index{component!shuffle}%
is a
topological disk with two $e$ punctures and two $w$ punctures, with
cyclic ordering $(e,w,e,w)$ around the boundary. (A shuffle component
is pictured in Figure~\ref{fig:shuffle1}, and an example where one
degenerates off is shown in Figure~\ref{fig:shuffle2}.) The map from a shuffle
component to $\RR\times Z$ has either one interior branch point or two
boundary branch points. Suppose that in counterclockwise order, the
punctures of a shuffle component are labeled
$(e,\rho_1),(w,\rho_2),(e,\rho_3),(w,\rho_4)$. Then the moduli space
is non-empty if and only if $\rho_1^+=\rho_2^+$, $\rho_2^-=\rho_3^-$,
$\rho_3^+=\rho_4^+$ and $\rho_4^-=\rho_1^-$.  
If a join curve admits a representative, exactly one of
the following cases occurs:
\begin{itemize}
\item $\{\rho_1,\rho_3\}$ is nested and
  $\{\rho_2, \rho_4\}$ is interleaved, or
\item $\{\rho_2, \rho_4\}$ is nested  and
  $\{\rho_1,\rho_3\}$ is interleaved.
\end{itemize}
(Recall the definition of interleaved and nested pairs of Reeb chords
from Definition~\ref{def:interleaved-nested}.)
In fact, we will see in Proposition~\ref{prop:gluing-shuffle} that the
two cases are \emph{not} symmetric for our purposes: when the first
occurs in a degeneration, it occurs an odd number of times. When the
second occurs, it occurs an even number of times. We will call the
first case an \emph{odd shuffle component} and the second case an
\emph{even shuffle component}, to remember the distinction. 
\index{component!shuffle!even}%
\index{component!shuffle!odd}%
\index{odd!shuffle component|see{component, shuffle, odd}}%
\index{even!shuffle component|see{component, shuffle, even}}%

A holomorphic curve in $\RR\times Z\times[0,1]\times\RR$ consisting of one shuffle component and some number of trivial components is called a \emph{shuffle curve}; a shuffle curve can be odd or even, depending on the type of its shuffle component.
\index{curve!shuffle}%
\index{curve!shuffle!even}%
\index{curve!shuffle!odd}%
\index{odd!shuffle curve|see{curve, shuffle, odd}}%
\index{even!shuffle curve|see{curve, shuffle, even}}%

\section{Compactifications via holomorphic combs}\label{sec:combs-compact}
Just as the compactifications of moduli spaces in Morse theory involve not just flow
lines but broken flow lines, in order to discuss the compactifications
of the moduli
spaces $\cM^B(\x,\y;\Source;P)$ we need to introduce a more general object than
holomorphic curves.  The idea is that, because of the non-compactness of the target
space $\Sigma\times[0,1]\times\RR$, holomorphic curves can degenerate into ``holomorphic
buildings'' with several stories, as described by Eliashberg, Givental,
and Hofer~\cite{EGH00:IntroductionSFT}. In our
setting, however, there are two different kinds of infinities: $\pm\infty$ and
$e\infty$. This leads us to the notion of \emph{holomorphic combs}.

By a \emph{simple holomorphic comb}
\index{holomorphic comb!simple}%
\index{simple!holomorphic comb|see{holomorphic comb, simple}}%
we mean a pair of holomorphic maps $(u,v)$ where
$u$ maps to $\Sigma\times[0,1]\times\RR$ and $v$ maps to
$\RR\times Z\times[0,1]\times\RR$, and such that the asymptotics
of $u$ at east $\infty$ match with the asymptotics of $v$ at west $\infty$. More
precisely:
\begin{definition}
  \label{def:SimpleComb}
  A \emph{simple holomorphic comb} is a pair $\gls*{uv}$ with
  $u\in\cM^B(\x,\penalty500\y\semico\Source)$ and $v\in\cN(\biSource)$ for some sources
  $\Source$ and $\biSource$, together with a one-to-one
  correspondence between $E(\Source)$ and
  $W(\biSource)$, preserving the labeling by Reeb chords, and such that
  $\ev(u)=\ev_w(v)$ inside $\RR^{E(\Source)}/\RR \cong
  \RR^{W(\biSource)}/\RR$.
\end{definition}

More generally, we also want to allow the components at east $\infty$
to degenerate further:
\begin{definition}
  A \emph{holomorphic story}
\index{holomorphic story}%
\index{story|see{holomorphic story}}%
  is a sequence
  $(u,v_1,\dots,v_k)$ (for
  some $k\ge 0$) where, for some $B$, $\Source$, and $\biSource_i$,
  \begin{itemize}
  \item $u\in\cM^B(\x,\y;\Source)$,
  \item $v_i\in\cN(\biSource_i)$,
  \item $(u,v_1)$ is a simple holomorphic comb (if $k \ge 1$),
  \item there is a correspondence between $E(\biSource_i)$ and
    $W(\biSource_{i+1})$ for $i=1,\dots,k-1$ which preserves the
    labelings by Reeb chords, and
  \item $\ev_e(v_i)=\ev_w(v_{i+1})$ in $\RR^{E(\biSource_i)}/\RR \cong
    \RR^{W(\biSource_{i+1})}/\RR$.
  \end{itemize}
  We refer to $(u,v_1,\dots,v_k)$ as the \emph{horizontal levels} of
  the holomorphic story.
  \index{holomorphic comb!horizontal level of}%
  \index{holomorphic story!horizontal level of}%
  \index{horizontal!level of holomorphic story}
\end{definition}

\begin{figure}
\includegraphics[scale=.83333]{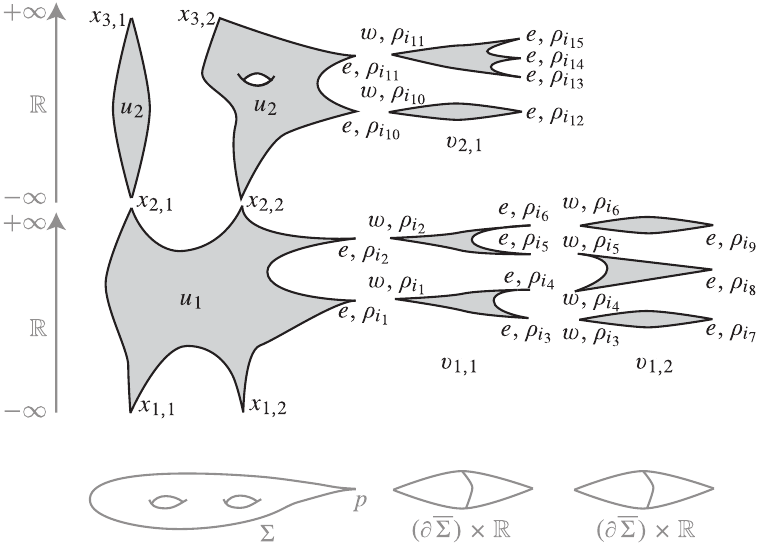}
\caption[Schematic of a height $2$ holomorphic comb]{\textbf{A schematic of a height $2$ holomorphic comb.}}
\label{fig:comb_schematic}
\end{figure}

We also need to allow degeneration at $\pm\infty$: 
\begin{definition}\label{def:comb}
  A \emph{holomorphic comb}
  \index{holomorphic comb}%
  \index{comb|see{holomorphic comb}}%
  \index{holomorphic comb!height}%
  of height~$N$ is a sequence
  $(u_j,v_{j,1},\allowbreak\dots,\allowbreak v_{j,k_j})$, for
  $j=1,\dots,N$, of holomorphic
  stories with $u_j$ a stable curve in
  $\cM^{B_j}(\x_j,\x_{j+1};\allowbreak\SourceSub{j})$ for some sequence of
  generalized generators $\x_1,\dots,\x_{N+1}$.  We admit the case
  $N=0$, the \emph{trivial holomorphic comb} from $\x_1$ to $\x_1$,
  which corresponds to a trivial (unstable) holomorphic curve.
  \index{holomorphic comb!trivial}\index{holomorphic comb!height!zero}%
  \index{trivial!holomorphic comb|see{holomorphic comb, trivial}}%
  \index{moduli space!of combs|see{holomorphic comb}}%
  We refer to the index $j$ as the \emph{vertical level}, so the comb
  $(u_j,v_{j,1},\dots,v_{j,k_j})$, $j=1,\dots,N$ has $N$ vertical
  levels (or \emph{stories}).
  \index{holomorphic comb!vertical level of}%
  \index{holomorphic comb!story of}%
  \index{vertical!level of comb}
\end{definition}
See Figure~\ref{fig:comb_schematic} for a schematic illustration of a
height 2 holomorphic comb.

We will call a holomorphic comb \emph{toothless} if it has no components at $e\infty$. The \emph{spine} of a holomorphic comb is the sub-comb of components mapped to $\Sigma\times[0,1]\times\RR$.
\index{spine! of comb}\index{toothless comb}\index{holomorphic comb!spine of}\index{holomorphic comb!toothless}%

Finally, to compactify our moduli spaces, we will also need to drop
Condition~(\ref{item:smooth_source}) of
Definition~\ref{def:decorated-source}, and allow our holomorphic
curves to have nodal sources. If we want to explicitly allow nodes, we
will refer to such combs as \emph{nodal holomorphic combs}; 
\index{nodal!holomorphic comb|see{holomorphic comb, nodal}}%
\index{holomorphic comb!nodal}%
if we want to
explicitly rule out nodal combs, we will refer to \emph{smooth combs}. When
the distinction is unimportant, we will simply talk about holomorphic
combs.
\index{smooth holomorphic comb|see{holomorphic comb, smooth}}%
\index{holomorphic comb!smooth}%

Note that a (nodal) holomorphic comb $\gls*{capU}=\{(u_j,v_{j,1},\dots,v_{j,k_j})\}_{j=1}^N$
naturally represents a homology class $B$ in $\pi_2(\x_1,\x_{N+1})$,
where $\x_1$ and $\x_{N+1}$ are the first and last generators
involved. 
\index{holomorphic comb!homology class of}%
\index{homology class!of comb}%
(Specifically,
$B= B_1*\cdots*B_N$, where $B_j$ is the domain of $u_j$. The
$v_{i,j}$ are irrelevant
to this homology class.)  Further, $U$ is asymptotic to a well-defined set of Reeb
chords at (far) east $\infty$: these are the asymptotics at east infinity
of the $v_{j,k_j}$.

Next, we turn to what it means for a sequence of holomorphic curves to
converge to a holomorphic comb or, more generally, for a sequence of holomorphic
combs to converge to another holomorphic comb. This is a simple adaptation of
definitions for holomorphic curves converging to holomorphic buildings
in~\cite{EGH00:IntroductionSFT,BEHWZ03:CompactnessInSFT}. 
\index{holomorphic building}
First, adapting the definition of convergence
in~\cite[Section 7.3]{BEHWZ03:CompactnessInSFT}, say, to the relative
case of holomorphic buildings in a symplectic manifold $W$ with
cylindrical ends with boundary on a Lagrangian $L\subset W$ with
cylindrical ends involves only trivial notational changes; we leave this to the reader. 
Next, call a component $C$ of $U$
\emph{$\Sigma$-stable} (respectively \emph{$\DD$-stable}) if
$(\pi_\Sigma\circ U)|_C$ (respectively $(\pi_\DD\circ U)|_C$) is a stable map (i.e., has a finite
automorphism group), and \emph{$\Sigma$-unstable} (respectively
\emph{$\DD$-unstable}) otherwise. (Here, we view trivial strips as
stable maps to $\DD$, i.e., we do not quotient by automorphisms of the
target, so every component of $U$ must be either $\Sigma$-stable or
$\DD$-stable or both.) With this terminology in hand we
are ready to define convergence to a holomorphic comb:
\begin{definition}\label{def:converge-to-comb}\index{convergence!to holomorphic comb}
  We say a sequence $u_n$ of holomorphic curves in
  $\Sigma\times[0,1]\times\RR$ \emph{converges} to a holomorphic comb $U$ if
  the following conditions hold:
  \begin{enumerate}[label=(C-\arabic*),ref=C-\arabic*]
  \item \index{(C-1)--(C-5)} Let $S_\Sigma$ be the result of
    collapsing all $\Sigma$-unstable components of $U$. Then
    $\{\pi_\Sigma\circ u_n\}$ converges to $\pi_\Sigma\circ U|_{S_\Sigma}$ as
    buildings (that is, in the sense of \cite[Section
    7.3]{BEHWZ03:CompactnessInSFT}). (Note that, for a non-split almost complex structure, the maps $\pi_\Sigma\circ u_n$ and $\pi_\Sigma\circ U$ are not holomorphic, but the topology defined in~\cite{BEHWZ03:CompactnessInSFT} makes sense for smooth buildings, not just holomorphic ones.) 
  \item Let $S_\DD$ be the result of collapsing all $\DD$-unstable
    components of $U$. Then $\{\pi_\DD\circ u_n\}$ converges to
    $\pi_\DD \circ U|_{S_\DD}$ as holomorphic buildings (i.e., in the sense
    of~\cite[Section 7.3]{BEHWZ03:CompactnessInSFT}).
  \item\label{item:c-infty-conv} Let $\tau_t\co
    \Sigma\times[0,1]\times\RR\to \Sigma\times[0,1]\times\RR$ denote
    translation by $t$ units in the $\RR$ direction. Then for each smooth
    point $q$ in the spine of $U$ there is a neighborhood $V\ni q$, a
    sequence of points $q_n$ in the source of $u_n$, neighborhoods
    $V_n\ni q_n$, diffeomorphisms $V_n\cong V$, and numbers $t_n\in
    \RR$ so that $\tau_{t_n}\circ u|_{V_n}$ converges in the
    $C^\infty_{\mathrm{loc}}$ topology to $u|_V$.
  \item All $u_n$ for $n$ sufficiently large represent the same
    homology class $B$, and $U$ represents $B$ as well.
  \end{enumerate}
  Convergence of a sequence of holomorphic combs to a holomorphic comb
  is defined similarly; we leave this extension to the reader.
\end{definition}


Given a simple holomorphic comb $(u,v)$ with $\Source$ and $\biSource$
the decorated sources of $u$ and $v$ respectively, there is a natural
way to (pre)glue $\Source$ and $\biSource$ to form a decorated source
$\Source\mathbin{\gls*{preglue}} \biSource$: the surface $S\glue T$ is obtained by
identifying small neighborhoods of the east punctures of $S$ with
neighborhoods of the corresponding west punctures of~$T$.  (Note that
this is a purely topological operation, not involving any differential
equations.  We use the term ``preglue'' to distinguish from gluing of
\index{pregluing sources}\index{gluing!holomorphic curves}%
holomorphic curves.)  The $e\infty$ punctures of $\biSource$ and
$\pm\infty$ punctures of $\Source$ carry over to decorate
$\Source\glue\biSource$. The pregluing of sources extends in an
obvious way to general smooth holomorphic combs: for any smooth
holomorphic comb there is a corresponding (preglued) smooth
surface. The construction also
extends to nodal holomorphic combs: there is a canonical way to deform
away the nodes in the source to produce a smooth surface.

From these constructions we can produce a compactification
of $\cM^B(\x,\y;\Source)$.
\index{moduli space!compactified}%
\index{compactified moduli space|see{moduli space, compactified}}%
\begin{definition}\label{def:ocM} We make the following definitions.
\begin{itemize}
\item 
  $\gls*{ModSpaceCpctCpct}$
  is the space of all
  (possibly nodal) holomorphic combs whose preglued surfaces are
  $\Source$, in the homology class $B$, with asymptotics $\x$ at
  $-\infty$, $\y$ at $+\infty$.
\item $\gls*{ModSpaceCpct}$
  is the closure of $\cM^B(\x,\y;\Source)$ in
  $\oocM^B(\x,\y;\Source)$.
\item For $p, q \in E(\Source)$,
  $\gls*{evpqcpct}\co\oocM^B(\x,\y,\Source)\to [-\infty,\infty]$
  is the extension of the map~$\ev_{p,q}$ (see Equation~\eqref{eq:DefEvPQ}).
\item $\gls*{ModSpaceUnorderedCpctCpct}$
  is the space of all
  holomorphic combs respecting the partition~$P$.  Formally, we set
  \[
  \oocM^B(\x,\y;\Source;P) \coloneqq
   \bigcap_{\substack{P_i \in P\\p,q\in P_i}}
   \widebar{\ev}_{p,q}^{-1}(0).
  \]
\item $\gls*{ModSpaceUnorderedCpct}$
  is the closure of
$\cM^B(\x,\y;\Source;P)$ in $\oocM^B(\x,\y;\Source;P)$.
\item $\gls*{ModSpaceOrderedCpct}$
  is the closure of
  $\cM^B(\x,\y;\Source;\vec{P})$ in $\ocM^B(\x,\y;\Source;P)$.
\end{itemize}
\end{definition}

\begin{figure}
\includegraphics[scale=.83333]{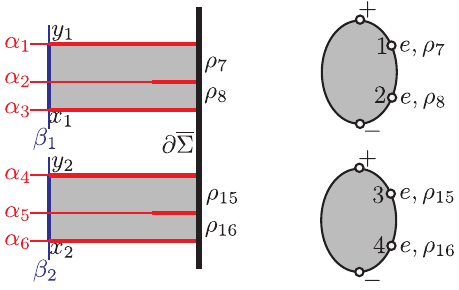}
\caption[$\ocM^B(\x,\y;\Source;P)$ may be a proper subset
  of $\oocM^B(\x,\y;\Source;P)$, I]{\textbf{An example illustrating that $\ocM^B(\x,\y;\Source;P)$ may be a proper subset
  of $\oocM^B(\x,\y;\Source;P)$. Left: a portion of a Heegaard diagram.} The
  domain of interest is shaded in gray. Right: the decorated source $\Source$.}\label{fig:moduli_proper_subset}
\end{figure}

\begin{figure}
\includegraphics[scale=.83333]{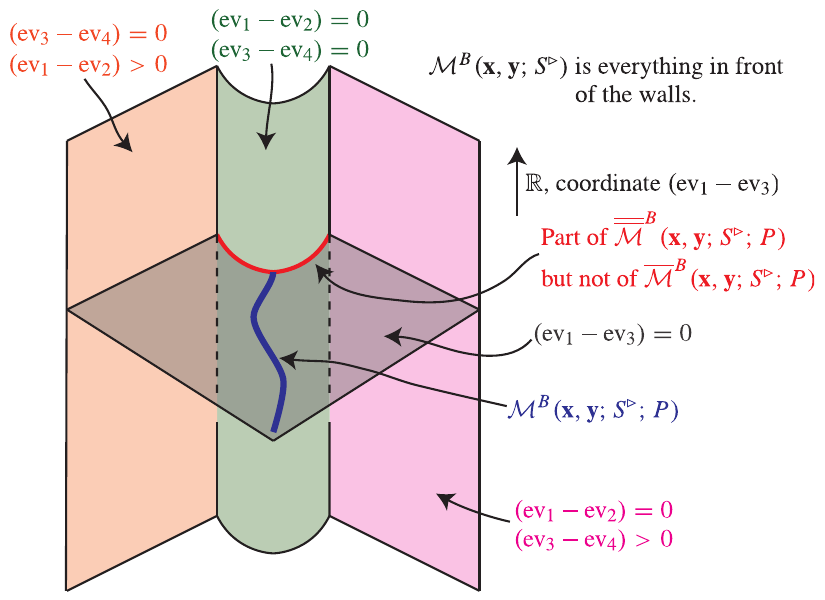}
\caption[$\ocM^B(\x,\y;\Source;P)$ may be a proper subset
  of $\oocM^B(\x,\y;\Source;P)$, II]{\textbf{An illustration of various moduli spaces for the Heegaard
  diagram from Figure~\ref{fig:moduli_proper_subset}.} The colored
  ``walls'' form the boundary of $\ocM^B(\x,\y;\Source)$, the open
  part of which one visualizes as the empty space in front of the
  walls. The subspace where $\ev_1=\ev_3$ is the gray horizontal
  plane. The moduli space $\ocM^B(\x,\y;\Source;P)$ is a thick,
  squiggly line segment; the space $\oocM^B(\x,\y;\Source;P)$ is a
  ``T'' formed by that line segment and an arc in
  $\bdy\ocM^B(\x,\y;\Source)$.}
\label{fig:moduli_proper_subset_2}
\end{figure}

The space $\ocM^B(\x,\y;\Source;P)$ is often a proper subset of
$\oocM^B(\x,\y;\Source;P)$, as the next example illustrates. (Because of
difficulties with transversality at east $\infty$,
$\ocM^B(\x,\y;\Source)$ may also be a proper subset of
$\oocM^B(\x,\y;\Source)$.)
\begin{example}Consider the portion of a Heegaard diagram shown in the left
  of Figure~\ref{fig:moduli_proper_subset}, and let $\Source$ be the decorated source shown
  on the right of Figure~\ref{fig:moduli_proper_subset}.  (We number
  the east punctures of $\Source$ for convenience in referring to them.)
Let $P$ denote the partition
  $\{\{1,3\},\{2,4\}\}$. The space $\widetilde{\cM}^B(\x,\y;\Source)$ is
    then homeomorphic to $\RR\times(0,\infty)\times\RR\times(0,\infty)$, via the map
    \[
    (\ev_1,\ev_1-\ev_2,\ev_3,\ev_3-\ev_4).
    \]
    The space $\widetilde{\cM}^B(\x,\y;\Source;P)$ is the subspace $x_1=x_3$, $x_2=x_4$
    of $\RR\times(0,\infty)\times\RR\times(0,\infty)$.
    
    Now, consider the compactification. In $\ocM^B(\x,\y\semico\Source)$, when
    $(\ev_1-\ev_2)\to 0$ the curve degenerates a split component
    at east $\infty$; similarly when $(\ev_3-\ev_4)\to 0$. If $(\ev_1-\ev_2)\to 0$ and
    $(\ev_3-\ev_4)\to 0$, there are two split components at east~$\infty$. There
    is a one-parameter family of these curves at east $\infty$,
    given by the relative $\RR$-coordinate of the branch point under
    $\pi_\Sigma$. Consequently the part of the space
    $\ocM^B(\x,\y;\Source)$ with $(\ev_1-\ev_2)$ and $(\ev_3-\ev_4)$ small has the form
    shown in Figure~\ref{fig:moduli_proper_subset_2}, and
    $\ocM^B(\x,\y;\Source;P)$ is the subspace pictured. By contrast,
    $\oocM^B(\x,\y;\Source;P)$ contains the entire stratum in
    which two split components have degenerated.
\end{example}
\colorused

Next, we turn to the technical results justifying the definitions of
the moduli spaces $\ocM^B(\x,\y\semico\Source\semico P)$.
\index{moduli space!compactified!is compact}%
\begin{proposition}\label{prop:compactness} The spaces
  $\ocM^B(\x,\y;\Source)$ are compact. That is, suppose that
  $\{U_n\}$ is a sequence of holomorphic combs in a fixed homology
  class, with a fixed preglued topological source. Then $\{U_n\}$ has
  a subsequence which converges to a (possibly nodal) holomorphic comb
  $U$, in the same homology class as the~$U_n$.  Similarly,
  $\ocM^B(\x,\y;\Source;P)$ and $\ocM^B(\x,\y;\Source;\vec{P})$ are compact.
\end{proposition}

We will deduce Proposition~\ref{prop:compactness} from two existing
compactness results. First, there is a local compactness result for
holomorphic curves:
\begin{citethm}\label{thm:MS-cpct}\cite[Theorem
  4.1.1]{MS04:HolomorphicCurvesSymplecticTopology} Let $(M,J)$ be a
  compact almost complex manifold, $L\subset M$ a compact totally real
  submanifold, and $J_\nu$ a sequence of almost complex structures on
  $M$ that converges in the $C^\infty$ topology to $J$. Let $(\Sigma,
  j_\Sigma, d\mathrm{vol}_\Sigma)$ be a (possibly non-compact) Riemann surface with boundary,
  $\Omega_\nu\subset \Sigma$ an increasing sequence of open sets that
  exhaust $\Sigma$, and $u_\nu\co
  (\Omega_\nu,\Omega_\nu\cap\bdy\Sigma)\to (M,L)$ a sequence of
  $J_\nu$-holomorphic curves such that
\[
\sup_\nu\|du_\nu\|_{L^\infty(K)}<\infty
\]
for every compact subset $K\subset \Sigma$. Then $u_\nu$ has a
subsequence which converges uniformly with all derivatives on compact
subsets of $\Sigma$ to a $J$-holomorphic curve $u\co (\Sigma,\bdy
\Sigma)\to (M,L)$.%
\end{citethm}

Recall from~\cite[Section~2]{BEHWZ03:CompactnessInSFT}, that if $V$ is
an odd-dimensional manifold, then an almost-complex structure $J$ on
$\RR\times V$ is called {\em cylindrical} if it is invariant under
translation in the $\RR$-factor and $\mathbf{R}\coloneqq
J(\frac{\partial}{\partial t})$ is tangent to $V$. In this case there is a one-form
$\lambda$ on $V$ which is characterized by the properties that
$\lambda|_{J TV\cap TV}\equiv 0$ and $\lambda(\mathbf{R})\equiv
1$. The cylindrical structure $J$ is called {\em symmetric} if the Lie
derivative of $\lambda$ in the direction $\mathbf{R}$ vanishes
identically. The almost complex structure $J$ is \emph{adjusted} to a
closed, maximal-rank $2$-form $\omega$ on $V$ if $\omega$ is
$J$-invariant; $\omega(v,Jv)$ is positive for all $v\in JTV \cap TV$;
and the Lie derivative of $\omega$ in the direction $\mathbf{R}$ is
zero. Given a map $u\co S\to \RR\times V$, the
\emph{$\omega$-energy}\index{energy!$\omega$}\index{$\omega$ energy}
of $u$ is $\int_S (\pi_V\circ u)^*\omega$; the
\emph{$\lambda$-energy}\index{energy!$\lambda$}\index{$\lambda$ energy}
of $u$ is $\sup_\phi \int_S [(\phi\circ t\circ u)dt]\wedge (\pi_V\circ
u)^*\lambda$, where the supremum is over all functions $\phi\co
\RR\to\RR_{\geq 0}$ with compact support and integral $1$; and the
\emph{energy}\index{energy} $E(u)$ of $u$ is the $\lambda$-energy plus the
$\omega$-energy~\cite[Section~5.3]{BEHWZ03:CompactnessInSFT}. 

The following degenerate case will play a role in the proof of Proposition~\ref{prop:compactness}:
\begin{example}\label{eg:1d-energy}
  If $V$ is $S^1$ then the $\omega$-energy (and, in fact, $\omega$)
  vanishes for dimension reasons. The $\lambda$-energy is given by the
  degree of the map $u\co S\to \RR\times S^1$. 
\end{example}

With these definitions in hand, we have the compactness theorem in
symplectic field theory:

\begin{citethm}\label{thm:BEHWZ-cyl-cpct}(\cite[Theorem
  10.1]{BEHWZ03:CompactnessInSFT}, see also~\cite{Abbas14:compactness})
  Let $(\RR\times V,J)$ be a symmetric cylindrical almost complex
  manifold. Suppose that the almost complex structure $J$ is adjusted
  to a $2$-form $\omega$. Then for every $E>0$, the
  space $\ocM_{g,\mu}(V)\cap \{E(F)\leq E\}$ is compact.
\end{citethm}
Here, $\ocM_{g,\mu}(V)$ denotes the moduli space of holomorphic
buildings in $\RR\times V$ with source a surface of genus $g$ with
$\mu$ marked points, and $E(F)$ denotes the energy of a building
$F$. (Punctures are viewed as marked points, but there may be
extra marked points as well.)

\begin{observation}\label{obs:BEHWZ-rel-cyl-cpct}
  Theorem~\ref{thm:BEHWZ-cyl-cpct} implies that any sequence of
  holomorphic maps
  \[
  u_n\co (S,\bdy S)\to ([0,1]\times\RR, \{0,1\}\times\RR)
  \]
  (with bounded energy) has a convergent subsequence: we can double
  the source $S$ across its boundary to obtain a sequence of
  holomorphic maps
  \[
  S\cup_\bdy (-S) \to S^1\times \RR
  \]
  and apply Theorem~\ref{thm:BEHWZ-cyl-cpct} to the result.
\end{observation}

Finally, we will use another relative version of symplectic field theory
compactness in the special case of maps between Riemann surfaces:
\begin{citethm}\label{thm:BEHWZ-2D-rel-cpct} Let $(W, j)$ be a
  punctured Riemann surface and $L\subset W$ a Lagrangian submanifold
  (union of curves) which is cylindrical near the punctures of
  $W$. Assume that $L$ is embedded away from finitely many transverse
  double points. Then for every $E>0$ the space
  $\ocM_{S,\mu}(W,L,J)\cap\{E(F)\leq E\}$ is compact, where
  $\ocM_{S,\mu}(W,L,J)\cap\{E(F)\leq E\}$ is the space of holomorphic buildings in
  $W$ with boundary on $L$, topological type $S$, $\mu$ marked
  points, and energy $\leq E$.
\end{citethm}
In this simple case, the energy of a map $F\co S\to \Sigma$ between Riemann surfaces is
the integral over $S$ of the pullback of a chosen area form on $\Sigma$.

See~\cite[Theorem 10.2]{BEHWZ03:CompactnessInSFT} for the analogous
result in the closed case, without the dimension restriction.  We will
not prove Theorem~\ref{thm:BEHWZ-2D-rel-cpct}, but there are at least
three ways one could go about doing so:
\begin{enumerate}
\item Imitate the arguments in~\cite{BEHWZ03:CompactnessInSFT} but in
  the relative case. This approach should yield
  Theorem~\ref{thm:BEHWZ-2D-rel-cpct}, where $W$ is an arbitrary
  symplectic manifold, not just a Riemann
  surface.
\item Apply Gromov's doubling trick~\cite[Paragraph 1.5.D${}_2$]{Gromov85} to deduce the relative case from
  the closed case. Again this would yield
  Theorem~\ref{thm:BEHWZ-cyl-cpct}, again without the dimension restriction,
  but with the restriction that the Lagrangians are real analytic,
  $J$ is integrable near the Lagrangians, and the
  Lagrangians are embedded. (This last restriction does not
  suit our purposes.)
\item Observe that a map between Riemann surfaces is determined by
  its branch points so, if we put marked points in the
  target (in addition to the punctures) to record the positions of
  branch points then the moduli space of maps is a quotient of a
  closed subspace of the Deligne-Mumford moduli space of 
  Riemann surfaces with boundary and marked points. (The quotient
  is because when branch points collide we remember less information
  than when marked points collide.)
  Theorem~\ref{thm:BEHWZ-2D-rel-cpct} then follows
  from the fact that the Deligne-Mumford moduli space of
  Riemann surfaces with boundary and marked points is compact.
\end{enumerate}
See also~\cite[Section 11.3]{BEHWZ03:CompactnessInSFT} for a little
more discussion of compactness in relative symplectic field theory
and~\cite{Abbas14:compactness} for further discussion, particularly
regarding the first approach. In the third point, compactness of the
moduli space of Riemann surfaces with boundary and marked points (and,
in fact, the description of the compactification itself) is induced
from compactness of the moduli space of closed Riemann surfaces with
marked points by doubling across the boundary (a simple case of the second strategy).

\begin{proof}[Proof of Proposition~\ref{prop:compactness}]
  In the case that the complex structure $J$ on
  $\Sigma\times[0,1]\times\RR$ is split, $J=j_\Sigma\times j_\DD$, the
  result follows by considering the projections $\pi_\DD\circ u_i$ and
  $\pi_\Sigma\circ u_i$ separately and applying
  Theorems~\ref{thm:BEHWZ-cyl-cpct} and~\ref{thm:BEHWZ-2D-rel-cpct},
  respectively. For general $J$, we use this projection argument in a
  neighborhood of east $\infty$ (using Property~(\ref{item:J4}) of
  $J$). Away from this neighborhood we use a combination of
  Theorem~\ref{thm:BEHWZ-cyl-cpct} applied to $\pi_\DD\circ u$ (which
  is holomorphic by Property~\ref{item:J2}) and
  Theorem~\ref{thm:MS-cpct}.  We spell this out in detail below; see
  also~\cite[Proposition 4.2.1]{Lipshitz06:BorderedHF} for a slightly
  different argument.

  To ease notation, we will restrict our attention to the case of a
  sequence $\{u_n\co S_n\to \Sigma\times[0,1]\times\RR\}_{n=1}^\infty$
  of holomorphic curves in $\cM^B(\x,\y;\Source)$, rather than a
  sequence of arbitrary holomorphic combs. It is routine to adapt the
  argument to the general case.

  \textbf{Step 1.} \emph{Extracting the vertical level structure on the
  limit.} Fix a point $p_r$ in each region $r$ of $\Sigma$. For
  generically chosen points $p_r$, the points $p_r$ are regular values
  of $\pi_\Sigma\circ u_n$ for all $n$.  Let
  $\{q_{r,i,n}\}=(\pi_\Sigma\circ u_n)^{-1}(p_r)$ be the preimages of
  $p_r$. We view the $q_{r,i,n}$ as extra marked points in $S_n$. Now,
  apply Theorem~\ref{thm:BEHWZ-cyl-cpct} (in light of
  Observation~\ref{obs:BEHWZ-rel-cyl-cpct}) to $\{\pi_\DD\circ u_n\}$
  (with these extra marked points) to extract a convergent
  subsequence. (Note that the energy bound required by
  Theorem~\ref{thm:BEHWZ-cyl-cpct} is trivially satisfied: the
  $\omega$-energy is $0$ (for dimension reasons) while the $\lambda$-energy is
  equal to $g$,
  the genus of $\Sigma$; compare Example~\ref{eg:1d-energy}.)  Relabel so that $\{u_n\}$ refers to this
  subsequence, and let $(\pi_\DD\circ U)$ denote the limit of this
  sequence, which is a (vertical) holomorphic building in
  $[0,1]\times\RR$. Let $S_\infty$ denote the source of $(\pi_\DD\circ
  U)$. (We have not yet defined $U$, just $\pi_\DD\circ U$.)

  Let $V_p\subset\Sigma$ be a closed disk neighborhood of east
  $\infty$ over which $J$ is split,
  $J|_{V_p\times[0,1]\times\RR}=j_\Sigma\times j_\DD$. Let $W_p$ be
  the complement of a
  slightly smaller closed disk around east $\infty$. We will extract
  convergent subsequences over $V_p$ and $W_p$
  separately, and then use unique continuation to conclude that they
  agree on the overlap.

  \textbf{Step 2.} \emph{Convergence over $V_p$.} Let $T_n = (\pi_\Sigma\circ
  u_n)^{-1}(V_p)$. We have holomorphic maps
  \[
  (\pi_\Sigma\circ u_n)|_{T_n}\co (T_n,\bdy T_n)\to (V_p, \bdy V_p\cup \alphas).
  \]
  By Theorem~\ref{thm:BEHWZ-2D-rel-cpct} the sequence
  $(\pi_\Sigma\circ u_n)|_{T_n}$ has a convergent subsequence. (Here,
  the energy bound comes from the fact that we are working in
  a fixed homology class $B$, so all holomorphic curves have the same area.)
  Replace $\{u_n\}$ by the corresponding subsequence, so
  $(\pi_\Sigma\circ u_n)|_{T_n}$ converges. Let $(\pi_\Sigma\circ
  U)|_{T}$ denote the limit of this subsequence; this is a
  (horizontal) holomorphic building in $V_p$ (with components
  extending into east $\infty$, $\RR\times Z$). Let
  \[
  T'=S_\infty\setminus [(\pi_\Sigma\circ U)|_T]^{-1}(W_p)
  \]
  be the subset of the limit surface $S_\infty$ which is mapped away from $e\infty$.

  \textbf{Step 3.} \emph{Convergence over $W_p$, i.e., to
  $T'$.} From Step 1, the sequence of decorated sources $S_n$ is
  converging (in the Deligne-Mumford moduli space) to $S_\infty$.
  This means that after passing to a subsequence we can view each of the $S_n$ as a
  complex structure on a single smooth surface (with boundary and
  punctures) $S$, and as $n\to\infty$ these complex structures
  converge in the $C^\infty$ topology away from some finite collection
  of pairwise-disjoint, embedded \emph{collapsing curves}. Near these
  collapsing curves, the complex structures are collapsing (or,
  equivalently, forming long necks).

  Let $q\in T'$ be a smooth point. Choose a neighborhood $X\ni
  q$ whose closure is compact (i.e., which does not touch the
  punctures in $S_\infty$) and contained in $T'$. We
  want to define $\pi_\DD\circ U$ on $X$. To this end, note that since
  \[
  \pi_2(\Sigma\times[0,1]\times\RR)=\pi_2(\Sigma\times[0,1]\times\RR,(\alphas\times\{1\}\times\RR)\cup(\betas\times\{0\}\times\RR))=0,
  \]
  Gromov's bubbling lemma (e.g., \cite[Lemma
  5.11]{BEHWZ03:CompactnessInSFT}) implies that $\|d(u_n)\|_{L^\infty}$ is
  bounded on $X$. The facts that the curves $\pi_\DD\circ u_n$
  converge and that the curves $u_n$ all represent the homology class
  $B$ gives a bound on the energy of the curves $u_n$. Consequently,
  Theorem~\ref{thm:MS-cpct} applies with $\Sigma=X$ to give a
  subsequence of the $u_n|_{X}$ converging to a holomorphic curve
  \[
  U|_{X}\co X\to W_p\times[0,1]\times\RR\subset \Sigma\times[0,1]\times\RR.
  \]
  (The reader concerned about the non-compactness of
  $W_p\times[0,1]\times\RR$ should note that by assumption the
  holomorphic curves $u_n|_{X}$ map to a compact subset $Y$ of
  $W_p\times[0,1]\times\RR$. Let $\gls*{Sigmae}$
  be the result of filling in the puncture $e$ to got a closed Riemann
  surface. One can embed $Y$ in $\Sigma_{\overline{e}}\times S^2$, the
  product of $\Sigma$ and the two-sphere, in such a way that the almost
  complex structure $J$ extends to an almost complex structure on all
  of $\Sigma_{\overline{e}}\times S^2$ which is tamed by the split
  symplectic form.)

  Choose a countable collection of such $X$ which covers the smooth
  part of $T'$ and take a diagonal subsequence. This subsequence then
  converges in $C^\infty_{\mathrm{loc}}$ away from the collapsing curves.

  \textbf{Step 4.} \emph{Nodes and punctures.}  The argument in Step 3
  defines $U$ on all of $T'$ except the nodes. The energy bound
  (coming from the fact that the curves $u_n$ all represent $B$) and
  the removable singularities theorem (see, e.g., \cite[Theorem
  4.1.2]{MS04:HolomorphicCurvesSymplecticTopology}) imply that $U$
  extends to the nodes, as well. To see that $U$ is continuous across
  the nodes (i.e., approaches the same value from both sides of the
  nodes) requires a slight further argument, for which we refer the
  reader to~\cite[Section
  4.7]{MS04:HolomorphicCurvesSymplecticTopology}).

  Next, we claim that near the punctures of $U|_{T'}$,
  $\pi_\Sigma\circ U$ converges to points in $\alphas\cap\betas$. This
  again follows from the fact that the maps $u_n$ all represent the
  homology class $B$.
  (Since $\pi_\Sigma\circ U$ is not
  holomorphic with respect to a single complex structure on $\Sigma$
  but rather with respect to a family parameterized by $[0,1]$, this
  claim uses a convergence result for strip-like ends, such
  as~\cite[Theorem 2]{Floer88:unregularized}.)

  \textbf{Step 5.} \emph{Piecing together $T$ and $T'$.} On $T\cap T'$ we
  have defined $U$ twice: once in Steps 1 and 2 and once in Steps 1, 3
  and 4, but it follows from the fact that $C^\infty_{\mathrm{loc}}$
  is Hausdorff that these definitions coincide. So, the maps $U|_T$ and $U|_{T'}$ glue to give
  a map on all of the components of $S_\infty$ mapped to
  $\Sigma\times[0,1]\times\RR$.

  Since the homology class of $u$ is determined by the local
  multiplicities at the $p_r$, it follows from the presence of the
  marked points we added in Step 1 that the curve $U$ represents the
  homology class $B$.

  \textbf{Step 6.} \emph{Components at east $\infty$.}  Finally, it remains
  to define $\pi_\DD\circ U$ on the components of $S_\infty$ mapped to
  east $\infty$ (i.e., $\RR\times Z\times[0,1]\times\RR$); these
  components of $S_\infty$, and the restriction of $\pi_\Sigma\circ U$
  from them to $\RR\times Z$, appeared in Step 2.  So, fix such a
  component $C$ of $S_\infty$. There are two cases: either the
  component $C$ appeared already in Step 1 or it did not. (Roughly, in
  the language of Definition~\ref{def:converge-to-comb}, these two
  cases correspond to $C$ being $\DD$-stable or $\DD$-unstable,
  respectively.) In the first case, Step 1 defines a map $C\to
  [0,1]\times\RR$. 
  Recall that the surface $S_\infty$ is obtained by collapsing some arcs in
  $S$, which we called collapsing curves. 
  Let $C_0$ be the preimage in $S$ of $C\subset S_\infty$. Then, since
  $C$ is mapped to $e\infty$, the boundary of $C_0$ consists of some
  collapsing arcs and some arcs in $S$ mapped by $\pi_\Sigma\circ u_n$
  to the $\alpha$-curves. It follows that $\pi_\DD\circ U$ maps $\bdy
  C$ to $\{1\}\times\RR$. The open mapping principle then implies that
  $(\pi_\DD\circ U)|_C$ is constant.
  %

  In the second case we have not yet defined $\pi_\DD\circ U|_C$. We
  can add extra marked points to the curves $\pi_\DD\circ u_n$ so that
  $C$ (with some extra marked points) does appear in the
  limit (cf.~\cite[Section 5.1]{MS04:HolomorphicCurvesSymplecticTopology}). Dropping these extra marked points makes $(\pi_\DD\circ
  U)|_C$ unstable (or else we would have been in the previous case),
  so in particular $(\pi_\DD\circ U)|_C$ is constant.
  
  We have now defined $U$ completely. Moreover, it is immediate from
  the construction that the sequence $u_n$ converges to $U$ in the
  sense of Definition~\ref{def:converge-to-comb}.
\end{proof}

\section{Gluing results for holomorphic combs}\label{sec:combs-gluing}
We next turn to gluing results, which show that the space
$\ocM^B(\x,\y;\Source;P)$ is well-behaved near some of the simplest
holomorphic combs.
\index{gluing!holomorphic curves}%
\begin{proposition}\label{prop:gluing_two_story}
  Let $(u_1,u_2)$ be a height $2$ holomorphic comb with
  $u_1\in\cM^{B_1}(\x,\penalty 500\y\semico\SourceSub{1}\semico P_1)$ and $u_2\in\cM^{B_2}(\y,\w\semico\SourceSub{2}\semico P_2)$. Then for
  sufficiently small open neighborhoods $U_1$ of $u_1$
  and $U_2$ of $u_2$, there is an open neighborhood of
  $(u_1,u_2)$ in
  $\ocM^{B_1*B_2}(\x,\w\semico\SourceSub{1}\glue\SourceSub{2}\semico P_1\cup P_2)$
  which is homeomorphic to $U_1\times U_2\times[0,1)$.
\end{proposition}

\begin{proposition}\label{prop:gluing_simple_comb} Let $(u,v)$ be a simple holomorphic comb with
  $u\in\cM^B(\x,\y\semico\Source)$ and $v\in\cN(\biSource\semico P_e)$.  Let $m =
  \abs{E(\Source)} = \abs{W(\biSource)}$.  Assume that the moduli spaces
  $\cM^B(\x,\y\semico\Source)$ and $\cN(\biSource\semico P_e)$ are transversally cut out at $u$ and
  $v$ respectively, and that $\ev\co\cM^B(\x,\y\semico\Source)\to\RR^m/\RR$ and
  $\ev_w\co\cN(\biSource\semico P_e)\to\RR^m/\RR$ are transverse at
  $(u,v)$. Then, for sufficiently small open neighborhoods $U_u$
  of~$u$ and $U_v$ of~$v$, there is an
  open neighborhood of $(u,v)$ in $\ocM^B(\x,\y\semico\Source\glue\biSource\semico P_e)$
  which is homeomorphic to $(U_u \times_{\ev}U_v)\times[0,1)$.
\end{proposition}
\begin{proof}[Proof of Propositions~\ref{prop:gluing_two_story} and~\ref{prop:gluing_simple_comb}]
  The results follow by standard arguments, dating back to Taubes’s
  work on gluing instantons~\cite{Taubes82:YM,
    Taubes84:YM} (see also~\cite[Chapter 7]{DonKron},~\cite[Chapter
  19]{KronheimerMrowka}). See, for instance,
  McDuff-Salamon \cite[Chapter
  10]{MS04:HolomorphicCurvesSymplecticTopology} or
  Pardon \cite[Appendix B]{Pardon:virtual} for a detailed explanation
  in the context of Gromov-Witten theory, and Pardon \cite[Appendix
  C]{Pardon:virtual}, \cite[Section 5]{Pardon:virtual-contact} for a
  detailed explanation for Hamiltonian Floer homology and contact
  homology, respectively. Since each puncture of a holomorphic curve
  in our setting is mapped either to $\pm\infty$ or to east or west
  $\infty$, the fact that we have two kinds of infinities is
  irrelevant to a local statement of this kind, and it follows from
  the corresponding result for the Morse-Bott case of relative
  symplectic field theory.

  In brief, to prove Proposition~\ref{prop:gluing_two_story}, say, one
  applies the implicit function theorem to define the gluing map
  \newcommand{\GlueMap}{\mathrm{Glue}}
  \[
  \GlueMap\co U_1\times U_2\times(0,\infty)\to \cM^{B_1*B_2}(\x,\w\semico\SourceSub{1}\glue\SourceSub{2}\semico P_1\cup P_2),
  \]
  (with $\infty$ corresponding to the broken curve, say).  The
  necessary estimates to apply the implicit function theorem can be
  obtained by adapting Bourgeois’s estimates~\cite[Section
  5.3]{Bourgeois02:MorseBott} to the relative case or Lipshitz’s
  estimates~\cite[Proposition A.1]{Lipshitz06:CylindricalHF} to the
  Morse-Bott case. One then shows that the gluing map is injective and
  (locally) surjective for sufficiently large gluing parameters.
  Local injectivity, i.e., injectivity of the restriction of
  $\GlueMap$ to subsets of the form $U_1\times
  U_2\times(t-\epsilon,t+\epsilon)$ for $\epsilon$ small, follows from
  the implicit function theorem.
  One way to deduce global injectivity is to mark points $p_1, p_2\in\Sigma$ such
  that $(\pi_\Sigma\circ u_i)^{-1}(p_i)$ is non-empty, choose a
  preimage $q_i$ of $p_i$ under $\pi_\Sigma\circ u_i$, and observe
  that for each $(u'_1,u'_2)\in U_1\times U_2$ the map
  $(0,\infty)\to(0,\infty)$,
  \[
  r \mapsto t(\GlueMap(u_1,u_2,r)(q_2))-t(\GlueMap(u_1,u_2,r)(q_2))
  \]
  is injective; compare~\cite[Lemma
  5.9]{Pardon:virtual-contact}. (Recall that $t$ is the
  $\RR$-coordinate on $\Sigma\times[0,1]\times\RR$.)  Surjectivity
  uses exponential convergence of holomorphic curves to Reeb
  chords, proved in the absolute case by Bourgeois~\cite[Section
  3.3]{Bourgeois02:MorseBott} (see also~\cite[Appendix
  A]{BEHWZ03:CompactnessInSFT}) and in the relative case by
  Abbas~\cite{Abbas04:strips,Abbas14:compactness}; see,
  e.g.,~\cite[Section C.11]{Pardon:virtual} for a nice explanation of
  surjectivity in a closely related context.
\end{proof}

We will also
use the following generalization of Proposition~\ref{prop:gluing_simple_comb} to
non-simple, height $1$ holomorphic combs in the proof of
Proposition~\ref{prop:gluing_degree_one}.
\index{gluing!holomorphic curves}%
\begin{proposition}\label{prop:gluing_height_1_comb}
  Consider a height~$1$ holomorphic
  comb $(u,v_1,\dots,v_k)$ with $u\in\cM^B(\x,\y\semico \Source)$, $v_i\in\cN(\biSource_i)$ for
  $i = 1,\dots,k-1$ and $v_k\in\cN(\biSource_k;P_e)$.  Let $m_i =
  \abs{W(\biSource_i)}$.  Assume that all the moduli spaces are
  transversally cut out, and that the map
  \begin{gather*}
\makebox[\textwidth][s]{$\displaystyle
    \cM^B(\x,\y;\Source)\times\cN(\biSource_1)\times\dots\times
        \cN(\biSource_k;P_e)\to
      \RR^{m_1}\times\RR^{m_1}\times\RR^{m_2}\times\RR^{m_2}\times
\dots\times
        \RR^{m_k}\times\RR^{m_k}$}\\
    (u,v_1,\dots,v_k)\mapsto
    (\ev(u),\ev_w(v_1),\ev_e(v_1),\ev_w(v_2),\dots,\ev_e(v_{k-1}),\ev_w(v_k))
  \end{gather*}
  is transverse to the diagonal $\{\,(x_1,x_1,\dots,x_k,x_k)\mid x_i \in
  \RR^{k_i}\,\}$.
  Then, for sufficiently small open neighborhoods $U_u$
  of~$u$ and $U_{v_i}$ of~$v_i$, there is an
  open neighborhood of $(u,v_1,\dots,v_k)$ in
  $\ocM^B(\x,\y\semico\Source\glue\biSource_1\glue\dots\glue\biSource_k\semico P_e)$
  which is homeomorphic to $(U_u
  \times_{\ev}U_{v_1}\times_{\ev}\dots\times_{\ev}U_{v_k})\times[0,1)^k$.
\end{proposition}
Again, the proof is essentially standard.

As mentioned at the end of Section~\ref{sec:moduli-overview}, we cannot prove that
the compactified moduli spaces are manifolds with corners in general.
We will get around this difficulty by proving that certain evaluation
maps are degree one, in an appropriate sense, at the corners; this
will be enough for the results in
Section~\ref{sec:degenerations-holomorphic-curves}. Before stating the
propositions, we make some definitions.

\begin{definition}\label{def:smeared-nhbd}Let $(u,v)$ be a simple
  holomorphic comb with $v\in\cN(\biSource)$. A \emph{smeared
    neighborhood of $(u,v)$ in
    $\ocM^B(\x,\y\semico\Source\semico P)$} is an open neighborhood of
  \index{smeared neighborhood}%
  \[
  \{\,(u,v') \mid v'\in\cN(\biSource), (u,v')\in\ocM^B(\x,\y;\Source;P)\,\}.
  \]
  There is an exactly analogous notion of a \emph{smeared neighborhood
  of $(u,v)$ in $\ocM^B(\x,\y\semico\Source\semico\vec{P})$}.
\end{definition}

\begin{remark}
  The notion of a smeared neighborhood of $(u,v)$ depends only on the
  combinatorics of (the decorated source of) $v$, not the map $v$ itself.
\end{remark}

\begin{definition}Given a continuous map $f\co X\to Y$ of topological
  spaces and a point $q\in Y$, we say $f$ is \emph{proper near $q$}
\index{proper near $q$}%
  if there is an open neighborhood $U\ni q$ such that $f|_{f^{-1}(U)}\co
  f^{-1}(U)\to U$ is proper.
\end{definition}

We will use the following weak notion of stratified spaces.

\begin{definition}
  \index{stratified!space}%
  \index{stratified!map}%
  An {\em $n$-dimensional stratified space} is a topological space $X$
  written as a union of strata $\{X_i\}_{i=0}^n$ where $X_i$ is a
  smooth $i$-dimensional manifold and the closure of $X_k$ is
  contained in $\bigcup_{i\leq k} X_i$. We typically suppress the
  strata from the notation.  Let $X$ and $Y$ be stratified spaces. A
  \emph{stratified map} $f\co X \to Y$ is a continuous map with the
  property that the preimage of any stratum in $Y_i\subset Y$ is a
  union of connected components of strata $X_j$ of $X$, and the
  restriction of $f$ to each stratum of $X_j$, thought of as a map
  into $Y_i$, is a smooth map.
\end{definition}

\begin{definition}
	\label{def:OddDegree}
  Let $X$ be a stratified space so that the top
  stratum is a smooth $m$-manifold, and let $f\co X\to\RR^m_+$ be a
  stratified
  map so that $f^{-1}\left((0,\infty)^m\right)$ is the top stratum of
  $X$. Let $q\in\RR^m_+$, and assume $f$ is proper near $q$. We say
  that $f$ is \emph{odd degree near $q$}
\index{degree near $q$, odd}\index{odd!degree near $q$}%
  if there is an open
  neighborhood $U$ of $q$ such that for any regular value $q'\in
  U\cap(0,\infty)^m$, $f^{-1}(q')$ consists of an odd number of
  points.
\end{definition}

\begin{lemma}\label{lemma:stratified-degree-one}Let $X$ be a
  stratified space such that
  the union of the top two strata of~$X$ forms an $m$-manifold
    with boundary.
  Let $f\co X\to \RR^m_+$ be a stratified map
  which
  is proper near $0$. Assume that $f^{-1}\left((0,\infty)^m\right)$ is
  the top stratum of~$X$. Then
  \begin{enumerate}
  \item\label{item:strata-lemma-1} If $q,q'\in(0,\infty)^m$ near $0$ are regular values of $f$
    then $\abs{f^{-1}(q)}\equiv \abs{f^{-1}(q')} \pmod 2$.
  \item\label{item:strata-lemma-2} If the restriction of $f$ to the
    preimage of some facet $\RR^{m-1}_+$ of $\RR^m_+$ is odd degree near
    $0\in \RR^{m-1}_+$ then $f$ is odd degree near $0\in\RR^m_+$.
  \end{enumerate}
\end{lemma}
\begin{proof}
  Let $U$ be a neighborhood of $0$ over which $f$ is proper. Let
  $V$ be $f^{-1}\left(U\cap\RR_+^m\right)$. Then $f|_V\co V\to
  \left(U\cap(0,\infty)^m\right)$ is a proper map. So
  part~(\ref{item:strata-lemma-1}) follows from standard degree theory.

  For part~(\ref{item:strata-lemma-2}), let $q\in \RR^{m-1}_+$ be a
  regular value of $f|_V$, and let $B
  \subset\RR^{m-1}_+$ be a ball neighborhood
  of $q$ small enough that $f|_{f^{-1}(B)}$ is a covering map. Write
  $f^{-1}(B)=B_1\amalg\dots\amalg B_{2\ell+1}$; each $B_i$ is an
  $(m-1)$-ball. For each $i$ choose an $(m-1)$-ball $B_i'$ whose interior lies in the top
  stratum of $X$ and with $\bdy B'_i=\bdy B_i$. Note that by our
  assumptions on $f$, $f(\Int(B'_i))$ is entirely contained in
  $(0,\infty)^m$. So $f(B_i\cup B_i')$ is a (singular) $(m-1)$-sphere
  in $\RR^m_+$, and for $q'\in(0,\infty)^m$ sufficiently close to $q$,
  the sphere $f(B_i\cup B_i')$ has winding number $\pm1$ around $q'$.
  It
  follows that $f$ has odd degree at regular values near
  $q$. Choosing $q$ close enough to $0$ gives the result.
\end{proof}

\begin{proposition}\label{prop:gluing_degree_one}
  Suppose that $(u,v)$ is a simple holomorphic comb in
  $\ocM^B(\x,\allowbreak\y\semico\Source\semico P)$. Assume that
  $v$ is a split curve and
  that there are two parts $P_1$
  and $P_2$ of $P$ such that, for each split component of $\biSource\!$,
  its bottom puncture belongs to $P_1$ and its
  top puncture belongs to~$P_2$.
  Assume also that $\ind(B,\Source,P)=2$ so $\dim \cM^B(\x,\y\semico
  \Source \semico P) = 1$.

  Let $q_1\in P_2$ and $q_2\in P_1$ be the top and bottom punctures,
  respectively, on one
  of the split components of~$\biSource$. Then, for generic
  $J$, there is a smeared neighborhood~$U$ of $(u,v)$ in
  $\ocM^B(\x,\y\semico\Source\semico P)$ so that the
  evaluation map $\widebar{\ev}_{q_1,q_2}\co U\to\RR_+$
  is proper near~$0$ and odd degree near~$0$.
\end{proposition}
\begin{proof}
  To see the $\widebar{\ev}$ is proper near~$0$ on~$U$, note that
  $\widebar{\ev}_{q_1,q_2}$ is clearly proper on the closure $\widebar{U}$
  of~$U$.  By examining the topology of $\Source$, we see that
  $\widebar{U} \setminus U$ is mapped to strictly positive values by
  $\widebar{\ev}_{q_1,q_2}$.

  Let the split components of $\biSource$ be $\biSource_i$ for $i =
  1,\dots,N$.  For notational convenience, we will assume that
  $\biSource$ has no trivial strips, so that the only parts of $P$ are
  $P_1$ and~$P_2$.  Let $P'$ be the partition~$P$ with $P_2$ split
  into the discrete partition; that is, $P'$ has $N+1$ parts, one of
  which is $P_1$ and all others consisting of a single puncture.  We
  will show that there is a smeared neighborhood $U'$ of $(u,v)$ in
  $\ocM^B(\x,\y\semico\Source\semico P')$ so that the evaluation map
  $\widebar{\ev} \co U' \to\RR_+^N$
  is odd degree near~$0$ in the sense of
  Definition~\ref{def:OddDegree}.  (We see that $\widebar{\ev}$ is
  proper near~$0$ as before.) This suffices, since for
  generic~$J$, the map $\ev$ from $\cM^B(\x,\y\semico\Source\semico P')$ is
  transverse to the diagonal in $(0,\infty)^N$ by
  Proposition~\ref{prop:transversality}.

  By Proposition~\ref{prop:gluing_height_1_comb}, we can make a smeared
  neighborhood~$U$ of $(u,v)$ inside $\ocM^B(\x,\y\semico\Source\semico P')$ into a
  stratified space of a particularly nice kind; in particular, the
  union of the top two strata is an $N$-manifold with boundary.
  The codimension~1 strata of~$U$ (of dimension $N-1$) consist of
  simple holomorphic combs $(u',v')$, where $v'$ is a split curve with
  no more than $N$ split components.  Higher codimension strata of~$U$
  are more complicated height $1$ holomorphic combs.

  We call such a neighborhood in an $N$-dimensional
  compactified moduli space $\ocM^B(\x,\y\semico\Source\semico P')$ a \emph{split
    neighborhood}: namely, a split neighborhood is an open set
  containing the entire stratum of simple combs $(u,v)$, where $v$ is
  a split curve with $N$ split components, with all bottom punctures
  of~$v$ belonging to the same part of~$P'$ and all top punctures
  belonging to different parts.

  We will show by induction on $N$ that $\widebar{\ev}$ on a split
  neighborhood has odd degree near~$0$.  For $N=1$, by
  Proposition~\ref{prop:gluing_simple_comb} we have a map from $[0,1)$
  to $\RR_+$ mapping $0$ to~$0$ which is proper near~$0$; such a map
  automatically has odd degree near~$0$.  For
  $N>1$, note that the preimage of the interior of a facet of
  $\RR_+^N$ consists of holomorphic curves $(u',v')$, where $v'$ has a
  single split component; thus, $u'$~lives in a split neighborhood~$U'$
  of a smaller moduli space, of dimension~$N-1$.  By induction,
  $\widebar{\ev}$ on $U'$ has odd degree near~$0$, and so by
  Lemma~\ref{lemma:stratified-degree-one}, $\widebar{\ev}$ on~$U$ also
  has odd degree near~$0$.
\end{proof}

In Chapter~\ref{chap:tensor-prod}, we will also need the following
generalization of Proposition~\ref{prop:gluing_degree_one}. To state
it, we need a little more terminology. A \emph{generalized split
  curve} is a union of holomorphic disks in
$\RR\times Z\times[0,1]\times\RR$ each of which has exactly one $w$
puncture. Each component of a generalized split curve is a
\emph{generalized split component.}
\index{component!generalized split}%
\index{generalized!split component|see{component, generalized split}}%
\index{curve!generalized split}%
\index{generalized!split curve|see{curve, generalized split}}%
\begin{proposition}\label{prop:generalized_gluing_degree_one}
  Let $(u,v)$ be a simple holomorphic comb in
  $\ocM^B(\x,\y\semico\Source\semico\vec{P})$, where $v$ is a generalized split
  curve. Write $\vec{P}=(P_1,\dots,P_k)$, label the components of~$v$
  as $T_1,\dots,T_m$ and label the east punctures of $T_i$ by
  $q_{i,1},\dots,q_{i,n_i}$. Assume that there are integers
  $0=\ell_1,\ell_2,\dots,\ell_k,\ell_{k+1}=m$ such that
  $P_j=\{q_{\ell_j+1,1},\dots,\allowbreak q_{\ell_{j+1},n_{\ell_{j+1}}}\}$. (That
  is, $P_j$ consists of the punctures on components
  $T_{\ell_j+1},\dots,T_{\ell_{j+1}}$.)
  
  Suppose that $\ind(B,\Source,P)=1$. 
  
  Then there is a smeared neighborhood $U$ of $(u,v)$ in $\ocM^B(\x,\y;\Source)$ so that the evaluation map
  \begin{multline*}
  \ev\co U\to\RR_+^{n_1-1}\times\RR\times\RR_+^{n_2-1}\times\RR\times\dots\times\RR_+^{n_{\ell_2}-1}\times\\
  \RR_+^{n_{\ell_2+1}-1}\times\RR\times\cdots\times\RR_+^{n_{\ell_3}-1}\times\cdots\cdots\times\RR_+^{n_m-1}
  \end{multline*}
  defined by 
  \begin{multline*}
    \ev(u)=(\overbrace{\ev_{q_{1,2}}(u)-\ev_{q_{1,1}}(u), \dots,
      \ev_{q_{1,n_1}}(u)-\ev_{q_{1,n_1-1}}(u)}^{\in\RR_+^{n_1-u}},\\
     \overbrace{\ev_{q_{2,1}}(u)-\ev_{q_1,n_1}(u)}^{\in\RR},
     \overbrace{\ev_{q_{2,2}}(u)-\ev_{q_{2,1}}(u), \dots, \ev_{q_{2,n_2}}(u)-\ev_{q_{2,n_2-1}}(u)}^{\in\RR_+^{n_2-1}},\\\overbrace{\ev_{q_{3,1}}(u)-\ev_{q_2,n_2}(u)}^{\in\RR},\dots)
  \end{multline*}
  is proper near $0$ and odd degree near $0$.
\end{proposition}
\begin{proof}
  The proof is similar to the proof of Proposition~\ref{prop:gluing_degree_one}. The moduli space of generalized split curves has ends where one degenerates off a split curve at (far) east or west infinity; the lowest-dimensional stratum consists of sequences of (ordinary) split curves. Using this, one can inductively reduce to the case of gluing on split curves. Details are left to the reader.
\end{proof}

\begin{remark}\label{remark:gluing-is-hard}It would be convenient to extend
  Proposition~\ref{prop:transversality} to the compactified moduli
  spaces, and in particular to assert that for height $1$ holomorphic
  combs, say, the extension of $\ev$ to $\ocM^B(\x,\y;\Source)$ is
  transverse to the partial diagonals $\Delta_P$. In particular, this
  would eliminate the need for the cumbersome
  Propositions \ref{prop:gluing_degree_one}, \ref{prop:generalized_gluing_degree_one}
  and~\ref{prop:gluing-shuffle}. Such a statement would
  require a smooth structure on the $\ocM^B(\x,\y;\Source)$, which we
  have not constructed. Moreover, the validity of the statement would
  depend strongly on the smooth structure: in the language of
  Hofer~\cite{Hofer06:FredholmPolyfolds}, it depends on the
  \emph{gluing profile.} (In fact, for the exponential gluing profile
  advocated there, the extension of $\ev$ would generally not be
  transverse to $\Delta_P$.)
  \index{gluing!profile}%
\end{remark}

We conclude this section with a result on shuffle curves. It states
that, when they occur, odd shuffle curves appear in the boundary of
$\ocM^B(\x,\y\semico\Source\semico P)$ (with $\ind(B,\Source,P)=2$) an odd number
of times, and even shuffle curves appear an even number of times.
\begin{proposition}\label{prop:gluing-shuffle}Fix a generic $J$. Let
  $(u,v)\in\ocM^B(\x,\y\semico\Source\semico P)$ be a simple holomorphic comb with
  $v$ a shuffle curve. Assume that the two $e$ punctures $q_i$ and
  $q_j$ of the shuffle component belong to the same part $P_i$ of
  $P$. Assume also that $\ind(B,\Source,P)=2$.
  \begin{enumerate}
  \item\label{item:odd}
    \index{curve!shuffle!odd}%
    If $v$ is an odd shuffle curve then there is a
    smeared neighborhood $U$ of $(u,v)$ in $\ocM^B(\x,\y\semico\Source\semico P)$
    so that $\overline{U}\setminus \bdy\ocM^B(\x,\y\semico\Source\semico P)$
    is homeomorphic to $\coprod_{i=1}^{2m-1}[0,1)$ (for some $m\in
    \mathbb{N}$).
  \item\label{item:even}
    \index{curve!shuffle!even}%
    If $v$ is an even shuffle curve then there is
    a smeared neighborhood $U$ of $(u,v)$ in $\ocM^B(\x,\y\semico\Source\semico P)$
    so that $\overline{U}\setminus \bdy\ocM^B(\x,\y\semico\Source\semico P)$
    is homeomorphic to $\coprod_{i=1}^{2m-2}[0,1)$ (for some $m\in
    \mathbb{N}$).
  \end{enumerate}
\end{proposition}
\begin{proof}
  First some notation: let
  $a_1,\dots,a_4$ be the endpoints of the Reeb chords on the shuffle
  component, ordered so that $a_1\lessdot\dots\lessdot a_4$. Let
  $\rho_{ij}$ (with $i < j$ and $i,j\in\{1,\dots,4\}$) denote the Reeb chord
  from $a_i$ to $a_j$.  Let $\biSource$ denote the source of $v$

  The two cases are similar; we first prove case~(\ref{item:odd}).
  In this case, the $e$ punctures of the shuffle component
  of~$\biSource$ are
  labeled by $\rho_{14}$ and $\rho_{23}$. The $w$ punctures are
  labeled by $\rho_{13}$ and $\rho_{24}$.
  The moduli space
  $\cN(\biSource)$ has a natural compactification
  $\ocN(\biSource)$ which is homeomorphic to $[0,1]$.  (The
  parameter is obtained by moving the branch point(s); see
  Figure~\ref{fig:shuffle-moduli} for a schematic illustration.)  By
  our transversality
  (Propositions~\ref{prop:transversality}
  and~\ref{prop:east_transversality}) and gluing
  (Proposition~\ref{prop:gluing_simple_comb}) results, any curve in
  $\cN(\biSource)$ can
  be glued to~$u$ to obtain a curve in $\cM^B(\x,\y;\Source)$.

  To see when there exist such curves in $\cM^B(\x,\y\semico\Source\semico P)$ we
  consider the ends of the moduli space $\ocN(\biSource)$. The two ends
  correspond to holomorphic combs $(u,v_1,v_2)$ and $(u,v_3,v_4)$. The
  curves $v_1,\dots,v_4$ are as follows:
  \begin{itemize}
  \item $v_1$ is a split curve with one split component, labeled by
    $(w,\rho_{13})$, $(e,\rho_{12})$ and $(e,\rho_{23})$, together
    with a trivial component labeled by $(e,\rho_{24})$ and $(w,\rho_{24})$.
  \item $v_2$ is a join curve, with a join component labeled by
    $(w,\rho_{24})$, $(w,\rho_{12})$ and $(e,\rho_{14})$, together
    with a trivial component labeled by $(e,\rho_{23})$ and $(w,\rho_{23})$.
  \item $v_3$ is a split curve with one split component, labeled by
    $(w,\rho_{24})$, $(e,\rho_{23})$ and $(e,\rho_{34})$, together
    with a trivial component labeled by $(e,\rho_{13})$ and $(w,\rho_{13})$.
  \item $v_2$ is a join curve, with a join component labeled by
    $(w,\rho_{34})$, $(w,\rho_{13})$ and $(e,\rho_{14})$, together
    with a trivial component labeled by $(e,\rho_{23})$ and $(w,\rho_{23})$.
  \end{itemize}
  See Figure~\ref{fig:shuffle-moduli}.

  Suppose one glues $v_1$ to $u$, with some gluing parameter
  $\epsilon_1$. Then in the resulting curve the Reeb chord $\rho_{23}$
  will be above $\rho_{12}$. If one then glues $v_2$ to the result,
  with some gluing parameter $\epsilon_2\ll\epsilon_1$, in the
  resulting curve the Reeb chord $\rho_{23}$ will be above
  $\rho_{14}$.

  By contrast, suppose one glues $v_3$ to $u$, with some gluing parameter
  $\epsilon_1$. Then in the resulting curve the Reeb chord $\rho_{34}$
  will be above $\rho_{23}$. If one then glues $v_2$ to the result,
  with some gluing parameter $\epsilon_2\ll\epsilon_1$, in the
  resulting curve the Reeb chord $\rho_{14}$ will be above
  $\rho_{23}$.
  
  So, since the two Reeb chords switched order as we moved from one
  end of $\ocN(\biSource)$ to the other, for given $\epsilon$
  there is algebraically one curve $v'\in\cN$ so that gluing $u$
  with~$v'$ with gluing parameter $\epsilon$ gives a curve with $\rho_{14}$
  and $\rho_{23}$ at the same height---i.e., an element of
  $\cM^B(\x,\y\semico\Source\semico P)$. Part~(\ref{item:odd}) of the proposition follows.

  \begin{figure}
    \includegraphics[scale=.55556]{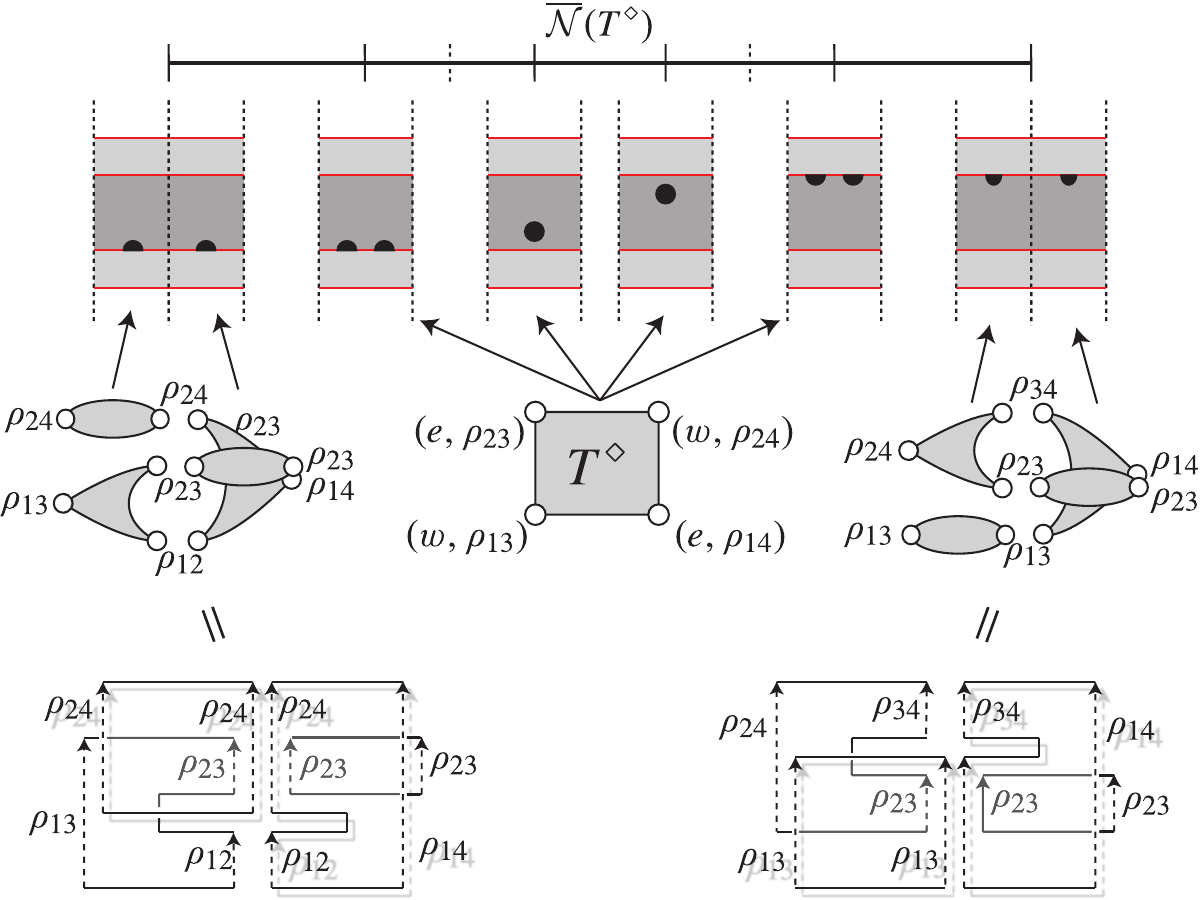}
    \caption[$1$-parameter family of maps $v$ from an odd
      shuffle component]{\textbf{The one-parameter family of maps $v$ from an odd
      shuffle component.} The dark dots denote interior branch points,
      and the dark semi-dots denote boundary branch points. The two
      ends of the moduli space are shown on the left and the right.
      At the bottom, the two ends are sketched in the
      style of Figure~\ref{fig:shuffle1}. The corresponding pictures for an even
      shuffle component are obtained by
      reflection.\index{drop shadow}}\label{fig:shuffle-moduli}
  \end{figure}
  \colorused

  Part~(\ref{item:even}) follows by a similar analysis. In this case,
  the $e$ punctures of the shuffle component are labeled by
  $\rho_{13}$ and~$\rho_{24}$, and the $w$ punctures are labeled by
  $\rho_{14}$ and~$\rho_{23}$.  A schematic illustration of
  $\ocN(\biSource)$ may be obtained by turning
  Figure~\ref{fig:shuffle-moduli} upside down.  An analysis as above
  shows that at both ends of $\ocN$, the Reeb chord $\rho_{24}$ is
  above $\rho_{13}$.  This implies that, algebraically, there
  are zero curves in $\cM^B(\x,\y\semico\Source\semico P)$ for a given gluing
  parameter.
\end{proof}

\section{Degenerations of holomorphic
  curves}\label{sec:degenerations-holomorphic-curves}

To construct our invariants we study degenerations of 1-dimensional moduli spaces.
The compactness theorem, Proposition~\ref{prop:compactness}, permits two kinds of degenerations to occur in $\ocM^B(\x,\y\semico\Source\semico P)$:
\begin{itemize}
\item\index{degeneration of holomorphic curves} The source $S$ can degenerate to a point in the boundary of the moduli space of Riemann surfaces. Such degenerations correspond to the conformal structure pinching along some circles and/or arcs, resulting in a surface with some nodes $q_i$. There are three possible behaviors of the limit near each $q_i$: 
  \begin{itemize}
  \item The $\RR$-coordinate of the map may converge to $\pm\infty$ at
    $q_i$, so $q_i$ corresponds to a level-splitting in the resulting
    height $N>1$ holomorphic comb.  This is a case of \emph{splitting at
      $\pm\infty$.}
    \index{degeneration of holomorphic curves!splitting at $\pm\infty$}%
    \index{splitting!at $\pm\infty$}%
  \item The holomorphic map may converge to the puncture $p$ in $\Sigma$ on one side of $q_i$, so the puncture $q_i$ corresponds to a splitting at east $\infty$; or the analogous behavior for combs with several east $\infty$ levels may occur. We call this \emph{splitting at east $\infty$.}
    \index{degeneration of holomorphic curves!splitting at east $\infty$}%
    \index{splitting!at east $\infty$}%
  \item The holomorphic map may extend continuously over $q_i$, sending $q_i$ to some point in $\Sigma\times[0,1]\times\RR$ (or perhaps $\RR\times Z\times[0,1]\times\RR$), resulting in a nodal holomorphic comb. We call this \emph{becoming nodal.}
    \index{degeneration of holomorphic curves!becoming nodal}\index{becoming nodal}\index{nodal!becoming}%
  \end{itemize}
\item The source $S$ may not degenerate, but the map could become
  singular. In this case, the derivative blows up at some point or
  points. (In the proof of Proposition~\ref{prop:compactness}, we
  added marked points in
  this case, so that it became a special case of the previous point.)
  There are three sub-cases.
  \begin{itemize}
  \item The derivative may blow up at a puncture of $S$. If the
    puncture is a $\pm$ puncture then this is the other case of
    \emph{splitting at $\pm\infty$}, where the curve that splits off
    has a non-stable source. If it were an $e$ puncture then
    this would be a case of splitting or east $\infty$, but there are
    no non-trivial curves at $e\infty$ with unstable sources, so this
    does not occur.
  \item  The derivative may blow up in the interior of $S$, which 
    results in \emph{bubbling a holomorphic sphere}.
    \index{bubble, sphere or disk}%
    \index{sphere bubble}%
  \item The derivative may blow up on the boundary of $S$, which
    results in \emph{bubbling a holomorphic disk}.
    \index{disk bubble}%
  \end{itemize}
\end{itemize}

The goal of this section is to restrict which of these degenerations
can occur in codimension $1$. Ultimately, we will rule out bubbling
disks and spheres and becoming nodal, and will restrict strongly which
curves can split off at $\pm\infty$ and $e\infty$. In
Section~\ref{sec:embedded-degen}, we show that a few more
degenerations can not occur for embedded curves. The remaining
degenerations will be incorporated in the algebra of our invariants,
in Chapters~\ref{chap:type-d-mod} and~\ref{chap:type-a-mod}.

\subsection{First restrictions on codimension-one degenerations}\label{sec:first-restrictions}
Let
\[
\gls*{BdyModSpaceOrderedCpct}
\coloneqq
\ocM^B(\x,\y\semico\Source\semico P)\setminus\cM^B(\x,\y\semico\Source\semico P).
\]

Some restrictions follow easily from the index formula and other simple considerations:
\begin{proposition}\label{prop:restrict_degens_1}
  Fix generalized generators $\x$, $\y$, $B\in\pi_2(\x,\y)$, decorated
  source $\Source$, and a partition $P$ so that
  $\cM^B(\x,\y\semico\Source\semico P)$ is $1$-dimensional. Then for
  generic $J$ every holomorphic comb in
  $\bdy\ocM^B(\x,\y\semico\Source\semico P)$ has one of the following forms:
\begin{enumerate}
\item a two-story holomorphic comb $(u_1,u_2)$;
\item \label{case:join_curve}a simple holomorphic comb $(u,v)$ where
  $v$ is a join curve;
\item \label{case:shuffle_curve}a simple holomorphic comb $(u,v)$
  where $v$ is a shuffle curve;
\item \label{case:split_curves}a height $1$ holomorphic comb
  $(u,v_1,\dots,v_k)$ such that each $v_i$ is a split curve, and,
  furthermore, the result
  $\biSource_1\glue\dots\glue\biSource_k$ of pregluing the sources
  of $v_1,\dots,v_k$ is also a split curve; or
\item \label{case:nodal}a nodal holomorphic curve, obtained by degenerating some arcs with boundary on $\bdy S$.
\end{enumerate}
\end{proposition}
(Recall the definitions of join curves, split curves and
shuffle curves from the end of
Section~\ref{sec:curves_at_east_infinity}; in particular, a join curve
or a shuffle curve has only one non-trivial component, while a
split curve may have arbitrarily many.)

\begin{proof}[Proof of Proposition~\ref{prop:restrict_degens_1}]
  We begin with some general remarks. Since $\pi_2(\Sigma)=0$,
  $\pi_2(\Sigma,\alphas)=0$ and $\pi_2(\Sigma,\betas)=0$, no disks or
  spheres can bubble off. As usual, Deligne-Mumford type
  degenerations, in which a closed circle or circles degenerate, are
  codimension $2$, and hence do not occur in these moduli spaces. By
  contrast, it is possible for an arc with boundary on $\bdy S$ to get
  pinched in codimension $1$ (and become a node); this is case~(\ref{case:nodal}).

  Next, suppose we have a simple holomorphic comb $(u,v)$ in
  $\bdy\ocM^B(\x,\y\semico\Source\semico P)$. Let $\Source'$ and $\biSource{}'$
  denote the sources of $u$ and $v$ respectively, and let $m$ be
  $\abs{W(\biSource{}')}$.

  The curve $v$ induces a partition $P'$ of $W(\biSource{}')$, by
  specifying that $q$ and $q'$ belong to the same part if they have the
  same $t$-coordinate in~$v$. Then, $u\in\cM^B(\x,\y\semico\Source'\semico P')$
  (using the identification of $E(\Source')$ and $W(\biSource{}')$).

  Since $t$ is constant on each component of $v$ and each component of
  $v$ with a $w$ puncture also has at least one $e$ puncture, $|P'|\leq|P|$. Also notice that
  $\chi(S')=\chi(S)-\chi(T')+m$. Since every component of $T'$ must
  have a west puncture, $\chi(T')\leq m$, with equality if and only
  $T'$ is a disjoint union of topological disks with one west puncture
  per component.

  Since, by hypothesis, the expected dimension of
  $\cM^B(\x,\y\semico\Source\semico P)$ is $1$, by Formula~(\ref{eq:Index}) we have
  \begin{align*}
  \ind(B,\Source,P)&=g-\chi(S)+2e(B)+|P|=2. \\
  \shortintertext{Now,}
  \ind(B,\Source',P')&=g-\chi(S')+2e(B)+|P'|\\
  &=g-\chi(S)+\chi(T')-m+2e(B)+|P'|\\
  &=2-|P|+\chi(T')-m+|P'|.
  \end{align*}
  Since $|P'|\leq|P|$ and $\chi(T')\leq m$, we therefore have
  $\ind(B,\Source',P')\leq 2$.

  Suppose that $\ind(B,\Source',P')=2$. Then $\chi(T')=m$ and
  $|P|=|P'|$. Consequently, $T'$ consists entirely of topological
  disks, each with a single puncture labeled~$w$. It follows from
  Propositions \ref{prop:east_transversality}
  and~\ref{prop:gluing_simple_comb} that $v$ can be glued to $u$,
  and so $\ind(B,\Source,P) > \ind(B,\Source',P')$, contradicting the
  hypothesis.

  Next, suppose that $\ind(B,\Source',P')=1$. There are now two cases:
  \begin{itemize}
  \item $|P|=|P'|$ and $\chi(T')=m-1$. In this case, either $T'$
    consists of $m-1$ topological disks and one topological annulus,
    or just $m-1$ topological disks. It is easy to see that the former
    does not occur: since $T'$ has $m$ $w$ punctures, in this case one
    boundary component of $T'$ would have no puncture. Thus, any curve
    close to this comb will still have a boundary component with no
    punctures.  The maximum modulus theorem forces $\pi_\DD\circ u$ to
    be constant on that component, a contradiction.  (See
    Figure~\ref{fig:annulus-degen}.)

    \begin{figure}
      \includegraphics[scale=.55556]{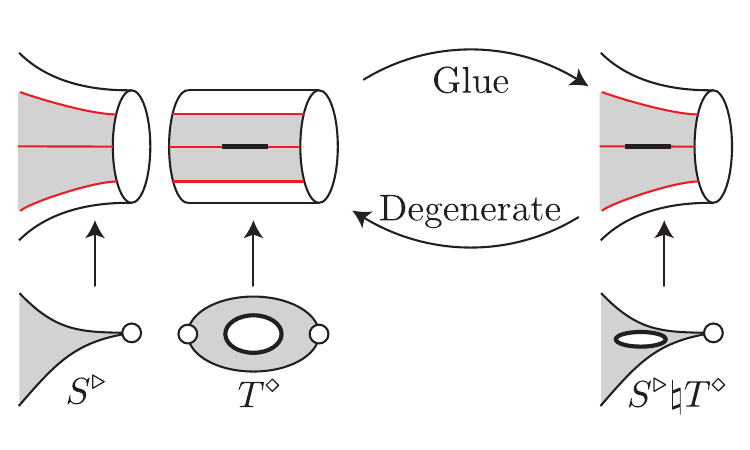}
      \caption[Annulus at east $\infty$ that does not exist]{\textbf{A putative annulus at east $\infty$.} Any near enough
        curve has a boundary component with no punctures, which
        violates the maximum modulus principle.}\label{fig:annulus-degen}
    \end{figure}
    \colorused

    By contrast, the case that $T'$ consists of $m-1$ topological
    disks can occur. In this case, all but one of these disks has one
    $w$ puncture, and one has two $w$ punctures. We claim that none of
    the disks has two consecutive $e$ punctures (not separated by a
    $w$ puncture). Indeed, any two consecutive $e$ punctures will have
    different heights in any curve with source
    $\Source'\glue\biSource{}'$. But this contradicts the fact that
    $|P|=|P'|$. It follows that $m-2$ of the disks in $\biSource{}'$
    are trivial components and the remaining disk is either a join
    component or a shuffle component.

  \item $|P'|=|P|-1$ and $\chi(T')=m$. In this case, $T'$ consists of
    $m$ topological disks, each with a single west puncture.
    On first glance, it appears each disk may have an arbitrary
    number of east punctures. However, by Propositions
    \ref{prop:transversality}, \ref{prop:east_transversality} and
    \ref{prop:gluing_simple_comb}, we can glue $u$ to $v$ in
    $\ocM(\x,\y\semico\Source)$ (though perhaps not in
    $\ocM(\x,\y\semico\Source\semico P)$). It is then easy to see that the east
    punctures on a single component of $T'$ will have different
    $t$-coordinates in all of the glued curves. It follows that if
    any component of $T'$ has more than two east punctures then
    $|P|-|P'|>1$, a contradiction.  Therefore $\biSource{}'$ is a
    split curve, as in
    Case~(\ref{case:split_curves}) of the statement.
  \end{itemize}

  If $\ind(B,\Source',P')\leq0$ then, by
  Proposition~\ref{prop:transversality}, the moduli space
  $\cM^B(\x,\y\semico\Source'\semico P')$ is empty.

  Next, consider a general height $1$ holomorphic comb
  $(u,v_1,\dots,v_\ell)$. The arguments above then imply that the
  bi-decorated source $\biSource_1\glue\dotsb\glue\biSource_\ell$
  obtained by pregluing the sources of $v_1,\dots,v_\ell$ is a
  join curve, split curve, or shuffle curve.

  Finally, consider a general holomorphic comb $U$ in
  $\bdy\ocM(\x,\y\semico\Source\semico P)$. It follows from
  Formula~\eqref{eq:Index} that each story of $U$ drops the expected
  dimension by~$1$, so there are at most two stories. Further, if
  there are two stories,
  each story individually has index~$1$, so by the arguments above
  cannot have any components at east~$\infty$.
\end{proof}

The first four kinds of degenerations will be involved in the
definitions of the invariants, in Chapters~\ref{chap:type-d-mod}
and~\ref{chap:type-a-mod}. So we now show that nodal degenerations
can not occur in certain parts of the moduli spaces. We give names to three different kinds of nodal degenerations:

\begin{definition}
  \label{def:BoundaryDegeneration}\index{boundary degeneration|see{degeneration, boundary}}\index{degeneration of holomorphic curves!boundary}
  A {\em boundary degeneration} is a nodal holomorphic comb with the property that some irreducible component
  of the source has no $+$ or $-$ punctures and is mapped via a
  non-constant map to $\Sigma$.
\end{definition}

\begin{definition}
  \label{def:BoundaryDoublePoints}\index{boundary double point|see{degeneration, boundary double point}}\index{double point, boundary|see{degeneration, boundary double point}}\index{degeneration of holomorphic curves!boundary double point}
  A {\em boundary double point} is a holomorphic comb whose source is
  a nodal curve and has a node $q$ on the boundary with the property that the projection to
  $[0,1]\times \RR$ is not constant near either preimage point $q_1$ or $q_2$ of $q$
  in the normalization of the nodal curve.
\end{definition}

\begin{definition}\index{component!ghost}\index{curve!haunted}\index{ghost component}\index{haunted curve}\label{def:ghost}
  A \emph{ghost component} of a holomorphic comb $u$ is a component $S_0$ of the source $S$ of $u$ such that $u|_{S_0}$ is constant. If $u$ has a ghost component then we say $u$ is \emph{haunted}.
\end{definition}

\begin{lemma}
  Every nodal degeneration results in a boundary degeneration, a curve with boundary double point, or a haunted curve.
\end{lemma}
\begin{proof}
  If a nodal curve has no ghost components and no boundary degenerations then $\pi_\DD\circ u$ is non-trivial on every component, and hence $u$ has a boundary double point.
\end{proof}

We rule out boundary degenerations immediately:

\begin{lemma}
  \label{lem:NoBoundaryDegenerations}
  Boundary degenerations do not exist.
  \index{degeneration of holomorphic curves!boundary!does not occur}%
\end{lemma}

\begin{proof}
  Consider an irreducible component $S$ of the source with no $+$ or $-$ punctures; and consider its corresponding map
  to $\Sigma\times [0,1]\times \RR$. Projecting further to $\RR$, we obtain a harmonic function on $S$ which is bounded
  on the boundary of $S$. Thus, this map is constant. It follows that the map $u$ maps into a single fiber $\Sigma$. Since we assumed that the map was non-trivial, its image gives a non-trivial relative homology class
  $B\in H_2(\Sigma,\alphas;\ZZ)$ with $n_z(B)=0$. Existence of such 
  a homology class violates homological linear independence of the $\alpha$-curves
  which our bordered diagrams are required to satisfy, Definition~\ref{def:BorderedDiagram}.
\end{proof}

\begin{remark}
  The non-existence of boundary degenerations is an artifact of our
  definition of bordered diagrams---in particular, the homological
  linear independence of the $\alpha$- (respectively $\beta$-)
  curves. In more general settings, one may wish to consider diagrams
  with an excess of $\alpha$- or $\beta$-curves. In those settings,
  boundary degenerations do play a role (and the definitions of the
  invariants need to be enhanced to take them into account).
\end{remark}

\begin{remark}
  In the proof of Proposition~\ref{prop:compactness} we observed that
  holomorphic curves at $e\infty$ can be viewed as ghost components of
  a curve in $\Sigma_{\overline e}\times[0,1]\times\RR$. Since
  Definition~\ref{def:ghost} is about curves in
  $\Sigma\times[0,1]\times\RR$, it does not include such ghosts.
\end{remark}
\subsection{More on strong boundary monotonicity}
A key tool in eliminating more nodal degenerations  will be strong
boundary monotonicity, to which we now turn.
The definition of strong boundary monotonicity in
Section~\ref{sec:curves-in-sigma} is a  condition on curves; 
but
we will show in Lemma~\ref{lem:MonotoneChords} that the condition in fact
depends only on the asymptotics of the curve.
We will also see, in Lemma~\ref{lem:MonotonicityClosed}, that the strong boundary monotonicity condition is closed. We start with a lemma about the boundaries of our holomorphic curves.

\begin{lemma}\label{lemma:arc_monotonicity}Let
  $u\in\cM(\x,\y\semico\Source\semico P)$. Let $a$ be an arc in $\bdy S$ mapped by
  $u$ to $\alpha_i\times\{1\}\times\RR$ for some $i$. Then $t\circ
  u|_a\co a\to\RR$ is either strictly monotone or constant.
\end{lemma}
\begin{proof}
The map $\pi_\DD\circ u$ is holomorphic
and has image
contained in $[0,1]\times\RR$, and its restriction to the arc $a$ is contained in $\{1\}\times \RR$.
Thus, if $t\circ u|_a$ is non-constant then it cannot have any critical points on the arc.
\end{proof}

\begin{definition}\label{def:strong-monotonicity-P}
  \index{boundary monotonicity!strong}\index{strong boundary monotonicity|see{boundary monotonicity, strong}}%
  Let $\SetS$ be a $k$-element multi-set
  of elements of $[2k]$ (i.e., a formal linear combination of
  elements of $[2k]$ with coefficients in $\NN\cup\{0\}$, so that the
  sum of the coefficients is $k$), and 
  $\gls*{ReebChordsSeqSets}=
  (\rhos_1,\dots,\rhos_n)$ be a sequence of non-empty multi-sets of Reeb
  chords. Define
  \[
  \gls*{osrhos}
  \coloneqq (\SetS \cup M(\vec\rhos^+)) \setminus M(\vec\rhos^-)
  \]
  where $\vec{\rhos}^+\!\coloneqq \bigcup_i \rhos_i^+$
  and the union and difference are interpreted in terms
  of multi-sets, i.e., as sums and differences of formal linear
  combinations over $\ZZ$.
  Let $\gls*{ReebChordsSeqSetsInt}$ 
  be the subsequence
  $(\rhos_i,\dots,\rhos_j)$ of~$\vec\rhos$.  (If $i > j$ the subsequence
  is empty.)
  We say that the pair $(\SetS,\vec{\rhos})$
  is \emph{strongly boundary monotone}
\index{boundary monotonicity!strong!for pair $(\SetS,\vec{\rhos})$}%
  if the
  following three conditions are satisfied:
  \begin{enumerate}[label=(SB-\arabic*),ref=SB-\arabic*]
\item\label{SB1}
\index{(SB-1), (SB-2), (SB-3)}%
  The multi-set $\SetS$ is actually a subset 
  of $[2k]$ (i.e., $\SetS$ has no repeated entries),
\item\label{SB2}
  $M(\rhos_{i+1}^-) \subset \gls*{osrhosint}$
  and
\item\label{SB3}
  $M(\rhos_{i+1}^+)$ is disjoint from $o(\SetS,\vec{\rhos}_{[1,i]})
  \setminus M(\rhos_{i+1}^-)$.
\end{enumerate}
This is equivalent to the following conditions:
\begin{enumerate}[label=(SB-\arabic*${}^\prime$),ref=N-\arabic*${}^\prime$]
\item\label{SB1p}
  For each $i=0,\dots,n$, the multi-set $o(\SetS,\vec\rhos_{[1,i]})$
  is a $k$-element
  subset of $[2k]$ with no repeated elements.
\index{(SB-$1'$), (SB-$2'$)}%
\item
  For each $i$, $M(\rhos_i^-)$ and $M(\rhos_i^+)$ (as multi-sets)
  have no elements with multiplicity bigger than~$1$.
\end{enumerate}
We also extend the definition to pairs  $(\x,\vec{\rhos})$, where $\x$
is a generalized generator, which we call
strongly boundary monotone if
$(o(\x),\vec{\rhos})$ is strongly boundary monotone.
\index{boundary monotonicity!strong!for pair $(\x,\vec{\rhos})$}%
\end{definition}

\begin{lemma}
  \label{lem:MonotoneChords}
  If $(\x,\vec{\rhos})$ is strongly boundary
  monotone and $u\in\cM^B(\x,\y\semico\Source\semico\vec{P})$ with
  $[\vec{P}] =\vec{\rhos}$, then
  $u$ is strongly boundary monotone. Conversely, if
  $\cM^B(\x,\y;\Source;\vec{P})$ contains a strongly boundary
  monotone holomorphic curve $u$ then $(\x,\vec{\rhos})$ is strongly
  boundary monotone.
\end{lemma}
(Recall that $[\vec{P}]$ is the sequence of multi-sets of Reeb chords
gotten by replacing the punctures in $\vec{P}$ by the Reeb chords
labeling them, and that given a term $P_i$ in $\vec{P}$, $[P_i]$ is
the corresponding multi-set of Reeb chords.)
\begin{proof}
  The fact that $(\x,\vec{\rhos})$ being strongly boundary monotone implies
  $u$ is strongly boundary monotone is immediate from
  Lemma~\ref{lemma:arc_monotonicity}.

  To see that $u\co S\to \Sigma\times[0,1]\times\RR$ strongly boundary
  monotone implies $(\x,\vec{\rhos})$ strongly boundary monotone, we also use
  the fact that there are no boundary degenerations, as
  follows. Condition~(\ref{SB1})
  is immediate from the fact that $\x$ is
  a generator. To prove Condition~(\ref{SB2}),
suppose that at some stage $M(\rhos_{i+1}^-)\not\in
  o(\SetS,\vec{\rhos}_{[1,i]})$, so there is some chord
  $\rho_j\in\rhos_{i+1}$ with $\rho_j^-\not\in
  o(\SetS,\vec{\rhos}_{[1,i]})$. Then there must also be a chord
  $\rho_\ell\in\rhos_{i+1}$ with $M(\rho_\ell^+)=M(\rho_j^-)$, and an arc
  $A$ on $\bdy S$ so that $u$ is asymptotic to $\rho_\ell$ at one end of
  $A$, $\rho_j$ at the other end of $A$, and $t\circ u|_A$ is
  constant. But by the open mapping principle applied to $\pi_\DD\circ
  u$, this implies that the component of $S$ containing $A$ is a
  boundary degeneration, violating Lemma~\ref{lem:NoBoundaryDegenerations}. 
  A similar argument implies
  Condition~(\ref{SB3}).
\end{proof}

In light of Lemma~\ref{lem:MonotoneChords}, we will say that a moduli space
$\cM^B(\x,\y;\Source;\vec{P})$ is strongly boundary monotone
if $(\x,[\vec{P}])$ is strongly boundary monotone.
\index{boundary monotonicity!strong!for moduli space}%

We next show that 
strong boundary monotonicity is a closed condition on the compactified moduli space. To make sense of this, we need to define strong boundary monotonicity for holomorphic combs. The following, somewhat weak definition will be sufficient for our purposes:
\begin{definition}
  A holomorphic comb $U$ is \emph{strongly boundary monotone} if, after deleting any ghost components, every story of the spine of $U$ is strongly boundary monotone.
\index{boundary monotonicity!strong!for comb}%
\end{definition}
\begin{lemma}\label{lem:MonotonicityClosed} The strong boundary monotonicity condition is closed. That is, if $\{U_j\}$ is a convergent sequence of strongly boundary monotone holomorphic combs then the limit $U$ of $\{U_j\}$ is strongly boundary monotone.
\end{lemma}
\begin{proof}
  This follows from the monotonicity of $t\circ u$ on arcs (Lemma~\ref{lemma:arc_monotonicity}) and the non-existence of boundary degenerations (Lemma~\ref{lem:NoBoundaryDegenerations}).
\end{proof}

\subsection{No nodal curves}
In this section, we rule out boundary double points and haunted curves
for moduli spaces satisfying the strong boundary monotonicity
condition. We start with boundary double points.

\begin{lemma}\label{lemma:NoBoundaryDoublePoints}
  If $(\x,\vec{\rhos})$ is strongly boundary monotone, then
  the moduli spaces $\ocM^B(\x,\y\semico\Source\semico\vec{P})$ with $[\vec{P}] =
  \vec{\rhos}$ do not contain any curves with boundary double points.
  \index{degeneration of holomorphic curves!boundary double point!does not occur}%
\end{lemma}
\begin{proof}
  Suppose that $u_j\in\cM^B(\x,\y\semico\Source\semico\vec{P})$
  converges to a holomorphic comb $u$ with source $S$ containing a
  boundary double point.  We have a node $q$ in the source $S$
  whose pre-images $q_1$ and $q_2$ in the normalization ${\widetilde{S}}$ of $S$ both
  get mapped non-trivially to $[0,1]\times\RR$. Since $q$ is on the boundary
  of $S$, it follows that $q$ is mapped to some $\alpha$- or $\beta$-curve.
  Now, since the projection $t\circ {\widetilde u}$ is a non-trivial harmonic function
  near both $q_1$ and $q_2$, by the open mapping principle, we can
  find disjoint neighborhoods $N_i$ of
  $q_i\in \partial{\widetilde S}$ (for $i=1,2$), with the property
  that both $t\circ {\widetilde u}(N_1)$ and
  $t\circ {\widetilde u}(N_2)$ are the same open interval $I$ around
  $t\circ u(q)$. 
  It follows that, for large enough $j$, $u_j$ is not strongly
  boundary monotone.
\end{proof}

Next we show that ghosts are vanishingly rare.
\index{haunted curve!rare}%
\begin{lemma}\label{lemma:no-ghosts}
  If $(\x,[\vec{P}])$ is strongly boundary monotone and
  $\ind(B,\Source,\vec{P})\leq 2,$ then $\ocM^B(\x,\y;\Source,\vec{P})$
  does not contain any haunted curves.  
\end{lemma}
\begin{proof}
  This follows from Lemma~\ref{lemma:NoBoundaryDoublePoints} and a
  standard argument, a version of
  which we give.  Suppose that arcs and circles $A_1,\dots,A_n$ in $S$
  collapse to form some ghost components. Let $S_1,\dots,S_m$ be the
  ghost components of the limit and $S_0$ the union of the other
  components of the limit, $\SourceSub{o}$ decorated source
  corresponding to $S_0$ and $u\co S_0\cup\dots\cup
  S_m\to\Sigma\times[0,1]\times\RR$ the limiting holomorphic
  curve. The arcs $A_i$ come in two types: arcs which separate a
  ghost-component from a non-ghost component, which we will call
  \emph{important arcs}; and arcs separating two ghost components, a
  ghost component from itself, two non-ghost components, or a
  non-ghost component from itself.

  Let $G$ be the \emph{ghost graph} with one vertex for each ghost
  component and one edge for each node connecting two ghost
  components. Let $G_1,\dots, G_p$ be the connected components of $G$
  and $T_i$, $i=1,\dots, p$, the union of the ghost components
  corresponding to $G_i$. We claim that each $T_i$ has exactly one
  node corresponding to an important arc. Certainly $T_i$ has at least
  one node corresponding to an important arc, and if some $T_i$ has
  more than one node corresponding to an important arc then $u|_{S_0}$
  would have a double point, violating
  Lemma~\ref{lemma:NoBoundaryDoublePoints}. This argument also implies
  that smoothing the nodes in $T_i$ gives a surface with a single
  boundary component.

  Each $S_i$ must be stable; in particular, no $S_i$ is a disk or
  sphere with $2$ or fewer nodes. Since exactly one node of $T_i$
  corresponds to an important arc, it follows that for each
  $i=1,\dots,p$, either $G_i$ is not a tree or some component of $T_i$
  is not a disk. Since the smoothing of $T_i$ is a surface with
  connected boundary, it follows that
  $\chi(S_0)\geq \chi(S)+2$.  Thus,
  \begin{align*}
    \ind(B,\SourceSub{0},\vec{P})&=g-\chi(S_0)+2e(B)+|\vec{P}|\\
    &\leq
    g-\chi(S)+2e(B)+|\vec{P}|-2\\
    &=\ind(B,\Source,\vec{P})-2\\
    &\leq 0.
  \end{align*}
  Hence, $u|_{S_0}$ belongs to a negative-dimensional moduli space,
  and so does not exist.
\end{proof}

\subsection{The codimension-one boundary of the moduli space of curves with a given source}\label{sec:codim-one-bound}
We have now ruled out bubbling and nodal degenerations, so the only
ends of the moduli spaces that remain are height $2$ combs, collisions of levels and degenerating join, split and shuffle curves at $e\infty$. 
At this point, we would like to assert that if strong boundary
monotonicity is satisfied then $\ocM^B(\x,\y\semico\Source\semico\vec{P})$ is a compact
$1$-manifold with boundary, and $\bdy\ocM^B(\x,\y\semico\Source\semico\vec{P})$ consists of
exactly those pieces. The
difficulty with this statement is the case of shuffle curves and
split curves with
many split components:
in these cases, the transversality needed by
Proposition~\ref{prop:gluing_simple_comb} is absent. One
would like to assert, for
instance, that simple combs $(u',v')$ in
$\ocM^B(\x,\y;\Source;\vec{P})\setminus\cM^B(\x,\y;\Source;\vec{P})$ are isolated. The
  maps $u'$ are, indeed, isolated, but this is not obvious for the components $v'$ at
  east $\infty$. (This is related to the fact that the evaluation map from
  $\ocM^B(\x,\y\semico\Source\semico\vec{P})$ collapses strata coming from non-trivial moduli at east
  $\infty$.)

  However, all we will need later is that the sum of the number of elements of the
  moduli spaces occurring as the boundaries of $\ocM^B(\x,\y\semico\Source\semico\vec{P})$ (when
  the latter space is one-dimensional) is zero modulo two. This weaker
  statement is
  Theorem~\ref{thm:master_equation}, below.

\begin{figure}
  \includegraphics[scale=.83333]{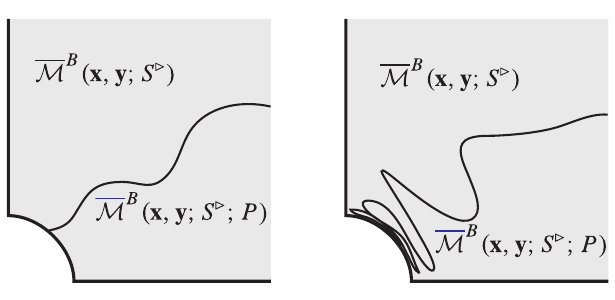}
  \caption[Expected and unexpected behavior of
    $\ocM(\x,\y\semico\Source\semico\vec{P})\subset \ocM(\x,\y;\Source)$]{\textbf{A picture of the expected behavior of
    $\ocM(\x,\y\semico\Source\semico\vec{P})\subset \ocM(\x,\y;\Source)$ (left)
    and a behavior of $\ocM(\x,\y\semico\Source\semico\vec{P})\subset
    \ocM(\x,\y;\Source)$ which we have not ruled out (right).} The
    reader might compare with
    Figure~\ref{fig:moduli_proper_subset_2}.}
\label{fig:moduli_bad_behavior}
\end{figure}
\colorused

\begin{figure}
  \includegraphics[scale=.4167]{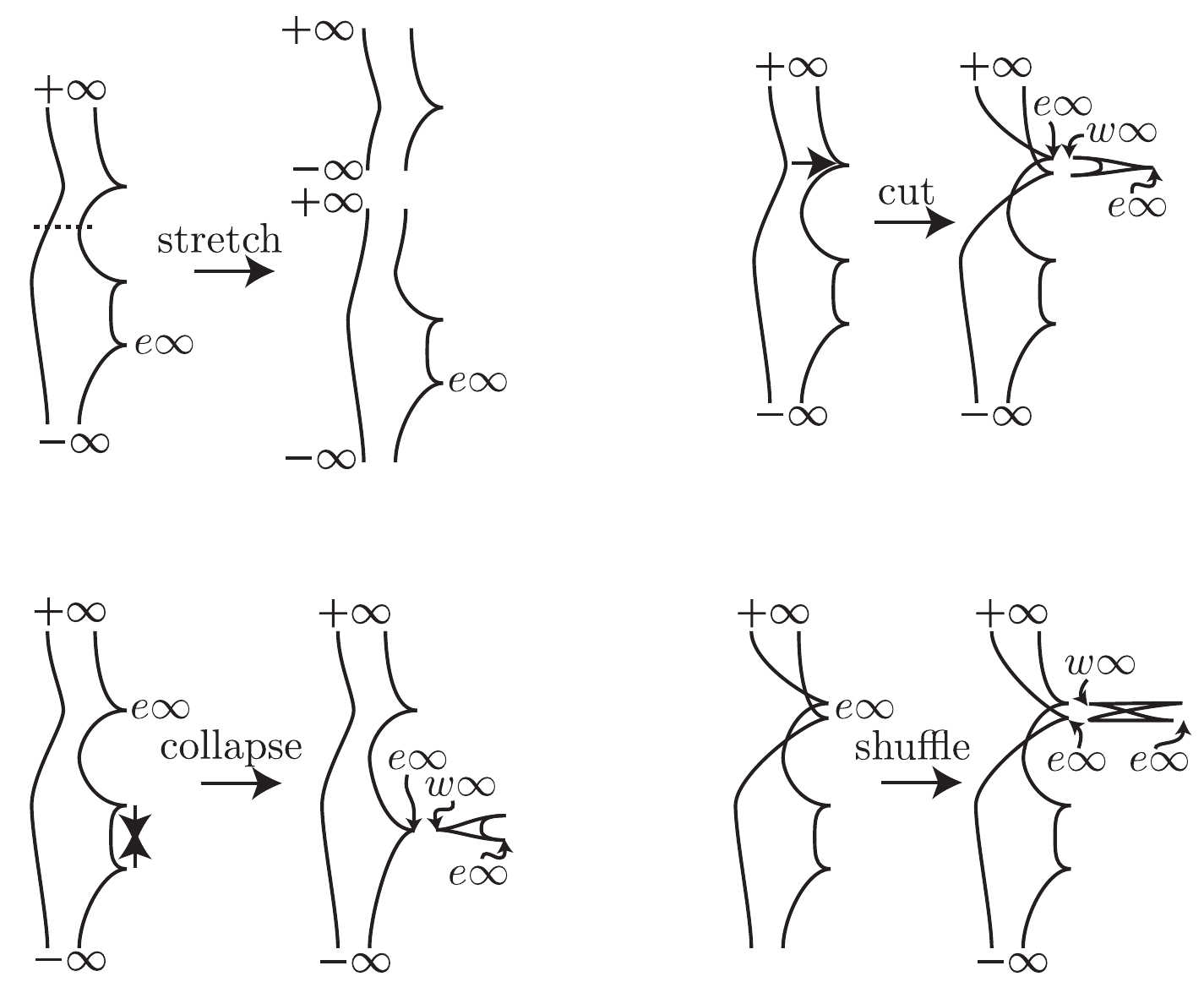}
  \caption[Degenerations
    discussed in Theorem~\ref{thm:master_equation}]{\textbf{Schematics and examples of the four kinds of degenerations
    of Theorem~\ref{thm:master_equation}.} Top left: degenerating
    into a height $2$ holomorphic comb. Top right: degenerating a
    join curve. Bottom left: degenerating a split curve. Bottom right:
    degenerating a shuffle curve.}
  \label{fig:schematics-degeneration}
\end{figure}

  Before stating this theorem, we define formally various moduli spaces that can
  appear at the end of one-dimensional spaces, as illustrated in
  Figure~\ref{fig:schematics-degeneration}, and place one more restriction on collisions of levels.  
\index{degeneration of holomorphic curves|seealso{end}}%
\index{end|seealso{degeneration}}%
\begin{definition}
  \label{def:ends-moduli}
  Fix a one dimensional moduli space $\cM^B(\x,\y\semico\Source\semico\vec{P})$
  (or $\cM$ for short)
  satisfying strong boundary monotonicity.

  \index{two-story end|see{end, two-story}}%
  \index{end!two-story}%
  A \emph{two-story end of $\cM$} is an element of
  $\cM^{B_1}(\x,\w\semico\SourceSub{1}\semico\vec{P}_1)\times
  \cM^{B_2}(\w,\y\semico\SourceSub{2}\semico\vec{P}_2)$, where
  $\w\in\S(\HD)$, $B_1 \in \pi_2(\x,\w)$ and $B_2 \in \pi_2(\w,\y)$
  are such that $B = B_1 * B_2$, and $\Source = \SourceSub{1} \glue
  \SourceSub{2}$ is a way of splitting $\Source$ in two which divides
  the ordered partition $\vec{P}$ as $\vec{P}_1 < \vec{P}_2$.

  \index{join!curve end|see{end, join curve}}%
  \index{end!join curve}%
  A \emph{join curve end of $\cM$ at level~$i$} is an element of
  $\cM^B(\x,\y\semico\Source'\semico\vec{P}')$, where $\Source'$ and~$\vec{P}'$
  are obtained in the following way.  Pick an east puncture~$q$
  of~$\Source$ in the $i\th$ part of~$\vec P$
  and a decomposition $\rho_q = \rho_a \uplus \rho_b$. Then
  the decorated source $\Source'$ is any source with a pair of
  punctures $a$ and $b$ labeled by $\rho_a$ and $\rho_b$, and such
  that $\Source$ is obtained from $\Source'$ by pregluing a
  join component to $\Source'$ at the punctures $a$ and~$b$. The
  partition $\vec{P}'$ is obtained from $\vec{P}$ by replacing $q$
  with $\{a,b\}$ in the $i\th$ part.

  \index{shuffle!curve end|see{end, shuffle curve}}%
  \index{end!shuffle curve}%
  An \emph{odd} (respectively \emph{even})
  \emph{shuffle curve end of $\cM$ at level~$i$} is an element of
  $\cM^B(\x,\y\semico\Source'\semico\vec{P}')$, where $\Source'$ and $\vec{P}'$
  are obtained in the following way.  Pick east punctures~$q_1$ and
  $q_2$ of~$\Source$ contained in the $i\th$ part of $\vec{P}$, and so
  that the corresponding Reeb chords are nested (respectively
  interleaved).  Then $\Source'$ is any source together with punctures
  $q_1'$ and $q_2'$ such that $\Source$ is obtained from $\Source'$ by
  pregluing an odd (respectively even) shuffle component $\biSource{}'$
  to $\Source'$
  at the punctures $q_1'$ and~$q_2'$ with $q_1$ and~$q_2$ labeling the
  punctures of $\biSource{}'$. The partition~$\vec{P}'$ is
  obtained from $\vec{P}$ by replacing $\{q_1,q_2\}$ with
  $\{q_1',q_2'\}$ in the $i\th$ part.

  \index{end!collision of levels}%
  \index{collision of levels|see{end, collision of levels}}%
  A \emph{collision of levels $i$ and $i+1$ in $\cM$} is an element of
  $\cM^B(\x,\y\semico\Source'\semico(P_1,\allowbreak\dots,\allowbreak
  P_i \uplus P_{i+1}, \dots, P_n))$
  with $i$ satisfying $1 \le i < n$ and $\Source'$ obtained from
  $\Source$ by contracting arcs on $\bdy\Source$ connecting punctures
  labeled by abutting pairs of Reeb chords from $P_i$ and
  $P_{i+1}$ (and replacing those pairs of Reeb chords $\rho_i,\rho_j$
  by their joins $\rho_i\uplus\rho_j$), and $P_i\uplus
  P_{i+1}$ denotes the Reeb chords in
  $\Source'$ coming from $P_i$ and $P_{i+1}$.  
  (The notions of abutting and $\uplus$ are given in
  Definition~\ref{def:Composable-one}.)

  We extend this definition inductively to a collision of levels
  $i,\dots,i+j$: a collision of levels $i,\dots,i+j$ is a collision of
  levels $i$ and $i+1$ in the result of a collision of levels $i+1,\dots,i+j$.
\end{definition}

We will see (Lemma~\ref{lemma:collision-weakly-composable}) that homological linear independence of the
$\alpha$-curves often prevents collisions of levels. First, a definition:
\begin{definition}\label{def:weakly-composable}
  \index{composable!sets of Reeb chords!weakly}%
  \index{weakly composable sets of chords}%
  We say that two sets of
  Reeb chords $\rhos_i$ and
  $\rhos_{i+1}$ are \emph{weakly composable} if, for all $\rho_j\in \rhos_i$
  and $\rho_k\in \rhos_{i+1}$, if $M(\rho_j^+)=M(\rho_k^-)$ then
  $\rho_j^+=\rho_k^-$. More generally, a sequence $\vec\rhos=(\rhos_1,\dots,\rhos_n)$ is said to be
  weakly composable if for all $i=1,\dots, n-1$, the sets
  $\rhos_1\uplus \dots\uplus \rhos_i$ and $\rhos_{i+1}$ are weakly composable.
\end{definition}

Note that if two sets of Reeb chords are composable in the sense of
Definition~\ref{def:Composable-sets} then they are weakly composable, but
the converse is false.

\begin{lemma}\label{lemma:collision-weakly-composable} Suppose that
  $\cM^B(\x,\y;\Source;\vec{P})$ is strongly boundary monotone and
  $\bdy\ocM^B(\x,\y\semico\Source\semico\vec{P})$ contains a point
  corresponding to a collision of levels
  $i$ and $i+1$ in $\vec{P}$. Then $[P_i]$ and $[P_{i+1}]$ are weakly
  composable.

  More generally, suppose that $\bdy\ocM^B(\x,\y;\Source;\vec{P})$
  contains a point corresponding to a collision of levels
  $i,i+1,\dots,j$ in $\vec{P}$. Then $([P_i],[P_{i+1}],\dots,[P_j])$
  is weakly composable. 
\end{lemma}
\begin{proof}
  Let $\{u_j\}$ be a sequence of curves in
  $\cM^B(\x,\y;\Source;\vec{P})$ whose spine converges to a collision
  of levels.  Fix chords
  $\rho\in [P_i]$ and $\sigma\in [P_{i+1}]$ with the property that
  $M(\rho^+)=M(\sigma^-)$.  Our goal is to show that
  $\rho^+=\sigma^-$.

  Consider the arc $A$ on $\partial S$ along which $u_j$ leaves $\rho$. The image
  $\pi_\Sigma\circ u_j(A)$ is contained in some curve
  $\alpha_i$. Similarly, for the arc $B$ on $\partial\Source$ along which $u_j$
  enters $\sigma$, $\pi_\Sigma\circ u_j(B)$ is contained in the same
  curve $\alpha_i$. By hypothesis, the sequences
  $t(u_j(\rho))$ and $t(u_j(\sigma))$ converge to the same point.

  We must have $A = B$ by strong boundary
  monotonicity on~$u_j$.  Therefore, by Lemma~\ref{lemma:arc_monotonicity}, the restriction of $t\circ u$ to $A$ is
  constant, and hence $u$ must be constant on the
  entire component of $S$ containing $A$.
  Lemmas~\ref{lem:NoBoundaryDegenerations} and~\ref{lemma:no-ghosts}
  imply that this component must be mapped to east $\infty$. It
  follows that the terminal point of $\rho$ coincides with the initial
  point of $\sigma$, as desired.

  The more general case is similar, and is left to the reader.
\end{proof}

We are now ready to state the theorem we will use to
prove all the algebraic properties of our modules.

\begin{theorem}\label{thm:master_equation}Suppose that
  $(\x,\vec{\rhos})$ is
  strongly boundary monotone. Fix $\y$, $B\in\pi_2(\x,\y)$,
  $\Source$, and $\vec{P}$ so that $[\vec{P}] = \vec{\rhos}$ and
  $\ind(B,\Source,P)=2$. Let $\cM = \cM^B(\x,\y\semico\Source\semico\vec{P})$.
  Then the
  total number of
  \begin{enumerate}[label=(ME-\arabic*),ref=ME-\arabic*]
  \item \index{(ME-1)--(ME-4)} \label{item:master:two-story} two-story ends of~$\cM$,
  \item \label{item:master:join} join curve ends of~$\cM$,
  \item \label{item:master:odd-shuffle} odd shuffle curve ends of~$\cM$ and
  \item \label{item:Collision} collision of levels $i$ and $i+1$ in $\cM$
  \end{enumerate}
  is even.
  Moreover, in Case~\eqref{item:Collision} the parts $P_i$ and $P_{i+1}$
	are weakly composable.
\end{theorem}
\begin{proof}
Let $U_{<\epsilon}$ denote the open subset of
$\ocM^B(\x,\y\semico\Source\semico\vec{P})$ where there are two parts of $\vec{P}$ with
height difference $<\epsilon$. Let $U_{\textrm{shuf}}$ denote the
union of a smeared neighborhood of each shuffle curve end
in $\ocM^B(\x,\y\semico\Source\semico\vec{P})$. It follows from the results above
that
\[
\ocM^B_{\textrm{cropped}}=\ocM^B(\x,\y\semico\Source\semico\vec{P})\setminus (U_{<\epsilon}\cup
U_{\textrm{shuf}})
\]
is a compact $1$-manifold with boundary
(smoothness following from Proposition~\ref{prop:transversality},
compactness from Proposition~\ref{prop:compactness}), the boundary of which
consists of the following pieces:
  \begin{itemize}
    \item Two-story holomorphic combs $(u_1,u_2)$. For $\epsilon$ and
      $U_{\textrm{shuf}}$ small
      enough, the number of these correspond exactly to
      the number of two-story ends of~$\cM$. (This is ensured
	by Proposition~\ref{prop:gluing_two_story}. The fact that splittings must occur at generators, not just generalized generators, follows from strong boundary monotonicity.)
    \item Simple holomorphic combs $(u,v)$ with $v$ a join curve. For
      $\epsilon$ and
      $U_{\textrm{shuf}}$ small enough the number of these
      boundary components is the number of join curve ends
      of $\cM$. (This is an application of Proposition~\ref{prop:gluing_simple_comb},
	with transversality hypotheses ensured by Propositions~\ref{prop:transversality}
	and~\ref{prop:east_transversality}.)
    \item An even number of points for each even
      shuffle curve end, and an odd number of points for each
      odd shuffle curve end, if $\epsilon$ and $U_{\textrm{shuf}}$ are
      small enough. (This is an application of 
	Proposition~\ref{prop:gluing-shuffle}.)
    \item The subspace of $\cM^B(\x,\y;\vec{P})$ where there are two
      parts which differ in
      height by $\epsilon$. By
      Propositions~\ref{prop:gluing_degree_one} (if at least one
      split component degenerates) and~\ref{prop:transversality} (if no
      split components degenerate), for $\epsilon$ small enough the number
      of boundary points of this form agrees, modulo~$2$, with the
      number of collisions of two levels from $\Source$. The fact that
      the levels must be weakly composable follows from Lemma~\ref{lemma:collision-weakly-composable}.
  \end{itemize}
  (In the first case, we take $\epsilon$ and $U_{\textrm{shuf}}$ small
  enough that no two-story combs occur in
  $\ocM^B\setminus\ocM^B_{\textrm{cropped}}$, and in the second case
  we take them small enough that no join curve degenerations occur in
  $\ocM^B\setminus\ocM^B_{\textrm{cropped}}$.  In the third case,
  $\epsilon$ should be small enough that $U_{<\epsilon}$ is disjoint
  from $U_{\textrm{shuf}}$. In the fourth case, $\epsilon$ should be
  small enough for Proposition~\ref{prop:gluing_degree_one} to apply.)

  We have accounted for all the boundary components, in view of
  Proposition~\ref{prop:restrict_degens_1} combined with
  Lemma~\ref{lemma:NoBoundaryDoublePoints}.
\end{proof}

\begin{remark}
  We will see in Lemma~\ref{lemma:collision-is-composable}
  that for moduli spaces of embedded curves, the only collisions which
  occur are ones in which the colliding parts are composable (as
  defined in Section~\ref{sec:reeb-chords-def}) and not merely weakly
  composable, as well as similar restrictions on the join curve and
  shuffle curve ends.
\end{remark}

\section{More on expected dimensions}\label{sec:expected-dimensions}

Recall from Proposition~\ref{Prop:Index}
that the expected dimension of $\cM^B(\x,\allowbreak\y\semico\Source\semico\vec{P})$ depends on the topology of the source curve. 
Our aim here is to study an index formula which depends on the source
curve only through its homology class~$B$ and its asymptotics~$[\vec{P}]$. We
establish two key properties of this index:
\begin{itemize}
\item In cases where we have an embedded curve representing the
  homology class, this formula agrees with our earlier index formula
  (Proposition~\ref{prop:asympt_gives_chi}). 
\item The index is additive under juxtapositions (Proposition~\ref{prop:indAdditive}). 
\end{itemize}
A fundamental consequence of the index formula is the following:
\begin{itemize}
\item For collisions of levels occurring in the boundaries of
  two-dimensional, strongly boundary monotone moduli spaces,
  the colliding levels are composable, rather than just weakly composable
(Lemma~\ref{lemma:collision-is-composable}).
\end{itemize}
This last property is the justification
for setting double crossings of strands to zero in the algebra
associated to a surface. There are also restrictions on join curve and odd shuffle curve ends,
proved in Section~\ref{sec:embedded-degen}.

\subsection{The index at an embedded curve}\label{sec:ind-at-emb}
\index{expected dimension!at embedded curve|(}
Before stating the index formula,
we recall and extend some definitions  from
Section~\ref{sec:pregrading}.  For $\alpha_1, \alpha_2 \in H_1(Z',\CircPts)$,
recall that $L(\alpha_1,\alpha_2) = m(\alpha_2,\bdy\alpha_1)$ is the
``linking'' of $\alpha_1$ and~$\alpha_2$.  Concretely, if the
$\alpha_i$ are represented by single Reeb chords~$\rho_i$, we have
\begin{equation}\label{eq:concrete-linking}
\glsadd{Linkingalpha}
L([\rho_1],[\rho_2])=
\begin{cases}
  1/2 & \textrm{if } \rho_1^+=\rho_2^-\textrm{ or } \rho_1^-\lessdot \rho_2^-\lessdot \rho_1^+=\rho_2^+ \textrm{ or } \rho_1^- = \rho_2^-\lessdot \rho_1^+\lessdot \rho_2^+\\
  -1/2 & \textrm{if } \rho_2^+=\rho_1^-\textrm{ or } \rho_2^-\lessdot \rho_1^-\lessdot \rho_1^+=\rho_2^+ \textrm{ or } \rho_1^- = \rho_2^-\lessdot \rho_2^+\lessdot \rho_1^+\\
  0 & \textrm{if $\rho_1\cap\rho_2=\emptyset$ or
  $\rho_1\subset\Int(\rho_2)$ or $\rho_2\subset\Int(\rho_1)$ or $\rho_1 = \rho_2$}\\
  1 & \textrm{if } \rho_1^-\lessdot
  \rho_2^-\lessdot \rho_1^+\lessdot \rho_2^+\\
  -1 & \textrm{if } \rho_2^-\lessdot
  \rho_1^-\lessdot \rho_2^+\lessdot \rho_1^+.
\end{cases}
\end{equation}
(The notation $\lessdot$ is defined in
Section~\ref{sec:matched-circles}.)
Recall that for $a\in\Alg(n,k)$ with starting idempotent~$S$, the quantity
$\iota(a) = \inv(a) - m([a],S)$ is the Maslov component of the grading
of~$a$.  We can again express this more concretely for in terms of Reeb chords.
\begin{lemma}\label{lem:iota-chords}
  For $\rhos$ a set of Reeb chords so that $a(\rhos) \ne 0$,
  \[
  \glsadd{iota}
  \iota(a(\rhos)) = -\!\!\sum_{\{\rho_1, \rho_2\}\subset \rhos}\!\!
      \bigl\lvert L([\rho_1],[\rho_2])\bigr\rvert - \frac{\abs{\rhos}}{2}.
  \]
\end{lemma}
\begin{proof}
  Recall from the proof of Proposition~\ref{prop:grading-descends} that
  $\iota(a)$ is unchanged under adding horizontal strands; thus we may
  pretend that the only strands in $a(\rhos)$ are those coming from the
  Reeb chords in~$\rhos$. Both terms in $\iota(a(\rhos))$, namely
  $\inv(a(\rhos))$ and
  $m([\rhos], \rhos^-)$, can be written as sums over pairs of Reeb chords or
  single chords.  A single chord contributes $-1/2$ to $m([\rhos],\rhos^-)$,
  while a straightforward case analysis shows that the contribution to
  the sum from a pair of chords is $-1$ if they are interleaved,
  $-1/2$ if the endpoints abut, and $0$ otherwise, as in the
  statement.
\end{proof}

We will shorten notation, writing $\gls*{Linkingrho}$ (resp. $\gls*{Linkingrhos}$)
for
$L([\rho_1],[\rho_2])$ (resp. $L([\rhos_1],[\rhos_2])$) and 
$\gls*{iotarhos}$
for $\iota(a(\rhos))$.  Also, for
a strongly boundary monotone
sequence of sets of Reeb chords
$\vec{\rhos}=(\rhos_1,\dots,\penalty 500\rhos_\ell)$, define
$\iota(\vec{\rhos})$ by
\begin{equation}\label{eq:def-iota}
\gls*{iotavecrhos}
\coloneqq\sum_{i}\iota(\rhos_i)+\sum_{i<j}
  L(\rhos_i,\rhos_j).
\end{equation}

\begin{lemma}\label{lem:G-n-products}
  The product of a sequence of elements in $\bigGroup(n)$ is given
  by
  \[
  (k_1,\alpha_1)\cdots(k_n,\alpha_n) = \Bigl(\sum_i k_i + \sum_{i < j}
  L(\alpha_i,\alpha_j),\sum_i \alpha_i\Bigr).
  \]
\end{lemma}

\begin{proof}
  As in the proof of Proposition~\ref{prop:G-n-group}, expand
  the right-associated product on the left-hand side of the equation
  and use bilinearity of~$L$.
\end{proof}

\begin{lemma}\label{lem:iota-grading}
  For $\vec\rhos$ a strongly boundary monotone sequence of Reeb
  chords, $\iota(\vec{\rhos})$ is the Maslov component of
  $\grb(a(\rhos_1))\cdots\grb(a(\rhos_\ell))$.
\end{lemma}
\begin{proof}
  Clear from the definitions and Lemma~\ref{lem:G-n-products}.
\end{proof}

For $x \in \alphas \cap \betas$ and $B$ a domain, define
$\gls*{nxB}$
to be the average of the local multiplicities of $B$ in the
four regions surrounding~$x$, and for $\x \in \S(\HD)$, define
$
\gls*{nxboldB}
\coloneqq \sum_{x \in \x}n_x(B)$.

\begin{definition}\label{def:emb-ind-emb-chi}
We say that a pair $(B,\vec{\rhos})$ with $B\in\pi_2(\x,\y)$ and
$\vec{\rhos}$ a
sequence of non-empty sets of Reeb chords is \emph{compatible}
\index{compatible pair!$(B,\vec{\rhos})$}%
if the homology
classes on the boundary agree, $[\vec \rhos]
= \bdy^\bdy B$, and $(\x,\vec{\rhos})$
is strongly boundary monotone.  For such a pair, define the \emph{embedded
  Euler characteristic}, \emph{embedded index}, and \emph{embedded
  moduli space} by
\index{embedded!Euler characteristic}%
\index{embedded!index|see{index, embedded}}%
\index{index|seealso{expected dimension}}%
\index{index!embedded}%
\index{embedded!moduli space}%
\index{moduli space!embedded}%
\begin{align*}
\gls*{chiemb}
&\coloneqq
g+e(B)-n_\x(B)-n_\y(B)-\iota(\vec{\rhos})\\
\gls*{indemb}
&\coloneqq e(B)+n_\x(B)+n_\y(B)+|\vec\rhos|+\iota(\vec{\rhos})\\
\gls*{ModSpaceEmbeddedOpen}
&\coloneqq
    \!\!\bigcup_
       {\substack{\chi(\Source)=\chi_{\emb}(B,\vec{\rhos})\\
           [\vec{P}] = \vec{\rhos}}}\!\!
    \Mod^B(\x, \y\semico\Source\semico\vec{P}).
\end{align*}
Define $\gls*{ModSpaceEmbeddedCpct}$
similarly.
\end{definition}

Justification for this terminology is given by the following propositions.

\begin{figure}
\includegraphics[scale=.8333]{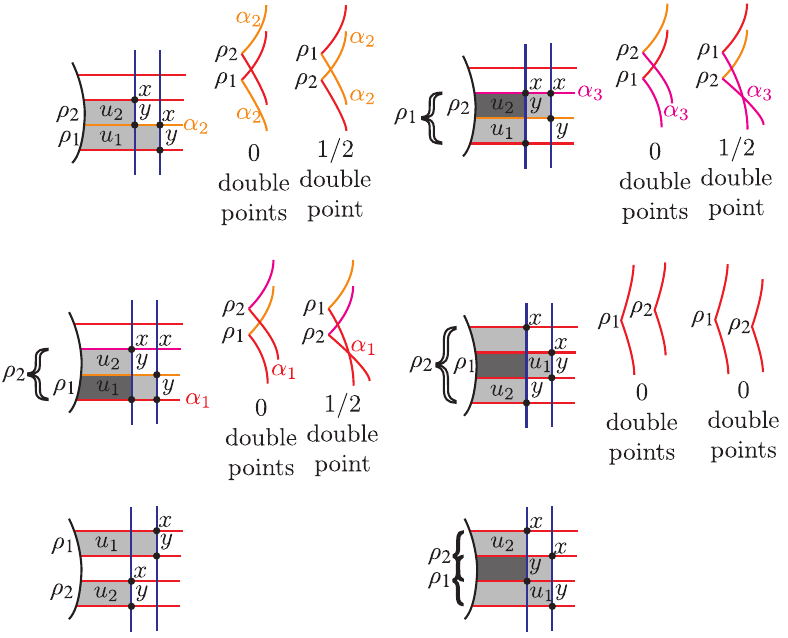}
\caption[Domains of model curves used to compute $\inv'$]{\textbf{The
    domains of model curves used to compute $\inv'$.} Top left:
  case~(\ref{case:ind-cross:1}). Case~(\ref{case:ind-cross:2}) is
  obtained by switching $\rho_1$ and $\rho_2$. Top right:
  case~(\ref{case:ind-cross:3}). Center left:
  case~(\ref{case:ind-cross:4}).  Center right:
  case~(\ref{case:ind-cross:7}). Bottom left:
  case~(\ref{case:ind-cross:8}). Bottom right:
  case~(\ref{case:ind-cross:9}). In the top four pictures, a schematic
  for $u_1(\partial S_1),u_2(\partial S_2)\subset
  \alphas\times\{1\}\times\RR$ is also shown.}\label{fig:IndexModels}
\end{figure}

\begin{proposition}\label{prop:asympt_gives_chi} Suppose that
  $\cM^B(\x,\y\semico\Source\semico\vec{P})$ admits a holomorphic
  representative~$u$. Then 
  \begin{equation}
    \chi(S)=\chi_\emb(B,[\vec{P}])\label{eq:chiDetermined}
  \end{equation}
  if and only if $u$ is embedded.  When there is such an embedded,
  holomorphic $u$, the expected dimension of
  $\cM^B(\x,\y\semico\Source\semico\vec{P})$ is given by
  \begin{equation}\label{eq:ComboIndex}
  \ind(B,\Source,\vec{P})=\ind(B,[\vec{P}]).
  \end{equation}
  If $\cM^B(\x,\y;\Source;\vec{P})$ is strongly boundary monotone and
  admits a non-embedded holomorphic representative then
  $\ind(B,\Source,\vec{P})\leq\ind(B,[\vec{P}])-2$.
\end{proposition}
\begin{proof}
  We imitate the proof of~\cite[Proposition 4.2]{Lipshitz06:CylindricalHF},
  to which the reader is referred for a more leisurely
  account.  We prove that if $u$ is embedded, the Euler characteristic
  is as stated; the other direction and Formula~\eqref{eq:ComboIndex} are
  immediate consequences; we will return to the inequality for
  non-embedded curves at the end of the proof.

  Fix an embedded holomorphic curve
  $u\in\cM^B(\x,\y\semico\Source\semico\vec{P})$.  Let $\br(u)$ denote the
  ramification number of $\pi_\Sigma\circ u$. (Here, boundary branch
  points contribute $1/2$.) By the Riemann-Hurwitz formula,
  \[
  \chi(S)=e(S)+g/2+\sum_i\abs{P_i}/2=e(B)-\br(u)+g/2+\sum_i\abs{P_i}/2.
  \]
  So, to compute $\chi(S)$ we only need to compute $\br(u)$.

  Let $\tau_R(u)$ denote a copy of $u$ translated by $R$ units in the
  $\RR$-direction.  We will prove the result by comparing
  $u\cdot\tau_R(u)$ for $R$ small and $R$ large.

  Assume now that the partition $\vec{P}$ is discrete, i.e., each
  part $P_i$ consists of a single puncture labeled by~$\rho_i$; we
  will return to the
  general case at the end of the proof.  
  Since $u$ is embedded, for
  small~$\epsilon$ the
  curves $u$ and $\tau_{\epsilon}(u)$ intersect only near where $u$ is
  tangent to $\partial/\partial t$, i.e., at branch points of $u$.
  Since $u$ and $\tau_\epsilon(u)$ are $J$-holomorphic, their
  algebraic intersection number is $\br(u)$.  (Here and later,
  intersections along the boundary count for $1/2$.)

  On the other hand, for $R$ large, from the vantage point of $u$,
  $\tau_R(u)$ looks like $g$ trivial strips $\x\times[0,1]\times\RR$;
  and from the vantage point of $\tau_R(u)$, $u$ looks like $g$
  trivial strips $\y\times[0,1]\times\RR$. Thus, the number of intersections of $u$
  and $\tau_{R}(u)$ is $u\cdot \tau_{R}(u)=n_\x(B)+n_\y(B)-g/2$; the
  $g/2$ comes from the $2g$ corners. We will show that
  \begin{equation}\label{eq:intChange}
    u\cdot\tau_{\epsilon}(u)-u\cdot\tau_{R}(u)=\sum_{i < j}L(\rho_i,\rho_j).
  \end{equation}
  and so
  \begin{align*}
    \chi(S) &= e(B) - \br(u) + g/2 + \abs{\vec{P}}/2\\
      &= e(B) - u\cdot \tau_\eps(u) + g/2 + \abs{\vec{P}}/2\\
      &= e(B) - u\cdot \tau_R(u) + g/2 - \sum_{i<j} L(\rho_i,\rho_j) +
        \abs{\vec{P}}/2\\
      &= g + e(B) - n_\x(B) - n_\y(B) -\iota([\vec{P}])
  \end{align*}
  as desired.
  
  It follows from a simple Schwartz reflection argument that the
  intersection number $u\cdot \tau_{r}(u)$ can only change when an
  $e$~puncture of $\tau_{r}(u)$, mapped to some Reeb chord~$\rho_1$, passes an
  $e$~puncture of $u$ mapped to some Reeb chord~$\rho_2$. Call this change
  $-L'(\rho_1,\rho_2)$; that is, $L'(\rho_1,\rho_2)$ is the number of
  double points which disappear when $\rho_1$ goes from below $\rho_2$ to
  above~$\rho_2$. We claim that $L'(\rho_1,\rho_2)$ is universal,
  depending only on $\rho_1$ and $\rho_2$, and not on $u$
  itself. Equation~(\ref{eq:intChange}) will then follow from a few
  model computations.

  To see that $L'(\rho_1,\rho_2)$ is universal we employ a doubling
  argument, which also computes $L'(\rho_1,\rho_2)$. That is, we
  construct a new Heegaard diagram with boundary
  $(\Sigma',\alphas',\betas',z)$, which we glue to $\Sigma$ along the
  boundary. We also construct holomorphic curves $u_2$ and $u_1$ which
  can be glued to $u$ and $\tau_r(u)$ respectively, yielding curves
  $uu_2$ and $\tau_r(u)u_1$, respectively. Since the number of
  intersections between $uu_2$ and $\tau_r(u)u_1$ is topological,
  remaining constant under deformations of $uu_2$ and $\tau_r(u)u_1$,
  the number $L'(\rho_1,\rho_2)$ is equal to the number of double
  points which appear when $u_1$ is slid up relative to
  $u_2$.

  In fact, it suffices to construct $\Sigma'$, $u_1$, $u_2$ and so on
  locally near $\rho_1\cup\rho_2$. The constructions of $\Sigma'$,
  $u_1$ and $u_2$ can then be made to depend only on $\rho_1$ and
  $\rho_2$, not on the curve $u$ itself. This implies that
  $L'(\rho_1,\rho_2)$ depends only on $\rho_1$ and $\rho_2$.

  Our model curves will have the domains shown in
  Figure~\ref{fig:IndexModels}. We distinguish 10 different cases of
  how $\rho_1$ and $\rho_2$ can overlap.
  \begin{enumerate}
  \item\label{case:ind-cross:1} $\rho_1^+=\rho_2^-$. In this case, $L'(\rho_1,\rho_2)=1/2$. We
    choose the curves $u_1$ and $u_2$ to consist of one strip each, as
    indicated in Figure~\ref{fig:IndexModels} (upper left). When
    sliding $u_1$ up past $u_2$, one boundary double point
    appears.
  \item\label{case:ind-cross:2} $\rho_1^-=\rho_2^+$. In this case,
    $L'(\rho_1,\rho_2)=-1/2$. This case is the same as case~(\ref{case:ind-cross:1}),
    with the roles of $u_1$ and $u_2$ exchanged.
  \item\label{case:ind-cross:3} $\rho_1^-\lessdot \rho_2^-\lessdot
    \rho_1^+=\rho_2^+$. In this case, $L'(\rho_1,\rho_2)=1/2$. We
    again choose $u_1$ and $u_2$ to consist of one strip each, though
    in this case they overlap, as indicated in
    Figure~\ref{fig:IndexModels} (top right). When sliding $u_1$ up past $u_2$,
    one boundary double point appears.  (The fact that $u_1\cup u_2$
    is not a valid holomorphic curve---its $-\infty$ asymptotics are
    not a generator---is irrelevant: we are only interested in its
    behavior near $\rho_1\cup \rho_2$.)
  \item\label{case:ind-cross:4} $\rho_1^- = \rho_2^-\lessdot
    \rho_1^+\lessdot \rho_2^+$. In this case,
    $L'(\rho_1,\rho_2)=1/2$. Similarly to
    case~(\ref{case:ind-cross:3}), we choose $u_1$ and $u_2$ to
    consist of one strip each, as indicated in
    Figure~\ref{fig:IndexModels} (center left). Again, when sliding $u_1$ up past
    $u_2$, one boundary double point appears.
  \item $\rho_2^-\lessdot \rho_1^-\lessdot \rho_1^+=\rho_2^+$. In this
    case $L'(\rho_1,\rho_2)=-1/2$. This case is the same as case~(\ref{case:ind-cross:4}) with
    the roles of $u_1$ and $u_2$ exchanged.
  \item\label{case:ind-cross:6} $\rho_1^- = \rho_2^-\lessdot \rho_2^+\lessdot \rho_1^+$. In this
    case $L'(\rho_1,\rho_2)=-1/2$. This case is the same as case~(\ref{case:ind-cross:3}) with
    the roles of $u_1$ and $u_2$ exchanged.
  \item\label{case:ind-cross:7} $\rho_1$ is nested inside $\rho_2$ or vice-versa. In this
    case, $L'=0$. Indeed, taking $u_1$ and~$u_2$ to consist of one
    strip each, as indicated in Figure~\ref{fig:IndexModels} (upper
    right), there are no intersection points between $u_1$ and $u_2$ no
    matter how they are slid. (A degenerate sub-case is when $\rho_1 =
    \rho_2$.)
  \item\label{case:ind-cross:8} $\rho_1$ and $\rho_2$ are disjoint. In this case,
    $L'(\rho_1,\rho_2)=0$. Taking $u_1$ and $u_2$ to each consist of
    a single strip, as in Figure~\ref{fig:IndexModels} (lower left)
    again gives no intersection points, no matter how they are slid.
  \item\label{case:ind-cross:9} $(\rho_1,\rho_2)$ is interleaved, with
    $\rho_1^-\lessdot \rho_2^-\lessdot \rho_1^+\lessdot \rho_2^+$. In
    this case, $L'(\rho_1,\rho_2)=1$. Take $u_1$ and $u_2$ to
    consist of a single strip each, as in Figure~\ref{fig:IndexModels}
    (lower right). Then if $u_2$ is slid much higher than $u_1$, there
    are no intersection points, while if $u_2$ is much lower than
    $u_1$, there is a single interior intersection point.
  \item $(\rho_2,\rho_1)$ is interleaved, with
    $\rho_2^-\lessdot \rho_1^-\lessdot \rho_2^+\lessdot \rho_1^+$. In
    this case, $L'=-1$. This is the same as case the previous case,
    with the roles of $u_1$ and $u_2$ reversed.
\end{enumerate}

This completes the proof in the case that $\vec{P}$ is a discrete
partition.

In the general case, it is not true that $\br(u)$ is equal
to the intersection number of $u$ with $\tau_\epsilon(u)$ for
$\epsilon$ small. That is, there may be intersection points of $u$ and
$\tau_\epsilon(u)$ which run off to east~$\infty$ as $\epsilon\to
0$. This phenomenon is a simple variant of the change in number of
double points studied above, and it follows from the model computations
used above that
\[
u\cdot\tau_\epsilon(u)-\br(u)=\sum_i\sum_{\{\rho_a,\rho_b\}\subset [P_i]}|L(\rho_a,\rho_b)|.
\]
Otherwise, the argument proceeds as before. The result follows
from Lemma~\ref{lem:iota-chords}.

Finally, if $\cM^B(\x,\y;\Source;\vec{P})$ admits a non-embedded
holomorphic representative then 
\[
u\cdot\tau_\epsilon(u)-\br(u)=2\sing(u)+\sum_i\sum_{\{\rho_a,\rho_b\}\subset [P_i]}|L(\rho_a,\rho_b)|,
\]
where $\sing(u)$ denotes the order of singularity of $u$ (i.e.,
number of double-points in an equivalent singularity;
compare~\cite{MicallefWhite95:intersection-positivity}). The rest of
the argument goes through, giving
\begin{align*}
  \chi(S)&=g+e(B)-n_\x(B)-n_\y(B)-\iota([\vec{P}])+2\sing(u)\\
  \ind(u)&=e(B)+n_\x(B)+n_\y(B)+|\vec{\rhos}|+\iota(\vec{\rhos})-2\sing(u)\\
  &=\ind(B,\vec{\rhos})-2\sing(u).
\end{align*}
Strong boundary monotonicity rules out boundary double points in $u$
(which contribute $1/2$ to $\sing(u)$), so $\sing(u)$ is a positive
integer. Thus, $\ind(u)\leq \ind(B,\vec{\rhos})-2$, as desired.
\end{proof}

\index{expected dimension!at embedded curve|)}
\colorused

\subsection{The embedded index formula is additive}
\label{sec:index-additive}
Since Proposition~\ref{prop:asympt_gives_chi} assumed the existence of a
holomorphic curve, it is not obvious that $\ind(B,\vec{\rhos})$ is
additive in cases when there is no holomorphic representative. To
prove this, we will adapt Sarkar's proof of the same fact
in the closed case \cite[Theorem
3.2]{Sarkar11:IndexTriangles} to our situation. We start with a
definition and some lemmas.

\begin{figure}
  \includegraphics[scale=.83333]{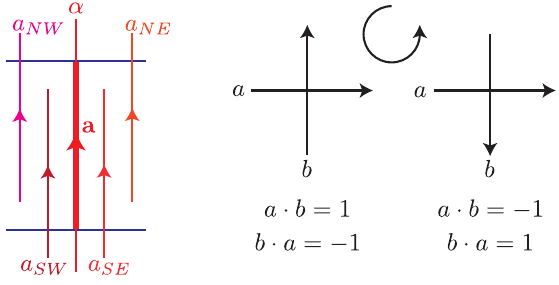}
  \caption[Jittered intersection number]{\textbf{Jittered intersection number.} Left: the four translates $a_{SW}$, $a_{SE}$, $a_{NW}$, and
    $a_{NE}$ of an arc $a$. Right: conventions for intersection
    numbers.}\label{fig:intersecting-arcs}
\end{figure}

Suppose $a$ is an oriented sub-arc of $\alphas$, with endpoints in
$\alphas\cap\betas$, and $b$ is an oriented sub-arc of $\betas$, with
endpoints in $\alphas\cap\betas$. We define the \emph{jittered
  intersection number $a\cdot b$} 
\index{jittered intersection number}\index{intersection number, jittered}%
as follows. 
\glsadd{aSW}%
Let $a_{SW}$, $a_{SE}$,
$a_{NW}$, and $a_{NE}$ denote the four translates of $a$ shown in
Figure~\ref{fig:intersecting-arcs}. The endpoints of these translates
are disjoint from $\betas$, and the translates do not intersect
$\alphas\cap\betas$. The translates inherit orientations from the
orientations of $a$ and $b$. Consequently, the intersection numbers
$a_{SW}\cdot b,\dots,a_{NE}\cdot b$ are well defined. Define
\[
\gls*{intnum}
=\frac{1}{4}\left(a_{SW}\cdot b+a_{SE}\cdot b+a_{NW}\cdot b+a_{NE}\cdot b\right).
\]
Note that instead of jittering $a$ we could have jittered $b$ with
the same result. We extend the definition of $a\cdot b$ bilinearly to
cellular $1$-chains $a$ contained in $\alphas$ and $b$ contained in
$\betas$. We define $b\cdot a$ analogously.  We also define $a\cdot
a'$ for $a$ and $a'$ both contained in $\alphas$ (or $\betas$) in
exactly the same way.  (This definition is analogous to the
definitions of $n_\x(B)$ for Proposition~\ref{prop:asympt_gives_chi}
and $m(\alpha,P)$ in Section~\ref{sec:pregrading}.)

\colorused

The following lemma is due to Sarkar \cite[Section
3]{Sarkar11:IndexTriangles}. For the reader's convenience we repeat
the proof here.

\begin{lemma}\label{lemma:jittered} The jittered intersection number
  $a\cdot b$ has the following properties:
\begin{enumerate}
\item \label{item:antisymmetric}$a\cdot b=-b\cdot a$.
\item \label{item:containedinalphas} If $a$ and $a'$ are both
  contained in $\alphas$ then $a\cdot a'=0$.
\end{enumerate}
For the remaining items, let $B \in \pi_2(\x,\y)$,
$a = \bdy^\alpha(B)$, $b = \bdy^\beta(B)$,
$B' \in \pi_2(\y,\w)$,
$a' = \bdy^\alpha(B')$ and $b' = \bdy^\beta(B')$.
\begin{enumerate}\setcounter{enumi}{2}
\item \label{item:nxny}
  $n_\x(B)-n_\y(B)=a\cdot b$.
\item \label{item:nxnw}
  $n_\x(B')-n_\y(B')=a'\cdot b$ and $n_\y(B)-n_\w(B)=a\cdot b'$.
\item \label{item:intsymmetry} If $B$ and $B'$ are
  provincial then $a\cdot b'+b\cdot a'=0$.
\end{enumerate}
\end{lemma}
\begin{proof}
  Properties (\ref{item:antisymmetric}) and
  (\ref{item:containedinalphas}) are obvious.

  To understand Property (\ref{item:nxny}), assume for simplicity that
  $a$ and $b$ are embedded arcs, and consider the path $b_{SE}$. The
  multiplicity of $B$ at the starting point of $b_{SE}$ is
  $n_{\y_{SE}}(B)$. Each time $b_{SE}$ crosses $a$, the multiplicity of $B$
  increases by $1$, for intersections contributing positively to
  $a\cdot b$, or decreases by $1$, for intersections contributing
  negatively to $a\cdot b$. Consequently, at the ending point of
  $b_{SE}$, the multiplicity $n_{\x_{SE}}$ is $n_{\y_{SE}}+a\cdot b_{SE}$
  and so $n_{\x_{SE}}(B) - n_{\y_{SE}}(B) = a\cdot b_{SE}$.  Averaging over the
  four possible directions gives the result.
  Property (\ref{item:nxnw})
  is proved similarly to Property (\ref{item:nxny}), using
  $b'$ instead of~$b$. See
  Figure~\ref{fig:nxnynw} for an example
  of these two properties.

  Finally, Property~(\ref{item:intsymmetry}) follows from the fact
  that $a+b$ and $a' + b'$ are null-homologous $1$-chains, so
  $(a+b)\cdot (a'+b')=0$. But by bilinearity and
  Property~(\ref{item:containedinalphas}), this reduces to
  $a\cdot b' + b\cdot a'=0$. (An alternate way to see Property~(\ref{item:intsymmetry}), and see
that the signs are right, is to use $\x - \y = \bdy(\bdy^\beta B)$ to
compute $n_\x(B') - n_\y(B') = B' \cdot (\partial(\bdy^\beta B)) =
(\partial B') \cdot \bdy^\beta B$, using an extension of the jittered
intersection number to 2-chains and 0-chains.)
\begin{figure}
\includegraphics[scale=.55556]{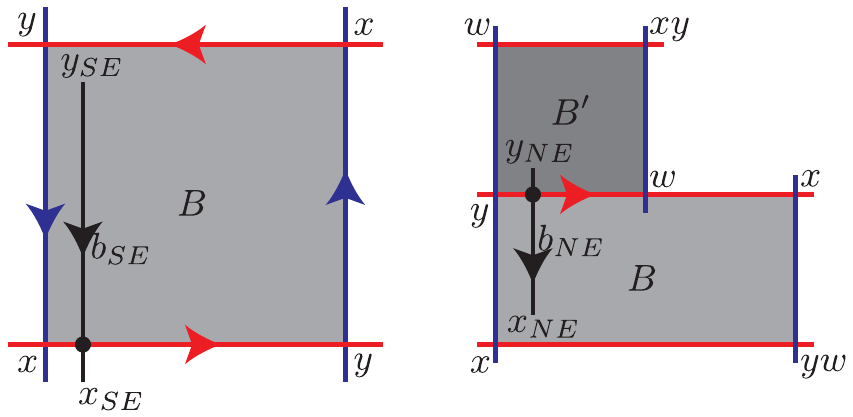}
\caption[Parts of the proof of Lemma~\ref{lemma:jittered} (index additivity)]{\textbf{Parts (\ref{item:nxny}) and
  (\ref{item:nxnw}) of the proof of Lemma~\ref{lemma:jittered}.} Left:
  part~(\ref{item:nxny}). The
  arrows show the orientations of $a$ and $b$ from $\partial
  B$. Notice that $a\cdot b_{SE}=-1$ and $n_{\y_{SE}}(B)=n_{\x_{SE}}(B)-(a\cdot
  b_{SE})$. Right: part (\ref{item:nxnw}). The horizontal arrow indicates the
  orientation induced by $\bdy B'$, the vertical arrow the orientation
  induced by $\bdy B$. Notice that $a'\cdot
  b_{NE}=-1=n_{\x_{NE}}(B')-n_{\y_{NE}}(B')$.}\label{fig:nxnynw}
\end{figure}
\end{proof}

We can extend Property~(\ref{item:intsymmetry}) of
Lemma~\ref{lemma:jittered} to non-provincial domains, as well.
\begin{lemma}\label{lem:jittered2}
  With notation as in Lemma~\ref{lemma:jittered}, for arbitrary $B$
  and~$B'$ in $\pi_2(\x,\y)$ and $\pi_2(\y,\w)$, respectively, we have $a \cdot b' + b \cdot a' = L(\bdy^\bdy B,
  \bdy^\bdy B')$.
\end{lemma}
\begin{proof}
  We now have that $a + b + \bdy^\bdy B$ and $a' + b' + \bdy^\bdy B'$
  are null-homologous 1-chains, so their cap product in homology
  is~$0$.  We can relate this to an extension of the jittered
  intersection product by, for instance, treating $\bdy \bSigma$ as an
  extra $\beta$-circle.  Computing this intersection number and dropping terms that always vanish gives:
  \[
  (a + b + \bdy^\bdy B) \cdot (a' + b' + \bdy^\bdy B') =
    a \cdot b' + b \cdot a' + a \cdot \bdy^\bdy B' +
    \bdy^\bdy B \cdot a' = 0.
  \]
  Clearly,  $a \cdot \bdy^\bdy B' = -\OneHalf m(\bdy^\bdy
  B', \bdy (\bdy^\bdy B))$ and $\bdy^\bdy B \cdot a' = \OneHalf
  m(\bdy^\bdy B, \bdy (\bdy^\bdy B'))$.  Both of these terms are equal
  to $-\OneHalf L(\bdy^\bdy B, \bdy^\bdy B')$, as desired.
\end{proof}

\begin{proposition}\label{prop:indAdditive} 
  \index{index!embedded!is additive}%
  The embedded index is
  additive in the following sense. Let $B\in\pi_2(\x,\y)$ and
  $B'\in\pi_2(\y,\w)$, and let $\vec{\rhos}$ and~$\vec{\rhos}'$ be
  sequences of non-empty sets of Reeb chords with $[\vec{\rhos}] = \bdy^\bdy B$
  and $[\vec{\rhos'}] = \bdy^\bdy B'$. Then
  $\ind(B,\vec{\rhos})+\ind(B',\vec{\rhos}')=\ind(B*B',(\vec{\rhos},\vec{\rhos}'))$
\end{proposition}
(Here, $\gls*{ConcatSeqs}$
denotes the sequence obtained by
simply concatenating $\vec{\rhos}$ and $\vec{\rhos}'$.)
\begin{proof}
  It is clear that terms $e(B)$ and $|\vec\rhos|$ of
  Definition~\ref{def:emb-ind-emb-chi} are additive. Thus, we must show that
\[
  n_\x(B)+n_\y(B)+n_\y(B')+n_\w(B')+\iota(\vec\rhos)+\iota(\vec{\rhos}')
  =n_\x(B+B')+n_\w(B+B')+\iota((\vec{\rhos},\vec{\rhos}')).
\]
On the other hand, we have 
\begin{align*}
[n_\x(B')-n_\y(B')]-[n_\y(B)-n_\w(B)]&=a'\cdot b-a\cdot b'\\
  &=a'\cdot b+b\cdot a' - L(\bdy^\bdy B, \bdy^\bdy B')\\
  &=\iota(\vec{\rhos}) + \iota(\vec{\rhos}') - \iota((\vec{\rhos},\vec{\rhos}')),
\end{align*}
where the first equality uses property~(\ref{item:nxnw}) of
Lemma \ref{lemma:jittered}, the second
uses Lemma {lem:jittered2}, and the third uses
property~(\ref{item:antisymmetric}) and the definition of~$\iota$.
This is
equivalent to what we were trying to show.
\end{proof}

\subsection{The codimension-one boundary of the moduli spaces of embedded curves}\label{sec:embedded-degen}
We finish this chapter with two lemmas relating to which holomorphic
curves occur in the boundaries of moduli spaces of embedded curves. 
We start by showing that when embedded holomorphic curves converge to
holomorphic combs, the stories are embedded.
\begin{lemma}\label{lem:splittings-embedded} For a generic almost
  complex structure~$J$ and any pair $(B,\vec{\rhos})$ with
  $\ind(B,\vec{\rhos})=2$, the height $2$
  combs which occur in the boundary of the embedded moduli space
  $\ocM^B(\x,\y\semico\vec{\rhos})$ (as in
  Theorem~\ref{thm:master_equation}, say) are
  embedded. Equivalently, if $(u_1,u_2)$ is such a height $2$
  comb, with $\SourceSub{i}$ the source of $u_i$, $B_i$ the
  homology class of $B_i$ and $\vec{P}_i$ the ordered partition
  associated to $u_i$, then:
  \begin{itemize}
  \item $\chi(\SourceSub{i})=\chi_{\emb}(B_i,\vec{P}_i)$
  \item $\ind(B_i,P_i)=1$.
  \end{itemize}
\end{lemma}
\begin{proof}
  This follows from the inequality $\ind(u)\leq\ind(B,\vec{\rhos})$
  when $u$ is non-embedded in Proposition~\ref{prop:asympt_gives_chi}
  and the fact that moduli spaces of negative expected dimension
  ($\ind\leq 0$) are empty (Proposition~\ref{prop:transversality}).
\end{proof}

The next lemma states that when two levels collide in a sequence of
embedded curves, the corresponding parts of a partition must be
composable (as in Section~\ref{sec:reeb-chords-def}), and not merely
weakly composable (Definition~\ref{def:weakly-composable}).  Similar
statements hold for splittings and shuffles.

\begin{lemma}\label{lemma:collision-is-composable}
  Let $(\x,\vec{\rhos})$ satisfy the strong boundary
  monotonicity condition, and $B\in\pi_2(\x,\y)$ be compatible
  with~$\vec\rhos$. Suppose that
  $\ind(B,\vec{\rhos})=2$.  Suppose there is a holomorphic comb
  $(u,v)$ in $\bdy\ocM^B(\x,\y\semico\vec\rhos)$, and let $\vec\rhos'$ be
  the sequence of sets of Reeb chords appearing as asymptotics of~$u$.  Then
  \begin{enumerate}
  \item if two levels $\rhos_i$ and $\rhos_{i+1}$ collide in~$u$,
    $\rhos_i$ and $\rhos_{i+1}$ are composable;
  \item if $u$ is a join curve end at level~$i$,
    $\rhos_i'$ is a splitting of~$\rhos_i$; and
  \item if $u$ is an odd shuffle curve end at level~$i$,
    $\rhos_i'$ is a shuffle of~$\rhos_i$.
  \end{enumerate}
\end{lemma}

As a preliminary, we first state a characterization of composability,
splittings, and shuffles in terms of gradings.  Recall that there is a partial order on
$\bigGroup(n)$ by comparing the Maslov components: $(k_1,\alpha_1) <
(k_2,\alpha_2)$ if $\alpha_1 = \alpha_2$ and $k_1 < k_2$.

\begin{lemma}\label{lem:composable-gr}
Let $\rhos$ and $\sigmas$ be sets of Reeb chords such that
$a(\rhos)$, $a(\sigmas)$ and $a(\rhos\uplus\sigmas)$ are non-zero.
Then:
\begin{enumerate}
\item The pair $(\rhos,\sigmas)$ is composable if and only if
  \[\grb(a(\rhos\uplus\sigmas)) = \grb(a(\rhos))\cdot\allowbreak\grb(a(\sigmas)).\]
\item A weak splitting $\rhos'$ of $\rhos$ is a splitting if and
  only if
  \[\grb(a(\rhos')) = \lambda^{-1}\grb(a(\rhos)).\]
\item A weak shuffle $\rhos'$ of $\rhos$ is a shuffle if and only if
  \[\grb(a(\rhos')) = \lambda^{-1}\grb(a(\rhos)).\]
\end{enumerate}
In all cases, if equality does not hold then the left-hand side has lower
grading.
\end{lemma}

\begin{proof}
  Clear.
\end{proof}

\begin{proof}[Proof of Lemma~\ref{lemma:collision-is-composable}]
  By Lemma~\ref{lem:splittings-embedded}, $u$ is embedded.

  For the first statement, by Theorem~\ref{thm:master_equation} we know that
  $\rhos_i$ and $\rhos_{i+1}$ are weakly composable.  By
  Lemma~\ref{lem:composable-gr}, if $\rhos_i$ and $\rhos_{i+1}$ are
  weakly composable but not composable, then
  $\grb(a(\rhos_i\uplus\rhos_{i+1})) <
  \grb(a(\rhos_i))\grb(a(\rhos_{i+1}))$, so, by
  Lemma~\ref{lem:iota-grading}, $\iota(\vec\rhos') <
  \iota(\vec\rhos)$.  Therefore $\ind(B,\vec\rhos') <
  \ind(B,\vec\rhos)+ (\abs{\vec\rhos'}-\abs{\vec\rhos}) = 1$. Thus
  the embedded moduli space $\cM^B(\x,\y;\vec\rhos')$
  has negative expected dimension and is therefore empty.

  Similarly, a join curve end at level~$i$ gives a weak splitting of
  $\rhos_i$ by definition.  By Lemma~\ref{lem:composable-gr}, for a
  weak splitting that is not a splitting we have $\grb(a(\rhos_i')) <
  \grb(a(\rhos_i)) - 1$.  This implies that $\ind(B,\vec\rhos') < 1$,
  so again $\cM^B(\x,\y;\vec\rhos')$ has negative expected dimension.

  The case of odd shuffle curve ends is similar.
\end{proof}

\begin{remark}
  A similar analysis shows that \emph{even} shuffle curve ends on the
  boundary of a 1-dimensional embedded moduli space all have negative
  expected dimension.  We do not use that, as we already showed in
  case~(\ref{item:even}) of Proposition~\ref{prop:gluing-shuffle} that
  such curves appear an even number of times on the boundary of any
  moduli space.
\end{remark}


\chapter{Type \textalt{$D$}{D} modules}
\label{chap:type-d-mod}

We now turn to defining the type $D$ module $\CFDa$ of a bordered
Heegaard diagram, using the moduli spaces developed in
Chapter~\ref{chap:structure-moduli}. The definition, in
Section~\ref{sec:def-CFD}, is fairly straightforward. More
work goes into proving that $\bdy^2=0$
(Section~\ref{sec:typed-d-sq-zero}) and that the module is
independent of the choices made in its definition
(Section~\ref{sec:typeD-invariance}). We also introduce a version of
$\CFDa$ with twisted coefficients, in
Section~\ref{sec:typeD-twisted}.
We defer the construction of the grading on $\CFDa$ to
Chapter~\ref{chap:gradings}.

\section{Definition of the type \textalt{$D$}{D} module}
\label{sec:def-CFD}\index{type $D$!invariant|(}
Let $\Heegaard$ be a bordered Heegaard diagram (in the sense of Definition~\ref{def:BorderedDiagram})
which is provincially admissible (in the sense of Definition~\ref{def:provincial-admissibility})
for a bordered three-manifold $Y$, and let 
$\PtdMatchCirc$ be the \emph{reverse} of the pointed matched circle
(in the sense of Definition~\ref{def:PointedMatchedCircle})
appearing on the boundary of~$\Heegaard=(\Sigma,\alphas,\betas,z)$.
That is, $\PtdMatchCirc$ is
$(Z,\widebar{\alphas}\cap\bdy\widebar{\Sigma},M,z)$,
where $Z = -\bdy\widebar{\Sigma}$ has the opposite of the boundary orientation. 
Let $\Alg(\PtdMatchCirc)$ be the associated algebra,
as in Chapter~\ref{chap:algebra}. We will abbreviate
$\Alg(\PtdMatchCirc)$ by~$\Alg$ when it causes no
confusion.\index{orientation reversal}

Fix a $\SpinC$ structure $\s$ over $Y$. The goal of this section is to define a left \dg module $\CFDa(\Heegaard,\spinc)$ over~$\Alg$.  

Let $\S(\HD)$ be the set of generators of $\HD$, as defined in
Definition~\ref{def:generator}. That is, $\Gen(\HD)$ consists of subsets of
$\alphas \cap \betas$ with
\begin{itemize}
\item $g$ elements in all,
\item exactly one element on each $\beta$ circle,
\item exactly one element on each $\alpha$ circle, and
\item at most one element on each $\alpha$ arc.
\end{itemize}
As in Section~\ref{sec:homology-classes-generators}, let $\S(\HD,\spinc)\subset \S(\HD)$ denote the subset
of generators representing the (absolute) $\SpinC$ structure $\spinc$.
Let 
$\gls*{XofHs}$
be the $\Field$-vector space spanned by $\S(\HD,\spinc)$ and
$\gls*{XofH}$
the $\Field$-vector spanned by $\S(\HD)$, so that
\[X(\HD)=\bigoplus_{\spinc\in\SpinC(Y)}X(\HD,\spinc).\]
Recall from Section~\ref{sec:homology-classes-generators} that, for
$\x \in \S(\HD)$, $o(\x) \subset [2k]$ denotes the set of $\alpha$-arcs
occupied by~$\x$.
Define 
$\gls*{IDofx}$
to be $I([2k]\setminus o(\x))$; that is, $I_D(\x)$ is the idempotent
corresponding to the complement of $o(\x)$.
Now, define an action on $X(\HD,\spinc)$ of
the sub-algebra of idempotents~$\Idem\subset\Alg$ via
\begin{equation}\label{eq:def-idem-DMod}
  I(\SetS) \cdot \x \coloneqq
  \begin{cases}
    \x&I(\SetS) = I_D(\x)\\
    0&\text{otherwise.}
  \end{cases}
\end{equation}
(Here, $\SetS$ is a $k$-element subset of $[2k]$, as in
Formula~\eqref{eq:idempotent-def-2}.)

As an $\Alg$-module, $\CFDa(\HD,\spinc)$ is defined by
\begin{equation*}
  \gls*{CFDspinc}
 \coloneqq \Alg \otimes_{\Idem} X(\HD,\spinc).
\end{equation*}
 So, the module structure on $\CFDa(\HD,\spinc)$ is quite simple; it is
given by
\begin{equation}
  \label{eq:def-DMod-mult}
  a\cdot(b\otimes\x) \coloneqq (ab) \otimes \x.
\end{equation}
In particular, the summands
$\Alg(\PtdMatchCirc,i)$ of $\Alg(\PtdMatchCirc)$ for $i\neq 0$ act trivially on
$\CFDa(\HD,\spinc)$.

By contrast, the differential on $\CFDa(\HD,\spinc)$ involves
counting moduli spaces of holomorphic curves $\Mod^B(\x,\y;\vec{\rho})$, where 
$\gls*{ReebChordSeq}$
is a sequence $(\rho_1,\dots,\rho_n)$ of Reeb chords.  
(Here we are
using $\vec\rho$ both for a sequence of Reeb chords and for the
corresponding sequence of one-element sets of Reeb chords.)

Recall from Section~\ref{sec:matched-circles} that $a(\rho_i)$
denotes the element of
$\Alg$ associated to $\rho_i$.
For $\vec\rho$ a sequence of Reeb chords, let 
$\gls*{aofvecrho}$
 be the
product $a(\rho_1)\cdots a(\rho_n)$.
(The reader is cautioned that this typically does not coincide
with the algebra element $a(\rhos)$ associated to a set of Reeb chords
$\rhos=\{\rho_1,\dots,\rho_n\}$.)
Let $-\vec{\rho}$ be the sequence $(-\rho_1,\dots,-\rho_n)$ of
chords with reversed orientation. (If $\rho$ is a Reeb chord in
$\bdy\Heegaard$, $-\rho$ is a Reeb chord
in~$\PtdMatchCirc=-\bdy\Heegaard$.)\index{orientation reversal}

\begin{definition}\label{def:Dmod-boundary}
  Fix an almost-complex structure on $\Sigma\times[0,1]\times\RR$
  which is admissible in the sense of
  Definition~\ref{def:admissible_J} and sufficiently generic in the
  sense of Definition~\ref{def:sufficiently-generic}.
  For $\x,\y\in\S(\HD,\spinc)$ and $B \in \pi_2(\x,\y)$, define\glsadd{aofminusvecrho}
  \[
  \gls*{aBxy}
  \coloneqq\!\!\sum_{\{\vec\rho\,\mid\,\ind(B,\vec{\rho}\,) = 1\}}\!\!\#\bigl(\Mod^B(\x, \y; \vec{\rho})\bigr)
  a(-\vec{\rho}).
  \]
  Here there is also an implicit condition that $(B,\vec\rho)$ is
  compatible (i.e., that $(\x,\vec\rho)$ is strongly boundary monotone
  and $\bdy^\bdy B = [\vec\rho\,]$); we will often omit such conditions from the notation.
  Compactness of the moduli spaces $\ocM^B(\x,\y;\vec{\rho})$
  implies that the sum defining $a^B_{\x,\y}$ is finite.

  Define the differential on $\CFDa(\HD,\spinc)$ by\index{differential!on $\CFDa$}
  \begin{equation*}
    \partial(\Unit\otimes \x) \coloneqq \sum_{\y\in\S(\HD,\spinc)}\sum_{B\in\pi_2(\x,\y)}a^B_{\x,\y} \otimes \y.
  \end{equation*}
  We will verify in Lemma~\ref{lem:finite-typeD} that if $\HD$ is
  provincially admissible then the sums defining~$\partial$ are
  finite.
  Extend $\partial$ by linearity and the Leibniz rule to all of $\CFDa(\HD,\spinc)$.

  Let 
$\gls*{CFD}
\coloneqq\bigoplus_{\spinc\in\SpinC(Y)}\CFDa(\HD,\spinc)$.
\end{definition}

We will often abbreviate elements $a\otimes \x\in \CFDa(\HD,\spinc)$ writing, simply $a\x$.
In particular, we blur the distinction between generators $\x\in\S(\HD,\spinc)$ 
and their associated elements $\Unit\otimes \x\in\CFDa(\HD,\spinc)$

\begin{remark}
  \label{rmk:UnderlyingTypeD}
  Consider the $\Idem$-module $X(\HD,\spinc)$, equipped with the
  map $\delta^1\co X(\HD,\spinc) \to \Alg\otimes_\Idem X(\HD,\spinc)$ given
  by
  \[
  \delta^1(\x) \coloneqq \partial(\Unit\otimes \x).
  \]
  This
  defines a type~$D$ structure over $\Alg$ with base ring $\Idem$ on
  $X(\HD,\spinc)$, in the sense of Definition~\ref{def:TypeD}.  (The
  compatibility relation is equivalent to
  Proposition~\ref{prop:typeD-d2} below.)
  $\CFDa(\HD)$, then, is the differential module associated
  to this type $D$ structure, as in Lemma~\ref{lemma:AssociatedTypeD}.
  Further, Theorem~\ref{thm:D-invariance}, below, and
  Lemma~\ref{lemma:ConverseAssociatedTypeD} imply that
  this type~$D$ structure is invariant up to homotopy as a type~$D$
  structure, not just as a module. 
\end{remark}

\begin{lemma}\label{lem:finite-typeD}
  If $\HD$ is provincially admissible, then the sum in the
  definition of $\partial\x$ on $\CFDa(\HD,\spinc)$ is finite for
  every $\x\in\S(\HD,\spinc)$.
  If $\HD$ is admissible, the map $\delta^1$ (for the underlying type~$D$ structure,
  as in Remark~\ref{rmk:UnderlyingTypeD}) is bounded
  in the sense of Definition~\ref{def:BoundedTypeD}.
\end{lemma}
\begin{proof}
By Lemma~\ref{lem:holo-has-pos-domain}, if $\cM^B(\x,\y;\vec{\rho})$ is non-empty then $B$ has a positive domain. From
Proposition~\ref{prop:provincial-admis-finiteness} we then see that if
$\HD$ is provincially admissible, then for a
given $\x$, $\y$ and
  $\vec{\rho}$ the union over $B$ of 0-dimensional moduli spaces
  $\Mod^B(\x,\y;\vec{\rho})$ is finite.  There
  are only finitely many possible generators~$\y$ and algebra
  elements~$a\in\Alg$; furthermore, for any given~$a$
  there are only finitely many ways to write it as a product of Reeb
  chords, and so only finitely many possible~$\vec\rho$.  The first
  statement then follows.

  If $\HD$ is admissible, then by
  Proposition~\ref{prop:admis-finiteness}, there are only finitely
  many $B\in\pi_2(\x,\y)$ with non-negative coefficients, and hence
  only finitely many possible terms in any $\delta^k(\x)$. More
  precisely, terms in $\delta^k(\x)$ count points in a product space
  $\prod_{i=1}^{k}\Mod^{B_i}(\x_i,\x_{i+1};\vec{\rho_i})$, where we
  sum over all $k+1$-tuples $\{\x_i\}_{i=1}^{k+1}$ in $\Gen$ with
  $\x_1 = \x$ and $\x_{k+1} = \y$,
  $k$-tuples $\{B_i\in\pi_2(\x_i,\x_{i+1})\}_{i=1}^k$ and compatible 
  sequences of Reeb chords $\{\vec{\rho_i}\}_{i=1}^k$.
  Thus, for each non-zero  $\delta^k(\x)$ we can find
  such sequences $\{\x_i\}_{i=1}^{k+1}$ and $\{B_i\in\pi_2(\x_i,\x_{i+1})\}_{i=1}^k$
  with the property that each $B_i$ is positive. Adding
  the $B_i$, we obtain a corresponding $B\in\pi_2(\x,\x_k)$ which is
  also positive. Moreover, if we define $|B|$ to
  be the sum of the coefficients of $B$ then $|B|\geq k$. Since there
  are only finitely many options for $B$, this gives an upper bound on $k$,
  as desired.
\end{proof}

\begin{remark}
\label{rmk:SpinCSplitting}
In the definition of the boundary operator, the requirement that
$\y\in\Gen(\HD,\s)$ is redundant: if $a^B_{\x,\y}\neq 0$, then
$B\in\pi_2(\x,\y)$ is non-empty and hence, by
Lemma~\ref{lem:SpinCStructures}, $\spinc_z(\x)=\spinc_z(\y)$. Thus, we
could have defined $\CFDa(\HD)$ without reference to $\SpinC$
structures, and then noted that
the complex splits as a direct sum over $\SpinC(Y)$.
\end{remark}
\index{type $D$!invariant|)}
\section{\textalt{$\bdy^2=0$}{Boundary-squared is zero}}
\label{sec:typed-d-sq-zero}
This section is devoted to proving that $\CFDa(\HD,\spinc)$ is a differential
module, i.e.:
\begin{proposition}
  \label{prop:typeD-d2}
  The boundary operator~$\partial$ on $\CFDa(\HD,\spinc)$
  satisfies $\partial^2 = 0$.
\end{proposition}

The proof of Proposition~\ref{prop:typeD-d2} involves a number of
pieces, some of which we give as lemmas.  The impatient reader is
encouraged to skip to the examples after the proof, which illustrate
the most important cases which arise.

In our next two lemmas, we look at which ordered lists of Reeb chords
$\vec{\rho}$ actually contribute to the boundary operator. The
second lemma follows from the first, which is more general than our
immediate needs, but will be useful in
Chapter~\ref{chap:tensor-prod}. The first lemma is phrased in terms
of sequences of sets of Reeb chords (rather than merely sequences of Reeb chords).
In particular, if $\vec{\rhos}=(\rhos_1,\dots,\rhos_n)$ is a sequence of sets of Reeb chords,
the notation 
$\gls{ReebChordsSeqSetsOrRev}$
means the corresponding sequence of sets of Reeb chords where 
each Reeb chord has its orientation reversed (but with the ordering
unchanged). Similarly, for a sequence
$\vec{\rho}=(\rho_1,\dots,\rho_n)$ of Reeb chords,
$\gls{ReebChordSeqOrRev}$
denotes $(-\rho_1,\dots,-\rho_n)$.

\begin{lemma}
  \label{lem:BiMonotonicity}
  Let $\SetS\subset [2k]$ be a subset
  and $\vec\rhos=(\rhos_1,\dots,\rhos_n)$ be a sequence of non-empty sets of Reeb chords.
  Let $o(\SetS,\vec{\rhos}_{[1,i]})$ be as in Definition~\ref{def:strong-monotonicity-P}.
  Suppose that for $i=0,\dots,n$, 
  $M(\rhos_{i+1}^-)\subset o(\SetS,\vec{\rhos}_{[1,i]})$
  and
  $M((-\rhos_{i+1})^-)\subset o([2k]\setminus\SetS,-\vec{\rhos}_{[1,i]})$.
  (That is, $(\SetS,\vec{\rhos})$ and $([2k]\setminus\SetS,-\vec{\rhos})$ satisfy
  Conditions~(\ref{SB1}) and~(\ref{SB2})
  of the strong boundary
  monotonicity condition, Definition~\ref{def:strong-monotonicity-P}.) Then:
  \begin{enumerate}
  \item 
    \label{conc:SBM}
    Both $(\SetS,\vec{\rhos})$ and $([2k]\setminus\SetS,-\vec{\rhos})$ are
    strongly boundary monotone.
  \item 
    \label{conc:DisjointSets}
    For $i=1,\dots,n$, the sets $o(\SetS,\vec{\rhos}_{[1,i]})$ and
    $o([2k]\setminus\SetS,-\vec{\rhos}_{[1,i]})$ are disjoint.
  \item 
    \label{conc:DisjointEndpoints}
    For $i=1,\dots,n$, $M(\rhos_i^-)$ and $M(\rhos_i^+)$
    are disjoint.
  \end{enumerate}
\end{lemma}

\begin{proof}
  We prove Conclusions~\eqref{conc:DisjointSets} and
  \eqref{conc:DisjointEndpoints} simultaneously by induction on $i$. The
  case where $i=0$ is the statement that $\SetS$ and $[2k]\setminus\SetS$ are
  disjoint. Assume the result for $i$.  By hypothesis,
  $M((-\rhos_{i+1})^-)\subset o([2k]\setminus\SetS,-\vec{\rhos}_{[1,i]})$, which, by
  the inductive hypothesis, is disjoint from
  $o(\SetS,\vec{\rhos}_{[1,i]})$, which in turn contains
  $M(\rhos_{i+1}^-)$. Since reversing
  orientation switches the sets of endpoints $\rhos_{i+1}^-$
  and\index{orientation reversal}
  $\rhos_{i+1}^+$, $M(-\rhos_{i+1}^-)=M(\rhos_{i+1}^+)$. Thus, we have
  verified Conclusion~\eqref{conc:DisjointEndpoints}.
  Thus, since $M({-\rhos}_{i+1}^+)=M(\rhos_{i+1}^-)$ and $M(\rhos_{i+1}^+)$
  are disjoint, it follows that $o(\SetS,\vec{\rhos}_{[1,i+1]})$
  and $o([2k]\setminus\SetS,-\vec{\rhos}_{[1,i+1]})$ are disjoint, as
  well, verifying Conclusion~(\ref{conc:DisjointSets}).

  Finally, we turn to Conclusion~(\ref{conc:SBM}).
  To see that $(\SetS,\vec{\rhos})$ is boundary monotone, it suffices to prove
  that $M(\rhos^+_{i+1})$ is disjoint from $o(\SetS,\vec{\rhos}_{[1,i]})\setminus M(\rhos^-_{i+1})$.
  But $M(\rhos^+_{i+1})=M(-\rhos^-_{i+1})\subset o([2k]\setminus\SetS,\vec{\rhos}_{[1,i]})$, which
  by Conclusion~\eqref{conc:DisjointSets} is disjoint from $o(\SetS,\vec{\rhos}_{[1,i]})$.
  By symmetry, $([2k]\setminus \SetS,-\vec{\rhos})$ is boundary monotone, as well.
\end{proof}

Lemma~\ref{lem:BiMonotonicity} has the following consequence which is
relevant for our present purposes:

\begin{lemma}\label{lemma:nonzero-implies-distinct}
  If $(\x,\vec\rho)$ is strongly boundary monotone and
  $I_D(\x)a(-\vec\rho)\neq 0$,  then for each $i$,
  $M(\rho_i^-) \ne M(\rho_i^+)$.
\end{lemma}

\begin{proof}
  The condition that $I_D(\x) a(-\vec\rho)\neq 0$ implies that
  $(I_D(\x),-\vec{\rho})=([2k]\setminus o(\x),-\vec{\rho})$ is strongly boundary monotone.
  Since by hypothesis, $(o(\x),\vec{\rho})$ is strongly boundary
  monotone as well, Lemma~\ref{lem:BiMonotonicity}
  applies, and gives the desired conclusion.
\end{proof}

The next two lemmas identify how some of the terms in
Theorem~\ref{thm:master_equation} relate to the algebraic
framework.  Specifically, a term in $\partial^2\x$ can come either
from $a_1\y$ appearing in $\partial\x$ and $a_2\z$ appearing in
$\partial\y$, contributing $a_1a_2\z$ to $\partial^2\x$, or from $a\y$
appearing in $\partial\x$ and $a'$ appearing in $\partial a$,
contributing $a'\y$ to $\partial^2\x$.  
Each type of term has its own interpretation in terms of
moduli spaces.

\begin{lemma}\label{lemma:first_term} For $\x, \w, \y \in \S(\HD,\spinc)$,
  $B_1 \in \pi_2(\x,\w)$ and $B_2\in\pi_2(\w,\y)$, we have
\begin{equation*}
  a^{B_1}_{\x,\w}a^{B_2}_{\w,\y} = \!\!\!\!\!\!\!
  \sum_{\substack{\vec\rho_1,\vec\rho_2\\\ind(B_1*B_2, (\vec\rho_1,\vec\rho_2))=2}}\!\!\!\!\!\!\!
\#\left(\Mod^{B_1}(\x,\w;\vec\rho_1)
\times\Mod^{B_2}(\w,\y;
  \vec\rho_2)\right)
a(-(\vec\rho_1,\vec{\rho_2})).
\end{equation*}
(Here, $(\vec{\rho}_1,\vec{\rho}_2)$ denotes the sequence obtained by
concatenating $\vec{\rho}_1$ and $\vec{\rho}_2$.)
\end{lemma}
\begin{proof}
Clear.
\end{proof}

\begin{lemma}\label{lemma:third_term} For $\x,\y\in\S(\HD,\spinc)$ and $B
  \in \pi_2(\x,\y)$, we have
\begin{equation*}
\bdy a^B_{\x,\y}=\!\!\!
\sum_{\substack{\vec\rho\\\ind(B,\vec{\rho})=2}}\,
\sum_{\substack{i\\\rho_i,\rho_{i+1}\textnormal{ abut}}}\!\!\!
\#\Mod^B(\x,\y;(\rho_1,\dots,\rho_i\uplus\rho_{i+1},\dots,\rho_n))
a(-\vec\rho).
\end{equation*}
\end{lemma}
\begin{proof}
By definition,
\[
a^B_{\x,\y}=
\sum_{\substack{\vec\rho\\\ind(B,\vec\rho)=1}}
\#\left(\Mod^B(\x,\y;\vec\rho)\right)
a(-\rho_1)\cdots a(-\rho_n).
\]
and so
\begin{align*}
\bdy a^B_{\x,\y}&=
\sum_{\substack{\vec\rho\\\ind(B,\vec\rho)=1}}\sum_i
\ \#\left(\Mod^B(\x,\y;\vec\rho)\right)
a(-\rho_1)\cdots\left(\bdy(a(-\rho_i))\right)\cdots
a(-\rho_n)\\
&=\!\sum_{\substack{\vec\rho\\\ind(B,\vec\rho)=1}}\!\sum_i\!\sum_{\substack{\rho_j,\rho_k\\\rho_i=\rho_j\uplus\rho_k}}\!
\#\left(\Mod^B(\x,\y;\vec\rho)\right)
a(-\rho_1)\cdots a(-\rho_j)a(-\rho_k)\cdots
a(-\rho_n).\\
\end{align*}
To see the last step, note that if $\rho_j$ and $\rho_k$ abut (in that order),
then $-\rho_j$ and
$-\rho_k$ do not abut; thus, for $\rho_i=\rho_j\uplus\rho_k$,
\[
\partial(a(-\rho_i))=
\sum_{\substack{\rho_j,\rho_k\\\rho_i=\rho_j\uplus\rho_k}}
a(\{-\rho_j,-\rho_k\}) =
\sum_{\substack{\rho_j,\rho_k\\\rho_i=\rho_j\uplus\rho_k}}
a(-\rho_j)a(-\rho_k).
\]
Now, from the definition of
$\ind(B,\vec{\rho})$ (Definition~\ref{def:emb-ind-emb-chi}),
\[
\ind(B,(\rho_1,\dots,\rho_j,\rho_k,\dots,\rho_n)) =
\ind(B,\vec{\rho})
+1 = 2.
\]
(By contrast,
$\ind(B,(\rho_1,\dots,\rho_k,\rho_j,\dots,\rho_n))=\ind(B,\vec\rho)$.)
The result follows by reindexing the sum to run over $\vec\rho\,' =
(\rho_1,\dots,\rho_j,\rho_k,\dots,\rho_n)$.
\end{proof}

\begin{proof}[Proof of Proposition~\ref{prop:typeD-d2}]

Roughly, the proof proceeds by considering the boundaries of the
index~$1$ moduli spaces. More precisely, we appeal to
Theorem~\ref{thm:master_equation}.

Observe that
\begin{align*}
\partial^2(a\x)&=\partial\biggl[(\partial a)\x+a\biggl(\sum_\w a_{\x,\w}\w\biggr)\biggr]\\
&=(\partial^2a)\x+2(\partial
a)(\partial\x)+a\biggl[\biggl(\sum_\w(\partial a_{\x,\w})\w\biggr)
+\biggl(\sum_\w\sum_\y a_{\x,\w}a_{\w,\y}\y\biggr)\biggr],
\end{align*}
where 
$
\gls*{axynoB}
= \sum_{B\in\pi_2(\x,\y)} a^B_{\x,\y}$.
So it suffices to show that for all $\x$ and $\y$,
\begin{equation}
  \label{eq:typeDequation}
\partial(a_{\x,\y})+\sum_\w a_{\x,\w}a_{\w,\y}=0.
\end{equation}
(This is, of course, equivalent to the compatibility condition for the underlying type $D$ structure, Definition~\ref{def:TypeD}.)

Fix generators $\x$ and $\y$. Fix also $B\in\pi_2(\x,\y)$, and a non-zero
algebra element~$a$. Fix an algebra element $a$ and a
sequence of Reeb chords $\vec\rho$ so that $\ind(B,\vec{\rho}) = 2$.
(The index condition depends only
on~$B$ and~$a$, not the particular sequence $\vec{\rho}$; see
Lemma~\ref{lem:iota-grading} and the definition of
$\ind(B,\vec{\rho})$.)
By Theorem~\ref{thm:master_equation}, applied to
$(\x,\y,B,\vec{\rho})$ and the union, over all sources $\Source$ with
$\chi(S)=\chi_{\emb}(B,\vec{\rho})$, of $\cM^B(\x,\y;\Source;\vec P)$,
the sum of the following terms
is equal to zero:
\begin{enumerate}
\item\label{case:dsquare-two-lev}
  The number of two-story ends, i.e., the number of elements of
\[
  \Mod^{B_1}(\x,\w; \vec\rho_1)\times\Mod^{B_2}(\w,\y; \vec\rho_2)
\]
where $B = B_1 * B_2$ and $\vec\rho = (\vec\rho_1, \vec\rho_2)$.
\item\label{case:dsquare-join}
  The number of join curve ends, i.e., the number of elements of
\[
\Mod^B(\x,\y;(\rho_1,\dots,\rho_{i-1},\{\rho_j,\rho_k\},\rho_{i+1},\dots,\rho_n))
\]
where $\rho_i = \rho_j \uplus \rho_k$.
\item\label{case:dsquare-split}
  The number of split curve ends (necessarily with one
  split component), i.e., the number of elements of
\[
\Mod^B(\x,\y;(\rho_1,\dots,\rho_{i-1},\rho_i\uplus\rho_{i+1},\rho_{i+2},\dots,\rho_n))
\]
 where $\rho_i^+=\rho_{i+1}^-$.
\item\label{case:dsquare-collide} The number of other collisions of levels,
  i.e., the number of elements of
\[
\Mod^B(\x,\y;(\rho_1,\dots,\rho_{i-1},\{\rho_{i},\rho_{i+1}\},\rho_{i+2},\dots,\rho_n))
\]
where $\rho_i^+\ne\rho_{i+1}^-$.
\end{enumerate}
(Both the third and fourth contributions come from Case~(\ref{item:Collision}) of
Theorem~\ref{thm:master_equation}. Shuffle curve ends cannot
happen when the partition is discrete.)

Next, we sum over all choices of vector $\vec{\rho}$  with $a=a(-\vec{\rho})$.
The sum of terms of
type~(\ref{case:dsquare-two-lev}) corresponds to the coefficient in $a$
of
$\sum_{\{\w,B_1,B_2\,\mid\,B_1*B_2=B\}} a^{B_1}_{\x,\w}a^{B_2}_{\w,\y}$, 
according to Lemma~\ref{lemma:first_term}; and the
sum of the  terms of type~(\ref{case:dsquare-split}) corresponds to
the coefficient of $a$ in $\bdy a^B_{\x,\y}$,
according to
Lemma~\ref{lemma:third_term}.
We will show that the remaining two types of terms cancel against each other.

For the terms of type~(\ref{case:dsquare-collide}), there are
several cases that cannot contribute:
\begin{itemize}
\item $M(\rho_i^-)=M(\rho_{i+1}^-)$ or
  $M(\rho_i^+)=M(\rho_{i+1}^+)$. This never occurs
  by Lemma~\ref{lemma:nonzero-implies-distinct} and strong boundary
  monotonicity.
\item $M(\rho_i^+)=M(\rho_{i+1}^-)$ and
  $\rho_i^+ \ne \rho_{i+1}^-$. This is ruled out by
  Lemma~\ref{lemma:collision-weakly-composable}.
\item $M(\rho_i^-) = M(\rho_{i+1}^+)$ and $\rho_i^- \ne
  \rho_{i+1}^+$.  In this case, $a(-\rho_i)a(-\rho_{i+1}) = 0$.
\item $(\rho_i,\rho_{i+1})$ are interleaved (in that order). 
  This degeneration does not occur in codimension one for embedded
  curves.  More precisely, then the sets $\{\rho_i\}$ and
  $\{\rho_{i+1}\}$ are not composable by
  Definition~\ref{def:Composable-sets}, so by
  Lemma~\ref{lemma:collision-is-composable} this degeneration cannot
  occur.
\item $(\rho_{i+1},\rho_i)$ are interleaved (in that order).  Then
  $a(-\rho_i)a(-\rho_{i+1}) = 0$ as there are double crossing strands.
\end{itemize}

The remaining possibilities are
\begin{enumerate}[label=(\ref*{case:dsquare-collide}\alph*),ref=\ref*{case:dsquare-collide}\alph*]
\item\label{case:dsquare-collide-join} $\rho_i^-=\rho_{i+1}^+$. These are exactly the moduli
  spaces which occur in Case~(\ref{case:dsquare-join}), for the
  factorization $a(-\rho_1)\cdots a(-(\rho_i\uplus\rho_{i+1}))\cdots a(-\rho_n)$
  of $a$. So, these two degenerations cancel.
\item\label{case:dsquare-collide-collide}
  $\{M(\rho_i^-),M(\rho_i^+)\}\cap\{M(\rho_{i+1}^-),M(\rho_{i+1}^+)\}=
  \emptyset$, with $\{\rho_i,\rho_{i+1}\}$ nested (in either order) or
  disjoint. In
  this case, the same degeneration also occurs
  for the factorization
  $a(-\rho_1)\cdots a(-\rho_{i+1})a(-\rho_i)\cdots a(-\rho_n)$
  of~$a$ (also in Case~(\ref{case:dsquare-collide-collide})). So,
  these two degenerations cancel.
\end{enumerate}

It follows that the coefficient of $a$ in $\partial a^B_{\x,\y}+\sum_{\{\w,B_1,B_2\,\mid\,B_1*B_2=B\}} a^{B_1}_{\x,\w}a^{B_2}_{\w,\y}$
vanishes. Summing over $B$ and noting that the choice of algebra element $a$ was arbitrary, Equation~\eqref{eq:typeDequation} follows.
\end{proof}

\begin{figure}
\includegraphics[scale=.83333]{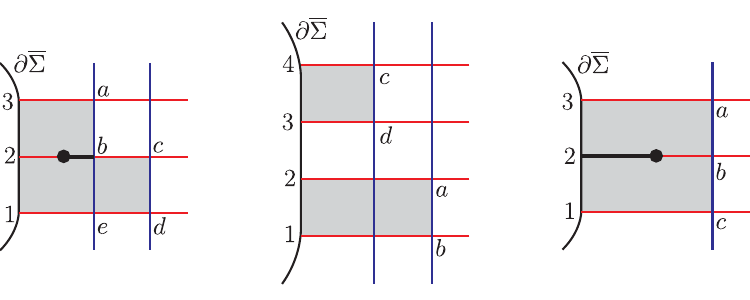}
\caption[Local illustrations of possible contributions to
    $\partial^2$ in $\CFDa$]{\textbf{Local illustrations of possible contributions to
    $\partial^2$ in the type $D$ module.} Left: a family of curves
  with one end a height $2$ comb and the other end a join curve
  end. Center: domain of a family of curves with a height $2$ comb at one end
  and a collision of two levels not resulting in any split curves at
  the other. Right: a family of curves with a height $2$ comb at
  one end and a collision of two levels, degenerating a split curve,
  at the other. (The thick lines indicate cuts; the boundary branch
  points at the ends of the cuts are shown as black dots.)}\label{fig:dsquare0D}
\end{figure}

The most interesting points in the proof of
Proposition~\ref{prop:typeD-d2} are seen in the
following examples, illustrated in Figure~\ref{fig:dsquare0D}; see also~\cite{LOT0}.  We put
the boundary of $\bSigma$ on the left to indicate visually that
the orientation of $\bdy\bSigma$ is reversed in the algebra.
\begin{example}\label{eg:dd1}On the left of Figure~\ref{fig:dsquare0D} is a piece of
  a diagram with four generators.  The complex $\CFDa$ of this piece
  satisfies the relations
  \begin{align*}
    \partial\{a,c\}&=\honestalg{1}{2}{}\{a,d\}+\honestalg{1}{3}{}\{c,e\}\\
    \partial\{c,e\}&=\{b,d\}\\
    \partial\{a,d\}&=\honestalg{2}{3}{}\{b,d\}.\\
  \end{align*} 
  The fact that $\partial^2=0$ follows from the relations in the algebra,
  more specifically
  \begin{equation*}
    \honestalg{1}{2}{}\cdot \honestalg{2}{3}{} = \honestalg{1}{3}{}.
  \end{equation*}
  This
  illustrates Case~(\ref{case:dsquare-two-lev}) in the proof of
  Proposition~\ref{prop:typeD-d2} (at the two two-story ends of the
  moduli space), as well as the
  cancellation of
  Case~(\ref{case:dsquare-join}) and Case~(\ref{case:dsquare-collide-join})
  (when the boundary branch point goes out to $\bdy\Sigma$). The
  reader may also want to compare with
  Figure~\ref{fig:degen_examples}.
\end{example}
\begin{example}\label{eg:dd2}
  In the center of Figure~\ref{fig:dsquare0D}, we have
  \begin{align*} 
    \partial\{a,c\}&=\honestalg{1}{2}{3}\{b,c\}
    +\honestalg{3}{4}{1}\{a,d\}\\
    \partial\{b,c\}&=\honestalg{3}{4}{2}\{b,d\}\\
    \partial\{a,d\}&=\honestalg{1}{2}{4}\{b,d\}.
  \end{align*}
  The fact that $\partial^2=0$ follows from the fact that
  $a\bigl(\reebchords{3}{4}\bigr)$ and
  $a\bigl(\reebchords{1}{2}\bigr)$
  commute in the algebra---or, more precisely, because
  \[
  \honestalg{1}{2}{3}\cdot\honestalg{3}{4}{2}
  =\honestalg{3}{4}{1}\cdot\honestalg{1}{2}{4}
  =\honestalg{1 & 3}{2 & 4}{}.
  \]
  This is the self-cancellation of
  Case~(\ref{case:dsquare-collide-collide}) of the proof of
  Proposition~\ref{prop:typeD-d2}.
\end{example}
\begin{example}\label{eg:dd3}
  On the right of Figure~\ref{fig:dsquare0D}, we have
  \begin{align*}
    \partial\{a\}&=\honestalg{2}{3}{1}\{b\}+\honestalg{1}{3}{2}\{c\}\\
    \partial\{b\}&=\honestalg{1}{2}{3}\{c\}.
  \end{align*}
  Here, $\partial^2=0$ because of a differential in the algebra:
  \[
  \partial\honestalg{1}{3}{2}=\honestalg{1 & 2}{2 & 3}{}
  =\honestalg{2}{3}{1}\cdot\honestalg{1}{2}{3}.
  \]
  This illustrates Case~(\ref{case:dsquare-split}) of the
  proof of Proposition~\ref{prop:typeD-d2}, or equivalently
  Lemma~\ref{lemma:third_term}.
\end{example}

The above examples are local: they are, of course, not the modules
associated to actual bordered Heegaard diagrams.  We refer the reader
wishing to see global examples to
Section~\ref{sec:surg-exact-triangle}, where we compute the
invariants of genus one handlebodies with certain borderings.

\section{Invariance}\label{sec:typeD-invariance}
\index{invariance!of $\CFDa$|(}%
\index{Heegaard moves!invariance under|see{invariance}}%
\begin{theorem}\label{thm:D-invariance}Up to chain homotopy equivalence, the
  differential module $\CFDa(\HD,\spinc)$ is independent of the choice
  of admissible, sufficiently generic almost complex structure and depends on
  the provincially admissible Heegaard diagram $\HD$ with $\PMC=-\partial \HD$
  only through its induced (equivalence class of)
  bordered three-manifold $(Y,\phi\co {-F(\PMC)}\to \partial Y)$.
  That is, if $\HD$ and $\HD'$ are provincially
  admissible bordered Heegaard diagrams for~$(Y, \phi\co {-F(\PMC)}\to
  \bdy Y)$; $J$ and $J'$ are admissible, sufficiently generic almost
  complex structures on $\Sigma\times[0,1]\times\RR$ and
  $\Sigma'\times[0,1]\times\RR$; 
  and $\spinc$ is a $\SpinC$ structure on $Y$\!, then
  $\CFDa(\HD,\spinc)$ (computed with respect to~$J$) and
  $\CFDa(\HD',\spinc)$ (computed with respect to~$J'$) are homotopy
  equivalent differential $\Alg(\PtdMatchCirc)$-modules.
\end{theorem}

The proof of Theorem~\ref{thm:D-invariance} is modeled on the
invariance proof for the closed case \cite{OS04:HolomorphicDisks}
(as modified for the cylindrical case
in~\cite{Lipshitz06:CylindricalHF}, say). In particular, the proof
proceeds by showing that each of the three kinds of Heegaard moves
from Proposition~\ref{prop:heegaard-moves} induces a homotopy
equivalence. (We also use
Proposition~\ref{prop:admis-achieve-maintain} to ensure that all of
the intermediate diagrams are provincially admissible.)

There are two essentially new issues. The first, and easier, is that
we must keep track of the algebra elements which occur as coefficients
in the maps. We will illustrate how this is done by defining the chain
map induced by deformation of the almost complex structure. The
second, more serious, issue is that the traditional proof of
handleslide invariance uses ``triangle maps,'' the general definitions
of which are subtle for bordered Heegaard diagrams. In
the special
case of handleslides the difficulties which come up in triangle maps
in general can be circumvented, however.
We illustrate this in a proof of invariance under handleslides
which carry one $\alpha$-arc over an $\alpha$-circle (the hardest
case). These two issues having been explained, the reader should have
no difficulty adapting the rest of the invariance proof to our
setting.  For instance, the proof of isotopy invariance is similar to
the proof of invariance under change of complex structure, but
notationally slightly more complicated.

\subsection{The chain map for change of complex structure}\label{sec:CFD-cx-str-change}
Recall that the moduli spaces $\cM^B(\x,\y;\Source;\vec{P})$ from
Chapter~\ref{chap:structure-moduli} depended on a choice of a generic,
admissible almost complex structure $J$ on
$\Sigma\times[0,1]\times\RR$, and consequently the differential module
$\CFDa(\HD,\spinc)$ also depends on this choice; to make this choice
explicit, we will write $\CFDa(\HD,\spinc; J)$ in this section.\glsadd{CFDJ}
As is usual in Floer homology theories,
one proves independence of the almost complex structure by
constructing ``continuation maps'', and proving that they are chain
homotopy equivalences; see, for example,~\cite[Theorem 4]{FloerContinuation}. 
More specifically, let $J_0$ and $J_1$ be two
generic, admissible almost complex structures on
$\Sigma\times[0,1]\times\RR$, and let $\partial_0$ and
$\partial_1$ denote the boundary operators on $\CFDa(\HD,\spinc; J_0)$ and
$\CFDa(\HD,\spinc; J_1)$, respectively. Call a (smooth)
path $\{J_r\mid r\in[0,1]\}$ of almost complex structures from $J_0$
to $J_1$ (constant near $0$ and $1$) \emph{admissible}
\index{complex structure!admissible!path of}
if each $J_r$ is admissible. To a
generic (in a sense to be made precise presently) admissible path~$J_r$
between $J_0$ and $J_1$ we will associate a ``continuation map''
\[
\gls*{FJr}
\co\CFDa(\HD,\spinc;J_0)\to
\CFDa(\HD,\spinc;J_1).
\]

To the path $J_r$ we can associate a single almost complex structure
$J$ on $\Sigma\times[0,1]\times\RR$ by
\begin{equation}\label{eq:non-cylind-J}
J|_{(x,s,t)}=\begin{cases}
J_1|_{(x,s,t)} & \text{if }t\geq 1\\
J_t|_{(x,s,t)} & \text{if }0\leq t\leq
1\\
J_0|_{(x,s,t)} & \text{if } t\leq 0.
\end{cases}
\end{equation}
(Here $J|_{(x,s,t)}$ means the map $J$ induces on the tangent space at
$x\in\Sigma$, $s\in[0,1]$, and $t\in\RR$.)

For $B\in\pi_2(\x,\y)$, let 
$\gls*{ModCxStrSource}$
denote the moduli spaces of holomorphic
curves defined in Definition~\ref{def:tcM} with respect to the almost
complex structure~$J$. (Note that here there is no $\RR$-action,
since $J$ is not $\RR$-invariant.) For generic~$J$, these moduli spaces are all
transversely cut out, as in
Proposition~\ref{prop:transversality}. This is the genericity
requirement we impose on $\{J_r\}$.
\index{generic|see{complex structure, sufficiently generic}}
\index{complex structure!sufficiently generic!for path}%
\index{sufficiently generic!path of almost complex structures}%
For $\vec\rhos$ a sequence of non-empty sets of Reeb chords so that
$(\x,\vec\rhos)$ is strongly boundary monotone and $B$ is compatible
with $\vec\rhos$, define the embedded moduli space by
\begin{equation}\label{eq:non-cyclind-embedded-M}
\gls*{ModCxStrEmb}
\coloneqq
    \!\!\bigcup_
       {\chi(\Source)=\chi_{\emb}(B,\vec{\rhos})}\!\!
  \cM^B(\x, \y;
  \Source; \vec{P};J).
\end{equation}
Since we do not divide by $\RR$-translation, we will
primarily be interested in
index~$0$ moduli spaces, not the index~$1$ ones.

Then define the map $F^{J_r}$ on $\x \in \Gen(\HD,\s)$ by
\begin{equation}
  \label{eq:ContinuationMap}
\gls*{FJr}
(\x)\coloneqq\sum_{\y\in\S(\HD,\spinc)}\,\sum_{B\in\pi_2(\x,\y)}\,
\sum_{\{\vec{\rho}\,\mid\,\ind(B,\vec\rho) = 0\}}\!
  \#\left(\Mod^B(\x,\y;\vec{\rho};J)\right)a(-\vec\rho) \y,
\end{equation}
and extend $F^{J_r}$ to all of $\CFDa(\HD,\spinc;J_0)$ by
$F^{J_r}(a\x)=aF^{J_r}(\x)$.  Compactness and provincial admissibility
imply that this sum is finite, as in Lemma~\ref{lem:finite-typeD}.

It is immediate from  Lemma~\ref{lem:SpinCStructures} that 
$F^{J_r}(\x)$ respects the decomposition into $\SpinC$ structures; see also 
Remark~\ref{rmk:SpinCSplitting}.

The next task is to show that $F^{J_r}$ is a chain map. We will use
the following analogue of
Theorem~\ref{thm:master_equation} for
almost complex structures which are not translation invariant.  The
statement is more general than we need for present purposes; we
will also use it to prove invariance of the type $A$ module in
Section~\ref{sec:A-invariance}.

\begin{proposition}\label{prop:non-cylindrical-master}
  Suppose that $(\x,\vec{\rhos})$ satisfies
  the strong boundary monotonicity condition. Fix $\y\in\Gen(\HD)$,
  $B\in\pi_2(\x,\y)$, $\Source$, and $\vec{P}$ such that
  $[\vec{P}] = \vec\rhos$ and $\ind(B,\Source,P)=1$. Consider the
  parameterized moduli space $\cM^B(\x,\y;\Source;\vec{P};J)$. Then the total
  number of the following ends of $\cM^B(\x,\y;\Source;\vec{P};J)$ is even:
  \begin{enumerate}
  \item two-story ends, with either
    \begin{enumerate}
    \item a $J_0$-holomorphic curve followed by
      a $J$-holomorphic curve or
    \item a $J$-holomorphic curve  followed by a
      $J_1$-holomorphic curve;
    \end{enumerate}
  \item join curve ends;
  \item odd shuffle curve ends; and
  \item collision of levels $P_i$ and $P_{i+1}$ from $\vec{P}$,
    where $[P_i]$ and $[P_{i+1}]$ are weakly composable.
  \end{enumerate}
\end{proposition}

(The definition of two-story ends, join curve ends,
odd or even shuffle curve ends, and collisions of
levels are given in Definition~\ref{def:ends-moduli}.  As indicated,
two-story ends come in two types: elements of
$\cM^{B_1}(\x,\w;\SourceSub{1};\vec{P}_1;J_0) \times
\cM^{B_2}(\w,\y;\SourceSub{2};J)$ and elements of
$\cM^{B_1}(\x,\w;\SourceSub{1};\vec{P}_1;J) \times
\cM^{B_2}(\w,\y;\SourceSub{2};J_1)$.)

\begin{proof}
The proof is essentially the same as the proof of
Theorem~\ref{thm:master_equation}.
\end{proof}

\begin{proposition}\label{prop:dffd}
  The map $F^{J_r}$ is a chain map, i.e.,
  $\partial_1\circ F^{J_r}+F^{J_r}\circ \partial_0=0$.
\end{proposition}
\begin{proof}
  The proof is essentially the same as the proof of
  Proposition~\ref{prop:typeD-d2}. First, observe that
\[
\partial_1(F^{J_r}(a\x))+F^{J_r}(\partial_0(a\x))=(\partial a)F^{J_r}(\x)+a\partial_1(F^{J_r}(\x))+aF^{J_r}(\partial_0(\x))+(\partial a)F^{J_r}(\x).
\]
So, it suffices to prove that $\partial_1(F^{J_r}(\x))+F^{J_r}(\partial_0(\x))=0$ for any generator $\x\in\S(\HD,\spinc)$.

Now, as in the proof of Proposition~\ref{prop:typeD-d2}, for given
$\x$, $\y$, $B$ and $\vec{\rho}=(\rho_1,\dots,\rho_n)$ with
$\ind(B,\vec{\rho}) = 1$,
Proposition~\ref{prop:non-cylindrical-master} implies that the sum of
the following terms is zero:
\begin{enumerate}
\item\label{case:dffd-split1}
  The number of height $2$ combs of the form
\[
  \Mod^{B_1}(\x,\w; \vec\rho_1; J_0)\times
    \Mod^{B_2}(\w,\y; \vec\rho_2; J)
\]
where $B = B_1 * B_2$ and $\vec\rho = (\vec\rho_1, \vec\rho_2)$.
\item\label{case:dffd-split2}
  The number of height $2$ combs of the form
\[
  \Mod^{B_1}(\x,\w; \vec\rho_1; J)\times
    \Mod^{B_2}(\w,\y; \vec\rho_2; J_1)
\]
where $B = B_1 * B_2$ and $\vec\rho = (\vec\rho_1, \vec\rho_2)$.
\item\label{case:dffd-join}
  The number of join curve ends, i.e., the number of elements of
\[
\Mod^B(\x,\y;(\rho_1,\dots,\rho_{i-1},\{\rho_j,\rho_k\},\rho_{i+1},\dots,\rho_n);
J)
\]
where $\rho_i = \rho_j \uplus \rho_k$.
\item\label{case:dffd-split}
  The number of split curve ends, i.e., the number of elements of
\[
\Mod^B(\x,\y;(\rho_1,\dots,\rho_{i-1},\rho_i\uplus\rho_{i+1},\rho_{i+2},\dots,\rho_n);
J)
\]
 where $\rho_i^+=\rho_{i+1}^-$.
\item\label{case:dffd-collide} The number of other collisions of levels,
  i.e., the number of elements of
\[
\Mod^B(\x,\y;(\rho_1,\dots,\rho_{i-1},\{\rho_{i},\rho_{i+1}\},\rho_{i+2}\dots,\rho_n); J)
\]
where $\rho_i^+\ne\rho_{i+1}^-$.
\end{enumerate}
The first term corresponds to $F^{J_t}(\partial_0(\x))$. The second
corresponds to the terms in $\partial_1(F^{J_t}(\x))$ coming from
applying $\partial_1$ to the resulting generators $\w\in\S(\HD,\spinc)$. The
fourth term corresponds to the terms in $\partial_1(F^{J_t}(\x))$
coming from the differential on the algebra $\Alg$. The third and
fifth terms cancel in pairs, as in the proof of
Proposition~\ref{prop:typeD-d2}.
\end{proof}

One can similarly use a path of almost complex structures from $J_1$ to
$J_0$ to define a chain map
$\CFDa(\HD,\spinc;J_1)\to \CFDa(\HD,\spinc;J_0)$. The
proof that these maps are mutually inverse chain homotopy equivalences
is obtained by adapting the standard arguments in analogous
manners to the above, and we leave this to the interested reader.
(Some more details are given in the type~$A$ case, in
Proposition~\ref{prop:change-cx-htpy-inv-A}.)

\subsection{Handlesliding an \textalt{$\alpha$}{alpha}-arc over an \textalt{$\alpha$}{alpha}-circle}\label{sec:CFD-handleslide}
Let $(\Sigma,\alphas,\betas,z)$ denote a bordered
Heegaard diagram, which we assume for
convenience to be admissible; this can be arranged by performing
isotopies of the $\alpha$-curves. Write
$\alphas=\{\alpha_1^a,\dots,\alpha_{2k}^a,\alpha_1^c,\dots,\alpha_{g-k}^c\}$. Provisionally,
let $\alpha_1^{a,H}$ denote the result of performing a handleslide of
$\alpha_1^a$ over $\alpha_1^c$, and
$\alphas^H=\{\alpha_1^{a,H},\alpha_2^a,\dots,\alpha_{g-k}^c\}$. Assume
that the Heegaard diagram $(\Sigma,\alphas^H,\betas,z)$ is also
admissible. We want to show that
$\CFDa(\Sigma,\alphas,\betas,z)$ and
$\CFDa(\Sigma,\alphas^H,\betas,z)$ are chain homotopy equivalent.

\subsubsection{Moduli spaces of triangles: basic definitions}
\begin{figure}
  \includegraphics[scale=.83333]{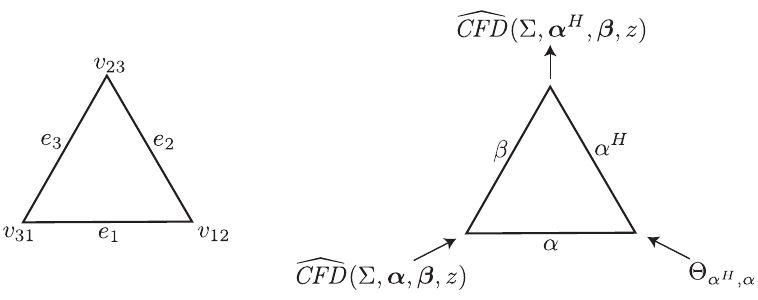}
  \caption[Labeling of triangle]{\textbf{The three-punctured disk (triangle) $\Delta$ (left), and the labeling inducing
    $F_{\alphas,\alphas^H,\betas}$ (right).} In counter-clockwise
  cyclic order, the edges are labeled $\alpha,\alpha^H,\beta$.}\label{fig:Triangle-Delta}
\end{figure}

As is traditional, we will produce a chain map
$\gls*{FaaHb}
\co \CFDa(\Sigma,\alphas,\betas,z)\to
\CFDa(\Sigma,\alphas^H,\betas,z)$ by counting holomorphic triangles
\index{holomorphic triangle}\index{triangle, holomorphic}%
in
$(\Sigma,\alphas,\alphas^H,\betas,z)$. That is, let 
$\gls*{TriangleDelta}$
denote a disk with three boundary punctures. Label the edges of $\Delta$ by
$e_1$, $e_2$ and $e_3$, counterclockwise, and let $v_{ij}$ denote the
\glsadd{TriangleEdges}\glsadd{TriangleVertices}%
puncture between $e_i$ and $e_j$; see
Figure~\ref{fig:Triangle-Delta}. To define
$F_{\alphas,\alphas^H,\betas}$ we will count holomorphic maps
\[
u\co (T,\bdy T)\to\left(\Sigma\times\Delta,(\alphas\times
e_1)\cup(\alphas^H\times e_2)\cup(\betas\times e_3)\right)
\]
which do not cover~$z$. 

\begin{figure}
  \includegraphics[scale=.83333]{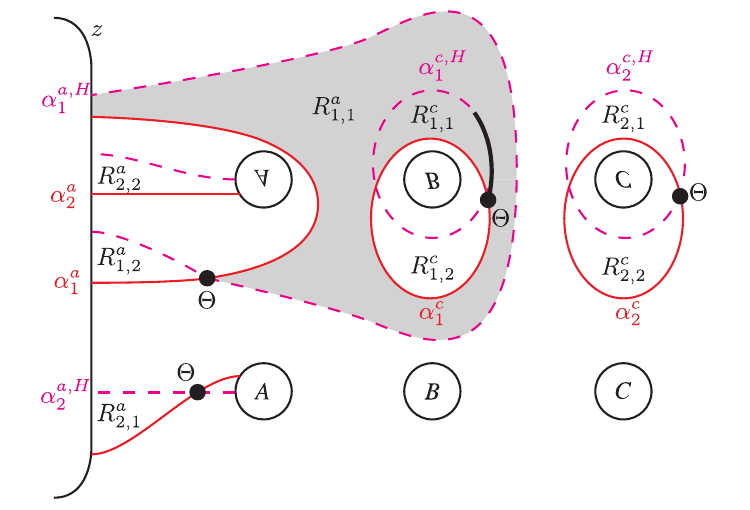}
  \caption[Curves $\alphas$ and $\alphas^H$, before and after a handleslide]{\textbf{The curves
      $\alphas$ and $\alphas^H$ in $\Sigma$.} The $\alpha$-curves are
    solid while the $\alpha^H$-curves are dashed. The intersection points
    which belong to generators $\Theta_o$ are marked.
    The annulus $R_{1,1}^a+R_{1,1}^c$
    is shaded. The branch cut corresponding to a
    possible holomorphic curve is also indicated by a thick
    arc.}\label{fig:alphasalphasH}
  \label{fig:AlphaAlphaH}
\end{figure}

The 
$\gls*{alphasH}$
we use are actually a slight modification of the
above, as shown in Figure~\ref{fig:alphasalphasH}.
Obtain
$\gls*{alphaicH}$
from $\alpha_i^c$ by performing
a small Hamiltonian perturbation, so that $\alpha_i^{c,H}$
intersects $\alpha_i^c$  transversely in two points, is
disjoint from all other $\alpha$-curves, and intersects
$\beta_j$ transversally close to $\beta_j\cap\alpha_i^c$. For
$i=2,\dots,2k$, 
$\gls*{alphaiaH}$
is an isotopic translate of
$\alpha_i^a$, intersecting $\alpha_i^a$ transversely in a single point, and such
that there are two short Reeb chords in $\bdy\widebar{\Sigma}$
running from $\alpha_i^{a,H}$ to $\alpha_i^a$ (rather than the other way). The isotopy is chosen
small enough that $\alpha_i^{a,H}$ is disjoint from all other
$\alpha$-curves and intersects $\beta_j$ transversally
close to $\beta_j\cap\alpha_i^a$. Obtain $\alpha_1^{a,H}$ by
handlesliding $\alpha_1^a$ over $\alpha_1^c$ and then performing
an isotopy so that $\alpha_1^{a,H}$ intersects $\alpha_1^a$ transversely in a
single point, and such that there are two short Reeb chords in
$\bdy\widebar{\Sigma}$ running from $\alpha_1^{a,H}$
to~$\alpha_1^a$. We also arrange that $\alpha_1^{a,H}$ is disjoint
from
all other $\alpha$-curves and intersects the $\beta_j$
transversally at points close to
$\beta_j\cap(\alpha_1^a\cup\alpha_1^c)$.
Also define
$\gls*{CircPtsH}$
to be the set
$\alphas^H\cap\bdy\bSigma$.  There is a natural bijection between
$\CircPts$ and~$\CircPts^H$.

Most of the notions from Chapter~\ref{chap:structure-moduli} generalize
in straightforward manners to the setting of triangles. 
Following Definition~\ref{def:admissible_J}, we
fix an almost complex structure $J$ on $\Sigma\times\Delta$ such that:
\begin{enumerate}
\item The projection map $\pi_\Delta\co\Sigma\times\Delta\to\Delta$ is
  $(J,j_\Delta)$-holomorphic.
\item The fibers of the projection map
  $\pi_\Sigma\co\Sigma\times\Delta\to\Sigma$ are $J$-holomorphic.
\item The almost complex structure is split near $p\times\Delta$.  (Recall
  that $p$ is the puncture of~$\Sigma$).
\item\label{item:tri-J-local-translation-inv} For each puncture
  $v_{i,j}$ of $\Delta$, choose a neighborhood $U_{i,j}$ which is
  biholomorphic to $[0,1]\times[0,\infty)$, and such a
  biholomorphism. This identification induces an action of $\RR_+$ on $\Sigma\times
  U_{i,j}$. We require that these actions be $J$-holomorphic.
\item\label{item:tri-J-global-translation-inv} Identify $\Delta$ with
  $([0,1]\times\RR)\setminus \{(1,0)\}$ so that $v_{12}$ is identified
  with $(1,0)$ (see Figure~\ref{fig:TriangleIsStrip}), and let
  $\tau_R\co ([0,1]\times\RR)\setminus \{(1,0),(1,-R)\}\to
  ([0,1]\times\RR)\setminus \{(1,R),(1,0)\}$ denote translation in the
  $\RR$-factor. We require that the map $\tau_R$ is
  $J$-holomorphic for each $R$.
\end{enumerate}
\index{complex structure!admissible!on $\Sigma\times\Delta$}%
We call such a $J$ \emph{admissible}. Notice that the restrictions of
an admissible almost complex structure $J$ on $\Sigma\times\Delta$ to
$\Sigma\times U_{i,j}$ induce three admissible structures $J_{i,j}$ on
$\Sigma\times[0,1]\times\RR$. Condition~\ref{item:tri-J-global-translation-inv}
can be thought of as a global version of
Condition~\ref{item:tri-J-local-translation-inv}, and in particular
implies Condition~\ref{item:tri-J-local-translation-inv} at the
punctures $v_{3,1}$ and $v_{2,3}$; together,
Conditions~\ref{item:tri-J-local-translation-inv} and
\ref{item:tri-J-global-translation-inv} imply that $J_{1,2}$ is split.
Condition~\ref{item:tri-J-global-translation-inv} is only needed for
the computation of the embedded index in
Proposition~\ref{prop:tri-emb-ind-general}.

Similarly, we can define decorated sources for maps to
\index{source!decorated!for map to $\Sigma\times\Delta$}%
$\Sigma\times\Delta$: a decorated source~$\gls*{TSource}$
is an
oriented surfaces with boundary and punctures, together with a
labeling of each puncture by an element of
$\{v_{1,2},v_{2,3},v_{3,1},e\}$ and a further labeling of each $e$
puncture either by a Reeb chord with endpoints in $\CircPts$ or
by a Reeb chord with endpoints in $\CircPts^H$.

Given a decorated source $\TSource$ we consider maps
\[
u\co (T,\bdy T)\to\left(\Sigma\times\Delta,(\alphas\times
e_1)\cup(\alphas^H\times e_2)\cup(\betas\times e_3)\right)
\]
satisfying the obvious analogues of
Conditions~(\ref{item:moduliFirst})--(\ref{item:moduliStrongBndyMntncty})
from Section~\ref{sec:curves-in-sigma}. Such a holomorphic curve is
asymptotic to a generator~$\x$ for $(\Sigma,\alphas,\betas,z)$
at~$v_{3,1}$ and a generator~$\y$ for $(\Sigma,\alphas^H,\betas,z)$
at~$v_{2,3}$. The
asymptotics at $v_{1,2}$, however, are more subtle. At each puncture of
$\TSource$ mapped to $v_{1,2}$, $u$ is either asymptotic to a point in
$\alphas^H\cap\alphas$ (in the interior of $\Sigma$), or is asymptotic
to a Reeb chord in $\bdy\widebar{\Sigma}$ running from $\CircPts^H$ to~$\CircPts$.  
In defining the map $F_{\alphas,\alphas^H,\betas}$ we
will use curves with only the former kind of asymptotics, but in
proving that $F_{\alphas,\alphas^H,\betas}$ is a chain map we will
need to consider curves with the latter asymptotics.
Such asymptotics present
certain technical issues, as we will discuss below.

For each $k$-element subset $o$ of $\{1,\dots,2k\}$ there is a
distinguished $g$-element subset of $\alphas^H\cap\alphas$,
denoted~$\gls*{Thetao}$.
The set $\Theta_o$ contains the unique intersection point
between $\alpha_i^a$ and $\alpha_i^{a,H}$ for $i\in o$, as well as the
higher-graded intersection point
between $\alpha_i^{c,H}$ and $\alpha_i^{c}$ for each
$i=1,\dots,g-k$. 
(The higher-graded intersection point~$x$
between $\alpha_i^{c,H}$ and $\alpha_i^c$ is the one such that there are two holomorphic maps
  $([0,1]\times\RR,\{1\}\times\RR,\{0\}\times\RR)\to
  (\Sigma,\alphas^H,\alphas)$ mapping $-\infty$ to~$x$.)
We are thinking of $\Theta_o$ as a generator
for counts of holomorphic curves in $(\Sigma, \alphas^H, \alphas, z)$,
although we will not define an object $\CFDa(\Sigma,\alphas^H,
\alphas, z)$.  

Like curves in $\Sigma\times[0,1]\times\RR$, curves
in $\Sigma\times\Delta$ carry (relative) homology classes, in a manner
analogous to Definition~\ref{def:pi2}. Let
$\gls*{pi2xyThetao}$
denote the set of homology classes
connecting $\x$, $\y$ and $\Theta_o$ (not covering~$z$, as before).
For a homology class $B$ in
$\pi_2(\x,\y,\Theta_o)$ we can consider its \emph{domain},
\index{domain!of triangle}%
the
multiplicities with which $\pi_\Sigma(B)$ covers each
component of $\Sigma\setminus(\alphas \cup \alphas^H \cup \betas)$. 
The map from homology classes to domains is injective.  Let
$\gls*{ModTriSource}$
denote the moduli space of
holomorphic curves asymptotic to $\x$ at $v_{3,1}$, $\y$ at $v_{2,3}$
and $\Theta_o$ at $v_{1,2}$, with decorated source $\TSource$ and in
the homology class~$B$.

Along $e_1$ (respectively $e_2$) a curve
$u\in\cM^B(\x,\y,\Theta_o;\TSource)$ is asymptotic to a collection of
Reeb chords in
$(\bdy\bSigma,\CircPts)$ (respectively $(\bdy\bSigma,\CircPts^H)$). 
The
orientation on $e_1$ (respectively $e_2$) induces a partial
order of these collections of Reeb chords. Thus, given ordered
partitions $\vec{P}$ and $\vec{P}^H$ of Reeb chords in 
$(\bdy\bSigma,\CircPts)$ and $(\bdy\bSigma,\CircPts^H)$
respectively, we may consider the subspace 
\[
\gls*{ModTriSourcePP}\subset \cM^B(\x,\y,\Theta_o;\TSource)
\]
consisting of those holomorphic curves so that the induced ordered
partitions of the Reeb chords along $e_1$ and $e_2$ are given by
$\vec{P}$ and $\vec{P}^H$ respectively.

\begin{lemma}\label{lem:tri-determine-o}
  If the moduli space
  $\cM^B(\x,\y,\Theta_o;\TSource;\vec{P},\vec{P}^H)$ 
  is non-empty then, in the notation of
  Definition~\ref{def:strong-monotonicity-P}, $o=o(\x,[\vec{P}])$.
\end{lemma}
\begin{proof}
  This is immediate from the definitions.
\end{proof}

  The arguments in Proposition~\ref{prop:compactness}---essentially,
  considering the maps $\pi_\Sigma
  \circ u$ and $\pi_\Delta \circ u$ separately---show that the moduli
  spaces $\cM^B(\x,\y,\Theta_o;\TSource;\vec P,
  \vec P^H)$ admit natural compactifications $\oocM^B(\x,\y,\Theta_o;\TSource;\vec P,
  \vec P^H)$.  Other than
  degenerations at $v_{12}$, which we will discuss presently, the
  limit curves are natural analogues of the holomorphic combs
  introduced in Section~\ref{sec:combs-compact}. As for bigons in
  Chapter~\ref{chap:structure-moduli}, we let 
  $\ocM^B(\x,\y,\Theta_o;\TSource;\vec P, \vec P^H)$
  denote the closure of $\cM^B(\x,\y,\Theta_o;\TSource;\vec P,
  \vec P^H)$ in $\oocM^B(\x,\y,\Theta_o;\TSource;\vec P,
  \vec P^H)$.
  \glsadd{ModTriSourcePPCpctCpct}\glsadd{ModTriSourcePPCpct}

\subsubsection{Homology classes and expected dimensions of triangles}
The expected dimension of $\cM^B(\x,\y,\Theta_o;\TSource;\vec{P},\vec{P}^H)$ is
\begin{equation}\label{eq:TriangleIndex}
\gls*{IndTri}
=\frac{g}{2}-\chi(T)+2e(B)+|P|+|P^H|.
\end{equation}
(The proof is similar to the proof of Proposition~\ref{Prop:Index}.
See \cite[p.\ 1018]{Lipshitz06:CylindricalHF} for the closed case.)
Like Formula~\eqref{eq:Index}, Formula~\eqref{eq:TriangleIndex}
depends on the topology of $\TSource$, and not just on the homology
class and asymptotics. However, as was the case for embedded curves in
$\Sigma\times[0,1]\times\RR$, the Euler characteristic of an embedded
curve in $\Sigma\times\Delta$ is determined by the homology class and
asymptotics.  To state the formula in a convenient form, first note
that for any generator $\x\in\S(\Sigma,\alphas,\betas)$, there is a
nearby generator 
$\gls*{xprime}
\in\S(\Sigma,\alphas^H,\betas)$ (but not vice
versa), as well as a canonical homology class~$\gls*{smallTriClass}$
in
$\pi_2(\x,\x',\Theta_{o(\x)})$, the union of a number of small
triangles (contained in the isotopy region) and possibly an annulus
with boundary on $\alpha_1^a$, $\alpha_1^{a,H}$ and $\alpha_1^c$.

\begin{lemma}\label{lem:triple-pi2-decomp}
  Any homology class $B\in\pi_2(\x,\y,\Theta_o)$ can be
  written uniquely as 
  \[
  (T_{\x} *_{12} B_{\alphas^H,\alphas}) *_{23} B_{\alphas^H,\betas}
  \]
  for some $B_{\alphas^H,\alphas}\in\pi_2(\Theta_{o(\x)},\Theta_o)$ and
  $B_{\alphas^H,\betas}\in\pi_2(\x',\y)$.
\end{lemma}
\glsadd{starij}%
Here $*_{12}$ is the natural operation of attaching a homology class in
$(\Sigma,\alphas^H,\alphas,z)$ to a homology class in
$(\Sigma,\alphas,\alphas^H,\betas,z)$ at the corner $v_{12}$
of~$\Delta$, and likewise $*_{23}$ is the operation of attaching a
homology class in $(\Sigma,\alphas^H,\betas,z)$ at the corner $v_{23}$
of~$\Delta$.
\begin{proof}
  The space
  $\pi_2(\x,\y,\Theta_o)$ has a free and transitive action by the space of
\index{triply-periodic domain}\index{domain!triply-periodic}%
  {\em triply-periodic domains}, the space of $2$-chains in $\Sigma$ which
  have local degree zero at~$z$ and whose boundary is a formal linear
  combination of elements of $\alphas\cup\alphas^H\cup \betas$.  
  The space of triply-periodic domains in
  $(\Sigma,\alphas,\alphas^H,\betas,z)$ (and therefore
  $\pi_2(\x,\y,\Theta_o)$) is isomorphic to the kernel of the
  map
  \[
  H_1(\betas)\oplus H_1(\alphas/\bdy\alphas)\oplus
  H_1(\alphas^{H}/\bdy\alphas^{H})\to H_1(\bSigma/\bdy\bSigma).
  \]
  Since $H_1(\alphas/\bdy\alphas)$ and $H_1(\alphas^H/\bdy\alphas^H)$ are the same subspace of $H_1(\bSigma/\bdy\bSigma)$, this kernel is isomorphic to 
  \begin{multline*}
    \ker\Bigl(H_1(\betas)\oplus H_1(\alphas/\bdy\alphas)\to
    H_1(\bSigma/\bdy\bSigma)\Bigr)\\ 
    \oplus 
    \ker\Bigl(H_1(\alphas/\bdy\alphas)\oplus
    H_1(\alphas^H/\bdy\alphas^H)\to H_1(\bSigma/\bdy\bSigma)\Bigr) \\ \cong \pi_2(\x,\x)\times\pi_2(\Theta_o,\Theta_o).
  \end{multline*}
  Thus, the action by $\pi_2(\x,\x)\times\pi_2(\Theta_o,\Theta_o)$ on 
  the set of triply-periodic domains is free and
  transitive. The result follows.
\end{proof}

For $B\in\pi_2(\x,\y,\Theta_o)$, let
$\gls*{bofB}
\in\pi_2(\x',\y)$ be the domain
$B_{\alphas^H,\betas}\in\pi{_2(\x',\y)}$ whose
existence is guaranteed
by Lemma~\ref{lem:triple-pi2-decomp}.

Next, we prove an analogue of Sarkar's index theorem for
triangles~\cite{Sarkar11:IndexTriangles}. Given a domain
$B\in\pi_2(\x,\y,\Theta_o)$, let $\bdy^{\alpha}(B)$ (respectively
$\bdy^{\alpha^H}(B)$) denote the part of $\bdy B$ contained in the
$\alpha$-curves (respectively $\alpha^H$-curves), oriented as $\bdy B$.
Let
$\bdy^{\alpha^H}(B)\cdot\bdy^\alpha(B)$ denote the jittered
intersection number of $\bdy^{\alpha^H}(B)$ and $\bdy^\alpha(B)$, as
in Section~\ref{sec:index-additive}. As for bigons in
Section~\ref{sec:expected-dimensions}, $n_\x(B)$ denotes the point
measure of $B$ at $\x$, and $e(B)$ denotes the Euler measure of
$B$. Also analogous to the bigon case, define
\begin{align*}
\iota(\vec\rhos, & \vec\rhos^H) \coloneqq \iota(\vec\rhos) + \iota(\vec\rhos^H)
  + L(\vec\rhos, \vec\rhos^H)\\
  &=-\sum_i\frac{\abs{\rhos_i}}{2}-\sum_i\frac{\abs{\rhos_i^H}}{2}
   -\sum_i\sum_{\{\rho_a,\rho_b\}\subset
    \rhos_i}\!\!|L(\rho_a,\rho_b)|
  -\sum_i\sum_{\{\rho_a,\rho_b\}\subset \rhos^H_i}\!\!|L(\rho_a,\rho_b)|\\[5pt]
  &\qquad+ \!\!\!\sum_{\substack{i < j\\
      \rho_i\in\rhos_i,\rho_j\in\rhos_j}}\!\!\!\!\!L(\rho_i,\rho_j)+
  \!\!\!\sum_{\substack{i < j\\
      \rho_i\in\rhos^H_i,\rho_j\in\rhos^H_j}}\!\!\!\!\!L(\rho_i,\rho_j)+
  \!\!\sum_{\rho_i\in\rhos_i,\rho_j\in\rhos^H_j}\!\!\!L(\rho_i,\rho_j).
\end{align*}
\begin{proposition}\label{prop:tri-emb-ind-general}
  If $\cM^B(\x,\y,\Theta_o;\TSource;\vec P, \vec P^H)$ contains an embedded
  holomorphic curve $u$ then the expected dimension of $\cM^B$ near $u$ is given
  by 
  \begin{equation}\label{eq:tri-emb-ind}
  \ind(u)=e(B)+n_\x(B)+n_\y(B)+\bdy^{\alpha^H}(B)\cdot\bdy^\alpha(B)-g/2+\iota([\vec
  P],[\vec P^H]).
  \end{equation}
\end{proposition}
\begin{proof}
  \begin{figure}
    \centering
    \includegraphics[scale=.83333]{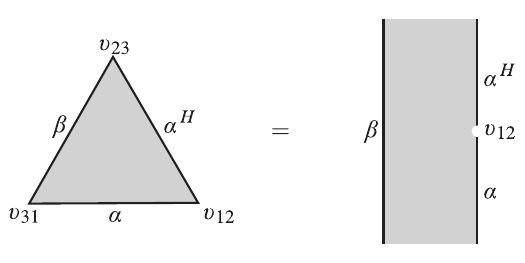}
    \caption[A triangle is a strip with a boundary puncture]{\textbf{A triangle is biholomorphic to a strip with a boundary puncture.} The edges are labeled by the corresponding arcs.}
    \label{fig:TriangleIsStrip}
  \end{figure}
  The proof is similar to the proof of
  Proposition~\ref{prop:asympt_gives_chi}. Identifying the triangle
  $\Delta$ with $[0,1]\times\RR\setminus \{(1,0)\}$, so that $v_{12}$
  is identified with $(1,0)$, we can view $u$
  as a map to $\Sigma\times[0,1]\times\RR$. (See
  Figure~\ref{fig:TriangleIsStrip}.)
  Let $S$ denote the source of $u$. Then
  \[
  \chi(S)=e(B)-\br(\pi_\Sigma\circ u)+\frac{3g}{4}
    +\sum_i\abs{P_i}/2 +\sum_i \abs{P_i^H}/2,
  \]
  where $\br$ denotes the order of branching and $\frac{3g}{4}$ comes
  from the corners. Let $\tau_r\co \Sigma\times[0,1]\times\RR\to
  \Sigma\times[0,1]\times\RR$ denote translation in the $\RR$ direction. Then for $\epsilon$
  sufficiently small,
  \[
   \br(\pi_\Sigma\circ u)=u\cdot \tau_\epsilon(u)-\sum_i\sum_{\{\rho_a,\rho_b\}\subset
     [P_i]}\!\!|L(\rho_a,\rho_b)|
   -\sum_i\sum_{\{\rho_a,\rho_b\}\subset [P^H_i]}\!\!|L(\rho_a,\rho_b)|.
  \]

  Translating farther, we see that for $R$ sufficiently large,
  \begin{multline}\label{eq:tri-index-intermediate-step}
    u\cdot \tau_\epsilon(u)=u\cdot
    \tau_R(u)+\bdy^{\alpha^H}(B)\cdot\bdy^{\alpha}(B)+\frac{g}{4}+\!\!\sum_{\substack{i
      < j\\ \rho_i\in[P_i],\rho_j\in[P_j]}}\!\!\!\!L(\rho_i,\rho_j)\\[5pt]
  + \!\!\!\sum_{\substack{i <
      j\\ \rho_i\in[P^H_i],\rho_j\in[P^H_j]}}\!\!\!\!L(\rho_i,\rho_j)+
    \!\!\sum_{\rho_i\in[P_i],\rho_j\in[P^H_j]}\!\!\!\!L(\rho_i,\rho_j).
  \end{multline}
  The $L(\cdot,\cdot)$ terms in
  Equation~\eqref{eq:tri-index-intermediate-step} arise exactly as in
  the proof of
  Proposition~\ref{prop:asympt_gives_chi}; note that this uses
  Condition~(\ref{item:tri-J-global-translation-inv}), so that
  $\tau_r(u)$ is holomorphic for any $r\in\RR$. The term
  $\bdy^{\alpha^H}(B)\cdot\bdy^{\alpha}(B)+g/4$ comes from the
  number of intersection points that appear (or disappear) at the
  boundary, as follows.
  For finitely many $r\in(0,\infty)$, the
  $\alpha^H$-boundary of $u$ intersects the $\alpha$-boundary of
  $\tau_r(u)$. Generically, these intersections will be transverse, so
  it suffices to consider that case. We also assume that the complex
  structure on $\Sigma\times\Delta$ is chosen so that $\alphas$ and
  $\alphas^H$ meet at right angles.
  Let  $\bdy^{\alpha}u = u(\bdy S)\cap \alphas$ and $\bdy^{\alpha^H}u
  = u(\bdy S)\cap \alphas^H$ denote the subset of $u(\bdy S)$.
  Consider a point $q\in
  (\bdy^{\alpha^H} u)\cap(\bdy^\alpha \tau_r(u))$, and let
  $q_1=u^{-1}(q)$, $q_2=(\tau_r(u))^{-1}(q)$. We can choose an
  identification of a neighborhood of $q$ with a neighborhood of
  $(0,1)$ in $\CC\times \{x+iy\mid x\leq 1\}$ and identifications of
  neighborhoods of $q_1$ and $q_2$ with neighborhoods of $0$ in
  $\HalfPlane$ so that, to first order, on these neighborhoods the
  maps $u$ and $\tau_r(u)$ are given by either
  \begin{align*}
    u(x+iy)&=(x+iy,1-y+ix)\\
    \tau_r(u)(x'+iy')&=(y'-ix',1-\lambda y'+i\lambda x')\\
    \shortintertext{or}
    u(x+iy)&=(x+iy,1-y+ix)\\
    \tau_r(u)(x'+iy')&=(-y'+ix',1-\lambda y'+i\lambda x'),
  \end{align*}
  depending on whether the sign of the intersection corresponding
  to~$q$ in $(\bdy^{\alpha^H} u)\cap(\bdy^\alpha
  \tau_r(u))$. (Again, we are using the fact that, by
  Condition~(\ref{item:tri-J-global-translation-inv}), $\tau_r(u)$ is
  holomorphic.)
  Here, $\lambda$ is a positive real number.
  Consider the first case. The points in the domains of $u$ and
  $\tau_{r+\delta}(u)$ that map to the same point are, respectively,
  \begin{align*}
    x + iy&=\frac{\lambda}{1+\lambda^2}(\delta + i\lambda\delta) \\
    x' + iy'&=\frac{\lambda}{1+\lambda^2}(-\lambda\delta + i\delta).
  \end{align*}
  These points lie in the domains of $u$ and $\tau_{r+\delta}(u)$ if
  $\delta>0$, and not if $\delta<0$. So, for this local model, one
  intersection point appears as we translate $u$ up. Similarly, in the
  second case, the points in the domain mapping to the intersection
  are
  \begin{align*}
    x + iy&=\frac{\lambda}{1+\lambda^2}(\delta - i\lambda\delta) \\
    x' + iy'&=\frac{\lambda}{1+\lambda^2}(-\lambda\delta - i\delta).
  \end{align*}
  These points lie in the domains of $u$ and $\tau_{r+\delta}(u)$ if
  $\delta<0$, and not if $\delta>0$. So, for this local model, one
  intersection point disappears as we translate $u$ up.  This explains
  the term $\bdy^{\alpha^H}(B)\cdot\bdy^{\alpha}(B)$, except for the
  part of $\bdy^{\alpha^H}(B)\cdot\bdy^{\alpha}(B)$ coming from
  intersections at the corners at $\Theta_o$. For contributions coming
  from a corner of $u$ intersecting an edge of $\tau_r(u)$ or
  vice-versa, the analysis is the same as before except that the
  intersection point that appears or disappears is on the boundary.
  The expression $\bdy^{\alpha^H}(B)\cdot\bdy^{\alpha}(B)$ also picks
  up $-1/4$ from each of the $g$ points in $\Theta_o$, coming from the
  intersections of the corners of the domain with themselves. These
  terms do not correspond to points appearing or disappearing; this is
  the reason for the $g/4$ term in
  Equation~\eqref{eq:tri-index-intermediate-step}.
  
  As in the bigon case, $u\cdot \tau_R(u)=n_\x(B)+n_\y(B)-g/2$. So, 
  \begin{align*}
    \chi(S)&=e(B)-\br(\pi_\Sigma\circ u)+\frac{3g}{4}
      +\sum_i\frac{\abs{P_i}}{2} +\sum_i \frac{\abs{P_i^H}}{2}\\
    &=e(B)+n_\x(B)+n_\y(B)+\bdy^{\alpha^H}(B)\cdot\bdy^{\alpha}(B)-g/2+\iota([\vec{P}],[\vec{P}^H]).
  \end{align*}
  (The last line uses Formula~\eqref{eq:TriangleIndex}). This proves
  the result.
\end{proof}

For notational convenience, let
\begin{align*}
  \ind(B,\vec{\rhos},\vec{\rhos}^H)&=e(B)+n_\x(B)+n_\y(B)+\bdy^{\alpha^H}(B)\cdot\bdy^\alpha(B)-g/2+\iota(\vec
  \rhos,\vec \rhos^H)\\
  \chi_\emb(B,\vec{\rhos},\vec{\rhos}^H)&=e(B)-n_\x(B)-n_\y(B)+g-\bdy^{\alpha^H}(B)\cdot\bdy^{\alpha}(B)-\iota(\vec{\rhos},\vec{\rhos}^H).
\end{align*}

\begin{lemma}\label{lem:tri-index-adds}
  The function $\ind(B,\vec{\rhos},\vec{\rhos}^H)$ is additive with respect to
  gluing domains at $v_{23}$ or $v_{31}$, in the following sense. If
  $B\in\pi_2(\x,\y,\Theta_o)$, $B_{\alpha,\beta}\in\pi_2(\x',\x)$ and
  $B_{\alpha^H,\beta}\in\pi_2(\y,\y')$, and
  $(\vec{\rhos},\vec{\rhos}^H)$, $\vec\rhos_{\alpha,\beta}$ and
  $\vec{\rhos}_{\alpha^H,\beta}$ are sequences of non-empty sets of Reeb chords
  compatible with $B$, $B_{\alpha,\beta}$ and
  $B_{\alpha^H,\beta}$, respectively, then
  \begin{multline*}
  \ind(B_{\alpha,\beta}*_{31}B*_{23}B_{\alpha^H,\beta},(\vec{\rhos}_{\alpha,\beta},\vec{\rhos}),(\vec{\rhos}^H,\vec{\rhos}_{\alpha^H,\beta}))\\=
  \ind(B_{\alpha,\beta},\vec{\rhos}_{\alpha,\beta})+
  \ind(B,\vec{\rhos},\vec{\rhos}^H)+
  \ind(B_{\alpha^H,\beta},\vec{\rhos}_{\alpha^H,\beta}).
  \end{multline*}
\end{lemma}
\begin{proof}
  The proof is essentially the same as in the closed case~\cite[Theorem
  3.3]{Sarkar11:IndexTriangles} (which also inspired the proof of
  Proposition~\ref{prop:indAdditive}).
\end{proof}

Next we prove an analogue of Lemma~\ref{lem:tri-index-adds} for gluing
at $v_{12}$:
\begin{lemma}\label{ind:tri-ind-move-past-v12}
  Let $B\in\pi_2(\x,\y,\Theta_o)$ and let
  $(\vec{\rhos},\vec{\rhos}^H)$ be a pair of sequences of non-empty sets of Reeb chords
  compatible with $B$. Write $\vec{\rhos}=(\rhos_1,\dots,\rhos_j)$ and
  $\vec{\rhos}^H=(\rhos^H_{j+1},\dots,\rhos^H_l)$. Let $\rhos^H_j$ be
  the set of $\alpha^H$-Reeb chords corresponding to $\rhos_j$. Let
  $B_{\alpha,\alpha^H}\in\pi_2(\Theta_{o'},\Theta_{o})$ (for some
  $o'$) be a domain so that $B*_{12} B_{\alpha,\alpha^H}$ is compatible with
  $((\rhos_1,\dots,\rhos_{j-1}),\penalty300(\rhos_j^H,\rhos_{j+1}^H,\dots,\rhos_l^H))$. Then
  \begin{equation}\label{eq:tri-ind-glue-v12}
  \ind(B*_{12}B_{\alpha,\alpha^H},(\rhos_1,\dots,\rhos_{j-1}),(\rhos_j^H,\rhos_{j+1}^H,\dots,\rhos_l^H))
  =\ind(B,\vec{\rhos},\vec{\rhos}^H).
  \end{equation}
\end{lemma}
\begin{proof}
  The difference between the two sides of
  Equation~\eqref{eq:tri-ind-glue-v12} is:
  \begin{equation}\label{eq:tri-difference}
    \begin{aligned}
      &e(B_{\alpha,\alpha^H})+n_\x(B_{\alpha,\alpha^H})+n_\y(B_{\alpha,\alpha^H})\\
      &+\bdy^{\alpha^H}(B)\cdot\bdy^\alpha(B_{\alpha,\alpha^H})
      +\bdy^{\alpha^H}(B_{\alpha,\alpha^H})\cdot\bdy^\alpha(B)
      +\bdy^{\alpha^H}(B_{\alpha,\alpha^H})\cdot\bdy^\alpha(B_{\alpha,\alpha^H})\\
      &+\sum_{i=1}^{j-1}\bigl(L(\rhos_i,\rhos_j^H)-L(\rhos_i,\rhos_j)\bigr)
      +\sum_{i=j+1}^{l}\bigl(L(\rhos_j^H,\rhos_i^H)-L(\rhos_j,\rhos_i^H)\bigr)\\
      &+\iota(\rhos_j^H)-\iota(\rhos_j).
    \end{aligned}
  \end{equation}

  By inspection, $\iota(\rhos_j^H)=\iota(\rhos_j)$. Since $\bdy^\bdy
  B_{\alpha,\alpha^H} = \rhos_j^H - \rhos_j$, by
  multi-linearity of $L$, 
  \begin{align*}
    \sum_{i=1}^{j-1}\bigl(L(\rhos_i,\rhos_j^H)-L(\rhos_i,\rhos_j)\bigr)&=L(\rhos_1+\dots+\rhos_{j-1},\bdy^\bdy
    B_{\alpha,\alpha^H})\\
    \sum_{i=j+1}^{l}\bigl(L(\rhos_j^H,\rhos_i^H)-L(\rhos_j,\rhos_i^H)\bigr)&=L(\bdy^\bdy
    B_{\alpha,\alpha^H},\rhos_{j+1}^H+\dots+\rhos_{l}^H).
  \end{align*}
  We next use the techniques of part~(\ref{item:nxnw}) of
  Lemma~\ref{lemma:jittered} and Lemma~\ref{lem:jittered2}.  We have
  $\Theta_o - \x = \partial(\bdy^\alpha B + \rhos_1 + \cdots +
  \rhos_j)$, so $n_{\Theta_o}(B_{\alpha,\alpha^H}) -
  n_\x(B_{\alpha,\alpha^H}) = \bdy B_{\alpha,\alpha^H} \cdot
  (\bdy^\alpha B + \rhos_1 + \cdots + \rhos_j)$, which yields
  \begin{align*}
    \bdy^{\alpha^H}(B_{\alpha,\alpha^H})\cdot\bdy^\alpha(B)
    +L(\rhos_1+\dots+\rhos_{j},\bdy^\bdy B_{\alpha,\alpha^H})
    +n_\x(B_{\alpha,\alpha^H})&=n_{\Theta_o}(B_{\alpha,\alpha^H}).\\
    \shortintertext{Similarly,}
    \bdy^{\alpha^H}(B)\cdot\bdy^\alpha(B_{\alpha,\alpha^H})
    +L(\bdy^\bdy B_{\alpha,\alpha^H},\rhos_{j+1}^H +\dots+\rhos_{l}^H)
    +n_\y(B_{\alpha,\alpha^H})
    &=n_{\Theta_o}(B_{\alpha,\alpha^H}).
  \end{align*}
  So, Formula~\eqref{eq:tri-difference} reduces to
  \[
  e(B_{\alpha,\alpha^H})+2n_{\Theta_o}(B_{\alpha,\alpha^H})+\bdy^{\alpha^H}(B_{\alpha,\alpha^H})\cdot\bdy^\alpha(B_{\alpha,\alpha^H})+L(\rhos_j,\bdy^\bdy
  B_{\alpha,\alpha^H}).
  \]
  Since $[\rhos_j^H]-[\rhos_j]=\bdy^\bdy B_{\alpha,\alpha^H}$, we have
  \[
  L(\rhos_j,\bdy^\bdy B_{\alpha,\alpha^H})=L(\rhos_j,\rhos_j^H).
  \]
  Similarly to part~(\ref{item:nxny}) of Lemma~\ref{lemma:jittered},
  \[
  n_{\Theta_o}(B_{\alpha,\alpha^H})+\bdy^{\alpha^H}(B_{\alpha,\alpha^H})\cdot\bdy^\alpha(B_{\alpha,\alpha^H})+L(\rhos_j,\rhos_j^H)=n_{\Theta_{o'}}(B_{\alpha,\alpha^H}).
  \]
  By inspection of the possible domains for $B_{\alpha,\alpha^H}$,
  \[
  n_{\Theta_o}(B_{\alpha,\alpha^H})+n_{\Theta_{o'}}(B_{\alpha,\alpha^H})+e(B_{\alpha,\alpha^H})=0.
  \]
  This proves the result.
\end{proof}

\begin{corollary}\label{cor:tri-emb-ind}
  If $\cM^B(\x,\y,\Theta_o;\TSource;\vec P, \vec P^H)$ contains an
  embedded curve $u$ then the expected dimension of $\cM^B$ near $u$
  is given by
  \[
  \ind(u)=\ind(b(B),\vec{\rhos}),
  \]
  where $\vec\rhos$ is the sequence of Reeb chords obtained
  from $([\vec P],[\vec P^H])$ by  by replacing
  each $\alpha^H$-chord with the
  corresponding $\alpha$-chord.
\end{corollary}
\begin{proof}
  This follows from Proposition~\ref{prop:tri-emb-ind-general}
  and Lemmas~\ref{lem:triple-pi2-decomp}, \ref{lem:tri-index-adds}
  and~\ref{ind:tri-ind-move-past-v12} and the observation that, in the
  notation of Lemma~\ref{lem:triple-pi2-decomp}, the
  small triangles $T_\x$ have index $0$.
\end{proof}

\subsubsection{Definition of the homomorphism associated to a handleslide}
For $\x\in\S(\Sigma,\penalty300\alphas,\betas)$,
$\y\in\S(\Sigma,\penalty300\alphas^H,\betas)$, $B\in\pi_2(\x',\y)$ and
$\vec\rhos$ a sequence of non-empty sets of Reeb chords so that
$(\x,\vec{\rhos})$ is strongly boundary monotone and $(B,\vec\rhos)$ is
compatible, define
  \begin{align*}
   \gls*{ModTriEmbRoneRtwo}
   &\coloneqq\!\!\!\!
      \bigcup_{\substack{
          [\vec P] = \vec\rhos_1, [\vec P^H] = \vec\rhos_2 \\
          b(B') = B\\
          \chi(\TSource)=\chi_{\emb}(B',\vec\rhos_1,\rhos_2)
        }}\!\!\!\!
      \Mod^{B'}(\x, \y,\Theta_{o};
      \TSource; \vec{P}_1,\vec{P}_2)\\
    \gls*{ModTriEmb}
    &\coloneqq
    \!\bigcup_{
        \substack{
          \vec{\rhos}=(\vec{\rhos}_1,\vec{\rhos}_2)
          }}\!
        \cM^B(\x,\y,\Theta_{o(\x,\vec{\rhos}_1)};\rhos_1,\rhos_2).
  \end{align*}
 The triangle map from $\CFDa(\Sigma,\alphas,\betas,z)$ to
$\CFDa(\Sigma,\alphas^H,\betas,z)$ is defined by
\begin{align}
  \label{eq:TriangleMap}
F_{\alphas,\alphas^H,\betas}(\x)&\coloneqq
\sum_{\y}\sum_{B\in\pi_2(\x',\y)}\,
\sum_{\{\vec{\rho}\,\mid\,\ind(B,\vec\rho) = 0\}}\!\!
\#\left(\Mod^B(\x,\y,\Theta;\vec{\rho})\right)a(-\vec\rho) \y\\
F_{\alphas,\alphas^H,\betas}(a\x)&\coloneqq a F_{\alphas,\alphas^H,\betas}(\x). \nonumber
\end{align}
\glsadd{FaaHb}%
As in the definition of the boundary operator, if both
$(\Sigma,\alphas,\betas,z)$ and
$(\Sigma,\alphas^H,\betas)$ are provincially admissible this sum is
finite, as the $0$-dimensional moduli spaces are compact and only
finitely many of them are non-empty.

The map $F_{\alphas,\alphas^H,\betas}$ is defined on all of
$\CFDa(\Sigma,\alphas,\betas,z)$, but
Lemma~\ref{lem:SpinCStructures} and the fact that
$\spinc_z(\x)=\spinc_z(\x')$ ensures that  restriction of $F_{\alpha,\alphas^H,\beta}$ to the summand
$\CFDa(\Sigma,\alphas,\betas,z;\spinc)$ of
$\CFDa(\Sigma,\alphas,\betas,z)$ representing a fixed $\SpinC$
structure $\s$  is
contained in the corresponding summand
$\CFDa(\Sigma,\alphas^H,\betas,z;\spinc)$ of
$\CFDa(\Sigma,\alphas^H,\betas,z)$.

Our next goal is to prove that $F_{\alphas,\alphas^H,\betas}$ is a
chain map. The new technical point relates to degenerations at
$v_{12}$, to which we turn next. 

\subsubsection{Holomorphic curves at \textalt{$v_{12}$}{v12}}
We start by completely analyzing what curves can degenerate
at $v_{12}$; but first some more notation.

There are $2(g+k)$ components of
$\Sigma\setminus(\alphas\cup\alphas^H)$ not adjacent to $z$.
\begin{itemize}
\item For each $\alpha$-circle $\alpha_i^c$ ($i=1,\dots,g-k$) there is
  a pair
  of bigons $R_{i,1}^c$ and $R_{i,2}^c$ contained entirely in the interior of
  $\overline{\Sigma}$, with boundaries in
  $\alpha_i^c\cup\alpha_i^{c,H}$. (The numbering $1$ and $2$ is
  arbitrary.)
\item For each $\alpha$-arc $\alpha_i^a$ with $i\geq 2$ there is a pair of
  bigons $R_{i,1}^a$ and $R_{i,2}^a$ with boundary in
  $\alpha_i^a\cup\alpha_i^{a,H}$.  (Again, the numbering $1$ and $2$
  is arbitrary.) These
  bigons are adjacent to $\bdy\overline{\Sigma}$; one of the two punctures
  of each $R_{i,j}^a$ corresponds to $e\infty$ in~$\Sigma$.
\item There is an annulus $R_{1,1}^a$ with boundary in
  $\alpha_1^a\cup\alpha_1^{a,H}\cup\alpha_1^c\cup\alpha_1^{c,H}$ and a
  bigon $R_{1,2}^a$ with boundary in
  $\alpha_1^a\cup\alpha_1^{a,H}$. Both $R_{i,1}^a$ and $R_{i,2}^a$ are
  adjacent to $\bdy\overline{\Sigma}$.
\end{itemize}
We call the regions $R_{i,1}^a$ and $R_{i,2}^a$ (respectively
$R_{i,1}^c$ and $R_{i,2}^c$) \emph{twin regions}.
\index{twin regions}\index{region!twin}%
See Figure~\ref{fig:alphasalphasH}.

\begin{lemma}\label{lem:holo-curves-v12}
  Let 
  \[
  u\co (S,\bdy S)\to (\Sigma\times[0,1]\times\RR,
  (\alphas\times\{0\}\times\RR)\cup(\alphas^H\times\{1\}\times\RR))
  \]
  be a holomorphic curve asymptotic to a generator $\Theta_o$ at
  $-\infty$, with multiplicity $0$ at~$z$. Then $S$ is a disjoint
  union of bigons and possibly an annulus. Moreover:
  \begin{enumerate}[label=(c-\arabic*),ref=c-\arabic*]
  \item \index{(c-1)--(c-4)}
    \label{item:Disjointness}
    If $S_i$ and $S_j$ are distinct components of $S$ then
    $\pi_\Sigma(u(S_i))\cap \pi_\Sigma(u(S_j))=\emptyset$.
  \item 
    \label{item:Bigons}
    If $S_i$ is a bigon then either $u|_{S_i}$ is a trivial strip
    (i.e., $(\pi_\Sigma\circ u)|_{S_i}$ is constant and $(\pi_\DD\circ
    u)|_{S_i}$ is a diffeomorphism) or $(\pi_\Sigma\circ u)|_{S_i}$ is
    a diffeomorphism onto its image, which is one of 
    the distinguished bigons $R_{i,j}^c$ or $R_{i,j}^a$.
  \item 
    \label{item:Annulus}
    The image of the annular component of $S$ (if it exists) is one of $R_{1,1}^a$,
    $R_{1,1}^a+R_{1,1}^c$, $R_{1,1}^a+R_{1,2}^c$,
    $R_{1,1}^a+R_{1,1}^c+R_{1,2}^c$, $R_{1,1}^a+2R_{1,1}^c$ or
    $R_{1,1}^a+2R_{1,2}^c$.
  \item 
    \label{item:AtMostOneCorner}
    If a region $R_{i,1}^a$ (respectively $R_{i,2}^a$, $R_{i,1}^c$,
    $R_{i,2}^c$) is in the image of $\pi_\Sigma\circ u$ then
    its twin region $R_{i,2}^a$ (respectively $R_{i,1}^a$,
    $R_{i,2}^c$, $R_{i,1}^c$) is not in the image of $\pi_\Sigma\circ u$,
    with the possible exception that both $R_{i,1}^c$ and $R_{i,2}^c$
    may be in the image of $\pi_\Sigma\circ u$ on an annular component. 
  \end{enumerate}
\end{lemma}
\begin{proof}
  Suppose first that $\pi_\Sigma\circ u$ does not cover the annulus
  $R_{1,1}^a$. Since $u$ is asymptotic to the generator $\Theta_o$,
  each bigon $R_{i,j}^a$ or $R_{i,j}^c$ appears with multiplicity $0$
  or~$1$, and at most one of each pair of twin bigons appears in the
  image of $\pi_\Sigma\circ u$. The result follows.

  Next, suppose that $\pi_\Sigma\circ u$ covers the annulus
  $R_{1,1}^a$. Then, as before, the domain of $u$ is a disjoint union
  of bigons away from a neighborhood of the annulus. The asymptotics
  at $-\infty$ imply that near the annulus, the domain has one of the
  forms $R_{1,1}^a$, $R_{1,1}^a+R_{1,1}^c$, $R_{1,1}^a+R_{1,2}^c$,
  $R_{1,1}^a+R_{1,1}^c+R_{1,2}^c$, $R_{1,1}^a+2R_{1,1}^c$ or
  $R_{1,1}^a+2R_{1,2}^c$. The result follows.
\end{proof}

\subsubsection{Gluing results at \textalt{$v_{12}$}{v12}}
Next we prove two gluing results for these curves at~$v_{12}$.  (See
Figure~\ref{fig:bigon-degen} for a typical example of this gluing.)
Again, we
need a little more terminology. The holomorphic curves in
Lemma~\ref{lem:holo-curves-v12} are asymptotic to a generator
$\Theta_o$ at $-\infty$. At $+\infty$, they are asymptotic to
$g$-tuples of the form
\[
\w^+_{\pd}=(\{w_1,\dots,w_\ell\}\times[0,1]\times\{+\infty\})\cup(\{\rho_{1},\dots,\rho_{g-\ell}\}\times[0,1]\times\{+\infty\})
\]
where each $w_i$ is a point in $\alphas^H\cap\alphas$ and each $\rho_i$ is a Reeb
chord in $\bdy\bSigma$ from some $\alpha_j^{a,H}$ to
some~$\alpha_k^a$. We call tuples $\w_\pd^+$ of this form
\emph{positive partially-diagonal generators}.
\index{generator!partially-diagonal}\index{generator!partially-diagonal!positive}%
\index{partially-diagonal generator|see{generator, partially-diagonal}}%
(In fact, for the curves in Lemma~\ref{lem:holo-curves-v12}, $j=k$,
but we do not require this condition for partially diagonal
generators.)
We will abuse notation and also refer to the set
$\w^+_\pd=\{w_1,\dots,w_\ell,\rho_1,\dots,\rho_{g-\ell}\}$ as a
positive partially-diagonal generator. 
\emph{Negative partially-diagonal generators} are defined
similarly, but with Reeb chords running from some $\alpha_j^a$ to some
$\alpha_k^{a,H}$, and are denoted by $\w_\pd^-$.
\index{generator!partially-diagonal!negative}
\glsadd{pdgen}

We collect the holomorphic curves at $v_{12}$ into moduli spaces. Specifically,
for each generator $\Theta_o$, positive
partially-diagonal generator $\w_\pd^+$ and domain
$B\in\pi_2^{\alphas^H,\alphas}(\Theta_o,\w_\pd^+)$, there is a
corresponding moduli space of (stable) holomorphic curves
$\cM^B(\Theta_o,\w_\pd^+)$. (As usual, we mod out by the action of $\RR$ by
translation.) Lemma~\ref{lem:holo-curves-v12} shows that the
topological type of the source of curves in $\cM^B(\Theta_o,\w_\pd^+)$
is determined by $B$ and $\w_\pd^+$, but when we want to name the
source explicitly we will also write this moduli space as
$\cM^B(\Theta_o,\w_\pd^+;\Source)$.
\glsadd{ModVonetwo}\glsadd{ModVonetwoSource}

It follows from Lemma~\ref{lem:holo-curves-v12} that the moduli spaces
$\cM^B(\Theta_o,\w_\pd^+)$ are typically high-dimensional:
\begin{corollary}\label{cor:small-ind-v12-curves} Fix a generic split
  almost complex structure on $\Sigma\times[0,1]\times\RR$.
  Let $\w_\pd^+=\{w_1,\dots,w_\ell,\rho_1,\dots,\rho_{g-\ell}\}$ be a
  positive partially-diagonal generator.  Then the moduli space
  $\cM^B(\Theta_o,\w_\pd^+)$ for generic $J$ is either empty or has dimension at least
  $\max\{g-\ell-1,0\}$. Further:
  \begin{enumerate}
  \item
    \label{case:Further1}
    If $g-\ell=0$ and $\cM^B(\Theta_o,\w_\pd^+)$ is non-empty and
    has dimension $0$ then $B$ is the bigon $R_{i,j}^c$ for some $c$
    (and $\w_\pd^+=\Theta_o$).
  \item
    \label{item:small-dim-w-chords} If $\cM^B(\Theta_o,\w_\pd^+)$ is non-empty and has dimension
    $g-\ell-1$ then every coefficient in $B$ is $0$ or $1$ and each
    connected component of the support of $B$ is one of $R_{i,j}^a$
    for some $i>1$;
    $R_{1,2}^a$; $R_{1,1}^a\cup R_{1,1}^c$; or
    $R_{1,1}^a\cup R_{1,2}^c$.
  \end{enumerate}
  Conversely, each bigon has a unique holomorphic representative, and
  the sum of the number of holomorphic representatives of
  $R_{1,1}^a+R_{1,1}^c$ and $R_{1,1}^a+R_{1,2}^c$ is odd.
\end{corollary}
\begin{proof}
  This is a straightforward analysis of the moduli spaces
  corresponding to the domains in Lemma~\ref{lem:holo-curves-v12}, as
  follows.
  
  Property~(\ref{item:Disjointness})
  ensures that the moduli space $\tcM^B(\Theta_o,\w_\pd^+)$
  (i.e., before we divide out by the $\RR$ action) splits as a product of
  moduli spaces associated to the component curves. By the
  classification from Lemma~\ref{lem:holo-curves-v12}, these factors
  correspond to moduli spaces representing the domains in the
  following table, listed with their number of chords and expected
  dimension (before dividing out by $\RR$):
  \begin{center}
  \begin{tabular}{ccc}
    \toprule
    & & Expected dimension \\
    Homology class & Number of chords & of moduli space \\
        \midrule 
    $0$ & $0$ & $0$ \\
    $R^c_{i,j}$ & $0$ & $1$ \\
    $R^a_{i,j}$, $i>1$ & $1$ & $1$ \\
    $R^a_{1,2}$ & $1$ & $1$ \\
    $R^a_{1,1}$ & $1$ & $0$ \\
    $R^a_{1,1}+R^c_{1,1}$ & $1$ & $1$ \\
    $R^a_{1,1}+R^c_{1,2}$ & $1$ & $1$ \\
    $R^a_{1,1}+2R^c_{1,1}$ & $1$ & $2$ \\
    $R^a_{1,1}+R^c_{1,1}+R^c_{1,2}$ & $1$ & $2$ \\
    \bottomrule
    \end{tabular}
    \end{center}
    The only homology class where the number of chords
    exceeds the dimension is $R^a_{1,1}$; but in this
    case, the dimension is $0$, so (for generic~$J$) the corresponding
    moduli space is empty,
    and hence cannot appear as a factor of a non-empty moduli space.
    The desired inequality now follows from 
    the fact that
    $\dim \cM^B(\Theta_o,\w^+_{\pd})$ is the sum of the dimensions
    of the constituent moduli spaces minus one (for the $\RR$ action), 
    while $g-\ell$ is the total
    number of chords.

    The condition that $g-\ell=0$ is equivalent to the condition that
    each constituent piece has no chords in it, in which case the homology
    class must be a sum of domains $R^c_{i,j}$. Moreover, each
    component contributes $1$ to the dimension (before dividing out by $\RR$).
    This establishes Point~(\ref{case:Further1}).
    
    The condition that $g-\ell-1$ coincides with the dimension ensures
    that each constituent piece must have exactly as many chords as
    its dimension: thus, each constituent piece must be one of
    $R^a_{i,j}$ for some $i>1$, $R^a_{1,2}$, $R^a_{1,1}+R^c_{1,1}$ or
    $R^{a}_{1,1}+R^c_{1,2}$. This establishes
    Point~(\ref{item:small-dim-w-chords}).
\end{proof}

Given a negative partially-diagonal generator $\w_\pd^-$ and generators
$\x\in\Gen(\alphas,\betas)$ and $\y\in\Gen(\alphas^H,\betas)$ we can
consider the moduli space
\[
\cM^B(\x,\y,\w_\pd^-;\TSource;\vec{P},\vec{P}^H)
\]
\glsadd{ModTriSourcePPpd}%
of holomorphic triangles asymptotic to the generators $\x$, $\y$ and
$\w_\pd^-$ at $v_{31}$, $v_{23}$ and $v_{12}$ respectively and with
ordered partitions of Reeb chords $\vec{P}$ and $\vec{P}^H$ along
$e_1$ and $e_2$, respectively. We can also consider the moduli space
\[
\cM^B(\x,\y,\w_\pd^-;\vec{\rhos},\vec{\rhos}^H)
\]
\glsadd{ModTriEmbpd}%
of embedded holomorphic triangles asymptotic to the sequences of non-empty sets
of Reeb chords $\vec{\rhos}$ and $\vec{\rhos}^H$ along $e_1$ and
$e_2$, respectively.

A pair of partially-diagonal generators
$\w^+_\pd=\{w^+_1,\dots,w^+_\ell,\rho_1,\dots\rho_{g-\ell}\}$ and
$\w^-_\pd=\{w^-_1,\dots,w^-_m,\rho'_1,\dots,\rho'_{g-m}\}$ of opposite polarity
are called
\emph{compatible} if $\{w^+_1,\dots,w^+_\ell\}=\{w^-_1,\dots,w^-_m\}$ (so,
in particular, $m=\ell$) and there is a bijection $\mu\co
\{\rho_1,\dots,\rho_{g-\ell}\}\to \{\rho'_1,\dots,\rho'_{g-m}\}$ so
that $\rho_i$ and $\mu(\rho_i)$ are
adjacent in the sense that either $\rho_i^+=\mu(\rho_i)^-$ or
$\rho_i^-=\mu(\rho_i)^+$.
\index{generators!partially-diagonal!compatible}%
We will call the map $\mu$ a \emph{matching} between $\w^+_\pd$ and $\w^-_\pd$.
\index{generators!partially-diagonal!matching of}
\glsadd{mumatching}

To a a compatible pair of partially-diagonal generators $\w^+_\pd$ and
$\w^-_\pd$ and a matching $\mu$ between them, we associate two
sets of Reeb chords $\rhos_{\alphas}(\w^+_\pd,\w^-_\pd,\mu)$ and
$\rhos_{\alphas^H}(\w^+_\pd,\w^-_\pd,\mu)$, given by
\begin{align*}
  \rhos_{\alphas^H}(\w^+_\pd,\w^-_\pd,\mu)&\coloneqq\bigcup_{\rho_i^+=\mu(\rho_i)^-}\rho_i\uplus\mu(\rho_i)
  \\
  \rhos_{\alphas}(\w^+_\pd,\w^-_\pd,\mu)&\coloneqq\bigcup_{\rho_i^-=\mu(\rho_i)^+}\mu(\rho_i)\uplus\rho_i.
\end{align*}
When the matching is implicit or determined by the rest of the data,
as will usually be the case below, we will drop it from the
notation, writing simply $\rhos_{\alphas^H}(\w^+_\pd,\w^-_\pd)$
and $\rhos_{\alphas}(\w^+_\pd,\w^-_\pd)$.
\glsadd{rhosHofpdgens}\glsadd{rhosofpdgens}

\begin{figure}
  \centering
  \includegraphics[scale=.83333]{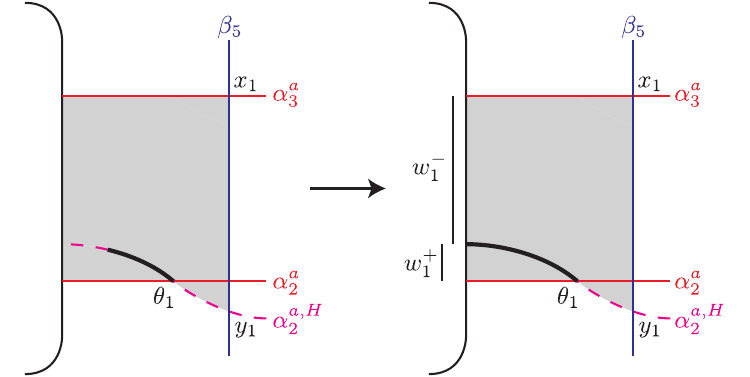}
  \caption[A degeneration at $v_{12}$]{\textbf{Example of a
      degeneration at $v_{12}$.} At the left is a $1$-parameter family of
    holomorphic curves connecting $x_1$, $y_1$ and $\theta_1$; the
    parameter is the length of the cut. One end of this moduli space
    corresponds to a degeneration at $v_{12}$ into a curve connecting
    $\theta_1$ and $w_1^+$ and a triangle connecting $x_1$, $y_1$ and $w_1^-$.}
  \label{fig:degen-at-v12}
\end{figure}

If $\w_\pd^+$ and $\w_\pd^-$ are compatible partially-diagonal
generators then the product of moduli spaces $\cM^B(\Theta_o,\w_\pd^+)\times
\cM^B(\x,\y,\w_\pd^-;\vec{\rhos},\vec{\rhos}^H)$ sits inside the
compactified moduli space
$\oocM^B(\x,\y,\Theta_o\semico\vec{\rhos}',\vec{\rhos}^{H,\prime}),$
where $\vec{\rhos}'$ is obtained from $\vec{\rhos}$ by adjoining the
set $\rhos_\alphas(\w_\pd^+,\w_\pd^-,\mu)$ at the end and
$\vec{\rhos}^{H,\prime}$ is obtained from $\vec{\rhos}$ by adjoining
the set $\rhos_{\alphas^H}(\w_\pd^+,\w_\pd^-,\mu)$ at the beginning.
(In principle, this construction depends on a choice of matching~$\mu$, but
boundary monotonicity of
$\cM^B(\x,\y,\w_\pd^-;\vec{\rhos},\vec{\rhos}^H)$ implies that this
matching is determined by the rest of the data.)  In other words,
there is an obvious notion of a sequence of curves in
$\cM^B(\x,\y,\Theta_o;\vec{\rhos}',\vec{\rhos}^{H,\prime})$
converging to a pair $\cM^B(\Theta_o,\w_\pd^+)\times
\cM^B(\x,\y,\w_\pd^-;\vec{\rhos},\vec{\rhos}^H)$. See
Figure~\ref{fig:degen-at-v12} for an example.

In order to understand the structure of
$\oocM^B(\x,\y,\Theta_o\semico\vec{\rhos}',\vec{\rhos}^{H,\prime})$
near the boundary stratum $\cM^B(\Theta_o,\w_\pd^+)\times
\cM^B(\x,\y,\w_\pd^-;\vec{\rhos},\vec{\rhos}^H)$, we will study some
auxiliary moduli spaces. Consider
$\cM^B(\Theta_o,\w_\pd^+\semico\Source)\times
\cM^B(\x,\y,\w_\pd^-\semico\TSource\semico\vec{P},\vec{P}^H)$. There is an obvious
way to preglue $\Source$ and $\TSource$ at the punctures corresponding
to $w_\pd^\pm$, giving a source $\Source\glue\TSource$. The
preglued surface $\Source\glue\TSource$ has punctures
$r_1,\dots,r_{a}$ labeled by the Reeb chords in
$\rhos_\alphas(\w_\pd^+,\w_\pd^-)$ and 
$r_1^H,\dots,r_{b}^H$ labeled by the Reeb chords in
$\rhos_{\alphas^H}(\w_\pd^+,\w_\pd^-)$, where $a+b=g-\ell$ is the
number of chords in $w_\pd^\pm$. For any permutations $\sigma_a\in S_a$
and $\sigma_b\in S_b$ we can consider the moduli space
\[
\cM^B(\x,\y,\Theta_o;\Source\glue\TSource;(\vec{P},(r_{\sigma_a(1)},\dots,r_{\sigma_a(a)})),((r^H_{\sigma_b(1)},\dots,r^H_{\sigma_b(b)})\vec{P}^H)).
\]
(Recall that for sequences $s$ and $t$, $(s,t)$ denotes the
concatenation of $s$ and $t$.)
There is an obvious notion for a sequence of curves in this moduli
space to converge to a curve in $\cM^B(\Theta_o,\w_\pd^+;\Source)\times
\cM^B(\x,\y,\w_\pd^-;\TSource;\vec{P},\vec{P}^H)$.

With this terminology in place, we have the following gluing results:
\begin{proposition}\label{prop:glue-v12} There is a neighborhood of
  the subspace
  $\cM^B(\Theta_o,\w_\pd^+;\Source)\times
  \cM^B(\x,\y,\w_\pd^-;\TSource;\vec{P},\vec{P}^H)$ in
  \[
  \bigcup_{\sigma_a\in S_a, \sigma_b\in S_b}
  \cM^B(\x,\y,\Theta_o;\Source\glue\TSource;
  (\vec{P},(r_{\sigma_a(1)},\dots,r_{\sigma_a(a)})),
  ((r^H_{\sigma_b(1)},\dots,r^H_{\sigma_b(b)})\vec{P}^H))
  \]
  that is homeomorphic to $\cM^B(\Theta_o,\w_\pd^+;\Source)\times
  \cM^B(\x,\y,\w_\pd^-;\TSource;\vec{P},\vec{P}^H)\times[0,\epsilon)$.
\end{proposition}
\begin{proof}
  At first glance, it seems that this requires a new kind of gluing
  argument. Specifically, suppose that
  $u_{\alpha^H,\alpha}\in\cM^B(\Theta_o,\w_\pd^+;\Source)$ and
  $u_{\alpha,\alpha^H,\beta}\in\cM^B(\x,\y,\w_\pd^-;\TSource;\vec{P},\vec{P}^H)$.
  Write
  $\w_\pd^+=\{w_1^+,\dots,w_\ell^+,\rho_1,\dots,\rho_{g-\ell}\}$ and
  $\w_\pd^-=\{w_1^-,\dots,w_\ell^-,\rho'_1,\dots,\rho'_{g-\ell}\}$. 
  At a puncture $q_i$ of $u_{\alpha^H,\alpha}$ corresponding to $\rho_i$,
  we have that
  $\lim_{q\to q_i}t(u_{\alpha^H,\alpha}(q))=\infty$ and $\lim_{q\to
    q_i}\pi_\Sigma(u_{\alpha^H,\alpha}(q))=p$, the puncture
  in~$\Sigma$ (so the curve goes to infinity in two
  directions at once).
  Similarly, at a puncture $q'_i$ of
  $u_{\alpha,\alpha^H,\beta}$ corresponding to $\rho'_i$, $\lim_{q\to
    q'_i}\pi_\Delta(u_{\alpha,\alpha^H,\beta}(q))=v_{12}$ and
  $\lim_{q\to q'_i}\pi_\Sigma(u_{\alpha,\alpha^H,\beta}(q))=p$. Such
  diagonal asymptotics seem not to have been treated in the
  literature.

  Because of the simple form of the curves in
  $\cM^B(\Theta_o,\w_\pd^+)$, we can avoid the new analysis by
  observing that our situation is biholomorphic to a familiar
  one, as follows.
  
  Consider first the case that $\ell=g-1$, so that each of
  $\w_\pd^\pm$ contains a single Reeb chord.  Recall that
  $\Sigma_{\widebar{e}}$ denotes the result of filling in the puncture
  on~$\Sigma$. The various $\widebar\alpha_i^a$ and
  $\widebar\alpha_i^{a,H}$ meet at the puncture~$p$ in
  $\Sigma_{\widebar e}$, in general with a corner. To fix notation,
  suppose that the chord $\rho_1$ runs from $\alpha_i^{a,H}$ to
  $\alpha_j^a$ and that $\rho'_1$ runs from $\alpha_j^a$ to
  $\alpha_k^{a,H}$; the case that $\rho_1'$ runs from $\alpha_k^a$
  to $\alpha_i^{a,H}$ is similar.
  Let $D_z$ be the component of
  $\Sigma\setminus(\alphas\cup\betas)$ adjacent to $z$ and let
  $\Sigma'_{\overline{e}}=\Sigma_{\overline e}\setminus D_z$. We can
  choose a complex structure $j'_\Sigma$ (and smooth structure) on
  $\Sigma'_{\overline e}$ so 
  that $j'_\Sigma$ agrees with the fixed complex structure $j_\Sigma$
  from Chapter~\ref{chap:structure-moduli} on $\Sigma'_{\overline
    e}\setminus\{e\}$ and so that the chord $\rho_1\uplus\rho'_1$
  corresponds to an $180^\circ$ angle with respect to
  $j'_\Sigma$. (See Figure~\ref{fig:straighten}. Note that while the
  total angle at the corner of $\Sigma'_{\overline e}$ may be more than $360^\circ$, the
  holomorphic curves are contained in a subset with total angle
  $180^\circ$.) 
  Let $J'$ be an
  almost complex structure on $\Sigma'_{\overline e}\times\Delta$
  which agrees with $J$ away from $p$ and with $j'_\Sigma\times
  j_\Delta$ near $p$. 

  The holomorphic curve $u_{\alpha,\alpha^H,\beta}$ has image in
  $\Sigma'_{\overline e}\times\Delta$ and is holomorphic with respect
  to $J'$.  Similarly, the holomorphic curve $u_{\alpha^H,\alpha}$ has
  image in $\Sigma'_{\overline e}\times[0,1]\times\RR$ and is
  holomorphic with respect to the complex structure on
  $\Sigma'_{\overline e}\times[0,1]\times\RR$ induced by
  $J'$. Moreover, the gluing problem for $u_{\alpha,\alpha^H,\beta}$
  and $u_{\alpha^H,\alpha}$ in $\Sigma'_{\overline
    e}\times[0,1]\times\RR$ has the same form as the usual gluing of
  holomorphic curves in Heegaard Floer theory at a generator. So, in
  this case the result follows from the usual gluing theorem,
  \cite[Proposition A.2]{Lipshitz06:CylindricalHF}, say.

  \begin{figure}
    \includegraphics[scale=.83333]{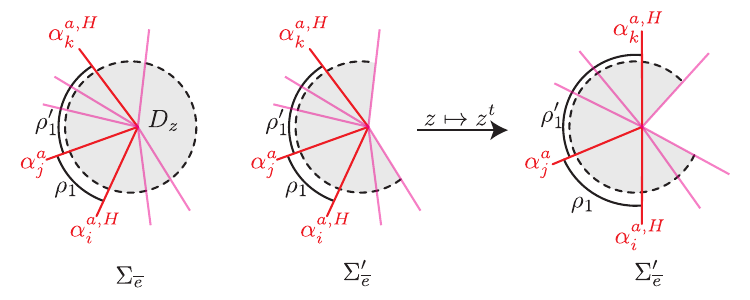}
    \caption[Straightening $\alpha$-arcs]{\textbf{Straightening-out
        $\alpha$-arcs.} The conformal map giving this straightening is
      locally the map $z\mapsto z^t$ for some appropriate
      $t>0$.}\label{fig:straighten}
  \end{figure}
  
  The general case, where $\w_\pd^\pm$ contains more than one Reeb
  chord, follows from the above argument with one more
  observation. The almost complex structure $J'$ that we used above depended on the
  endpoints of the chords $\rho_1$ and $\rho'_1$ being glued; it may not
  be possible to choose a single almost complex structure which
  straightens out all of the relevant arcs. However, we may choose the
  complex structure $J'$ on $\Sigma\times \Delta$ to depend on the point
  in the source of $u_{\alpha,\alpha^H,\beta}$. (In other words, for
  curves with source $S$, we are
  now considering holomorphic sections of a topologically trivial bundle
  $\Sigma\times\Delta\times S\to S$, with a complex structure on
  $\Sigma\times\Delta\times S$ given by $J'_s\times j_S$, where $j_S$
  is the complex structure on $S$ and $J_s$ is an almost complex
  structure on $\Sigma\times\Delta$; compare~\cite[Section
  3]{Seidel08:DehnTwistLectures}.)
  In particular, we can choose
  $J'=J'(s,d,p)$ ($s\in\Sigma$, $d\in\Delta$, $p\in S$ the source of
  $u_{\alpha,\alpha^H,\beta}$) so that for $p$ near the puncture $p'_i$,
  the angle corresponding to $\rho_i\uplus\rho'_i$ or
  $\rho'_i\uplus\rho_i$ (as appropriate) is $180^\circ$. The result then
  follows as before.
 \end{proof}

It follows from Proposition~\ref{prop:glue-v12} that if
$\cM^B(\x,\y,\Theta_o;\vec{\rho})$ is a $1$-dimensional moduli space
of holomorphic curves, where $\vec{\rho}$ is the discrete partition,
then, in degenerations at $v_{12}$, the resulting curve for
$(\alphas^H,\alphas)$ has index $1$, and hence at most one non-trivial
component. In particular, by Lemma~\ref{lem:holo-curves-v12}, at most
one $\alpha$-$\alpha^H$-Reeb chord is involved.

\begin{proposition}\label{prop:tri-master}
  Suppose that $B\in\pi_2(\x,\y)$ and sequences of non-empty sets of Reeb chords
  $\vec{\rhos}=(\rhos_1,\dots,\rhos_m)$ and
  $\vec{\rhos}^H=(\rhos^H_1,\dots,\rhos^H_n)$ so that
  $\ind(B,(\vec{\rhos},\vec{\rhos}^H))=1$. For a sufficiently generic,
  admissible almost complex structure $J$, the boundary of the moduli space
  $\ocM^B(\x,\y,\Theta_o;\vec{\rhos},\vec{\rhos}^H)$ (where
  $o=o(\x,\vec{\rhos})$) consists of curves of the following types:
  \begin{enumerate}[label=(MT-\arabic*),ref=MT-\arabic*]
  \item\label{item:tri-2story-31}\index{(MT-1)--(MT-9)} Two-story buildings corresponding to breaks at $v_{31}$, i.e., curves in 
    \[
    \cM^{B_1}(\x,\mathbf{u};(\rhos_1,\dots,\rhos_i))\times \cM^{B_2}(\mathbf{u},\y,\Theta_o;(\rhos_{i+1},\dots,\rhos_m),\vec{\rhos}^H)
    \]
    where $B_1\in\pi_2^{\alphas,\betas}(\x,\mathbf{u})$ and
    $B_2\in\pi_2(\mathbf{u},\y)$ are such that $B_1*B_2=B$, and $1\leq i\leq
    m$.
  \item\label{item:tri-2story-23}  Two-story buildings corresponding to breaks at $v_{23}$, i.e., curves in 
    \[
    \cM^{B_1}(\x,\mathbf{u},\Theta_o;\vec{\rhos},(\rhos_{1}^H,\dots,\rhos_i^H))
    \times \cM^{B_2}(\mathbf{u},\y;(\rhos_{i+1}^H,\dots,\rhos_n^H))
    \]
    where $B_1\in\pi_2(\x,\mathbf{u})$ and
    $B_2\in\pi_2^{\alphas^H,\betas}(\mathbf{u},\y)$ are such that
    $B_1*B_2=B$, and $1\leq i\leq n$.
  \item\label{item:tri-join} Join curve ends.
  \item\label{item:tri-split} Split curve ends.
  \item\label{item:tri-collision} Other collisions of levels among $\vec{\rhos}$ or $\vec{\rhos}^H$.
  \item\label{item:tri-shuffle} Shuffle curve ends.
  \item\label{item:tri-v12-zero} Degenerations at $v_{12}$ where the
    resulting partially diagonal generator involves only intersection
    points (not Reeb chords), i.e., curves in 
    \[
    \cM^{B_1}(\Theta_o,\w_\pd^+)\times \cM^{B_2}(\x,\y,\w_\pd^-;\vec{\rhos},\vec{\rhos}^{H})
    \]
    where
    $\w_\pd^+=\w_\pd^-=\{w_1,\dots,w_g\}\subset\alphas\cap\alphas^H$. Moreover,
    the homology class $B_1$ (and generator $\w_\pd^+$) is such that
    $\cM^{B_1}(\Theta_o,\w_\pd^+)$ is $0$-dimensional.
  \item\label{item:tri-v12-one} Degenerations where $\rhos_m$ approaches $v_{12}$, i.e.,
    curves in 
    \[
    \cM^{B_1}(\Theta_o,\w_\pd^+)\times \cM^{B_2}(\x,\y,\w_\pd^-;(\rhos_1,\dots,\rhos_{m-1}),\vec{\rhos}^{H})
    \]
    where $\w_\pd^\pm$ are compatible partially-diagonal generators,
    $\rhos_\alphas(\w^+_\pd,\w^-_\pd)=\rhos_m$ and
    $\rhos_{\alphas^H}(\w^+_\pd,\w^-_\pd)=\emptyset$.  Moreover, the
    homology class $B_1$ (and generator $\w_\pd^+$) is such that moduli space
    $\cM^{B_1}(\Theta_o,\w_\pd^+)$ is $(|{\rhos}_m|-1)$-dimensional
    (compare Corollary~\ref{cor:small-ind-v12-curves}).
  \item\label{item:tri-v12-two} Degenerations where $\rhos_1^H$ approaches $v_{12}$, i.e.,
    curves in 
    \[
    \cM^{B_1}(\Theta_o,\w_\pd^+)\times \cM^{B_2}(\x,\y,\w_\pd^-;\vec{\rhos},(\rhos^H_2,\dots,\rhos^H_{n}))
    \]
    where $\w_\pd^\pm$ are compatible partially-diagonal generators,
    $\rhos_\alphas(\w^+_\pd,\w^-_\pd)=\emptyset$ and
    $\rhos_{\alphas^H}(\w^+_\pd,\w^-_\pd)=\rhos^H_1$.  Moreover, the
    homology class $B_1$ (and generator $w_\pd^+$) is such that moduli space
    $\cM^{B_1}(\Theta_o,\w_\pd^+)$ is $(|{\rhos}^H_1|-1)$-dimensional
    (compare Corollary~\ref{cor:small-ind-v12-curves}).
  \end{enumerate}
  Conversely, each curve of type (\ref{item:tri-2story-31}),
  (\ref{item:tri-2story-23}), (\ref{item:tri-join}),
  (\ref{item:tri-split}), (\ref{item:tri-collision}),
  (\ref{item:tri-v12-zero}) and the odd case of
  (\ref{item:tri-shuffle}) corresponds to an odd number of ends of
  $\cM^B(\x,\y,\Theta_o\semico\vec{\rhos},\vec{\rhos}^H)$. 
  (By contrast, we will see in the proof of
  Proposition~\ref{prop:v12-ends-match} that not all curves of types
  (\ref{item:tri-v12-one}) and (\ref{item:tri-v12-two}) correspond to
  ends of $\cM^B(\x,\y,\Theta_o;\vec{\rhos},\vec{\rhos}^H)$.)
\end{proposition}
\begin{proof}
  That any degeneration has one of forms
  (\ref{item:tri-2story-31})--(\ref{item:tri-shuffle}) or
  corresponds to a degeneration at $v_{12}$ follows from an argument
  analogous to the proof of Theorem~\ref{thm:master_equation}. That
  the degenerations at $v_{12}$ have the specified forms follows from
  Proposition~\ref{prop:glue-v12}. 
  For example, for
  case~(\ref{item:tri-v12-one}), Proposition~\ref{prop:glue-v12}
  implies that without the height constraints imposed by $\rhos_m$,
  the glued moduli space has dimension 
  \[
  \dim(\cM^{B_1}(\Theta_o,\w_\pd^+))+\dim(\cM^{B_2}(\x,\y,\w_\pd^-;(\rhos_1,\dots,\rhos_{m-1}),\vec{\rhos}^{H}))+1.
  \]
  So, the dimension of
  $\cM^B(\x,\y,\Theta_o;\vec{\rhos},\vec{\rhos}^H)$ is
  \[
  \dim(\cM^{B_1}(\Theta_o,\w_\pd^+))+\dim(\cM^{B_2}(\x,\y,\w_\pd^-;(\rhos_1,\dots,\rhos_{m-1}),\vec{\rhos}^{H}))+1-|\rhos_m|+1.
  \]
  By Corollary~\ref{cor:small-ind-v12-curves}, this is at least 
  \[
  |\rhos_m|-1+0+1-|\rhos_m|+1=1,
  \]
  with equality only if $\dim(\cM^{B_1}(\Theta_o,\w_\pd^+))=|\rhos_m|-1$.

  The converse follows from gluing
  arguments similar to those used to prove
  Theorem~\ref{thm:master_equation} (and Proposition~\ref{prop:glue-v12}).
\end{proof}

As a warm-up to our next proposition, we have the following:
\begin{lemma}\label{lem:tri-v12-zero-cancel} For any $B$,
  $\vec{\rhos}$ and $\vec{\rhos}^H$ so that
  $\ind(B,(\vec{\rhos},\vec{\rhos}^H))=1$, the moduli space
  $\cM^B(\x,\y,\penalty 500\Theta_o\semico\vec{\rhos},\vec{\rhos}^H)$ has an even number
  of ends of type~(\ref{item:tri-v12-zero}).
\end{lemma}
\begin{proof}
  By Corollary~\ref{cor:small-ind-v12-curves}, each end of
  $\cM^B(\x,\y,\Theta_o;\vec{\rhos},\vec{\rhos}^H)$ of
  type~(\ref{item:tri-v12-zero}) corresponds to a pair of curves
  $(u_{\alpha^H,\alpha},u_{\alpha,\alpha^H,\beta})$ where
  $u_{\alpha^H,\alpha}$ is a bigon with domain some $R_{i,j}^c$. There
  is a corresponding end
  $(u'_{\alpha^H,\alpha},u_{\alpha,\alpha^H,\beta})$ where
  $u'_{\alpha^H,\alpha}$ is a bigon with domain $R_{i,j\pm
    1}^c$. The result follows.
\end{proof}

The following proposition is similar in spirit to
Lemma~\ref{lem:tri-v12-zero-cancel}, but for the bigons with domains
$R_{i,j}^a$ and the annulus with domain $R_{1,1}^a+R_{1,1}^c$ or $R_{1,1}^a+R_{1,2}^c$.
\begin{proposition}\label{prop:v12-ends-match}
  Fix $B\in\pi_2(\x,\y)$ and a sequence of non-empty sets of Reeb chords
  $(\rhos_1,\dots,\rhos_m)$ such that
  $\ind(B,(\rhos_1,\dots,\rhos_n))=1$. Let
  $o_i=o(\vec{\x},(\rhos_1,\dots,\rhos_i))$.  For a sufficiently
  generic almost complex structure $J$, the number of ends of the moduli space
  \[
  \cM^B_{i,i+1}\coloneqq\cM^B(\x,\y,\Theta_{o_i};(\rhos_1,\dots,\rhos_i),(\rhos_{i+1}^H,\dots,\rhos_n^H))
  \]
  where $\rhos_i$ approaches~$v_{12}$ has the same parity as the
  number of ends of
  \[
    \cM^B_{i-1,i}\coloneqq\cM^B(\x,\y,\Theta_{o_{i-1}};(\rhos_1,\dots,\rhos_{i-1}),(\rhos_{i}^H,\dots,\rhos_n^H))
  \]
  where $\rhos_i^H$ approaches~$v_{12}$.
\end{proposition}
(Since we have fairly weak control over the behavior of $\cM^B$ near
the boundary, the parity of the number of ends might deserve some
explanation. For each $\epsilon>0$, consider
$\{u\in\cM^B_{i,i+1}\mid
d(\pi_\Delta\circ u(P_i),v_{12})\geq \epsilon\}$, where $d$ denotes
the distance in $\Delta$, viewed as the unit disk in $\CC$ with three
boundary punctures, and $P_i$ denotes the punctures of the source of
$u$ mapped to $\rhos_i$. For generic $\epsilon$, this is a manifold
with boundary. The number of boundary components of this manifold
modulo $2$ is
independent of $\epsilon$ if $\epsilon$ is generic and sufficiently
small. This is the (modulo $2$) number of ends of
$\cM^B_{i,i+1}$ where $\rhos_i$ approaches $v_{12}$.)

\begin{proof}
  We start with the case that $\rhos_i$ consists of a single Reeb
  chord $\rho$; the argument in this case is illustrated schematically
  in Figure~\ref{fig:bigon-degen}. By Proposition~\ref{prop:tri-master}, ends of
  $\cM^B_{i,i+1}$ at which $\rhos_i$ approaches $v_{12}$ correspond to
  a pair of curves
  $(u_{\alpha^H,\alpha},u_{\alpha,\alpha^H,\beta})$. Proposition~\ref{prop:tri-master}
  and Corollary~\ref{cor:small-ind-v12-curves} imply that the curve
  $u_{\alpha^H,\alpha}$ is either a bigon with domain $R_{i,j}^a$ (for
  some $i,j$) or an annulus with domain $R_{1,1}^a+R_{1,1}^c$ or
  $R_{1,1}^a+R_{1,2}^c$. In either case, there is algebraically one
  such curve $u_{\alpha^H,\alpha}$ in the moduli space
  $\cM^B(\Theta_o,\w_\pd^+)$; this is obvious in the case of bigons,
  and follows from \cite[Lemma~9.4]{OS04:HolomorphicDisks} in the case
  of annuli. Which domain $R_{i,j}^a$ occurs in the domain is determined by
  the chord $\rho$: the region $R_{i,j}^a$ has the point $\rho^+$ in
  its boundary. Proposition~\ref{prop:glue-v12} implies that near
  $(u_{\alpha^H,\alpha},u_{\alpha,\alpha^H,\beta})$, the compactified
  moduli space $\ocM^B_{i,i+1}$ is homeomorphic to $[0,\epsilon)$.
  
  \begin{figure}
    \includegraphics[scale=.83333]{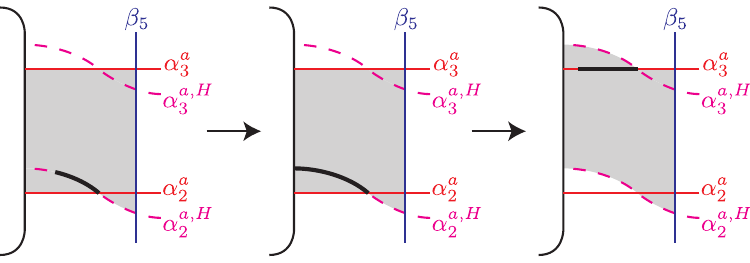}
    \caption[Illustration of a boundary bigon degenerating, in proof
    of invariance]{\textbf{An illustration of a boundary bigon degenerating, and a
      different one being glued on the other
      side.} The domain (shaded) changes slightly in the process.}\label{fig:bigon-degen}
  \end{figure}

  There is a corresponding curve $u'_{\alpha^H,\alpha}$, again a bigon
  $R_{i,j}$ or annulus $R_{1,1}^a+R_{1,1}^c$ or $R_{1,1}^a+R_{1,2}^c$,
  with $\rho^-$ in the boundary of the domain of
  $u'_{\alpha^H,\alpha}$. The pair
  $(u_{\alpha^H,\alpha},u_{\alpha,\alpha^H,\beta})$ is a point in
  $\bdy\ocM^B_{i-i,i}$. Again, Proposition~\ref{prop:glue-v12} implies
  that near $(u_{\alpha^H,\alpha},u_{\alpha,\alpha^H,\beta})$ the
  moduli space $\ocM^B_{i-1,i}$ is homeomorphic to $[0,\epsilon)$. The
  result, in the case that $\rhos_i$ consists of a single chord, follows.

  As a warm-up to the general case, 
  we consider the case that $\rhos_i$ contains exactly two Reeb
  chords, $\rho_1$ and $\rho_2$. Our argument is analogous
  to the proof of Proposition~\ref{prop:gluing_degree_one}. By
  Proposition~\ref{prop:tri-master}, when $\rhos_i$ approaches the
  puncture, the curve degenerates a curve
  $u_{\alpha^H,\alpha}$
  in $(\Sigma,\alphas^H,\alphas,z)$ with two non-trivial
  components. By Corollary~\ref{cor:small-ind-v12-curves}, each of
  these components maps either to the bigon 
  $R_{i,j}^a$ or to the annulus $R_{1,1}^a+R_{1,1}^c$ or
  $R_{1,1}^a+R_{1,2}^c$.  There is also a leftover
  component~$u_{\alpha,\alpha^H,\beta}$ in
  $\Sigma\times\Delta$. Let $\cM(u_{\alpha^H,\alpha})$ denote the
  moduli space of curves in
  $(\Sigma\times[0,1]\times\RR,\alphas\times\{0\}\times\RR\cup\alphas^H\times\{1\}\times\RR)$
  of which $u_{\alpha^H,\alpha}$ is a member. The moduli space
  $\cM(u_{\alpha^H,\alpha})$ is homeomorphic to (an odd number of
  copies of) $\RR$; different
  points are gotten by sliding the two non-trivial components of
  $u_{\alpha^H,\alpha}$ relative to each other.

  By Proposition~\ref{prop:glue-v12}, a neighborhood $\cU$ of
  $\cM(u_{\alpha^H,\alpha})\times \{u_{\alpha,\alpha^H,\beta}\}$ in 
  \begin{multline*}
    \ocM^B(\x,\y,\Theta_{o_i};(\rhos_1,\dots,\rhos_{i-1},\{\rho_1\},\{\rho_2\}),(\rhos_{i+1}^H,\dots,\rhos_n^H))
    \\
    \cup
    \ocM^B(\x,\y,\Theta_{o_i};(\rhos_1,\dots,\rhos_{i-1},\{\rho_2\},\{\rho_1\}),(\rhos_{i+1}^H,\dots,\rhos_n^H))
  \end{multline*}
  is homeomorphic to $\RR\times [0,\epsilon)$.
  Choose this homeomorphism so that $\delta\in(0,\epsilon)$ is the distance in
  $\Delta$ between $v_{12}$ and the closer of $\pi_\Delta(u(P_1))$ and
  $\pi_\Delta(u(P_2))$, where $P_1$ and $P_2$ are the punctures
  corresponding to $\rho_1$ and $\rho_2$, respectively.
  Given
  $\delta\in(0,\epsilon)$, for $t$ sufficiently small the point
  $(t,\delta)$ lies in
  $\ocM^B(\x,\y,\Theta_{o_i}\semico(\rhos_1,\dots,\rhos_{i-1},\{\rho_1\},\{\rho_2\}),\penalty400(\rhos_{i+1}^H,\dots,\rhos_n^H))$,
  say, while for $t$ sufficiently large the point $(t,\delta)$ lies in
  $\ocM^B(\x,\y,\Theta_{o_i}\semico(\rhos_1,\dots,\rhos_{i-1},\{\rho_2\},\{\rho_1\}),\penalty400(\rhos_{i+1}^H,\dots,\rhos_n^H))$.
  It follows that $\cM^B_{i,i+1}$ has an odd number of ends in
  $\cM(u_{\alpha^H,\alpha})\times \{u_{\alpha,\alpha^H,\beta}\}$.

  For each of the components of $u_{\alpha^H,\alpha}$ there is a
  matching bigon or annulus, as in the single Reeb chord case; denote the union
  of the matching components by $u'_{\alpha^H,\alpha}$. (That is,
  $u'_{\alpha^H,\alpha}$ has two non-trivial components
  whose domains contain the points $\rho_1^-$
  and $\rho_2^-$, respectively.) A similar argument to the
  previous paragraph shows that $\cM^B_{i-1,i}$ has an odd number of
  ends in $\cM(u'_{\alpha^H,\alpha})\times
  \{u_{\alpha,\alpha^H,\beta}\}$. The result follows.

  For the general case, write $\rhos_i=\{\rho_1,\dots,\rho_m\}$.
  Consider a pair of curves
  $(u_{\alpha^H,\alpha},u_{\alpha,\alpha^H,\beta})\in
  \bdy\ocM^B_{i,i+1},$ as above. Let $\Source$ denote the source of
  $u_{\alpha^H,\alpha}$ and $\TSource$ the source of
  $u_{\alpha,\alpha^H,\beta}$. Let $\vec{P},\vec{P}^H$ be the ordered
  partitions of punctures on $u_{\alpha,\alpha^H,\beta}$ along $e_1$
  and $e_2$ induced by $(\rhos_1,\dots,\rhos_{i-1})$ and
  $(\rhos_{i+1}^H,\dots,\rhos_n^H)$. Let $r_1,\dots,r_m$ denote the
  punctures of $\Source\glue\TSource$ mapped to
  $\rho_1,\dots,\rho_m$.
  Consider the moduli space
  \[
  \bigcup_{\sigma\in S_m}
  \cM^B(\x,\y,\Theta_o;\Source\glue\TSource;
  (\vec{P},(r_{\sigma(1)},\dots,r_{\sigma(m)})),\vec{P}^H)
  \]
  There is an evaluation map
  \[
  \ev\co \bigcup_{\sigma\in S_m}
  \cM^B(\x,\y,\Theta_o;\Source\glue\TSource;
  (\vec{P},(r_{\sigma(1)},\dots,r_{\sigma(m)})),\vec{P}^H)\to (e_1)^m
  \]
  given by
  \[
  \ev(u)=(\pi_\Delta\circ u)(r_1,\dots,r_m).
  \]
  The moduli space $\cM^B_{i,i+1}$ is
  $\ev^{-1}\bigl(\{(x,x,\dots,x)\mid x\in e_1\}\bigr)$.
  We can also evaluate at only some of the punctures; let 
  \[
  \ev_{[1,j]}(u)=(\pi_\Delta\circ u)(r_1,\dots,r_j)
  \]

  Let
  \begin{multline*}
  \Gamma\co \cM^B(\Theta_o,\w_\pd^+;\Source)\times
  \cM^B(\x,\y,\w_\pd^-;\TSource;\vec{P},\vec{P}^H)\times[0,\epsilon)\\
  \to 
  \bigcup_{\sigma\in S_m}
  \ocM^B(\x,\y,\Theta_o;\Source\glue\TSource;
  (\vec{P},(r_{\sigma(1)},\dots,r_{\sigma(m)})),\vec{P}^H)
  \end{multline*}
  denote the gluing map from Proposition~\ref{prop:glue-v12}. Recall
  from Corollary~\ref{cor:small-ind-v12-curves} that there is a
  bijection between the nontrivial components of $u_{\alpha^H,\alpha}$ and the
  punctures $r_i$. 
  Gluing at the punctures $r_1,\dots,r_j$ (i.e.,
  translating the components corresponding to $r_{j+1},\dots,r_m$ to
  $-\infty$ and then gluing at $v_{12}$) gives a gluing map 
  \begin{multline*}
  \Gamma_{[1,j]}\co V_{[1,j]}\times
  \cM^B(\x,\y,\w_\pd^-;\TSource;\vec{P},\vec{P}^H)\times[0,\epsilon)\\
  \to  \bigcup_{\sigma\in S_m}
  \ocM^B(\x,\y,\Theta_o;\Source\glue\TSource;
  (\vec{P},(r_{\sigma(1)},\dots,r_{\sigma(m)})),\vec{P}^H)
  \end{multline*}
  for the corresponding stratum $V_{[1,j]}$ of 
  $\cM^B(\Theta_o,\w_\pd^+;\Source)$ (where the curve at $v_{12}$ has
  at least two stories, and the punctures $r_1,\dots,r_j$ are on a
  higher story than the punctures $r_{j+1},\dots,r_m$).
  Write $V_{[1,j]}=W_{[1,j]}\times W_{[j+1,m]}$ where $W_{[1,j]}$ is the moduli
  space corresponding to the punctures $r_1,\dots,r_j$ and
  $W_{[j+1,m]}$ is the moduli space corresponding to the punctures
  $r_{j+1},\dots,r_{m}$.

  We will prove the following statement
  by induction on $m$:
  \begin{itemize}
  \item Given any $(u_{\alpha^H,\alpha},u_{\alpha,\alpha^H,\beta})$ as above
    there is an open neighborhood $U$ in $(e_1)^j$ containing
    $\{(v_{12},v_{12},\dots,v_{12})\}$ and an
    $\epsilon>0$ so that 
    \[
    \ev_{[1,j]}\circ \Gamma_{[1,j]}\co W_{[1,j]}\times 
    \{u_{\alpha,\alpha^H,\beta}\}\times(0,\epsilon)\to (e_1)^j
    \]
    maps with odd degree onto $U$.
  \end{itemize}
  Given this statement for the case $j=m$, taking the preimage of the
  diagonal implies the proposition.

  The case $j=1$ is clear from (the proof of) Proposition~\ref{prop:glue-v12}.
  The inductive step follows from
  Lemma~\ref{lemma:stratified-degree-one} with $X=W_{[1,j]}$ and
  $f=\ev_{[1,j]}\circ \Gamma_{[1,j]}$. The fact that the union of the
  top two strata of $W_{[1,j]}$ is a manifold with boundary is clear
  from the structure of $W_{[1,j]}$ (analyzed completely in
  Corollary~\ref{cor:small-ind-v12-curves}) or, alternatively, follows
  from gluing theory. That $f$ is proper near $0$ follows from the
  gluing construction. The
  inductive hypothesis gives the fact that the restriction of $f$ to
  the (codimension $1$) facet where $r_1$ is on the bottom story and $r_2,\dots,r_{m}$
  are on the top story (adjacent to $v_{12}$), with image in
  $(e_1)^{m-1}$, has odd degree near $0$.
\end{proof}

\subsubsection{The map associated to a handleslide is a chain map}\label{sec:handleslide-is-chain-map}
\begin{proposition}\label{Prop:FaaHbChainMap}$F_{\alphas,\alphas^H,\betas}$
  is a chain map, i.e., $F\circ\bdy+\bdy\circ F=0$.
\end{proposition}

\begin{proof}
  This follows from Proposition~\ref{prop:tri-master},
  Lemma~\ref{lem:tri-v12-zero-cancel} and
  Proposition~\ref{prop:v12-ends-match}, similarly to the proof of
  Proposition~\ref{prop:typeD-d2}. Ends of
  type~(\ref{item:tri-2story-31}) correspond to
  $F_{\alphas,\alphas^H,\betas}\circ\bdy(\x)$. Ends of
  type~(\ref{item:tri-2story-23}) correspond to part of $\bdy\circ
  F_{\alphas,\alphas^H,\betas}(\x)$. Ends of type~(\ref{item:tri-split})
  correspond to the rest of $\bdy\circ
  F_{\alphas,\alphas^H,\betas}(\x)$. Ends of
  type~(\ref{item:tri-collision}) cancel with each other and with ends
  of type~(\ref{item:tri-join}), as in the proof of
  Proposition~\ref{prop:typeD-d2}. Ends of
  type~(\ref{item:tri-v12-zero}) cancel in pairs by
  Lemma~\ref{lem:tri-v12-zero-cancel}. Ends of
  type~(\ref{item:tri-v12-one}) cancel with ends of
  type~(\ref{item:tri-v12-two}) by
  Proposition~\ref{prop:v12-ends-match}. The result follows.
\end{proof}

\begin{figure}
  \includegraphics[scale=.83333]{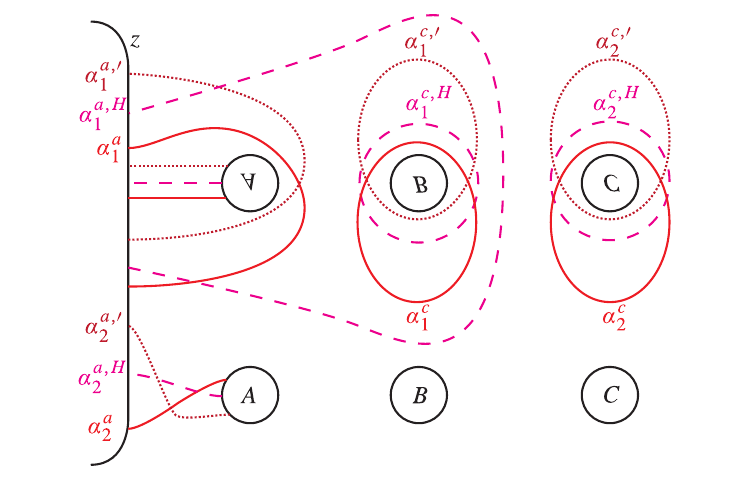}
  \caption[Curves $\alphas$, $\alphas^H$ and $\alphas'$ in
    $\Sigma$]{\textbf{The curves $\alphas$, $\alphas^H$ and $\alphas'$ in
    $\Sigma$.}  The curves in $\alphas$ are solid, those in $\alphas^H$
  are dashed, and those in $\alphas'$ are dotted.}
\label{fig:alphasalphasHalphasP}
\end{figure}

To show that $F_{\alphas,\alphas^H,\betas}$ is an isomorphism, we wish
to construct a chain map the other direction and a chain homotopy
between their composition and the
identity map.  A count of triangles in $(\Sigma,\alphas^H,\alphas,\betas)$
does not work as desired, as the small bigons $R_{i,j}^a$ (and
annuli $R_{1,1}^a+R_{1,1}^c$ and $R_{1,1}^a+R_{1,2}^c$) are not
holomorphic: their orientations are consistent with
$(\alphas,\alphas^H,\betas)$, not $(\alphas^H,\alphas,\betas)$.
Instead we define a third set of circles, which we denote
$\gls*{alphasprime}$.
Obtain
$\gls*{alphaicprime}$
from $\alpha_i^c$ by
performing a small Hamiltonian perturbation, so that $\alpha_i^{c,\prime}$
intersects each of $\alpha_i^{c}$ and $\alpha_i^{c,H}$
transversely in two points and is disjoint from all other curves in
$\alphas\cup\alphas^H$. Let 
$\gls*{alphaiaprime}$
be an isotopic translate
of $\alpha_i^a$, intersecting each of $\alpha_i^a$ and
$\alpha_i^{a,H}$ in a single point, and such that there are two short
Reeb chords in $\bdy\widebar{\Sigma}$ running from
$\alpha_i^{a,\prime}$ to $\alpha_i^{a,H}$ (and hence two slightly
longer Reeb chords running from $\alpha_i^{a,\prime}$ to
$\alpha_i^a$). The isotopy is chosen small enough that
$\alpha_i^{a,\prime}$ is disjoint from all other curves in
$\alphas\cup\alphas^H$. See
Figure~\ref{fig:alphasalphasHalphasP}. 
Now there is a chain map $F_{\alphas^H,\alphas',\betas}$ defined
using $(\Sigma,\alphas^H,\alphas',\betas,z)$, in exactly the same way
that $F_{\alphas,\alphas^H,\betas}$ was defined using
$(\Sigma,\alphas,\alphas^H,\betas,z)$.

\subsubsection{Completion of the proof of handleslide invariance}\label{sec:end-handleslide-pf}
To prove handleslide invariance, we also need to consider the map
$F_{\alphas,\alphas',\betas}$ defined in the same
way as $F_{\alphas,\alphas^H,\betas}$ but with $\alphas'$ in place of~$\alphas^H$.
\begin{proposition}\label{Prop:FaapbChainIso}The map
  $F_{\alphas,\alphas',\betas}$ is an isomorphism of differential
  $\Alg$-modules.
\end{proposition}
\begin{proof}
The proof that $F_{\alphas,\alphas',\betas}$ is a chain map is
parallel to, but marginally easier than, the proof of
Proposition~\ref{Prop:FaaHbChainMap}, the only difference being that
the discussion of annuli is irrelevant here.

Next we argue that $F_{\alphas,\alphas',\betas}$ is in fact an
isomorphism. The proof, which makes use of the
\index{energy filtration}%
energy filtration, is a straightforward analogue of the proof
for closed 3-manifolds~\cite[Proposition
9.8]{OS04:HolomorphicDisks}. We will use the
following standard lemma.
\begin{lemma}\label{lemma:filtrations-good}
  Let $F\co A\to B$ be a map of filtered groups which is
  decomposed as a sum $F=F_0+\ell$ where $F_0$ is a
  filtration-preserving isomorphism and $\ell$ has strictly lower
  order than $F_0$. Suppose that the filtration on $B$ is bounded
  below. Then $F$ is an isomorphism of groups.
\end{lemma}

Choose an area form $\Area$ on $\Sigma$ such that if $D$ is a
periodic domain in $(\Sigma,\alphas,\betas)$ then
$\Area(D)=0$. (The existence of $\Area$ is guaranteed by
Lemma~\ref{lemma:admiss-reform}.) Arrange also that periodic domains in
$(\Sigma,\alphas',\alphas)$ have area~$0$, and that the
boundary bigons all have the same area. It follows that
for $D$ a periodic domain in $(\Sigma,\alphas',\betas)$,
$\Area(D)=0$ as well.

For each $\spin^c$ structure $\s$ on $Y$, we define a map
$\mathcal{F}$ from
$\S(\Sigma,\alphas,\betas,\s)$ to~$\RR$. Choose one generator
$\x_0\in\S(\Sigma,\alphas,\betas,\s)$ and declare $\mathcal{F}(\x_0)=0$.
Then, for any other $\x\in\S(\Sigma,\alphas,\betas,\s)$ pick
$A_{\x_0,\x}\in\pi_2(\x_0,\x)$ and define $\mathcal{F}(\x)=-\mathord{\Area}(A_{\x_0,\x})$. Since
periodic domains have area $0$, $\mathcal{F}(\x)$ is independent of
the choice of $A_{\x_0,\x}$. More generally, for $a\in\Alg$ such that
$a\x\neq0$, define $\mathcal{F}(a\x)=\mathcal{F}(\x)$. Since
holomorphic curves have positive domains, the map $\mathcal{F}$ induces a
filtration on $\CFDa(\Sigma,\alphas,\betas)$.

There is an obvious identification between
$\S(\Sigma,\alphas,\betas,\s)$ and
$\S(\Sigma,\alphas',\betas,\s)$. As before, let $\x'$ denote the
generator corresponding to~$\x$, and let
$T_{\x}\in\pi_2(\x,\x',\Theta_{o(\x)})$ be the canonical small
triangle.
Define a map $\mathcal{F}'_0\co
\S(\Sigma,\alphas',\betas,\s)\to\RR$ in the same way as $\mathcal{F}$,
still using the area form $\Area$ but using $\x_0'$ in place of
$\x_0$. Define $\mathcal{F}'\co
\S(\Sigma,\alphas',\betas,\s)\to\RR$ by
$\mathcal{F}'=\mathcal{F}'_0-\Area(T_{\x_0})$. By extending $\mathcal{F}'$
by $\mathcal{F}'(a\y)=\mathcal{F}'(\y)$, we get
a filtration~$\mathcal{F}'$ on $\CFDa(\Sigma,\alphas',\betas)$.

To see that $F_{\alphas,\alphas',\betas}$ respects the filtrations
$\mathcal{F}$ and $\mathcal{F}'$, suppose there is a term $a\y$ in
$F_{\alphas,\alphas',\betas}(\x)$, and thus a positive domain
$B \in \pi_2(\x,\y',\Theta_o)$, and pick a domain
$A_{\x_0,\x}\in\pi_2(\x_0,\x)$.  Then, for an appropriate $o' =
o(\x_0)$ and domain
$A_{o,o'}\in\pi_2(\Theta_o,\Theta_{o'})$, the domain
$A_{\x_0,\x}+B-T_{\x_0}-A_{o,o'}$ connects $\x_0'$ and $\y'$.  Our
assumptions guarantee that $\Area(A_{o,o'}) = 0$, so
\begin{align*}
\cF'(\y') &=-\Area(A_{\x_0,\x}+B-T_{\x_0}-A_{o,o'})-\Area(T_{\x_0})\\
&= -\Area(A_{\x_0,\x}) - \Area(B)\\ & < \cF(\x),
\end{align*}
as desired.

  Now, if we choose $\Area$ so that for all
  $\x$ the triangle
  $T_\x$ is the unique triangle of minimal area connecting $\x$, $\y'$
  and $\Theta_o$ for any $\y'$ and $o$ (this is easily accomplished),
  then the top order part
  of $F_{\alphas,\alphas',\betas}$ with respect to $\mathcal{F}$ and
  $\mathcal{F}'$ is simply the map $\x\mapsto\x'$. In particular, the
  top order part of $F_{\alphas,\alphas',\betas}$ is a group isomorphism.

  Consequently, it follows from Lemma~\ref{lemma:filtrations-good}
  that the $\Alg$-module homomorphism $F_{\alphas,\alphas',\betas}$ is
  an $\Field$ vector space isomorphism
  $\CFDa(\Sigma,\alphas,\betas,z)\to\CFDa(\Sigma,\alphas',\betas,z)$,
  and hence also a $\Alg$-module isomorphism.
\end{proof}

\begin{remark}
  We could weaken our assumptions and prove handleslide invariance for
  Heegaard diagrams~$\HD$ that are provincially admissible (rather
  than admissible) by being slightly more clever in the choice of
  filtration~$\cF$, as follows. Pick a map $f\co H_1(\bdy\bSigma,\CircPts) \to
  \RR$ so that for any periodic domain~$B$ for
  $(\Sigma,\alphas,\betas)$,
  $f(\bdy^\bdy B) +
  \Area(B) = 0$.  (This is possible for any provincially admissible
  diagram.)  Then define $\cF(\x) = \Area(A_{\x_0,\x})
  + f(\bdy^\bdy A_{\x_0,\x})$, and make similar other adjustments.
  Alternatively, we can do an isotopy to make $\HD$ admissible, do the
  handleslide, and then do an inverse isotopy.
\end{remark}

\begin{figure}
    \includegraphics[scale=.82]{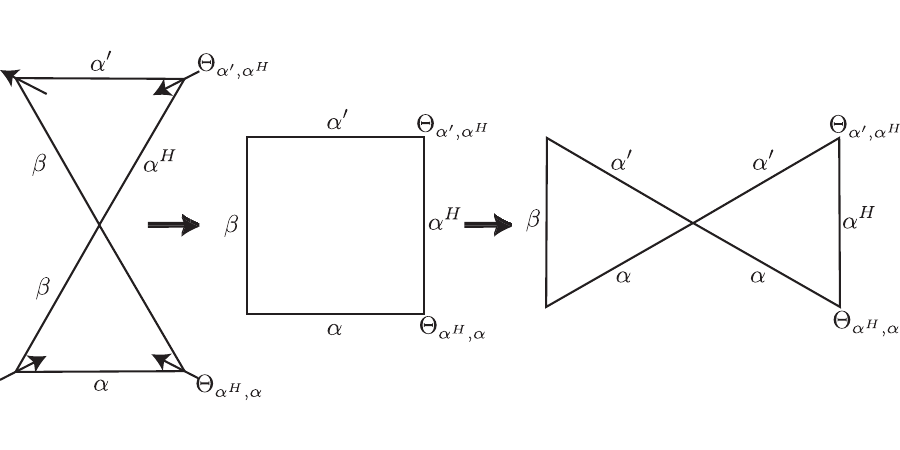}
    \caption[Two compositions of triangles, as degenerations of a
    rectangle]{\textbf{Two compositions of triangles, as degenerations of a
    rectangle.} The bold arrows indicate a path in the compactified moduli
  space of rectangles, starting at a configuration giving $F_{\alphas^H,\alphas',\betas}\circ
  F_{\alphas,\alphas^H,\betas}$ and
  ending at a configuration giving
  $F_{\alphas,\alphas',\betas}$.}\label{fig:rectangle-degen}
  \end{figure}

\begin{proposition}\label{prop:triangle-assoc}The map $F_{\alphas,\alphas',\betas}$ is chain
  homotopic to $F_{\alphas^H,\alphas',\betas}\circ
  F_{\alphas,\alphas^H,\betas}$.
\end{proposition}
\begin{proof}
  The proof is in two steps.  First, the composition
  $F_{\alphas^H,\alphas',\betas}\circ F_{\alphas,\alphas^H,\betas}$
  is a count of pairs of triangles, in
  $(\Sigma,\alphas,\alphas^H,\betas,z)$ and
  $(\Sigma,\alphas^H,\alphas',\betas,z)$. This map is chain homotopic to
  a count of pairs of triangles in $(\Sigma,\alphas,\alphas',\betas,z)$
  and $(\Sigma,\alphas,\alphas^H,\alphas',z)$, via a chain homotopy
  counting (index $-1$) holomorphic curves in a path of almost complex
  structures on the product of $\Sigma$ and a rectangle.
  (See
  Figure~\ref{fig:rectangle-degen} for a schematic illustration.) The
  reader unfamiliar with such arguments is directed to any
  of~\cite[Theorem~8.16]{OS04:HolomorphicDisks},
  \cite[Proposition~10.29]{Lipshitz06:CylindricalHF} or
  \cite[Lemma~7.2.5]{Lipshitz06:BorderedHF}.
  The only new types of degeneration are curves in either
  $(\Sigma,\alphas^H,\alphas,z)$ or $(\Sigma,\alphas',\alphas^H,z)$ touching
  $\bdy\bSigma$, either small bigons or annuli.  As in
  Proposition~\ref{Prop:FaaHbChainMap}, these degenerations
  cancel in pairs.

  Second, by direct computation, the only rigid triangles in
  $(\Sigma,\alphas,\alphas^H,\alphas',z)$ asymptotic to
  $\Theta_{o,\alphas^H,\alphas}$ and $\Theta_{o,\alphas',\alphas^H}$ at two
  of the corners are asymptotic to $\Theta_{o,\alphas',\alphas}$ at the
  third corner---and that there is algebraically a single such
  triangle. Again, the reader is directed to~\cite[Lemma
  9.7]{OS04:HolomorphicDisks} for details---the proof there applies
  without change in our case.
\end{proof}

\begin{proposition}
  \label{prop:HandleslideInvariance}
  The $\Alg$-modules $\CFDa(\Sigma,\alphas,\betas,z)$ and
  $\CFDa(\Sigma,\alphas^H,\betas,z)$ are homotopy equivalent as \dg
  modules.
\end{proposition}
\begin{proof}
If we arrange that $\alphas'$ and $\alphas$ are close enough then the
obvious identification of generators of
$\CFDa(\Sigma,\alphas,\betas,z)$ and $\CFDa(\Sigma,\alphas',\betas,z)$
induces an isomorphism of $\Alg$-modules (because all of the moduli
spaces counted in the definitions of $\partial$ are the same). Hence, by
Propositions~\ref{Prop:FaapbChainIso} and~\ref{prop:triangle-assoc},
the composition of $F_{\alphas^H,\alphas',\betas}\circ
F_{\alphas,\alphas^H,\betas}$ and the evident identification of
generators is chain homotopic to an automorphism of
$\CFDa(\Sigma,\alphas,\betas,z)$.

Choose a fourth set of curves $\alphas^{H,\prime}$ isotopic to
$\alphas^H$ and define another map
$F_{\alphas',\alphas^{H,\prime},\betas}\co\CFDa(\Sigma,\alphas',\betas,z)\to\CFDa(\Sigma,\alphas^{H,\prime},\betas,z)$.
The same argument shows that the composition of
$F_{\alphas',\alphas^{H,\prime},\betas}\circ
F_{\alphas^H,\alphas',\betas}$ and the evident identification of generators
is chain homotopic to an automorphism of
$\CFDa(\Sigma,\alphas^H,\betas,z)$. The result follows.
\end{proof}

\subsection{Completion of the invariance proof}\label{sec:all-together}

\begin{proof}[Proof of Theorem~\ref{thm:D-invariance}]
  In Section~\ref{sec:CFD-cx-str-change}, we saw that $\CFDa(\HD,\spinc;J)$ is independent
  of the choice of generic almost-complex structure $J$ in its definition.
  Suppose now that $\HD$ and $\HD'$ are two Heegaard diagrams which
  represent equivalent bordered three-manifolds (and in particular $-\partial
  \HD=\PMC=-\partial \HD'$). Then, according to
  Proposition~\ref{prop:heegaard-moves}, we can pass from $\HD$ to
  $\HD'$ by a sequence of Heegaard moves.  Invariance under isotopies
  has not been spelled out, but this can be done by writing down the
  continuation maps for isotopies modeled on the continuation maps for
  variation of almost-complex structures as in
  Equation~\eqref{eq:ContinuationMap}.
  Invariance under handleslides of an arc over a circle was verified
  in Proposition~\ref{prop:HandleslideInvariance} in the case that the
  diagrams are admissible; in the general case, first do an isotopy to
  make the diagrams admissible, then a handleslide, and then undo the
  isotopy. The case of
  handleslides of a circle over a circle is similar but easier, and is
  left to the reader.
  Invariance under stabilization is straightforward: by stabilizing
  near the basepoint, 
  one can see that the chain complexes before and after the
  stabilization are isomorphic 
  (just as in the closed case for $\HFa$ \cite[Theorem~10.1]{OS04:HolomorphicDisks}).
\end{proof}

In view of Theorem~\ref{thm:D-invariance}, we can now refer to $\CFDa$
without reference to a Heegaard diagram, defining $\CFDa(Y,\phi,\s)$
(where $\phi\co F(-\PMC) \to \partial Y$) to be $\CFDa(\HD,\s)$, where $\HD$ is any
bordered Heegaard diagram representing the bordered three\hyp manifold
$(Y,\phi)$. We will sometimes drop the
map $\phi$ from the notation, writing 
$\gls*{CFDspincY}$
to denote
$\CFDa(Y,\phi\co {F(-\PMC)}\to {\bdy Y},\spinc)$.

To obtain a topological invariant of the bordered three-manifold 
without reference to a $\SpinC$ structure, we simply sum over all $\SpinC$ structures, defining
$$
\gls*{CFDY}
=\bigoplus_{\s\in\SpinC(Y)} \CFDa(Y,\s).
$$ 

\begin{proof}[Proof of Theorem~\ref{intro:D-invariance}]
  With the above definition, this follows immediately from Theorem~\ref{thm:D-invariance}.
\end{proof}
\index{invariance!of $\CFDa$|)}%
\section{Twisted coefficients}
\label{sec:typeD-twisted}
As with the original Heegaard Floer homology \cite[Section
8]{OS04:HolDiskProperties}, there are versions of $\CFDa(\HD,\spinc)$
with twisted coefficients.  
\index{twisted coefficients!for $\CFDa$}%
We describe the (maximal) twist by
$H_2(Y,\bdy Y)$.  (Recall from Lemma~\ref{lem:pi2-h2} that, for any
$\x \in \S(\HD,\spinc)$, $\pi_2(\x,\x) \cong H_2(Y,\bdy Y)$.)  Besides
its intrinsic interest, we explain this version because, as we shall
see in Chapter~\ref{chap:gradings},  the grading on the
twisted theory is easier to understand than
the grading on the untwisted theory.

Fix a $\SpinC$ structure $\s$ over $Y$, and pick a base generator $\x_0=\x_0(\spinc) \in \gls*{generatorsS}$
representing $\s$.  (For $\spin^c$
structures~$\s$ without any generators, the twisted chain complex
is~$0$.)  Let $\gls*{generatorsUs}$ be the set of elements of the
form $e^{B_0}\x$, where $\x\in\S(\HD, \s)$,
$B_0 \in \pi_2(\x_0, \x)$, and 
$\gls*{formale}$
is a formal variable.
Let 
$\gls*{XofHUs}$
be the $\Field$ vector space spanned by
$\underline{\S}(\HD,\x_0)$. There is an action of $\Idem(\PtdMatchCirc)$
on $\underline{X}(\HD,\x_0)$, where
\begin{equation*}
  I(\SetS) \cdot e^{B_0}\x \coloneqq
  \begin{cases}
    e^{B_0}\x&I(\SetS) = I_D(\x)\\
    0&\text{otherwise,}
  \end{cases}
\end{equation*}
as in Equation~\eqref{eq:def-idem-DMod}.  

As a module
over $\Alg=\Alg(\PMC)$ (where $\PMC=-\partial\HD$), define the twisted
chain complex 
$\gls*{tCFD}$
to be $\Alg \otimes_\Idem \underline{X}(\HD,\x_0)$. 

There is a natural action of $H_2(Y,\bdy Y)$ on
$\tCFDa(\HD,\x_0)$, by composition with the corresponding
periodic domain: for $\alpha \in H_2(Y, \bdy Y)$ corresponding to $B_\alpha
\in\pi_2(\x_0,\x_0)$, define
\[
\alpha * e^{B_0}\x \coloneqq e^{B_\alpha* B_0}\x.
\]
As in the untwisted case, the action of $\Alg$ on $\tCFDa(\HD, \s)$
is given by left multiplication:
\[
a\cdot(b \otimes e^{B_0}\x) \coloneqq (ab) \otimes e^{B_0}\x.
\]
For the boundary operator, 
the contributions of the holomorphic curves connecting $\x$ to $\y$
are separated according to their homology classes:
\[
\partial(e^{B_0} \x) \coloneqq
\sum_{\y\in\S(\HD,\s)}\sum_{B\in\pi_2(\x,\y)}a^B_{\x,\y} e^{B_0* B}\y.
\]

It is convenient to think of the action by $H_2(Y,\partial Y)$
as a right action of the group-ring $\Field[H_2(Y,\partial Y)]$
on $\tCFDa(\HD,\x_0,\spinc)$.
Writing the group-ring as sums
$\sum_{\alpha\in H_2(Y,\partial Y)} n_\alpha e^\alpha$,
where $n_\alpha\in\Zmod{2}$ is zero for all but finitely
$\alpha\in H_2(Y,\partial Y)$, we let
$$e^{B_0}\x \cdot \Bigl(\sum_{\alpha} n_\alpha \cdot e^{\alpha}\Bigr)
= \sum_\alpha n_\alpha e^{B_\alpha * B_0} \x.$$

\begin{proposition}
  $\tCFDa(\HD,\x_0)$ is a differential 
  $\Alg$--$\Field[H_2(Y,\partial Y)]$ bimodule.
\end{proposition}

\begin{proof}
  Here, we are thinking of $\Field[H_2(Y,\partial Y)]$ as a \dg algebra with trivial
  differential, so that the statement  means that
  \begin{itemize}
  \item  $\tCFDa(\HD,\x,\spinc)$ is a chain complex, and
  \item for $a\in \Alg$, $m\in \tCFDa(\HD,\x_0,\spinc)$, 
    and $c\in\Field[H_2(Y,\partial Y)]$, 
    \begin{align*}
      (a \cdot m)\cdot c &= a\cdot (m\cdot c) \\
      \partial (a\cdot m \cdot c)&=(\partial a) m \cdot c + a\cdot (\partial m)\cdot c.
    \end{align*}
  \end{itemize}
  
  The proof of Proposition~\ref{prop:typeD-d2} readily adapts to
  show that $\tCFDa(\HD,\x_0,\spinc)$ is a chain complex.
  The action by $\Field[H_2(Y,\partial Y)]$ commutes both with the differential
  and with the action by $\Alg$. The result follows.
\end{proof}

\begin{theorem}\label{thm:tD-invariance}Up to chain homotopy equivalence, the
  differential module $\tCFDa(\HD,\x_0)$ depends on the generator $\x_0$ only through
  its induced $\SpinC$ structure $\spinc_z(\x)$ and on the Heegaard diagram
  only through its (equivalence class of) 
  underlying  bordered three-manifold~$(Y,\phi\co
  {-F(\PMC)}\to \partial Y)$
  (with $\PMC=-\partial \HD$). That
  is, if $(\HD,\x_0)$ and $(\HD',\x_0')$ are provincially
  admissible bordered Heegaard diagram for~$Y$ with $-\partial \HD=\PMC=-\partial \HD'$, and $\x_0\in\Gen(\HD,\spinc)$ and $\x_0'\in\Gen(\HD',\spinc)$ are
  base generators representing the same $\SpinC$ structure over $Y$,
  then
  $\tCFDa(\HD,\x_0)$ and
  $\tCFDa(\HD',\x_0')$ are homotopy equivalent
  $\Alg(\PtdMatchCirc)$--$\Field[H_2(Y,\partial Y)]$ bimodules.
\end{theorem}

\begin{proof}
  First, we address the dependence of $\tCFDa(\HD,\x_0)$ on the base
  generator~$\x_0$.  Suppose that $\HD$ is an admissible diagram, and that
  $\x_0$ and $\x_1$ are two different generators which represent the
  same $\SpinC$ structure over $Y$. Then, by
  Lemma~\ref{lem:SpinCStructures}, we can find some
  $b\in\pi_2(\x_0,\x_1)$. We use $b$ to induce a map
  $\underline{X}(\HD,\x_1)\to \underline{X}(\HD,\x_0)$ sending the
  generator $e^B \x$ to $e^{b*B}\x$. To see that this respects the
  action by $H_2(Y,\partial Y)$, note that for $\alpha\in
  H_2(Y,\partial Y)$, $b*B_\alpha=B_\alpha*b$ (the map $\alpha\to
  B_\alpha$ from $H_2(Y;\ZZ)$ to $\pi_2(\x_0,\x_0)$ gives a domain
  which does not depend on the basepoint $\x_0$). 
  The map $\underline{X}(\HD,\x_1)\to \underline{X}(\HD,\x_0)$ 
  extends to a
  bimodule map $\Phi_b\co
  \tCFDa(\HD,\x_1)\to\tCFDa(\HD,\x_0)$. If $\x_0=\x_1$ and
  $b\in\pi_2(\x_0,\x_1)$ is the trivial homology class then $\Phi_b$ is the
  identity. Moreover, given $b_0\in \pi_2(\x_0,\x_1)$ and
  $b_1\in\pi_2(\x_1,\x_2)$, $\Phi_{b_0*b_1}=\Phi_{b_2}\circ
  \Phi_{b_1}$. It follows that for any $b\in\pi_2(\x_0,\x_1)$,
  $\Phi_b$ is an isomorphism.

  To adapt the  proof of
  Theorem~\ref{thm:D-invariance} to the twisted setting, we must
  give twisted analogues of the
  continuation and triangle maps.  For example, for a suitable path of
  almost complex structures $J_r$, we define
  $${\underline F}^{J_r}\co\tCFDa(\HD,\spinc;J_0)\to\tCFDa(\HD,J_1)$$
  by adapting the map defined in Equation~\eqref{eq:ContinuationMap}:
  \[{\underline F}^{J_r}(e^{B_0}\x)=
  \sum_{\y\in\S(\HD,\spinc)}\,\sum_{B\in\pi_2(\x,\y)}\,
  \sum_{\{\vec{\rho}\,\mid\,\ind(B,\vec\rho) = 0\}}\!
  \#\left(\Mod^B(\x,\y;\vec{\rho};J)\right)a(-\vec\rho) e^{B_0*B}\y.
  \]
  We modify the triangle map from Equation~\eqref{eq:TriangleMap}
  similarly, writing
  $$
  {\underline F}_{\alphas,\alphas^H,\betas}(e^{B_0}\x)\coloneqq
  \sum_{\y}\sum_{B\in\pi_2(\x',\y)}\,
  \sum_{\{\vec{\rho}\,\mid\,\ind(B,\vec\rho) = 0\}}
  \#\left(\Mod^B(\x,\y,\Theta;\vec{\rho})\right)a(-\vec\rho) e^{B_0'* B}\y,
  $$
  where $B_0'\in\pi_2(\x_0',\x)$ is determined 
  (via Lemma~\ref{lem:triple-pi2-decomp}) by
  $$T_{\x}*_{13} B_0=T_{\x_0}*_{12} B_{\alphas,\alphas^{H}}
  *_{23} B_0'.$$
  
  With the maps in place, the proof Theorem~\ref{thm:D-invariance}
  adapts easily to the twisted setting.
\end{proof}


\chapter{Type \textalt{$A$}{A} modules}
\label{chap:type-a-mod}
We now turn to defining the type $A$ module $\CFAa$ associated
to a (provincially admissible) bordered Heegaard diagram~$\HD$. This uses more of the
moduli spaces from Chapter~\ref{chap:structure-moduli} than the
type~$D$ module does; in particular, we must consider moduli spaces of curves
where the (ordered) partition of Reeb chords has more than one chord
per part.

\section{Definition of the type \textalt{$A$}{A} module}\label{sec:def-CFA}
\index{type $A$ invariant|(}
Fix a bordered Heegaard diagram $\HD=(\Sigma,\alphas,\betas,z)$,
satisfying the provincial admissibility criterion of
Definition~\ref{def:provincial-admissibility}, and a $\SpinC$ structure $\s$ over $Y\!$.
The module $\CFAa(\HD,\s)$
is a right ($\Ainf$)
$\Alg(\PtdMatchCirc)$-module, where $\PMC=\bdy\HD$. Unlike $\CFDa$,
which is essentially free (projective) over $\Alg(\PMC)$, much of the
data of the Heegaard diagram is encoded in the multiplication on
$\CFAa(\HD,\s)$, as well as in certain higher products.

The module $\gls*{CFAspinc}$
is generated
over $\Field$ by the set of
generators $\S(\HD,\spinc)$. Given $\x\in\S(\HD,\spinc)$, define 
$\gls*{IAofx}\coloneqq
I(o(\x))$.  (Recall that $o(\x) \subset [2k]$ is the set of $\alpha$-arcs
occupied by~$\x$.)
We define a right action of
$\Idem(\PMC)$ on $\CFAa(\HD)$ by
\begin{equation}
  \label{eq:typeAactionI}
  \x\cdot I(\SetS) \coloneqq
  \begin{cases}
    \x&I_A(\x) = I(\SetS)\\
    0&\text{otherwise.}
  \end{cases}
\end{equation}
Again, the summands $\Alg(\PtdMatchCirc,i)$ of
$\Alg(\PtdMatchCirc)$ act trivially on $\CFAa(\HD,\s)$ for $i\neq 0$.

Recall from Section~\ref{sec:reeb-chords-def} that to a set $\rhos$ of
Reeb chords we associate an algebra element $a(\rhos)$.

\begin{lemma}\label{lem:strong-monotone-tensor}
  For $\x \in \S(\HD,\s)$ and $\vec\rhos$ a sequence of non-empty sets of Reeb
  chords, $(\x,\vec\rhos)$ is strongly boundary monotone
  (Definition~\ref{def:strong-monotonicity-P}) if and only if the
  tensor product
  \[
  I_A(\x) \otimes_\Idem a(\rhos_1) \otimes_\Idem\dots\otimes_\Idem
  a(\rhos_n)
  \]
  is not zero.
\end{lemma}
\begin{proof}
  In the notation of Definition~\ref{def:strong-monotonicity-P}, if $
  I_A(\x) \otimes_\Idem a(\rhos_1) \otimes_\Idem\dots\otimes_\Idem
  a(\rhos_i)\neq 0 $ then the right idempotent of $I_A(\x) \otimes_\Idem
  a(\rhos_1) \otimes_\Idem\dots\otimes_\Idem a(\rhos_i)$ is
  $I(o(o(\x),\vec\rhos_{[1,i]}))$, i.e.,
  \[
  (I_A(\x) \otimes_\Idem a(\rhos_1) \otimes_\Idem\dots\otimes_\Idem
  a(\rhos_i))I(o(o(\x),\vec\rhos_{[1,i]}))=I_A(\x) \otimes_\Idem a(\rhos_1) \otimes_\Idem\dots\otimes_\Idem
  a(\rhos_i),
  \]
  and then $I(o(o(\x),\vec\rhos_{[1,i]})) a(\rhos_{i+1})\neq 0$ if
  and only if $M(\rhos_{i+1}^-)\subset o(o(\x),\vec\rhos_{[1,i]})$ and
  $M(\rhos_{i+1}^+)$ is disjoint from
  $o(o(\x),\vec\rhos_{[1,i]})\setminus M(\rhos_{i+1}^-)$, which is
  exactly the condition of strong boundary monotonicity.
\end{proof}

As explained in Definition~\ref{def:ainf-module},
an $\Ainf$ module
structure on $\CFAa(\HD,\spinc)$ is a family of maps
\[
m_{n+1}\co
\CFAa(\HD,\spinc)\otimes_\Idem\Alg\otimes_\Idem\dots\otimes_\Idem\Alg\to \CFAa(\HD,\spinc)
\]
satisfying Equation~\ref{eq:Ainf-mod-compat}.
In particular, in view of Lemma~\ref{lem:strong-monotone-tensor},
since the
tensor products are over $\Idem(\PMC)$, it suffices to define
\[
m_{n+1}(\x,a(\rhos_1),\dots,a(\rhos_n))
\]
when $(\x,\rhos_1,\dots,\rhos_n)$ is strongly boundary monotone.  (In
the following sections we will use commas to separate the tensor
factors in the argument to $m_{n+1}$ and similar maps.)

\begin{definition}\label{def:Amod-mult}
  Fix an admissible (in the sense of Definition~\ref{def:admissible_J}),
  sufficiently generic (in the sense of Definition~\ref{def:sufficiently-generic})
  almost-complex structure on $\Sigma\times[0,1]\times\RR$.
  For $(\x,\vec\rhos)$ strongly boundary monotone (and
  $\rhos_i\neq\emptyset$ for each $i$), define
\begin{align}
\gls*{monCFA}
(\x,a(\rhos_1),\dots,a(\rhos_n))&\coloneqq
  \sum_{\y\in\S(\HD,\s)}
  \sum_{\substack{B\in\pi_2(\x,\y)\\
       \ind(B,\vec\rhos)= 1}}
    \#\left(\Mod^B(\x,\y;\rhos_1,\dots,\rhos_n)\right)\y\\
  m_2(\x,\Unit) &\coloneqq \x \label{eq:CFA-unital-1}\\
  m_{n+1}(\x,\dots,\Unit,\dots) &\coloneqq 0,\quad\textrm{$n>1$}\label{eq:CFA-unital-2},
\end{align}
\glsadd{CFAspinc}%
where $\gls*{Unit}$
is the unit in $\Alg(\PMC)$. (Equations~(\ref{eq:CFA-unital-1}) and~(\ref{eq:CFA-unital-2}) state that $\CFAa(\HD)$ is strictly unital.)\index{strictly unital!$\CFAa$ is}
\end{definition}

We will sometimes use the alternate notation 
$m_1(\x)=\gls*{bdyCFA}$
and
$m_2(\x,a)=\gls*{multCFA}$.

Let 
$\gls*{CFA}
=\bigoplus_{\spinc\in\SpinC(Y)}\CFAa(\HD,\spinc)$.

\begin{lemma}
  \label{lem:finite-typeA}
  If $\HD$ is
  provincially admissible, then for any fixed generators $\x$ and $\y$ and
  a consistent sequence $\vec\rhos$ of non-empty sets of Reeb
  chords, 
  there are at most finitely many $B\in\pi_2(\x,\y)$ such that
  $\Mod^B(\x,\y;\vec\rhos)$ is non-empty; in particular, the sums appearing in the
  definition of $m_k$ are finite.
  Moreover, if $\HD$ is admissible, then there are only finitely many
  non-zero maps $m_k$ (i.e., in the language of
  Section~\ref{sec:AinfModules}, $(\CFAa(\HD),\{m_k\})$ is bounded).
  \index{bounded!$\Ainf$ module}%
\end{lemma}
\begin{proof}
  The first statement is immediate from
  Proposition~\ref{prop:provincial-admis-finiteness} and the fact that
  if a domain admits a holomorphic representative then all of its
  coefficients must be non-negative (Lemma~\ref{lem:holo-has-pos-domain}).  For the second, observe that if
  we let $|B|$ be the sum of all local multiplicities, then $m_{n+1}$
  counts domains with $|B|\geq n$. Thus, Proposition~\ref{prop:admis-finiteness}
  gives the desired result.
\end{proof}

\begin{remark}Notice that the operation $m_1\co\CFAa\to\CFAa$ is the
  differential obtained by counting only provincial holomorphic curves,
  i.e., curves which do not approach~$\bdy\widebar{\Sigma}$.\index{differential!on $\CFAa$}
\end{remark}

\begin{remark}
  The value of $m_{n+1}$ when one of the arguments is $\Unit$ can be
  made more natural by allowing the corresponding $\rhos_i$ to be
  empty (since $a(\emptyset) = 1$).  The corresponding moduli spaces
  should be interpreted as having a choice of height $t_i$, and no
  other conditions at $t_i$. Except in the case of trivial strips, these
  moduli spaces cannot be $0$-dimensional, and hence they 
  do not contribute to $m_{n+1}$.
\end{remark}

\begin{remark}In Chapter~\ref{chap:nice-diagrams} we will show
  that for special kinds of diagrams (the obvious analogues of ``nice
  diagrams'' from \cite{SarkarWang07:ComputingHFhat}) the $m_i$ vanish
  for $i>2$, so $\CFAa$ is an ordinary (not $\Ainf$) differential
  module. However, there are reasons for considering more general
  diagrams. First, it is not clear how to prove invariance if one
  restricts only to nice diagrams.  Further, in spite of the added
  analytic and algebraic complications, it is often easier to compute
  with a diagram which is not nice, as general diagrams often have a
  much smaller set of generators $\S(\HD)$.
\end{remark}
\index{type $A$ invariant|)}
\section{Compatibility with algebra}
In this section we will prove that $\CFAa(\HD,\s)$ is an
$\Ainf$ module over $\Alg(\PMC)$. Recall from
Definition~\ref{def:ainf-module} that the basic relation for an
$\Ainf$ module (over an $\Ainf$ algebra over $\Field$) is
\begin{align*}
    0 &= \sum_{i+j=n+1}\!m_i(m_j(\x, a_1, \dots, a_{j-1})
        , \dots, a_{n-1})\\
    &\quad+\!\sum_{i+j=n+1}\sum_{\ell=1}^{n-j}m_i(\x, a_1, \dots, a_{\ell-1}, \mu_j(a_\ell, \dots, a_{\ell+j-1}), \dots, a_{n-1}).
\end{align*}
For the algebra $\Alg(\PMC)$, the $\mu_i$ vanish for
$i>2$. Therefore, we must check that
\begin{equation}\label{eq:A-Ainf-rel}
  \begin{split}
0 &=
\sum_{i+j=n+1}m_i(m_j(\x,a_1,\dots,a_{j-1}),\dots,a_{n-1})\\
&\quad+\sum_{\ell=1}^{n-1} m_n(\x,a_1,\dots,\partial a_\ell,\dots,a_{n-1})\\
&\quad+\sum_{\ell=1}^{n-2} m_{n-1}(\x,a_1,\dots,a_\ell a_{\ell+1},\dots,a_{n-1}).
  \end{split}
\end{equation}

\begin{proposition}\label{prop:A-module-defined}The data
$(\CFAa(\HD,\s), \{m_i\}_{i=1}^\infty)$ forms an $\Ainf$ module over
$\Alg(\PMC)$.
\end{proposition}

\begin{figure}
\includegraphics[scale=.83333]{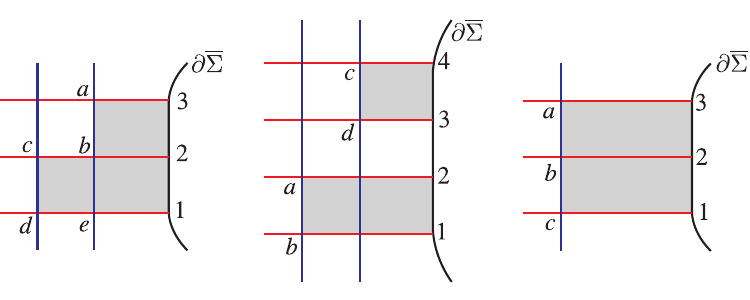}
\caption[Local illustrations of terms occurring the
    $\Ainf$ relation for $\CFAa$]{\textbf{Local illustrations of terms occurring the
    $\Ainf$ relation for $\CFAa$.} The moduli spaces are
  the same as the moduli spaces in Figure~\ref{fig:dsquare0D} in
  Section~\ref{sec:def-CFD}, but the algebraic interpretations of the
  phenomena are different for $\CFAa$ than for
  $\CFDa$.}\label{fig:AASquared}
\end{figure}

We give the proof shortly.  The most interesting points in
the proof of Proposition~\ref{prop:A-module-defined} can be seen in
the following examples. (The reader may wish to compare
Examples~\ref{eg:a1}--\ref{eg:a3} with
Examples~\ref{eg:dd1}--\ref{eg:dd3}.)
\begin{example}\label{eg:a1}On the left of Figure~\ref{fig:AASquared}
  is a piece of a diagram with four generators. The non-trivial
  multiplications in the corresponding module $\CFAa$ are
  \begin{align*}
   \bdy\{b,d\} &= \{c,e\}\\
   \{a,d\}\cdot\honestalg{1}{2}{3} &= \{a,c\}\\
   \{c,e\}\cdot\honestalg{1}{3}{2} &= \{a,c\}\\
   \{b,d\}\cdot\honestalg{2}{3}{1} &= \{a,d\}\\
   \{b,d\}\cdot\honestalg{1&2}{2&3}{} &= \{a,c\} 
  \end{align*}
  There is an $\Ainf$ compatibility equation 
  \[ (\partial \{b,d\})\cdot \honestalg{1}{3}{2} + 
  \{b,d\}\cdot \partial \honestalg{1}{3}{2}=0.\] 
  (Recall that
  $\bdy\honestalg{1}{3}{2}=\honestalg{1 & 2}{2 & 3}{}$.)
  Geometrically, this corresponds to considering the two-dimensional
  moduli space of curves from $\{b,d\}$ to $\{a,c\}$ with a Reeb chord
  from $1$ to $3$. This moduli space has two ends, one of which is
  a height $2$ comb, and the other of which is a join curve end.
\end{example}

\begin{example}\label{eg:a2}
  In the center of Figure~\ref{fig:AASquared} is a piece
  of a diagram with four generators. The non-trivial products are
  \begin{align*}
    \{b,d\}\cdot\honestalg{1}{2}{3}&=\{a,d\}\\
    \{b,c\}\cdot\honestalg{1}{2}{4}&=\{a,c\}\\
    \{a,d\}\cdot\honestalg{3}{4}{2}&=\{a,c\}\\
    \{b,d\}\cdot\honestalg{3}{4}{1}&=\{b,c\}\\
    \{b,d\}\cdot\honestalg{1 & 3}{2 & 4}{} &=\{a,c\}.
  \end{align*}
  Here, associativity follows from the fact that
  \[
  \bigl(\{b,d\} \cdot \honestalg{1}{2}{3}\bigr)\cdot\honestalg{3}{4}{2}=
  \{b,d\} \cdot \honestalg{1 & 3}{2 & 4}{}=
  \bigl(\{b,d\} \cdot \honestalg{3}{4}{1}\bigr)\cdot\honestalg{1}{2}{4}.
  \]
  The second equality above can be seen geometrically as follows.
  Consider  the two-dimensional moduli
  space connecting $\{b,d\}$ to $\{a,c\}$ where the Reeb chord from $3$ to $4$
  is encountered before the one from $1$ to $2$. Then this moduli space
  has two ends, one of which is a height $2$ comb, 
  while at the other end the heights of two Reeb chords are the same.
\end{example}

\begin{example}\label{eg:a3}
  On the right in Figure~\ref{fig:AASquared} is a piece
  of a diagram with three generators. The non-trivial products are
  \begin{align*}
    \{c\}\cdot\honestalg{1}{2}{} &= \{b\}\\
    \{c\}\cdot\honestalg{1}{3}{} &= \{a\}\\
    \{b\}\cdot\honestalg{2}{3}{} &= \{a\}.
  \end{align*}
  Associativity follows from the fact that
  \[
  \bigl(\{c\} \cdot \honestalg{1}{2}{}\bigr)\cdot\honestalg{2}{3}{}=
  \{c\} \cdot \honestalg{1}{3}{}.
  \]
  Geometrically, at one end of the moduli space is a height $2$ comb. At the other end, two levels collide and a split
  curve degenerates.
\end{example}

\begin{figure}
\includegraphics[scale=.83333]{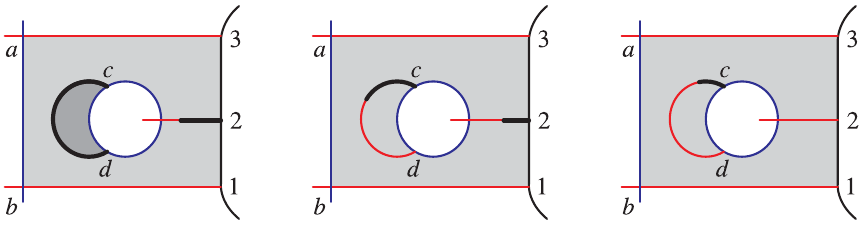}
\caption[Domain where $m_3$ is necessary for associativity]{\textbf{A domain where $m_3$ is necessary for associativity.}
  The domain is shaded; the dark arcs denote cuts, whose ends are
  boundary branched points of $\pi_\Sigma\circ u$. A typical curve in the one-parameter
  family is indicated in the center. On the left is a
  two-story end
  of the moduli space. On the right is
  the other end, corresponding to
  degenerating a split curve.}\label{fig:AHigherProduct}
\end{figure}
\begin{example}\label{eg:a4}
  Figure~\ref{fig:AHigherProduct} shows an example in which the higher
  product $m_3$ is needed for $\Ainf$ associativity. The diagram has
  four generators, $\{a,c\}$, $\{a,d\}$, $\{b,c\}$ and $\{b,d\}$. The
  differential is given by $\bdy\{a,c\}=\{a,d\}$ and
  $\bdy\{b,c\}=\{b,d\}$. For an appropriate choice of complex
  structure $J$, there is a non-trivial product ($m_2$) given by
  $\{b,c\}\cdot\honestalg{1}{3}{}=\{a,c\}$; this corresponds to the
  domain and branch cuts indicated in the right of the figure. One can
  factor
  $\honestalg{1}{3}{}=\honestalg{1}{2}{}\cdot\honestalg{2}{3}{}$. Obviously
  $\{b,c\}\cdot\honestalg{1}{2}{}=0$, so ordinary associativity
  fails. There is, however, for the same choice of complex structure
  $J$ as above, a higher product
  \[
  m_3\left(\{b,d\},\honestalg{1}{2}{},\honestalg{2}{3}{}\right)=\{a,c\},
  \]
  with domain given by the lighter-colored region in the middle of
  Figure~\ref{fig:AHigherProduct}.
  Thus, the $\Ainf$ associativity condition is
  satisfied:
  \[
  m_3\left((\bdy \{b,c\}),\honestalg{1}{2}{},\honestalg{2}{3}{}\right)
  +\left(\{b,c\}\cdot\honestalg{1}{2}{}\right)\cdot\honestalg{2}{3}{}
  +\{b,c\}\cdot\honestalg{1}{3}{}=0.
  \]
  Geometrically, this corresponds to the one-parameter moduli space of
  curves asymptotic to two Reeb chords (at different heights), shown
  at the left of Figure~\ref{fig:AHigherProduct}. One end (pictured in
  the center) is a height $2$ comb where one of the stories
  contains both Reeb chords. At the other end (pictured at the right),
  the two Reeb chords collide and a split curve degenerates.
\end{example}

\begin{proof}[Proof of Proposition~\ref{prop:A-module-defined}]
  Morally, the proof proceeds by considering the ends of the index $2$
  moduli spaces. As in the type $D$ module case, however, we will instead
  appeal to Theorem~\ref{thm:master_equation}.

  Fix generators $\x$ and $\y$, a domain $B\in\pi_2(\x,\y)$
  and a sequence of non-empty sets of Reeb chords
  $\vec{\rhos}=(\rhos_1,\dots,\rhos_n)$ so that $(\x,\vec\rhos)$ is
  strongly boundary monotone, $B$ is compatible with~$\vec\rhos$, and
  $\ind(B,\vec\rhos) = 2$. Set
  $\rhos_i=\{\rho_{i,j}\}$. 
  Theorem~\ref{thm:master_equation} applied to the data
  $(\x,\y,B,\vec{\rhos})$ (for all sources $\Source$ with
  $\chi(\Source)=\chi_{\emb}(B,\vec{\rhos})$) shows that the sum of the
  following numbers is zero:
  \begin{enumerate}
  \item\label{item:AMod-two-lev}The number of two-story ends, i.e., the
    number of elements of
    \[
    \cM^{B_1}(\x,\w;\vec\rhos_1) \times \cM^{B_2}(\w,\y;\vec\rhos_2)
    \]
    where $B = B_1 * B_2$ and $\vec\rhos = (\vec\rhos_1,\vec\rhos_2)$.
  \item\label{item:AMod-join}The number of join curve ends,
    i.e., the number of elements of
    \[
    \cM^B(\x,\y;(\rhos_1,\dots,\rhos_i^{a,b},\dots,\rhos_n))
    \]
    where $\rho_{i,j} = \rho_a \uplus \rho_b$ and
    $\rhos_i^{a,b}$ is obtained from $\rhos_i$ by replacing $\rho_{i,j}$
    with $\{\rho_a,\rho_b\}$.
  \item\label{item:AMod-shuffle}The number of odd shuffle
    curve ends, i.e., the number of elements of
    \[
    \cM^B(\x,\y;\vec\rhos')
    \]
    where $\vec\rhos'$ is obtained from $\vec\rhos$ by performing a
    weak shuffle on one $\rhos_i$.
  \item\label{item:AMod-collide} The number of collisions of levels,
    i.e., the number of elements of
      \[
        \cM^{B}(\x,\y;\rhos_1,\dots,\rhos_i\uplus\rhos_{i+1},\dots,\rhos_n)
      \]
      where $\rhos_i$ and $\rhos_{i+1}$ are weakly composable.
  \end{enumerate}
  (We are also using Lemma~\ref{lem:splittings-embedded} to ensure
  that the limiting curves are embedded.)

  The two-story ends correspond to the first term in
  Formula~\eqref{eq:A-Ainf-rel}.

  Next we discuss the second term in Formula~(\ref{eq:A-Ainf-rel}).
  Lemma~\ref{lem:reeb-differential} implies that, for a consistent
  set $\rhos$ of Reeb chords,
  \[
  \partial a(\rhos) =
  \sum_{\substack{\rhos'\textrm{ a splitting}\\
                   \textrm{of $\rhos$}}} a(\rhos')
  +\sum_{\substack{\rhos'\textrm{ a shuffle}\\
                   \textrm{of $\rhos$}}} a(\rhos').
  \]
  On the other hand, sums~(\ref{item:AMod-join})
  and~(\ref{item:AMod-shuffle}) above correspond, respectively, to the
  sum over all weak splittings and weak shuffles one of
  the sets of Reeb chords in~$\vec\rhos$.  By
  Lemma~\ref{lemma:collision-is-composable}, only honest (not just weak)
  splittings or shuffles contribute, and therefore the second term in
  Formula~\eqref{eq:A-Ainf-rel} corresponds to
  sums~(\ref{item:AMod-join}) and~(\ref{item:AMod-shuffle}) above.

  Finally, sum~(\ref{item:AMod-collide}) corresponds to the third term in
  Formula~\eqref{eq:A-Ainf-rel}, as follows.  Since, by
  Lemma~\ref{lem:reeb-product},
  $a(\rhos_i)a(\rhos_{i+1})$ is $a(\rhos_i\uplus\rhos_{i+1})$ if $\rhos_i$
  and $\rhos_{i+1}$ are composable, and is~$0$ otherwise, the
  coefficient of $\y$ in
  $m_{n}\left(\x,a(\rhos_1),\dots,a(\rhos_i)a(\rhos_{i+1}),\dots,a(\rhos_n)\right)$
  is equal to sum~(\ref{item:AMod-collide}) if $\rhos_i$ and $\rhos_{i+1}$ are
  composable, and is~$0$ otherwise.  But by
  Lemma~\ref{lemma:collision-is-composable},
  sum~(\ref{item:AMod-collide}) is also~$0$ if $\rhos_i$ and $\rhos_{i+1}$
  are not composable, and otherwise gives the third term of
  Formula~\eqref{eq:A-Ainf-rel}.
\end{proof}

\section{Invariance}
\label{sec:A-invariance}\index{invariance!of $\CFAa$|(}%
\begin{theorem}\label{thm:A-invariance}Up to $\Ainf$ homotopy equivalence, the
  $\Ainf$ module $\CFAa(\HD,\s)$ is independent of the almost complex
  structure used to define it, and depends on the provincially
  admissible Heegaard diagram $\HD$ with $\PMC=\partial\HD$ only
  through its induced equivalence class of bordered
  three-manifold~$(Y,\phi\co F(\PMC)\to\partial Y)$. That is, if $\HD$
  and $\HD'$ are provincially admissible bordered Heegaard diagram
  for~$(Y,\phi\co F(\PMC)\to\bdy Y)$ and $\s$ is a $\SpinC$ structure
  on $Y$ then $\CFAa(\HD,\s)$ and $\CFAa(\HD',\s)$ are $\Ainf$
  homotopy equivalent $\Ainf$ $\Alg(\PtdMatchCirc)$-modules.
\end{theorem}

As for invariance of $\CFDa$, the proof of
Theorem~\ref{thm:A-invariance} differs from the standard arguments in
essentially two ways: the additional algebraic complications in
dealing with $\Ainf$ maps and homotopies, and the additional analytic
complications in the case of handleslide invariance because of the
presence of east $\infty$. Both of these issues are slightly more
involved for $\CFAa$ than for $\CFDa$. We will take a similar approach
as for $\CFDa$: first we will explain succinctly the proof of
invariance under change in almost complex structure, and then we will
discuss the (mild) new complications in the triangle maps used to
prove handleslide invariance.

\subsection{Change of complex structure}\label{sec:A-invariance-complex}
In this section, we address invariance of $\CFAa$ under change in
almost complex structure. 
In the course of this analysis, we will relate the
algebra of $\Ainf$ maps and
their $\Ainf$ homotopies 
to holomorphic curve counts; because this is somewhat more involved
than the algebra of differential graded algebras required in
Section~\ref{sec:CFD-cx-str-change}, we will explain more steps in the
proof.

Fix a Heegaard diagram $\HD=(\Sigma,\alphas,\betas,z)$ and sufficiently
generic,
admissible almost complex structures $J_0$ and $J_1$. The complex
structures $J_0$ and $J_1$ lead to potentially different moduli
spaces, and hence to two potentially different $\Ainf$ modules
$\CFAa(\HD,\s;J_0)$ and $\CFAa(\HD,\s;J_1)$. We will show that these two
$\Ainf$ modules are $\Ainf$ homotopy equivalent.

As in Section~\ref{sec:CFD-cx-str-change}, call a (smooth) path
$\{J_r\mid r\in[0,1]\}$ of almost complex structures from $J_0$ to
$J_1$ \emph{admissible} if each $J_r$ is admissible. To the path $J_r$
we can associate a single almost complex structure $J$ on
$\Sigma\times[0,1]\times\RR$ by Formula~\eqref{eq:non-cylind-J},
and define embedded moduli spaces $\cM(\x,\y;\vec\rhos;J)$ by
Formula~\eqref{eq:non-cyclind-embedded-M}.

Recall from Definition~\ref{def:AinfMorphism} that an $\Ainf$ morphism
between $\Ainf$ modules $\cModule_1$ and~$\cModule_2$ is a compatible family of maps
\[
f_n\co
M_1\otimes\overbrace{\Alg\otimes\dots\otimes\Alg}^{n-1}\to M_2,
\]
for $n\in\NN$.  Define an $\Ainf$ map $f^{J_r}=\{f^{J_r}_n\}_{n=1}^{\infty}$ from
$\CFAa(\HD,\s;J_0)$ to $\CFAa(\HD,\s;J_1)$ by setting
\begin{align*}
\gls*{CFAJContMap}
(\x,a(\rhos_1),\dots,a(\rhos_{n-1}))&\coloneqq
  \sum_{\y\in\S(\HD,\s)}
  \sum_{\substack{B\in\pi_2(\x,\y)\\\ind(B,\vec\rhos) = 0}}
    \#\left(\Mod^B(\x,\y;\rhos_1,\dots,\rhos_{n-1};J)\right)\y \\
f^{J_r}_n(\x,\dots,\Unit,\dots) &\coloneqq 0
\end{align*}
and extending $f^{J_r}$ to all of $\CFAa(\HD,\s;J_0)\otimes\Alg^{\otimes(n-1)}$ multilinearly over
$\Idem(\PMC)$. The above map is defined over all of $\CFAa(\HD;J_0)$; but,
in view of Lemma~\ref{lem:SpinCStructures}, it clearly
maps $\CFAa(\HD,\s;J_0)$ to $\CFAa(\HD,\s;J_0)$.

Our first task is to show that $f^{J_r}$ is an
$\Ainf$ morphism. We again use
Proposition~\ref{prop:non-cylindrical-master}.

\begin{lemma}\label{lem:complex-contiuation-A}The map
  $f^{J_r}=\{f^{J_r}_n\}_{n\in\NN}\co\CFAa(\HD,\s;J_0)\to\CFAa(\HD,\s;J_1)$ is
  a strictly unital morphism of $\Ainf$ modules, as defined in
  Definition~\ref{def:AinfMorphism}.
\end{lemma}
\begin{proof}
  The map
  $f^{J_r}$ is strictly unital by definition. It remains
  to check the $\Ainf$ relation for morphisms, which we recall from
  Definition~\ref{def:AinfMorphism} (specialized for $\mu_i = 0$ for
  $i > 2$):
  \begin{align*}
    0 &=\sum_{i+j=n+1}\!m_i^{J_1}(f_j^{J_r}(\x, a_1, \dots, a_{j-1}), \dots, a_{n-1})\\
    &\quad+\!\sum_{i+j=n+1}\!f_i^{J_r}(m_j^{J_0}(\x, a_1, \dots, a_{j-1}), \dots, a_{n-1}) \\
    &\quad+\sum_{\ell=1}^{n-1}f^{J_r}_n(\x, a_1, \dots,
    a_{\ell-1}, \partial a_\ell, a_{\ell+1}, \dots, a_{n-1})\\
    &\quad+\sum_{\ell=1}^{n-2}f^{J_r}_{n-1}(\x, a_1, \dots,
    a_{\ell-1}, a_\ell a_{\ell+1}, a_{\ell+2}, \dots, a_{n-1}).
  \end{align*}
This follows in a similar manner to
Proposition~\ref{prop:A-module-defined}. The first term corresponds to
the two-story ends with a $J$-holomorphic curve followed by a
$J_1$-holomorphic curve. The second term corresponds to two-story
ends with a $J_0$-holomorphic curve followed by a
$J$-holomorphic curve. The third term corresponds to
join curve ends and odd shuffle curve ends. The fourth term corresponds to
collisions of two levels.
\end{proof}

Similarly, fix a generic path $J'_r$ of admissible almost complex
structures with $J'_0=J_1$ and $J'_1=J_0$. Let $J'$ denote the
associated complex structure on $\Sigma\times[0,1]\times\RR$. Then
$J'_r$ defines an $\Ainf$ morphism $f^{J'_r}\co
\CFAa(\HD,\s;J_1)\to\CFAa(\HD,\s;J_0)$. We will show:

\begin{proposition}\label{prop:change-cx-htpy-inv-A}
  The maps $f^{J_r}$ and $f^{J'_r}$ are homotopy inverses to each
  other. In particular, $\CFAa(\HD,\s;J_0)$ is $\Ainf$ homotopy
  equivalent to $\CFAa(\HD,\s;J_1)$.
\end{proposition}
The proof is again essentially standard, but as the algebra may be
unfamiliar we will give a few details.

For $R\gg 0$, let $J\glue_R J'$ denote the almost complex structure on
$\Sigma\times[0,1]\times\RR$ given by
\[
J\glue_R J'|_{(x,s,t)}\coloneqq\begin{cases}
  J_0|_{(x,s,t)} & \text{if } R+1\leq t\\
  J'_{t-R}|_{(x,s,t-R)} & \text{if }R\leq
  t\leq R+1\\
  J_1|_{(x,s,t)} & \text{if }1\leq t\leq R\\
  J_t|_{(x,s,t)} & \text{if }0\leq t\leq
  1\\
  J_0|_{(x,s,t)} & \text{if } t\leq 0.
  \end{cases}
\]
Associated to $J\glue_R J'$ is an $\Ainf$ morphism $f^{J\glue_R J'}\co
\CFAa(\HD,\s;J_0)\to \CFAa(\HD,\s;J_0)$.
\begin{lemma}\label{lemma:AMapsCompose}For $R$ sufficiently large the map $f^{J\glue_R J'}$
  is the $\Ainf$ composite $f^{J'_r}\circ f^{J_r}$.
\end{lemma}
\begin{proof}
Recall from Section~\ref{sec:AinfModules} that $f^{J'_r}\circ
f^{J_r}$ has
\[
(f^{J'_r}\circ
f^{J_r})_n(\x,a_1,\dots,a_{n-1})=\sum_{i+j=n+1}f^{J'_r}_j\left(f^{J_r}_i(\x,a_1,\dots,a_{i-1}),\dots, a_{n-1}\right).
\]
The coefficient of $\y$ in
\[
f^{J'_r}_j\left(f^{J_r}_i(\x,a(\rhos_1),\dots,a(\rhos_{i-1})),\dots, a(\rhos_{n-1}))\right)
\]
counts holomorphic curves in
$
\Mod^{B_1}(\x,\w;\vec\rhos_{[1,{i-1}]};J)\times
 \Mod^{B_2}(\x,\w;\vec\rhos_{[i,{n-1}]};J')
$
where 
$\w \in \S(\HD,\s)$, $B_1\in\pi_2(\x,\w)$, $B_2\in\pi_2(\w,\x)$,
$\ind(B_1,\vec\rhos_{[1,i-1]}) = 0$ and
$\ind(B_2,\allowbreak\vec\rhos_{[i,n-1]}) = 0$.
It follows from obvious non-$\RR$-invariant analogues of
Propositions~\ref{prop:compactness} and~\ref{prop:gluing_two_story}
(compactness and gluing) that for $R$ large enough such curves
correspond to curves in
\[
\bigcup_{B\in\pi_2(\x,\y)}\cM^{B}(\x,\y;\rhos_1,\dots,\rhos_{n-1};J\glue_R J').
\]
The result follows.
\end{proof}

We next define an $\Ainf$ homotopy $h$ between $f^{J\glue_R J'}$
and the identity map. Fix a path of almost complex structures $\{J^w\}_{w\in [0,1]}$ on
$\Sigma\times[0,1]\times\RR$ with $J^0=J\glue_R J'$ and $J^1=
J_0$, and such that each $J^w$ satisfies
conditions (\ref{item:J1}), (\ref{item:J2}) and~(\ref{item:J4}) of
Definition~\ref{def:admissible_J}---i.e., $J^w$ is admissible except
for not being $\RR$-invariant. For a generic path $\{J^w\}_{w\in
  [0,1]}$ and any
$B,\Source,P$ with $\ind(B,\Source,P) = -1$, there are
finitely many $w\in(0,1)$ for which $\cM^B(\x,\y;\Source;P;J^w)$ is
non-empty. So set
\begin{align*}
\Mod^B(\x,\y;\vec\rhos;\{J^w\})&\coloneqq\bigcup_{w\in[0,1]}
  \cM^B(\x,\y;\vec{\rhos};J^w)\\
h_n(\x,a(\rhos_1),\dots,a(\rhos_{n-1}))&\coloneqq
  \sum_{\y\in\S(\HD,\s)}
  \sum_{\substack{B\in\pi_2(\x,\y)\\
                  \ind(B,\vec\rhos) = -1}}
\#\left(\Mod^B(\x,\y;\vec\rhos;\{J^w\})\right)\y.
\end{align*}

In this context, the analogue of
Theorem~\ref{thm:master_equation} is:
\begin{proposition}\label{prop:master_equation_A3}Suppose that
  $(\x,\vec{\rhos})$ satisfies the strong boundary monotonicity
  condition. Fix $\y$, $B\in\pi_2(\x,\y)$, $\Source$, and $\vec{P}$
  such that both $[\vec{P}] = \vec\rhos$ and
  $\ind(B,\Source,\vec{P})=0$. Then, for a generic path
  $\{J^w\}_{w\in[0,1]}$, in the moduli space
  $\bigcup_{w\in[0,1]}\cM^B(\x,\y;\Source;\vec{P};J^w)$, the
  total number of all
  \begin{enumerate}
  \item\label{item:A3-1} elements of $\cM^B(\x,\y;\Source;\vec{P};J^0)$;
  \item\label{item:A3-2} elements of $\cM^B(\x,\y;\Source;\vec{P};J^1)$;
  \item two-story ends, with either
    \begin{enumerate}
    \item a $J_0$-holomorphic curve followed by a
      $J^w$-holomorphic curve or
    \item a $J^w$-holomorphic curve followed by
      a $J_0$-holomorphic curve;
    \end{enumerate}
  \item join curve ends;
  \item odd shuffle curve ends; and
  \item collisions of two levels $P_i$ and $P_{i+1}$ from $\vec{P}$,
    where $[P_i]$ and $[P_{i+1}]$ are composable,
  \end{enumerate}
  is even.
\end{proposition}

\begin{lemma}\label{lemma:AHomotopy}The map $h$ is an $\Ainf$ homotopy between
  $f^{J\glue_R J'}$ and the identity map on $\CFAa(\HD, J_0)$.
\end{lemma}
\begin{proof}
We must check that 
\begin{equation}\label{eq:A-homotopy1}
f^{J\glue_R J'}_1(\x)-\x=h_1(m^{J_0}_1(\x))-m_1^{J_0}(h_1(\x))
\end{equation}
and that for $n>1$,
\begin{equation}\label{eq:A-homotopy2}
  \begin{split}
  f^{J\glue_R J'}_n(\x, a_1 ,\dots, a_{n-1})&=
    \sum_{i+j=n+1}\!h_i(m^{J_0}_j(\x, a_1, \dots, a_{j-1}), \dots, a_{n-1})
    \\
    &\quad+\!\sum_{i+j=n+1}\!m^{J_0}_i(h_j(\x, a_1,
    \dots, a_{j-1}), \dots, a_{n-1})\\
    &\quad+\!\sum_{\ell=1}^{n-1}h_n(\x, a_1,\dots,
    \bdy a_\ell,\dots, a_{n-1})\\
    &\quad+\!\sum_{\ell=1}^{n-2}h_{n-1}(\x, a_1,\dots,
    a_\ell a_{\ell+1},\dots, a_{n-1}).
  \end{split}
\end{equation}
As expected, this follows from
Proposition~\ref{prop:master_equation_A3}. The $f^{J\glue_R J'}$
terms in Formulas~(\ref{eq:A-homotopy1}) and~(\ref{eq:A-homotopy2})
correspond to item~(\ref{item:A3-1}) of
Proposition~\ref{prop:master_equation_A3}. The $-\x$ term of
Formula~\eqref{eq:A-homotopy1} corresponds
to item~(\ref{item:A3-2}), since $J^1 = J_0$ and for an
$\RR$-invariant complex
structure the only index $0$ holomorphic curves are trivial strips. As
usual, two-story ends correspond to the remaining two terms in
Formula~(\ref{eq:A-homotopy1}) as well as the first two sums in
Formula~(\ref{eq:A-homotopy2}).  The join curve and
shuffle curve ends account for the third sum in
Formula~(\ref{eq:A-homotopy2}). Collisions of two levels account for
the final sum in Formula~(\ref{eq:A-homotopy2}).
\end{proof}

\begin{proof}[Proof of Proposition~\ref{prop:change-cx-htpy-inv-A}]
  This is immediate from Lemmas~\ref{lemma:AMapsCompose}
  and~\ref{lemma:AHomotopy}, together with 
  the obvious analogues for the composition
  $f^{J_r}\circ f^{J'_r}$.
\end{proof}

\subsection{Handlesliding an \textalt{$\alpha$}{alpha}-arc over an \textalt{$\alpha$}{alpha}-circle}
\label{sec:A-Handleslides}
In the present subsection we prove handleslide
invariance:
\begin{proposition}\label{prop:handleslide-CFA}
  Let $(\Sigma,\alphas,\betas,z)$ be a bordered Heegaard diagram, and
  let
  $(\Sigma,\alphas^H,\betas^H,z)$ be the result of performing a handleslide
  of one of the $\alpha$-curves over an $\alpha$-circle or one
  $\beta$-circle over another $\beta$-circle. Then there is
  a strictly unital $\Ainf$ homotopy equivalence between
  $\CFAa(\Sigma,\alphas,\betas,z)$ and $\CFAa(\Sigma,\alphas^H,\betas,z)$.
\end{proposition}

The most subtle case of handleslide invariance is handlesliding an
$\alpha$-arc over an $\alpha$-circle; we focus on this case. We will
be brief, as most of the argument is the same as the type $D$ case, in 
Section~\ref{sec:CFD-handleslide}.

Fix a Heegaard triple diagram $(\Sigma,\alphas,\alphas^H,\betas)$ as
in Section~\ref{sec:CFD-handleslide}. There are associated moduli
spaces of holomorphic curves in $\Sigma\times\Delta$
$$\Mod^B(\x,\y,\Theta;\rhos_1,\dots,\rhos_{n})
\coloneqq\bigcup_{i=0}^n
\Mod^B(\x,\y,\Theta_{o_i};(\rhos_1,\dots,\rhos_i),(\rhos_{i+1}^H,\dots,
\rhos_{n}^H))\y.$$
Define the $\Ainf$ morphism $f^{\alphas,\alphas^H,\betas}=\{f^{\alphas,\alphas^H,\betas}_n\}_{n=1}^{\infty}$
associated to the handleslide by
\begin{equation}
  \label{eq:TriangleMapA}
f_n^{\alphas,\alphas^H,\betas}(\x,a(\rhos_1),\dots,a(\rhos_{n-1})
)\coloneqq
\sum_{\y}\sum_{\substack{B\in\pi_2(\x',\y)\\\ind(B,\vec\rhos)=0}}\,
\#\left(\Mod^B(\x,\y,\Theta;\rhos_1,\dots,\rhos_{n-1}\right))\y
\end{equation}
and extending to all of $\CFAa(\HD,\s,J_0)$ by multi-linearity over $\Idem(\PMC)$.

\begin{proof}[Proof of Proposition~\ref{prop:handleslide-CFA}]
  This follows along the same lines as the proof of
  Proposition~\ref{prop:HandleslideInvariance}.  That
  $f^{\alphas,\alphas^H,\betas}$ is a strictly unital $\Ainf$ morphism
  follows from Propositions~\ref{prop:tri-master}
  and~\ref{prop:v12-ends-match}, similarly to how
  Proposition~\ref{prop:A-module-defined} follows from
  Theorem~\ref{thm:master_equation}. To prove that
  $f^{\alphas,\alphas^H,\betas}$ is invertible we consider its
  composition with another triangle map,
  $f^{\alphas^H,\alphas',\betas}$, where $\alphas'$ is a small
  Hamiltonian perturbation of $\alphas$, as in
  Section~\ref{sec:handleslide-is-chain-map}.  Associativity for
  holomorphic triangle maps, the type $A$ analogue of
  Proposition~\ref{prop:triangle-assoc}, shows that
  $f^{\alphas^H,\alphas',\betas}\circ f^{\alphas,\alphas^H,\betas}$ is
  homotopic to $f^{\alphas,\alphas',\betas}$. It follows from
  Lemma~\ref{lemma:filtrations-good} that
  $f^{\alphas,\alphas',\betas}$ is an isomorphism (noting that an
  $\Ainf$ homomorphism $\{f_i\}$ is an isomorphism if and only if
  $f_1$ is an isomorphism). Thus,
  $f^{\alphas,\alphas^H,\betas}$ is injective on homology (or, better, monic in the homotopy category).
  That $f^{\alphas,\alphas^H,\betas}$ is surjective on homology (or, better,
  epic in the homotopy category) follows from a
  similar argument using a fourth set of circles $\alphas^{H,\prime}$,
  as in the proof of Proposition~\ref{prop:triangle-assoc}.
\end{proof}

\subsection{Completion of the invariance proof}

\begin{proof}[Proof of Theorem~\ref{thm:A-invariance}]
  Independence from the complex structure was verified in Section~\ref{sec:A-invariance-complex}.
  Handleslide invariance was outlined in Section~\ref{sec:A-Handleslides}.
  Stabilization invariance is obvious, by stabilizing in the region
  containing the basepoint.
  Isotopy invariance is left to the reader. In light of
  Proposition~\ref{prop:heegaard-moves} this verifies topological
  invariance.
\end{proof}

In view of Theorem~\ref{thm:A-invariance}, if $(Y,\phi\co F(\PMC)\to \partial Y)$ is a bordered
three-manifold and $\s$ is a $\SpinC$ structure over $Y$, define
$\CFAa(Y,\phi\co F(\PMC)\to\partial Y)$  to be $\CFAa(\HD,\s)$, where $\HD$ is any
bordered diagram representing the bordered three-manifold.
Let
$$\CFAa(Y)\coloneqq \bigoplus_{\s\in\SpinC(Y)}\CFAa(Y,\s).$$

\begin{proof}[Proof of Theorem~\ref{intro:A-invariance}]
  This follows at once from Theorem~\ref{thm:A-invariance}.
\end{proof}
\index{invariance!of $\CFAa$|)}%
\section{Twisted coefficients}
We can also define a twisted theory $\tCFAa(\HD, \x_0)$, for
$\x_0$ a generator representing the $\s\in\SpinC(Y)$, in an analogous
way to the definition of $\tCFDa(\HD,\x_0)$ in
Section~\ref{sec:typeD-twisted}.

Consider the generating set 
$\gls*{generatorsUs}$,
consisting of elements of the
form $e^{B_0}\x$, where $\x\in\S(\HD,\s)$ and $B_0 \in \pi_2(\x_0, \x)$.
$\gls*{tCFA}$
is the $\Field$ vector space spanned by
$\underline{\S}(\HD,\s)$. It is given an action of $\Idem(\PtdMatchCirc)$
generalizing Equation~\eqref{eq:typeAactionI}, via
\begin{equation*}
  I(\SetS) \cdot e^{B_0}\x \coloneqq
  \begin{cases}
    \x&I(\SetS) = I_A(\x)\\
    0&\text{otherwise}
  \end{cases}
\end{equation*}

There is an action of $H_2(Y,\bdy Y)$ on
$\underline{\CFAa}(\HD, \s)$, by composition with the corresponding
periodic domain: For $\alpha \in H_2(Y, \bdy Y)$ corresponding to $B
\in\pi_2(\x_0,\x_0)$, define
\[
\alpha * e^{B_0}\x \coloneqq e^{B*B_0}\x.
\]
We can view this as a left action by the group ring $\Field[H^2(Y,\partial Y)]$.

The $\Ainf$ action by $\Alg(\PMC)$ is specified by:
\begin{align*}
m_{n+1}(e^{B_0}\x,a(\rhos_1),\dots,a(\rhos_n))&\coloneqq
  \sum_{\y\in\S(\HD)}
  \sum_{\substack{B\in\pi_2(\x,\y)\\
       \ind(B,\vec\rhos)= 1}}
    \!\!\!\!\#\left(\Mod^B(\x,\y;\rhos_1,\dots,\rhos_n)\right) e^{B_0*B} \y\\
  m_2(e^{B_0}\x,\Unit) &\coloneqq \x \\
  m_{n+1}(e^{B_0}\x,\dots,\Unit,\dots) &\coloneqq 0,\quad\textrm{$n>1$},
\end{align*}
where here $\Unit$ is the unit in $\Alg(\PMC)$. Extend this to an
$\Ainf$ bimodule action of $\Field[H^2(Y,\partial Y)]$ and $\Alg(\PMC)$ by setting
\begin{align*}
  m_{1,1,0}(\alpha,e^{B_0}\x)&\coloneqq\alpha*e^{B_0}\x\\
  m_{m,1,n}(\alpha_1,\dots,\alpha_m,e^{B_0}\x,a_1,\dots,a_n)&\coloneqq 0
  &&\text{if }nm>0\text{ or }m>1.
\end{align*}

\begin{proposition}\label{prop:tA-module-defined}The data
$(\tCFAa(\HD,\s), \{m_i\}_{i=1}^\infty)$ forms an $\Ainf$ bimodule
over the rings
$\Field[H_2(Y,\partial Y)]$ and $\Alg(\PMC)$.
\end{proposition}

\begin{proof}
  This is a simple adaptation of the proof of
  Proposition~\ref{prop:A-module-defined}.
\end{proof}

\begin{theorem}\label{thm:tA-invariance}Up to $\Ainf$ homotopy
  equivalence, the $\Ainf$ bimodule
  $\tCFAa(\HD,\spinc)$ over $\Field[H_2(Y,\penalty 600\partial Y)]$ and
  $\Alg(\PMC)$ is independent of the almost complex structure
  used to define it, and depends on the Heegaard diagram $\HD$
  only through its (equivalence class of) 
  underlying  bordered three-manifold~$(Y,\phi\co {F(\PMC)}\to \partial Y)$ 
  and on the generator $\x_0$ only through
  its induced $\SpinC$ structure $\spinc_z(\x_0)$.
\end{theorem}

\begin{proof}
  The proof of Theorem~\ref{thm:A-invariance} can be adapted easily to the twisted context,
  following the model of Theorem~\ref{thm:tD-invariance}. Note once again that all the maps
  commute with the $\Field[H_2(Y,\partial)]$-action (rather than merely commuting up to homotopy).
\end{proof}


\chapter{Pairing theorem via nice diagrams}
\label{chap:nice-diagrams}

We now give a proof of the pairing
theorem (Theorem~\ref{thm:TensorPairing}),
adapting methods from Sarkar and
Wang~\cite{SarkarWang07:ComputingHFhat}.

\begin{definition}
  \label{def:NiceDiagram}
  A bordered Heegaard diagram $(\Sigma, \alphas,
  \betas, z)$ is called \emph{nice}\index{nice diagram}\index{Heegaard diagram!nice}
  if every region in
  $\widebar{\Sigma}$ except for the region 
  $\gls*{Dz}$
  adjacent to $z$ is a (topological) disk with at
  most $4$ corners. That
  is, each region in the interior of $\Sigma$ is either a bigon or a
  quadrilateral and
  any region at the boundary of $\Sigma$ except for $D_z$ has two $\alpha$-arcs, one
  $\beta$-arc and one arc of $\partial\widebar{\Sigma}$ in its boundary.
\end{definition}

\begin{proposition}\label{prop:nice-diagram}
  Any bordered Heegaard diagram $(\Sigma, \alphas,
  \betas, z)$ can be turned into a nice bordered Heegaard diagram
  $(\Sigma, \alphas, \betas', z)$ via
  \begin{itemize}
    \item isotopies of the $\beta_i$ not crossing $\partial \Sigma$ and    
    \item handleslides among the $\beta_i$.
  \end{itemize}
\end{proposition}

\begin{figure}
  \includegraphics[scale=.4167]{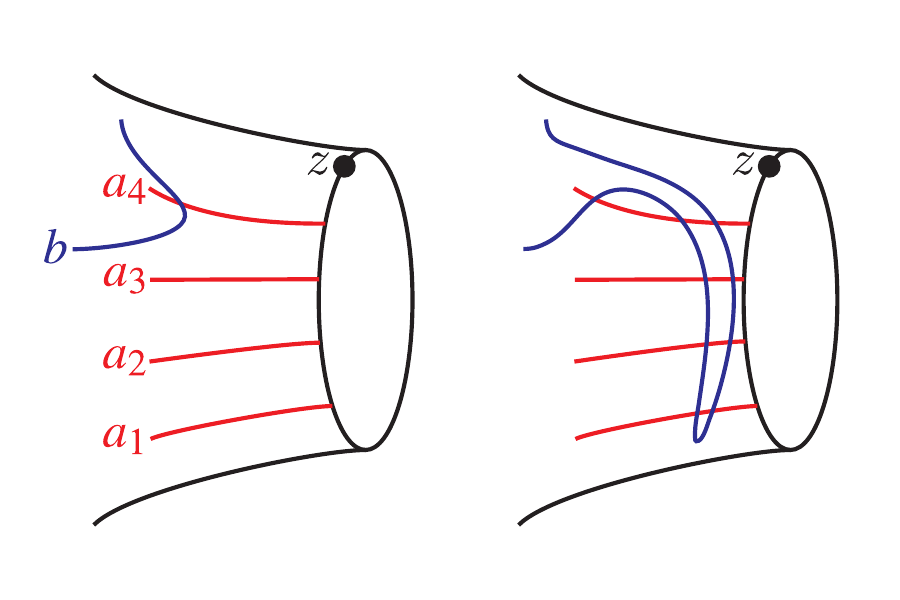}
  \caption[First step in making a bordered Heegaard diagram
  nice]{\textbf{First step in making a bordered Heegaard diagram
      nice:} perform a finger-move to separate off the boundary.}\label{Fig:MakingNice}
\end{figure}

\begin{proof}
  We start by doing a finger move as in
  Figure~\ref{Fig:MakingNice}. That is, let $a_1, \dots, a_{4k}$
  denote the
  ends of the $\alpha$-arcs intersecting~$\partial \Sigma$, enumerated
  (say) counter-clockwise
  around $\partial \Sigma$, with $a_1$ and $a_{4k}$ forming part of $\partial
  D_z$. Let $b$ be the $\beta$-arc on the boundary of $D_z$ which intersects
  $a_1$ closest to $\bdy\bSigma$. Do a finger move pushing $b$ consecutively through $a_1$, \dots,
  $a_{4k}$. In the resulting diagram, the regions $R_1,\dots,R_{4k-1}$ adjacent to
  $\partial\widebar{\Sigma}$ other than $D_z$ are all rectangles.
  
  Now apply the Sarkar-Wang
  algorithm~\cite{SarkarWang07:ComputingHFhat}, which we
  recall briefly, to the resulting diagram. Assign to each region $R$ its
  \emph{distance}, which is the minimum number of $\beta$-arcs which a
  path in $\Sigma\setminus\alphas$ connecting $R$ and $D_z$ must
  cross. A region is called \emph{bad} if it is not a bigon or
  rectangle.  The algorithm is inductive. At each step, choose a bad
  region $R$ of maximal distance and with minimal badness (i.e., maximal
  Euler measure) among such
  regions. Let $b_R$ be a $\beta$-arc separating $R$ from a region of
  smaller distance. Do a finger move, pushing $b_R$ through one of the
  $\alpha$-edges of $R$, and continue the finger move until reaching
  either a lower distance region, a bigon, another bad region, or the
  region $R$ again. In the last case, the algorithm tells us to either
  choose a different edge through which to push $b_R$ or perform an
  appropriate handleslide; otherwise, one repeats the process.  The
  crucial point for us is that the finger move will never reach any of
  the regions $R_1,\dots,R_{2n-1}$. The proof in the closed case thus
  implies that the process eventually yields a nice bordered Heegaard diagram.
\end{proof}

Nice diagrams are automatically admissible:
\begin{lemma}\label{lem:nice-admissible}
  If $\HD=(\Sigma,\alphas,\betas,z)$ is nice then $\HD$ is
  admissible.
\end{lemma}
\begin{proof}
The proof is exactly the same as in the closed case (\cite[Corollary
3.2]{LMW06:CombinatorialCobordismMaps}). One shows that the region
containing $z$ occurs both to the left and to the right of each
$\beta$-circle: if all the regions to the left (say) of some $\beta_i$ are
bigons or rectangles then one sees that the $\beta$-circles are
homologically linearly dependent.
\end{proof}

Sarkar and Wang apply the index formula
to show that for nice Heegaard diagrams (in the closed case),
the differential counts only bigons and rectangles. This fact has
the following analogue for nice bordered Heegaard
diagrams. 

\begin{proposition}\label{Prop:CurvesInNiceDiagrams}Let $\HD$ be a
  nice bordered Heegaard diagram.
  Suppose that $\ind(B,\vec{\rhos})=1$. Then any holomorphic
  curve $u\in\cM^B(\x,\y;\vec{\rhos})$ has one of the following three
  forms:
\begin{enumerate}
\item\label{Item:ManyBigons} The source $S$ of $u$ consists of $g$
  bigons. Each bigon is either a trivial strip or has a unique
  puncture mapped to east $\infty$. Further, $\vec{\rhos}$ has only one part.
\item\label{Item:OneBigon} The source $S$ consists of $g$ bigons. All
  but one of the bigons are trivial strips, and the other is mapped by
  $\pi_\Sigma\circ u$ to the interior of $\Sigma$.
\item\label{Item:OneQuadrilateral} The source $S$ consists of $g-1$
  bigons and one quadrilateral. All of the bigons are trivial strips,
  and the quadrilateral is mapped by $\pi_\Sigma\circ u$ to the
  interior of~$\Sigma$.
\end{enumerate}

Conversely, suppose the homology class $B$ and asymptotic data
$\vec{\rhos}$ admit a topological map satisfying the condition of
Case~(\ref{Item:ManyBigons}), (\ref{Item:OneBigon}), or
(\ref{Item:OneQuadrilateral}), and
$\ind(B,\vec{\rhos})=1$.  Then there is a
holomorphic map in $\cM^B(\x,\y;\vec{\rhos})$,
unique up to translation.
\end{proposition}
\begin{proof}
By Formula~\eqref{eq:Index}, for a rigid holomorphic curve we have
\[
1 = g-\chi(S)+2e(D(u))+|P|.
\]
Since $S$ has at most $g$ connected components, $\chi(S)\leq g$. Since
the Heegaard diagram is nice, $e(D(u))\geq 0$. Consequently, one of the
following three cases holds:
\begin{itemize}
\item $g=\chi(S)$, $e(D(u))=0$ and $|P|=1$. This corresponds to Case~\ref{Item:ManyBigons}, above.
\item $g=\chi(S)$, $e(D(u))=1/2$ and $|P|=0$. This corresponds to Case~\ref{Item:OneBigon}, above.
\item $g-1=\chi(S)$, $e(D(u))=0$ and $|P|=0$. This corresponds to Case~\ref{Item:OneQuadrilateral}, above.
\end{itemize}

It is clear that in each case the holomorphic representative is
unique. The converse, that such domains admit holomorphic
representatives, follows from the Riemann mapping theorem.
\end{proof}

\begin{remark}\label{rmk:sarkar-nice}
  In the closed case, Sarkar and Wang show that for a rigid
  holomorphic curve $u$ in a nice diagram, the image of
  $\pi_\Sigma\circ u$ is an empty, embedded bigon or
  rectangle~\cite[Theorem 3.3]{SarkarWang07:ComputingHFhat}. The
  analogue in the bordered case is that the domain of a rigid curve is
  an empty embedded bigon, rectangle or half strip (rather than, say,
  merely immersed).  This follows
  immediately from the Proposition~\ref{Prop:CurvesInNiceDiagrams} and
  the result in the closed case, by doubling the Heegaard diagram
  along the boundary and considering the doubled domain.

  We do not need this stronger result.
  Proposition~\ref{Prop:CurvesInNiceDiagrams} already gives an
  algorithm for counting holomorphic curves in a nice diagram: one
  checks which domains have index $1$, and are represented by
  orientation-preserving maps of a bigon, rectangle or half-strip.
  (This is analogous to the approach taken in
  \cite{LMW06:CombinatorialCobordismMaps} for holomorphic triangles.)
  Proposition~\ref{Prop:CurvesInNiceDiagrams} also suffices to prove
  the pairing theorem, our main goal.
\end{remark}

\begin{corollary}\label{Cor:NiceIsHonest}Let $\HD$ be a nice Heegaard diagram with
  boundary, and $\CFAa=\CFAa(\HD)$ the associated
  type $A$ module. Then the higher $\Ainf$ operations $m_i$ ($i>2$)
  vanish on $\CFAa$.
\end{corollary}
\begin{proof}
  The higher multiplications $m_i$ ($i>2$) count rigid holomorphic
  curves in moduli spaces $\cM^B(\x,\y;\vec{\rhos})$ with
  $|\vec{\rhos}|>1$. But by
  Proposition~\ref{Prop:CurvesInNiceDiagrams},
  in a nice diagram there are no such curves.
\end{proof}

As a consequence of these remarks, we have a quick proof of
Theorem~\ref{thm:TensorPairing}.

\begin{proof}[Proof of Theorem~\ref{thm:TensorPairing}.]
  Let $\HD_1=(\Sigma_1,\alphas_1,\betas_1,z)$ and
  $\HD_2=(\Sigma_2,\alphas_2,\betas_2,z)$ be nice Heegaard diagrams
  with boundary for $Y_1$ and $Y_2$ respectively.  Their
  union $\HD=({\overline\Sigma}_1\cup_{\partial}{\overline\Sigma}_2,
  \balphas_1\cup_\bdy\balphas_2,\betas_1\cup\betas_2,z)$ is a Heegaard
  diagram for $Y=Y_1\cup_{F} Y_2$. We want to show that
  $\CFa(\HD)$ is chain homotopy equivalent to $\CFAa(\HD_1)\DTP\CFDa(\HD_2)$.

  It is easy to see that $\HD$ is a nice diagram, since $\HD_1$ and
  $\HD_2$ are. Sarkar and Wang~\cite{SarkarWang07:ComputingHFhat} showed
  that the differential for
  $\CFa(\HD)$ counts index $1$ rectangles and bigons. (This is also a
  special case of Proposition~\ref{Prop:CurvesInNiceDiagrams}.)
  We compare this complex with $\CFAa(\HD_1)\DT\CFDa(\HD_2)$ from
  Section~\ref{sec:DT}, which by
  Proposition~\ref{prop:IdentifyDT} serves as a model for the derived
  tensor product~$\DTP$. (Note, in particular, that since $\HD_2$ is admissible,
  $\CFDa(\HD_2)$ is bounded.)

  There is an obvious identification of generators of
  $\CFAa(\HD_1)\DT\CFDa(\HD_2)$ and $\CFa(\HD)$. Let
  $\x_1\otimes\x_2\in
  \CFAa(\HD)\DT\CFDa(\HD_2)$. (Recall that the tensor products
  $\otimes$ are over the ground ring
  $\Ground=\Idem(\PMC)$, the subring of idempotents in $\Alg(\PMC)$.)

  By Proposition~\ref{Prop:CurvesInNiceDiagrams}, the
  differential
  $\partial(\x_2)$ can be written as $\sum_{\y_2}a_{\x_2,\y_2}\otimes
  \y_2$, where each term in
  $a_{\x_2,\y_2}$ is either $1$---coming from a provincial rectangle
  or bigon---or an algebra element in which a single strand moves.
  Applying Corollary~\ref{Cor:NiceIsHonest}, the differential 
  in $\CFAa(\HD_1)\DT\CFDa(\HD_2)$ is given by
  \[
  \partial^{\DT}(\x_1\otimes\x_2)=(m_1 (\x_1))\otimes\x_2
  +\sum_{\y_2}m_2(\x_1,a_{\x_2,\y_2})\otimes\y_2.
  \]
  The first term in the sum counts holomorphic bigons and rectangles
  contained entirely in $(\Sigma_1,\alphas_1,\betas_1,z)$. The second
  sum counts both rectangles and bigons contained entirely in
  $(\Sigma_2,\alphas_2,\betas_2,z)$ (corresponding to terms where
  $a_{\x_2,\y_2}$ is~$1$) and rectangles in
  $(\widebar{\Sigma}_1\cup_\bdy\widebar{\Sigma}_2,\alphas_1\cup\alphas_2,\betas_1\cup\betas_2,z)$
  crossing $\bdy\widebar{\Sigma}_1$ (corresponding to terms with a single
  moving strand in $a_{\x_2,\y_2}$). But, since $\HD$ is a nice diagram,
  this is exactly the differential on
  $\CFa(\HD)$.
\end{proof}

\begin{remark}
The idea to use the Sarkar-Wang algorithm to prove a gluing result is
partly inspired by work of Juh{\'a}sz~\cite{Juhasz08:SuturedDecomp}.
\end{remark}

\begin{remark}
  Since for $\HD_1$ a nice diagram $\CFAa(\HD_1)$ is an ordinary
  $\Alg$-module, we can also consider the ordinary tensor product
  $\CFAa(\HD_1)\otimes\CFDa(\HD_2)$.  It is easy to see that this
  agrees, on the one hand, with $\CFAa(\HD_1)\DT\CFDa(\HD_2)$.  On the
  other hand, it is straightforward to
  check that $\CFDa(\HD_2)$ is a projective $\Alg$-module, and thus
  $\CFAa(\HD_1)\otimes\CFDa(\HD_2)$ is chain homotopy equivalent to
  $\CFAa(\HD_1)\DTP\CFDa(\HD_2)$.  This gives an alternative to the use
  of Proposition~\ref{prop:IdentifyDT}.
\end{remark}

With a little extra work, one can establish a graded version of the
pairing theorem along these lines; see Section~\ref{sec:GradedPairingThm}.


\chapter{Pairing theorem via time dilation}
\label{chap:tensor-prod}

The aim of the present section is to give an alternate proof of
Theorem~\ref{thm:TensorPairing}, which describes $\HFa$ of a
three-manifold which is separated by a surface in terms of the
bordered invariants of the two pieces. Whereas the proof from
Chapter~\ref{chap:nice-diagrams} used special Heegaard diagrams, the
current proof starts from any (suitably admissible) diagrams for the two pieces.  Our
motivation for including this second proof is
twofold. First, the proof illuminates the connection between the
analysis of holomorphic curves and the algebra underlying the
definitions of the bordered Floer modules and their $\Ainf$
tensor products.
Second, the proof applies to other contexts,
such as counting higher polygon maps; see~\cite{LOT:DCov2}.

In slightly more detail,
let $F$ be a parameterized surface, and let $Y_1$ and $Y_2$ be
three-manifolds with compatibly parameterized boundaries, $\partial Y_1=F(\PMC)=-\partial Y_2=-F(\PMC)$. Let $Y = Y_1 \cup_F Y_2$. 
Any two (suitably admissible) diagrams~$\HD_1$ for~$Y_1$ and~$\HD_2$
for~$Y_2$ can be glued together to form an admissible Heegaard diagram~$\HD$
for~$Y$. Standard gluing theory describes the moduli spaces counted in the differential
for $\CFa(\HD)$ in terms of fibered products of moduli spaces from~$\HD_1$ and~$\HD_2$.
Starting from this description, the differential for $\CFa(\HD)$ is then deformed to the differential on the
derived tensor product by a procedure we call \emph{time dilation}.

More precisely, we construct a family of chain complexes $\CFa(T;\HD_1,
\HD_2)$ from fibered products of the moduli spaces for $\HD_1$ and
$\HD_2$, indexed by a real parameter~$T>0$
which, informally, speeds up the time parameter~$t$
on the type~$D$
side (i.e., $\HD_2$), as compared to time on the type~$A$ side ($\HD_1$).
Standard gluing
theory identifies $\CFa(1;\HD_1, \HD_2)$ with
$\CFa(\HD_1 \cup_\bdy \HD_2)$.  A continuation argument shows that the chain
homotopy type of the complex $\CFa(T;\HD_1,\HD_2)$ is independent of~$T$.
A further compactness and gluing argument shows that for all sufficiently large $T$, the
behavior of the complex stabilizes. In the large $T$ limit,
the complex is then identified with the $\DT$~product, in the sense of
Definition~\ref{def:DT}, of the chain complexes $\CFAa(\HD_1)$ and
$\CFDa(\HD_2)$.  Since this in turn is identified with the derived
tensor product (Proposition~\ref{prop:IdentifyDT}), 
Theorem~\ref{thm:TensorPairing} follows.

This chapter is organized as follows. In
Section~\ref{sec:moduli-matched-pairs} we explain how, via a
neck-stretching argument, holomorphic curves for $\HD_1\cup_\bdy\HD_2$
correspond to certain pairs of holomorphic curves for $\HD_1$ and
$\HD_2$ (for suitable almost complex structures). In
Section~\ref{sec:dilating-time} we introduce the chain complexes
$\CFa(T;\HD_1,\HD_2)$ and show that they are all chain homotopy
equivalent. In Section~\ref{sec:dilating-infinity} we study the
behavior of holomorphic curves as the dilation parameter $T$ goes to
$\infty$. This is the most technically involved part of the
chapter. In Section~\ref{sec:pairing-proof} we combine the results
from the previous sections to prove
Theorem~\ref{thm:TensorPairing}. Finally, in
Section~\ref{sec:tensor-prod-eg} we give a local example of the steps
in the argument; the reader may want to read this simple example
before delving into the details of the proof.

\section{Moduli of matched pairs}
\label{sec:moduli-matched-pairs}

\glsadd{HDsubi}\glsadd{pairingSigma}
Let $\HD_1 = (\Sigma_1,\alphas_1,\betas_1,z)$ and
$\HD_2 = (\Sigma_2,\alphas_2,\betas_2,z)$ be provincially admissible
Heegaard diagrams for $Y_1$ and~$Y_2$ of genus $g_1$ and~$g_2$, respectively, inducing
compatible markings of~$F$; i.e., $\bdy \HD_1 = -\bdy \HD_2$. Let
$\Sigma=\overline{\Sigma}_1\cup_\bdy\overline{\Sigma}_2$,
$\alphas=\balphas_1\cup_\bdy \balphas_2$, and $\betas =
\betas_1\cup\betas_2$.  Assume that one of $\HD_1$ or $\HD_2$
is admissible.  Then, by Lemma~\ref{lem:closed-admissible},
$(\Sigma,\alphas,\betas,z)$
is an admissible Heegaard diagram (denoted~$\HD$) for~$Y$. Recall from
Section~\ref{sec:GluingDiagrams} that generators $\Gen(\HD)$ of
$\CFa(\HD)$ correspond to certain elements
$(\x_1,\x_2)$ of $\Gen(\HD_1)\times\Gen(\HD_2)$, called compatible
pairs of generators; these are the pairs
\index{generator!compatible pair of}\index{compatible pair!of generators}%
$(\x_1,\x_2)$ so that $I_A(\x_1)=I_D(\x_2)$.
Moreover, by
Lemma~\ref{lem:HoClassFibProd}, $\pi_2(\x,\y)$
is naturally a subset of $\pi_2(\x_1,\y_1)\times\pi_2(\x_2,\y_2)$.

\begin{definition}
  \label{def:Compatible}\index{compatible pair!of decorated sources}%
  By a \emph{compatible pair} of decorated sources we mean decorated
  sources $\SourceSub{1}$ and $\SourceSub{2}$ together with a
  bijection $\varphi$ between $E(\SourceSub{1})$ and
  $E(\SourceSub{2})$, such that for each puncture
  $q_i$ of $\SourceSub{1}$, the Reeb chord labeling $\varphi(q_i)$ is the
  orientation reverse of the Reeb chord labeling $q_i$.
  A compatible pair of sources $\SourceSub{1}$ and $\SourceSub{2}$
  can be {\em preglued}\index{pregluing sources!giving closed source}
  to form a surface
  $
  \gls*{closedSource}
  \coloneqq\SourceSub{1}\mathbin{\gls*{preglue}}\SourceSub{2}$, 
  by identifying a
  neighborhood of $q_i$ with a neighborhood of $\phi(q_i)$ for each $i$.
\end{definition}

\begin{definition}
  \label{def:weak-matched}
  Fix two compatible pairs of generators $\x_i,\y_i\in\Gen(\HD_i)$ ($i=1,2$),
  and homology classes $B_1\in\pi_2(\x_1,\y_1)$ and
  $B_2\in\pi_2(\x_2,\y_2)$ inducing a homology class
  $B\in\pi_2(\x,\y)$. Fix a compatible pair of sources $\SourceSub{1}$
  and $\SourceSub{2}$ connecting $\x_1$ to $\y_1$ and $\x_2$ to $\y_2$
  respectively.  Then define the \emph{moduli
    space of matched pairs}
  \index{moduli space!of matched pairs}%
  \index{matched pairs|see{moduli space of matched pairs}}%
  to be the fibered product
  $$
  \gls*{MatchMod}
  \coloneqq
  \tcM^{B_1}(\x_1,\y_1;\SourceSub{1})\times_{\ev_1=\ev_2}
  \tcM^{B_2}(\x_2,\y_2;\SourceSub{2}).$$
  That is, the moduli space of matched pairs consists of pairs
  $(u_1,u_2)$ with $u_i\in\tcM^{B_i}(\x_i,\y_i;\SourceSub{i})$ such
  that $\ev(u_1)=\ev(u_2)$, under the correspondence induced
  by~$\varphi$.
\end{definition}

Define the index of a matched pair by
\index{moduli space!of matched pairs!expected dimension}%
\index{index!of matched pair}%
\begin{equation}
  \label{eq:IndexMatchedPair}
  \begin{split}
    \gls*{MatchIndSrc}
     &\coloneqq \ind(B_1,\SourceSub{1},P)+
         \ind(B_2,\SourceSub{2},P)-m\\
      &= g_1+g_2+2e(B_1)+2e(B_2)-\chi(S_1)-\chi(S_2)+m,
  \end{split}
\end{equation}
where $P$ is the discrete partition on
$E(\SourceSub{1})$ (or equivalently $E(\SourceSub{2})$) and $m =
\abs{P} = \abs{E(\SourceSub{1})}$.
This is the expected dimension of the moduli space of matched pairs coming from
its description as a fibered product of two moduli spaces over $\RR^m$.

\begin{lemma}\label{lemma:Matched-Exp-Dim}
  For generic, admissible almost complex structures on
  $\Sigma_i\times[0,1]\times\RR$, the moduli spaces of matched
  pairs
  $\tcMM^{B_1\glue B_2}(\x_1,\y_1,\SourceSub{1}\semico
  \x_2,\y_2,\SourceSub{2})$ are 
  transversely cut out by the $\dbar$-equation and the evaluation map.
\index{transversality!for matched pair}%
  Thus, the moduli space is a manifold whose dimension is given by
  $\ind(B_1,\SourceSub{1}\semico B_2,\SourceSub{2})$.
\end{lemma}

\begin{proof}
  For the transversality statement, by Proposition~\ref{prop:transversality}, for a generic almost complex structure $J_1$ on $\Sigma_1$, the moduli spaces of holomorphic curves 
$\tcM^{B_1}(\x_1,\y_1;\SourceSub{1})$ are transversally cut out for all (countably many) choices of $\x_1$, $\y_1$, $B_1$ and $\SourceSub{1}$. Again by Proposition~\ref{prop:transversality}, for a generic almost complex structure $J_2$ on $\Sigma_2$, the moduli spaces of holomorphic curves $\tcM^{B_2}(\x_2,\y_2;\SourceSub{2})$, for all choices of $B_2$, $\x_2$, $\y_2$ and $\SourceSub{2}$, are transversally cut out and the evaluation maps $\ev\co\tcM^{B_2}(\x_2,\y_2;\SourceSub{2})\to\RR^m$ are transverse to the images $\ev(\tcM^{B_1}(\x_1,\y_1;\SourceSub{1}))$ for all choices of $B_1$, $\x_1$, $\y_1$ and $\SourceSub{1}$.

The statement about the expected dimension is immediate from the definitions and the index formula, Formula~\eqref{eq:Index}.
\end{proof}

\begin{lemma}
  \label{lem:ModuliSpacesAreMonotone}
  Fix a compatible pair of generators $\x_1$ and $\x_2$ with $\x_i\in\Gen(\HD_i)$.
  Let $(u_1,u_2)$ be any matched pair, i.e., $u_i\in\tcM^{B_i}(\x_i,\y_i;\SourceSub{i})$ such
  that $\ev(u_1)=\ev(u_2)$, under the correspondence induced
  by~$\varphi$. Then, both $u_1$ and $u_2$ are strongly boundary monotone.
\end{lemma}
\begin{proof}
  The condition that $\x_1$ and $\x_2$ can be combined to give a generator of the glued
  diagram is equivalent to the condition that $o(\x_1)\amalg o(\x_2)=[2k]$.
  Now, let $\vec{\rhos}=(\{\rho_1\},\dots,\{\rho_n\})$ denote the sequence of (singleton sets of) Reeb chords
  in order on the $u_1$ side. Of course, we encounter the same sequence of chords,
  only with orientations reversed, in the same order on the $u_2$ side.
  Continuity of $u_1$ and $u_2$ ensures that 
  $M(\rhos_{i+1}^-)\subset o(\x_1,\vec{\rhos}_{[1,i]})$ and
  $M(-\rhos_{i+1}^-)\subset o(\x_2,-\vec{\rhos}_{[1,i]})$. So,
  Lemma~\ref{lem:BiMonotonicity} implies the result.
\end{proof}

Translation induces an $\RR$-action on
$\tcMM{}^B(\x_1,\y_1,\SourceSub{1}\semico
\x_2,\y_2,\SourceSub{2})$.\index{$\RR$-action!on matched pairs}
This action is free except
where both sides of the matching are trivial strips.  When the action
is free we say the curves in the moduli space are
\emph{stable},\index{stable!matched pair}
and
denote the 
quotient $\tcMM{}^B(\x_1,\y_1,\SourceSub{1}\semico
\x_2,\y_2,\SourceSub{2})/\RR$
by
$\gls*{MatchModUnparam}$.

Given a compatible pair $\SourceSub{1}$ and $\SourceSub{2}$ with
$\closedSource=\SourceSub{1}\glue \SourceSub{2}$, 
consider the moduli space of holomorphic curves
$\cM^B(\x,\y;\closedSource)$ for the diagram~$\HD$. The expected
dimension of $\cM^B(\x,\y;\closedSource)$ is given by
$$\ind(B,\closedSource)\coloneqq g-\chi(S)+2e(B).$$
(This is~\cite[Formula (6)]{Lipshitz06:CylindricalHF}, which is a
special case of Equation~\eqref{eq:Index}.)  If
$\abs{E(\SourceSub{1})} = m$, we have $\chi(S_1\glue
S_2)=\chi(S_1)+\chi(S_2)-m$. It follows from
Formula~\eqref{eq:IndexMatchedPair} that the index of the moduli
space of curves from $\x$ to $\y$ represented by the pre-glued
source~$\Source$ coincides with the index of the moduli space of
matched
pairs, i.e.,
$$\ind(B_1,\SourceSub{1};B_2,\SourceSub{2})
=\ind(B_1\glue B_2,\SourceSub{1}\glue\SourceSub{2}).$$

A stronger identification occurs at the level of moduli spaces.

\begin{proposition}
  \label{prop:Gluing1}
  There are generic admissible complex structures on
  $\Sigma\times[0,1]\times \RR$ and on $\Sigma_i\times[0,1]\times\RR$
  for $i=1,2$
  with the property that for all $\x,\y\in\Gen(\HD)$,
  homology classes $B\in\pi_2(\x,\y)$, and sources $\closedSource$ with
  $\ind(B,\closedSource) = 1$,
  the number of elements in $\cM^B(\x,\y,\closedSource)$ is
  equal (modulo $2$) to the number of elements in
  \[\bigcup_{\closedSource=\SourceSub{1}\glue\SourceSub{2}}
  \cMM^B(\x_1,\y_1,\SourceSub{1};\x_2,\y_2,\SourceSub{2}),\]
  where 
  $\x=\x_1\cup\x_2$ and $\y=\y_1\cup\y_2$.
\end{proposition}

\begin{proof}
  The result follows from standard compactness and gluing techniques,
  in a completely analogous way to Propositions~\ref{prop:compactness}
  (but using~\cite[Theorem 10.3]{BEHWZ03:CompactnessInSFT} in place
  of~\cite[Theorem 10.2]{BEHWZ03:CompactnessInSFT})
  and~\ref{prop:gluing_simple_comb} (using a Morse-Bott analogue
  of~\cite[Proposition A.2]{Lipshitz06:CylindricalHF}, say).
\end{proof}

As always, the curves of primary importance to us are the embedded
ones.
For a Heegaard diagram $\HD$ for a closed manifold with genus~$g$,
$\x,\y\in\Gen(\HD)$, and
$B\in\pi_2(\x,\y)$, define
\begin{align*}
  \chi_\emb(B)&\coloneqq g+e(B)-n_\x(B)-n_\y(B) \\
  \ind(B)&\coloneqq e(B)+n_\x(B)+n_\y(B).
\end{align*}
(This is a special case of Definition~\ref{def:emb-ind-emb-chi}.)

\begin{definition}\label{def:emb-matched}
  The \emph{embedded matched moduli space}
\index{moduli space!of matched pairs!embedded}%
  $\gls*{MatchModEmb}$
 is the union of
  $\tcMM{}^B(\x_1,\y_1,\SourceSub{1}\semico\x_2,\y_2,\SourceSub{2})$ over all
  pairs $\SourceSub{1}$, $\SourceSub{2}$ with
  $\chi(\SourceSub{1}\glue\SourceSub{2}) = \chi_\emb(B)$.  (Note there
  are only finitely many terms in this union.)  Also define
  \gls*{MatchModEmbUnparam}
  to be the quotient of
  $\tcMM{}^B(\x_1,\y_1\semico\x_2,\y_2)$ by the $\RR$-action (assuming
  $B\neq 0$, so the $\RR$-action is free).
\end{definition}

It follows from Lemma~\ref{lemma:Matched-Exp-Dim} that for generic
$J_1$ and $J_2$, $\cMM^B(\x_1,\y_1\semico\x_2,\y_2)$ is a manifold of
dimension $\ind(B)-1$.

\begin{lemma}\label{lem:matched-embedded}
  A matched pair of curves $(u_1,u_2)$ is in the corresponding
  embedded matched moduli space $\tcMM{}^B(\x_1,\y_1\semico\x_2,\y_2)$ if and
  only if both $u_i$ are embedded and every pair of Reeb chords
  appearing at the same height in the $u_i$ are either nested or disjoint.
\end{lemma}

\begin{proof}
  Let $(u_1,u_2)$ be any matched pair of curves.
  Let $\SourceSub{i}$ be the source of $u_i$, and let $\vec\rhos$ be
  sequence of sets of Reeb chords on $u_1$.  Let $m =
  \abs{E(\SourceSub{1})}$ and $B_1 \glue B_2 = B$.  By 
  Proposition~\ref{prop:asympt_gives_chi},
  $\chi(S_1)\geq\chi_\emb(B_1, \vec\rhos)$ and
  $\chi(S_2)\geq\chi_\emb(B_2, -\vec\rhos)$.

  We claim that
  \[\iota(\vec\rhos) + \iota(-\vec\rhos) \le -m,\]
  as follows.  The $L(\rhos_i,\rhos_j)$ terms in the definition of $\iota$
  (Formula~\eqref{eq:def-iota}) contribute oppositely to
  $\iota(\vec\rhos)$ and $\iota(-\vec\rhos)$, as the orientation on the
  circle is reversed, so the only contributions to $\iota(\vec\rhos) +
  \iota(-\vec\rhos)$ are from the
  $\iota(\rhos_i)$ terms.  From Lemma~\ref{lem:iota-chords} we see
  that each chord contributes $-1$, and each pair of
  chords $\rho_1$, $\rho_2$ at the same level has a further negative contribution of
  $-2\abs{L(\rho_1,\rho_2)}$, from which the claim follows.

  Thus
  \begin{align*}
  \chi(\SourceSub{1}\glue\SourceSub{2}) &= \chi(S_1) + \chi(S_2) - m\\
    &\geq \chi_\emb(B_1,\vec\rhos) + \chi_\emb(B_2,-\vec\rhos) - m\\
    &= g + e(B) - n_\x(B) - n_\y(B) - \iota(\vec\rhos) -
       \iota(-\vec\rhos) - m\\
    &\ge \chi_\emb(B).
  \end{align*}
  The pair $(u_1,u_2)$ is in the embedded moduli space if and only if
  we have equality at each step, which is equivalent to both $u_1$ and $u_2$ being
  embedded and all $\abs{L(\rho_1,\rho_2)}$ terms vanishing.
  These last terms vanish if and only if for all $\rho_1$ and $\rho_2$
  on the same level, $\rho_1$ and $\rho_2$ are
  nested or disjoint, as desired.
\end{proof}

\begin{definition}\label{def:paired-complex}
Let $\CFa(\HD_1,\HD_2)$ be the $\Field$-vector space with basis those
$\x_1\times\x_2\in\Gen(\HD_1)\times\Gen(\HD_2)$ so that $\x_1\cup\x_2$
is in $\Gen(\HD_1\cup_\bdy\HD_2)$.  Define a boundary operator
$\bdy_1$ on $\CFa(\HD_1,\HD_2)$ by
\begin{equation*}
  \bdy_1(\x_1\times\x_2) \coloneqq
    \sum_{\y_1,\y_2} \sum_{\{B\mid\ind(B)=1\}}\!\!
    \#\bigl(\cMM^B(\x_1,\y_1;\x_2,\y_2)\bigr) \cdot (\y_1\times\y_2).
\end{equation*}
Since, by Lemma~\ref{lem:closed-admissible}, $\HD_1\cup_\bdy\HD_2$ is
admissible, the argument from Lemma~\ref{lem:finite-typeD} shows that
this sum is finite.
\end{definition}

\begin{theorem}
  \label{thm:PrimitivePairing}
  $\CFa(\HD_1,\HD_2)$ is a chain complex isomorphic to the complex
$\CFa(\HD_1\cup_\bdy\HD_2)$, for suitable choices of almost-complex
structures.
\end{theorem}

\begin{proof}
  According to \cite[Theorem~2]{Lipshitz06:CylindricalHF},
  $\CFa(\HD_1\cup_\bdy\HD_2)$ can be calculated by the chain complex
  whose differential
  counts points in $\cM^B(\x,\y,\closedSource)$, where we take the union
  over all choices of $B\in\pi_2(\x,\y)$ and compatible $\closedSource$ so
  that $\ind(\closedSource,B)=\ind(B)=1$.
  According to Proposition~\ref{prop:Gluing1}, these counts
  agree with~$\bdy_1$.
\end{proof}

\begin{remark}
  Theorem~\ref{thm:PrimitivePairing} shows, in particular, that
  $\bdy_1^2=0$. This also follows as a special case of
  Proposition~\ref{prop:DSquaredZero}, below.
\end{remark}

\section{Dilating time}
\label{sec:dilating-time}

We now generalize the notion of matched pairs by inserting a real
parameter~$\gls*{Tdilation}$. 

\begin{definition}
  \label{def:TMatched}
  Fix a real number $T>0$.  Fix two compatible pairs of generators
  $\x_i,\y_i\in\Gen(\HD_i)$ ($i=1,2$). Fix compatible
  sources $\SourceSub{1}$ and $\SourceSub{2}$ connecting $\x_1$ to
  $\y_1$ and $\x_2$ to $\y_2$ respectively.  The \emph{moduli space of
    $T$-matched pairs}
  \index{moduli space!of T-matched pairs@of $T$-matched pairs}%
  $\gls*{TMatchMod}$
  is defined to be the fibered product
  $$\tcM^{B_1}(\x_1,\y_1;\SourceSub{1})\times_{T\cdot\ev_1=\ev_2}
  \tcM^{B_2}(\x_2,\y_2;\SourceSub{2}).$$
  That is, the moduli space of $T$-matched pairs consists of pairs $(u_1,u_2)$ with
  $u_i\in\tcM^{B_i}(\x_i,\y_i\semico\SourceSub{i})$
  and $T\cdot\ev(u_1)=\ev(u_2)$. Let
  \[
  \gls*{TMatchModEmb}
  \coloneqq\!\!\!
  \bigcup_{\{\SourceSub{1},\SourceSub{2}\mid
    \chi(\SourceSub{1}\glue\SourceSub{2})=\chi_\emb(B)\}}
  \!\!\!\tcMM^B(T\semico\x_1,\y_1,\SourceSub{1}\semico
  \x_2,\y_2,\SourceSub{2})
  \]
  denote the embedded moduli space of $T$-matched pairs.
\end{definition}

The index of a $T$-matched pair does not depend on the
$T$-parameter, and so is still given by
$\ind(B_1,\SourceSub{1}\semico B_2,\SourceSub{2})$ as defined in 
Equation~\eqref{eq:IndexMatchedPair}.

\begin{definition}
  The moduli space of $T$-matched pairs has an $\RR$-action given by
  $\tau_{t}(u_1,u_2)\coloneqq(\tau_{T\cdot t}(u_1),\tau_{t}(u_2))$.
  \index{$\RR$-action!on matched pairs}%
  Let
  \begin{align*}
    \gls*{TMatchModUnparam}
    &\coloneqq\tcMM^B(T\semico\x_1,\y_1,\SourceSub{1}\semico
    \x_2,\y_2,\SourceSub{2})/\RR\\
    \gls*{TMatchModEmbUnparam}
    &\coloneqq\tcMM^B(T\semico\x_1,\y_1\semico\x_2,\y_2)/\RR.
  \end{align*}
\end{definition}

\begin{lemma}
  \label{lem:MatchedCurveTransversality}
  For generic admissible almost complex structures $J_i$ on
  $\Sigma_i\times[0,1]\times\RR$, for all choices of
  $\x_i,\y_i\in\Gen(\HD_i)$ and $B
  \in\pi_2(\x_1\cup\x_2,\y_1\cup\y_2)$ so that
  $B \ne 0$,
  \begin{itemize}
    \item for generic values of~$T$, $\cMM^{B}(T\semico\x_1,\y_1\semico
      \x_2,\y_2)$ is a manifold of dimension
      $\ind(B)-1$, and
    \item $\bigcup_{T>0}\cMM^{B}(T\semico\x_1,\y_1\semico
      \x_2,\y_2)$ is a manifold of dimension
      $\ind(B)$.
    \end{itemize}
\end{lemma}
\begin{proof}
  As in Lemma~\ref{lemma:Matched-Exp-Dim}, the first statement follows
  from Proposition~\ref{prop:transversality} and
  Equation~\eqref{eq:Index}. The proof of the second
  statement proceeds in a standard way; see, for instance, \cite[Section
  3.4]{MS04:HolomorphicCurvesSymplecticTopology} for a nice
  explanation of this type of argument.
\end{proof}

\begin{definition}
  For $T\in (0,\infty)$, let $\CFa(T;\HD_1,\HD_2)$ be the vector space
  with the same basis as $\CFa(\HD_1,\HD_2)$ from
  Definition~\ref{def:paired-complex} and with boundary operator
  $$\partial_T \co \CFa(T;\HD_1,\HD_2)\to \CFa(T;\HD_1,\HD_2)$$
  defined by 
  \begin{equation*}
    \gls*{Tdiff}
    \coloneqq
    \sum_{\y_1,\y_2} \sum_{\{B\mid\ind(B)=1\}}\!\!
    \#\bigl(\cMM^B(T;\x_1,\y_1;\x_2,\y_2)\bigr) \cdot (\y_1\times\y_2).
  \end{equation*}%
\glsadd{CFTimeDilated}%
\end{definition}

The next goal is to show that $(\partial_T)^2=0$. 
As usual, this requires considering the ends of
$\cMM^B(T\semico\x_1,\y_1\semico\x_2,\y_2)$ when $\ind(B)=2$.  To this end, we
consider a compactification of
$\cMM^B(T\semico\x_1,\y_1,\SourceSub{1}\semico\x_2,\y_2,\SourceSub{2})$ by
holomorphic combs, following Section~\ref{sec:combs-compact}.

\begin{definition}\label{def:T-matched-comb}
  Given $\x_i,\y_i\in\Gen(\HD_i)$ for $i=1,2$, a \emph{$T$-matched
    story}
  \index{T-matched@$T$-matched!story}%
  \index{holomorphic story!T-matched@$T$-matched}%
  from $\x=(\x_1,\x_2)$ to $\y=(\y_1,\y_2)$ is a sequence
  $(u_1,v_1,\dots,v_k,u_2)$ where 
  \begin{itemize}
  \item 
    $u_1\in\cM^{B_1}(\x_1,\y_1\semico\SourceSub{1})$,
    $u_2\in\cM^{B_2}(\x_2,\y_2\semico\SourceSub{2})$, 
    and $v_j\in\cN(\biSource_j)$ for
    some sources $\SourceSub{i}$ and $\biSource_j$;
  \item the sources are equipped with
    bijective correspondences between $E(\SourceSub{1})$ and
    $W(\biSource_1)$, $E(\biSource_i)$ and $W(\biSource_{i+1})$, and
    $E(\biSource_k)$ and $E(\SourceSub{2})$;
  \item the correspondences
    preserve the labelings by Reeb chords, with orientation reversal
    between $E(\biSource_k)$ and $E(\SourceSub{2})$; and
  \item the following compatibility conditions hold:
    \begin{align*}
      \ev(u_1) &= \ev_w(v_1),\\
      \ev_e(v_i) &= \ev_w(v_{i+1})\textrm{ and}\\
      T\cdot\ev_e(v_k) &= \ev (u_2).
    \end{align*}
  \end{itemize}
  We call $u_1$ the \emph{west-most level}\index{west-most level}
  of the $T$-matched story and $u_2$
  the \emph{east-most level}\index{east-most level}
  of the $T$-matched story.
  A $T$-matched story is \emph{stable}
  \index{stable!T-matched story@$T$-matched story}%
  if either $u_1$ or
  $u_2$ is stable.  (In particular, if $\abs{E(\SourceSub{1})} > 0$,
  the comb is automatically stable.)

  A \emph{$T$-matched comb}
  \index{T-matched@$T$-matched!comb}%
  \index{holomorphic comb!T-matched@$T$-matched}%
  of height~$N$ is a sequence of stable
  $T$-matched stories running from $\x^j=(\x^j_1,\x^j_2)$ to
  $\x^{j+1}=(\x^{j+1}_1,\x^{j+1}_2)$ for some sequences of generators $\x^j_1$ in
  $\HD_1$ and
  $\x^j_2$ in $\HD_2$ (for $j=1,\dots,N$).
\end{definition}

As in Section~\ref{sec:combs-compact}, there is a natural
compactification
$\glsadd{MatchModCpct}\ocMM^B(\x_1,\y_1,\SourceSub{1}\semico\x_2,\y_2,\penalty700\SourceSub{2})$ 
of
$\cMM^B(\x_1,\y_1,\SourceSub{1}\semico\x_2,\y_2,\SourceSub{2})$ inside the
space of stable $T$-matched holomorphic combs.
The following is the analogue of
Proposition~\ref{prop:restrict_degens_1}, and the proof is similar.

\begin{proposition}
  \label{prop:EndsOfTwoDimModuliSpacesT}
  Suppose that
  $\ind(B_1,\SourceSub{1}\semico B_2,\SourceSub{2})=2$. Then 
  for generic $J$, every $T$-matched comb in
  $\bdy\ocMM^B(T\semico\x_1,\y_1,\SourceSub{1}\semico
  \x_2,\y_2,\SourceSub{2})$
  has one of the following
  forms:
  \begin{enumerate}
  \item a height $2$ $T$-matched comb;
  \item a $T$-matched comb $(u_1,v,u_2)$ where the only non-trivial
    component of~$v$ is a join component; or
  \item a $T$-matched comb $(u_1,v,u_2)$ where the only non-trivial
    component of~$v$ is a split component.
  \end{enumerate}
\end{proposition}

\begin{proof}
  The proof of Lemma~\ref{lem:ModuliSpacesAreMonotone} applies to show
  that the components of a $T$-matched pair are strongly monotone, as
  well. Thus, boundary double points cannot appear on the boundaries
  of the $T$-matched moduli spaces, owing to the strong boundary
  monotonicity of the moduli spaces, Lemma~\ref{lemma:NoBoundaryDoublePoints}.
  Boundary degenerations are ruled out by
  Lemma~\ref{lem:NoBoundaryDegenerations}, and ghost curves have
  codimension $2$ as in Lemma~\ref{lemma:no-ghosts}.

  Now suppose that we have a height $1$ $T$-matched comb
  $(u_1,v_1,\dots,v_\ell,u_2)$ in the boundary of the moduli space, with
  source
  $(\SourceSub{1}',\biSource_1,\dots,\biSource_\ell,\SourceSub{2}')$.  Let
  $\biSource = \biSource_1\glue\dots\glue\biSource_\ell$.  By the
  index hypothesis,
  \[
    2 = g - \chi(\SourceSub{1}) - \chi(\SourceSub 2) + 2 e(B) + m_0,
  \]
  where $m_0 = \abs{E(\SourceSub{1})} = \abs{E(\SourceSub{2})}$.  We
  also have $\SourceSub{1} \glue \SourceSub{2} \cong
  \SourceSub{1}'\glue\biSource\glue\SourceSub{2}'$ and so
  \[
  \chi(\SourceSub{1}) + \chi(\SourceSub{2}) - m_0 =
  \chi(\SourceSub{1}') + \chi(\biSource) + \chi(\SourceSub{2}')
  -m_1 - m_2
  \]
  where $m_1 = \abs{E(\SourceSub{1}')}$ and $m_2 = \abs{E(\SourceSub{2}')}$.
  Let $k$ be the number of components of $\biSource$.
  Let $P_i$ be the partition of $E(\SourceSub{i}')$ according to the
  components of $\biSource$, so that $\abs{P_i} = k$.
  Then an upper bound for the dimension of the
  space of limit curves for generic~$J$ is
  \begin{multline*}
    \ind(B_1,\SourceSub{1}',P_1) + 
    \ind(B_2,\SourceSub{2}',P_2) - k -1\\
    \begin{aligned}
    &= g - \chi(\SourceSub{1}') - \chi(\SourceSub{2}') + 2 e(B) + k-1\\
    &= \chi(\SourceSub{1}) - \chi(\SourceSub{1}') +
      \chi(\SourceSub{2}) - \chi(\SourceSub{2}') - m_0+k+1\\
    &= (\chi(\biSource) - k) + (k - m_1) + (k-m_2) + 1.
    \end{aligned}
  \end{multline*}
  All three terms in parentheses are non-positive, so if
  the space of limit curves is non-empty, at most one can be negative.
  If all are equal to~$0$, $\biSource$ consists of trivial strips and
  there is no degeneration. If
  $\chi(\biSource) < k$, one component of $\biSource$ is an annulus,
  which is ruled out as in Proposition~\ref{prop:restrict_degens_1}.
  If $m_1 = k+1$ or $m_2 = k+1$, then $\biSource$ is a join or split
  curve, respectively.

  The arguments above also show that the index of a story in the
  degenerate curve is at least~$1$, and that if the index is~$1$ there
  are no components at east infinity.  Thus the only other possibility
  for a degeneration is a height $2$ matched curve.
\end{proof}

In our setting, the matched combs with components at east infinity
cancel in pairs, as boundary branch points move between $\Sigma_1$
and~$\Sigma_2$. Roughly speaking, the next proposition states
that rigid $T$-matched combs $(u_1,v,u_2)$ appear with 
multiplicity one in the boundaries of their corresponding
one-dimensional moduli spaces.

\begin{proposition}
  \label{prop:CombsCancel} Let $(u_1,v,u_2)$ be a $T$-matched comb of
  index two with source $(\SourceSub{1}$, $\biSource$,
  $\SourceSub{2}$), where $\biSource$ has one non-trivial component,
  which is either a join or split component. Then there are arbitrarily small
  open
  neighborhoods $U_1$ of $(u_1,v,u_2)$ in
  \[
  \ocM^{B_1}(\x_1,\y_1,\SourceSub{1}\glue\biSource)\times
  \ocM^{B_2}(\x_2,\y_2,\SourceSub{2})
  \]
  with the property that $\partial {\overline{U_1}}$  meets $\cMM^{B}(T\semico
  \x_1,\y_1,\SourceSub{1}\glue\biSource\semico\x_2,\y_2,\SourceSub{2})$
  in an odd number of points.
  Similarly,
  there are arbitrarily small open neighborhoods $U_2$ of $(u_1,v,u_2)$ in
  \[
  \ocM^{B_1}(\x_1,\y_1,\SourceSub{1})\times
  \ocM^{B_2}(\x_2,\y_2,\biSource\glue\SourceSub{2})
  \]
  with the property that $\partial \overline {U_2}$ meets 
  $\cMM^{B}(T\semico\x_1,\y_1,\SourceSub{1}\semico\x_2,\y_2,\biSource\glue\SourceSub{2})$
  in an odd number of points.
\end{proposition}

\begin{proof}
  Suppose that $v$ is a split curve.  We label the east
  punctures of $\SourceSub{1}\glue\biSource$ by $\{e_i\}_{i=0}^{\ell}$ so that $e_0$
  and $e_1$ are the two east punctures which are in the same component in the split
  curve.  Let $\{w_i\}_{i=0}^{\ell}$ be the corresponding punctures of
  $\SourceSub{2}$.
  There are open neighborhoods $N_1\subset \ocM^{B_1}(\x_1,\y_1;\SourceSub{1}\glue \biSource)$ of $(u_1,v)$ and $N_2\subset \ocM^{B_2}(\x_2,\y_2;\SourceSub{2})$ of $u_2$ which admit evaluation maps
  \begin{align*} 
   \ev^1&\co N_1
   \to \RR^{\ell} &
   \ev^2& \co
  N_2 \to \RR^{\ell},
  \end{align*} 
  where the $i\th$ coordinate of $\ev^2$
  is $\ev_{w_i,w_0}$
  (as in Equation~\eqref{eq:DefEvPQ}), and the $i\th$ coordinate
  of $\ev^1$ is $T\cdot\ev_{e_i,e_0}$.
  Let
  $$f_j\co N_j\to \RR$$
  be the first component of $\ev^j$ for $j=1,2$, i.e., $f_2=\ev_{w_1}-\ev_{w_0}$
  and $f_1=\ev_{e_1}-\ev_{e_0}$.

  By the assumptions on the index, $\dim N_1 + \dim N_2 = \ell+1$, so
  if we knew that the maps $\ev^1$ and $\ev^2$ were smooth (and
  transverse) it would
  follow that the images intersect in a 1-manifold, which would prove
  the result.  The difficulty is that we do not know that $\ev^1$ is
  smooth on the boundary. Instead, we will use linking number
  considerations.

  Let $U$ be a small neighborhood of $u_1$ inside
  $\ocM^{B_1}(\x_1,\y_1;\SourceSub{1})$. Gluing
  (Proposition~\ref{prop:gluing_simple_comb})
  gives a map 
  $$\gamma
  \co U\times [0,\epsilon] \to
N_1\subset \ocM^{B_1}(\x_1,\y_1;\SourceSub{1}\glue\biSource)$$
which is a
homeomorphism onto a neighborhood of $(u_1,v)$. Shrinking $U$
slightly, we may assume that $\gamma$ is defined on
$\overline{U}\times[0,\epsilon]$. Note that
$\gamma(\overline{U}\times\{0\})$ consists of pairs $(u'_1,v')$ where
$v'$ is a split curve whose height is determined by $u'_1$.
Choose $U$ so that
\begin{enumerate}
\item $\overline{U}$ is a ball contained in
  $\cM^{B_1}(\x_1,\y_1;\SourceSub{1})\subset\ocM^{B_1}(\x_1,\y_1;\SourceSub{1})$ and
\item\label{item:combs-U-small} $U$ is sufficiently small that
$\ev^1(\gamma({\overline U}\times \{0\}))\cap\ev^2(N_2)$
consists of the single point $\ev^1((u_1,v))=\ev^1(\gamma(u_1\times0))$.
\end{enumerate}
To see that the second requirement can be met, note that $\ev^1(\gamma({\overline U}\times
\{0\}))=\{0\}\times\ev^1(\overline{U})$; and by
Proposition~\ref{prop:transversality} we may assume that $\ev^2$ is
transverse to this submanifold (compare
Lemma~\ref{lemma:Matched-Exp-Dim}).

The rest of the argument is similar to the proof of
Proposition~\ref{prop:gluing_degree_one}. 
The intersection $\ev^1(\gamma(\partial{\overline
  U}\times \{0\}))\cap\ev^2(N_2)$ is empty and hence, for sufficiently
small~$\epsilon$:
\begin{enumerate}[resume]
\item\label{item:combs-eps-small} The intersection $\ev^1(\gamma (\partial {\overline
    U}\times[0,\epsilon]))\cap\ev^2(N_2)$ is empty.
\end{enumerate}
Since $\gamma(\overline{U}\times\{\epsilon\})$ is contained in
$\cM^{B_1}(\x_1,\y_1;\SourceSub{1}\glue\biSource)\subset
\ocM^{B_1}(\x_1,\y_1;\SourceSub{1}\glue\biSource)$,
$f_1$ is positive on $\gamma(\overline{U}\times\{\epsilon\})$.  
Therefore:
\begin{enumerate}[resume]
\item\label{item:comps-delta-small} For sufficiently small $\delta>0$,
  $\ev^2(f_2^{-1}((-\infty,\delta]))$ is disjoint from
  $\ev^1(\gamma(U\times\{\epsilon\}))$.
\end{enumerate}
Consider now the ball
$D_1=\overline{U}\times[0,\epsilon]$ and choose a closed ball $D_2 \subset
N_2$ which is a neighborhood of $u_2$ and contained in
$f_2^{-1}((-\infty,\delta])$.  Combining \eqref{item:combs-U-small},
\eqref{item:combs-eps-small} and \eqref{item:comps-delta-small}, $\ev^1(\bdy D_1)$ intersects
$\ev^2(D_2)$ only in the single point $\ev^1(u_1,v)$, which is in the image
$\ev^2(D_2^\circ)$ of the interior $D_2^\circ$ of $D_2$.  Therefore linking number
considerations show that $\ev^1(D_1)\cap\ev^2(\bdy D_2)$ consists of
an odd number of points.  We can therefore take $U_1$ in the statement
of the proposition to be $\gamma(U\times[0,\epsilon)) \times
D_2^\circ$: then, the intersection of $\bdy\overline{U_1}$ with $\ocMM^{B}(T\semico
  \x_1,\y_1,\SourceSub{1}\glue\biSource\semico\x_2,\y_2,\SourceSub{2})$ consists of the point
$(u_1,v,u_2)$, from $(\bdy\gamma(\overline{U}\times[0,\epsilon))) \times
D_2^\circ$,
which is not in $\cMM^{B}(T\semico
  \x_1,\y_1,\SourceSub{1}\glue\biSource\semico\x_2,\y_2,\SourceSub{2})$, and an odd number
of other points, from $\gamma(U\times[0,\epsilon)) \times
(\bdy D_2)$.

A similar analysis holds for the other moduli space or if $v$ has a
join component.
\end{proof}

\begin{proposition}
  \label{prop:DSquaredZero}
  The map $\partial_T$ is a differential.
\end{proposition}

\begin{proof}
  Fix $\x_i,\y_i\in\Gen(\HD_i)$ and a homology class
  $B\in\pi_2(\x_1\cup\x_2,\y_1\cup\y_2)$
  with $\ind(B)=2$. Consider the ends of the moduli space
  $\ocM^B(T\semico\x_1,\x_2\semico\y_1,\y_2)$ of embedded
  $T$-matched holomorphic curves representing $B$. According to
  Proposition~\ref{prop:EndsOfTwoDimModuliSpacesT}, these
  ends correspond either to simple matched combs, which cancel in pairs
  according to Proposition~\ref{prop:CombsCancel}, or to height $2$
  holomorphic buildings, with each story an index one
  $T$-matched curve. Indeed, each possible height $2$ building
  appears as an end an odd number of times, according to an extension of
  Proposition~\ref{prop:gluing_two_story} along the lines of
  Proposition~\ref{prop:CombsCancel}.  Hence, the number of such
  height $2$ buildings must be even. But the number of
  such terms (over all homology classes $B$ with index two)
  is the $\y_1\times\y_2$ coefficient of $(\partial_T)^2(\x_1\times\x_2)$.
\end{proof}

The chain complex $\CFa(1;\HD_1,\HD_2)$ is $\CFa(\HD_1\cup_\bdy\HD_2)$
by Theorem~\ref{thm:PrimitivePairing}.
We will next show that the homotopy type of $\CFa(T;\HD_1,\HD_2)$ is
independent of the choice of~$T$.

\begin{definition}
  \label{def:PhiMatched} Let $T_1$ and $T_2$ be positive real
  numbers, and fix a smooth function $\psi\co \RR\to
  \RR$ with positive derivative so that
  $$
  \gls*{psidilation}
  =
  \begin{cases}
    T_1\cdot t & t\leq -1\\
    T_2\cdot t & t\geq 1.
  \end{cases}
  $$ 
  There are induced maps $\psi^m\co \RR^m\to \RR^m$
  defined by $\psi^m(t_1,\dots,t_m)=(\psi(t_1),\dots,\allowbreak\psi(t_m))$.
  Fix two compatible pairs of generators $\x_i, \y_i\in\Gen(\HD_i)$ ($i=1,2$).
  Fix compatible sources $\SourceSub{1}$ and $\SourceSub{2}$
  connecting $\x_1$ to $\y_1$ and $\x_2$ to $\y_2$ respectively.  The
  {\em moduli space of $\psi$-matched pairs},
  \index{moduli space!of $\psi$-matched pairs}%
  denoted
  \[
  \gls*{psiMatchMod}
  \]
  is defined as the fibered product
  $$\cM^{B_1}(\x_1,\y_1;\SourceSub{1})\times_{\psi^m\circ\ev_1=\ev_2}
  \cM^{B_1}(\x_2,\y_2;\SourceSub{2}),$$
  where $m = \abs{E(\SourceSub{1})}$.
  That is, the moduli space consists of pairs
  $(u_1,u_2)$ with $u_i\in\cM^{B_i}(\x_i,\y_i\semico\SourceSub{i})$ such that
  $\psi^m\circ\ev(u_1)=\ev(u_2)$.
  Let 
  \[
  \gls*{psiMatchModEmb}\coloneqq\!\!\!
  \bigcup_{{\{\SourceSub{1},\SourceSub{2}\mid
      \chi(\SourceSub{1}\glue\SourceSub{2})=\chi_\emb(B)\}}}\!\!\!
  \gls*{psiMatchMod}
  \]
  denote the space of embedded $\psi$-matched curves.

  The space of $\psi$-matched pairs has a natural compactification by
  \emph{$\psi$-matched combs}, defined analogously to
  Definition~\ref{def:T-matched-comb}, except using the map induced by
  $\psi$ rather than simply rescaling on the left
  side.\index{holomorphic comb!$\psi$-matched}
\end{definition}

The moduli spaces of $\psi$-matched pairs do not usually have a
natural action by~$\RR$.\index{$\RR$-action!on matched pairs}

For $\psi$ as above, define 
$$F_{T_1,T_2}\co \CFa(T_1;\HD_1,\HD_2)\to
\CFa(T_2;\HD_1,\HD_2)$$ 
by 
\begin{equation*}
  \gls*{Tcont}
  (\x_1\times\x_2)
  =
  \sum_{\y_1,\y_2}
  \sum_{\{B\mid\ind(B) = 0\}}\!\!
    \#\bigl(\cMM^B(\psi;\x_1,\y_1;\x_2,\y_2)\bigr)
  \cdot (\y_1\times\y_2).
\end{equation*}

We sketch now familiar arguments (generalizing the above proof of
Proposition~\ref{prop:DSquaredZero}) which  show that $F_{T_1,T_2}$ is a
chain map.

\begin{lemma}
  \label{lem:EndsOfTwoDimModuliSpacesPsi}
  Suppose that $\cM^B(\psi\semico\x_1,\y_1,\SourceSub{1}\semico
  \x_2,\y_2,\SourceSub{2})$ is one-dimensional. Then for a
  generic almost complex structure~$J$, every $\psi$-matched comb in
  $\bdy\ocM^B(\psi\semico\x_1,\y_1,\penalty700\SourceSub{1}\semico\x_2,\y_2,\SourceSub{2})$
  has one of the following
  forms:
  \begin{enumerate}
  \item a height $2$ matched curve $(u_1',u_2')*(u_1'',u_2'')$, where
    $(u_1',u_2')$ is a $\psi$-matched curve and
    $(u_1'',u_2'')$ is a $T_2$-matched curve;
  \item a height $2$ matched curve $(u_1',u_2')*(u_1'',u_2'')$, where
    $(u_1',u_2')$ is a $T_1$-matched curve and
    $(u_1'',u_2'')$ is a $\psi$-matched curve;
  \item a $\psi$-matched comb $(u_1,v,u_2)$ where the only non-trivial
    component of~$v$ is a join component; or
  \item a $\psi$-matched comb $(u_1,v,u_2)$ where the only non-trivial
    component of~$v$ is a split component.
  \end{enumerate}
\end{lemma}

\begin{proof}
  This follows along the same lines of
  Proposition~\ref{prop:EndsOfTwoDimModuliSpacesT} above.
  (Observe, once again, that the proof of Lemma~\ref{lem:ModuliSpacesAreMonotone}
  applies to show that $\psi$-matched combs also consist of boundary monotone factors.)
\end{proof}

\begin{proposition}
  \label{prop:VaryTChain}
  $F_{T_1,T_2}\colon \CFa(T_1;\HD_1, \HD_2)\to \CFa(T_2\semico\HD_1,\HD_2)$
  is a chain homotopy equivalence. In particular, for any $T>0$, $\CFa(T;\HD_1,\HD_2)$ is
  chain homotopy equivalent to $\CFa(\HD)$.
\end{proposition}

\begin{proof}
  For fixed $\x_i,\y_i\in\Gen(\HD_i)$ and homology class $B$ with
  $\ind(B) = 1$, consider the moduli space $\cM^B(\psi\semico\x_1,\y_1\semico
  \x_2,\y_2)$.
  Lemma~\ref{lem:EndsOfTwoDimModuliSpacesPsi} accounts for the ends of
  these moduli spaces. A straightforward adaptation of
  Proposition~\ref{prop:CombsCancel} shows that the ends enumerated in
  Lemma~\ref{lem:EndsOfTwoDimModuliSpacesPsi} which do not {\em a
    priori} cancel in pairs are the height $2$
  $\psi$-matched curves. Gluing as in the proof of
  Proposition~\ref{prop:DSquaredZero} shows that those ends are in
  one-to-one correspondence with the types of height $2$ buildings
  enumerated in Lemma~\ref{lem:EndsOfTwoDimModuliSpacesPsi}.
  But these counts, which must add up to zero, are precisely the
  coefficients of $\y_1\times\y_2$ in $\partial_{T_2} \circ
  F_{T_1,T_2}(\x_1\times\x_2)+F_{T_1,T_2}\circ \partial_{T_1}(\x_1\times\x_2)$.
  
  Showing that $F_{T_1,T_2}$ is a chain homotopy equivalence also
  follows along familiar lines in Floer homology. The homotopy inverse
  is provided by $F_{T_2,T_1}$ defined with a suitable function
  $\psi_{21}$, and the chain homotopy to the identity
  is provided by counting curves matched by a
  suitable family of functions, depending on a parameter~$c$,
  interpolating between
  multiplication by $T_1$ and a function which agrees with $\psi(t)$
  for $t$ sufficiently small, and $\psi_{21}(-t+c)$ for $t$
  sufficiently large.

  The second statement is immediate from this and
  Theorem~\ref{thm:PrimitivePairing}.
\end{proof}

\section{Dilating to infinity}
\label{sec:dilating-infinity}
Now that we have established suitable independence of
$\CFa(T\semico\HD_1,\HD_2)$ from $T$, we turn to the large $T$
behavior of this complex. The work, of course, is in understanding the
behavior of $T$-matched curves as $T\to\infty$. In
Section~\ref{sec:ideal-cpctness} we prove a compactness result,
showing that $T$-matched curves converge to so-called ideal-matched
curves as $T\to\infty$. While general ideal-matched curves can be
fairly complicated, we show in Section~\ref{sec:simple-ideal-curves}
that for rigid moduli spaces they are quite simple; this special class
of ideal-matched curves are called simple ideal-matched curves. It
turns out that simple ideal-matched curves contain too much
information at east $\infty$ to correspond one-to-one with $T$-matched
curves. In Section~\ref{sec:trimmed-gluing} we trim this extra
information and prove the desired gluing result.

\subsection{Ideal limits of \textalt{$T$}{T}-matched curves}\label{sec:ideal-cpctness}
\begin{definition}\label{def:time-parameter-space}
  Let $\overline{\RR} \coloneqq [-\infty,+\infty]$ be the standard
  compactification of~$\RR$.  For $U$ a holomorphic comb of
  height~$N$, define its \emph{time-parameter
    space}\index{time-parameter space}
  \gls*{TimeSpace}
  to be the
  union of $N$ copies of $\overline{\RR}$ modulo identifying $+\infty$
  in the $i\th$ copy with $-\infty$ in the $(i+1)\st$ copy:
  \[
  T(U) \coloneqq \biggl(\bigcup_{i=1}^N \overline{\RR}_i\biggr)\bigg/
    \bigl((+\infty)_i \sim (-\infty)_{i+1}\bigr).
  \]
  When $U$ is the trivial comb ($N=0$), let $T(U)$ be a single point.
\end{definition}

\begin{figure}
  \centering
  \includegraphics[scale=.83333]{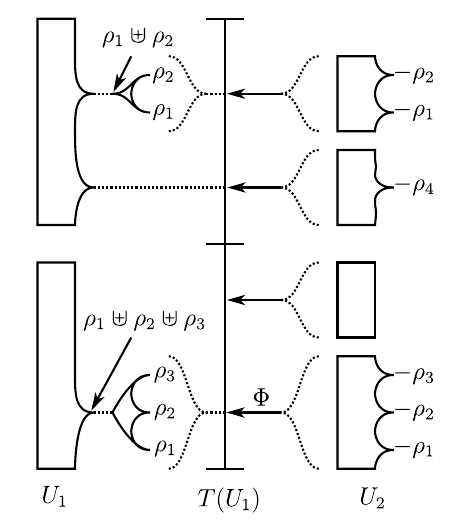}
  \caption[Schematic of an ideal-matched comb]{\textbf{Schematic
      illustration of an ideal-matched comb.} On the left is a height
    $2$ comb $U_1$. In the center is the time-parameter space
    $T(U_1)$. On the right is a height $4$ comb $U_2$. The arrows indicate the map $\Phi$ in the definition of an ideal-matched comb.}
  \label{fig:ideal-comb-schematic}
\end{figure}

\begin{definition}\label{def:ideal-matched-comb}
  An \emph{ideal-matched holomorphic comb}
  \index{holomorphic comb!ideal-matched}%
  \index{ideal-matched holomorphic comb}%
  consists of
  \begin{itemize}
  \item a pair $(U_1,U_2)$ of strongly boundary monotone, stable holomorphic combs,
    for $\HD_1$ and $\HD_2$, respectively;
  \item an order-preserving map~$\Phi$ from the set of stories of
    $U_2$ to $T(U_1)$; and
  \item a correspondence~$\gls*{phicorresp}$
    from $E(U_2)$ to $E(U_1)$,
  \end{itemize}
  so
  that, for every east-most puncture~$p$ of $U_2$,
  \begin{itemize}
  \item if $p$ is labeled by the Reeb chord $\rho$, then $\varphi(p)$
    is labeled by $-\rho$ and
  \item if $p$ is a puncture on the story~$C$ of $U_2$,
    the $t$-coordinate of $\varphi(p)$ is $\Phi(C)$.
  \end{itemize}
  See Figure~\ref{fig:ideal-comb-schematic} for a schematic illustration of an ideal-matched comb.
  As usual, we will often
  suppress $\varphi$ from the notation.  Let the \emph{total west
  source}~$\SourceSub{1}$\index{total west source} and \emph{total
  east source}~$\SourceSub{2}$\index{total east source}
  of $(U_1,U_2)$ be the result of pregluing all the sources of $U_1$
  and $U_2$, respectively.  We view such an ideal-matched comb as a point in
  \[ \ocM^{B_1}(\x_1,\y_1;\SourceSub{1};\vec P_1) \times
  \ocM^{B_2}(\x_2,\y_2;\SourceSub{2};P_2) \] where $\x_i$, $\y_i$, and
  $B_i$ are as usual, $\vec P_1$ is the induced partition on
  $E(\SourceSub{1})$, and $P_2$ is the discrete partition on
  $E(\SourceSub{2})$.

  Two ideal-matched combs are called \emph{equivalent} if they differ
  only by translating some of the components of $U_2$ and/or some of
  the components of $U_1$ together with the map $\Phi$.

  Let
  \[
  \ocMM^B(\infty;\x_1,\y_1,\SourceSub{1};\x_2,\y_2;\SourceSub{2})
  \]
  denote the moduli space of equivalence classes of ideal-matched
  combs in the homology class~$B$ connecting $\x_1,\x_2$ to
  $\y_1,\y_2$ and with total west and
  east sources $\SourceSub{1}$ and~$\SourceSub{2}$. Let
  \[
  \ocMM^B(\infty;\x_1,\y_1;\x_2,\y_2)\coloneqq\!\!
  \bigcup_{\substack{\SourceSub{1},\SourceSub{2}\\
    \chi(\SourceSub{1}\glue\SourceSub{2})=\chi_\emb(B)}}\!\!
  \ocMM^B(\infty;\x_1,\y_1,\SourceSub{1};\x_2,\y_2;\SourceSub{2}).
  \]
\end{definition}
In the definition, either $U_1$ or $U_2$ may be of height~$0$, in which
case the other is an ordinary provincial comb.

\begin{definition}\label{def:converge-to-ideal}
  Let $\{U^j\}_{j=1}^\infty$ be a sequence of $T_j$-matched combs,
  where $T_j\to \infty$.  Let $(U_1,U_2)$ be an ideal-matched
  holomorphic comb.  From each $T_j$-matched comb $U^j$ extract
  a pair of combs $u_1^j$ and $u_2^j$ by first taking $u_2^j$ to consist of the east-most
  levels of each story and $u_1^j$ to be the rest of each
  story and then throwing out any stories consisting
  of trivial strips.  We say that the sequence $\{U^j\}$
  \emph{converges}\index{convergence!to ideal-matched comb}
to $(U_1,U_2)$ if the following conditions are satisfied:
  \begin{itemize}
  \item The sequences $\{u_1^j\}$ and $\{u_2^j\}$ converge to $U_1$
    and $U_2$, respectively.
  \item For sufficiently large $j$
the identification between the east-most punctures
of~$u_1^j$ and the east punctures of $u_2^j$ 
corresponds to the identification between the 
east-most punctures of $U_1$ and the east-most punctures of~$U_2$.
\item For each story $C$ of $U_2$, $\Phi(C)\in T(U_1)$ is
  determined as follows.
Place marked points $p_j$ on the source of $u^2_j$ 
so that the sequence $\{p_j\}$ converges to a point on the source of~$C$.
Choose marked points~$q_j$ in the corresponding story of the source
of~$u^1_j$ so that
$$T_j\cdot t(u^1_j(q_j))= t(u^2_j(p_j)).$$
Then $\Phi(C)$ is the limit of $t(u^1_j(q_j))$ (which may be
$\pm\infty$).
\end{itemize}
\end{definition}
It is clear that if the third condition of
Definition~\ref{def:converge-to-ideal} is satisfied for one choice of
sequences $p_j$ and~$q_j$ for the story~$C$ then it holds for any
choice of sequences $p_j$ and~$q_j$. If $C$ has an $e\infty$
puncture then $\Phi(C)$ is given by the $\RR$-coordinate of the
corresponding $e\infty$ puncture of $U_1$, so $\Phi$ is determined by
the rest of the data in this case. By contrast, if $C$ is provincial
then $\Phi(C)$ is not determined by the rest of the data. Also,
$\Phi(C)=\infty$ is only possible for provincial stories $C$ of $U_2$.

\begin{proposition}
\label{prop:ideal-compactness}
Fix decorated sources $\SourceSub{1}$ and $\SourceSub{2}$ and homology
classes $B_1$ and $B_2$.  Let $\{(U_1^j,U_2^j)\}_{j=1}^\infty$ be a
sequence of $T_j$-matched holomorphic combs with total source
$(\SourceSub{1},\SourceSub{2})$ in the homology class $B_1\glue
B_2$. Suppose that $T_j\to\infty$. Then
$\{(U_1^j,U_2^j)\}_{j=1}^\infty$ has a subsequence which converges to
an ideal-matched holomorphic comb, with total source
$(\SourceSub{1},\SourceSub{2})$ and in the homology class $B_1\glue
B_2$.
\end{proposition}

\begin{proof}
  The key point is that, by Proposition~\ref{prop:compactness},
  the sequence $\{(U_1^j,U_2^j)\}_{j=1}^\infty$ limits to some pair of
  holomorphic combs. Strong boundary monotonicity follows from Lemma~\ref{lem:MonotonicityClosed}. 
  Moreover, since the evaluation maps are
  continuous on the compactified moduli space, the limit is an
  $\infty$-matched (i.e., ideal-matched) holomorphic comb.  More
  details follow.
  
  We prove the claim in the case that each $U_i^j$ is a holomorphic
  curve; the general case is similar, but the notation more
  complicated.  So, let $\{(u_1^j,u_2^j)\}$ be a sequence of
  $T_j$-matched holomorphic curves connecting $\x_1\otimes\x_2$ to
  $\y_1\otimes\y_2$ with west source~$\SourceSub{1}$ and east
  source~$\SourceSub{2}$. Replacing the sequence with a subsequence if
  necessary, we may assume that there is an ordering
  $q_{i,1},\dots,q_{i,n}$ of the punctures of $\SourceSub{i}$ so that
  if $a<b$ then for any $j$, we have
  $\ev_{q_{1,b},q_{1,a}}(u_1^j)\geq0$ (and so also
  $\ev_{q_{2,b},q_{2,a}}(u_2^j)\geq0$).

  Consider the space
  \[
  \ocM^{B_1}(\x_1,\y_1;\SourceSub{1};\vec{P}_1)\times \ocM^{B_2}(\x_2,\y_2;\SourceSub{2};\vec{P}_2)\times[1,\infty].
  \]
  where $\vec{P}_i$ is the sequence of sets of punctures $(\{q_{i,1}\},\dots,\{q_{i,n}\})$.
  By Proposition~\ref{prop:compactness}, this space is compact.
  On the open part, there are evaluation maps $\ev^i\co \cM^{B_i}(\x_i,\y_i\semico\SourceSub{i}\semico\vec{P}_i)\to
  \RR^{E(\SourceSub{i})}/\RR$, which we wish to extend to
  $\ocM^{B_i}(\x_i,\y_i;\SourceSub{i};\vec{P}_i)$. To this end, identify
  $\RR^{E(\SourceSub{i})}/\RR$ with $\RR^{n-1}$ by using the ordering
  on the east punctures $q_{i,1},\dots,q_{i,n}$ and considering
  differences of successive coordinates. Thus the map $\ev^i$ is
  given by
  \begin{multline*}
    \ev^i(u_i)\\
    \begin{aligned}
    &=(\ev_{q_{i,2}}(u_i)-\ev_{q_{i,1}}(u_i),\ev_{q_{i,3}}(u_i)-\ev_{q_{i,2}}(u_i),\dots,\ev_{q_{i,n}}(u_i)-\ev_{q_{i,n-1}}(u_i))\\
    &=(\ev_{q_{i,2},q_{i,1}}(u_i),\ev_{q_{i,3},q_{i,2}}(u_i),\dots,\ev_{q_{i,n},q_{i,n-1}}(u_i)).
    \end{aligned}
  \end{multline*}
  The image of $\ev^i$ lies inside $\RR_+^{n-1}=[0,\infty)^{n-1}$, and
  the map $\ev^i$ has a continuous extension to a map
  $\overline{\ev}^i\co \ocM^{B_i}(\x_i,\y_i;\SourceSub{i};\vec{P}_i)\to [0,\infty]^{n-1}.$
  
  Now, consider the map 
  \[
  f\co \ocM^{B_1}(\x_1,\y_1;\SourceSub{1};\vec{P}_1)\times
  \ocM^{B_2}(\x_2,\y_2;\SourceSub{2};\vec{P}_2)\times[1,\infty]\to [0,\infty]^{n-1}
  \]
  defined by
  \[
  f(u_1,u_2,T)=T\cdot \overline{\ev}^1(u_1)-\overline{\ev}^2(u_2).
  \]
  The preimage $f^{-1}(0)$ is the union over all $T\in[1,\infty)$ of the moduli
  space of $T$-matched holomorphic combs (with total sources
  $\SourceSub{1}$ and $\SourceSub{2}$, in the homology class
  $B_1\glue B_2$, with the partial order on the Reeb chords
  consistent with the order $q_{i,1},\dots,q_{i,n}$), together with
  the fiber over $T=\infty$, which we will discuss presently.

  The space $f^{-1}(0)$ is a closed subspace of a compact space, and
  hence compact. Consequently, the sequence of points
  $\{(u_1^j,u_2^j,T_j)\}_{j=1}^\infty$ in $f^{-1}(0)$ has a convergent
  subsequence. The limit $(U_1,U_2)$ of this subsequence lies in
  $f^{-1}(0)\cap\{T=\infty\}$. The limit $(U_1,U_2)$ can be endowed
  with the structure of an ideal-matched holomorphic comb as
  follows. The correspondence $\varphi$ is inherited from the
  $\varphi$ on each $(U_1^j,U_2^j)$. For each story of $U_2$ with at
  least one $e$ puncture, the map $\Phi$ is induced from the rest of
  the structure. For each story $u_{i}$ of $U_2$ with no $e$ puncture,
  place a marked point $p_i$ on the source of $u_i$; a sequence of
  marked points $p_{i,j}$ on the source of $u_2^j$ converging
  to~$p_i$; and a sequence of points~$q_{i,j}$ on $u_1^j$ with
  $T_j\cdot t(u_1^j(q_{i,j}))=t(u_2^j(p_{i,j}))$. Since the
  time-parameter space~$T(U_1)$ is compact, we
  may take a further subsequence of the $\{(u_1^j,u_2^j)\}$ so that for
  each $i$, the sequence $\{t(u_1^j(q_{i,j}))\}$ converges. Then $\Phi$ takes
  the $i\th$ story of $U_2$ to $\lim_{j\to\infty}t(q_{i,j})$. It is
  immediate from the definitions that our subsequence of
  $\{(u_1^j,u_2^j)\}$ converges to $(U_1,U_2,\phi,\Phi)$.
\end{proof}

\begin{remark}
  We have not claimed that the topology corresponding to
  Definition~\ref{def:converge-to-ideal} is Hausdorff.
\end{remark}

\subsection{Simple ideal-matched curves, the limits of rigid \textalt{$T$}{T}-matched curves}\label{sec:simple-ideal-curves}
We will see in this section that there are strong restrictions on which ideal-matched combs come from rigid $T$-matched curves.
\begin{definition}\index{simple ideal-matched curve|see{curve, simple ideal-matched}}\index{curve!ideal-matched!simple}
  A \emph{simple ideal-matched curve} is an ideal-matched comb
  $(U_1,U_2)$ in which
  \begin{itemize}
  \item $U_1$ has at most one story;
  \item $U_2$ has no components at east infinity;
  \item every component of $U_1$ at east infinity is a disk with a
    single west puncture; and
  \item exactly one of the following two conditions holds:
  \begin{enumerate}
  \item\label{case:ProvincialFlow} $E(\SourceSub{1}) = E(\SourceSub{2}) = \emptyset$, one of
    $U_1$ or $U_2$ is the trivial comb of height~$0$ and the other of
    $U_1$ or $U_2$ is a height $1$ holomorphic comb with index $1$, or
  \item\label{case:StandardFlow} each story of $U_2$ has at least one
    $e$~puncture and has index~$1$
    (with the discrete partition), the map $\Phi$ is injective, and $U_1$ has index $1$ (with the
    induced partition).
  \end{enumerate}
  \end{itemize}
\end{definition}
In the above definition, the east infinity components of $U_1$ are
like split components, except that there may be any number of east
punctures; in the terminology of
Section~\ref{sec:combs-gluing} they are
generalized split components.
\index{component!generalized split}%
Note that curves with generalized split components that are not split
components are not transversally cut with respect to the induced
partition.

\begin{lemma}\label{lem:sim-ind-1}
  Simple ideal-matched curves have
$\ind(B_1,\SourceSub{1};B_2,\SourceSub{2})=1$.
\end{lemma}
\begin{proof}
  In Case~\eqref{case:ProvincialFlow} in the definition of simple
  ideal-matched curves the statement is obvious.  To prove
  Case~\eqref{case:StandardFlow}, let
  $\vec{P}=(P_1,\dots,P_n)$ denote the induced ordered partition of the $e$
  punctures of $\SourceSub{1}$ (by their $\RR$-coordinates), and let
  $Q^i$ denote the discrete partition of the $e$ punctures of
  $\SourceSub{i}$. By definition,
  \[
  \ind(B_1,\SourceSub{1},\vec{P})=1 \qquad\qquad \ind(B_2,\SourceSub{2},Q^2)=n.
  \]
  Note that $|Q^1|=|Q^2|=\sum |P_i|$.
  By Formula~\eqref{eq:IndexMatchedPair},
  \begin{align*}
    \ind(B_1,\SourceSub{1};B_2,\SourceSub{2})&=\ind(B_1,\SourceSub{1},Q^1)+\ind(B_2,\SourceSub{2},Q^2)-|Q^1|\\
    &=\ind(B_1,\SourceSub{1},\vec{P})+\ind(B_2,\SourceSub{2},Q^2)-n\\
    &=1+n-n=1.\qedhere
  \end{align*}
\end{proof}

\begin{lemma}
  \label{lem:IdealMatchedTransversality}
  With respect to generic almost complex structures $J_1$ and $J_2$,
  if $\ind(B_1,\SourceSub{1}\semico B_2,\SourceSub{2}) < 1$ then
  $\ocMM^B(\infty;\x_1,\y_1,\SourceSub{1};\x_2,\y_2;\SourceSub{2})$ is
  empty, and if $\ind(B_1,\SourceSub{1}\semico B_2,\SourceSub{2}) = 1$
  then
  $\ocMM^B(\infty;\x_1,\y_1,\SourceSub{1};\x_2,\y_2;\SourceSub{2})$
  consists of simple ideal-matched curves $(U_1,U_2)$.
\end{lemma}

\begin{proof}
  This is similar to the proofs of
  Propositions~\ref{prop:restrict_degens_1}
  and~\ref{prop:EndsOfTwoDimModuliSpacesT}.  The case with no east
  punctures is immediate.  Otherwise, suppose
  first that $U_1$ has a single story, and for notational simplicity
  assume that $U_2$ also has a single story, so that all
  east punctures of $U_1$ appear at a single height.  Let
  $\SourceSub{1}'$ and $\SourceSub{2}'$ be the sources of the main
  components of $U_1$ and $U_2$, respectively, and similarly let
  $\biSource_1$ and $\biSource_2$ be the result of gluing together the
  horizontal levels of the sources of the east
  components of $U_1$ and~$U_2$.  Let $k_2$ be the number of
  components of $\biSource_2$, $m_0 =
  \abs{E(\SourceSub{1})}=\abs{E(\SourceSub{2})}$, and
  $m_i=\abs{E(\SourceSub{i}')}$ for $i=1,2$.  Let $P_i'$ be the
  induced partition on $E(\SourceSub{i}')$.  For generic~$J$, an upper
  bound for the dimension
  of $\ocMM^B(\infty;\x_1,\y_1,\SourceSub{1};\x_2,\y_2;\SourceSub{2})$ is
  \begin{multline*}
    \ind(B_1,\SourceSub{1}',P_1') + \ind(B_2,\SourceSub{2}',P_2') -2\\
    \begin{aligned}
      &= g + 2e(B) - \chi(\SourceSub{1}') - \chi(\SourceSub{2}') + k_2
         - 1\\
      &= g + 2e(B) - \bigl(\chi(\SourceSub{1}) - \chi(\biSource_1) + m_1\bigr)
           - \bigl(\chi(\SourceSub{2}) - \chi(\biSource_2) + m_2\bigr) + k_2-1\\
      &= \bigl(\chi(\biSource_1) - m_1\bigr) + \bigl(\chi(\biSource_2) - m_0\bigr) + \bigl(k_2
      - m_2\bigr) + \bigl(\ind(B_1,\SourceSub{1};B_2,\SourceSub{2})-1\bigr).
    \end{aligned}
  \end{multline*}
  When $\ind(B_1,\SourceSub{1}\semico B_2,\SourceSub{2})<1$, the sum is
  negative and so there are no ideal-matched combs.  When
  $\ind(B_1,\SourceSub{1}\semico B_2,\SourceSub{2}) = 1$, the three other terms in
  parentheses must be~$0$, which implies that $(U_1,U_2)$ is a simple
  ideal-matched curve as stated.

  The case when $U_2$ has multiple stories is similar.
  If $U_1$ has more than one story, each story must contribute at
  least one to $\ind(B_1, \SourceSub{1}; B_2, \SourceSub{2})$ by the
  above analysis, contradicting the index
  assumption.
\end{proof}

\subsection{Trimmed curves and gluing}\label{sec:trimmed-gluing}
We showed in Section~\ref{sec:ideal-cpctness} that $T$-matched curves
converge to ideal-matched curves as $T\to\infty$. It is not true,
however, that every ideal-matched curve arises this way. For instance, the split
curves that normally occur at $e\infty$ for a simple ideal-matched
curve have moduli, so the moduli spaces of simple ideal-matched curves
of expected dimension $0$ are not, in fact, $0$-dimensional. The
philosophically correct solution is probably to remember some
microscopic information at the curves in $e\infty$ (e.g., maps to the
upper half-plane rather than $[0,1]\times\RR$). For our purposes,
however, it is enough to simply trim off the components at $e\infty$:

\begin{definition}
  \label{def:TrimmedSimpleIdealMatchedCurve}
  A \emph{trimmed simple ideal-matched
    curve}\index{curve!ideal-matched!trimmed simple}
  is 
  a pair of holomorphic combs $(w_1,w_2)$ connecting two compatible pairs of
  generators $\x_i,\y_i\in\Gen(\HD_i)$ ($i=1,2$)
  such that either:
  \begin{enumerate}[label=(T-\Alph*),ref=T-\Alph*,leftmargin=*]
  \item\label{TSIC:provincial}\index{(T-A)--(T-B)} One of $w_1$ or $w_2$ is trivial and the other is a
    index $1$ holomorphic curve with no $e$ punctures or
  \item\label{TSIC:cosmopolitan} $(w_1,w_2)$ has the following properties:
  \end{enumerate}
  \begin{enumerate}[label=\hspace{1em}(T-B\arabic*),ref=T-B\arabic*,leftmargin=*]
  \item 
    \label{TSIC:WestCurve}\index{(T-B1)--(T-B7)}
    The comb $w_1$ is a holomorphic curve for $\HD_1$ (height $1$ building with no
    components at east infinity)
    which is asymptotic to the non-trivial sequence of non-empty sets of Reeb chords
    $\vec{\rhos}=(\rhos_1,\dots,\rhos_n)$.
  \item 
    \label{TSIC:WestIndex}
      The curve $w_1$ has index $1$ with respect to $\vec{\rhos}$.
  \item
    \label{TSIC:EastBrokenFlow}
    The comb $w_2$ is a height $n$ holomorphic building for $\HD_2$
    with no components at east infinity.
  \item
    \label{TSIC:EastIndex}
    Each story of $w_2$ has index one.
  \item \label{TSIC:monotone} Each of $w_1$ and $w_2$ is strongly boundary monotone.
  \item
    \label{TSIC:EastPunctures}
    For each $i=1,\dots,n$, the east punctures of the $i\th$
    story of $w_2$ are labeled, in order, by a non-empty sequence
    of Reeb chords
    $(-\rho^i_1,\dots,-\rho^i_{\ell_i})$,
    which have the property that
    the sequence of singleton sets of chords
    $\vec{\rho}\,^i=(\{\rho_1^i\},\dots,\{\rho_{\ell_i}^i\})$ is
    composable. We will let $\vec{\rho}~$ be the concatenation of the
    sequences $\vec{\rho}\,^1,\dots,\vec{\rho}\,^n$.
  \item
    \label{TSIC:WestPunctures}
    The composition of the sequence of singleton sets of Reeb chords
    $\vec{\rho}\,^i$ on the
    $i\th$ story of $w_2$ (with reversed orientation)
    coincides with the $i\th$ set of Reeb chords $\rhos_i$,
    in the partition for $w_1$; i.e.,
    \[
    \rhos_i=\biguplus_{j=1}^{\ell_i} \{\rho^i_j\}.
    \]
  \end{enumerate}
\end{definition}

\begin{lemma}
  \label{lem:ChangeTSIC}
  When Conditions~\eqref{TSIC:WestCurve} through \eqref{TSIC:monotone} hold,
  then Conditions~\eqref{TSIC:EastPunctures} and~\eqref{TSIC:WestPunctures}
  hold if and only if the following holds:
{\rm \begin{enumerate}[label=\hspace{1em}(T-B6$^{\prime}$),ref=T-B6$^\prime$,leftmargin=*]
  \item
        \label{TSICp}\index{(T-B6')}
  For all $i=1,\dots,n$, the east punctures of the $i\th$ story of $w_2$ are
  labeled, in order, by a non-empty sequence of Reeb chords $(-\rho_1^i,\dots -\rho_{\ell_i}^i)$
  which has the property that
  \begin{equation}
    \label{eq:ProductNonZero}
    I(o(\x_1,\vec{\rhos}_{[1,i-1]}))
    a(\rhos_i)=I(o(\x_1,\vec{\rhos}_{[1,i-1]})) \prod_{j=1}^{\ell_i}
    a(\rho^i_j)\neq 0.
  \end{equation}
  Here, $o(\x_1,\vec{\rhos}_{[1,i-1]})$ denotes the set of
  $\alpha$-arcs occupied by $w_1$ between $\rhos_{i-1}$ and $\rhos_i$;
  see Definition~\ref{def:strong-monotonicity-P}.
  \end{enumerate}}
\end{lemma}
\begin{proof}
  Suppose that $(w_1,w_2)$ is a trimmed simple ideal-matched comb.
  Let $\y^i$ denote the initial point of the $i\th$ story of $w_2$
  (so that $\y^1=\x_2$).
  
  \begin{claims}
  \claim\label{claim:TSIC-initial-disjoint} \emph{The sets $o(\x_1,\vec{\rhos}_{[1,i-1]})$ and
  $o(\y^i)$ are disjoint.}  In the case where $i=1$, this is the fact
  that $\x_1$ and $\y_1$ are a compatible pair of generators.  The
  case where $i>1$ follows from Conclusion~(\ref{conc:DisjointSets})
  of Lemma~\ref{lem:BiMonotonicity}, applied to $\vec\rhos$, using
  Condition~(\ref{TSIC:WestPunctures}).

  Let $o_i=o(\x_1,\vec{\rhos}_{[1,i-1]})$.

  \claim\label{claim:TSIC-product} \emph{For each $i=1,\dots,n$,
  \begin{equation}\label{eq:want-this-for-a-not-a-naught}
  I(o_i) a_0(\rhos_i)=I(o_i)\prod_{j=1}^{\ell_i} a_0(\rho^i_j)\neq 0.
  \end{equation}}
  This follows from Claim~\ref{claim:TSIC-initial-disjoint}, the
  composability of the~$\vec{\rho}\,^i$ (Condition~(\ref{TSIC:EastPunctures})) and
  Lemma~\ref{lem:reeb-product}.

  We wish to conclude that
  Equation~\eqref{eq:want-this-for-a-not-a-naught} holds with $a$ in
  place of $a_0$.

  For the next three claims, pick some $p,q\in \CircPts$ with
  $M(p)=M(q)$ and $p\neq q$. Let $\{\rho^i_{j_1},\dots,\rho^i_{j_m}\}$
  be the subsequence of chords on the $i\th$ level which start or end
  in $\{p,q\}$.

  \claim\label{claim:TSIC-no-period} \emph{There is no $\rho^i_j$ such that
  $(\rho^i_j)^-=p$ and $(\rho^i_j)^+=q$.} We argue this as follows. Let
  $\rho^i_{j_k}$ be the first $\rho^i_j$ such that
  $(\rho^i_{j_k})^-=p$ and $(\rho^i_{j_k})^+=q$.
  By strong boundary monotonicity of~$w_2$ (Condition~(\ref{TSIC:monotone})), for $\ell=1,\dots, k-1$,
  the chords $\rho^i_{j_\ell}$ alternate between starting at one of
  $\{p,q\}$ and terminating at one of $\{p,q\}$. Similarly, since
  $\rho^i_{j_k}$ terminates at $q$, strong boundary
  monotonicity of $w_2$ implies that $\rho^i_{j_{k-1}}$ starts at one of $\{p,q\}$.
  Now, there are two cases, according to the parity of~$k$. If $k$ is even,
  then there are at least two strands starting in $\{p, q\}$ in $\biguplus_{j=1}^{j_k} \{\rho^i_j\}$,
  and hence at least two strands starting at $p$ and $q$ in $\rhos_i$,
  violating boundary monotonicity of $w_1$
  (Condition~(\ref{TSIC:monotone})). If $k$ is odd, there is at least one
  strand starting in $\{p, q\}$ in $\biguplus_{j=1}^{j_k} \{\rho^i_j\}$, and hence
  $\{p,q\}$ appears in $o_i$.
  But in this case, $\rho^i_{j_1}$~terminates at one of $\{p,q\}$, and hence $\{p,q\}$ appears in 
  $o(\y^i)$. This violates Claim~\ref{claim:TSIC-initial-disjoint}.

  \claim\label{claim:TSIC-alternate} \emph{For $k=1,\dots,m$, $\rho^i_{j_k}$ alternates between
  starting at one of $\{p,q\}$
  and ending at one of $\{p,q\}$.} This follows immediately from
  boundary monotonicity of $w_2$, together with Claim~\ref{claim:TSIC-no-period}.

  \claim\label{claim:TSIC-match} \emph{There is no $k$ such that
  $(\rho^i_{j_k})^+$ and $(\rho^i_{j_{k+1}})^-$ are both in $\{p,q\}$
  but not equal to each other.}
  Suppose that there is such a~$k$.
  From Claim~\ref{claim:TSIC-alternate}, it follows that if $k$ is
  odd then $\rhos_i$ has a chord starting at one of $\{p, q\}$, and
  $\rho^i_{j_1}$ is a chord terminating at one of $\{p, q\}$; hence $\{p,q\}$ is
  in the initial idempotent for both $w_1$ and $w_2$, contradicting 
  Conclusion~(\ref{conc:DisjointSets})
  from Proposition~\ref{lem:BiMonotonicity}. If $k$ is even, there are at least
  two strands leaving one of $p$ or~$q$ in~$\rhos_i$, contradicting
  boundary monotonicity of~$w_1$.
  \end{claims}

  Claim~\ref{claim:TSIC-alternate} ensures that for $k=1,\dots,m$, $\rho^i_{j_k}$ and
  $\rho^i_{j_{k+1}}$ do not have matched initial or terminal 
  points. Together with
  Claim~\ref{claim:TSIC-match} (using a straightforward inductive application of
  Lemma~\ref{lem:VanishingProduct}), this ensures that not only is
  $I(o_i)\prod_{j=1}^{\ell_i} a_0(\rho^i_j) \neq 0$, but in fact
  $I(o_i)\prod_{j=1}^{\ell_i} a(\rho^i_j) \neq 0$.

  Conversely, suppose that Conditions~\eqref{TSIC:WestCurve}--\eqref{TSIC:monotone} hold,
  along with Condition~\eqref{TSICp}.
  Then, the non-vanishing of
  $I(o_i)\prod_{j=1}^{\ell_i} a(\rho^i_j)$ implies the non-vanishing of 
  $I(o_i)\prod_{j=1}^{\ell_i} a_0(\rho^i_j)$, which, by Lemma~\ref{lem:reeb-product},
  ensures both the composability of ${\vec\rho}\,^i$ (i.e., Condition~\eqref{TSIC:EastPunctures})
  and the property that
  $$I(o_i)\prod_{j=1}^{\ell}a_0(\rho^i_j)=I(o_i)a_0\left(\biguplus_{j=1}^{\ell_i} \{\rho^i_j\}\right).$$
  Equation~\eqref{eq:ProductNonZero} now gives
  Conclusion~\eqref{TSIC:WestPunctures}.  
\end{proof}

Let 
$\gls*{tsicMod}$
denote the moduli space of trimmed simple ideal-matched curves
connecting $\x_1,\x_2$ to $\y_1,\y_2$ with homology class $B$ and east/west 
sources $\SourceSub{1}$ and $\SourceSub{2}$, modulo
$\RR$-translation of $w_1$ and each story of $w_2$.\index{$\RR$-action!on matched pairs} Let
\[
\gls*{tsicModEmb}
\coloneqq\!\!\!\bigcup_{
  \substack{
    \chi(\SourceSub{1})=\chi_\emb(B_1,\vec{\rhos})\\
    \chi(\SourceSub{2})=\chi_\emb(B_2,-\vec{\rho})}
}\!\!\!
\cMM^B_\tsic(\x_1,\y_1,\SourceSub{1};\x_2,\y_2;\SourceSub{2}).
\]

\begin{lemma}\label{lem:trimming-gives-trimmed}
  The spine\index{spine!of ideal-matched curve} of a simple
  ideal-matched curve 
  is a trimmed simple
  ideal-matched curve. That is, if $(U_1,U_2)$ is a simple
  ideal-matched curve, where $U_1=(w_1,v_1)$ (with $v_1$ a generalized
  split curve) and $U_2$ is toothless comb~$w_2$\index{toothless comb},
  then $(w_1,w_2)$ is an embedded trimmed simple ideal-matched curve.
\end{lemma}
\begin{proof}
  Case~(\ref{TSIC:provincial}) is trivial. In case~(\ref{TSIC:cosmopolitan}),
  properties~(\ref{TSIC:WestCurve}),
  (\ref{TSIC:WestIndex}), (\ref{TSIC:EastBrokenFlow}),
  (\ref{TSIC:EastIndex}), and (\ref{TSIC:monotone}) are all immediate from the fact that $(U_1,U_2)$
  is a simple ideal-matched curve. 

  Properties~(\ref{TSIC:EastPunctures}) and~(\ref{TSIC:WestPunctures})
  are verified as follows. 

  First, weak composability follows from
  Property~(\ref{TSIC:monotone}) and
  Lemma~\ref{lem:MonotoneChords}.

  Next, let $-\vec{\rho}$ denote the sequence of Reeb chords coming from
  $w_2$, and $\vec{\rhos}$ be the sequence of sets of Reeb chords on $w_1$.
  Then $\vec{\rhos}$ is the partition of Reeb chords at west infinity
  for the generalized split curve $v_1$, while $\vec{\rho}$ is the
  collection of Reeb chords at east infinity of $v_1$, with a discrete
  partition and an ordering compatible with $\vec{\rhos}$. It follows
  that $\vec{\rhos}$ is obtained from $\vec{\rho}$ by contracting
  various consecutive sequences of Reeb chords.

  The stronger condition of composability follows from the following
  index considerations (akin to
  Lemma~\ref{lemma:collision-is-composable}).
  By Lemma~\ref{lem:composable-gr}, we have an inequality:
  $$
  \prod_{i=1}^{|\vec{\rhos}|} \grb(a(\rhos_i)) \leq
  \prod_{i=1}^{|\vec{\rho}|} \grb(a(\rho_i)),
  $$
  with strict inequality if and only if some sequence in $\vec{\rho}$
  which is contracted to get an element $\rhos_i$ is not
  composable. Equivalently, by Lemma~\ref{lem:iota-grading}, if $\iota$
  is the function from Equation~\eqref{eq:def-iota} then
  $\iota(\vec{\rhos})\leq \iota(\vec{\rho})$ with strict inequality if
  and only if there is a non-composable sequence which is contracted.

  It follows from Definition~\ref{def:emb-ind-emb-chi} that:
  \begin{align*} 
    \ind(B_1,\vec{\rhos})
    &= \ind(B_1,\vec{\rho})+|\vec{\rhos}|-|\vec{\rho}| + (\iota(\vec{\rhos})-\iota(\vec{\rho}))\\
    &= \ind(B_1,\vec{\rho})+\ind(B_2,\vec{\rho})-|\vec{\rho}|+ (\iota(\vec{\rhos})-\iota(\vec{\rho}))+(|\vec{\rhos}|-\ind(B_2,\vec{\rho}))\\
    &= \ind(B_1\glue B_2) +  (\iota(\vec{\rhos})-\iota(\vec{\rho}))+(|\vec{\rhos}|-\ind(B_2,\vec{\rho})).
  \end{align*}
  By Lemma~\ref{lem:sim-ind-1}, $\ind(B_1\glue B_2)=1$, and we just verified that
  $(\iota(\vec{\rhos})-\iota(\vec{\rho}))\leq 0$. Since $|\vec\rhos|$
  is equal to the number of stories of $U_2$ and each story has index $1$,
  $(|\vec{\rhos}|-\ind(B_2,\vec{\rho}))=0$. Since there is a
  holomorphic curve, $\ind(B_1,\vec{\rhos})\geq 1$ by
  Proposition~\ref{prop:asympt_gives_chi}, so
  $\iota(\vec{\rhos})=\iota(\vec{\rho})$ and
  $|\vec{\rhos}|=\ind(B_2,\vec{\rho})$. It follows that $w_1$ and
  $w_2$ are embedded and the Reeb chords in $\vec{\rho}$ which are
  collapsed in $\vec{\rhos}$ are composable. 
\end{proof}

It turns out that there is a combinatorially unique way to un-trim
curves. To this end, it helps to have the following:

\begin{definition}
  Let $\vec\rho=(\rho_1,\dots,\rho_n)$ be a sequence of Reeb chords.
  A subsequence $(\rho_{j_1},\dots,\rho_{j_m})$ is called
  \emph{abutting} if $\rho_{j_k}^+=\rho_{j_{k+1}}^-$ for each
  $k=1,\dots,m-1$.  An \emph{amalgamation} of $\vec\rho$ is a
  partition of $\vec\rho$ into abutting sub-sequences. Given an
  amalgamation of $\vec\rho$, there is an associated set of Reeb
  chords $\rhos$ whose elements are formed by joining together
  the chords in the abutting sub-sequences specified by the amalgamation. We call
  such a set $\rhos$ an \emph{amalgam} of~$\vec\rho$.
  \index{amalgamation of chords}\index{Reeb chords!amalgamation}
  \index{abutting!sequence of Reeb chords}\index{Reeb chords!abutting}
\end{definition}

\begin{lemma}\label{lem:amalgam-unique}
  If $\vec\rho=(\rho_1,\dots,\rho_n)$ is a sequence of Reeb chords with the property that
  $a_0(\rho_1)\cdots a_0(\rho_n) = a_0(\rhos) \neq 0$, then there is a unique
  amalgamation of~$\vec{\rho}$ with resulting amalgam~$\rhos$.
\end{lemma}
\begin{proof}
  For existence of the amalgamation, let $j_1=1$, let $j_2$ be the first index
  after $j_1$ so that $\rho_{j_2}^-=\rho_{j_1}^+$, and so on. This
  gives an abutting subsequence
  $(\rho_{j_1},\dots,\rho_{j_m})$.  Since $\rhos = \biguplus_j
  \{\rho_j\}$, one of the chords in $\rhos$ is $\rho_{j_1}\uplus\dots\uplus\rho_{j_m}$.  We
  take $(\rho_{j_1},\dots,\rho_{j_m})$ to be one part of the desired
  amalgamation.  The other parts of the amalgamation
  are gotten by iteratively repeating this process with~$\vec{\rho}$ replaced by $\vec{\rho}\setminus
  (\rho_{j_1},\dots,\rho_{j_m}).$  We call the resulting amalgamation the
  \emph{maximal amalgamation} of~$\vec\rho$; by construction the
  corresponding amalgam is $\rhos$.

  For uniqueness, we argue as follows.  By the non-vanishing of $a_0(\rho_1)\cdots
  a_0(\rho_n)$, given $p\in\PMC$, the subsequence of Reeb chords in
  $\vec\rho$ which start or end at $p$ alternates between starting
  at $p$ and ending at $p$.  In the maximal amalgamation, a chord ending
  at~$p$ always is joined to the next chord starting at~$p$. For an
  amalgamation which does not always join these pairs, the number of
  chords starting or ending at~$p$ is greater than for the maximal
  amalgamation. It follows that the only amalgamation whose
  amalgam is~$\rhos$ is the maximal amalgamation: any other
  amalgamation results in an amalgam with more chords in it.
\end{proof}

\begin{lemma}\label{lem:grow-stubble}
  Let $(w_1,w_2)$ be a trimmed simple ideal-matched curve. Then there
  is a generalized split curve $v\co\biSource\to
  \bdy\bSigma\times\RR\times[0,1]\times\RR$ such that $(w_1,v,w_2)$ is
  a simple ideal-matched curve. Moreover, the combinatorial type of
  $\biSource$ is determined by $(w_1,w_2)$.
\end{lemma}
\begin{proof}
  Again Case~\eqref{TSIC:provincial} is trivial, so we consider
  Case~\eqref{TSIC:cosmopolitan}.
  Observe that each generalized split curve~$v$ with
  $e$~punctures labeled by~$\vec{\rho}$ gives an amalgamation
  of~$\vec{\rho}$, depending only on the combinatorial type of the
  source~$\biSource$: the parts of the amalgamation correspond to the
  generalized split components making up~$v$.  Furthermore, the
  $w$~punctures of~$v$
  are labeled by the corresponding amalgam of~$\vec{\rho}$.  In
  addition, in any simple ideal-matched curve, $v$ must match the
  $i\th$ level of~$w_1$ with the $i\th$ story of~$w_2$.  So it suffices to
  show that for a trimmed simple ideal-matched curve with $w_1$ asymptotic to
  $\vec{\rhos}$ and the $i\th$ level of $w_2$ asymptotic to
  $-\vec\rho\,^i$, there is a unique amalgamation
  of~$\vec\rho\,^i$ whose amalgam is~$\rhos_i$.  But this follows
  from Property~\eqref{TSICp} and Lemma~\ref{lem:amalgam-unique}.
\end{proof}

\begin{lemma}\label{lem:TSIC-finite}
  For a generic choice of almost complex structures $J_1$ and $J_2$
  and a homology class $B=B_1\glue B_2$, the moduli space
  $\cMM^B_\tsic(\x_1,\y_1,\SourceSub{1};\x_2,\y_2;\SourceSub{2})$
  consists of a finite number of points.
\end{lemma}
\begin{proof}
  Consider a curve $(w_1,w_2)\in
  \cMM^B_\tsic(\x_1,\y_1,\SourceSub{1};\x_2,\y_2;\SourceSub{2})$.
  There are only finitely many sequences of sets of Reeb chords
  $\vec\rhos$ compatible with $B_1$. Thus, by
  Condition~\eqref{TSIC:WestIndex}, $w_1$ lies in one of finitely many
  index $1$ moduli spaces which, by
  Propositions~\ref{prop:transversality} and~\ref{prop:compactness}
  and Theorem~\ref{thm:master_equation} (transversality, compactness,
  and the analysis of the boundary strata), is a compact $0$-manifold.
  Similarly, there are a finite number of ways to decompose $B_2$ as a
  sum $B_2=B_{2,1}+\dots+B_{2,n}$ with each $B_{2,i}$ a positive
  domain, and a finite number of sequences of sets of Reeb chords
  compatible with each $B_{2,j}$. Thus, each story of $w_2$ lies in
  one of finitely many index $1$ moduli spaces, each of which is a
  compact $0$-manifold.
\end{proof}

\begin{proposition}
  \label{prop:LargeTLimit}
  For $\x_i,\y_i\in\Gen(\HD_i)$;
  $B_1$, $B_2$ so that $\ind(B_1\glue B_2)=1$; generic $J_1$ and~$J_2$;
  and generic $T$ sufficiently large,
  \[
  \#\cMM^{B_1\glue B_2}(T\semico
  \x_1,\y_1\semico\x_2,\y_2)=\#\cMM^{B_1\glue B_2}_\tsic(\x_1,\y_1;\x_2,\y_2).
  \]
\end{proposition}
\begin{proof}
  Let $B=B_1\glue B_2$. By Lemma~\ref{lem:TSIC-finite}, the moduli space
  $
  \cMM^{B}_\tsic(\x_1,\y_1\semico\x_2,\y_2)
  $
  consists of finitely many points, which we label
  $(w_1^i,w_2^i).$
  Let $(w_1^i,v^i,w_2^i)$ be an un-trimming of $(w_1^i,w_2^i)$, as
  given by Lemma~\ref{lem:grow-stubble}.
 
  Consider the union of moduli spaces
  $$
  \gls*{geqTMatchMod}
  \coloneqq
  \bigcup_{T\geq T_0} (T,\cMM^B(T,\x_1,\y_1,\SourceSub{1};
  \x_2,\y_2,\SourceSub{2})).
  $$
  Our goal is to find disjoint
  neighborhoods $\mathcal{W}^i$ of each $(w_1^i,v^i,w_2^i)$ in
  $\ocMM^{B_1\glue B_2}(\geq T_0\semico\x_1,\y_1\semico\x_2,\y_2)$ so that, for $T$
  sufficiently large, $\cMM^{B}(T\semico
  \x_1,\y_1\semico\x_2,\y_2)\subset \bigcup_i\mathcal{W}^i$ and, for
  each~$i$, $\cMM^{B}(T\semico
  \x_1,\y_1\semico\x_2,\y_2)$ intersects $\mathcal{W}^i$ in an odd
  number of points. The neighborhoods $\mathcal{W}^i$ will have the
  form
  \[
  \mathcal{W}^i=(\cU^i_1\times\cU^i_2)\cap \ocMM^{B_1\glue B_2}(\geq
  T_0;\x_1,\y_1;\x_2,\y_2)
  \]
  where $\cU^i_1\subset \ocM^{B_1}(\x_1,\y_1;\SourceSub{i,1}\glue
  \biSource_i)$ is a smeared neighborhood
  (Definition~\ref{def:smeared-nhbd}) of $(w_1^i,v^i)$ and
  $\cU^i_2\subset \ocM^{B_2}(\x_2,\y_2;\SourceSub{i,2})$ is a
  neighborhood of $w_2^i$.

  The proof involves two steps, of course: compactness and gluing.
  We start with the easier step, compactness. Fix any collection
  of neighborhoods $\mathcal{W}^i$ as above. By
  Proposition~\ref{prop:ideal-compactness}, any sequence of
  holomorphic curves in $\cMM^{B}(T_j;\x_1,\y_1;\x_2,\y_2)$
  with $T_j\to\infty$ converges to an ideal-matched holomorphic comb,
  which must be a simple
  ideal-matched holomorphic curve, by Lemma~\ref{lem:IdealMatchedTransversality}. By
  Lemma~\ref{lem:trimming-gives-trimmed} and the uniqueness part of
  Lemma~\ref{lem:grow-stubble}, this limit lies in some
  $\mathcal{W}^i$. It follows that
  $\cMM^{B}(T\semico\x_1,\y_1\semico\x_2,\y_2)\subset
  \bigcup_i\mathcal{W}^i$ for $T$ sufficiently large.

  We turn next to the gluing step. Fix a trimmed simple ideal-matched
  curve $(w_1,w_2)\in \cMM^{B}_\tsic(\x_1,\y_1;\x_2,\y_2)$, and let
  $(w_1,v,w_2)$ be its un-trimming. Let $\SourceSub{i}$
  denote the source of $w_i$ and $\biSource$ the source of $v$. Let
  $\vec{P}=(p_1,\dots,p_n)$ be the sequence of $e$~punctures of~$w_2$,
  ordered by the story in which they occur and, within each story, by
  their $\RR$-coordinates. Via the correspondence $\phi$ (from
  Definition~\ref{def:ideal-matched-comb}), we can also view $\vec{P}$ as a
  list of the $e$ punctures of $v$.
  We can also group the $p_i$
  according to the stories of $w_2$ or, equivalently, according to
  their heights in $v$; write this grouping as
  \[
  \vec{Q}=(\overbrace{\{p_1,\dots,p_{\ell_1}\}}^{Q_1},\overbrace{\{p_{\ell_1+1},\dots,p_{\ell_2}\}}^{Q_2},\dots,\overbrace{\{p_{\ell_{m-1}+1},\dots,p_{\ell_m}\}}^{Q_m}).
  \]
  Thus, $n$ is the total number of $e$ punctures of
  $w_2$, $m$ is the number of stories of $w_2$ and $\ell_i - \ell_{i-1}$ is the
  number of $e$ punctures on the $i\th$ story of $m$, so
  $n=\ell_m$.

  Observe that, by the proof of Lemma~\ref{lem:grow-stubble}, the
  ordering of the elements of $Q_i$ induced by $P$ is compatible with
  the induced ordering around the boundary of each generalized
  split-component of $v$.

  View $(w_1,v)$ as an element of
  $\ocM^{B_1}(\x_1,\y_1\semico \SourceSub{1}\glue\biSource\semico \vec{P})$ and
  $w_2$ as an element of
  $\ocM^{B_2}(\x_2,\y_2\semico \SourceSub{2}\semico \vec{P})$. There are evaluation
  maps
  \begin{align*}
    \ev^1&\co
    \cM^{B_1}(\x_1,\y_1;\SourceSub{1}\glue\biSource;\vec{P})\to
    \RR_+^{n-1} &
    \ev^2&\co \cM^{B_2}(\x_2,\y_2;\SourceSub{2};\vec{P})\to \RR_+^{n-1}
  \end{align*}
  given by 
  \[
  \ev^i(u_i)=(\ev_{p_{2}}(u_i)-\ev_{p_{1}}(u_i),\ev_{p_{3}}(u_i)-\ev_{p_{2}}(u_i),\dots, \ev_{p_{n}}(u_i)-\ev_{p_{n-1}}(u_i)).
  \]
  Define $\pi_1\co \RR^{n-1}\to \RR^{n-m}$ and $\pi_2\co
  \RR^{n-1}\to\RR^{m-1}$ by
  \begin{align*}
    \pi_1(t_1,\dots,t_n)&= (t_1,\dots,t_{\ell_1-1},t_{\ell_1+1},\dots,t_{\ell_2-1},t_{\ell_2+1},\dots\dots,t_{\ell_m-1},t_{\ell_m})\\
    \pi_2(t_1,\dots,t_n)&= (t_{\ell_1},t_{\ell_2},\dots,t_{\ell_{m-1}}).
  \end{align*}

  We can find a smeared neighborhood $\cU_1$ of $(w_1,v)$ and a
  neighborhood $\cU_2$ of $w_2$ with the following properties:
  \begin{enumerate}[label=(i-\arabic*),ref=(i-\arabic*),leftmargin=*]
  \item
    \label{item:pi2zero}\index{(i-1)--(i-6)}
    The map $\pi_2\circ\ev^1$ sends
    $\cU_1$  to a neighborhood of some
    vector $y_0$ in the interior of $\RR^{m-1}_+$ (disjoint from $\bdy\RR^{m-1}_+$). This
    uses continuity of the relative evaluation maps on the
    compactified moduli space. The vector $y_0$ records the height
    differences of the punctures on $w_1$. 
  \item 
    \label{item:pi1zero}
    The map $\pi_1\circ \ev^2$ sends $\cU_2$ to a neighborhood of some
    vector $x_0$ in the interior of $\RR^{n-m}_+$. Again, this uses
    continuity of the evaluation maps: the vector $x_0$ records the
    height differences of the punctures within each of the stories of
    $w_2$.
  \item 
    \label{item:pi1one}
    For all sufficiently
    small~$\delta$, $\pi_1\circ \ev^1$ maps $\cU_1$ properly and with
    odd degree to $\{x\in\RR^{n-m}_+\mid \|x\|<\delta\}$. This uses Proposition~\ref{prop:generalized_gluing_degree_one},
    and our observation about the ordering of the points in
    each~$Q_i$ (to see that the image of the punctures on each
    generalized
    split component lands in $\RR^{n-m}_+\subset\RR^{n-m}$).
  \item 
    \label{item:pi2one}
    For all sufficiently
      large~$R$, $\pi_2\circ \ev^2$ maps $\cU_2$ properly and with odd degree onto
    $\{x\in\RR^{m-1}_+\mid \|x\|>R\}\cup\{\infty\}$ (so that
    $w_2$ is mapped to $\infty$). 
    This uses
    Proposition~\ref{prop:gluing_two_story} (or rather, a height $n$
    analogue thereof).
  \item 
    \label{item:disjointU1}
    The closure of $\cU_1$ intersects $\{(w_1^i,v^i)\}$ in the one point $(w_1,v)$.
  \item
    \label{item:disjointU2}
    The closure of $\cU_2$ intersects $\{w_2^i\}$ in the one point $w_2$.
  \end{enumerate}
  It follows from these observations and an intersection theory
  argument that 
  for any generic, large choice of $T$, the fibered product
  \begin{equation}\label{eq:odd-fiber-prod}
  \cU_1\times\cU_2\cap \ocMM^{B_1\glue B_2}(T;\x_1,\y_1;\x_2,\y_2) =
  \cU_1\sos{T\ev^1}{\times}{\ev^2} \cU_2
  \end{equation}
  consists of an odd number of points. We spell this out presently.
  First, we rephrase Conditions~\ref{item:pi2zero}
  and~\ref{item:pi1zero} more concretely:
  \begin{enumerate}[label=(i-\arabic*$^\prime$),ref=(i-\arabic*$^\prime$),leftmargin=*]
  \item\index{(i-1')--(i-2')}
    \label{item:pi2zeroP}
    There is a $y_0\in\RR^{m-1}_+$ and  $\epsilon>0$ so that
    $\pi_2\circ\ev^1$ sends $\cU_1$ into $B_2\coloneqq B_{\epsilon}(y_0)\subset \RR^{m-1}_+$.
  \item 
    \label{item:pi1zeroP}
    There is an $x_0\in\RR^{n-m}_+$ and $\eta>0$ so that 
    $\pi_1\circ\ev^2$ sends $\cU_2$ into $B_1\coloneqq B_\eta(x_0)\subset \RR^{n-m}_+$.
  \end{enumerate}
  See Figure~\ref{fig:FiberProd}.

  \begin{figure}
    \centering
    \includegraphics[scale=.83333]{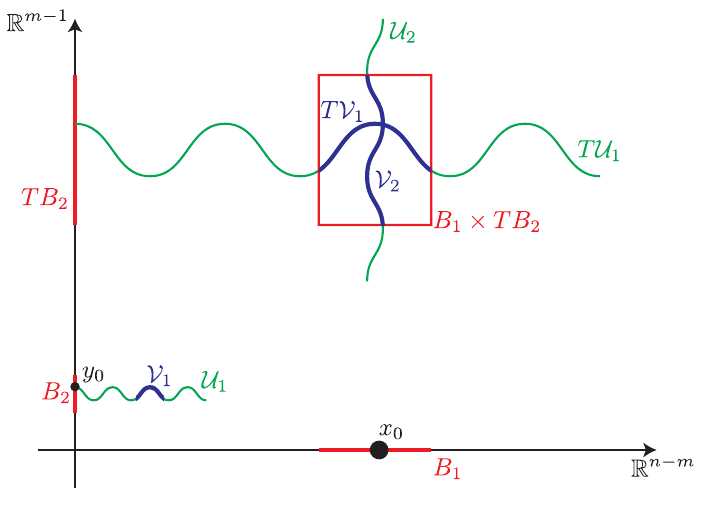}
    \caption[Illustration of the fibered product in the proof                                                            
        of Proposition~\ref{prop:LargeTLimit}]{\textbf{Illustration of the fibered product in the proof
        of Proposition~\ref{prop:LargeTLimit}.}}
    \label{fig:FiberProd}
  \end{figure}
  \colorused

  Think of the fibered product $\cU_1\sos{T\ev^1}{\times}{\ev^2}
  \cU_2$ as a subset of $\RR^{n-m}_+\times \RR^{m-1}_+$. By the above
  two conditions, it follows that the fibered product is
  contained in the subset $B_1 \times TB_2$.
  Choose $T$ large enough so that $T^{-1}B_1$ is contained inside the ball of
  radius $\delta$ about the origin in $\RR_+^{n-m}$, and so
  that $TB_2$ lies outside a ball of radius~$R$ around the
  origin in $\RR_+^{m-1}$, where $\delta$ and~$R$ are constants from
  Conditions~\ref{item:pi1one} and~\ref{item:pi2one}, respectively.  Then let $\overline{\cV}_1$
  denote the preimage under $T\ev^1$ of ${\overline B}_1\times
  T{\overline B}_2$, and let $\overline{\cV}_2$ denote the preimage
  under
  $\ev^2$ of ${\overline B}_1\times T{\overline B}_2$.  Both are
  compact, because $\pi_1 \circ \ev^1$ is proper on the preimage of
  $T^{-1}B_1$ and $\pi_2 \circ \ev^2$ is proper on the preimage of~$B_2$.

  Now, Condition~\ref{item:pi2zeroP} shows that $T\ev^1\co
  \overline{\cV}_1 \to
  \overline{B}_1 \times T\overline{B}_2$ can be viewed as a relative
  cycle, representing some homology class in $H_{n-m}({\overline
    B}_1\times T B_2, (\partial {\overline B}_1)\times T
  B_2;\Field)\cong \Field.$ Indeed, Condition~\ref{item:pi1one}
  ensures that $[\overline{\cV}_1]$ is the non-trivial
  class: the degree of $\pi_1 \circ T\ev^1\colon
  \overline{\cV}_1\rightarrow \overline{B}_1$ can be interpreted
  as the
  intersection pairing of $[{\overline\cV}_1]$ with an element of
  $H_{m-1}(B_1\times T{\overline B}_2, B_1\times (\partial T{\overline
    B}_2);\Field)\cong \Field.$

  Symmetrically, by Conditions~\ref{item:pi1zeroP} and~\ref{item:pi2one},
  $\ev^2\co {\overline \cV}_2 \to \overline{B}_1 \times T\overline{B}_2$
  represents a generator of $H_{m-1}(B_1\times T{\overline B}_2,
  B_1\times (\partial T{\overline B}_2);\Field)\cong \Field.$ The
  count of points in the fibered product
  $\overline{\cV}_1\sos{T\ev^1}{\times}{\ev^2} \overline{\cV}_2=
  \cU_1\sos{T\ev^1}{\times}{\ev^2} \cU_2$ can now be interpreted as
  the intersection pairing of $[{\overline\cV}_1]$ with $[{\overline
    \cV}_2]$ under the non-trivial intersection pairing
\begin{align*}
    H_{n-m}&({\overline B}_1\times T B_2,
    (\partial {\overline B}_1)\times T B_2)
    \otimes 
    H_{m-1}(B_1\times T {\overline B}_2,
    B_1\times (\partial T {\overline B}_2)) \\
    &\rightarrow H_0(B_1\times TB_2)\cong\Field,
\end{align*}
as desired.
\end{proof}

\section{Completion of the proof of the pairing
  theorem}\label{sec:pairing-proof}
We are now in a position to give our second proof of the pairing
theorem, Theorem~\ref{thm:TensorPairing}.

\begin{proof} [Proof of Theorem~\ref{thm:TensorPairing}]
  Fix bordered Heegaard diagrams $\HD_1$ and $\HD_2$ for
  $(Y_1,\phi_1\co F(\PMC)\to\bdy Y_1)$ and $(Y_2,\phi_2\co
  {-F(\PMC)}\to\bdy Y_2)$. By
  Proposition~\ref{prop:admis-achieve-maintain}, we may assume that
  both $\HD_1$ and $\HD_2$ are provincially admissible, and at
  least one of $\HD_1$ or $\HD_2$ is admissible. (We could
  assume both diagrams are admissible, but we will not need this.) By
  Lemma~\ref{lem:closed-admissible}, this implies that
  $\HD=\HD_1\cup_\bdy\HD_2$ is a Heegaard diagram for~$Y$ which is
  weakly admissible for all $\SpinC$-structures. By
  Lemmas~\ref{lem:finite-typeD} and~\ref{lem:finite-typeA}, our
  assumption also implies that at least one of $\CFAa(\HD_1)$ or
  $\CFDa(\HD_2)$ is bounded.

  By Proposition~\ref{prop:IdentifyDT}, the $\DT$-product
  $\CFAa(\HD_1)\DT\CFDa(\HD_2)$ (Definition~\ref{def:DT}) is a model
  for $\CFAa(Y_1)\DTP\CFDa(Y_2)$. By Proposition~\ref{prop:VaryTChain}
  (which uses Theorem~\ref{thm:PrimitivePairing}), $\CFa(Y)$ is
  homotopy equivalent to $\CFa(T;\HD_1,\HD_2)$ for any $T>0$. So, it
  remains to identify $\CFAa(\HD_1)\DT\CFDa(\HD_2)$ with
  $\CFa(T;\HD_1,\HD_2)$ for $T$ sufficiently large.

  On the level of $\Field$-vector spaces, the identification is
  straightforward. As a vector space, $\CFAa(\HD_1)\DT\CFDa(\HD_2)$ is
  given by $X(\HD_1)\otimes_{\Idem(\PMC)}X(\HD_2)$, which has a basis 
  $\{(\x_1,\x_2)\in\Gen(\HD_1)\times\Gen(\HD_2)\mid o(\x_1)\cap
  o(\x_2)=\emptyset\}$. This is exactly the set of generators for
  $\CFa(T;\HD_1,\HD_2)$ (and for $\CFa(\HD)$).

  Next we identify the differentials. It follows from the
  admissibility criteria that there are only finitely many homology
  classes $B$ which can contribute to the differential on
  $\CFa(T;\HD_1,\HD_2)$ (compare Lemmas~\ref{lem:finite-typeD}
  and~\ref{lem:finite-typeA}). So, for large enough~$T$,
  Proposition~\ref{prop:LargeTLimit}
  implies that the differential on $\CFa(T;\HD_1,\HD_2)$ counts
  trimmed simple ideal-matched curves.

  Let $\{a_i\}$ be the set of basic generators for the algebra
  $\Alg(\PMC)$. Define operators $D_i\co X(\HD_2)\to X(\HD_2)$ by
  $$\bdy_2(\x_2)=\sum_{i} a_i\otimes D_i(\x_2)$$
  as in Example~\ref{ex:explicitDT}.
  Then, for large $T$, the differential on $\CFa(T;\HD_1,\HD_2)$ is
  given by
  \begin{equation}
    \partial_T(\x_1\otimes \x_2) = \sum m_{n+1}(\x_1,a_{i_1},\dots
    a_{i_{n}}) (D_{i_n}\circ\dots\circ D_{i_1})(\x_2), \label{eq:pair-pf-diff}
  \end{equation}
  where the sum is taken over all finite sequences
  $a_{i_1},\dots,a_{i_n}$ of basic generators for $\Alg(\PMC)$
  (including the empty sequence, where $n=0$). To see this, consider a
  trimmed simple ideal-matched curve $(w_1,w_2)$ connecting
  $\x_1\otimes \x_2$ to $\y_1\otimes \y_2$. The holomorphic building
  $w_2$ corresponds to a sequence of differentials
  $D_{i_n}\circ\dots\circ D_{i_1}$ on $\CFDa(\HD_2)$, where $n$ is the
  number of stories of $w_2$. (This uses one direction of
  Lemma~\ref{lem:ChangeTSIC} to guarantee that the algebra elements
  coming from $w_2$ are nonzero and agree with the algebra elements
  coming from $w_1$.) The curve $w_1$ corresponds to a higher
  product $m_{n+1}(\x_1,a_{i_1},\dots,a_{i_n})$ on $\CFAa(\HD_1)$ of
  $\x_1$ with the corresponding algebra elements
  $a_{i_1},\dots,a_{i_n}$. Every non-zero term on the right of
  Equation~\eqref{eq:pair-pf-diff} arises this way. (This uses the
  other direction of Lemma~\ref{lem:ChangeTSIC}, to guarantee that the
  corresponding curves $w_1$ and $w_2$ satisfy
  conditions~(\ref{TSIC:EastPunctures}) and~(\ref{TSIC:WestPunctures})
  of Definition~\ref{def:TrimmedSimpleIdealMatchedCurve}.)
  
  Equation~\eqref{eq:pair-pf-diff} is exactly the differential on
  $\CFAa(\HD_1)\DT\CFDa(\HD_2)$ (compare
  Equation~\eqref{eq:ExplicitDT}). This concludes the proof.
\end{proof}

\section{A twisted pairing theorem}
Next, we give a version of the pairing theorem for the
twisted-coefficient variants of $\CFDa$ and $\CFAa$.  Loosely
speaking, it says that the derived tensor product of twisted bordered
invariants computes the twisted Heegaard Floer complex of the glued-up
three-manifold. Before stating the precise version, we introduce a
little algebra.

\begin{lemma}\label{lem:twisted-pairing-coefficients}
  The left action on $\tCFAa(\HD_1,\spinc_1)$ by $\Field[H_2(Y_1,\partial Y_1)]$ and the right action
  of $\tCFDa(\HD_2,\spinc_2)$ by $\Field[H_2(Y_2,\partial Y_2)]$ induce 
  a left $\Field[H_2(Y,F)]$-action on the tensor product $\tCFDa(\HD_1,\spinc_1)\DT \tCFAa(\HD_2,\spinc_2)$.
\end{lemma}

\begin{proof}
   The two actions give $\tCFDa(\HD_1,\spinc_1)\DT \tCFAa(\HD_2,\spinc_2)$ the structure
   of a bimodule over $\Field[H_2(Y_1,F)]$ and $\Field[H_2(Y_2,F)]$; since
   $\Field[H_2(Y_2,F)]$ is commutative, the right action by $\Field[Y_2,F]$ can be
   viewed as a left action, giving the product complex the structure of a left
   $\Field[H_2(Y_1,F)]\otimes \Field[H_2(Y_2,F)]\cong \Field[H_2(Y,F)]$-module.
\end{proof}

The quotient map $H_2(Y)\to H_2(Y,F)$ induces a ring map $i\colon
\Field[H_2(Y)]\to \Field[H_2(Y,F)]$. Note that $i$ is injective.
As in~\cite[Section 8]{OS04:HolDiskProperties}, the twisted Heegaard
Floer complex $\tCFa(Y)$ is a chain
complex of free $\Field[H_2(Y)]$-modules, and hence induces a chain
complex of free $\Field[H_2(Y,F)]$-modules
$i_*(\tCFa(Y))=\tCFa(Y)\otimes_{\Field[H_2(Y)]}\Field[H_2(Y,F)]$.

\begin{theorem}
\index{pairing theorem!twisted coefficients}%
\index{twisted coefficients!pairing theorem with}%
  There is a homotopy equivalence of $\Field[H_2(Y,F)]$-modules.
  $$\tCFAa(\HD_1,\s_1)\DT \tCFDa(\HD_2,\s_2)\simeq i_*\biggl(\bigoplus_{\substack{\spinc\in\SpinC(Y)\\
    \spinc|_{Y_i}=\spinc_i, i=1,2}}\!\! \tCFa(\HD,\spinc)\biggr).$$
\end{theorem}

\begin{proof}
  Lemma~\ref{lem:twisted-pairing-coefficients} implies that as modules
  the two sides agree. It remains to show that the differentials agree
  (up to homotopy). For this we can use either proof of the pairing
  theorem, and only need to check that the group-ring coefficients on
  the two sides agree. Let $\x_i$, $i=1,2$, be the chosen base
  generator for $\HD_i$ in the $\SpinC$-structure $\spinc_i$. Fix
  generators $\y=\y_1\otimes\y_2$ and $\w=\w_1\otimes\w_2$ and domains
  $B_{i}\in\pi_2(\x_i,\y_i)$. Given domains $C_i\in\pi_2(\y_i,\w_i)$,
  the group-ring coefficient of $\w_1\otimes\w_2$ in the differential
  of $e^{B_1}\x_1\otimes \x_2 e^{B_2}$ coming from $(C_1,C_2)$ is
  $e^{B_1*C_1}\x_1\otimes\x_2 e^{B_2*C_2}$. The group-ring coefficient
  of $\w$ in $\bdy(\x)$ coming from $C_1\glue C_2$ is
  $e^{(B_1\glue B_2)*(C_1\glue C_2)}$. With respect to the
  identification from Lemma~\ref{lem:twisted-pairing-coefficients},
  these coefficients agree.
\end{proof}

\begin{remark}
  We can recover $\tCFa(Y)$ from $i_*(\tCFa(Y))$ as follows. Since
  $H_1(F)$ is a free abelian group, so is the cokernel of the
  projection map $i\co H_2(Y)\to H_2(Y,F)$. So, we can choose a splitting
  $\ell\co H_2(Y,F)\to H_2(Y)$, with $\ell\circ i=\Id_{H_2(F)}$. Then
  $\tCFa(Y)\cong \ell_*(i_*(\tCFa(Y)))$.
\end{remark}

\section{An example}\label{sec:tensor-prod-eg}
\begin{figure}
\centerline{\input{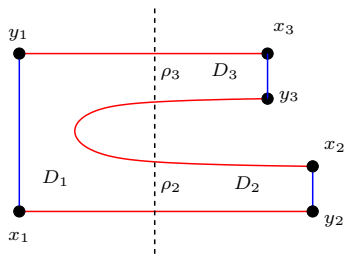}}
\caption[Degenerating a hexagon]{\label{fig:Hexagon}
  {\bf{Degenerating a hexagon.}} 
  A hexagon in the Heegaard diagram 
  (giving a flow from $\x=\{x_1,x_2,x_3\}$ to $\y=\{y_1,y_2,y_3\}$)
  is divided into three pieces $D_1$, $D_2$, and $D_3$,
  grouped as $D_1$ and $D_2\cup D_3$.
  }
\end{figure}

We conclude this chapter with a local example.
Consider a hexagon connecting the generator $\x=\{x_1,x_2,x_3,\dots\}$ to
$\y=\{y_1,y_2,y_3,\dots\}$ in a Heegaard diagram~$\HD$ for
some three-manifold~$Y$. (Points not listed are the same in all
generators, and will be dropped from the notation.) Here, the hexagon
is a domain $D$ in
$\Sigma$ whose boundary, as we traverse it, contains arcs
in three distinct $\alpha$-curves and arcs in three distinct
$\beta$-curves, so that if the arc $\partial D\cap \alpha_i$ (with its
boundary orientation) goes from $x_i$ to~$y_i$.  This domain
corresponds to an element $B\in \pi_2(\x,\y)$ with $\Mas(B)=1$. In
fact, it is not difficult to see that the domain $B$ always
contributes $\y$ to the differential of~$\x$, i.e.,
$\#\cM(\x,\y;B)=1$; see, for instance,~\cite[Lemma 9.11]{Rasmussen03:Knots}.
We wish to see this contribution from the point of
view of time dilation in a bordered setting where
the dividing curve $Z$, which separates $\HD$ into two bordered diagrams,
crosses the hexagon in two arcs.
Such a curve divides the hexagon into three regions which we denote $\{D_i\}_{i=1}^3$
with the convention that $D_i$ connects $x_i$ to $y_i$. We label so that
$D_1$ meets both arcs in $Z\cap D$, whereas $D_2$ and $D_3$ each meet only one apiece.
We label the Reeb chords corresponding to these arcs $\rho_2$ and $\rho_3$, numbered
so that $\rho_i$ meets $D_i$ for $i=2,3$.
See Figure~\ref{fig:Hexagon} for an illustration.

We start, as a warm-up, by considering the problem from the 
point of view of matched
pairs (i.e., using the complex $\CFa(\HD_1,\HD_2)$ from 
Theorem~\ref{thm:PrimitivePairing}). We verify that
\[
\#\bigl(\cMM^B(\{x_1\},\{y_1\};\{x_2,x_3\},\{y_2,y_3\})\bigr)=1.
\]
Now, $\cM(\{x_1\},\{y_1\};D_1)$ consists of a single curve $u$. Then
$\ev_{\rho_3,\rho_2}(u)=\ev_{\rho_3}(u)-\ev_{\rho_2}(u)$ is some fixed
real number $t_0$, the height difference between the Reeb chords
$\rho_3$ and $\rho_2$. The number $t_0$ is necessarily positive, but
the precise value of $t_0$ depends on the conformal parameters in the
picture.

The moduli space on the other side,
$\cM(\{x_2,x_3\},\{y_2,y_3\};D_2+D_3)$ consists of a pair of disjoint
disks, which can be moved relative to each other subject to the
constraint that $\rho_2$ is below $\rho_3$. (The constraint comes from
strong boundary monotonicity.) In particular, the
evaluation map $\ev_{\rho_3,\rho_2}$ sends this moduli space
homeomorphically to $\RR_{\geq 0}$.

It follows that there is a unique point in the fibered product, as claimed.

We compare this against two different computations which arise by time dilation.
The two different computations arise by thinking of $D_1$ as the part on the $A$ side or on the $D$ side.

We start by placing $D_1$ on the $D$ side. 
In this case, the type $A$ module for $D_2\cup D_3$ has four
generators $\gxx=\{x_2, x_3\}$, $\gxy=\{x_2, y_3\}$, $\gyx=\{y_2, x_3\}$, and $\gyy=\{y_2, y_3\}$ and the following operations:
\begin{align*}
  m_2(\gxx,\rho_2)&=\gyx \\
  m_2(\gxy,\rho_2)&=\gyy \\
  m_2(\gxx,\rho_3)&=\gxy \\
  m_2(\gyx,\rho_3)&=\gyy \\
  m_2(\gxx,\rho_2\rho_3)&=\gyy.
\end{align*}
The first four of these are gotten by looking the union of a bigon
with a trivial strip.  The last is a union of the two non-trivial
bigons $D_2+D_3$, and comes from the element of
$\cM(\{x_2,x_3\},\{y_2,y_3\},D_2+D_3)$ with the property that the
$\ev_{\rho_2,\rho_3}=0$.  We have already seen that there is a unique
such element.

The type $D$ module of $D_1$ has two generators $x_1$ and $y_1$, and
the bigon $D_1$ is interpreted as giving a differential
$$\partial x_1 = \rho_2 \rho_3 \cdot y_1.$$

The $\DT$~product of these two complexes has two generators, 
$\x=\gxx\otimes x_1$ and $\y=\gyy\otimes y_1$, and a  differential
single from $\x$ to $\y$, consistent our earlier computation.
Algebraically, this differential pairs the operation
$m_2(\gxx,\rho_2\cdot \rho_3)=\gyy$ with the (unbroken) differential
on the type $D$ side. In the graphical notation of
Chapter~\ref{chap:ainfinity}, this operation is given by
\[
\mathcenter{\begin{tikzpicture}
 \node at (0,0) (tl) {$\gxx$};
 \node at (2,0) (tr) {$x_1$};
 \node at (2,-1) (delta) {$\delta^1$};
 \node at (0,-2) (m) {$m_2$};
 \node at (0,-3) (bl) {$\gyy$};
 \node at (2,-3) (br) {$y_1$};
 \draw[amar] (tl) to (m);
 \draw [amar] (m) to (bl);
 \draw[dmar] (tr) to (delta);
 \draw[dmar] (delta) to (br);
 \draw[aar] (delta) to node[above,sloped] {\lab{\rho_2\rho_3}} (m);
\end{tikzpicture}}.
\]
Analytically, the moduli spaces $\cMM^{B}(T;\gxx,\gyy;x_1,y_1)$
consist of the curve $u$ in $\HD_1$ and the pair of disks in $\HD_2$
with punctures at height difference $t_0/T$. Taking $T\to\infty$, the
separation between these heights goes to $0$; this corresponds to the
stated $m_2$.

Next, place $D_1$ on the $A$ side. Now, the type $A$ module has 
two generators $x_1$ and $y_1$, and a single non-trivial operation
$m_3(x_1,\rho_2,\rho_3)=y_1$. 
The type $D$ module for $D_2\cup D_3$ is generated by the four
generators $\gxx$, $\gxy$, $\gyx$, and $\gyy$, with differentials
\begin{align*}
  \partial \gxx& =\rho_2\cdot \gyx + \rho_3\cdot \gxy \\
  \partial \gxy&= \rho_2\cdot \gyy \\
  \partial \gyx&= \rho_3\cdot \gyy.
\end{align*}

On the $\DT$~product of these two complexes, we once again get two
generators and a single differential. The generators are now
labeled $x_1\otimes \gxx$ and $y_1\otimes\gyy$.  Algebraically, the
differential in this case is gotten by pairing the $m_3$ on the type
$A$ side with the differential from $\gxx$ to $\rho_2\cdot \gyx$
followed by the differential from $\gyx$ to $\rho_3\cdot \gyy$.  In
the graphical notation of Chapter~\ref{chap:ainfinity}, this operation
is given by
\[
\mathcenter{\begin{tikzpicture}
 \node at (0,0) (tl) {$x_1$};
 \node at (2,0) (tr) {$\gxx$};
 \node at (2,-1) (delta1) {$\delta^1$};
 \node at (2,-2) (delta2) {$\delta^1$};
 \node at (0,-3) (m) {$m_3$};
 \node at (0,-4) (bl) {$y_1$};
 \node at (2,-4) (br) {$\gyy$};
 \draw[amar] (tl) to (m);
 \draw [amar] (m) to (bl);
 \draw[dmar] (tr) to (delta1);
 \draw[dmar] (delta1) to (delta2);
 \draw[dmar] (delta2) to (br);
 \draw[aar] (delta1) to[pos=0.25] node[above,sloped]{\lab{\rho_2}} (m);
 \draw[aar] (delta2) to[pos=0.25] node[above,sloped]{\lab{\rho_3}} (m);
\end{tikzpicture}}.
\]
Analytically, the moduli spaces $\cMM^{B}(T;\x_1,\y_1;\gxx,\gyy)$
consist of the curve $u$ in $\HD_1$ and the pair of disks in $\HD_2$
with punctures at height difference $T\cdot t_0$. Taking $T\to\infty$,
the separation between these heights goes to $\infty$; this
corresponds to the pair of $\delta^1$'s.


\chapter{Gradings}
\label{chap:gradings}

We now turn to gradings on the structures we have already constructed
in Chapters \ref{chap:type-d-mod}--\ref{chap:tensor-prod}:
Type~$D$ invariants $\CFDa(\HD)$; Type~$A$ invariants $\CFAa(\HD)$;
and their tensor products $\CFAa(\HD_1) \DT \CFDa(\HD_2)$.  This
extends the gradings on the algebra from
Section~\ref{sec:gradings-algebra}, and we will use the material on
the dimension of the moduli spaces $\tcM^B(\x,\y;\vec{\rhos})$ from
Section~\ref{sec:expected-dimensions}.

As described in
Definition~\ref{def:module-nc-grading}, the gradings $\grb$ on
$\CFDa(\HD,\s)$ and $\CFAa(\HD,\s)$ take
values in sets $\DBigGrSet(\HD, \s)$ and $\ABigGrSet(\HD, \s)$ (depending on the $\spin^c$
structure~$\s$) with actions of the grading
group~$\bigGroup(4k)$; a left action in the case of $\DBigGrSet(\HD,
\s)$ (Section~\ref{sec:typeD-gradings}) and a
right action in the case of $\ABigGrSet(\HD,\s)$ (Section~\ref{sec:typeA-gradings}).
We will also define refinements $\gr$ of these
grading taking values in sets with actions of the smaller group
$\smallGroup(\PtdMatchCirc)$ (Section~\ref{sec:refined-gradings-2}).  As discussed in
Section~\ref{sec:GradedPairingThm}, the tensor product takes values in the
product over $\smallGroup(\PMC)$ of these two grading sets. Before
constructing these gradings, we collect some properties of the
grading groups and relationships with the expected dimensions of
moduli spaces from earlier chapters
(Section~\ref{sec:grading-alg-sum}) and give a useful
grading on domains (Section~\ref{sec:domains}).

\section{Algebra review}
\label{sec:grading-alg-sum}

For reference, we recall briefly the relevant material from
Sections~\ref{sec:gradings-algebra} and~\ref{sec:expected-dimensions}.

In Section~\ref{sec:gradings-algebra} we introduced two grading groups:
\begin{itemize}
\item $\bigGroup(4k)$, which is a central extension of $H_1(Z',
  \CircPts)$ by $\ZZ$, and
\item $\smallGroup(\PMC)$, which is a subgroup of $\bigGroup(4k)$
  isomorphic to a central extension of $H_1(F(\PMC))$ by~$\ZZ$.
\end{itemize}
We review $\bigGroup(4k)$ here, and delay reviewing $\smallGroup(\PMC)$
until Section~\ref{sec:refined-gradings-2} where we construct the refined
gradings.

\glsadd{bigGroup}
Elements of $\bigGroup(4k)$ are written 
$\gls*{jalpha}$, where $j \in
\OneHalf\ZZ$ and $\alpha \in H_1(Z',\CircPts)$.
We impose the additional restriction that $j \equiv
\gls*{epsilonalpha} \pmod{1}$, where $\epsilon(\alpha)$ is one quarter
the number of parity changes in~$\alpha$.
The multiplication on
$\bigGroup(4k)$ is given by
\glsadd{biggpprod}
\begin{align*}
  (j_1,\alpha_1)\cdot(j_2,\alpha_2) &\coloneqq (j_1 + j_2 +
    L(\alpha_1, \alpha_2), \alpha_1 + \alpha_2)
\end{align*}
where 
$\gls*{Linkingalpha}$
is the linking of the boundaries of $\alpha_1$ and
$\alpha_2$, given by
\[
L(\alpha_1, \alpha_2) \coloneqq m(\alpha_2, \bdy\alpha_1)
\]
or, concretely by Formula~\eqref{eq:concrete-linking}.
We call $j$ the \emph{Maslov component} of $(j,\alpha)$ and $\alpha$ the
\emph{$\SpinC$ component} of $(j,\alpha)$.
\index{spinc component of grading@$\SpinC$ component of grading}%
\index{grading!Maslov component}%
\index{grading!spinc component@$\SpinC$ component}%
The distinguished central element of $\bigGroup(4k)$ is the element
$\lambda = (1,0)$; thus the differential lowers the first component of
the grading
by~$1$.

More generally, we have (Lemma~\ref{lem:G-n-products})
  \[
  (j_1,\alpha_1)\cdots(j_n,\alpha_n) = \Bigl(\sum_i j_i + \sum_{i < j}
  L(\alpha_i,\alpha_j),\sum_i \alpha_i\Bigr).
  \]

\index{grading!on $\Alg(\PMC)$|(}%
The grading $\grb(a)$ of an element $a\in\Alg(\PMC)$
with starting idempotent~$S$ is given by
\begin{align*}
  \gls*{iota}
  &\coloneqq \inv(a) - m([a], S)\\
  \gls*{grprime}(a)
  &\coloneqq (\iota(a), [a])
\end{align*}
(Definition~\ref{def:grading-algebra}), where $\inv(a)$ denotes the
number of inversions of $a$.  
For convenience, if $\vec\rho$ is a sequence of Reeb chords, let
$\gls*{grprimevecrho}$ 
be the
product $\prod_{i=1}^n\grb(a(\rho_i))$.
(If $a(\vec\rho) \ne 0$, then
$\grb(\vec\rho) = \grb(a(\vec\rho))$.)  Similarly set
$\gls*{grprimevecrhos}$ 
to be $\prod_{i=1}^n \gr'(a(\rhos_i))$.

The function $\iota$ can also be written, for
$\rhos$ a set of Reeb chords with $a(\rhos) \ne 0$, as
\[
  \iota(a(\rhos)) = -\!\sum_{\{\rho_1, \rho_2\}\subset \rhos}\!
      \abs{L(\rho_1,\rho_2)} - \frac{\abs{\rhos}}{2}
\]
(Lemma~\ref{lem:iota-chords}).  We also define
\begin{align*}
  \gls*{iotarhos}&\coloneqq \iota(a(\rhos))\\
  \gls*{iotavecrhos}
&\coloneqq\sum_{i}\iota(\rhos_i)+\sum_{i<j}
  L(\rhos_i,\rhos_j)
\end{align*}
so that $\iota(\vec{\rhos})$ is the Maslov component of
$\grb(\vec{\rhos})$
(Lemma~\ref{lem:iota-grading}).
\index{grading!on $\Alg(\PMC)$|)}

These quantities are related to the dimensions of the embedded moduli
space as follows.  For $B \in \pi_2(\x, \y)$ and $\vec{\rhos}$ a compatible
sequence of sets of Reeb chords, $\tcM^B(\x,\y;\vec{\rhos})$ has
expected dimension
\[
\gls*{indemb}= e(B)+n_\x(B)+n_\y(B)+|\vec\rhos|+\iota(\vec{\rhos})
\]
(Definition~\ref{def:emb-ind-emb-chi} and
Proposition~\ref{prop:asympt_gives_chi}).

\section{Domains}
\label{sec:domains}

Before defining a grading on modules, we define a grading on domains.

\begin{definition}\label{def:def-g-B}
  For $\x, \y \in \S(\HD)$ and $B \in \pi_2(\x,\y)$, define
  $\gb(B) \in \OneHalf\ZZ \times H_1(Z', \CircPts)$ by
  \begin{equation}
    \label{eq:Largeg}
    \gls*{gbig}
    \coloneqq (-e(B) - n_\x(B) - n_\y(B), \bdy^\bdy(B)).
  \end{equation}
\end{definition}

\begin{lemma}\label{lem:g-B-integral}
  For any $B \in \pi_2(\x,\y)$, the element $g'(B)$ is in $G'(4k)$.
\end{lemma}

A special case of Lemma~\ref{lem:g-B-integral} is that for provincial
domains~$B$, the expected dimension $e(B) + n_\x(B) + n_\y(B)$, which
is \emph{a priori} in $\OneQuart\ZZ$, is actually in~$\ZZ$.

\begin{proof}
  We must show that $-e(B) -n_\x(B) -n_\y(B) \equiv
  \epsilon(\bdy^\bdy(B)) \pmod{1}$.
  To do so, we follow \cite[Lemma
  2.17]{OS04:HolomorphicDisks} to construct a surface related
  to~$B$.  First, consider $\tB \coloneqq \ell[\Sigma] + B$ where
  $\ell$ is large enough that all
  local multiplicities of $\tB$ are non-negative.  Here,
  $\tB$ is an element of $\piBig(\x,\y)\glsadd{pitwoBig}$,
  the set of homology classes of curves connecting $\x$ to $\y$ which
  are allowed to cross the basepoint~$z$.  (This is denoted
  $\pi_2(\x,\y)$ in~\cite{OS04:HolomorphicDisks}, but in this book we
  have reserved $\pi_2$ for classes which have local multiplicity zero
  at~$z$.)  Since $e([\Sigma]) + n_\x([\Sigma]) + n_\y([\Sigma]) = 1$,
  we have
  $e(B) + n_\x(B) + n_\y(B)\equiv e(\tB)+n_\x(\tB)+n_\y(\tB)$ modulo~$1$.

  Label the components of $\Sigma \setminus (\alphas \cup \betas)$ by
  $R_i$ for $i = 1,\dots, N$, and write $\widetilde{B} = \sum_i m(R_i) R_i$
  with $i \ge 0$.  Then we will construct our surface by
  identifying copies of $R_i$:
  \[
  F \coloneqq \biggl(\coprod_{i=1}^N \coprod_{j=1}^{m(R_i)} R_i^{(j)}\biggr)\bigg/\sim
  \]
  where each $R_i^{(j)}$ is a diffeomorphic copy of $R_i$.  To describe the
  identification~$\sim$, first note that the $\alpha$-curves are
  divided by the $\beta$-circles into shorter pieces, which we call
  $\alpha$-segments, and similarly the $\beta$-circles are divided by
  the $\alpha$-curves into $\beta$-segments.  For each $\alpha$- or
  $\beta$-segment~$c$, let $R_1(c)$ and $R_2(c)$ be the two regions
  adjoining $c$, with $m(R_1(c)) \le m(R_2(c))$.  For $x \in c$, let
  $x_1^{(j)}$ and $x_2^{(j)}$ be the lifts of $x$ to
  $R_1^{(j)}(c)$ and $R_2^{(j)}(c)$, respectively.  Let $\delta_c \coloneqq
  m(R_2(c)) - m(R_1(c))$.
  Then the equivalence relation~$\sim$
  is generated by
  \begin{alignat*}{2}
    x_1^{(j)} &\sim x_2^{(j+\delta_a)}
      &\quad&\text{for $x$ in an $\alpha$-segment $a$, $j = 1,\dots,m(R_1(a))$}\\
    x_1^{(j)} &\sim x_2^{(j)}
      &&\text{for $x$ in an $\beta$-segment $b$, $j = 1,\dots,m(R_1(b))$.}
  \end{alignat*}
  Let $\Phi$ be the evident projection from $F$ to~$\Sigma$, so $\tB =
  \Phi_*[F]$.
  Topologically, $F$ is a manifold with boundary; requiring the map
  $\Phi$ to be a local isometry makes $F$ into a cone manifold.

  We next analyze the local behavior of~$F$ around a point $x \in
  \alpha_i \cap \beta_j$.
  Let the four regions around~$x$ be $R_1(x), \dots, R_4(x)$ in cyclic
  order, arranged so that $m(R_1(x)) \le m(R_i(x))$ for $i = 2,3,4$.
  The fact that $\widetilde{B} \in \piBig(\x,\y)$ means that
  \[
  m(R_1(x)) + m(R_3(x)) = m(R_2(x)) + m(R_4(x)) + C(x),
  \]
  where $C(x) \in \{0,-1,1\}$.  Specifically, $C(x) = 0$ when $x$ is
  in neither $\x$ nor~$\y$ or when $x$ is in both $\x$ and~$\y$.  For
  $x\in\Sigma$ let
  $D(x)$ be a small disk neighborhood of $x$ and let
  $N(x)=\Phi^{-1}(D(x))$. The geometry of $N(x)$ depends on $C(x)$ as follows:
  \begin{itemize}
  \item If $C(x)=0$ then $N(x)$ is a union of $m(R_1(x))$ disks and some
    half-disks, with the boundary of the half-disks lying above either
    $\alpha_i$ or $\beta_j$; in this case, $N(x)$ is a manifold with
    boundary but no corners.
  \item If $C(x) = 1$ then one component of $N(x)$ has a cone point at~$x$
    with total angle $2\pi m(R_1(x)) +
    \pi/2$.  The other components of $N(x)$ are half-disks as in the
    $C=0$ case.
  \item If $C(x) = -1$, then one component of $N(x)$ has a cone point at~$x$
    with total angle $2\pi (m(R_1(x))+1) - \pi/2$.  Again, the other
    components of $N(x)$ are half-disks.
  \end{itemize}
  Let $n_+\coloneqq\#\{x\in\alphas \cap \betas \mid C(x)=1\}$
  and $n_-\coloneqq\#\{x\in\alphas \cap \betas \mid C(x)=-1\}$, and
  let $h(x)$ be the number of half-disks in $N(x)$.

  Recall that $\bdy^\bdy \tB$ denotes the part of $\bdy \tB$ contained in
  $\bdy\bSigma$.  For each $a_i \in \CircPts$, let $\delta_i\coloneqq
  m(\bdy^\bdy \tB, a_i+\epsilon) - m(\bdy^\bdy \tB, a_i - \epsilon)$
  be the difference in multiplicities with which $\bdy^\bdy \tB$
  covers the regions on the two sides of $a_i$.
  The surface $F$ has corners at intersection
  points with $C(x)=\pm 1$ and at points above $\bdy \Sigma$, and the number of corners of
  $F$ over $\bdy\bSigma$ is $\OneQuart \sum_i \abs{\delta_i}$.
  Thus,
  \[
  e(\tB) = e(F) \equiv \frac{1}{4}\biggl(- n_+ + n_- +
  \sum_{i=1}^{4k} \abs{\delta_i}\biggr)
  \pmod 1.
  \]
  On the other hand,
  \begin{align*}
      n_\x(\tB) + n_\y(\tB)=n_{\x+\y}(\tB)
       &\equiv 
       \sum_{x \in \x+\y}
         \bigl(\OneQuart C(x) + \OneHalf h(x)\bigr)\\
       &\equiv \OneQuart (n_+ - n_-) +
         \frac{1}{2}\sum_{x \in [(\x + \y) \cap \alphas^a]} h(x).
  \end{align*}
  For the second congruence, observe that on each $\beta$-circle $\beta_j$ there
  is one $x_j \in \x$ and one $y_j \in \y$, and $h(x_j) = h(y_j)$.
  A similar statement holds for each $\alpha$-circle, so it suffices
  to sum $h(x)$ only when $x$ is in an $\alpha$-arc.

  It follows that
  \[
  e(\tB)+n_\x(\tB)+n_\y(\tB)\equiv \frac{1}{4}\sum_{i=1}^{4k} \abs{\delta_i}
  +\frac{1}{2}\sum_{x \in [(\x + \y) \cap \alphas^a]} h(x).
  \]
  It remains to show that the right-hand side is equal to
  $\epsilon(\bdy^\bdy \tB)$.  Consider the contribution to the
  right-hand side from a particular $\alpha$-arc $\alpha_\ell^a$ and its
  endpoints $a_i$ and~$a_j$.
  \begin{itemize}
  \item If $\ell$ is not in $o(\x)$ or $o(\y)$, the set of
    $\alpha$-arcs occupied by $\x$ and $\y$, respectively, 
    then there is no
    contribution to the $h(x)$ term and $\delta_i = -\delta_j$, so the
    net contribution is $\OneHalf$ if $\delta_i$ is odd and $0$
    otherwise, agreeing with the contribution to $\epsilon(\bdy^\bdy \tB)$.
  \item If $\ell$ is in both $o(\x)$ and $o(\y)$, then we have points
    $x_\ell \in \x \cap \alpha_\ell^a$ and $y_\ell \in \y \cap
    \alpha_\ell^a$.  We have $h(x_\ell) = h(y_\ell)$, so the same
    argument as in the previous case applies.
  \item If $\ell$ is in exactly one of $o(\x)$ and $o(\y)$, we have
    $\delta_j = \delta_i \pm 1$.  Assume without loss of generality
    that $\abs{\delta_j} = \abs{\delta_i} + 1$.  Let $x_\ell \in (\x +
    \y) \cap \alpha_\ell^a$ be the unique corner on this arc.  Then
    $h(x_\ell) = \abs{\delta_i}$, so
    \[
    \OneQuart(\abs{\delta_i} + \abs{\delta_j}) + \OneHalf h(x_\ell)
      \equiv \OneQuart(2\abs{\delta_i} + 1) + \OneHalf \abs{\delta_i}
      \equiv \OneQuart.
    \]
    Since exactly one of $\delta_i$ and~$\delta_j$ is odd, the
    contribution of $a_i$ and $a_j$ to $\epsilon(\bdy^\bdy \tB)$ is
    also $\OneQuart$, as desired. \qedhere
  \end{itemize}
\end{proof}

\begin{lemma}\label{lem:gB-mult}
  For $B_1 \in \pi_2(\x,\y)$ and $B_2 \in \pi_2(\y,\w)$, we have
  \begin{align*}
  \gb(B_1) \gb(B_2) &= \gb(B_1\ast B_2).
  \end{align*}
\end{lemma}
\begin{proof}
  This is obvious for the $\SpinC$ component.  For the Maslov
  component, the proof is similar to
  Proposition~\ref{prop:indAdditive}.
  If we set $a_i = \bdy^\alpha(B_i)$ and
  $b_i = \bdy^\alpha(B_i)$, the
  Maslov component of $\gb(B_1\ast B_2)$ is
{\setlength\multlinegap{0pt} 
  \begin{multline*}
    (-n_\x - n_\w - e)(B_1 * B_2)\\
    \begin{aligned}
      &= (-n_\x - n_\y - e)(B_1)+ (-n_\y - n_\w - e)(B_2)
      + (n_\y - n_\w)(B_1) + (n_\y - n_\x)(B_2)\\
      &= (-n_\x - n_\y - e)(B_1)+ (-n_\y - n_\w - e)(B_2)
      + a_1 \cdot b_2 - a_2 \cdot b_1\\
      &= (-n_\x - n_\y - e)(B_1)+ (-n_\y - n_\w - e)(B_2) + L_{\bdy
        \HD}(\bdy^\bdy(B_1), \bdy^\bdy(B_2)).
    \end{aligned}
  \end{multline*}
}
  where the second equality uses Lemma~\ref{lemma:jittered}
  and the third equality uses Lemma \penalty700\ref{lem:jittered2}.  The last formula is
  the Maslov component of $\gb(B_1)\gb(B_2)$, as desired.
\end{proof}

\begin{definition}
For a generator $\x \in \S(\HD)$, let 
$\gls*{Pprimeofx}$
be the subset of $\bigGroup(4k)$
given by $\{\gb(B)\mid B \in \pi_2(\x,\x)\}$.  
\end{definition}

\begin{corollary}
  $P'(\x)$ is a subgroup of $\bigGroup(4k)$.
\end{corollary}

This grading on domains is closely related to $\ind(B, \vec{\rhos})$.

\begin{lemma}\label{lem:index-vs-gB-2}
  For a compatible pair $(B, \vec{\rhos})$ where $\vec{\rhos} =
  (\rhos_1, \rhos_2, \dots, \rhos_\ell)$, we have 
  \[
  \gls*{gbig}
  = \lambda^{\ell-\ind(B,\vec{\rhos})} \grb(\vec{\rhos}).
  \]
\end{lemma}

\begin{proof}
  Applying Lemma~\ref{lem:iota-grading},
  Definition~\ref{def:emb-ind-emb-chi} and
  Definition~\ref{def:def-g-B}, in turn, gives
  \begin{align*}
    \grb(\vec{\rhos})
      &= (\iota(\vec{\rhos}), \bdy^\bdy B)\\
      &= (\ind(B, \vec{\rhos}) - \ell - e(B) - n_\x(B) - n_\y(B), \bdy^\bdy B)\\
      &= \lambda^{\ind(B,\vec{\rhos})-\ell} \gb(B).\qedhere
  \end{align*}
\end{proof}

\section{Type \textalt{$A$}{A} structures}
\label{sec:typeA-gradings}
In this section we describe the (unrefined) grading on $\CFAa(\HD)$.
Per Definition~\ref{def:module-nc-grading}, we are looking
for a right $\bigGroup(4k)$-set $\ABigGrSet(\HD)$ and a map $\grb\co
\Gen(\HD) \to \ABigGrSet(\HD)$ so that
\begin{align*}
\grb(\partial \x) &= \lambda^{-1}\grb(\x)\\
\grb(\x a) &= \grb(\x)\gr(a)\\
\shortintertext{and similarly, for higher products,}
\grb(m_{\ell+1}(\x,a_1,\dots,a_\ell))&=\lambda^{-1+\ell}\grb(\x)\gr(a_1)\cdots\gr(a_\ell).
\end{align*}
Since the
different $\SpinC$ structures do not interact in any algebra
action, we will construct one grading set $\ABigGrSet(\HD,\s)$ for
each $\SpinC$ structure $\s$ on~$Y$; then
$\ABigGrSet(\HD)$ is defined as
$\coprod_{\spinc\in\SpinC(Y)}\ABigGrSet(\HD,\spinc)$.

We first give the grading on the module with twisted
coefficients.

\begin{definition}
  To define the grading on the module $\tCFAa(\HD,\spinc)$, define
  $\tABigGrSet(\HD,\spinc)$ to be $\bigGroup(4k)$ (independent of
  $\HD$ and~$\spinc$), and set 
  $\gls*{tgrprime}(e^{B_0}\x) \coloneqq \gb(B_0)$.
\end{definition}

\begin{proposition}\label{prop:typeA-graded-tw}
  The map $\tgrb$
  defines a grading on $\tCFAa(\HD,\s)$.  Furthermore, for $B \in
  \pi_2(\x,\y)$, $\tgrb$ satisfies\index{grading!on $\tCFAa$}
  \begin{equation}\label{eq:AB-fund-tw}
  \tgrb(e^{B_0}\x)\cdot \gb(B) = \tgrb(e^{B_0 * B}\y),
  \end{equation}
  and any other grading $\widetilde{\gr}$ with this property is
  related to $\tgrb$ by
  left translation, in the sense that there is an element $g_0 \in
  \bigGroup(4k)$ so that for all~$\x$,
  \[
  \widetilde{\gr}(\x) = g_0 \cdot \tgrb(\x).
  \]
\end{proposition}

\begin{proof}
  Equation~\eqref{eq:AB-fund-tw} is immediate from
  Lemma~\ref{lem:gB-mult}:
  \[
  \tgrb(e^{B_0 * B}\y) = \gb(B_0 * B) = \gb(B_0)\gb(B) =
  \tgrb(e^{B_0}\x)\cdot \gb(B).
  \]

  To check that $\tgrb$ is a grading, from
  Definition~\ref{def:ainf-module} (in the form in
  Equation~\eqref{eq:DegreeAinfAction}), we want to show that for
  $\vec{\rhos} = (\rhos_1, \dots, \rhos_\ell)$ any
  sequence of sets of Reeb chords,
  $m_{\ell+1}(e^{B_0}\x,\allowbreak a(\rhos_1),\dots,a(\rhos_\ell))$
  is homogeneous of degree
  $\lambda^{\ell-1}\tgrb(e^{B_0}\x)\allowbreak\grb(\vec{\rhos})$.
  Let $B\in \pi_2(\x,\y)$ and let $\vec{\rhos} =
  (\rhos_1,\dots,\rhos_\ell)$ be a sequence of sets of Reeb chords compatible
  with~$B$.
  If $e^{B_0 * B}\y$ appears as a term in
  $m_{\ell+1}(e^{B_0}\x,\allowbreak a(\rhos_1),\dots,a(\rhos_\ell))$, then
  $\ind(B,\vec{\rhos}) = 1$, so by Equation~\eqref{eq:AB-fund-tw} and
  Lemma~\ref{lem:index-vs-gB-2},
  \[
   \tgrb(e^{B_0 * B}\y) = \tgrb(e^{B_0}\x) \cdot \gb(B) = \lambda^{\ell-1}
   \tgrb(e^{B_0}\x) \cdot \grb(\vec{\rhos}).
  \]

  For any other grading
  $\widetilde{\gr}$ satisfying Equation~\eqref{eq:AB-fund-tw}, let
  $g_0 = \widetilde{\gr}(e^{I_0}\x_0)$, where $I_0 \in
  \pi_2(\x_0,\x_0)$ is the trivial domain.  Then for any generator
  $\x$ and any $B_0\in\pi_2(\x_0,\x)$,
  \[
  \widetilde{\gr}(e^{B_0}\x) = \widetilde{\gr}(e^{I_0}\x)\cdot
  \gb(B_0) = g_0 \cdot \tgrb(e^{B_0}\x). \qedhere
  \]
\end{proof}

For the untwisted theory, the presence of periodic domains means that
we can not just define the
relative grading of $\x$ and $\y$ to be $\gb(B)\in \bigGroup(4k)$ for some $B \in
\pi_2(\x,\y)$, as the analogue of Equation~\eqref{eq:AB-fund-tw} would
fail for any periodic domain~$B$ for which $\gb(B)$ is non-trivial.
Instead, we look for a right $\bigGroup(4k)$-set $\ABigGrSet(\HD,\s)$ satisfying
an analogue of Equation~\eqref{eq:AB-fund-tw}.  The $\bigGroup(4k)$ action on
$\ABigGrSet(\HD,\s)$ will be transitive but not in general free; this
generalizes the fact that for
closed 3-manifolds the gradings in a non-torsion $\spin^c$ structure
take values in a quotient of~$\ZZ$.

\begin{definition}\label{def:typeA-grading}
Suppose that $\Gen(\HD,\spinc)\neq\emptyset$.
Pick any $\x_0 \in
  \Gen(\HD, \spinc)$ and 
define $\ABigGrSet(\HD, \s)$ to be the coset space
\begin{equation}\label{eq:GA-def}
\gls*{ABigGrSet}
\coloneqq  P'(\x_0) \backslash \bigGroup(4k).
\end{equation}
Define $\gls*{grprime}\co \S(\HD,\s) \to \ABigGrSet(\HD,\s)$ on $\x$ by
picking any $B_0\in\pi_2(\x_0,\x)$ and setting
\begin{equation}
\grb(\x) \coloneqq P'(\x_0)\cdot \gb(B_0).\label{eq:gr-typeA-def}
\end{equation}
The value of $\grb(\x)$ is independent of~$B_0$: If $B_1\in\pi_2(\x_0, \x)$ is
any other element, then $\gb(B_1 * (B_0)^{-1})\in P'(\x_0)$. Thus,
$\grb(\x) = P'(\x_0)\cdot \gb(B_1)$.

For completeness, suppose next that $\Gen(\HD,\spinc)=\emptyset$.
Choose a diagram $\HD'$ isotopic to $\HD$ and so that
$\Gen(\HD',\spinc)\neq\emptyset$ (by winding transverse to the
$\beta$-circles; compare~\cite[Lemma
5.2]{OS04:HolomorphicDisks}). Define
$\ABigGrSet(\HD,\s)=\ABigGrSet(\HD',\s)$, and let $\grb\co
\Gen(\HD,\spinc)\to \ABigGrSet(\HD,\spinc)$ be the unique map from the
empty set to $\ABigGrSet(\HD',\spinc)$.
\end{definition}

\begin{lemma}\label{lem:gr-fundamental-A}
  The function $\grb: \S(\HD,\s) \to \ABigGrSet(\HD,\s)$ above
  satisfies, for $\x, \y \in \S(\HD, \s)$ and $B \in \pi_2(\x,\y)$,
  \begin{equation}\label{eq:AB-fund}
  \grb(\x)\cdot\gb(B) = \grb(\y).
  \end{equation}
  Furthermore, if $\Gen(\HD,\s)\neq\emptyset$ then $\ABigGrSet(\HD,\s)$ is universal, in the sense that
  for any other $\bigGroup(4k)$-set $\widetilde{S}$ and grading
  function $\widetilde{\gr}$ on $\Gen(\HD,\s)$
  satisfying~\eqref{eq:AB-fund}, there is a unique $\bigGroup(4k)$-set
  map $f \co \ABigGrSet(\HD,
  \s) \to \widetilde{S}$ so that for any generator~$\x$, $\widetilde{\gr}(\x) = f(\grb(\x))$.
\end{lemma}

\begin{proof}
To check Equation~\eqref{eq:AB-fund}, pick $B_0 \in \pi_2(\x_0,\x)$.
Then $B_0 * B \in \pi_2(\x_0, \y)$, so we have
\[
\grb(\y) = P'(\x_0) \cdot \gb(B_0)\cdot\gb(B) = \grb(\x)\cdot\gb(B).
\]
(Note that we used the freedom to pick an arbitrary element of
$\pi_2(\x_0,\y)$ to define $\grb(\y)$.)

For universality, let $v_0$ be $\widetilde{\gr}(\x_0)$.  Then
Equation~\eqref{eq:AB-fund} implies that $P'(\x_0)$ stabilizes $v_0$.
We can therefore set, for $h \in \bigGroup(4k)$,
\[
f(P'(\x_0)\cdot h) \coloneqq v_0 \cdot h.
\]
It is immediate that $\widetilde{\gr}(\x) = f(\grb(\x))$.
\end{proof}
Observe that any right $\bigGroup(4k)$-set map from $\bigGroup(4k)$ to
itself is left translation as in Proposition~\ref{prop:typeA-graded-tw}.

\begin{corollary}\label{cor:change-base-gen-CFA}
If we use a different base generator
$\x_0'\in\S(\HD,\s)$ to define $\ABigGrSet(\HD, \s)$, the grading of an element $\x$ with respect to $\x_0'$ is
$\gb(B_0) \cdot \grb(\x)$ for any $B_0 \in \pi_2(\x_0',\x_0)$.
\end{corollary}

\begin{proof}
  This follows from noticing that $P'(\x_0') = \gb(B_0)\cdot
  P'(\x_0)\cdot\gb(B_0)^{-1}$.
\end{proof}
Left translation by
$\gb(B_0)$ is an isomorphism of right $\bigGroup(4k)$-sets 
\[P'(\x_0)\backslash\bigGroup(4k)\cong P'(\x'_0)\backslash\bigGroup(4k).\]

\begin{proposition}\label{prop:typeA-graded}
  For the grading $\grb$ of Definition~\ref{def:typeA-grading},
  $\CFAa(\HD,\s)$ is a graded right
  $\Ainf$ module in the sense of
  Definition~\ref{def:module-nc-grading}.
  \index{grading!on $\CFAa$}%
\end{proposition}
\begin{proof}
  This follows from Lemma~\ref{lem:gr-fundamental-A} as in the proof of
  Proposition~\ref{prop:typeA-graded-tw}.
\end{proof}

We now have the following graded version of
Theorems~\ref{thm:A-invariance} and~\ref{thm:tA-invariance}
(and indeed Theorem~\ref{intro:A-invariance}):

\begin{theorem}
  \label{thm:A-invariance-graded}\index{invariance!of $\CFAa$!graded}
  Up to graded $\Ainf$ homotopy equivalence, the
  $\bigGroup(4k)$-graded $\Ainf$ module 
  $\CFAa(\HD,\spinc)$ over $\Alg(\bdy\HD)$ is independent of the choices
  made in its definition. That is, for any choices of provincially
  admissible Heegaard diagrams $\HD_1$ and $\HD_2$, sufficiently
  generic admissible almost complex structures $J_i$ for $\HD_i$ and
  base generators $\x_0^i\in\Gen(\HD_i,\spinc)$, there is an isomorphism of
  $\bigGroup(4k)$-sets 
  \[
  \phi\co \ABigGrSet(\HD_1,\x_0^1)\to \ABigGrSet(\HD_2,\x_0^2)
  \]
  and an $\Ainf$ homotopy equivalence of $\Ainf$ modules
  \[
  f=\{f_{1+i}\co \CFAa(\HD_1,J_1)\otimes\Alg(\bdy\HD)^{\otimes i}\to \CFAa(\HD_2,J_2)\}
  \]
  such that, for $\x\in\Gen(\HD_1)$ and homogeneous elements
  $a_1,\dots,a_\ell\in\Alg(\bdy\HD)$, 
  \[
  \grb(f_{1+\ell}(\x,a_1,\dots,a_\ell))=\phi(\grb(\x))\grb(a_1)\dots\grb(a_\ell)\lambda^{\ell}.
  \]

  The graded $\Ainf$
  bimodule $\tCFAa(\HD,\s)$
  over $\Field[H_2(Y,\partial Y)]$ and $\Alg(\bdy\HD)$ is similarly
  invariant.
\end{theorem}

\begin{proof}
  We prove invariance of the untwisted theory; the proof for the
  twisted theory is similar.

  Independence of the grading from the choice of base generator
  (assuming $\HD_1=\HD_2$ and $J_1=J_2$) is
  immediate from Corollary~\ref{cor:change-base-gen-CFA}: the map $g$
  is the identity map (i.e., $g_1=\Id$ and $g_i=0$ for $i>1$) and the
  map $\phi$ is given by $\phi(s)=g'(B_0)s$, where
  $B_0\in\pi_2(\x_0^2,x_0^1)$. (Left translation is an
  isomorphism of right $\bigGroup(4k)$-sets.)

  The rest of the proof follows the proof of the ungraded version,
  Theorem~\ref{thm:A-invariance}, noting that the
  identifications induced by changes of complex structure, isotopies,
  handleslides, and stabilizations all respect gradings.

  Stabilization invariance is immediate: stabilizing near the basepoint, the $\Ainf$ modules before
  and after stabilization are identified, along with all the grading data. For changes of complex
  structure, the fact that the continuation
  maps preserve gradings is a straightforward adaptation of
  proof that $\CFAa$ is graded
  (Proposition~\ref{prop:typeA-graded}). The key point is that
  now, if $\y$ appears in
  $f^{J_r}_{1+\ell}(\x,a(\rhos_1),\dots,a(\rhos_\ell))$, 
  then there is a domain $B\in\pi_2(\x,\y)$ with
  $\ind(B,\vec{\rhos})=0$ (rather than $1$, as it would be for $m_{1+\ell}$ in
  place of $f^{J_r}_{1+\ell}$), and hence
  $$\grb(\y)=\grb(\x)\cdot \grb(a_1)\cdots
  \grb(a_\ell)\lambda^{\ell}.$$
  This is precisely what is needed for
  $f^{J_r}=\{f^{J_r}_n\}_{n=1}^{\infty}$ to be a graded $\Ainf$
  homomorphism. 

  The proof of isotopy invariance is the same as the proof of
  invariance under change of complex structure except that the base
  generator may disappear during the isotopy. If there are no
  generators in $\Gen(\HD_1,\spinc)$ or $\Gen(\HD_2,\spinc)$ then the
  result is obvious. Otherwise, we can choose the isotopy so that
  there is an intermediate Heegaard diagram $\HD_{3/2}$ so that there
  is some generator $\x^1_0$ in both $\Gen(\HD_1,\spinc)$ and
  $\Gen(\HD_{3/2},\spinc)$, and a generator $\x^2_0$ in both
  $\Gen(\HD_{3/2},\spinc)$ and $\Gen(\HD_{2},\spinc)$. The same
  argument as for independence of complex structures gives a graded
  homotopy equivalence between $\CFAa(\HD_1,\spinc)$ and
  $\CFAa(\HD_{3/2},\spinc)$, using $\x^1_0$ to define the grading, and
  a graded homotopy equivalence between $\CFAa(\HD_{3/2},\spinc)$ and
  $\CFAa(\HD_{2},\spinc)$, using $\x^2_0$ to define the
  grading. Together with the independence of the grading from the base
  generator, which we already verified, this gives the invariance
  result for isotopies.

  For handleslides, we must take into account the fact that the
  grading sets for the two diagrams are different.  Specifically,
  suppose that $\HD_2$ is gotten from $\HD_1$ by a handleslide; as in
  Section~\ref{sec:A-Handleslides}, suppose this is a handleslide of
  an $\alpha$-arc over an $\alpha$-circle, and write
  $\HD_2=(\Sigma,\alphas^H,\betas,z)$,
  $\HD_1=(\Sigma,\alphas,\betas,z)$. Since we have already verified
  independence from the basepoint, we may choose $\x_0^2$ to be the
  generator connected to $\x_0^1$ by a union of small triangles $T_\x$
  (see Section~\ref{sec:CFD-handleslide}).  Then, the grading set
  for $\HD_1$ is $P'(\x_0^1) \backslash \bigGroup(4k)$, whereas the
  grading set for $\HD_2$ is $P'(\x_0^2) \backslash \bigGroup(4k)$.
  There is a map $\Phi\colon P'(\x_0)\to P'(\x_0')$ defined by sending
  $B\in \pi_2(\x_0^1,\x_0^1)$ to the element
  $\Phi(B)\in\pi_2(\x_0^2,\x_0^2)$ uniquely determined by
   $$T_\x *_{13} B=T_\x *_{12} * B_{\alpha,\alpha^H} *_{23} \Phi(B),$$
   using Lemma~\ref{lem:triple-pi2-decomp}. (Here,
   $B_{\alpha,\alpha^H}\in\pi_2(\Theta_o,\Theta_{o'})$ for appropriate
   $o$ and $o'$.)  By inspection, $\gb(\Phi(B))=\gb(B)$. So,
   $P'(\x_0^1)=P'(\x_0^2)$, and
   $\ABigGrSet(\HD_1,\x_0^1)=\ABigGrSet(\HD_2,\x_0^2)$.

   Next we check that the handleslide map $f^{\alphas,\alphas^H,\betas}$
   respects this identification of grading sets. By
   Corollary~\ref{cor:tri-emb-ind} (and the definition of the map
   $f^{\alphas,\alphas^H,\betas}$, Equation~\eqref{eq:TriangleMapA}),
   \[
   f^{\alphas,\alphas^H,\betas}_{1+\ell}(\x,a(\rhos_1),\dots,a(\rhos_\ell))
   \]
   is given by a count of triangles in homology classes $b$ with
   $\ind(b(B),(\rhos_1,\dots,\rhos_\ell))=0$. So, as for the map
   associated to changes of almost complex structure, if $\y$ occurs
   in
   $f^{\alphas,\alphas^H,\betas}_{1+\ell}(\x,a(\rhos_1),\dots,a(\rhos_\ell))$
   then 
  \[
  \grb(\y)=\grb(\x)\cdot \grb(a_1)\cdots \grb(a_\ell)\lambda^{\ell}.
  \]
  This completes the proof.
\end{proof}

\section{Type \textalt{$D$}{D} structures}
\label{sec:typeD-gradings}

We turn now to defining a grading on $\CFDa(\HD)$.
One point that merits attention is that the relevant grading groups
$\bigGroup(4k)$ and $\smallGroup(\PMC)$ depend on the orientation
of~$\PMC$.  In\index{grading!group!dependence on orientation}
particular, in the definition of the product on $\bigGroup(4k)$,
\[
  (j_1,\alpha_1)\cdot(j_2,\alpha_2) \coloneqq (j_1 + j_2 +
    L(\alpha_1, \alpha_2), \alpha_1 + \alpha_2)
\]
(from Equation~\eqref{eq:Gn-mult-def}), the correction term
$L(\cdot,\cdot)$ depends on whether we are in $\bdy \HD$ or $-\bdy
\HD$.  To avoid confusion, let 
$
\gls*{reverseOr}\co Z \to -Z$
be the identity map, which is orientation-reversing,\index{orientation reversal} and let
$
\gls*{reverseOrGp}
\co \bigGroup(4k) \to
\bigGroup(4k)$ or $\smallGroup(\PMC) \to \smallGroup(-\PMC)$ be the
induced map:
\begin{equation}\label{eq:def-R}
R(j,\alpha) = (j, r_* (\alpha)).
\end{equation}
It follows from the form of the product that $R$ is a group
anti-homomorphism, so, for instance,
\[
  R(\gb(B_1\ast B_2)) = R(\gb(B_2)) R(\gb(B_1)).
\]

We have an analogue of Lemma~\ref{lem:index-vs-gB-2} for the moduli
spaces relevant for $\CFDa$.

\begin{lemma}\label{lem:index-vs-gB}
  For a compatible pair $(B,\vec\rho)$,  we have
  \[\grb(-\vec{\rho})R(\gb(B)) = (-\mathord{\ind}(B,\vec{\rho}),0).\]
\end{lemma}
\begin{proof}
  This can be deduced from Lemma~\ref{lem:index-vs-gB-2}, but we find
  it clearer to give a direct proof.
  The $\SpinC$ component of $R(\gb(B))$ is the negative of the
  $\SpinC$ component of $\grb(-\vec{\rho})$, since all the boundary
  intervals are reversed in the latter. Thus the $\SpinC$
  component of the result is~$0$, and Maslov component of the
  product is obtained by adding the Maslov components of the factors. The Maslov
  component of $R(\gb(B))$ is, by definition, $-e(B)-n_\x(B)-n_\y(B)$.  On
  the other hand, if $\vec \rho = (\rho_1,\dots,\rho_n)$ then
  \begin{equation}\label{eq:gr-rho-seq}
    \begin{aligned}
  \grb(-\vec{\rho}) &= (-\OneHalf, -[\rho_1])\cdots(-\OneHalf, -[\rho_n])\\
    &= \Bigl(-{\textstyle \frac{n}{2}} +
           \sum_{i < j} L(r_*(-[\rho_i]), r_*(-[\rho_j]),
           -\sum_i[\rho_i]\Bigr)\\
    &= \Bigl(-{\textstyle \frac{n}{2}} -
           \sum_{i < j} L([\rho_i], [\rho_j]),
           -\sum_i[\rho_i]\Bigr)\\
    &= \Bigl(-n-\iota(\vec{\rho}),
           -\sum_i[\rho_i]\Bigr),
    \end{aligned}
  \end{equation}
  where we use Lemma~\ref{lem:G-n-products}, bilinearity of~$L$, and
  the definition of~$\iota$ (Equation~\eqref{eq:def-iota}).  Thus, by
  Definition~\ref{def:emb-ind-emb-chi}, the Maslov component of the
  product is $-\mathord{\ind}(B,\vec{\rho})$.
\end{proof}

\index{grading!on $\tCFDa$}
\begin{definition}
  For $\tCFDa(\HD,\spinc)$ as for $\tCFAa(\HD,\spinc)$, define
  $\tDBigGrSet(\HD,\spinc)$ to be $\bigGroup(4k)$, and set
  $\gls*{tgrprime}(e^{B_0}\x) \coloneqq R(\gb(B_0))$.
\end{definition}

\begin{proposition}\label{prop:grading-typeD-tw}
  The map $\tgrb$ above
  defines a grading on $\tCFDa(\HD,\s)$.  Furthermore, for $B \in
  \pi_2(\x,\y)$, $\tgrb$ satisfies
  \begin{equation}\label{eq:BD-fund-tw}
  R(\gb(B)) \cdot\tgrb(e^{B_0}\x) = \tgrb(e^{B_0 * B}\y),
  \end{equation}
  and any other grading $\widetilde{\gr}$ with this property is
  related to $\tgrb$ by right translation, in the sense that there is
  an element $g_0 \in
  \bigGroup(4k)$ so that for all~$\x$,
  \[
  \widetilde{\gr}(\x) = \tgrb(\x)\cdot g_0.
  \]
\end{proposition}

\begin{proof}
  If there is a term of the form  $a(-\vec\rho)e^{B_0 * B}\y$ occurring in
  $\partial(e^{B_0}\x)$ then $\ind(B,
  \vec{\rho}) = 1$, and so
  \begin{align*}
    \tgrb(a(-\vec\rho)e^{B_0 * B}\y)
      &= \grb(a(-\vec\rho)) R(\gb(B)) \tgrb(e^{B_0}\x)\\
      &= \lambda^{-1} \tgrb(e^{B_0}\x).
  \end{align*}
  The universality statement is analogous to the case for $\tCFAa$,
  Proposition~\ref{prop:typeA-graded-tw}.
\end{proof}

Again, for the untwisted theory, we need to find an appropriate left
$\bigGroup(4k)$-set satisfying a version of
Equation~\eqref{eq:BD-fund-tw}.

\begin{definition}\label{def:typeD-grading}
  To define the grading on $\CFDa(\HD,\spinc)$, 
take
  $\DBigGrSet(\HD,\s)$ to have the same underlying set as
  $\ABigGrSet(\HD,\s)$, with the left action of $g$ on
  $\DBigGrSet(\HD,\s)$ being given by the right action of $R(g)$ on
  $\ABigGrSet(\HD,\s)$.  Concretely, if
  $\Gen(\HD,\spinc)\neq\emptyset$,
  \[
  \DBigGrSet(\HD,\s) \coloneqq \bigGroup(4k)/R(P'(\x_0)),
  \]
and define $\gls*{grprime}\co \S(\HD,\s) \to \DBigGrSet(\HD,\s)$ on $\x$ by
picking any $B_0\in\pi_2(\x_0,\x)$ and setting
\begin{equation}
\grb(\x) \coloneqq R(\gb(B_0)) \cdot R(P'(\x_0)).\label{eq:gr-typeD-def}
\end{equation}
(If $\Gen(\HD,\spinc)=\emptyset$ we proceed as in this degenerate case of Definition~\ref{def:typeA-grading}.)
\end{definition}

\begin{lemma}\label{lem:gr-fundamental-D}
  The function $\grb: \S(\HD,\s) \to \DBigGrSet(\HD,\s)$ above
  satisfies, for $\x, \y \in \S(\HD, \s)$ and $B \in \pi_2(\x,\y)$,
  \begin{equation}\label{eq:BD-fund}
  R(\gb(B))\cdot\grb(\x) = \grb(\y).
  \end{equation}
  Furthermore, $\DBigGrSet(\HD,\s)$ is universal, in the sense that
  for any other $\bigGroup(4k)$-set $\widetilde{S}$ and grading
  function $\widetilde{\gr}$ on $\Gen(\HD,\s)$
  satisfying~\eqref{eq:BD-fund}, there is a unique $\bigGroup(4k)$-set
  map $f \co \DBigGrSet(\HD,
  \s) \to \widetilde{S}$ so that for any generator~$\x$, $\widetilde{\gr}(\x) = f(\grb(\x))$.
\end{lemma}

\begin{proof}
  This follows from Lemma~\ref{lem:gr-fundamental-A} by reversal of
  all products.
\end{proof}

\index{grading!on $\CFDa$}
\begin{proposition}\label{prop:grading-typeD}
  For the grading $\grb$ given by Definition~\ref{def:typeD-grading},
  $\CFDa(\HD,\s)$ is a graded left differential module
\end{proposition}
\begin{proof}
  This is similar to the proof of
  Proposition~\ref{prop:grading-typeD-tw}, and is left to the reader.
\end{proof}

\begin{theorem}
  \label{thm:D-invariance-graded}\index{invariance!of $\CFDa$!graded}
  Up to graded homotopy equivalence, the $\bigGroup(4k)$-set
  graded differential module $\CFDa(\HD,\spinc)$ over $\Alg(-\bdy\HD)$ is
  independent of the choices made in its definition. That is, for any
  choices of provincially admissible Heegaard diagrams $\HD_1$ and
  $\HD_2$, sufficiently generic admissible almost complex structures
  $J_i$ for $\HD_i$ and base generators $\x_0^i\in\Gen(\HD_i,\spinc)$,
  there is an isomorphism of $\bigGroup(4k)$-sets
  \[
  \phi\co \DBigGrSet(\HD_1,\x_0^1)\to \DBigGrSet(\HD_2,\x_0^2)
  \]
  and a homotopy equivalence of differential modules
  \[
  f\co \CFDa(\HD_1,J_1)\to \CFDa(\HD_2,J_2)
  \]
  such that, for each $\x\in\Gen(\HD_1)$ and each term $a\y$ in
  $f(\x)$ (where $\y\in\Gen(\HD_2,\spinc)$ and $a\in\Alg(-\bdy\HD)$), 
  \[
  \grb(a)\grb(\y)=\phi(\grb(\x)).
  \]

  The \dg bimodule $\tCFDa(\HD,\s)$ over $\Field[H_2(Y,\partial
  Y)]$ and $\Alg(-\bdy\HD)$ is similarly invariant.
\end{theorem}
\begin{proof}
  This is similar to the proof of
  Theorem~\ref{thm:A-invariance-graded}, and is left to the reader.
\end{proof}

\section{Refined gradings}
\label{sec:refined-gradings-2}
\index{grading!on $\Alg(\PMC)$!refined|(}
We turn now to a refined grading $\gr$ on these modules, building 
on the refined grading of the algebra.
\index{grading!group!refined}
We start with a review of the small grading group
$\smallGroup(\PMC)$.  The group
$\gls*{smallGroup}$
can be viewed either
abstractly as a central extension or concretely as a subgroup of
$\bigGroup(4k)$.  Abstractly, it is the central extension of
$H_1(F(\PMC))$, with a sequence of groups
\[
\ZZ \overset{\lambda}{\hookrightarrow} \smallGroup(\PtdMatchCirc)
  \overset{[\cdot]}{\twoheadrightarrow} H_1(F)
\]
satisfying the relation
\begin{equation*}
  gh = hg\lambda^{2([g] \cap [h])}
\end{equation*}
(Equation~\eqref{eq:commutators}).
Concretely, the group 
$\gls*{smallGroup}$
is the subgroup of
$\bigGroup(4k)$ consisting of those elements $(j, \alpha)$ where
$M_*\bdy(\alpha)=0$, that is, for each pair $p,q$ of points
in~$\CircPts$ identified by the matching~$M$ defining~$\PMC$, the sum
of the multiplicities of
$\bdy\alpha$ at~$p$ and at~$q$ is~$0$.

The refined grading on $\Alg(\PMC,0)$ is built from grading refinement
data: a choice of a base idempotent $I(\SetS_0)\in\Alg(\PtdMatchCirc,0)$
and, for every other minimal idempotent $I(\SetS)$, an element
\glsadd{refinementdata}
\index{grading!refinement data}
$\psi(\SetS)=\psi(I(\SetS))\in \bigGroup(4k)$ satisfying
$M_*\bdy([\psi(\SetS)])=\SetS-\SetS_0$.  The grading on the algebra
is then defined by
\[
  \gls*{smallgr}
  \bigl(I(\SetS)a(\rhos)I(\SetT)\bigr) \coloneqq \psi(\SetS)\grb(a(\rhos))\psi(\SetT)^{-1}
\]
(Equation~\eqref{eq:small-grading}).
\index{grading!on $\Alg(\PMC)$!refined|)}

We now turn to refined gradings on domains.
Recall first that every generator $\x \in \Gen(\HD)$ has an associated
set $\gls*{ox}$ of occupied $\alpha$-arcs (Definition~\ref{def:generator})
and associated idempotents $\gls*{IAofx}
= I(o(\x))$ and $\gls*{IDofx} = I([2k] \setminus o(\x))$.
\index{occupied $\alpha$-arcs}%
Then,
for $\x,\y \in \Gen(\HD)$ and $B\in\pi_2(\x,\y)$, define
\begin{align}\label{eq:DSmallg}
\gls*{gsmall}
&\coloneqq \psi(I_A(\x))\gb(B)\psi(I_A(\y))^{-1}\\
\gls*{Pofx}&\coloneqq \{\gs(B)\mid B \in \pi_2(\x,\x)\}.
\end{align}
\begin{lemma}\label{lem:gs-fund}
  The function $\gs$ satisfies $\gs(B) \in \smallGroup(\PMC)$ and
  $\gs(B_1*B_2)=\gs(B_1)\cdot\gs(B_2)$.
\end{lemma}
\begin{proof}
  To verify that $\gs(B)\in\smallGroup(\PMC)$, compute
  \begin{align*}
    M_*\bdy([\gs(B)]) &= M_*\bdy([\psi(I_A(\x))]) +
        M_*\bdy([\gb(B)]) - M_*(\bdy([\psi(I_A(\y))]))\\
      &= (o(\x) - \SetS_0) + (o(\y) - o(\x)) - (o(\y) - \SetS_0)\\
      &= 0.
  \end{align*}
  The fact that $\gs(B_1*B_2)=\gs(B_1)\cdot\gs(B_2)$ is immediate from the
  definitions and Lemma~\ref{lem:gB-mult}.
\end{proof}
In particular, $P(\x)$ is a subgroup of $G(\PMC)$.

\begin{definition}\label{def:refined-gr-CFA}
  To define the refined grading on $\tCFAa(\HD,\spinc)$, define
  $\tASmallGrSet(\HD,\spinc)$ to be $\smallGroup(\bdy\HD)$ (independent of
  $\HD$ and~$\spinc$), and set 
  $\gls*{tgr}(e^{B_0}\x) \coloneqq \gs(B_0)$.
\index{grading!on $\tCFAa$!refined}%

  To define the refined grading on $\CFAa(\HD,\spinc)$, 
  suppose first that $\Gen(\HD,\spinc)\neq\emptyset$.
  Pick any $\x_0 \in \Gen(\HD, \spinc)$.
  Define $\ASmallGrSet(\HD,\spinc)$ by
  \begin{align*}
    \gls*{ASmallGrSet}
    (\HD,\s)&\coloneqq P(\x_0)\backslash \smallGroup(\bdy\HD)\\
\intertext{and define $\gr(\x)$ by picking any $B_0 \in
  \pi_2(\x_0,\x)$ and setting}
    \gls*{smallgr}
    (\x) &\coloneqq P(\x_0) \cdot \gs(B).
  \end{align*}
  \index{grading!on $\CFAa$!refined}%
  For completeness, if $\Gen(\HD,\spinc)=\emptyset$ choose a diagram
  $\HD'$ isotopic to $\HD$ and so that
  $\Gen(\HD',\spinc)\neq\emptyset$, let
  $\ASmallGrSet(\HD,\s)=\ASmallGrSet(\HD',\s)$ and let $\gr\co
  \Gen(\HD,\spinc)\to \ASmallGrSet(\HD,\spinc)$ be the unique map from
  the empty set to $\ASmallGrSet(\HD',\spinc)$.
\end{definition}

\begin{proposition}
  The function $\tgr \co\tGen(\HD,\spinc) \to
  \tASmallGrSet(\HD,\spinc)$ above
  defines a grading on $\tCFAa(\HD,\s)$.
  Likewise, $\gr \co \Gen(\HD,\spinc) \to
  \ASmallGrSet(\HD,\spinc)$
  defines a grading on $\CFAa(\HD,\s)$.
\end{proposition}

\begin{proof}
  The proof follows the proofs of
  Propositions~\ref{prop:typeA-graded-tw}
  and~\ref{prop:typeA-graded}, noting that if
  $m_{n+1}(\x,a_1,\dots,a_n) \ne 0$, the idempotent of $\x$ agrees
  with the left idempotent of $a_1$ and the right idempotent of $a_i$
  agrees with the left idempotent of $a_{i+1}$, so the corresponding
  $\psi$'s cancel.
\end{proof}

To define the refined gradings on $\tCFDa$ and $\CFDa$, we first need to fix
grading refinement data for $\Alg(-\bdy \HD,0)$. We choose it as follows.
Fix a choice of $(\SetS_0)_{\bdy\HD}$
and $\psi_{\bdy\HD}$ determining a grading refinement for
$\Alg(\bdy\HD,0)$ and define
\begin{align*}
  (\SetS_0)_{-\bdy\HD} &\coloneqq [2k] \setminus (\SetS_0)_{\bdy\HD}\\
  \psi_{-\bdy\HD}(\SetS) &\coloneqq R(\psi([2k] \setminus \SetS))^{-1}.
\end{align*}
With this choice,
\[
M_*\bdy([\psi_{-\bdy\HD}(\SetS)])
  = -M_*\bdy([\psi_{\bdy\HD}([2k]\setminus\SetS)])
  = -(([2k] - \SetS) - (\SetS_0)_{\bdy\HD})
  = \SetS - (\SetS_0)_{-\bdy\HD},
\]
so $(\SetS_0)_{-\bdy\HD}$ and $\psi_{-\bdy\HD}(\SetS)$ define a
grading refinement on $\Alg(-\bdy\HD,0)$; we call this the
\emph{reverse} of the grading refinement $(\SetS_0)_{\bdy \HD}$ and
$\psi_{\bdy\HD}$.
Then, for $B \in \pi_2(\x,\y)$,
\begin{equation}
  \begin{aligned}
  R(g(B)) &= R\bigl(\psi_{\bdy\HD}(I_A(\x)) \gb(B) \psi_{\bdy\HD}(I_A(\y))^{-1}\bigr)\\
   &= \psi_{-\bdy\HD}([2k] \setminus I_A(\y)) R(\gb(B))
     \psi_{-\bdy\HD}([2k] \setminus I_A(\x))^{-1}\\
   &= \psi_{-\bdy\HD}(I_D(\y)) R(\gb(B)) \psi_{-\bdy\HD}(I_D(\x))^{-1}.
  \end{aligned}
\label{eq:R-small-g}
\end{equation}

\begin{definition}\label{def:refined-gr-CFD}
  To define the refined grading on $\tCFDa(\HD,\spinc)$, define
  $\tDSmallGrSet(\HD,\spinc)$ to be $\smallGroup(-\bdy\HD)$ (independent of
  $\HD$ and~$\spinc$), and set
  $\gls*{tgr}(e^{B_0}\x) \coloneqq R(\gs(B_0))$.
  \index{grading!on $\tCFDa$!refined}%

  To define the refined grading on $\CFDa(\HD,\spinc)$, 
  suppose first that $\Gen(\HD,\spinc)\neq\emptyset$.
  Pick any $\x_0 \in \Gen(\HD, \spinc)$.
  Define
  $\DSmallGrSet(\HD,\spinc)$ by
  \begin{align*}
    \gls*{DSmallGrSet}
    (\HD,\s)&\coloneqq\smallGroup(-\bdy\HD)/R(P(\x_0))\\
\intertext{and define $\gr(\x)$ by picking any $B_0 \in
  \pi_2(\x_0,\x)$ and setting}
    \gls*{smallgr}
    (\x)&\coloneqq R(\gs(B)) \cdot R(P(\x_0)).
  \end{align*}
  \index{grading!on $\CFDa$!refined}%

  For completeness, if $\Gen(\HD,\spinc)=\emptyset$ choose a diagram
  $\HD'$ isotopic to $\HD$ and so that
  $\Gen(\HD',\spinc)\neq\emptyset$, let
  $\DSmallGrSet(\HD,\s)=\DSmallGrSet(\HD',\s)$ and let $\gr\co
  \Gen(\HD,\spinc)\to \DSmallGrSet(\HD,\spinc)$ be the unique map from
  the empty set to $\DSmallGrSet(\HD',\spinc)$.
\end{definition}

\begin{proposition}
  The function $\tgr \co\tGen(\HD,\spinc) \to
  \tDSmallGrSet(\HD,\spinc)$ above
  defines a grading on $\tCFDa(\HD,\s)$.
  Likewise, $\gr \co \Gen(\HD,\spinc) \to
  \DSmallGrSet(\HD,\spinc)$
  defines a grading on $\CFDa(\HD,\s)$.
\end{proposition}
\begin{proof}
  The proof follows the proofs of
  Propositions~\ref{prop:grading-typeD-tw}
  and~\ref{prop:grading-typeD}, noting that if
  $a\x \ne 0$, the right idempotent of $a$ agrees
  with the idempotent of $\x$, so the intermediate
  $\psi$ cancels.
\end{proof}

We have the following analogue of the graded invariance theorems
(Theorems~\ref{thm:D-invariance-graded}
and~\ref{thm:A-invariance-graded}) with respect to the refined
gradings.
\begin{theorem}\label{thm:A-D-invariance-small-gr}
  Up to graded $\Ainf$ homotopy equivalence, the
  $\smallGroup(\bdy\HD)$-set graded $\Ainf$ module $\CFAa(\HD,\spinc)$
  over $\Alg(\bdy\HD)$ is independent of the choices of provincially
  admissible Heegaard diagram, sufficiently generic admissible almost
  complex structure and base generator used to define it; and up to
  graded homotopy equivalence, the $\smallGroup(-\bdy\HD)$-set graded
  differential module $\CFDa(\HD,\spinc)$ over $\Alg(-\bdy\HD)$ is
  independent of the choices of provincially admissible Heegaard
  diagram, sufficiently generic admissible almost complex structure
  and base generator used to define it. Similar statements hold for
  $\Field[H_2(Y,\partial Y)]$--$\Alg(\bdy\HD)$ and $\Field[H_2(Y,\partial
  Y)]$--$\Alg(-\bdy\HD)$ bimodules $\tCFAa(\HD,\spinc)$ and
  $\tCFDa(\HD,\spinc)$.
\end{theorem}
\begin{proof}
  This is immediate from the definitions and
  Theorems~\ref{thm:D-invariance-graded}
  and~\ref{thm:A-invariance-graded}.
\end{proof}

\begin{remark}
  Theorem~\ref{thm:A-D-invariance-small-gr} does not assert that the
  refined gradings are independent of the grading refinement data. For
  some discussion of this point, see
  Remark~\ref{rem:refined-grading-change}.
\end{remark}
\section{Tensor product}
\label{sec:GradedPairingThm}

We now turn our attention to a graded version of the
pairing theorem. Let $(Y_1,\phi_1\co F(\PMC)\to \bdy Y_1)$ and
$(Y_2,\phi_2\co F(-\PMC)\to\bdy Y_2)$ be
bordered three-manifolds which agree along their boundary. As
discussed in Section~\ref{sec:GluingDiagrams}, given bordered
Heegaard diagrams $\HD_i$ for the $(Y_i,\phi_i)$, gluing $\HD_1$ and
$\HD_2$ along their boundary gives a Heegaard diagram
$\HD=\HD_1\cup_\bdy \HD_2$ for $Y=Y_1\cup_\bdy Y_2$.

For each $\SpinC$ structure $\spinc_1$ over $Y_1$,
the right module $\CFAa(\HD_1,\spinc_1)$ is graded by the set
$P_1(\x_1)\backslash \smallGroup(\PtdMatchCirc)$,
where $\x_1\in\Gen(\HD_1)$ is a chosen base generator representing
$\spinc_1$,
and $P_1(\x_1)$ denotes the image of $\pi_2(\x_1,\x_1)$ in $\smallGroup(\PMC)$.
Similarly, for each $\SpinC$ structure $\spinc_2$ over $Y_2$,
the left module $\CFDa(\HD_2,\spinc_2)$ is graded by the set
$\smallGroup(\PtdMatchCirc)\slash R(P_2(\x_2))$,
where $\x_2\in\Gen(\HD_2)$ is a chosen base generator representing
$\spinc_2$,
and $P_2(\x_2)$ denotes the image of $\pi_2(\x_2,\x_2)$ in
$\smallGroup(-\PMC)$.
(Here the grading refinement on $G(-\PMC)$ is the reverse
of the grading refinement on $G(\PMC)$, so that the grading
refinement used in the definition of the grading of
$\CFAa(\HD_1,\s_1)$ is the same as the grading refinement used in the
definition of $\CFDa(\HD_2,\s_2)$.)
Correspondingly, the tensor product 
$\CFAa(\HD_1,\spinc_1)\DT\CFDa(\HD_2,\spinc_2)$ inherits
a grading 
$\gls*{grBox}$
with values in  the double-coset space
$$P_1(\x_1)\backslash \smallGroup(\PtdMatchCirc)\slash R(P_2(\x_2)).$$
We choose the base generators $\x_1$ and $\x_2$ to occupy complementary $\alpha$-arcs,
i.e., so that $\x_1\times\x_2\in\Gen(\HD)$. Let $\x \coloneqq \x_1
\times \x_2$.

The above double-coset space inherits an action by
$\ZZ$ (acting on the Maslov component). The quotient by this action is
naturally identified with the quotient of $H_1(F)$ by the image of
$H_2(Y_1,F)\oplus H_2(Y_2,F)$ under the boundary homomorphism, or,
equivalently, with the image of $H_1(F)$ in $H_1(Y)$.
More succinctly, we have
\[
\xymatrix{
  P_1(\x_1)\backslash \smallGroup\slash R(P_2(\x_2)) \ar[d] \ar[dr]^{[\cdot]}& \\
  \bigl(P_1(\x_1)\backslash \smallGroup\slash R(P_2(\x_2))\bigr)/\ZZ
  \ar[r]_(.53)\cong & 
  \Image\bigl(H_1(F) \rightarrow H_1(Y)\bigr),
}
\]
where $\gamma\mapsto [\gamma]$ denotes the map from the double-coset
space to $H_1(Y)$.
A more invariant version is
\[
\xymatrix{
  P_1(\x_1)\backslash \smallGroup\slash R(P_2(\x_2)) \ar[d]\ar[dr]^{\Pi} & \\
  \bigl(P_1(\x_1)\backslash \smallGroup\slash R(P_2(\x_2))\bigr)/\ZZ
  \ar[r]_(.44)\cong & 
  \bigl\{\spinc\in\SpinC(Y)\mathbin{\big|}
  \spinc\vert_{Y_i}\cong\spinc_i, i=1,2\bigr\},
}
\]
where the map~$\Pi$ is defined by
$
\gls*{PiOfGamma}\coloneqq\spinc(\x)+[\gamma].
$
We are implicitly using Poincar\'e duality to think of
$\SpinC(Y)$ as an affine space over $H_1(Y)$ (rather than the usual
$H^2(Y)$).

\glsadd{divis}%
Let $\divis(c_1(\spinc'))$ be the divisibility of the first
Chern class of $\spinc'$, i.e.,
\[
H^1(Y;\ZZ)\cup c_1(\spinc') =\divis(c_1(\spinc'))\cdot H^3(Y;\ZZ).
\]

\begin{lemma}\label{lem:gradings-tensor}
For $\s' \in \SpinC(Y)$, there is an isomorphism of $\ZZ$-sets $\Pi^{-1}(\s')\cong\ZZ/\!\divis(c_1(\spinc'))$.
\end{lemma}

\begin{proof}
Suppose that $\gamma\in G(\PMC)$ is a group element with the property that
$\gamma$ and $\lambda^t\cdot \gamma$ represent the same double coset. This means
that there are $B_i\in\pi_2(\x_i,\x_i)$ with the property that
\[\gamma = \lambda^t\cdot \gs(B_1)\cdot \gamma \cdot R(\gs(B_2)).\]
Let $\psi_0 = \psi(I_A(\x_1)) = \psi(I_D(\x_2))$.  Then, by
\eqref{eq:DSmallg} and~\eqref{eq:R-small-g},
\begin{align}
\gamma &= 
\lambda^t\cdot \psi_0 \gb(B_1) \psi_0^{-1} \cdot \gamma \cdot
  \psi_0 R(\gb(B_2)) \psi_0^{-1}\nonumber\\
\intertext{or, setting $\gamma' = \psi_0^{-1}\gamma\psi_0$,}
\gamma' &= \lambda^t\cdot \gb(B_1) \cdot \gamma' \cdot
  R(\gb(B_2)).\label{eq:tens-prod-lem-step}
\end{align}
Note that $[\gamma]$ and
$[\gamma']$ are the same element of $H_1(F)$.

This implies that $\bdy^\bdy B_1+\bdy^\bdy B_2=0$, so $B_1$ and $B_2$
glue to a domain
$B=B_1\glue B_2\in\pi_2(\x,\x)$ (Lemma~\ref{lem:HoClassFibProd}).
Using the formula $e(B)+2n_{\x}(B)=\langle c_1(\spinc(\x)),[B]\rangle$
(\cite[Proposition~7.4]{OS04:HolDiskProperties}),
Equation~\eqref{eq:tens-prod-lem-step} implies
\begin{align*}
1&=\lambda^t\cdot \lambda^{-e(B)-2n_{\x}(B)}(0,\bdy^\bdy B_1)\cdot \gamma'\cdot 
(0,-\bdy^\bdy B_1)\cdot \gamma'{}^{-1} \\
&=\lambda^{t-\langle
  c_1(\spinc(\x)),[B]\rangle -2[\gamma]\cap[\bdy^\bdy B_1]} \\
&=\lambda^{t-\langle
  c_1(\spinc(\x)),[B]\rangle -2[\gamma]\cap[B]}\\
&=\lambda^{t-\langle c_1(\Pi(\gamma)),[B]\rangle},
\end{align*}
where $[\gamma]\cap[\bdy^\bdy B_1]$ denotes the algebraic intersection
number of  $[\gamma]$ and $[\bdy^\bdy B_1]$ in $H_1(F)$, while
$[\gamma]\cap[B]$
denotes the algebraic intersection number of $[\gamma]\in H_1(Y)$ and
$[B]\in H_2(Y)$. The second equality uses 
Equation~\eqref{eq:commutators} and the last uses the identity
$c_1(\Pi(\gamma))=c_1(\spinc(\x))+2[\gamma]$.
From this (and the fact that any periodic domain $B$ on $Y$
can be decomposed into periodic domains $B_1$ and
$B_2$ (Lemma~\ref{lem:HoClassFibProd}), the result follows.
\end{proof}

\begin{theorem}
\label{thm:GradedPairing}
\index{pairing theorem!graded}\index{grading!for pairing theorem}%
Fix $\spinc_i\in\SpinC(Y_i)$.  There is a
homotopy equivalence
$$\Phi\co \CFAa(\HD_1,\spinc_1)\DT\CFDa(\HD_2,\spinc_2)
\to\!\! \bigoplus_{\substack{\spinc\in\SpinC(Y)\\
    \spinc|_{Y_i}=\spinc_i, i=1,2}}\!\! \CFa(\HD,\spinc)
$$
which respects the identification between grading sets in the following sense:
\begin{enumerate}
\item
\label{item:SpinCComponentOfDoubleCoset}
If  $m_1\DT m_2$ is a homogeneous 
element with grading 
$\grT(m_1\DT m_2)\in P_1(\x_1)\backslash \smallGroup\slash P_2(\x_2)$,
then $\Phi(m_1\DT m_2)$ lies in the summand of $\CFa(Y)$ with grading
$\Pi(\grT(m_1\DT m_2))\in\SpinC(Y)$.
\item 
\label{item:MaslovComponentOfDoubleCoset}
If $m_1\DT m_2$ and $n_1\DT n_2$ are
homogeneous elements whose gradings are related by
$$\grT(m_1\DT m_2)=\lambda^t \grT(n_1\DT n_2),$$ 
so that $\Pi(\grT(m_1\DT m_2))=\Pi(\grT(n_1\DT n_2))=\spinc'\in\SpinC(Y)$ 
then $\Phi(m_1\DT m_2)$ and $\Phi(n_1\DT n_2)$ are elements whose relative
Maslov grading is given by $t$ modulo $\divis(c_1(\spinc'))$.
\end{enumerate}
\end{theorem}

\begin{proof}
The map $\Phi$ is the homotopy equivalence from
Theorem~\ref{thm:TensorPairing}. We must prove that $\Phi$ respects
the grading sets.

To define the identification of grading sets we chose complementary
generators $\x_1$ and~$\x_2$ with $\s(\x_i)=\s_i$. (This is always
possible after modifying $\HD_1$ and $\HD_2$ by an isotopy, which changes
$\CFAa(\HD_1,\s_1)$ and $\CFDa(\HD_2,\s_2)$ by a graded homotopy
equivalence.)  Assume for
simplicity that $I_A(\x_1) = I_D(\x_2)$ is~$\SetS_0$, the base
idempotent used to define the grading refinement.  (Again, we may need
to modify $\HD_1$ and~$\HD_2$ by an isotopy to achieve this.)
In this case, some formulas simplify: $P(\x_i) =
P'(\x_i)$, and, for $\y_i \in \Gen(\HD_i,\s_i)$,
\begin{equation*}
\gr(\y_1) = \grb(\y_1) \cdot \psi(I_A(\y_1))^{-1}\qquad\text{and}\qquad
\gr(\y_2) =  \psi(I_D(\y_2)) \cdot \grb(\y_2).
\end{equation*}

Let $\y_1\in\S(\HD_1)$ and $\y_2\in\S(\HD_2)$ be generators with
$\s(\y_i)=\s_i$ and such that
$\y=\y_1\times\y_2$ is a generator for $\HD$.
Let $B_i\in\pi_2(\x_i,\y_i)$.
Then
\begin{align*}
\gr(\y_1)&=P_1(\x_1)\cdot\gs(B_1) \in \ASmallGrSet(\HD_1)=P_1(\x_1)\backslash\smallGroup (\PMC)\\
\gr(\y_2)&=R(\gs(B_2))\cdot R(P_2(\x_2)) \in\DSmallGrSet(\HD_2)=\smallGroup(\PMC)/R(P_2(\x_2)).
\end{align*}
Consequently,
\begin{align*}
\gr(\y_1\times\y_2)&=
\grb(\y_1)\cdot \psi(I_A(\y_1))^{-1}\cdot \psi(I_D(\y_2))\cdot\grb(\y_2) \\
&=
\grb(\y_1)\cdot\grb(\y_2),
\end{align*}
where we have used
$I_A(\y_1)=I_D(\y_2)$.
We wish to use this to calculate the
difference element between the two $\SpinC$ structures associated to
$\x=\x_1\times\x_2$ and $\y=\y_1\times\y_2$. The $\SpinC$ component
of $\gr'(\y_1)\cdot \gr'(\y_2)$ is $\bdy^\bdy B_1 + \bdy^\bdy B_2$.
This represents a homology class $\gamma$ in $H_1(F)$. Moreover, the image
of $\gamma$ in $H_1(Y)$ agrees with the difference element
\glsadd{epsilonxy}%
$\epsilon(\x_1\times\x_2, \y_1\times \y_2)$ defined
in~\cite[Section~2.6]{OS04:HolomorphicDisks} (see also Lemma~\ref{lem:SpinCStructures}), whose Poincar{\'e} dual
is $\spinc_z(\x_1\times \x_2)-\spinc_z(\y_1\times \y_2)$.
It follows  that
$\Pi(\grb(\y_1)\cdot\grb(\y_2))$ represents the same $\SpinC$
structure as $\y=\y_1\times\y_2$.  This verifies
Part~\ref{item:SpinCComponentOfDoubleCoset}.

For Part~\ref{item:MaslovComponentOfDoubleCoset}, 
let $\y'_1\in\S(\HD_1)$ and $\y'_2\in\S(\HD_2)$ be another pair of generators with
$\s(\y_i)=\s_i$ and such that
$\y'=\y'_1\times\y'_2$ is a generator for~$\HD$ in the same $\SpinC$
structure as~$\y$.
Fix some
$B'\in\pi_2(\y,\y')$, where $B'=B'_1\glue B'_2$.  Thus
$B_i*B_i'\in\pi_2(\x_i,\y_i')$. 
By
Lemma~\ref{lem:gB-mult} and the fact that $R$ is an anti-homomorphism, the grading of $\y'$ is
given by
\begin{align*}
\gr(\y') &= P_1(\x_1)\cdot\gb(B_1*B_1') R(\gb(B_2*B_2'))\cdot P_2(\x_2) \\
&=P_1(\x_1)\cdot\gb(B_1)\gb(B_1')R(\gb(B_2'))R(\gb(B_2))
  \cdot P_2(\x_2)\\
&=P_1(\x_1)\cdot\gb(B_1)\lambda^{-e(B')-n_{\y}(B')-n_{\y'}(B')}
   R(\gb(B_2))\cdot P_2(\x_2) \\
&=\lambda^{-e(B')-n_{\y}(B')-n_{\y'}(B')}\cdot\gr(\y).
\end{align*}
This proves Part~\ref{item:MaslovComponentOfDoubleCoset}.
\end{proof}

\begin{remark}
  To get the correct grading set on the tensor product, we need to use
  the refined gradings: the double-coset space $P'_1(\x_1) \backslash \bigGroup(4k) /
  R(P'_2(X_2))$ is larger than the double-coset space $P_1(\x_1) \backslash \smallGroup(\PMC) /
  R(P_2(X_2))$ in Theorem~\ref{thm:GradedPairing}.
\end{remark}

\begin{remark}
  An alternative to the refined grading would be to consider a
  ``grading groupoid'', with objects corresponding to the different
  idempotents $I(\SetS)$ and $\Mor(\SetS_1, \SetS_2)$ being those
  elements $g \in \bigGroup(4k)$ with $M_* \bdy([g]) = \SetS_1 -
  \SetS_2$.  Then $\ASmallGrSet$ and $\DSmallGrSet$ are equipped with
  right and left actions of this grading groupoid, and the grading set
  for the tensor product is the (suitably defined) product over the
  grading groupoid.
\end{remark}


\newcommand\eG{\widetilde G}
\chapter{Bordered manifolds with torus boundary}
\label{chap:TorusBoundary}

In this chapter, we specialize to the case of bordered manifolds with
torus boundary.  In particular, we show that the bordered Floer
invariant of a manifold with torus boundary is closely related to its
knot Floer homology.  As a warm-up, after introducing some notation in
Section~\ref{sec:torus-algebra}, we use the pairing theorem to give a
quick proof of a version of the surgery exact triangle for Heegaard
Floer homology, in Section~\ref{sec:surg-exact-triangle}. We then turn
to the relationship between bordered Floer homology for manifolds with
torus boundary and knot Floer homology. In Section~\ref{sec:PrelimHFK}
we recall the conventions on knot Floer homology. The easier task of
extracting (associated graded) knot Floer homology from bordered Floer
homology is addressed in Section~\ref{sec:CFD-to-HFK}. This follows
quickly from a suitable adaptation of the pairing theorem and a simple
model calculation for the solid torus.

Following this, Section~\ref{sec:CFKm-to-CFDa-statement} explains
how to extract the type $D$ module for a knot complement in $S^3$ from
the (filtered) knot Floer homology.  After further developing the holomorphic
curve machinery in Section~\ref{sec:BoundaryDegenerations}, this
result is proved in Sections~\ref{sec:CFK-to-CFD}
and~\ref{sec:HFKtoHFDproof}.  Finally, in
Section~\ref{sec:CablesAgain}, we use bordered Floer homology to
study satellite knots. 
(The knot Floer homology groups of various families of satellite
knots have also been studied using a more direct
analysis.  See~\cite{Eftekhary05:LongitudeWhitehead},~\cite{HeddenWhitehead},
and~\cite{Hedden}.)

\section{Torus algebra}
\label{sec:torus-algebra}

\begin{figure}
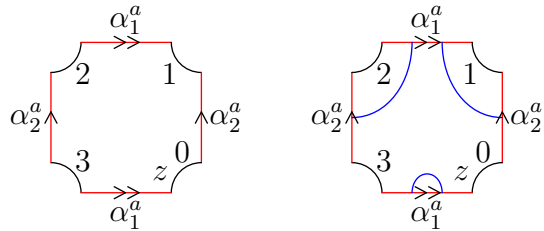

  \[
  \mfigb{torus-0}\qquad \mfigb{torus-19}
  \]
  \caption[Bordered Heegaard diagrams with torus boundary]{\textbf{Bordered Heegaard diagrams with torus boundary.}  Left: the
    necessary arrangement of the $\alpha$-arcs.  Opposite sides are
  identified as indicated.  There may be additional handles attached
  with corresponding $\alpha$-circles. Right: a $\beta$-circle to
  create an admissible diagram for a solid torus.}
  \label{fig:torus-diagrams}
\end{figure}

\begin{figure}
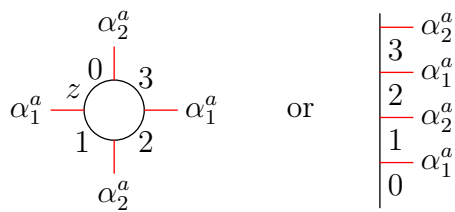

\[
\mfigb{torus-1}\qquad\text{or}\qquad\mfigb{torus-2}
\]
\caption[Labeling of regions around
    $\bdy\widebar{\Sigma}$ (torus case)]{\textbf{Labeling of the regions around
    $\bdy\widebar{\Sigma}$.} With respect to the boundary orientation,
they are labeled $0,3,2,1$ with $z$ in the region labeled $0$.}\label{fig:torus-region-labeling}
\end{figure}

A genus $g$ bordered Heegaard diagram~$\HD$ for a 3-manifold $Y$ with
torus boundary has two $\alpha$-arcs, $g-1$ $\alpha$-circles, and $g$
$\beta$-circles.  The $\alpha$-arcs are necessarily arranged as in
Figure~\ref{fig:torus-diagrams}.  As a result, the parametrization of
$\bdy Y$ is specified by a pair of a meridian and a longitude,
intersecting once, and the matching on $\bdy\HD$ is always the
same. Let 
$
\gls*{AlgTorus}
=\Alg(\Torus,0)$ denote the $i=0$ part of the
corresponding algebra.  Since we will mostly be interested in
$\CFDa(\HD$), we label the regions between the arcs in the opposite
order to the induced orientation on the boundary, as in
Figure~\ref{fig:torus-region-labeling}. (From the type $A$ side, the
regions would be labeled in the opposite order around the puncture.)
The point $z$ is in the region labeled~$0$.  Label the $\alpha$-arcs
so that regions $0$ and $1$ are separated from regions $2$ and $3$ by
$\alpha_2^a$.\glsadd{torusalphaarcs}

The algebra $\Alg(\Torus)$ has
two minimal idempotents $\iota_0$ and $\iota_1$,\glsadd{torusidempotents}
corresponding respectively to
$\alpha_1^a$ and
$\alpha_2^a$ being occupied on the type $A$ side, and 6 other basic generators,
denoted graphically as follows: \glsadd{torusreebgens}
\begin{align*}
  \rho_1&:\mfigb{torus-112}  & \rho_2&:\mfigb{torus-123} &
  \rho_3&:\mfigb{torus-134} \\[5pt]
 \rho_{12}&:\mfigb{torus-113} &
  \rho_{23}&:\mfigb{torus-124} &
  \rho_{123}&:\mfigb{torus-114}
\end{align*}
The differential is zero, and the non-zero products are\index{torus!algebra}
\[
  \rho_1\rho_2 = \rho_{12} \qquad \rho_2\rho_3 = \rho_{23} \qquad
  \rho_1\rho_{23} = \rho_{123} \qquad \rho_{12}\rho_{3} = \rho_{123}.
\]
(All other products of two non-idempotent basic generators vanish identically.)  There
are also compatibility conditions with the idempotents:
\begin{align*}
  \rho_1&=\iota_0\rho_1 \iota_1&
  \rho_2&=\iota_1\rho_2 \iota_0&
  \rho_3&=\iota_0\rho_3 \iota_1 \\
  \rho_{12}&=\iota_0\rho_{12} \iota_0&
  \rho_{23}&=\iota_1\rho_{23} \iota_1&
  \rho_{123}&=\iota_0\rho_{123} \iota_1.
\end{align*}

The unrefined grading takes values in the group $\bigGroup(\Torus)$,
which we denote $\bigGroup$\glsadd{torusbiggrgp} in this chapter.  It is
generated by
quadruples $(j;a,b,c)$ with $j \in \OneHalf\ZZ$, $a,b,c \in \ZZ$ and
$j$ is an integer if all of $a,b,c$ are even or if $b$ is
even and $a$ and $c$ are odd, and is a half-integer otherwise.
The group
law is\glsadd{torusbiggrgpelt}\index{torus!unrefined grading group}
\begin{multline*}
(j_1; a_1,b_1,c_1) \cdot (j_2; a_2,b_2,c_2)
= \biggl(j_1+j_2+
\frac{1}{2}
\begin{vmatrix}
a_1 & b_1 \\
a_2 & b_2 
\end{vmatrix}
+ 
\frac{1}{2}
\begin{vmatrix}
b_1 & c_1 \\
b_2 & c_2 
\end{vmatrix};\\ a_1+a_2, b_1+b_2, c_1+c_2\biggr).
\end{multline*}
The $\bigGroup$-grading on the algebra is given by
\begin{align*}
  \grb(\rho_1)&=\left(-\textstyle\frac{1}{2};1,0,0\right) \\
  \grb(\rho_2)&=\left(-\textstyle\frac{1}{2};0,1,0\right) \\
  \grb(\rho_3)&=\left(-\textstyle\frac{1}{2};0,0,1\right).
\end{align*}
The distinguished central element is
$\lambda = (1;0,0,0)$.

There is also a refined grading, as in
Section~\ref{sec:refined-grading}, with values in the
group $G(\Torus)$.  We will instead give a grading with values in a somewhat
larger group~$\smallGroup$\glsadd{torussmallgrgp} given explicitly by triples
$(j;p,q)$ where $j,p,q\in \OneHalf\ZZ$ and $p+q \in \ZZ$.
The group law\glsadd{torussmallgrgpelt}\index{torus!refined grading group}
is
$$(j_1; p_1,q_1) \cdot (j_2; p_2,q_2)
= \biggl(j_1+j_2+
\begin{vmatrix}
p_1 & q_1 \\
p_2 & q_2 
\end{vmatrix}; p_1+p_2, q_1+q_2\biggr)$$
and the distinguished central element is $\lambda=(1;0,0)$.
(The group $G(\Torus)$ from Section~\ref{sec:refined-grading} would
have $p,q \in \ZZ$ and $2j \equiv p+q \pmod{2}$.)
Let $G'_\QQ\supset\bigGroup$ be the group defined in the same way as
$\bigGroup$ but
where $j,a,b,c\in\QQ$ (so $G'_\QQ$ is a $\QQ$-central extension of
$H_1(Z,\mathbf{a};\QQ)$).
There is a homomorphism $G\to G'_\QQ$ defined by $(j;p,q) \mapsto (j;p,p+q,q)$.
To define the refined grading, choose $\iota_0$ as
the base idempotent and set $\psi(\iota_1) \coloneqq
\bigl(0;\OneHalf,0,\OneHalf\bigr)\in G'_\QQ$.
From
Equation~\eqref{eq:small-grading} we then find
\begin{equation}\label{eq:grading-torus-alg}
   \begin{split}
  \gr(\rho_1) &= \bigl(-\OneHalf;\OneHalf,-\OneHalf\bigr)\\
  \gr(\rho_2) &= \bigl(-\OneHalf;\OneHalf,\OneHalf\bigr)\\
  \gr(\rho_3) &= \bigl(-\OneHalf;-\OneHalf,\OneHalf\bigr).
   \end{split}
\end{equation}
Note that for any algebra element $a$, the gradings $\grb(a)$ and
$\gr(a)$ are related by the homomorphism from $\bigGroup$ to
$\smallGroup$ defined by
\begin{equation}
  \label{eq:grading-relation}
  (j;a,b,c)\longmapsto \biggl(j;\frac{a+b-c}{2},\frac{-a+b+c}{2}\biggr).
\end{equation}

We find it convenient to have the following concrete reformulation of
type $D$ structures over the torus algebra. (This will be helpful in
Section~\ref{sec:CFK-to-CFD}, where we understand how to construct the
type $D$ module of a knot complement from the knot Floer complex.)

\begin{definition}
  \label{def:CoeffMaps}
  Let $V=V^0\oplus V^1$ be a $\Zmod{2}$-graded vector space.
  A collection of (torus) \emph{coefficient maps}
  \index{coefficient maps}%
  \index{module!type $D$!coefficient maps|see{coefficient maps}}%
  \index{torus!coefficient maps|see{coefficient maps}}%
  is a collection of maps
  $$
  \gls*{coeffmap}
  \co V^{[i_0-1]} \to V^{[i_n]}
  $$
  indexed by
  increasing sequences of consecutive integers $I=\{i_0,\cdots,
  i_n\}\subset \{1,2,3\}$, where $[i]$ denotes the reduction of
  $i\pmod{2}$, satisfying the compatibility equation stated below.
  (The empty sequence also has a corresponding map $D_{\emptyset}$ or just~$D$,
  which maps $V^i$ to $V^i$ for
  $i=0,1$.) For any increasing sequence of consecutive integers $I\subset
  \{1,2,3\}$ we require that
  \begin{equation}
    \label{eq:Compatibility}
    \sum_{\{I=J\cup K\mid J<K\}} \!\!D_{K}\circ D_{J}=0,
  \end{equation}
  where
  $J<K$ means that each element of $J$ is smaller than every element
  of~$K$.
\end{definition}

Explicitly, Equation~\eqref{eq:Compatibility} in the case where
$I=\emptyset$ says that $D$ is a differential; similarly,
$I=\{1\}$ says that $D_1$ is a
chain map with respect to the differential $D$; 
$I=\{1,2\}$ says that
$D_{12}$ is a null-homotopy of the composite $D_2\circ
D_1$;
$I=\{123\}$ gives the following:
$$D \circ D_{123}+D_3\circ D_{12}+ D_{23}\circ D_1 + D_{123}\circ D = 0.$$

\begin{lemma}
  A type~$D$ structure in the sense of Definition~\ref{def:TypeD}
  over $\Alg(\Torus)$ with base ring $\Idem(\Torus)$ corresponds
  to a pair of $\Zmod{2}$-graded vector spaces equipped with coefficient maps 
  in the sense of Definition~\ref{def:CoeffMaps}.
\end{lemma}

\begin{proof}
  Let $V$ be a vector space underlying a type~$D$ structure over
  $\Alg(\Torus)$.  The projection
  $V\to V^0$ is induced by multiplication by the
  idempotent $\iota_0$, while the projection $V\to V^1$ is
  induced by multiplication by the idempotent $\iota_1$.  The
  coefficient maps can be extracted from the map $\delta^1$
  by the formula
  $$\delta^1 = \Unit \otimes D + \sum_{i} \rho_i\otimes D_i +
  \!\!\sum_{\{i,j\mid j=i+1\}}\!\!\rho_{ij} \otimes D_{i j} +
  \rho_{123} \otimes D_{123}.$$
  Equation~\eqref{eq:Compatibility} gives the components of the
  compatibility condition for $\delta^1$.
\end{proof}

In particular, thinking of $\CFDa(Y)$ as a type $D$ structure, as in
Remark~\ref{rmk:UnderlyingTypeD}, we can talk about its coefficient
maps.

\begin{remark}
  The coefficient maps depend on the type $D$ structure itself, not just
  on its corresponding differential module.
\end{remark}

\section{Surgery exact triangle}
\label{sec:surg-exact-triangle}
Recall that Heegaard Floer homology admits a surgery exact
triangle~\cite{OS04:HolDiskProperties}.  Specifically, for a pair
$(M,K)$ of a 3-manifold~$M$ and a framed knot~$K$ in $M$, there
is an exact triangle\index{surgery exact triangle}\index{exact triangle, surgery}
\begin{equation}
  \label{eq:surgery-exact-triangle}
\mathcenter{\begin{tikzpicture}[x=2.3cm,y=48pt]
  \node at (0,0) (m0) {$\HFa(M_{-1})$} ;
  \node at (2,0) (m1) {$\HFa(M_0)$} ;
  \node at (1,-1) (minf) {$\HFa(M_\infty)$} ;
  \draw[->] (m0) to (m1) ;
  \draw[->] (m1) to (minf) ;
  \draw[->] (minf) to node[auto] {} (m0) ;
\end{tikzpicture}}
\end{equation}
where $M_{-1}$, $M_0$, and $M_\infty$ are $-1$, $0$, and
infinity surgery on $K$,  respectively.  As a simple application of
bordered Floer theory, we prove a version of this result.
A similar computation was carried out in
\cite[Section 5.3]{Lipshitz06:BorderedHF}.
We do not, however,
verify that the maps in our present triangle agree with those
in the original version; see~\cite{LOT:DCov2}
or~\cite{LOTCobordisms}.

We prove this exact triangle by constructing three provincially
admissible bordered Heegaard diagrams $\HD_{-1}$, $\HD_0$, and
$\HD_\infty$ for
\index{solid torus, Heegaard diagram for}%
\index{Heegaard diagram!bordered!for solid torus}%
the three solid tori filled at the corresponding slopes, and
exhibiting a short exact sequence relating the chain complexes
$\CFDa(\HD_\bullet)$.  For any knot complement with (admissible) bordered
Heegaard diagram $\HD'$, we can then take the derived tensor product
of $\CFAa(\HD')$ with this short exact sequence.
The exact triangle in Equation~\eqref{eq:surgery-exact-triangle}
is the induced exact triangle on homology (see Proposition~\ref{prop:InducedLongExactSequence}).

The three diagrams we use are
\begin{equation}
  \label{eq:torus-diagrams}
  \HD_\infty:\mfigb{torus-30}\qquad
  \HD_{-1}:\mfigb{torus-20} \qquad
  \HD_0:\mfigb{torus-10}
\end{equation}%
\glsadd{HDSolidTori}%
A generator for $\CFDa(\HD_\bullet)$ consists of a single
intersection point between the $\beta$-circle in $\HD_\bullet$ and an
$\alpha$-arc. These intersections are labeled above.

The boundary operators on the $\CFDa(\HD_\bullet)$ (and the support
of the relevant domains) are given by
\newcommand{\tdiag}[2]{\underset{\includegraphics[scale=0.5]{torus-#2}}{\strut #1}}
\begin{equation*}
  \partial r = \tdiag{\rho_{23}r}{31}
  \qquad\qquad
  \partial a = \tdiag{\rho_3b}{21} + \tdiag{\rho_1b}{22}
  \qquad\qquad
  \partial b = \tdiag{0}{255}
  \qquad\qquad
  \partial n = \tdiag{\rho_{12}n}{11}
\end{equation*}
There is a short exact sequence
\[
0 \longrightarrow \CFDa(\HD_\infty)\overset{\varphi}{\longrightarrow} \CFDa(\HD_{-1})
  \overset{\psi}{\longrightarrow} \CFDa(\HD_0) \longrightarrow 0
\]
where the maps $\varphi$ and $\psi$ are given by
\begin{equation*}
  \varphi(r) = \tdiag{b}{41} + \tdiag{\rho_{2}a}{42} \qquad\qquad
  \begin{aligned}[t]
    \psi(a) &= \tdiag{n}{46}\\
    \psi(b) &= \tdiag{\rho_{2}n}{47}.
  \end{aligned}
\end{equation*}
Indeed, one can guess the maps $\varphi$ and $\psi$ by considering
holomorphic triangles, as indicated above. It is
straightforward to verify that $\psi \circ \varphi = 0$, that
$\varphi$ is injective, that $\psi$ is surjective, and that $\ker \psi
= \im \varphi$.

\section{Preliminaries on knot Floer homology}
\label{sec:PrelimHFK}

We focus now on the relationship between Heegaard Floer invariants of
three-manifolds with torus boundary and knot Floer homology. Before
delving into the details, we recall some of the basics of knot Floer
homology, with a special emphasis on the bigradings.  For simplicity
of notation, we restrict attention to the case where the ambient
three-manifold is~$S^3$.  More detailed accounts may be found in the
original papers~\cite{OS04:Knots, Rasmussen03:Knots}.

The knot Floer invariant of a knot is the filtered chain
homotopy type of a filtered
chain complex over $\Field[\gls*{U}]$. Specifically, a knot in $S^3$ can be
specified by a suitable Heegaard diagram
\index{Heegaard diagram!doubly pointed for knot}%
$\HD_K=\gls*{DblyPtdHD}
$, where $w$ and $z$ are a pair of
basepoints in the complement of the $\alpha$- and $\beta$-circles.
The \emph{knot Floer complex}
\index{knot Floer!complex}%
is defined as a chain complex 
$\gls*{knotCxminus}$
generated over $\Field[U]$ by 
$\gls*{knotGen}$,
the usual generators of Heegaard
Floer homology with respect to the Heegaard diagram $\HD_K$, endowed
with the differential
\begin{equation}
  \label{eq:DiffCFKm}
\partial^-(\x)\coloneqq\sum_{\y\in\Gen_K}\,
  \sum_{\substack{B\in\piBig(\x,\y)\\ \Mas(B)=1}}
\#\bigl(\Mod^B(\x,\y)\bigr) U^{n_{w}(B)}\cdot \y.
\end{equation}
In the above equation, 
$\gls*{pitwoBig}$
refers to the set of homology
classes of curves connecting $\x$ to $\y$
which are allowed to cross both $w$ and $z$; this is denoted
$\pi_2(\x,\y)$ in~\cite{OS04:Knots, Rasmussen03:Knots},
but in this book we have reserved $\pi_2$ for classes which have local multiplicity
zero at $z$.

The complex is endowed with an \emph{Alexander filtration}.
\index{Alexander filtration}\index{filtration, Alexander}%
The Alexander depth (filtration degree)~$\gls*{Alexander}$
of the generators is characterized, up to an
overall translation, by
\begin{equation}
  \label{eq:AlexanderUpToTranslation}
  \begin{aligned}
  A(\x)-A(\y)&=n_{z}(B)-n_{w}(B)\\
  A(U\cdot\x) &= A(\x) - 1
  \end{aligned}
\end{equation}
for any $\x,\y\in\Gen_K$ and homology class $B\in\piBig(\x,\y)$.
The boundary operators do not increase~$A$.
Furthermore, the boundary operators decrease (by one) the \emph{Maslov
  grading}~$\gls*{Maslov}$,
\index{grading!Maslov!on $\CFK$}%
defined up to overall translation by
\begin{equation}
  \label{eq:MaslovGrading} 
\begin{aligned}
    M(\x)-M(\y)&=\Mas(B)-2n_{w}(B) \\
    M(U\cdot\x) &= M(\x)-2. 
\end{aligned}
\end{equation}
(Note that the knot enters the definition only through the Alexander filtration.)
We will explain how to remove the indeterminacy in the Alexander and Maslov gradings 
below, in Equations~\eqref{eq:AlexanderSymmetry} and 
in Equations~\eqref{eq:RenormalizedMaslov}
and~\eqref{eq:FixAlgebraicGrading}, respectively. 

Sometimes it is natural to consider a different grading~$\gls*{Naslov}$,
the 
{\em $z$-normalized Maslov grading},  
\index{z-normalized Maslov grading@$z$-normalized Maslov grading}%
\index{grading!Maslov!z-normalized@$z$-normalized}%
characterized up to an overall translation, by
\begin{equation}
\begin{aligned}
  N(\x)-N(\y)&=\ind(B)-2n_{z}(B) \\
  N(U\cdot\x) &= N(\x).
\end{aligned}
\end{equation}
Any two of $A$, $M$, and $N$ determine the remaining one, by the
relation
\begin{equation}
  \label{eq:RenormalizedMaslov}
  N=M-2A.
\end{equation}

Let
$\gls*{gknotCxminus}$
denote the associated graded object with respect to the Alexander filtration.
Explicitly, this is a bigraded complex with a differential
$$\partial^-(\x)\coloneqq\sum_{\y\in\Gen_K}\,
  \sum_{\substack{B\in\piBig(\x,\y)\\ \Mas(B)=1, n_z(B)=0}}
\#\bigl(\Mod^B(\x,\y)\bigr) U^{n_{w}(B)}\cdot \y.$$
The functions $A$ and $M$ induce Alexander and Maslov gradings on $\gCFKm(K)$; 
the Alexander grading is preserved by the differential, and the Maslov grading is decreased by
one by the differential.

The homology groups of $\gCFKm(K)$ are the {\em knot Floer homology
  groups}\index{knot Floer!homology}
of~$K$,
$$
\gls*{knotHminus}
=\bigoplus_{r,d}\HFKm_d(K,r).$$
The index $r$ is the Alexander grading 
and $d$ is the Maslov grading.

The above constructions can be specialized to $U=0$, where we have a
different version of knot Floer homology. The $U=0$ specialization of
$\CFKm(K)$ is a filtered chain complex over $\Field$ denoted
$\gls*{knotCxHat}$.
We denote the associated graded complex to $\CFKa(K)$ by $\gCFKa(K)$.
Explicitly, 
$\gls*{gknotCxHat}$
is a bigraded complex 
generated over $\Field$ by $\Gen_K$, and endowed with the differential
$${\widehat\partial}(\x)\coloneqq\sum_{\y\in\Gen_K}\,
  \sum_{\substack{B\in\piBig(\x,\y)\\ \Mas(B)=1, n_z(B)=n_w(B)=0}}
\#\bigl(\Mod^B(\x,\y)\bigr) \y.$$ 
(Sometimes, when we wish to call attention to a fixed doubly-pointed Heegaard diagram, we denote
$\gCFKa(K)$ by $\gCFKa(\HD,w,z)$.)
The homology of $\gCFKa(K)$ is denoted
$$
\gls*{knotHhat}
=\bigoplus_{r,d}\HFKa_d(K,r).
$$

The indeterminacy (up to an overall translation) of the Alexander grading can be 
removed by requiring that the function $A$ have the following symmetry property
for all $r\in\ZZ$:
\begin{equation}
  \label{eq:AlexanderSymmetry}
  \#\{\x\in\Gen_K \mid A(\x)=r \} \equiv \#\{\x\in\Gen_K \mid A(\x)=-r\}\pmod{2}.
\end{equation}
The graded Euler characteristic of
$\HFKa$ coincides with the (symmetrized) Alexander polynomial $\Delta_K$
in the sense that
$$
\sum_{d,r} (-1)^d \rank \HFKa_d(K,r) \cdot T^r \doteq \Delta_K(T), 
$$
where $\doteq$ indicates that the two Laurent polynomials agree up to multiplication by
$\pm T^k$ (for some $k$). 
The symmetry property from Equation~\eqref{eq:AlexanderSymmetry}, which 
pins down the indeterminacy in the Alexander grading, can 
be reformulated as the condition that
$$\sum_{d,r} (-1)^d \rank \HFKa_d(K,r) \cdot T^r = \pm \Delta_K(T). $$

The indeterminacy in the Maslov grading
can be removed by first removing the indeterminacy in the $z$-normalized Maslov grading,
and demanding that Equation~\eqref{eq:RenormalizedMaslov} continue to
hold. In turn, to pin down the algebraic grading $N$, consider the chain complex
$\gCFKm(K)/(U=1)$. Note that neither the Alexander nor the Maslov
gradings descend to this complex; but the $z$-normalized Maslov grading does. The
homology of the resulting complex has rank one. The algebraic grading
is then normalized by requiring that this homology is supported in
grading equal to zero; i.e.,
\begin{equation}
  \label{eq:FixAlgebraicGrading}
  H_*\left(\frac{\gCFKm(K)}{U=1}\right)= \begin{cases}
      \Field & {*=0} \\
      0 & {\text{otherwise,}}
\end{cases}
\end{equation}
where $*$ denotes the algebraic grading $N$.

It is more traditional to define the Maslov gradings in
terms of holomorphic disks crossing the basepoint $z$. This is
normalized by considering
$\CFa(S^3)$, defined as holomorphic disks which do not cross the
basepoint $w$ (i.e., set $U=0$ in the chain complex defined in
Equation~\eqref{eq:DiffCFKm} and ignore its Alexander filtration), and
require that the resulting homology is supported in
degree zero with respect to the Maslov grading. Our present
formulation defines gradings in terms of
disks which never cross $z$. The equivalence of the two points of view
follows from the fact that $\gCFKm(K)/(U=1)$ is
the chain complex $\CFa(S^3)$, using the
basepoint $z$ in place of $w$.

Knot Floer homology naturally gives rise to an
integral valued concordance invariant~$\gls*{tau}$
for
knots~\cite{OS03:4BallGenus}.  One
construction of~$\tau$ comes
from considering the filtered chain homotopy type of
$\CFKa(K)$. From this
point of view, $\tau(K)$ is the minimal $s$ for which the
generator of $H_*(\CFa(S^3))$ can be represented as a sum of
generators in Alexander grading less than or equal to
$s$. In an alternative formulation, $\tau$ has an
interpretation from the associated graded object $\HFKm(K)$,
namely,
\begin{equation}
  \label{eq:CharacterizeTau}
  \tau(K)=-\max\{s \mid \forall d\geq 0, U^{d}\cdot \HFKm(K,s)\neq 0\}.
\end{equation}
(See, for example, \cite[Lemma~A.2]{OST}.)

The above definition of the Alexander grading carries over with little
change to the more general case where $K\subset Y$ is a
null-homologous knot in an arbitrary closed, oriented three-manifold
(rather than just $S^3$). The Maslov grading is a little more subtle;
but that is because Maslov gradings on closed manifolds are more
complicated. Our present applications, however, primarily concern the case where $Y=S^3$.

\section{From \textalt{$\CFDa$}{CFD\textasciicircum} to \textalt{$\HFKm$}{HFK-}}
\label{sec:CFD-to-HFK}

Let $K$ be a null-homologous knot in a $3$-manifold $Y$. Recall from
Section~\ref{sec:Bestiary} that we can construct a bordered Heegaard
diagram for $Y\setminus \nbd(K)$ as follows. Fix a Heegaard diagram
$(\Sigma_g,\alpha_1^c,\dots,\alpha_{g-1}^c,\beta_1,\dots,\beta_g)$
for $Y\setminus \nbd(K)$, in the classical
sense. Let $\mu\subset \Sigma$ and $\lambda\subset \Sigma$ denote a
meridian and longitude for $K$, respectively. Arrange that $\lambda$
and $\mu$ intersect in a single point $p$, and that they are disjoint
from all $\alpha_i^c$. Set $\alpha_1^a \coloneqq \lambda\setminus\{p\}$
and $\alpha_2^a \coloneqq \mu\setminus\{p\}$. Then
\[
(\Sigma\setminus \{p\},\alpha_1^a,\alpha_2^a, \alpha_1^c,\dots,\alpha_{g-1}^c,\beta_1,\dots,\beta_g)
\]
is a bordered Heegaard diagram for $Y\setminus\nbd(K)$.

Number the four regions in $(\Sigma,\alphas,\betas)$ adjacent to $p$,
in counter-clockwise cyclic order (i.e., the opposite of the induced
boundary orientation), so that the meridian $\alpha_2^a$ separates
regions $0$ and~$3$..
Let $z$ and $w$ be basepoints in regions $0$ and $2$,
respectively. (See Figure~\ref{fig:Twisting}.) Then 
\[
(\Sigma,\mu,\alpha_1^c,\dots,\alpha_{g-1}^c,\beta_1,\dots,\beta_g)
\]
is a doubly-pointed Heegaard diagram for $K$ in $Y$.  So, with
notation as in Section~\ref{sec:torus-algebra}, for the Heegaard
diagram $(\Sigma,\alphas,\betas,z)$ it is immediate that
$(V^0,D)$ is the knot Floer chain complex $\gCFKa(Y,K)$. (The chain
complex $(V^1,D)$ is the \emph{longitude Floer
complex}\index{longitude Floer complex} studied by Eftekhary~\cite{Eftekhary05:LongitudeWhitehead}.)

It follows that for any bordered Heegaard diagram
$(\Sigma,\alphas,\betas,z)$ for $Y\setminus\nbd(K)$, with framing such
that $\alpha_2^a$ is a meridian of $K$, the complex $(V^0,D)$ is
homotopy equivalent to $\gCFKa(Y,K)$. To see this, define an $\Alg$-module
$W_0$ to be one-dimensional over~$\Field$ with
$\iota_0$ acting as the identity and
all other basic generators of $\Alg$ acting by $0$. Then the tensor product
$W_0\DT\CFDa(\Sigma,\alphas,\betas,z)$ is exactly $(V^0,D)$. By
Lemma~\ref{lem:ThetaComplex2}, the homotopy type of this
tensor product depends only on the homotopy type of
$\CFDa(\Sigma,\alphas,\betas,z)$. Thus, we have the following.
\begin{proposition}\label{prop:CFD-to-CFKhat}Let
  $(\HD,z)$ be any bordered Heegaard diagram for
  $Y\setminus\nbd(K)$, with framing such that $\alpha_2^a$ is a
  meridian of $K$. Then, with notation as in
  Section~\ref{sec:torus-algebra}, $(V^0,D)$ is homotopy equivalent to
  $\gCFKa(Y,K)$. In particular, $\gCFKa(Y,K)$ is determined by
  $\CFDa(\HD,z)$.
\end{proposition}

The relative Maslov grading on $\gCFKa(Y,K)$ within each $\SpinC$ structure (on the
knot complement) is induced
from the grading on $\CFDa(\HD,z)$ as follows.  Two generators represent the
same $\SpinC$-structure on $Y\setminus K$ if and only if they can be
connected by a provincial domain. Consequently, the gradings of two
such generators differ by a power of $\lambda$, and this
$\lambda$-power is their Maslov grading difference.
One can also determine the Alexander grading, the relative Maslov grading
between generators in different Alexander gradings and, in the case that
$Y=S^3$, the absolute Maslov grading; we return to these points later.

There is a simple geometric interpretation of $W_0$, which we explain
after a more general discussion.

A {\em doubly-pointed bordered Heegaard diagram}
\index{Heegaard diagram!bordered!doubly pointed}%
is a bordered
Heegaard diagram $(\Sigma,\alphas,\betas)$ for $(Y,\bdy Y)$, together
with basepoints $z$ and $w$ in $\bSigma\setminus(\balphas\cup\betas)$.
We further assume $z \in \bdy\bSigma$.
A doubly-pointed bordered Heegaard diagram specifies a knot $K$ in
$Y$: let $\gamma_\alpha$ (respectively $\gamma_\beta$) denote a path
in $\Sigma\setminus\alphas$ (respectively $\Sigma\setminus\betas$)
from $z$ to~$w$. Let $\tilde{\gamma}_\alpha$ denote the result of
pushing the interior of $\gamma_{\alpha}$ into the $\alpha$-handlebody
and $\tilde{\gamma}_\beta$ the result of pushing the interior of
$\gamma_{\beta}$ into the $\beta$-handlebody. Then
$K=\tilde{\gamma}_\alpha\cup\tilde{\gamma}_\beta$. We
orient $K$ as $\tilde{\gamma}_\alpha-\tilde{\gamma}_\beta$.

\begin{figure}
  \centering
  \includegraphics[scale=.83333]{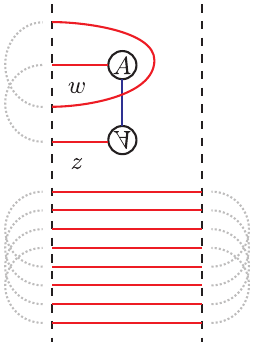}
  \caption[From bordered diagrams to doubly-pointed diagrams]{\textbf{Turning a bordered diagram for $Y'=Y\setminus\nbd(K\cup
      A)$ into a doubly-pointed bordered diagram for $(Y,K)$.} This
    diagram can be thought of as a bordered Heegaard diagram with two
    basepoints for the
    standard $1$-handle attachment cobordism $M$ from a surface of
    genus $g$ to a surface of genus $g+1$. Gluing the left boundary of
    this diagram, or an obvious analogue, to a bordered Heegaard
    diagram for $(Y',\bdy Y'=T^2\connectsum\Sigma_g)$, so that the top portion
    of this diagram is glued to the $T^2$ summand in $\bdy Y'$, gives a doubly-pointed
    Heegaard diagram for $(Y'\cup_\bdy M, K)$, where $K$ is the core
    of the $1$-handle in $M$.}
  \label{fig:bordered-to-2ptd}
\end{figure}
Every knot in a bordered three-manifold can be realized by a
doubly-pointed bordered Heegaard diagram. Specifically, let $Y$ be a
bordered three-manifold with genus $g$ boundary, and let $K\subset Y$
be a knot. By connecting $K$ to the boundary by an arc~$A$ and
deleting a regular neighborhood of $K\cup A$, we can write $Y$ as the
union of a three-manifold $Y'$ with genus $g+1$ boundary and a
standard cobordism $M$ (a two-handle attachment) from the surface of
genus $g+1$ to a surface with genus $g$. Gluing the (partial) diagram
shown in Figure~\ref{fig:bordered-to-2ptd} to a bordered Heegaard
diagram for $Y'$ gives a doubly-pointed bordered Heegaard diagram for $K$
in $Y$.

To a doubly-pointed bordered Heegaard diagram
$(\HD,z,w)$, we can associate type~$D$ and
type~$A$ modules as before, now working over a ground ring of
$\Field[U]$, with the understanding that holomorphic disks crossing
$w$ with multiplicity $n$ contribute a factor of $U^n$.  (Holomorphic
disks crossing~$z$ are still forbidden.)
We denote the
resulting modules by 
$\gls*{knotBorderedmA}$ and $\gls*{knotBorderedmD}$.
Setting $U=1$
recaptures $\CFAa(\HD,z)$ and
$\CFDa(\HD,z)$, while setting $U=0$ gives a
pair of modules denoted
$\gls*{knotBorderedaA}$ and $\gls*{knotBorderedaD}$
where we count only those
holomorphic disks with multiplicity zero at $w$.

\begin{remark}
  These new modules $\CFDm(\HD,z,w)$, etc., are not invariants of the
  pair $(Y,K)$.  Because the proofs of invariance in Sections
  \ref{sec:typeD-invariance} and~\ref{sec:A-invariance} require one
  basepoint, $z$, to remain on the boundary in the course of the
  handleslides and isotopies, the module $\CFDm(\HD,z,w)$ is only an
  invariant of the graph in $(Y,\bdy Y)$ obtained by attaching a point
  in~$K$ to $\bdy Y$; see
  also~\cite[Example~\ref*{LOT2:rmk:bdy-dehn-twist}]{LOT2}.
\end{remark}

For a doubly-pointed bordered Heegaard diagram with boundary~$\PMC$,
we can enhance the
grading set to include an Alexander grading.
\index{Alexander grading|see{grading, Alexander}}%
\index{grading!Alexander!bordered case}%
Specifically, the
enhanced grading group has the form 
$
\gls*{AGrading}
=\smallGroup(\PMC)\times
\ZZ$. (This is a direct product.) We call the new $\ZZ$ summand
the {\em Alexander factor}.\index{Alexander factor (of grading)}
Given $\x,\y\in\S(\HD)$ and
$B\in\piBig(\x,\y)$, define
$$
\gls*{Ags}
(B)\coloneqq (g(B),n_{w}(B)-n_{z}(B)),$$
where $g(B)$
is the grading from Equation~\eqref{eq:DSmallg}.  We can use this
enhanced grading group to define an enhanced grading on
$\CFD^-({\mathcal H},z,w)$. Following the discussion from
Section~\ref{sec:typeD-gradings}, a $\SpinC$ structure
$\spinc$ on $Y$  gives rise to a grading on the
type $D$ module with values in the coset space
$$
\gls*{AGradingD}
({\mathcal H},\spinc)\coloneqq (\AGrading/\APeriodics(\x_0)),$$
where $\x_0$ is any intersection point representing $\spinc$,
and 
$\gls*{APeriodics}$
denotes the image of the space of periodic domains
in $\AGrading$ under the homomorphism $\Ags$.  Namely, we
define 
$\gls*{Agr}(\x)\coloneqq \Ags(B)\cdot \APeriodics(\x_0)$, 
where $B\in\pi_2(\x,\x_0)$.
Note that when we extend this grading to 
$\CFDm({\mathcal H},z,w)$, we take the convention that multiplication by 
$U$ drops the Alexander factor of the grading by one (and preserves the 
other components of the grading).
\index{grading!multiplication by $U$ on}%
\glsadd{AgrUx}%

Suppose next that $({\mathcal H}_1,z)$ is a bordered Heegaard
diagram for $(Y_1,\partial Y_1=F)$ and that $({\mathcal H}_2,z,w)$ is a
doubly-pointed bordered Heegaard diagram for $(Y_2,\partial Y_2=-F)$,
equipped with the knot~$K$. Then, the grading set of the tensor
product $\CFAa(\mathcal{H}_1,z)\DT\CFDm(\mathcal{H}_2,z,w)$ is
naturally the double-coset space $\APeriodics(\x_1)\backslash
{\AGrading}/\APeriodics(\x_2)$. (In the definition of $\APeriodics(\x_1)$, we take the convention that $n_w(B)=0$ for any domain $B$
since
$w$ is thought of as located on the $Y_2$ side of the Heegaard diagram.)

Projecting onto the Alexander factor in the double-coset space induces
a map from the grading set to a transitive $\ZZ$-set.  
Indeed, the transitive $\ZZ$-set is isomorphic to
$\ZZ/n_{w}(\pi_2(\x,\x))$ for any generator $\x$ of $\HD_1 \cup_Z \HD_2$.
For any $B \in \pi_2(\x,\x)$ representing the homology class $A\in
H_2(Y_1 \cup Y_2)$, the quantity $n_w(B)$ is the intersection of~$K$ with~$A$.
In particular, if
$K$ is null-homologous, the $\ZZ$-set is (non-canonically) isomorphic to~$\ZZ$.
In this way,  when $Y$ is an integer homology sphere, the tensor product
$\CFAa(\HD_1,z)\DT \CFDm(\HD_2,z,w)$ inherits a (relative) $\ZZ$ grading,
which we call its {\em Alexander grading}. Corresponding remarks also apply to $\CFAm(\HD_2,z,w)\DT \CFDa(\HD_1,z)$.

We now have a pairing theorem for knot Floer homology.

\begin{theorem}
  \label{thm:PairingKnot}
  \index{pairing theorem!for knot Floer}%
  \index{knot Floer!pairing theorem}%
  Let $(\HD_1,z,w) = (\Sigma,\alphas_1,\betas_1,z,w)$ be a
  doubly-pointed bordered Heegaard
  diagram for $(Y_1,\partial Y_1=F,K)$ and
  $(\HD_2,z)=(\Sigma_2,\alphas_2,\betas_2,z)$ be a bordered
  diagram for $(Y_2,\partial Y_2=-F)$.
  Then, there are homotopy equivalences
  \begin{align*}
    \gCFKm(Y_1\cup_F Y_2,K)&\simeq \CFAm(\HD_1,w,z)\DT\CFDa(\HD_2,z)\\
    &\simeq \CFAa(\HD_2,z)\DT\CFDm(\HD_1,w,z) \\
    \gCFKa(Y_1\cup_F Y_2,K)&\simeq \CFAa(\HD_1,w,z)\DT\CFDa(\HD_2,z)\\
    &\simeq \CFAa(\HD_2,z)\DT\CFDa(\HD_1,w,z),
  \end{align*}
  which respect the gradings in an obvious sense. 
  In particular,
  when $K$ is null-homologous in $Y_1\cup_{F} Y_2$, 
  the above homotopy equivalences
  respect the relative $z$-normalized Maslov and Alexander gradings.
  (The Alexander and $z$-normalized Maslov gradings on
  the knot Floer homology appearing on the left-hand side
  are as recalled in Section~\ref{sec:PrelimHFK}; 
  the relative Alexander grading on  the right-hand side is
  defined immediately above, and the relative Maslov grading is induced by
  projecting onto the Maslov component of the double-coset space).
\end{theorem}

\begin{proof}
  The identification of complexes follows from an extension of
  either of the two
  proofs of Theorem~\ref{thm:TensorPairing} to the doubly-pointed
  context, without any significant change. The statement about
  Alexander gradings now follows
  from Lemma~\ref{lem:HoClassFibProd}.
\end{proof}

\begin{figure}
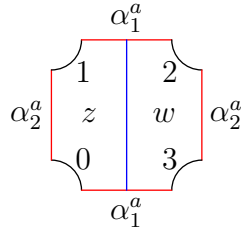

$\mfigb{torus-14}$
  \caption[Diagram for the longitude of a solid torus]{{\bf The
      doubly-pointed Heegaard diagram $\HDst$ for the longitude in
      the solid torus.} This is a genus~$1$ Heegaard diagram with
    boundary, with $\alpha$-arcs $\alpha_1^a$ and $\alpha_2^a$, and a
    $\beta$-circle $\beta_1$ intersecting $\alpha_1^a$ in a single
    point and disjoint from~$\alpha_2^a$. Note that the boundary markings
    follow the conventions of a type $A$ module.}
  \label{fig:Solid-Torus-with-Longitude}
\end{figure}

Proposition~\ref{prop:CFD-to-CFKhat} can be seen as a special case of
Theorem~\ref{thm:PairingKnot}, as follows.  Consider the
doubly-pointed Heegaard diagram $\HDst$ with
boundary for the solid torus shown in
Figure~\ref{fig:Solid-Torus-with-Longitude}.  If we restrict to
holomorphic curves which cross neither $z$ nor $w$, the resulting
type~$A$ module $\CFAa(\HDst)$ has a single generator over $\Field$, no differential,
and
indeed is isomorphic to the module $W_0$ described before the statement of
Proposition~\ref{prop:CFD-to-CFKhat}.  Thus,
Proposition~\ref{prop:CFD-to-CFKhat} follows from the $U=0$ version of
Theorem~\ref{thm:PairingKnot} (with the solid torus on the type $A$
side).

To reconstruct $\gCFKm$ of a knot from $\CFDa$, we must describe
$\CFAm(\HDst)$. But this is easily done.

\begin{lemma}
        \label{lemma:SolidTorus}\index{type $A$ invariant!for solid torus}
        The type $A$ module $\CFAm(\HDst)$ for $\HDst$, representing
        the solid torus
        equipped with a longitudinal unknot, has a single generator $\x_0$
        and higher products given by
\begin{align*}
m_{3+i}(\x_0,\rho_3,\overbrace{\rho_{23},\dots,\rho_{23}}^{i},\rho_2)
=U^{i+1}\cdot \x_0
\end{align*}
for all $i\geq 0$.
\end{lemma}

\begin{proof}
  The only curves which can contribute are disks covering the
  region containing~$w$ a total of $k$ times.  There is exactly one
  such disk for each~$k$, and it contributes the term in the statement
  with $k=i+1$.
\end{proof}

Thus, according to Theorem~\ref{thm:PairingKnot}, $\HFKm(Y,K)$ can be
calculated as a derived tensor product of $\CFDa(Y\setminus K)$ with
$\CFAm(\HDst)$, generalizing
Proposition~\ref{prop:CFD-to-CFKhat}.

In a different direction, Theorem~\ref{thm:PairingKnot} gives a
technique for studying the knot Floer homology of satellite knots.  We
return to this point in Section~\ref{sec:CablesAgain} (where we
give also some examples), contenting ourselves for now with
the following corollary.

\begin{corollary}\label{cor:gives-satellite}
  Suppose that $K_1$ and $K_2$ are knots in $S^3$ such
  that $\CFDa(S^3\setminus K_1)\simeq\CFDa(S^3\setminus K_2)$ (with
  respect to the $0$-framing on $S^3\setminus K_1$ and $S^3\setminus
  K_2$, say). Let $K_1^C$ (respectively $K_2^C$) denote the satellite
  of $K_1$ (respectively $K_2$) with pattern $C$. Then
  $\HFKm(K_1^C)\cong\HFKm(K_2^C)$.
\end{corollary}

\begin{proof}
  According to Theorem~\ref{thm:PairingKnot},
  the knot Floer homology of the satellite knot is obtained
  as the derived tensor product of the type $D$ module of the companion
  knot ($K_1$ or $K_2$) with the type $A$ module of a doubly-pointed
  Heegaard diagram for the pattern knot in the solid torus.
\end{proof}

In light of Theorem~\ref{THM:HFKTOHFD} (and
Proposition~\ref{prop:HotopyEquivalentDStructures}),
Corollary~\ref{cor:gives-satellite} implies that if $\CFKm(K_1)$ is
filtered homotopy equivalent to $\CFKm(K_2)$ then $\HFKm(K_1^C)$
is isomorphic to $\HFKm(K_2^C)$.

\section{From \textalt{$\CFKm$}{CFK-} to \textalt{$\CFDa$}{CFD\textasciicircum}: Statement of results}
\label{sec:CFKm-to-CFDa-statement}

The rest of this chapter concerns the opposite direction, computing $\CFDa(S^3\setminus
K)$ from $\CFKm(K)$. The main results are Theorems~\ref{THM:HFKTOHFD}
and~\ref{thm:HFKtoHFD2} below. The statement of Theorem~\ref{THM:HFKTOHFD} will 
refer to some general algebraic properties of the chain complex $\CFKm(K)$,
which is a 
free, $\ZZ$-filtered, $\ZZ$-graded chain complex over
$\Field[U]$; and in particular, it refers to certain convenient bases
for this chain complex. After setting up the needed algebra, we give the statement 
of the theorem.

\begin{definition}
  \index{basis!filtered}\index{filtered!basis}%
  \index{filtered!linearly-independent}%
  Let $\cdots\subset V_{i}\subset V_{i+1}\subset \cdots$ be an
  (ascending) filtration
  of a vector space $V$, so $V=\bigcup_{i=-\infty}^\infty V_i$. Assume
  that $\bigcap_{i=-\infty}^\infty V_i=0$. Given $v\in V$, define
  $A(v)\coloneqq\inf\{i\in\ZZ\mid v\in V_i\}$; we call $A(v)$ the \emph{filtration
    level of $v$}. 
  Let $\gr(V)$ denote the associated graded vector
  space, i.e., 
  \begin{align*}
  \gr(V) &\coloneqq \bigoplus\nolimits_{i=-\infty}^\infty \gr_i(V) &
    \gr_i(V) &\coloneqq V_i/V_{i-1}.
  \end{align*}
  Composing the projection $V_i\to V_i/V_{i-1}$ and the inclusion
  $V_i/V_{i-1}\to \gr(V)$ gives a map $\iota_i\co V_i\to \gr(V)$.
  For $v\in V$ define $[v]\coloneqq\iota_{A(v)}(v)$.

  A \emph{filtered basis} of a filtered vector space $V$ is a basis
  $\{v_1,\dots,v_n\}$ for $V$ such that $\{[v_1],\dots,[v_n]\}$ is a
  basis for $\gr(V)$. Similarly, we say $\{v_1,\dots,v_m\}$ is
  \emph{filtered linearly\hyp independent} if $\{[v_1],\dots,[v_m]\}$ is
  linearly independent.
\end{definition}
There is a filtered analogue of the basis extension theorem: if
$\{v_1,\dots,v_m\}$ is a filtered linearly-independent set then there is a
filtered basis for $V$ containing $\{v_1,\dots,v_m\}$.

We now look at filtered $\Field[U]$ modules in general.  Let $C$ be a
finitely-generated free $\Field[U]$
module equipped with
\begin{itemize}\index{complex!$\ZZ$-graded, $\ZZ$-filtered over $\Field[U]$}
\item an integer-valued
grading, which we refer to as the {\em Maslov
  grading}~$M$,\index{grading!Maslov}
\item an integer-valued ascending filtration, the {\em Alexander
    filtration}~$A$,\index{Alexander filtration}\index{filtration, Alexander} and
\item a differential, denoted $\partial$. 
\end{itemize}
Assume these data satisfy the following compatibility conditions:
\begin{itemize}
\item The differential $\partial$ drops the Maslov grading by one:
  for homogeneous $\xi\in C$, $M(\partial \xi) = M(\xi) - 1$.
\item The differential $\partial$ respects the Alexander
  filtration: $A(\partial\xi) \le A(\xi)$.
\item  Multiplication by $U$ is compatible with the Maslov grading, in
  the sense that $M(U\cdot \xi) = M(\xi)-2$.
\item  Multiplication by $U$ drops the Alexander filtration by at
  least~$1$, in the sense that $A(U\cdot \xi)
  \le A(\xi)-1$.
\end{itemize}

Given such a $C$, there are two naturally associated chain complexes.
One complex, which we call the {\em vertical
  complex}
\index{vertical!complex}\index{complex!vertical} 
$
\gls*{Cvert}
\coloneqq C/(U\cdot C)$, is a filtered complex which inherits
the Alexander filtration.
Its homology is called the {\em
  vertical homology},
\index{vertical!homology}%
denoted~$\gls*{Hvert}(C)$.
Next define
$C^{\infty} \coloneqq C\otimes_{\Field[U]} \Field[U,U^{-1}]$.  This has
two filtrations, namely the Alexander filtration and the filtration by
powers of $U^{-1}$.
(Precisely, the $s\th$ filtered part of the second filtration contains
generators $U^t \x$, where $t \ge -s$ and $\x$ is a generator of~$C$.)
Then the \emph{horizontal
  complex}~$\gls*{Chor}$ 
\index{horizontal!complex}\index{complex!horizontal}%
is
defined to be the degree~0 part of the associated graded complex
to~$C^{\infty}$ with
respect to the Alexander filtration.  (Note that $\Cvert$ is the
degree 0 part of the associated graded space to $C^{\infty}$ with respect to the
filtration by $U$ powers.)  The complex $\Chor$ still has the filtration by powers
of~$U^{-1}$ and
inherits a differential, whose homology
is called the {\em horizontal homology},\index{horizontal!homology}
denoted 
$\gls*{Hhor}(C)$.

We illustrate $C^{\infty}$ graphically by picking a filtered basis $\xi_i$
over $\Field[U]$ of $C$ and
plotting a generator
over $\Field$ of $C\otimes\Field[U,U^{-1}]$ of the form $U^{-x} \cdot
\xi_i$ with Alexander
depth~$y$ on the plane at the
position $(x,y)$. Then the
differential of a generator at $(x,y)$ can be graphically represented
by arrows
connecting the point at $(x,y)$ with the coordinates of other
generators.  These arrows necessarily point (non-strictly) to the left and
down.  The complex $\Cvert$ is obtained by restricting to one vertical
column, and $\Chor$ is obtained by restricting to one horizontal row. (See the left-hand side of Figure~\ref{fig:HFKtoHFD}.)

\begin{definition}
\label{def:Simplified}
We call $C$ {\em reduced}\index{complex!reduced}\index{reduced!complex}
 if $\partial$ strictly drops either the Alexander
or $U$ power filtration;
i.e., if for all $\x$, we have $\partial\x=\y_1 + U\cdot\y_2$ where $A(\y_1)<A(\x)$.

We call a filtered basis $\{\xi_i\}$ over $\Field[U]$ for $C$ {\em vertically simplified}
\index{vertically simplified basis}%
\index{basis!vertically simplified}%
if for
each basis vector~$\xi_i$, either $\partial \xi_i\in U\cdot C$ or
$\partial \xi_i\equiv \xi_{i+1} \pmod{U\cdot C}$. In the
latter case we say that there is a {\em vertical arrow from $\xi_i$
  to $\xi_{i+1}$}; so, in a vertically simplified basis, there is at
most one vertical arrow starting or ending at each basis
element.  (The fact that $\bdy^2=0$ says that there cannot be a vertical arrow
both starting and ending at a single element.)
\index{vertical!arrow}\index{vertical!arrow!length}%
The {\em length}\index{length!of arrow}\index{arrow!vertical}\index{arrow!length}
of the vertical arrow is the difference
$A(\xi_i)-A(\xi_{i+1})$.

Similarly, we call a filtered basis $\{\eta_i\}$ over $\Field[U]$ for $C$ {\em horizontally
  simplified}
\index{horizontally simplified basis}%
\index{basis!horizontally simplified}%
if for each basis vector~$\eta_i$, either
$A(\partial\eta_i)<A(\eta_i)$ or there is an $m$ so that
$\partial \eta_i= U^m\cdot \eta_{i+1} +
\eps$ where $A(\eta_i)=A(U^m\cdot \eta_{i+1})$ and
$A(\eps)<A(\eta_i)$.
In the latter case we say that there is a
\emph{horizontal arrow from $\eta_i$ to $\eta_{i+1}$}; so, in a
horizontally simplified basis, there is at
most one horizontal arrow starting or ending at each basis
element. 
\index{horizontal!arrow}\index{horizontal!arrow!length}%
\index{arrow!horizontal}\index{arrow!length}%
The {\em length}\index{length of arrow} of
the horizontal arrow is the integer~$m$.
\end{definition}

\begin{remark}
  These notions could alternatively be formulated without reference to
  a basis. Vertical (respectively horizontal) arrows correspond to non-vanishing
  differentials in the spectral sequence belonging to the filtered
  complex $\Cvert$ (respectively~$\Chor$).
\end{remark}

Every finitely generated chain complex $C$ as above is homotopy equivalent
to one which is reduced; moreover, for a reduced complex, we can find a vertically or horizontally simplified
basis.
See Proposition~\ref{prop:SimplifyComplex}.

In particular, these constructions apply to knot Floer homology
$\CFKm(K)$.  The complexes $\Cvert(\CFKm(K))$ and $\Chor(\CFKm(K))$
(ignoring the filtration) both represent $\CFa(S^3)$
and so if $K$ is a knot in $S^3$, then 
\[\Hvert(\CFKm(K))\cong
\Hhor(\CFKm(K))\cong\HFa(S^3)\cong\Field.\]  As such, after passing to 
a reduced model for $\CFKm(K)$, a vertically
(respectively horizontally) simplified basis $\xi_i$ (respectively $\eta_i$)
has a preferred element which
(after reordering) we label $\xi_0$ (respectively $\eta_0$), characterized
by the property that it has no in-coming
or out-going vertical (respectively horizontal) arrows.

Observe that
\begin{equation}
  \label{eq:GradingsOnXiEta}
  A(\xi_0)=\tau(K) \qquad M(\xi_0)=0 \qquad A(\eta_0)=-\tau(K) \qquad
  M(\eta_0)=-2\tau(K).
\end{equation}
The first equation follows from the definition of $\tau(K)$; the
second from the fact that $\HFa(S^3)$ is supported in degree $M=0$;
the third follows from the reformulation of $\tau(K)$
(Equation~\eqref{eq:CharacterizeTau}); and the fourth from the
normalization of $N=M-2A$.

\begin{theorem}\label{THM:HFKTOHFD}
  \label{thm:HFKtoHFD} Let $\CFKm(K)$ be a model for a chain complex
  for a knot $K \subset S^3$ which is reduced.  Let $Y$ be the bordered
  three-manifold $S^3\setminus\nbd{K}$. Given any sufficiently large integer $n$,
  we consider $Y$ with framing $-n$. The associated type $D$ module
  $\CFDa(Y)$ can be extracted from $\CFKm(K)$ by the following
  procedure.

  In the idempotent $\iota_0$, $\CFDa(Y)$ is isomorphic as an
  $\Field$-module to $\HFKa(K)$.\footnote{That is, the submodule
    $\iota_0\CFDa(Y)$ is identified with $\HFKa(K)$ as an
    $\Field$-module.}
  Let $\{\xi_i\}$ be a vertically simplified basis
  for $\CFKm(K)$.  To each vertical arrow of length $\ell$ from
  $\xi_i$ to $\xi_{i+1}$, we introduce a string of basis elements
  $\kappa^{i}_1,\dots,\kappa^{i}_\ell$
  for $\iota_1\CFDa(Y)$ and differentials
  $$
  \xi_i 
  \stackrel{D_{1}}\longrightarrow
  \kappa^{i}_1
  \stackrel{D_{23}}\longleftarrow
  \cdots
  \stackrel{D_{23}}\longleftarrow
  \kappa^{i}_{k}
  \stackrel{D_{23}}\longleftarrow
  \kappa^{i}_{k+1}
  \stackrel{D_{23}}\longleftarrow
  \cdots
  \stackrel{D_{23}}\longleftarrow
  \kappa^{i}_{\ell}
  \stackrel{D_{123}}\longleftarrow
  \xi_{i+1}
  $$
  where $\xi_i$ and $\xi_{i+1}$ are viewed as elements of $\CFDa(Y)$.
  (The reader is cautioned to note the
  directions of the above arrows.)
  
  Similarly, endow $\CFKm(K)$ with a horizontally simplified
  basis $\{\eta_i\}$. 
  To each horizontal arrow from $\eta_i$ to $\eta_{i+1}$ of length
  $\ell$, we introduce a string of basis elements
  $\lambda^{i}_1,\dots,\lambda^{i}_\ell$
  connected by a string of
  differentials
  $$
  \eta_i 
  \stackrel{D_{3}}\longrightarrow
  {\lambda}^{i}_1
  \stackrel{D_{23}}\longrightarrow
  \cdots
  \stackrel{D_{23}}\longrightarrow
  {\lambda}^{i}_{k}
  \stackrel{D_{23}}\longrightarrow
  {\lambda}^{i}_{k+1}
  \stackrel{D_{23}}\longrightarrow
  \cdots
  \stackrel{D_{23}}\longrightarrow
  {\lambda}^{i}_\ell
  \stackrel{D_{2}}\longrightarrow
  \eta_{i+1}.
  $$
  (As before, $\eta_i$ and $\eta_{i+1}$ are viewed as elements of $\CFDa(Y)$,
  via the identification of $\iota_0\CFDa(Y)$ with $\CFKm(K)$.)
  
  There is a final string, the {\em unstable chain}, of generators
  $\mu_1,\dots,\mu_m$ connecting $\xi_0$ and $\eta_0$, connected by
  a string of differentials, 
  $$\xi_0\stackrel{D_1}{\longrightarrow}\mu_1\stackrel{D_{23}}\longleftarrow
  \mu_2\stackrel{D_{23}}\longleftarrow
  \cdots
  \stackrel{D_{23}}\longleftarrow
  \mu_m
  \stackrel{D_3}\longleftarrow
  \eta_0.$$
  \glsadd{coeffmap}
  The length~$m$ of this string is $n+2\tau(K)$.

  The gradings are determined as follows:
  \begin{itemize}
  \item The grading set is $G/\lambda^{-1} \gr(\rho_{23})^n
    \gr(\rho_{12})^{-1} = G/\bigl(\frac{n-1}{2};-1,n\bigr)$.
  \item Recall that any element $\x_0$ of $V^0$ is represented by a generator
    $\x$ for the knot Floer complex. The grading of $\x_0$ in the above grading
    set is determined by the Maslov grading $M$
    and Alexander grading $A$ of $\x$ by the formula
    \begin{equation}\label{eq:gr-on-CFK-to-CFD}
    \gr(\x_0)=\lambda^{M(\x)-2A(\x)}\left(\gr(\rho_{23})\right)^{-A(\x)}
    = (M-{\textstyle\frac{3}{2}}A;0,-A).
    \end{equation}
  \end{itemize}
\end{theorem}

\begin{remark}
  \label{rmk:Normalization}
  The normalization in Formula~\ref{eq:gr-on-CFK-to-CFD} is chosen so
  that if one tensors the module with $\CFAa(\HD_\infty)$, the
  invariant of the $\infty$-framed solid torus, and declares that the
  generator $r$ for the solid torus has grading $0$, then the
  generator $r\otimes \eta_0$ of $\HFa(S^3)$ has (Maslov) grading $0$, as well.
  Although we do not specify the
  gradings of the elements of $V^1$, these can be readily worked out,
  since they are connected to elements of $V^0$ by coefficient maps.
  See Equation~\eqref{eq:TrefoilCFDGradings} for an example.
\end{remark}

As stated, Theorem~\ref{thm:HFKtoHFD} only computes 
$\CFDa$ for sufficiently negative framing
parameters.  In~\cite{LOT2}, we give bimodules which allow us to
deduce the type $D$ module for arbitrary framings; this is summarized
in Appendix~\ref{app:Bimodules}.  It turns out that most of
the statement of Theorem~\ref{thm:HFKtoHFD} goes through for arbitrary
integral framings; the only part which depends on the framing is the unstable chain,
connecting $\xi_0$ and~$\eta_0$, which can be explicitly described.
In particular, for framing $2\tau(K)$, this entire chain can be
replaced by a single map
\begin{equation}
  \label{eq:UnstableChain}
  \xi_0 \stackrel{D_{12}}\longrightarrow \eta_0.
\end{equation}
See Theorem~\ref{thm:HFKtoHFDframed} for the complete statement.

The procedure for translating knot Floer complexes
to type $D$ modules given by Theorem~\ref{thm:HFKtoHFD} is illustrated in
Figure~\ref{fig:HFKtoHFD}, with one proviso: to keep the pictures
small, we have illustrated type $D$ modules with framing parameter
$n=2\tau(K)$, for which one needs to refer to
Theorem~\ref{thm:HFKtoHFDframed}.

\begin{figure}
  \input{HFKtoHFD}
  \caption[Examples of reconstructing $\CFDa$ from $\CFK$]{\label{fig:HFKtoHFD}
    {\bf{$\CFKm$ and $\CFDa$.}}  The top, middle, and
    bottom rows give data for the left-handed $(3,4)$ torus knot, the
    $(2,-1)$ cable of the left-handed trefoil, and the figure eight
    knot respectively.  Knot Floer
    complexes are on the left, with vertical (resp. horizontal)
    arrows specifying  $\bdy_z$ (resp. $\bdy_w$); arrow lengths
    encode  their Alexander grading changes.
    Type $D$ modules with framing parameter $2\tau(K)$ ($-6$,
    $-4$ and $0$, respectively) are on the right: 
    generators in the idempotent $\iota_0$ (resp. $\iota_1$) are 
    indicated by 
    black
    (resp. white) dots.
    The dotted arrow is
    the unstable chain.}
\end{figure}

It is interesting to note that the above construction of $\CFDa$ uses
only the vertical and horizontal arrows, while $\CFKm$ in general can
contain diagonal arrows.
\index{arrows!diagonal}%
In concrete terms, the vertical arrows
correspond to counting holomorphic disks~$\phi$ with $n_{w}(\phi)=0$, while
horizontal ones correspond to counts of holomorphic disks with
$n_z(\phi)=0$.  
Diagonal arrows come from disks crossing both $z$ and $w$.

We will deduce Theorem~\ref{thm:HFKtoHFD} from a more
complicated statement, Theorem~\ref{thm:HFKtoHFD2},
which does not make reference to a choice of basis, and which comes
more directly from a Heegaard diagram. This basis-free version is
stated and proved in Section~\ref{sec:CFK-to-CFD}, and we deduce
Theorem~\ref{thm:HFKtoHFD} from it in
Section~\ref{sec:HFKtoHFDproof}.  Before turning giving these
proofs, we discuss some further counts of holomorphic curves, in
Section~\ref{sec:BoundaryDegenerations}.

\section{Generalized coefficient maps and boundary degenerations}
\label{sec:BoundaryDegenerations}

So far, the holomorphic curves we count have not covered the
basepoint, and in particular have not involved the Reeb chords
that contain the basepoint.
In this section, we explain a certain further structure on the type $D$
module coming from counting certain holomorphic curves which contain
$\rho_0$ in their boundary as well. These additional
maps will be used in Section~\ref{sec:CFK-to-CFD} to construct
$\CFDa$ from $\CFKm$. This further structure exists for arbitrary
(connected) boundary, but for notational simplicity we will restrict
our attention to the torus boundary case.

Fix a bordered Heegaard diagram $(\HD,z)$ for
$(Y^3,T^2)$, and label the regions around $\bdy\widebar{\Sigma}$ as in
Section~\ref{sec:torus-algebra}. Recall that for $I=\{i_0,\dots,i_n\}$
an interval in $\{1,2,3\}$ there is an associated map $D_I\co
V^{[i_0-1]}\to V^{[i_n]}$. This is induced by the differential on
$\CFDa$, which in turn is defined by counting holomorphic curves in
$\Sigma\times[0,1]\times\RR$, cf.\ Chapter~\ref{chap:type-d-mod}. By moving
the basepoint $z$ from the region labeled $0$ to one of the regions
labeled $1$, $2$ or~$3$, we obtain maps $D_I\co V\to V$ for intervals
$I$ in $\{0,1,2,3\}$ with respect to the cyclic ordering, of length
at most $3$. (For instance, putting $z$ in region $1$ the new maps we
obtain are $D_0$, $D_{30}$ and $D_{230}$.) We will call these new maps
\emph{generalized coefficient maps.}
\index{coefficient maps!generalized}%
\index{generalized!coefficient maps|see{coefficient maps, generalized}}

\glsadd{coeffmap}
We define four more generalized coefficient maps, $D_{0123}$,
$D_{1230}$, $D_{2301}$ and $D_{3012}$ as follows. By dropping
Condition~(\ref{item:moduliPenUlt}) on holomorphic maps from
Section~\ref{sec:curves-in-sigma}, we can consider moduli spaces of
holomorphic curves asymptotic to Reeb chords $\rho_0$, $\rho_{01}$ and
so on. Collecting these into moduli spaces as in
Section~\ref{sec:def-CFD}, for $\x$ a generator of
$\CFDa(\HD,z)$ in idempotent $\iota_0$, we set
\[
D_{0123}(\x)\coloneqq\!\!\!\sum_{\substack{\y\in\S(\HD)\\
    B\in\piBig(\x,\y)\\ \vec{\rho}\,\mid\, \ind(B,\vec{\rho})=1}}
  \!\!\!\#\left(\Mod^B(\x,\y;\vec{\rho})\right)\y
\]
where the sum over $\vec\rho$ is over
$\{(\rho_0,\rho_1,\rho_2,\rho_3),(\rho_{012},\rho_3),(\rho_{0},\rho_{123})\}$.
(These are the three sequences of Reeb chords
for which $\prod
a(-\rho_\alpha) = a(-\rho_{0123})$ in a suitable generalization of
$\Alg$.)%
\footnote{See the Errata at the end for corrections to
  this definition and the following proof.}
The maps $D_{1230}$, $D_{2301}$ and $D_{3012}$ are defined similarly.

By moving the basepoint $z$, it follows from the discussion in
Section~\ref{sec:torus-algebra} that these maps satisfy the
compatibility conditions
\[
  \sum_{\{I=J\cup K\,\mid\,J<K\}} D_{K}\circ D_{J}=0
\]
where $<$ now denotes the cyclic ordering on $\{0,1,2,3\}$ and $I$ is
any \emph{proper} cyclic interval in $\{0,1,2,3\}$.  The compatibility
equations involving all four Reeb chords are, perhaps, more
surprising, as there is now a right hand side.
\begin{proposition}
  \label{prop:MatrixFactorization}
 \index{coefficient maps!generalized}
  The maps $D_{0123}$ and $D_{2301}$ satisfy
  \begin{align}
    \label{eq:Homotopy0123}
    D\circ D_{0123}+D_3\circ D_{012} + D_{23}\circ D_{01} + D_{123} \circ D_0
  + D_{0123}\circ D &=\Id\co V^1\to V^1\\
  \label{eq:Homotopy2301}
    D\circ D_{2301}+D_1\circ D_{230} + D_{01}\circ D_{23} + D_{301} \circ D_2
    + D_{2301}\circ D &=\Id\co V^1\to V^1.
  \end{align}
  The maps $D_{1230}$ and $D_{3012}$ satisfy similar formulas.
\end{proposition}
(Here, $\Id$ denotes the identity map.)

\begin{remark}
  The grading shifts of the new
  coefficient maps may be computed directly. Alternatively, they are
  uniquely
  specified by the grading shifts for the original grading shift
  maps, together with the above formulas, which are homogeneous.
\end{remark}

\begin{figure}
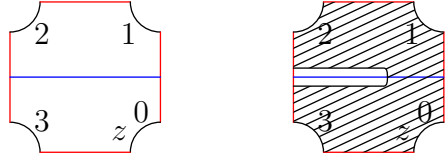

  \[
  \mfigb{torus-60}\qquad\qquad\mfigb{torus-61}
  \]
  \caption[Boundary degeneration as end of a one-dimensional
      moduli space]{\textbf{Boundary degeneration as end of a one-dimensional
      moduli space.} Left: a genus $1$ Heegaard diagram
    $(\Sigma,\alphas,\betas)$ for the solid torus. Here, $\beta_1$ is
    disjoint from $\alpha_1^a$ and intersects $\alpha_2^a$ in a single
    point. Right: one point in a one-parameter family of holomorphic
    curves in $\Sigma\times[0,1]\times\RR$, with one end a height $2$ 
    comb and the other end a boundary
    degeneration.}\label{fig:simplest-torus-diagram}
\end{figure}

The reason that the right hand sides of
Equations~\eqref{eq:Homotopy0123} and~\eqref{eq:Homotopy2301} are the
identity map rather than zero is the existence of certain boundary
degenerations.
\index{boundary degeneration}%
\index{degeneration of holomorphic curves!boundary}%
We can see these boundary degenerations explicitly in
simple cases. For example, let $(\Sigma,\alphas,\betas,z)$ be a genus
$1$ Heegaard diagram for a solid torus---for instance, the diagram in
Figure~\ref{fig:simplest-torus-diagram}. Then $\Sigma\setminus\alphas$
is a rectangle (disk with four boundary punctures). There is a
holomorphic map $u\co S\to\Sigma\times[0,1]\times\RR$, where $S$ is
the union of a trivial strip and a rectangle, and the rectangle is
mapped in the obvious way to $\Sigma\setminus\alphas$, and by a
constant map to $[0,1]\times\RR$.

Such curves look quite strange in the cylindrical setting. (In
particular, it is not clear when they are transversely cut out.)
We therefore invoke, for the first and only time in this book, the
\emph{tautological correspondence}\index{tautological correspondence}
between holomorphic curves in
$\Sigma\times[0,1]\times\RR$ and holomorphic disks in $\Sym^g(\Sigma)$. That is,
recall that there is a bijective correspondence
\[
  \begin{tabular}{c}
    Holomorphic maps\\
    $\phi\co[0,1]\times\RR\to\Sym^g(\Sigma)$\\
    with $\phi(\{0\}\times\RR)\subset\prod\beta_i$\\
    and $\phi(\{1\}\times\RR)\subset\prod\alpha_i$\\ \ 
  \end{tabular}\;\;\longleftrightarrow
  \begin{tabular}{c}
    Holomorphic curves\\
    $u$ in $\Sigma\times[0,1]\times\RR$\\
    with boundary in $\betas\times\{0\}\times\RR$\\
    and in $\alphas\times\{1\}\times\RR$\\
    so that $\pi_\DD\circ u$ is a $g$-fold branched cover.
  \end{tabular}
\]
(In the simplest case, one uses a complex structure $\Sym^g(j_\Sigma)$
on  $\Sym^g(\Sigma)$ and a split complex structure $j_\Sigma\times
j_\DD$ on $\Sigma\times[0,1]\times\RR$. More generally, admissible
almost complex structures $J$ on $\Sigma\times[0,1]\times\RR$ correspond to
working with paths of complex structures on $\Sym^g(\Sigma)$.)
See, for instance, \cite[Section 13]{Lipshitz06:CylindricalHF} for a
detailed explanation. 

To extend this to the bordered setting, recall that
$\Sigma_{\overline{e}}$ denotes the result of filling-in the puncture
$p$ of $\Sigma$. Viewing $\Sym^g(\Sigma)$ as the complement of the
divisor $p\times\Sym^{g-1}(\Sigma_{\overline e})$ in
$\Sym^g(\Sigma_{\overline e})$, we see that $\Sym^g(\Sigma)$ is a symplectic
manifold with a cylindrical end, in a way compatible with the complex
structure
$\Sym^g(j_\Sigma)$. The end of $\Sym^g(\Sigma)$ is modeled on the unit
normal bundle $U\nu(\Sym^{g-1}(\Sigma_{\overline e}))$ to
$p\times\Sym^{g-1}(\Sigma_{\overline e})$, an $S^1$-bundle over
$\Sym^{g-1}(\Sigma_{\overline e})$; compare~\cite[Example 2.2 and
Remark 5.9]{BEHWZ03:CompactnessInSFT}. There is a Reeb-like vector
field $\vec{R}$ on $U\nu(\Sym^{g-1}(\Sigma_{\overline e}))$ tangent to
the $S^1$-fibers.

The Lagrangian subspace
$\alpha_i^a\times\alpha_1^c\times\dots\times\alpha_g^c$ ($i=1,2$)
intersects the ideal boundary of $\Sym^{g}(\Sigma)$ in the subspace
$(\bdy\overline{\alpha}_i^a)\times
\alpha_1^c\times\dots\times\alpha_g^c$. So, for instance, the Reeb
chord $\rho_1$ corresponds to a $(g-1)$-dimensional family of
$\vec{R}$-chords between $\bdy\overline{\alpha}_1^a$ and
$\bdy\overline{\alpha}_2^a$, parameterized by the points in
$\alpha_1^c\times\dots\times\alpha_g^c$; and similarly for the other
Reeb chords.

The tautological correspondence extends to a correspondence between
holomorphic maps $u\co S\to\Sigma\times[0,1]\times\RR$ asymptotic to a
sequence of Reeb chords $(\rho_1,\dots,\rho_n)$ at $e\infty$ and
generators $\x$ and $\y$ at $\pm\infty$ and holomorphic maps
\[
\phi\co (([0,1]\times\RR)\setminus\{(1,t_1),\dots,(1,t_n)\})\to \Sym^g(\Sigma)
\]
asymptotic to a chord of the form $\rho_i\times
\{(a_{i,1},\dots,a_{i,g-1})\}$ at the puncture $(1,t_i)$ and to $\x$ and
$\y$ at $\pm\infty$.  Here, the values of the $t_i$ and the points
$a_{i,j}\in \alpha_j$ are allowed to vary, with the restriction that
$t_1 < t_2 < \dots < t_n$. In particular, the
$\RR$-coordinates of the $e$ punctures of $u$ correspond to the values
of the $t_i$.

With respect to this correspondence, the boundary degeneration
mentioned above comes from a disk in $\Sym^g(\Sigma)$ with boundary
entirely in the $\alpha$-tori; we will view this as a map from the
upper half-plane $\HalfPlane$ with some boundary punctures. There is
nothing strange about this disk. In particular, for generic $J$ on the
symmetric product $\Sym^g(\Sigma)$, the moduli spaces of these disks
will be transversely cut out.

We collect the boundary degenerations in $\Sym^g(\Sigma)$ into moduli
spaces:
\begin{definition}
  Given $i=1$ or $2$,
  $\x\in\alpha^a_{i}\times\alpha_1^c\times\dots\times\alpha_{g-1}^c$,
  a sequence $\vec{\rho}=(\rho_1,\dots,\rho_n)$ of Reeb chords, and an
  almost complex structure $J$ on $\Sym^g(\Sigma)$, let
  $\cM^{[\Sigma]}(\x;\vec{\rho};J)$ denote the moduli space of
  $J$-holomorphic maps
  \[
  \phi\co \HalfPlane\setminus\{(t_1,\dots,t_n)\}\to \Sym^g(\Sigma),
  \]
  where the $t_i$ are real numbers (points in $\bdy\HalfPlane$) which
  are allowed to vary with the given order, such that the
  boundary of $\HalfPlane$ is
  mapped into $\{\alpha_1^a \cup \alpha_2^a\}\times\alpha_2^c\times
  \dots\times \alpha_{g-1}^c$, the homology class of $\phi$ is
  $[\Sigma]$, $\phi$ is asymptotic to $\x$ at $\infty$ and to a chord
  of the form $\rho_i\times\{(a_{i,1},\dots,a_{i,g-1})\}$ at $t_i$.
\end{definition}

\begin{proposition}\label{prop:bdy-degen-exists} 
  There is a non-empty open subset of nearly-symmetric almost complex
  structures $J$ on $\Sym^g(\Sigma)$, as in \cite[Section
  3]{OS04:HolomorphicDisks}, such that the moduli spaces
  $\cM^{[\Sigma]}(\x;\vec{\rho};J)$, where $\vec{\rho}$ is a cyclic
  permutation of $(\rho_0,\rho_1,\rho_2,\rho_3)$ and $\x\in\S(\HD)$, are
  transversely cut out and have an odd number of points. Further, for
  these almost complex structures, all
  $\cM^{[\Sigma]}(\x;\vec{\rho};J)$ for other $\vec{\rho}$ are empty.
\end{proposition}
\begin{proof}
  As already discussed, it is standard that
  $\cM^{[\Sigma]}(\x;\vec{\rho};J)$ is transversely cut out for
  generic $J$.  To prove the moduli space has an odd number of points,
  we start by studying a model case. Consider the Heegaard diagram
  $(\Sigma_1,\alphas_1,\betas_1)$ shown in
  Figure~\ref{fig:simplest-torus-diagram}. In this diagram, there is a
  single generator $\x$. 
  With respect to any almost complex structure
  $j$ on $\Sigma_1$, there is a
  one-parameter family of holomorphic disks in $\Sigma_1=\Sym^{g_1}(\Sigma_1)$, 
  in
  the homology class $[\Sigma_1]$, with asymptotics
  $(\rho_0,\rho_1,\rho_2,\rho_3)$. One end of this moduli space is a
  broken holomorphic curve and the other end is the boundary
  degeneration under discussion. 
  Since the domain of these curves has multiplicity $1$ somewhere (in
  fact, everywhere), for a generic $j$ on $\Sigma$, this moduli space
  of $\Sym^g(j)$-holomorphic
  curves is transversely cut out, including at its boundary.
  So, the count of boundary degenerations in
  $\Sym^g(\Sigma_1)=\Sigma_1$ is $1$.

  Now, for the general case, let $(\Sigma,\alphas,\betas)$ be any
  bordered Heegaard diagram with torus boundary. Then
  $(\Sigma,\alphas)$ can be decomposed as a connect sum
  $(\Sigma_1,\alphas_1)\connectsum(\Sigma_2,\alphas_2)$ where
  $(\Sigma_1,\alphas_1)$ is the $\alpha$-curves of the diagram discussed above and
  $\alphas_2$ consists of closed curves.  A simple adaptation of the
  arguments showing stabilization invariance of $\HF$ of closed
  3-manifolds \cite[Section 10]{OS04:HolomorphicDisks} shows that for
  almost complex structures with a long neck between
  $(\Sigma_1,\alphas_1)$ and $(\Sigma_2,\alphas_2)$, there are an odd
  number of holomorphic curves in
  $\cM^{[\Sigma]}(\x\semico(\rho_0,\rho_1,\rho_2,\rho_3);J)$. Similar
  arguments apply for other cyclic orderings of
  $(\rho_0,\rho_1,\rho_2,\rho_3)$. This proves the first claim in the
  proposition.

  Curves with other asymptotics (e.g., $(\rho_0,\rho_{123})$) cannot
  occur in $(\Sigma_1,\alphas_1,\betas_1)$ by inspection. Hence, if
  there is a long neck between $(\Sigma_1,\alphas_1)$ and
  $(\Sigma_2,\alphas_2)$, they do not occur for
  $(\Sigma,\alphas)$. This proves the second claim in the proposition.
\end{proof}

\begin{proof}[Proof of Proposition~\ref{prop:MatrixFactorization}]
  We will prove Relation~(\ref{eq:Homotopy0123}); the other relations
  follow similarly. Consider the moduli
  space
  \[
  \bigcup_{\ind(B,(\rho_0,\rho_1,\rho_2,\rho_3))=2}\!\!\!\!\Mod^B(\x,\y;(\rho_0,\rho_1,\rho_2,\rho_3)).
  \]
  For $\x\neq\y$, the ends of this moduli space all correspond to
  terms on the left side of Formula~(\ref{eq:Homotopy0123}). By
  contrast, for $\x=\y$ it follows from
  Proposition~\ref{prop:bdy-degen-exists}, together with a standard
  gluing argument (compare~\cite[Theorem 5.1]{OS05:HFL}) and the tautological
  correspondence, that an odd number of ends of the moduli space
  correspond to boundary degenerations. The result is immediate.
\end{proof}

\section{From \textalt{$\CFKm$}{CFK-} to \textalt{$\CFDa$}{CFD\textasciicircum}:
  Basis-free version}
\label{sec:CFK-to-CFD}

In this section we state and prove a basis-free version of
Theorem~\ref{thm:PairingKnot}; this is Theorem~\ref{thm:HFKtoHFD2} below. We
deduce Theorem~\ref{thm:PairingKnot} from Theorem~\ref{thm:HFKtoHFD2}
in the next section.  Before stating the basis-free version, we
introduce a little more notation for knot Floer homology.

Fix a doubly-pointed
Heegaard diagram~$\HD$ for $K$, with basepoints $w$ and~$z$.
We write
$\gls*{knotCxAssocGr}$
for the chain complex $\gCFKa(\HD,w,z)$
computing $\HFKa(K)$
and 
$\gls*{knotCxAssocGrSummand}$
for its summand in Alexander grading~$r$.
There are maps 
$\gls*{bdyi}\co C(r)\rightarrow C(r+i)$
defined by counting holomorphic disks $\phi$ with
$n_{z}(\phi)=-i$ and $n_{w}(\phi)=0$ if $i\leq 0$, and $n_{z}(\phi)=0$
and $n_{w}(\phi)=i$ if $i\geq 0$. Set 
$\gls*{bdyw} \coloneqq \sum_{i\geq 0} \partial^i$
and
$\gls*{bdyz}\coloneqq\sum_{i\leq 0} \partial^i$.
The maps $\partial_w$ and $\partial_z$ each give differentials on~$C$.
Furthermore, $C$ has two filtrations,
\[
\gls*{knotCxAssocGrPos}\coloneqq\bigoplus_{r\geq s} C(r)\qquad\text{and}\qquad
\gls*{knotCxAssocGrNeg}\coloneqq\bigoplus_{r\leq s} C(r),
\]
preserved by $\partial_w$ and $\partial_z$, respectively.  We will
principally use not the complex $(C(\geqo s),\partial_w)$, but the
complex $(C(\leqo s),\partial_w)$, defined to be the quotient of
$(C,\partial_w)$ by the subcomplex $(C(\geq s+1),\partial_w)$.  We
likewise need the quotient complex $(C(\geqo s),\partial_z)$, with a
similar definition.

The above constructions can be interpreted using the notions from 
Section~\ref{sec:CFKm-to-CFDa-statement}. The complex $\CFKm(K)$ is 
a finitely generated, free, $\ZZ$-filtered, $\ZZ$-graded chain complex
over $\Field[U]$
satisfying compatibility conditions spelled out in that section.
The complex $(C,\partial_w)$ is the horizontal complex associated to $\CFKm(K)$,
while $(C,\partial_z)$ is the vertical complex.
In particular, $(C(\geqo s),\partial_w)$ is identified with the
$s\th$ summand of $\gCFKm(\HD,s)$ with respect to the Alexander grading,
while $(C(\leqo s), \partial_z)$ is identified with the $s\th$
filtered complex of $\CFKa(\HD)$.  Thus $(C(\geqo s), \partial_w)$ and
$(C(\geqo s),\partial_z)$ for
$s$ sufficiently small both have homology isomorphic to $\Field$;
likewise for $(C(\leqo s), \partial_z)$ and
$(C(\leqo s),\partial_w)$ for $s$ sufficiently large.

\begin{theorem}
  \label{thm:HFKtoHFD2}
  Let $K\subset S^3$ be a knot with meridian~$\mu$ and 0-framed
  longitude~$\lambda$, and let $C$, $\partial_w$, and $\partial_z$ be
  the data associated as above. Given a
  positive integer~$\gls*{nframing}$,
  let 
  $$
  \gls*{VectZero}
  \coloneqq\bigoplus_{s\in\ZZ} V^0_s
  \qquad\text{and}\qquad 
  \gls*{VectOne}
  \coloneqq\bigoplus_{s\in\ZZ+\frac{n+1}{2}}
  V^1_s
  $$
  where
  \begin{align*}
    \gls*{VectZeroS}
    &\coloneqq C(s) \\
    \gls*{VectOneS}
    &\coloneqq \begin{cases}
        C\bigl(\leq s+\frac{n-1}{2}\bigr) & s\leq -\frac{n}{4} \\
        \Field & |s|< \frac{n}{4} \\
        C\bigl(\geq s-\frac{n-1}{2}\bigr) & s\geq \frac{n}{4}.
      \end{cases}
  \end{align*}
  Let 
  $\gls*{Vect}\coloneqq V^0\oplus V^1$.
  Then, for $n$ large enough, the following maps~$D_I$ satisfy the
  compatibility condition~(\ref{eq:Compatibility}):
  \begin{itemize}
  \item\glsadd{coeffmap} The differential $D$ restricted to $V^0_s=C(s)$ is the
    differential $\partial^0$ on the knot complex; for $s\leq
    -\frac{n}{4}$, $D$ restricted to
    $V^1_s =C\left(\leq s +\frac{n-1}{2}\right)$ is the
    differential~$\partial_w$; for $s\geq \frac{n}{4}$, $D$ restricted to
    $V^1_s = C\left(\geq s-\frac{n-1}{2}\right)$ is
    the differential $\partial_z$; for $|s|< \frac{n}{4}$, $D$
    restricted to $V^1_s = \Field$ is identically zero.
  \item The map
    $$D_1\co V^0_s=C(s)\to V^1_{s+\frac{n-1}{2}}=C(\geqo s)$$
    is the obvious inclusion of the subcomplex.
  \item The map
    $D_2$ is nonzero only on $V^1_s$ for $s\leq -\frac{n}{4}$, in which case
    $$D_2\co V^1_{s}=C\left(\leq s+\textstyle\frac{n-1}{2}\right) \to
    V^0_{s+\frac{n+1}{2}}=C\left(s+\textstyle\frac{n+1}{2}\right)$$
    is given by
    $$\pi\circ \partial_w\co C\left(\leq s +
      \textstyle\frac{n-1}{2}\right)\to 
    C\left(s+\textstyle\frac{n+1}{2}\right)$$
    where $\pi\co C\to
    C\left(s+\frac{n+1}{2}\right)$ is the projection.
  \item
    The map
    $$D_{3}\co V^0_s=C(s)\to V^1_{s-\frac{n-1}{2}}=C(\leqo s)$$
    is the obvious inclusion of the subcomplex.
  \item
    The map $D_{12}$ is identically zero.
  \item
    For $s< -\frac{n}{4}$, the map
    $$D_{23}\co V^1_s = C\left(\leq s+\textstyle\frac{n-1}{2}\right) \to
    V^1_{s+1}=C\left(\leq s+\textstyle\frac{n+1}{2}\right)$$
    is the obvious inclusion map;
    for the $s$ with $s \le - \frac{n}{4} < s+1$,
    $$D_{23}\co V^1_s=(C,\partial_w) \to V^1_{s+1}=\Field$$
    is a chain map inducing an isomorphism in homology;
    for $\frac{-n+2}{4}<s<\frac{n-2}{4}$,
    $D_{23}\co V^1_s = \Field\to V^1_{s+1} = \Field$
    is the isomorphism;
    for the $s$ with $s < \frac{n}{4} \le s+1$,
    $$D_{23}\co V^1_s =\Field\to V^1_{s+1} = (C,\partial_z)$$
    is a chain map inducing an isomorphism on homology;
    for
    $s>\frac{n}{4}$,
    $$D_{23}\co V^1_s=C\left(\geq s-\textstyle\frac{n-1}{2}\right)\to
    V^1_{s+1}= C\left(\geq s-\textstyle\frac{n+1}{2}\right)$$
    is the obvious projection map.
  \item\glsadd{coeffmap} The map
    $$D_{123}\co V^0_s\to V^1_{s+\frac{n+1}{2}}$$
    is
    the composite of one component of $\partial_w$,
    $$\partial^1_w\co C(s)\to C(s+1),$$
    with the inclusion of $C(s+1)$ in $(C(\geq
    s+1),\partial_z)=V^1_{s+\frac{n+1}{2}}$.
  \end{itemize}
  Moreover, the associated
  type $D$ module is homotopy equivalent to the bordered invariant
  $\CFDa(S^3\setminus\nbd(K))$ where the boundary is
  marked by the curves $\mu$ and $-n\cdot\mu+\lambda$.
  
  Let $\DSmallGrSet(K)$ denote the grading set for
  $\CFDa(S^3\setminus\nbd(K))$, so $\DSmallGrSet(K)$ is the quotient
  $\smallGroup/\langle R(g(P))\rangle$ where $P$ is a primitive periodic
  domain. The type $D$ structure $V^0\oplus V^1$ has a grading
  by $\DSmallGrSet(K)$ so that the homotopy equivalence with
  $\CFDa(S^3\setminus\nbd(K))$ respects the gradings. Further,
  there is a projection map $p\co \smallGroup\to
  (\OneHalf\ZZ)^2$ given by $p(m;i,j)=(i,j)$. The map $p$ induces a
  map $\DSmallGrSet(K)\to (\OneHalf\ZZ)^2/\langle p(R(g(P)))\rangle
  \to
  \OneHalf\ZZ$, giving a grading on $\CFDa(S^3\setminus\nbd(K))$ by
  $\OneHalf\ZZ$. (Explicitly, this map sends $(m;i,j)$ to $-ni-j$.) Then the homotopy equivalence identifies this
  grading with the $s$ grading above.
\end{theorem}

\begin{figure}
\begin{center}
\includegraphics{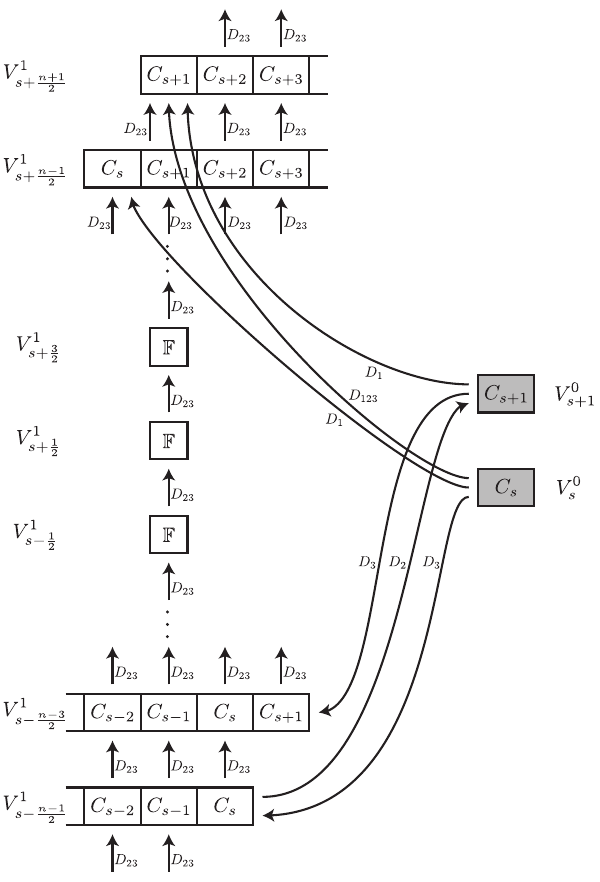}
\caption[Sketch of some of 
the maps from Theorem~\ref{thm:HFKtoHFD2}]{\textbf{A sketch of some of 
the maps from Theorem~\ref{thm:HFKtoHFD2}.} The vertical
coordinate represents the grading by relative $\spin^c$ structures.
The shaded squares correspond to $V^0$, while the unshaded ones represent
$V^1$.}\label{fig:HFKtoHFDschematic}
\end{center}
\end{figure}

Figure~\ref{fig:HFKtoHFDschematic} sketches
some of the maps from Theorem~\ref{thm:HFKtoHFD2}. 

\begin{remark}
  The complexes $C(\geqo s)$ appear in Hedden's cabling and Whitehead
  doubling formulas~\cite{Hedden, HeddenWhitehead}.
\end{remark}

In effect, Theorem~\ref{thm:HFKtoHFD2} describes how to associate a 
type $D$ module to a suitable chain complex
$C$, equipped with boundary operators $\bdy_w$ and $\bdy_z$ and
filtrations. In that statement, we apply it to $C=\CFKm(K)$; but the
association is functorial, in the following sense:

\begin{proposition}
  \label{prop:HotopyEquivalentDStructures}
  Up to homotopy equivalence, for sufficiently large~$n$ the type $D$
  module associated in
  Theorem~\ref{thm:HFKtoHFD2} depends only on the $\ZZ$-filtered,
  $\ZZ$-graded chain homotopy type of the $\Field[U]$-chain complex
  $\CFKm(K)$.
\end{proposition}
\begin{proof}
  A homomorphism $\phi$ of type $D$ structures (over $\Alg(\Torus)$)
  can be thought of as a collection of coefficient maps $\Phi_I$,
  $I\subset\{1,2,3\}$, satisfying the following straightforward
  analogue of Equation~\ref{eq:Compatibility}:
  \begin{equation}\label{eq:ModCoeffMapCompat}
  \sum_{\{I=J\cup K\mid J<K\}} \!\!D'_{K}\circ \Phi_{J}+\Phi_K\circ D_J=0.
  \end{equation}

  Let $\phi\co C\to C'$ be a $\ZZ$-filtered, $\ZZ$-graded homotopy
  equivalence of $\Field[U]$-chain complexes. The map $\phi$ induces
  homotopy equivalences on the associated complexes $C(s)$, $C(\geqo s)$, and
  $C(\leqo s)$, which in turn induces homotopy equivalences between
  the associated $V^0$ and $V^1$'s. Call these homotopy equivalences
  $\Phi_\emptyset$. The map $\phi$ also induces maps $C_{\leq s}\to
  C'_{s+1}$ which fit together to give maps $\Phi_{\{2\}}\co
  V^1_{s-\frac{n-1}{2}}\to V^{\prime\,0}_{s+1}$. Define all other
  $\Phi_I$ to be~$0$. It is straightforward to verify that $\Phi$
  satisfies
  Equation~\eqref{eq:ModCoeffMapCompat}. A homotopy of $\ZZ$-filtered,
  $\ZZ$-graded maps induces a homotopy of maps of type $D$ structures
  in the same way. This proves the result.%
\end{proof}


\begin{figure}
\begin{center}
\input{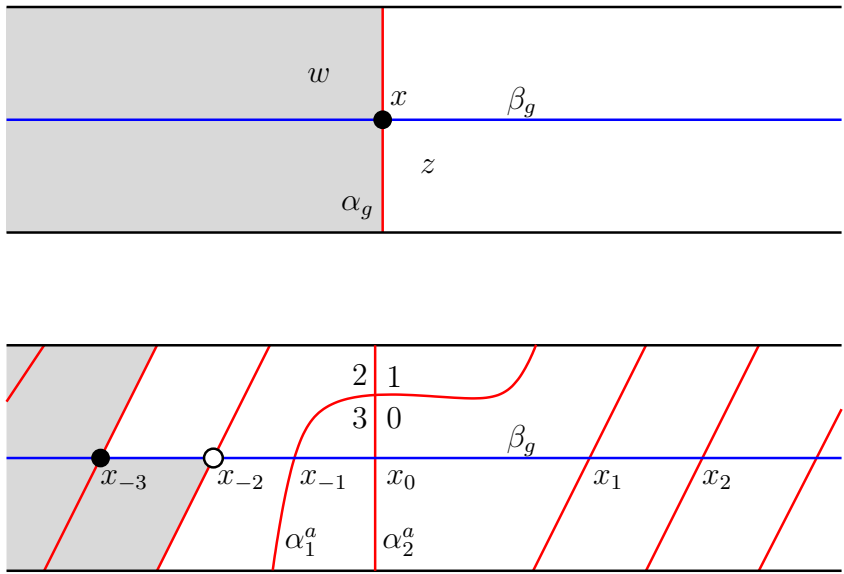}
\end{center}
\caption[Twisting, in proof of Theorem~\ref{thm:HFKtoHFD2}] {\label{fig:Twisting} {\bf{Twisting.}} At the top, we have
  pictured a portion of a doubly-pointed Heegaard diagram for a knot.
  Below it, we have a bordered diagram for the knot complement, with
  suitably twisted framing. The gray shading represents a domain from
  ${\mathbf a}\times x$ to ${\mathbf b}\times x$ (for some ${\mathbf
    a}$ and ${\mathbf b}$) which, in the top diagram, crosses $w$ with
  multiplicity one and $z$ with multiplicity zero. There is a corresponding
  (provincial) domain below connecting ${\mathbf a}\times x_{i}$ to
  ${\mathbf b}\times x_{i+1}$, provided that $i\leq -2$.}
\end{figure}

We prove Theorem~\ref{thm:HFKtoHFD2} by considering a bordered
Heegaard diagram~$\HD(n)$ for the
knot complement with sufficiently large negative framing~$-n$, in a way
analogous to the surgery formula for Heegaard Floer
homology~\cite{OS04:Knots}. For an analogous choice of bordered
Heegaard diagram~$\HD$, the corresponding type $D$ module,
$\CFDa(\HD)=W^0\oplus W^1$,\glsadd{Wi}
has
generators which correspond to the generators of $V=V^0\oplus V^1$ as
stated in the theorem, except for $W^1_{s}$ with $|s|<
\frac{n}{4}$.  Moreover, again except for $W^1_s$ with
$|s|<\frac{n}{4}$, the coefficient maps
are as stated in the theorem.  We then construct a homotopy equivalence
which fixes
the cases with $\abs{s} < \frac{n}{4}$.  We give the proof after some
preliminary lemmas.

To describe our bordered Heegaard diagram $\HD(n)$ we start with a
doubly-pointed Heegaard diagram
$\gls*{HDforK}$
for the knot $K$. We stabilize
this diagram to obtain a new diagram 
$\gls*{HDforKstab}$.
Here, $\Sigma$ is gotten by
attaching a two-dimensional one-handle to $\Sigma_0$,
with feet near the markings $z$ and $w$; 
$\alphas$ is gotten from $\alphas_0$
by introducing a new circle~$\gls*{alphag}$, 
contained in the new one-handle; and
$\betas$ is gotten from $\betas_0$ by introducing a new circle
$\gls*{betag}$,
which meets 
$\alpha_g$ in a single transverse intersection
point, and is disjoint from the other $\beta$-circles. In particular,
the circle $\alpha_g$ is a meridian for the knot.

Identify a neighborhood $\Winding$ of $\alpha_g$ in $\Sigma$ with
$[-1,1]\times S^1$ so that $\alpha_g$ is identified with $\{0\}\times
S^1$ and $\beta_g\cap \Winding$ is identified with
$[-1,1]\times\{e^{-\pi i/2}\}$. Let $\lambda$ denote a $0$-framed
longitude of $K$ in $\Sigma$, disjoint from the $\alpha_i$ for $i\neq
g$, and with $\lambda\cap \Winding=[-1,1]\times\{1\}$. Let
\[
\overline{\alpha}_1^a=(\lambda\cap(\Sigma\setminus \Winding))\cup\{(t,e^{\pi
  in(t+1)})\in \Winding\mid t\in[-1,1]\}.
\]
That is, $\overline{\alpha}_1^a$ is obtained from the $0$-framed longitude
$\lambda$ by winding $n$ times around the meridian $\alpha_g$
inside~$\Winding$. We
refer to 
$\gls*{Winding}$
as the \emph{winding region.}\index{winding region}\index{region!winding} Note that
$\overline{\alpha}_1^a$ intersects $\beta_g$ in $n$ points in $\Winding$, and
intersects $\alpha_g$ in a single point $p$. Define
$\alpha_1^a\coloneqq\overline{\alpha}_1^a\setminus\{p\}$
and $\alpha_2^a
\coloneqq \alpha_g\setminus\{p\}$.  Then, with notation as in
Section~\ref{sec:Bestiary},
\[
\gls*{HDofn}
\]
is a bordered Heegaard diagram for $S^3\setminus\nbd(K)$, with
framing~$-n$.  We will also write $\HD$ for $\HD(n)$ when the framing
is clear.
We label the four quadrants around $p$ by $0,\dots,3$, arranged
in a counterclockwise order (i.e., the order induced by
$\bdy\nbd(p)$), so that, if you forget $\alpha_1^a$, $0$ and $1$ are in
the region which contains $z$ in $\HD_K$. As
usual, we place the basepoint in region~$0$.  See
Figure~\ref{fig:Twisting} for an illustration. For convenience of
notation, we will
henceforth assume that $n$ is divisible by~$4$.

We now turn to the grading on $\CFDa(\HD)$.  Following
Section~\ref{sec:typeD-gradings}, the grading takes values in the
coset space $\smallGroup/P(\x_0)$, where $P(\x_0)$ is the span of
$g(B)$ over all periodic domains $B \in \pi_2(\x_0,\x_0)$.  For a
space with the same homology groups as a solid torus (such as a knot complement), the group of periodic
domains is isomorphic to $\ZZ$. Let 
$\gls*{BsubZero}$
be a generator of
$\pi_2(\x_0,\x_0)$.
Since the knot has framing~$-n$, we can
choose the sign of $B_0$ so that $B_0$ has local multiplicities
$(0,1,1-n,-n)$ in the regions
$(0,1,2,3)$, respectively.  Thus $P(\x_0)\subset\smallGroup$ is
$\langle(v;1,-n)\rangle$ for some $v \in \OneHalf\ZZ$ (depending on $\x_0$).
For present purposes it suffices to consider only the $\SpinC$
component of the grading on $\CFDa(\HD)$, which we can recover with
the homomorphism from $\smallGroup$ to $\OneHalf\ZZ$ given by
\begin{equation}
  \label{eq:GradingMap}
  (j;p,q)\longmapsto -np-q.
\end{equation}
(Note that $g(P(\x_0))$ is in the kernel of this map.)

More generally, let 
$
\gls*{AlexGrFun}\co \Gen(\HD)\rightarrow \OneHalf\ZZ$ be
a function characterized up to overall translation by the formula
\begin{equation}
    \label{eq:DefOfS}
  S(\x)-S(\y)=
   \left(\textstyle\frac{n+1}{2}\right)\cdot n_{0}(B)
  -\left(\textstyle\frac{n-1}{2}\right)\cdot n_{1}(B)
  -\left(\textstyle\frac{n+1}{2}\right)\cdot n_{2}(B)
  +\left(\textstyle\frac{n-1}{2}\right)\cdot n_{3}(B),
\end{equation}
where $B$ is any homology class in $\piBig(\x,\y)$ (i.e., classes possibly covering
the basepoint~$z$).

\begin{lemma}
  \label{lem:DModDegrees}
  Let $Y$ be a homology solid torus, i.e., a three-manifold with torus
  boundary and
  $H_1(Y;\ZZ)\cong\ZZ$.
  Pick $\mu,\lambda \in H_1(\bdy Y)$ so that
  $\lambda$ generates the kernel of the map on $H_1$ induced by
  inclusion of the boundary while the image of $\mu$
  generates $H_1(Y)$.  Let $\CFDa(Y) \cong W^0\oplus W^1$ be
  the decomposition by idempotents of the type~$D$ invariant of~$Y$,
  with boundary marked by $\mu$ and
  $-n\cdot\mu+\lambda$.  Then a function $S$ satisfying~\eqref{eq:DefOfS} gives a
  decomposition\glsadd{WiDecomp}
  $$
  W^i \cong\bigoplus_{s\in\frac{1}{2}\ZZ}W^i_s,
  $$
  well-defined up to overall translation.
  The coefficient maps respect this grading in the following sense:
  \begin{align*}
    D_{1}\co W^0_s&\to W^1_{s+\frac{n-1}{2}} &
    D_{2}\co W^1_s&\to W^0_{s+\frac{n+1}{2}} &
    D_{3}\co W^0_s&\to W^1_{s-\frac{n-1}{2}} \\
    D_{12}\co W^0_s&\to W^0_{s+n} &
    D_{23}\co W^1_s&\to W^1_{s+1} &
    D_{123}\co W^0_s&\to W^1_{s+\frac{n+1}{2}}.
  \end{align*}
\end{lemma}

\begin{proof}
  We first check that there is
  a function~$S$ satisfying Equation~\eqref{eq:DefOfS}. 
  To see this, observe that there is a natural splitting
  $\piBig(\x,\y)=\pi_2(\x,\y)\oplus \ZZ[\Sigma]$.
  Now, the right-hand side of Equation~\eqref{eq:DefOfS} remains unchanged when 
  we add a copy of $[\Sigma]$ to the domain $B$; so it suffices to show
  that the right-hand side of Equation~\eqref{eq:DefOfS} depends only
  on the generators $\x$ and $\y$, when $B$ is chosen to lie in $\pi_2(\x,\y)$
  (i.e., when $B$ satisfies $n_0(B)=0$).
  For such $B$, $S(B)$ is obtained by taking $g(B)$ and applying first
  the map in
  Formula~\eqref{eq:grading-relation} followed by the map in
  Formula~\eqref{eq:GradingMap}.  Since Formula~\eqref{eq:GradingMap}
  was chosen to annihilate periodic domains, it follows that the
  result is independent of tho choice of~$B$.

  The fact that the coefficient maps change grading as stated is
  immediate from the definition of $S$, Formula~\eqref{eq:DefOfS}.
\end{proof}

Let 
$\gls*{GenK}$
denote the set of generators for the knot Floer complex
$\CFKa(\HD_K)$,
and $\Gen(n)$\glsadd{Genn} denote the generators for
$\CFDa(\Heegaard(n))$.
Order the intersection points of $(\alpha_1\cup\alpha_2)$ with $\beta_g$ in
the winding region in the order they are encountered along~$\beta_g$,
$\gls*{xiWind}$,
so that $x_0$ is the
intersection point of $\alpha_2$ with~$\beta_g$.
(Recall that we assume for notational convenience that $n$ is
divisible by $4$.)
For each generator
$\x\in\Gen_K$ there are $n+1$ corresponding generators in $\Gen(n)$, gotten
by adding a point in $(\alpha_1\cup\alpha_2)\cap\beta_g$ in the
winding region. That is, we let 
$\gls*{xiWindBold}$
denote the
sequence of generators with the property that the
$(\alpha_1\cup\alpha_2)\cap\beta_g$-component of $\x_i$ is $x_i$ and,
outside the winding region, $\x_i$ agrees with~$\x$.

\begin{lemma}
  \label{lem:NormalizeS}
  For a given $n$, we can choose a function~$\gls*{AlexGrFun}$ 
  on $\Gen(n)$ as above so
  that it is related to the Alexander grading $A\co \Gen_K
  \to \ZZ$\index{grading!Alexander!relationship to $S$} of elements in $\Gen_K$ by the formula
  
  \begin{equation}
    S(\x_k)=A(\x)-k+\left(\frac{(n+1)\cdot \sgn(k)}{2}\right),\label{eq:SNormalization}
  \end{equation}
  where $\sgn(k)=-1$, $0$, or $1$ if $k<0$, $k=0$, or $k>0$ respectively.
\end{lemma}

\begin{proof}
  The region $3$ is a domain
  connecting $\x_0$ to $\x_{-1}$, and for $k < 0$ there is a
  domain connecting $\x_k$ to $\x_{k-1}$ with multiplicities
  $(0,0,1,1)$ in regions $(0,1,2,3)$.  (See Figure~\ref{fig:Twisting} and the middle of
  Figure~\ref{fig:DegenerateWinding}.)  Similarly, region~$1$ connects
  $\x_0$ to $\x_1$ and, for $k>0$, the generators $\x_k$ and $\x_{k+1}$ are
  connected by a domain with multiplicities $(1,1,0,0)$.  It follows
  from Equation~\eqref{eq:DefOfS} that
  $$
  S(\x_k)=S(\x_0)-k+\left(\frac{(n+1)\cdot \sgn(k)}{2}\right).
  $$
  According to Equation~\eqref{eq:AlexanderUpToTranslation},
  the Alexander grading on knot Floer
  homology is characterized (up to translation) by
  the property that if $\x,\y \in \Gen_K$ and $\phi\in\piBig(\x,\y)$ is any domain in $\HD_K$, then
  $A(\x)-A(\y)=n_{z}(\phi)-n_{w}(\phi)$. We can convert $\phi$ into
  a domain $\phi'$ in $\HD(n)$ that connects $\x_0$ to $\y_0$ with
  $n_{0}(\phi)=n_{1}(\phi)=n_z(\phi)$ and
  $n_{2}(\phi)=n_{3}(\phi)=n_w(\phi)$ (as is
  evident from Figure~\ref{fig:Twisting}). 
  So,
  $$
  S(\x_0)=A(\x)+c.
  $$
  for some constant independent of the intersection point $\x\in
  \Gen_K$. We can modify $S$ by an overall translation to take $c=0$.
\end{proof}

We now have a collection of generators in $\Gen(n)$
which have corresponding elements of $\Gen_K$. The remaining
generators (coming from intersections of $\alpha_2^a$ with the
$\beta$-curves outside of the winding region) are
called {\em exterior generators}.\index{exterior generator}\index{generator!exterior}
Note that the non-exterior generators are the generators of the form
$\x_i$ for some $\x\in\Gen_K$, $i\in \{-n/2,-n/2+1,\dots,n/2\}.$

\begin{lemma}
  \label{lem:Exterior}
  For any doubly-pointed Heegaard diagram $\HD_K$ of a knot in~$S^3$,
  there is a constant $c$ (independent of the framing parameter $-n$)
  with the following property.  For the function $S\co \Gen(n)\to\QQ$
  as defined in Lemma~\ref{lem:NormalizeS}, the generators $\y\in\Gen(n)$
  with $|S(\y)|\geq c$ satisfy the following properties:
  \begin{enumerate}[label=(c-\arabic*),ref=c-\arabic*,leftmargin=*]
  \item\index{(c-1)--(c-2)}
    \label{item:Exterior}
    $\y$ is not an exterior generator, so $\y=\x_i$ for some
    $\x\in\Gen_K$, $i\in\ZZ$; and
  \item
    \label{item:TypeDetermined}
    the sign of $S(\y)$ agrees with the sign of $i$.
  \end{enumerate}
\end{lemma}

\begin{proof}
  First, we find $c_1$ so that generators~$\y$ with $\abs{S(\y)} \ge c_1$ satisfy
  Property~\eqref{item:Exterior}.
  Note that the set of exterior generators is independent of the
  framing parameter~$-n$. Start with some initial framing
  parameter~$-m$, fix some $\x\in\Gen_K$ and consider
  $\x_{-1}\in\Gen(m)$. For each
  exterior generator $\y$, we can find a domain
  $B_0\in\pi_{2}(\y,\x_{-1})$ with $n_{0}(B_0)=n_{1}(B_0)=0$ and
  $n_2(B_0)=n_3(B_0)=c(\y)$.  (Start with an arbitrary domain in
  $\piBig(\y,\x_{-1})$ and add copies of $[\Sigma]$ and the periodic domain to make $n_0 = n_1 =
  0$; since $\y$ and $\x_{-1}$ are in the same idempotent, we will
  then have $n_2 = n_3$.)  Now, in $\HD(n)$, the domain
  $B_{0}$ can be used to construct a new domain
  $B_{n}\in\piBig(\y,\x_{\frac{m-n}{2}-1})$ with
  $n_0(B_n)=n_1(B_n)=0$ and $n_2(B_n)=n_3(B_n)=c(\y)$.  Since there
  are only
  finitely many possible exterior generators, we see that
  $\{\,|S(\y)-S(\x_{\frac{m-n}{2}})|\,\mid \y\text{ an exterior generator}\}$
  is bounded independently of $n$. It follows from Lemma~\ref{lem:NormalizeS}
  that $|S(\y)|$ for exterior generators~$y$ is bounded by some
  constant~$c_1$ independent of~$n$, giving Property~(\ref{item:Exterior}).

  Next, we claim that there is a constant $c_2$ with the property that
  if $S(\x_i)>c_2$, then $i>0$. Specifically, choose $c_2\geq
  \max_{\x\in\Gen_K} A(\x)$.
  Then Equation~\eqref{eq:SNormalization} shows that
  $S(\x_0)\leq c_2$. Similarly,
  if $\x_{-i}$ is a non-exterior generator with $i>0$
  then, since $i\leq \frac{n}{2}$,
  Equation~\eqref{eq:SNormalization} guarantees that
  $$S(\x_{-i})=A(\x)+i-\frac{n+1}{2}\leq A(\x)-\frac{1}{2}\leq
  c_2-\frac{1}{2}<c_2.$$ 
  This guarantees Property~\eqref{item:TypeDetermined} for
  $S(\y)>0$. A symmetric argument ensures the property also
  in the case where $S(\y)<0$.
\end{proof}

Lemmas~\ref{lem:OrderGens}--\ref{lem:D23-big-s} below state that for an
appropriate class of almost
complex structures, the coefficient maps have the desired forms from
Theorem~\ref{thm:HFKtoHFD2}.  The
boundary of the winding region~$\Winding$ consists of two circles
$C_1$ and~$C_2$ parallel to the meridional $\alpha$-circle; we label
these so that $C_1$ is on the $2,3$ side of $\Winding$ and $C_2$ is on
the $0,1$ side of winding; see Figure~\ref{fig:PinchingWindingWinding}.
Fix an admissible
complex structure $J_0$ on $\Sigma\times[0,1]\times\RR$. We now
consider a two-parameter family of complex structures $J_{T_1,T_2}$ on
$\Sigma\times[0,1]\times\RR$ gotten by inserting a necks of lengths
$T_1$ and $T_2$ along $C_1$ and $C_2$, respectively. We say that a
result \emph{holds for sufficiently pinched almost complex structures}
if there is a $K_1$ so that for any $T_1>K_1$ there is a $K_2$ so that
for any $T_2>K_2$ the result holds for all generic almost complex
structures of the form $J_{T_1,T_2}$.%
\index{complex structure!sufficiently pinched}%
\index{sufficiently pinched complex structure}

\begin{lemma}
  \label{lem:OrderGens}
  Let $K\subset S^3$ be a knot with Heegaard diagram~$\HD_K$, and let
  $c$ be the constant from Lemma~\ref{lem:Exterior}.
  For any sufficiently large~$n$ the following holds on~$\HD(n)$. Given any
  sufficiently pinched almost complex
  structure, for each generator $\x\in\Gen_K$, with corresponding
  generators $\{\x_k\}$ in $\Gen(n)$, we have:
  \begin{itemize}
  \item $D_{23}(\x_k)=\x_{k-1}$ if $k<0$ and $S(\x_k) \leq  -c-1$;
  \item $D_{01}(\x_k)=\x_{k+1}$ if $k>0$ and $S(\x_k) \geq  c+1$;
  \item $D_3(\x_0)=\x_{-1}$; and
  \item $D_{1}(\x_0)=\x_1$.
  \end{itemize}
  That is, for each $\x\in\Gen_K$, we have a string
  of generators connected by coefficient maps as follows:
  $$
  \dots 
  \stackrel{D_{23}}\longleftarrow
  \x_{-3} 
  \stackrel{D_{23}}\longleftarrow
  \x_{-2} 
  \stackrel{D_{23}}\longleftarrow
  \x_{-1} 
  \stackrel{D_{3}}\longleftarrow
  \x_0
  \stackrel{D_1}\longrightarrow
  \x_1
  \stackrel{D_{01}}\longrightarrow
  \x_2 
  \stackrel{D_{01}}\longrightarrow
  \x_3
  \stackrel{D_{01}}\longrightarrow
  \dots
  $$
\end{lemma}

\begin{proof}  
  The domains realizing the maps displayed in the statement can be found inside the
  winding region. For example, the domain realizing $D_{3}$ carrying
  $\x_0$ to $\x_{-1}$ is the
  region marked by $3$ in Figure~\ref{fig:Twisting}. Domains from
  $\x_{-k}$ to $\x_{-k-1}$ are annuli~$\gls*{phik}$
  which have
  a boundary component in $\alpha_1^a$ as in the middle of
  Figure~\ref{fig:DegenerateWinding}.

  \textbf{The case of $D_{23}$.} We analyze the map $D_{23}$, verifying it has the
  form stated, in the following eight steps.

  \begin{steps}
  \step\label{step:D23-xk1-occurs}
  \emph{The term $\x_{-k-1}$ occurs in $D_{23}(\x_{-k})$.}
  By cutting $\phi_K$ along $\alpha_1^a$ at
  $x_{-k}$, we obtain a bigon realizing a non-zero coefficient of
  $\x_{-k-1}$ in $D_{23}(\x_{-k})$.
  The rest of the proof (for the operator $D_{23}$) consists in showing that this is the only term in $D_{23}(\x_{-k})$.

  \step\label{step:D23-no-yk1}
  \emph{If $\y_{-k-1}$ occurs as a term in $D_{23}(\x_{-k})$
    then $\y=\x$.}
  Suppose that $\y_{-k-1}$ occurs in $D_{23}(\x_{-k})$.
  Considering the corners shows that any domain
  $\phi\in\piBig(\x_{-k},\y_{-k-1})$ is a disjoint union of
  $\phi_k$ with a different positive domain $\phi'$.  As such, the
  corresponding moduli space naturally has an $(\RR\oplus \RR)$-action
  (translating the two components independently). If this
  moduli space is generic and one-dimensional, it follows that the
  other component $\phi'$ is trivial, so $\y_{-k-1}=\x_{-k-1}$.

  \step\label{step:D23-no-ext}
  \emph{No exterior generators occur in $D_{23}(\x_{-k})$.}
  If an exterior generator~$\y$ appears in $D_{23}(\x_{-k})$, then $S(\y)=S(\x_{-k})+1\leq -c$,
  in view of Lemma~\ref{lem:DModDegrees} and our assumption on
  $S(\x_{-k})$. By Lemma~\ref{lem:Exterior} (and specifically
  Property~\eqref{item:Exterior}), it follows that $\y$ is not an
  exterior generator.

  \step\label{step:D23-no-positive}
  \emph{No terms of the form $\y_{j}$ where $j>0$
    occur in $D_{23}(\x_{-k})$.}
  This is same argument as in Step~\ref{step:D23-no-ext}, using
  Property~\eqref{item:TypeDetermined}) this time.

  In view of Steps~\ref{step:D23-no-ext}
  and~\ref{step:D23-no-positive}, only terms of the form $\y_{-k+i}$
  with
  $-k+i<0$ can appear in $D_{23}(\x_{-k})$. Suppose now that such a
  term does, indeed, appear, and let $B\in\piBig(\x_{-k},\y_{-k+i})$
  be a corresponding homology class with a pseudo-holomorphic
  representative. Our goal now is to study properties of this domain.
  
  \step\label{step:D23-pinching}
  \emph{Pinching $\Sigma$ induces a decomposition of the domain $B$.}

  The boundary of the winding region~$\Winding$
  consists of two circles $C_1$ and $C_2$ parallel to the meridional
  $\alpha$-circle, where $C_1$ is 
  the component on the side of the regions $2$
  and~$3$.  Fix a complex structure $J$ and, as in the definition of
  ``sufficiently pinched,'' consider degenerating $J$ (and so
  $\Sigma$) along $\bdy\Winding$. 
  Suppose that there is a
  $\varphi\in\ModFlow(\x_{-k},\y_{-k+i})$ with multiplicity $+1$ in regions
  $2$ and $3$ and $0$ in regions $0$ and $1$.  We will show
  by positivity
  considerations that $i\geq -1$; in fact, we will see that
  the homology class~$B$ can be constructed as a kind of connected sum
  of a domain $B_0\in\piBig(\x,\y)$ with $n_{w}(B_0)=i+1$ and
  $n_{z}(B_0)=0$ with a homology class $B_1\in\piBig(x_{-k},x_{-k+i})$.

  More precisely, let $B\in\piBig(\x_{-k},\y_{-k+i})$ be a homology
  class giving a term in $D_{23}(\x_{-k})$.
  The domain~$B$ has multiplicity zero at the regions $0$ and~$1$, no corners
  outside of the winding region on $\beta_g$ or~$\alpha_1^a$, and in fact,
  both of its corners on~$\beta_g$ (and on~$\alpha_1^a$)
  occur in the negative half of the winding region.
  (Specifically the corners occur at $x_{-k}$ and $x_{-k+i}$, where both
  $-k$ and $-k+i$ are negative; here we are using
  Step~\ref{step:D23-no-positive}.)
  It follows that  
  $B$ induces a two-chain $B_0$ on
  the destabilized surface $\Sigma_0$, having local multiplicity zero
  at $z$, and some multiplicity~$j$
  at the point~$w$ (which in turn
  corresponds to $C_1$, shrunk to a point). 

\begin{figure}
    \begin{center}
      \includegraphics[scale=.83333]{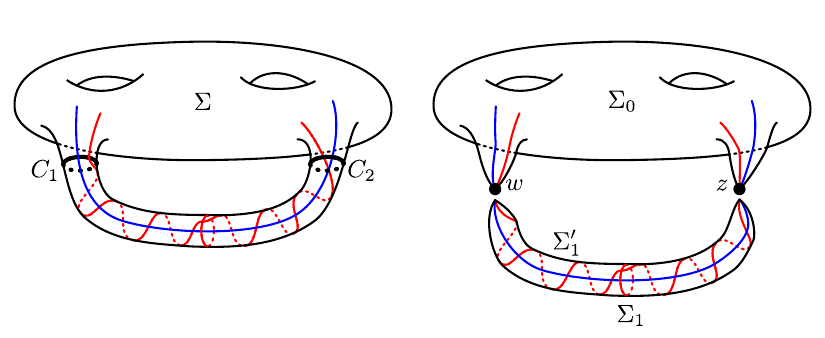}
    \end{center}
    \caption[Pinching off the winding region]
    {\label{fig:PinchingWindingWinding} \textbf{Pinching off the
        winding region.} The Heegaard surface is degenerated at the curves
      $C_1$ and $C_2$ to give a nodal surface with two components: one of these is
      the de-stabilized knot diagram
      $\Sigma_0$ (equipped with a pair of arcs $\alpha_g$ and $\beta_g$ connecting the
      $w$ and $z$ punctures)  and the other as a sphere, which contains the winding region.}
  \end{figure}

  \begin{figure}
    \begin{center}
      \input{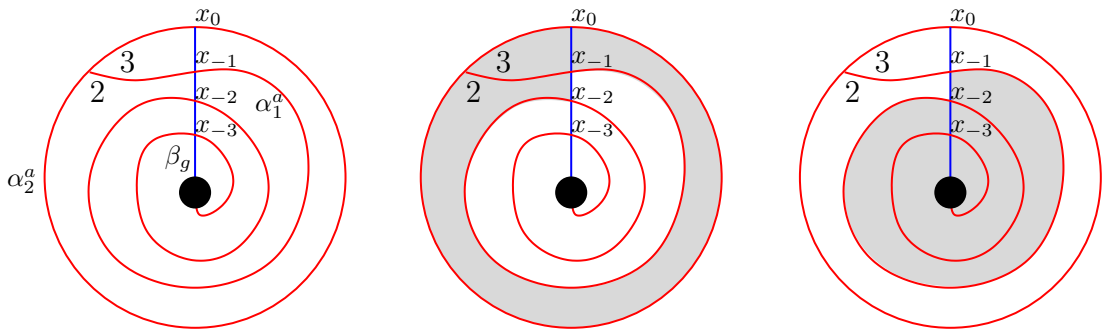}
    \end{center}
    \caption[Homology classes in the
        winding region]
    {\label{fig:DegenerateWinding} \textbf{Homology classes in the
        winding region.} The winding region to the left of
      the meridian, from Figure~\ref{fig:Twisting}, is drawn inverted so that
      the outer boundary of the winding region corresponds to the dark
      marked point, and the meridional arc corresponds to the outer
      boundary of the picture.  In the center and right we illustrate two
      homology classes: $\phi_1\in\piBig(\x_{-1},\x_{-2})$ and its
      complement, $\psi_1\in \piBig(\x_{-2},\x_{-1})$ from the proof
      of Lemma~\ref{lem:OrderGens}. The regions are labeled $R_i$ in
      order of increasing distance from the region containing~$\rho_3$.}
  \end{figure}

  The domain at
  the winding region induces another two-chain, supported in a
  disk.  Specifically, let $\Winding_L$ be the part of the winding
  region to the left of $\alpha_2^a$, so $\Winding_L$ contains the
  regions $2$ and $3$. As we degenerate our curve $C_1$ to a point,
  $\Winding_L$ becomes a disk $\Sigma'_1$ whose outer boundary is $\alpha^a_2$,
  the meridian of
  the knot from Figure~\ref{fig:Twisting}, with a preferred central
  point~$\delta_1$, the node obtained by degenerating $C_1$. See
  Figure~\ref{fig:DegenerateWinding}.  The
  two-chain~$B_1$ induced from~$B$ must have local multiplicity~$j$
  at~$\delta_1$.
  There are two distinguished arcs in~$\Sigma'_1$ connecting $\delta_1$ to the
  boundary. One is a portion of~$\beta_g$. The other is a portion
  of~$\alpha^a_1$
  which winds around, meeting
  $\beta_g$ in a sequence of intersection points
  $\{x_k\}_{k=-n/2}^0$.  Recall that for each $k$ there is a homology class
  $\gls*{phik}$
  connecting $x_k$ to $x_{k-1}$, with multiplicity~$+1$ at the
  regions $2$ and $3$, and multiplicity zero at~$\delta_1$. There is
  also a homology class
  $\gls*{psik}$
  connecting
  $x_{k-1}$ to $x_k$, which can be thought of as the complement of
  $\phi_k$. It is easy to see that a positive domain in this disk
  region which connects intersections $x_k$ and $x_{k+i}$ of the
  $\alpha$- and $\beta$-arcs and has multiplicity one in regions $2$
  and~$3$ necessarily decomposes as a composite $\phi_k * \psi_k * \psi_{k+1} *
  \dots * \psi_{k+i}$.  In particular, for our initial
  $B\in\piBig(\x_{-k},\y_{-k+i})$, we get
  $B_1\in\piBig(x_{-k},x_{-k+i})$, which we decompose as
  $\phi_k * \psi_k * \psi_{k+1} * \dots * \psi_{k+i}$. Consequently,
  the multiplicities of $B_1$ near the boundary of the winding
  region are all equal (and indeed equal to~$i+1$). Thus, the homology
  class $B_0$ on $\Sigma_0$ has multiplicity $j = i+1$ at~$w$.

  \step\label{D23:degenerate}
  \emph{Holomorphic curves have constrained behavior under the degeneration.}
  We turn next to the behavior of holomorphic curves as we stretch the
  neck along $\bdy\mathcal{W}$. By an
  argument analogous to the proof of Proposition~\ref{prop:Gluing1}, a
  sequence of curves $u_r$ holomorphic with respect to
  neck lengths $r\to\infty$ converges to a pair $(u_0,u_1)$ consisting of a
  holomorphic curve (or comb) $u_0$ in $\Sigma_0\times[0,1]\times\RR$
  and a holomorphic curve (or comb) $u_1$ in
  $\mathcal{W}\times[0,1]\times\RR$, where we view $\mathcal{W}$ as a
  sphere with three punctures, corresponding to the two circles where
  $\mathcal{W}$ is joined to the rest of $\Sigma$ and the puncture $p$
  in $\mathcal{W}$ (i.e., $e\infty$ of $\Sigma$). (For the purposes of some of the arguments
  here, it is useful to bear in mind that $\Sigma_0$ is equipped with a little extra structure now,
  in the form of two arcs, remnants of $\alpha_g$ and $\beta_g$, which connect the punctures
  $w$ and~$z$.)
  For notational
  convenience, write $\Sigma_1$ to denote $\mathcal{W}$ viewed as a
  three-punctured sphere.  Let $\delta_1$ denote the puncture of $\Sigma_1$
  corresponding to the curve $C_1$.

  The maps $u_i$ are analogous to the matched pairs discussed in
  Section~\ref{sec:moduli-matched-pairs}
  (Definition~\ref{def:weak-matched}), but can have Reeb orbits in
  addition to Reeb chords. More precisely, in the case that $u_0$ and
  $u_1$ are curves rather than combs,
  \begin{itemize}
  \item $u_0$ is a holomorphic map $S_0\to
    \Sigma_0\times[0,1]\times\RR$ (with respect to an appropriate
    almost complex structure $J$ on $\Sigma_0\times[0,1]\times\RR$), where
    $S_0$ is a surface with boundary, boundary punctures, and interior
    punctures. At each boundary puncture, $u_0$ is asymptotic to one
    of:
    \begin{itemize}
    \item a chord in $\x\times[0,1]$ at $-\infty$ or $\y\times[0,1]$
      at $+\infty$,
    \item a chord in $C_1$ starting and ending on
      $\alpha_1^a\cup C_1$,
    \item a chord in $C_1$ starting and ending on
      $\beta_g\cap C_1$, or
    \item a chord in $C_1$ starting on $\beta_g$ and ending on
      $\alpha_1^a$ or starting on $\alpha_1^a$ and ending on $\beta_g$.
    \end{itemize}
    We call punctures of these types \emph{$\pm\infty$ punctures},
    \emph{$\alpha$ boundary punctures}, \emph{$\beta$ boundary
      punctures} and \emph{mixed boundary punctures}, respectively. 

    At each interior puncture, $u_0$ is asymptotic to the circle
    $C_1$ at some point $(x,y)\in[0,1]\times\RR$ (i.e.,
    $\pi_\Sigma\circ u\to \delta_1$ and $\pi_\DD\circ u\to (x,y)$); or
    possibly to a multiple cover of $C_1$.
  \item $u_1$ is a map $S_1\to \Sigma_1\times[0,1]\times\RR$ (with
    respect to an appropriate almost complex structure $J$ on
    $\Sigma_1\times[0,1]\times\RR$), where
    $S_1$ is a surface with boundary, boundary punctures, and interior
    punctures. At each boundary puncture, $u_1$ is asymptotic to one
    of:
    \begin{itemize}
    \item a chord in $x_i\times[0,1]$ at $\pm\infty$,
    \item a chord $\rho$ at $e\infty$ (at some point $(1,t)\in[0,1]\times\RR$).
    \item a chord in $C_1$ starting and ending on
      $\alpha_1^a\cup C_1$,
    \item a chord in $C_1$ starting and ending on
      $\beta_g\cap C_1$, or
    \item a chord in $C_1$ starting on $\beta_g$ and ending on
      $\alpha_1^a$ or starting on $\alpha_1^a$ and ending on $\beta_g$.
    \end{itemize}
    We call punctures of these types \emph{$\pm\infty$ punctures},
    $e\infty$ punctures, \emph{$\alpha$ boundary punctures},
    \emph{$\beta$ boundary punctures} and \emph{mixed boundary
      punctures}, respectively.

    At each interior puncture, $u_0$ is asymptotic to the circle
    $C_1$ at some point $(x,y)\in[0,1]\times\RR$ (i.e.,
    $\pi_\Sigma\circ u\to \delta_1$ and $\pi_\DD\circ u\to (x,y)$); or
    possibly to a multiple cover of $C_1$.
  \item There is a bijection $\phi_a$ between the $\alpha$ boundary
    punctures of $S_0$ and the $\alpha$ boundary punctures of $S_1$ so
    that for each boundary puncture $q_i^a$ of $S_0$, $u_0$
    $t(u_0(q_i^a))=t(u_1(\phi_a(q_i^a)))$; and similarly for the
    $\beta$ boundary punctures.
  \item There is a bijection $\phi_c$ between the interior punctures
    of $S_0$ and the interior punctures of $S_1$ so that for each
    interior puncture $q_i^c$ of $S_0$,
    $\pi_\DD(u_0(q_i^c))=\pi_\DD(u_1(\phi_c(q_i^c)))$.
  \end{itemize}
  More generally, $u_0$ and $u_1$ can be combs---the obvious analogues
  of combs from Section~\ref{sec:combs-compact}, but allowing Reeb
  orbits. (Indeed, in the presence of mixed boundary punctures, $u_1$
  is forced to be a holomorphic comb with at least two vertical levels.)

  \step\label{step:D23-no-Reeb}
  \emph{Holomorphic curves with Reeb chords do not occur in
    the limit.}
  We show next that the limiting curves $(u_0,u_1)$ have no mixed
  boundary punctures, $\alpha$ boundary punctures or $\beta$ boundary
  punctures. 

  If $u_1$ has a mixed boundary puncture then $u_1$ necessarily has at
  least two stories. Let $u_{1,1}$ be the first story of $u_1$
  that contains a mixed boundary puncture, so $u_{1,1}$ is asymptotic
  to some $x_{-i}\times[0,1]$ at $-\infty$ and to a chord $\gamma_1$
  connecting $\alpha_1^a$ to $\beta_g$ at $+\infty$. 
  Label the regions in the winding region $\{R_i\}_{i=0}^{n/2+1}$, in
  order of increasing distance from the region containing the Reeb
  chord~$\rho_3$, as on the left of Figure~\ref{fig:DegenerateWinding}, and 
  let the local multiplicity of $u_{1,1}$ in $R_i$ be~$a_i$.  Note
  that $a_0$ and $a_1$
  coincide with the local multiplicities of $u_{1,1}$ around~$\rho_3$
  and $\rho_2$, respectively. In particular, $a_0,a_1\in\{0,1\}$. We
  rule out the four possibilities for $a_0$ and $a_1$ in turn.

  If $a_0=1$ and $a_1=0$ then the domain $u_{1,1}$ consists entirely
  of the region $3$. This contradicts the fact that $u_{1,1}$ is
  asymptotic to a mixed chord at $+\infty$. If $a_0=0$ and $a_1=1$
  then $a_j=j$ for $j=1,\dots, i$ and $a_{j}=i$ for $j>i$ (since
  $u_{1,1}$ has a $-\infty$ corner at $x_{-i}$ and no other corners in
  the interior of $\Winding$). But the fact that $u_{1,1}$ is
  asymptotic to a mixed chord at $+\infty$ implies that
  $a_{n/2}=a_{n/2-1}-1$, a contradiction.  If $a_0=a_1=0$ then
  $0=a_0=a_1=\dots=a_{i}$. Since $x_{-i}$ is an initial corner, it
  follows that $a_{i+1}=-1$, a contradiction to the fact that
  $u_{1,1}$ is holomorphic. Finally, if $a_0=a_1=1$ then $a_j=1$ for
  $j\leq i$, $a_{i+1}=0$ and $a_{i+2}=-1$, which again contradicts the
  fact that $u_{1,1}$ is holomorphic.  (Here we use the assumption
  that $n$ is sufficiently large, to guarantee that $k+2 < n/2$ from
  the assumption on $S(\x_i)$.) It follows that $u$ has no mixed
  chords.

  Next, suppose $u_0$ has an $\alpha$ boundary puncture. Then the
  source $S_0$ of $u_0$ has some boundary component which is mapped to
  $\alpha^a_1$. Since $\x$ and $\y$ have no components in
  $\alpha_1^a\cap\Sigma_0$, $t\circ u_0$ remains bounded on this
  component of $\bdy S_0$. By the maximum principle, it follows that
  $\pi_\DD\circ u_0$ is constant on this component of $S_0$. In
  particular, the boundary of this component of $S_0$ is contained
  entirely in $\alphas$. 
  Thus, $\pi_{\Sigma_0} \circ u_0$ gives a homological relation between $\alphas\cap \Sigma_0$.
  Since these curves and arcs are homologically linearly independent, and the two-chain has multiplicity $0$ at $z\in \Sigma_0$
  (as established in Step~\ref{step:D23-pinching}), it follows that
  all its local multiplicities vanish.
  (See also
  Lemma~\ref{lem:NoBoundaryDegenerations} for a similar argument.)
  Thus, $u_0$ (and hence $u_1$) has no $\alpha$ boundary punctures. A
  similar argument shows that $u_0$ and $u_1$ have no $\beta$ boundary
  punctures.
 
  \step\label{step:D23-dim-constraint}
  \emph{Dimension constraints on holomorphic curves with
    Reeb orbits rule out all holomorphic curves which are not supported in the winding region.}  In sum, the
  holomorphic curves $u_0$ and $u_1$ are asymptotic to some collection
  of Reeb orbits (copies of $C_1\times(x,y)$, or multiple covers of
  these) at $\delta_1$. Generically, all of these Reeb orbits will be
  simple (not multiply covered), so it suffices to consider this case.
  As in Step~\ref{step:D23-pinching}, we decompose $B_1$ as 
  $B_1=\phi_k*\psi_k*\cdots*\psi_{k+i-1}$, where now $i$
  is the number
  of Reeb orbits of $u_0$ and $u_1$. (In the case  $i=0$ there are no
  $\psi_j$ in the decomposition and no Reeb orbits.) The (expected)
  dimension of the
  moduli space of curves containing $u_1$ is $\ind(u_1)=2i+1$, and the
  (expected) dimension of the
  moduli space of curves near $u_0$ is given by $\ind(B_0)$. Finally,
  each Reeb orbit contributes a two-dimensional matching condition
  between $u_0$ and $u_1$, so
  \begin{align*}
    \ind(B,(\rho_2,\rho_3))&=\ind(u_0)+\ind(u_1)-2i\\
    &=\ind(B_0)+1.
  \end{align*}
  Thus, since $\ind(B,(\rho_2,\rho_3))=1$, $\ind(B_0)=0$, so $u_0$ is
  a trivial holomorphic curve. It follows that $i=0$, and $B$ is
  exactly the class $\phi_k$, showing that the only possible curve is
  the one analyzed in Step~\ref{step:D23-xk1-occurs}.
  \end{steps}

  This completes the proof for $D_{23}$.

  {\bf{The case of $D_{01}$.}} This case is symmetric to the case of $D_{23}$,
  via a rotation of the winding region (which also reverses the sign of the indexing
  of the $x_i$ and the sign of the function $S$).

  {\bf{The case of $D_{3}$.}} Suppose that $\phi\in\piBig(\x_0,\y)$ is
  a positive domain with local multiplicity one at the region marked
  with $3$ and zero at the regions marked by $0$, $1$ and $2$. 
  Considering local
  multiplicities around the intersection point $x_{-1}$, we see that
  $x_{-1}$ is forced to be a corner point, and indeed that the
  intersection of $\phi$ with the winding region $\Winding$
  is forced to be the bigon domain from $x_0$ to
  $x_{-1}$. It follows (as in Step~\ref{step:D23-no-yk1} of the
  analysis of $D_{23}$) that
  for generic almost-complex structures, the only positive, index one
  $\phi\in\piBig(\x_0,\y)$ which admits a pseudo-holomorphic
  representative is this bigon, connecting $\x_0$ to $\y=\x_{-1}$.

  {\bf{The case of $D_{1}$.}} This is symmetric to the case of $D_3$
  via a rotation of the winding region.
\end{proof}

\begin{lemma}
  \label{lem:VanishingHomotopies}
  For any sufficiently large $n$, if $\y\in\Gen(n)$ satisfies
  $|S(\y)|<\frac{n}{4}$ then
  $D_{230}(\y)=D_{2}(\y)=D_{0}(\y)=D_{012}(\y)=0$.
\end{lemma}

\begin{proof}
  By the same argument as used to prove Lemma~\ref{lem:DModDegrees}, the grading shifts of these maps are given by
  \begin{align*}
    D_{0}\co W^1_s&\to W^0_{s-\frac{n+1}{2}} &
    D_{2}\co W^1_s&\to W^0_{s+\frac{n+1}{2}} \\
    D_{230}\co W^1_s&\to W^0_{s-\frac{n-1}{2}} &
    D_{012}\co W^1_s&\to W^0_{s+\frac{n-1}{2}}.
  \end{align*}
  Thus, 
  if $|S(\y)|<\frac{n}{4}$ then, writing $D_*$ for any of the above maps, 
  $|S(D_*(\y))|\geq \frac{n-2}{4}$.
  By Equation~\eqref{lem:NormalizeS}, the $S$-grading of any generator $\x_0$ of $W^0$ coincides with
  the Alexander grading of the corresponding $\x\in\Gen_K$. 
  Thus, taking $n\geq 4\max_{\x\in\x_S} |A(\x)|+2$, $D_*(\y)$ lies in
  the zero group.
\end{proof}

\begin{lemma}
  \label{lem:ProvincialDomains}
  The following holds for any sufficiently pinched almost complex
  structure on $\Sigma\times[0,1]\times\RR$.
  For each positive domain $B\in\piBig(\x,\y)$ in $\HD_K$ with
  $\x,\y\in\Gen_K$, $n_w(B)=k\geq 0$, and $n_z(B)=0$, there is a
  corresponding sequence of domains $B_i\in\piBig(\x_{i-k}, \y_{i})$
  for $k-\frac{n}{2}\leq i\leq 0$ such that:
  \begin{itemize}
  \item $B_i$ is positive.
  \item If $i<0$ or $k=0$ then $B_i$ is provincial, and if $i=0$ and
    $k>0$ then $B_i$ crosses region~$2$ with multiplicity one and has
    multiplicity $0$ at regions $0$, $1$ and $3$.
  \item If $i<0$ or $k=0$ then $\ind(B_i,())=\ind(B)$, and if
    $i=0$ and $k>0$ then $\ind(B_i,(\{\rho_2\}))=\ind(B)$.
  \item Outside the winding region, $B_i$ agrees with $B$ (in the
    obvious sense).
  \end{itemize}
  In the other direction, suppose $B'\in\piBig(\x_{i-k},\y_i)$
  satisfies $n_0(B')=n_1(B')=n_3(B')=0$ and $n_2(B')\in\{0,1\}$. Let
  $\vec{\rho}(B')$ be the empty sequence of sets of Reeb chords if
  $n_2(B')=0$ and the sequence $(\{\rho_2\})$ if $n_2(B')=1$. 
  Suppose that
  \begin{align*}
    \eDim(B', \vec\rho(B'))&=1 \\
    \UnparModFlow(B', \vec\rho(B'))&\neq \emptyset.
  \end{align*}
  Then $B'=B_i$ for some $B\in\piBig(\x,\y)$, and
  $\#\UnparModFlow(B',\vec\rho(B'))=\#\UnparModFlow(B)$.

  Similarly, for each positive domain $B\in\piBig(\x,\y)$ in $\HD_K$ with
  $\x,\y\in\Gen_K$, $n_z(B)=k\geq 0$, and $n_w(B)=0$, there is a
  corresponding sequence of domains $B_i\in\piBig(\x_{i+k}, \y_{i})$
  for $0\leq i\leq \frac{n}{2}-k$ such that:
  \begin{itemize}
  \item $B_i$ is positive.
  \item If $i>0$ or $k=0$ then $B_i$ is provincial, and if $i=0$ and
    $k>0$ then $B_i$ crosses region~$0$ with multiplicity one and has
    multiplicity $0$ at regions $1$, $2$ and $3$.
  \item If $i>0$ or $k=0$ then $\ind(B_i,())=\ind(B)$, and if
    $i=0$ and $k>0$ then $\ind(B_i,(\{\rho_0\}))=\ind(B)$.
  \item Outside the winding region, $B_i$ agrees with $B$ (in the
    obvious sense).
  \end{itemize}
  In the other direction, suppose $B'\in\piBig(\x_{i+k},\y_i)$
  satisfies $n_1(B')=n_2(B')=n_3(B')=0$ and $n_0(B')\in\{0,1\}$. Let
  $\vec{\rho}(B')$ be the empty sequence of sets of Reeb chords if
  $n_0(B')=0$ and the sequence $(\{\rho_0\})$ if $n_0(B')=1$. 
  Suppose that
  \begin{align*}
    \eDim(B', \vec\rho(B'))&=1 \\
    \UnparModFlow(B', \vec\rho(B'))&\neq \emptyset.
  \end{align*}
  Then $B'=B_i$ for some $B\in\piBig(\x,\y)$, and
  $\#\UnparModFlow(B',\vec\rho(B'))=\#\UnparModFlow(B)$.
\end{lemma}
\begin{proof}
  Given $B$, the domain $B_i$ is constructed using the procedure
  illustrated by the shaded region in Figure~\ref{fig:Twisting}; $B_i$
  is uniquely determined by its corners $x_{i\pm k}$ and $x_i$, its
  local multiplicities at the regions $0,1,2,3$, and the
  condition that $B_i$ agrees with $B$ outside the winding region. It
  is immediate from the index formulas that $B_i$ has the same index
  as $B$.
 
  Suppose next that $B'\in\piBig(\x_{-s},\y_{-s+t})$ for $s, t \ge 0$
  is an
  index one homology class with $n_0(B')=n_1(B')=n_3(B')=0$ which has
  a holomorphic representative. We
  claim that $B'=B_{-s}$ for some $B\in\piBig(\x,\y)$ with
  $n_{w}(B)=k$, $n_{z}(B)=0$, and
  $\#\UnparModFlow(B',\vec\rho(B'))=\#\UnparModFlow(B)$.
  
  By stretching along the boundary of the winding region, as in the proof of
  Lemma~\ref{lem:OrderGens}, we see that $B'$ can be written as a
  connected sum of a homology class $B_0$ with a standard chain
  $\psi_{s,t}\in\piBig(x_{-s},x_{-s+t})$. We have that
  $\eDim(\psi_{s,t}, \vec\rho(B'))=2t$, and
  $\eDim(B')=\eDim(B_0)+\eDim(\psi_{s,t})-2t$. Thus,
  $\eDim(B_0)=1$.

  We claim the moduli space is obtained as a fibered product
  $\ModFlow(\varphi_0)\times_{\Sym^t(\DD)} \ModFlow(\psi_{s,t})$.  As
  in the proof of Lemma~\ref{lem:OrderGens}, this follows from
  compactness and gluing results, together with the following argument
  to rule out curves asymptotic to Reeb chords at the boundary of
  $\Winding$ (so the only asymptotics at $\bdy\Winding$ are Reeb
  orbits). Chords with both endpoints on $\alpha$ or both endpoints on
  $\beta$ are ruled out exactly as in the proof of
  Lemma~\ref{lem:OrderGens}: such chords force components of the curve
  in $\Sigma_0$ to be constant. Mixed Reeb chords are ruled out by
  considering the local multiplicities at the regions in $\Winding$,
  similarly to the proof of Lemma~\ref{lem:OrderGens}; for the
  bottom-most story containing a mixed puncture the local
  multiplicities $a_i$ are all $1$ for $i$ sufficiently large, which
  is inconsistent with being asymptotic to a mixed Reeb chord at $+\infty$.

  Having identified our moduli space with the
  fiber product $\ModFlow(B_0)\times_{\Sym^t(\DD)}
  \ModFlow(\psi_{s,t})$, we next claim that the map
  $\ModFlow(\psi_{s,t})\to \Sym^t(\DD)$, gotten by projecting
  the preimage of the connect sum point~$q$, has degree one (onto
  $\Sym^t(\DD)$), so that the moduli space has the same count modulo
  two as $\ModFlow(B_0)$, as desired.

  In the case where $t=1$, this can be seen directly. Consider the
  moduli space~$\psi_{s,1}$, which, in the notation of
  the proof of Lemma~\ref{lem:OrderGens}, is $\psi_{s-1}$; see also the right of
  Figure~\ref{fig:DegenerateWinding}. Our claim is that the map of
  $\ModFlow(\psi_{s-1})$ to $[0,1]\times \RR$ has degree one or,
  equivalently, after we divide out the space $\ModFlow(\psi_{s-1})$ by
  translations, the preimage of~$q$ is mapped with degree
  one onto the interval $[0,1]$. But the quotient of
  $\ModFlow(\psi_{s-1})$ by translations can be parameterized by depth of
  the cut at $x_{-s}$. In the limit as the cut goes out to the
  marked point along the $\alpha$-circle, the preimage of~$q$ is
  mapped (under projection to $[0,1]$) to $\{1\}$; if the cut
  goes out along the $\beta$-circle, then the preimage of~$q$ is
  mapped to $\{0\}$. Hence the evaluation map has degree one
  when $t=1$, as stated.

  For $t>1$ we show that the evaluation map has degree one as
  follows. By a straightforward analysis of domains in the winding
  region (see Step~\ref{step:D23-pinching} in the proof of
  Lemma~\ref{lem:OrderGens}),
  $\psi_{s,t}$ can be written as a
  juxtaposition of flowlines with Maslov index two. In fact, any
  decomposition of $\psi_{s,t}$ consisting of non-negative and
  non-constant homology classes must necessarily be a decomposition
  into homology classes each of which has non-zero multiplicity
  at~$q$.  It follows that the map
  $\ModFlow(\psi_{s,t})\to \Sym^t(\DD)$ is proper. To
  calculate its degree, consider the fibers of points in
  $\Sym^t(\DD)$ where one of the components has $\RR$ coordinate much
  larger than all the others. Such holomorphic disks are in
  the image of the gluing map of a moduli space with multiplicity $1$
  with another moduli space with multiplicity $t-1$. By induction, it
  follows that the evaluation map has degree one.
\end{proof}

\begin{lemma}\label{lem:D23-big-s} For any sufficiently pinched almost
  complex structure on $\Sigma\times[0,1]\times\RR$, 
  for $s\geq \frac{n}{4}$ the map
  $
  D_{23}\co W^1_{s}\to W^1_{s+1}
  $ 
   is as given by Theorem~\ref{thm:HFKtoHFD2}.
\end{lemma}
\begin{proof}  
  We introduce some preliminary notation. A curve $w$ representing $D_{23}$ has local multiplicity
  zero at the regions marked by $0$ and $1$. Thus, if we let $\Omega \subset \Sigma$ denote the complement in $\Sigma$
  of the regions marked by $0$ and $1$, the source of $w$ maps to $\Omega \times [0,1]\times \RR$.

  We prove the lemma in nine steps. 

  \begin{steps}
  \step \emph{$D_{23}(\x_1)=0$ for any $\x\in\Gen_K$, since there
  are no domains which can represent such a map.}  On the one hand, any
  domain with initial corner at $x_1$ and vanishing local
  multiplicities at $0$ and $1$ must in fact vanish at all four
  regions near $x_1$. On the other hand, a domain which represents
  $D_{23}$ must have a cut (along $\alpha_1^a$) which goes out the boundary (separating
  regions $2$ and $3$). Thus, it cannot have vanishing multiplicities around $x_1$.

  \step \emph{Degenerate the surface to extract a boundary degeneration.}
  Consider $D_{23}(\x_i)$ for some $1<i$. We are interested in
  terms $\y_j$ in $D_{23}(\x_i)$ where $j>0$. We use a degeneration
  argument similar to the proof of Lemma~\ref{lem:OrderGens}. As in
  the definition of sufficiently pinched almost complex structures let
  $C_2$ be the boundary component of the winding region on the $0,1$
  side (i.e., the right side as shown in
  Figure~\ref{fig:PinchingWindingWinding}). Pinch the surface $\Omega\subset \Sigma$ along $C_2$.
  The resulting nodal surface has two irreducible components: a
  topological disk which we denote $\Winding_R$
  corresponding to part of the winding region; and another component
  $\Sigma'$ corresponding to most of the diagram and
  containing the regions $2$ and $3$. See
  Figure~\ref{fig:PinchingHalfWinding}. Let $\delta_2$ denote the node (or
  puncture) corresponding to $C_2$ after the degeneration. Holomorphic
  curves giving terms in $D_{23}(\x_i)$ degenerate to pairs of
  holomorphic curves $(u,v)$ where $u$ is a curve in
  $\Winding_R\times[0,1]\times\RR$ and $v$ is a curve whose interior
  lies in $\Sigma'\times[0,1]\times\RR$.

  \begin{figure}
    \begin{center}
      \includegraphics[scale=.83333]{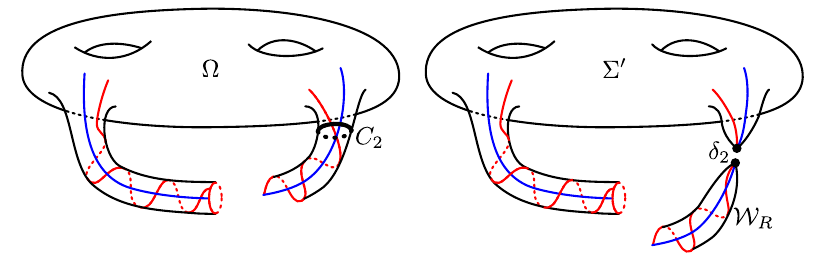}
    \end{center}
    \caption[Pinching $\Omega$ along $C_2$]
    {\label{fig:PinchingHalfWinding} \textbf{Pinching $\Omega$ along $C_2$.} 
      This is half the pinching from Figure~\ref{fig:PinchingWindingWinding}
      (done on $\Omega\subset \Sigma$).}
  \end{figure}

  By hypothesis, the components of $\x_i$ and $\y_j$ on $\alpha_1^a$
  both lie in $\Winding_R$, and $v$ has punctures mapped to the Reeb
  chords $\rho_2$ and $\rho_3$. Let $S_0$ be the component of (the
  source of) $v$ part of the boundary of which is mapped to
  $\alpha_1^a$. Then one component of $\bdy \overline{S_0}$
  is mapped entirely to $\alpha_1^a$, and the $t$-coordinate of $v$
  remains bounded on this component. Consequently, as in the proof of
  Lemma~\ref{lemma:arc_monotonicity}, $\pi_\DD\circ v$ is constant on
  $S_0$, so $v|_{S_0}$ is a boundary degeneration.

  Write the source of $v$ as $S=S_0\amalg S_1$.

    \step \emph{The curve $v$ consists of the boundary degeneration and
  trivial strips.} By Proposition~\ref{prop:asympt_gives_chi},
  $D_{23}$ counts domains $B$ with 
  \begin{equation}
    \label{eq:MaslovIndexOne}
    1/2=e(B)+n_{\x}(B)+n_{\y}(B).
  \end{equation}
  Let $B_2$ be the domain in $\Winding_R$ determined by $u$ and $B_1$
  the domain in $\Sigma'$ determined by $v|_{S_1}$. Let
  $n_{\delta_2}(B_i)$ denote the local multiplicity of $B_i$ at
  $\delta_2$. Then,
  \begin{equation}
    \begin{split}
      n_{\delta_2}(B_2)&=n_{\delta_2}(B_1)+1\\
      n_{x_i}(B_2)+n_{\x}(B_1)& =n_{\x}(B)-(g-1)\\
      n_{x_j}(B_2)+n_{\y}(B_1)&= n_{\y}(B)-(g-1) \\
      e(B_2)+e(B_1)&= e(B)-5/2+2g+2n_{\delta_2}(B_2).
    \end{split}\label{eq:Bone-versus-B}
  \end{equation}
  (In these formulas, for the purposes of the Euler measure $e$, the
  point $\delta_2$ is treated as a node, not a puncture.)
  Combining Equations~\eqref{eq:Bone-versus-B}
  and~\eqref{eq:MaslovIndexOne} gives
  \begin{equation}
    e(B_1)+n_{\x}(B_1)+n_{\y}(B_1)+e(B_2)+n_{x_i}(B_2)+n_{x_j}(B_2)=2n_{\delta_2}(B_1)+2.\label{eq:another-D23-big-s-eq}
  \end{equation}  
  A direct analysis of the domains in the winding region shows that
  \[
  e(B_2)+n_{x_i}(B_2)+n_{x_j}(B_2)=2n_{\delta_2}(B_2).
  \]
  It follows that $e(B_1)+n_\x(B_1)+n_\y(B_1)=0$. This is the expected
  dimension of the parameterized moduli space of holomorphic curves in
  the domain $B_1$. It follows that $v|_{S_1}$ consists of trivial
  strips, so $B_1=0$ and $n_{\delta_2}(B_2)=1$. In particular, $\y=\x$.

  \step \emph{Analyzing $u$.} The homology class of $u$ is
  uniquely determined by the properties that its local multiplicities
  at the regions $0$ and $1$ are zero, $x_i$ is a corner, and its
  local multiplicity $\delta_2$ is $1$. It
  follows that $u$ is connects $x_i$ to $x_{i-1}$.  Combining this
  with the previous step, it follows that $\y_j=\x_{i-1}$. It also follows
  from the above considerations that there is a unique holomorphic
  curve $u$ in this homology class.

  At this point, we have shown that either $D_{23}(\x_i)=\x_{i-1}$ or
  $D_{23}(\x_i)=0$. It remains to exclude the latter case; this
  involves a gluing argument.

  \step \emph{Relaxing the almost complex structure.} We would
  like to assert that there is a unique holomorphic curve in the
  moduli space of $v$ with the Reeb chord in $C_2$ at a given
  $t$-coordinate, and that we can glue this curve $v$ to $u$ to get a
  holomorphic curve in $\Sigma$. Unfortunately, if $g>1$ then the
  curve $v$ is not transversally cut out for any choice of admissible
  almost complex structure in the sense of Definition~\ref{def:admissible_J}.

  Instead, we relax Condition~(\ref{item:J1}) of the definition of
  admissible complex structures, to~\cite[Condition
  (J5$^\prime$)]{Lipshitz06:CylindricalHF}, which we restate here as: 
  \begin{enumerate}[label=(J$^\prime$-\arabic*),ref=J$^\prime$-\arabic*]
  \item\label{Jprime1}\index{(J'-1)} There is a 2-plane distribution \(\xi\) on
    \(\Sigma\times[0,1]\) such that the restriction of \(\omega\) to
    \(\xi\) is non-degenerate, \(J\) preserves \(\xi\), and the
    restriction of \(J\) to \(\xi\) is compatible with \(\omega\).  We
    further assume that \(\xi\) is tangent to \(\Sigma\) near
    \((\alphas\cup\betas)\times[0,1]\).
  \end{enumerate}
  This is sufficient for \cite{BEHWZ03:CompactnessInSFT} to give
  compactness for the moduli spaces of holomorphic curves. 
  We call such almost-complex structures {\em relaxed}.
  
  Given an almost complex structure $J$ on
  $\Sigma\times[0,1]\times\RR$, let $\cM(v;\betas;J)$ denote the
  moduli space of $J$-holomorphic curves containing $v$.  For a
  generic almost complex structure satisfying (\ref{Jprime1}) rather
  than (\ref{item:J1}), the moduli spaces $\cM(v;\betas;J)$ is
  transversally cut out (and, in particular, does not contain any
  boundary degenerations).

  \step\label{step:big-s-independence}
  \emph{Independence of $\cM(v,\betas;J)$ from the
  $\beta$-curves.}  Let $\betas'$ be another set of $\beta$-curves that
  agree with $\betas$ near $\x$ and let $\cM(v;\betas';J)$ denote the
  moduli space of holomorphic curves with respect to $\betas'$
  corresponding to $\cM(v;\betas;J)$.
  We claim that for $J$ a generic almost complex
  structure satisfying condition (\ref{Jprime1}) close enough to an almost
  complex structure $J_0$ satisfying condition (\ref{item:J1}),
  $\#\cM(v;\betas';J)=\#\cM(v;\betas;J)$. Indeed, taking $J\to J_0$,
  elements of $\cM(v;\betas';J)$ converge to the unions of
  $\alpha$-boundary degenerations and trivial strips. Let $U_\x$ be a
  small neighborhood of $\x$ so that $\betas\cap U_\x=\betas'\cap
  U_\x$. It follows that for $J$ close enough to $J_0$ and $w\in
  \cM(v;\betas';J)$, $w(\bdy S)\subset \alphas\cup (\betas\cap U_\x)$;
  and similarly for elements of $\cM(v;\betas;J)$. Thus, for $J$ close
  enough to $J_0$, $\cM(v;\betas;J)=\cM(v;\betas';J)$.

  \step\label{step:big-s-model}
  \emph{A model computation when $g=1$.} In the special case
  where $g=1$, the moduli space $\cM(v;\betas;j_\Sigma\times j_\DD)$ is
  transversally cut out, and consists of a single point. This is
  analogous to Proposition~\ref{prop:east_transversality}.

  \step\label{step:big-s-stabilize}
  \emph{Stabilizing the model computation.} Let
  $\betas'=\{\beta'_1,\dots,\beta'_{g-1}\}$ be a
  $(g-1)$-tuple of curves so that $\beta'_i\cap\alpha_i$ is a single
  point and $\beta'_i\cap \alpha_j=\emptyset$ for $i\neq j$ (and the
  $\beta'_j$ are disjoint from $\alpha_1^a$ and $\alpha_2^a$). Choose
  a curve $\gamma$ so that $\Sigma'\setminus \gamma$ has two connected
  components, one containing $\betas'$ and the other containing
  $\alpha_1^a$. Stretching the neck along $\gamma$, $\cM(v;\betas';J)$
  breaks up as a fibered product of two moduli spaces over a disk. One is the
  moduli space from the model computation in Step~\ref{step:big-s-model}, which has a single, transversally cut-out
  point. The other maps degree~$1$ onto the disk, by the same argument
  used to prove stabilization invariance in the cylindrical setting
  \cite[Proposition A.3]{Lipshitz06:CylindricalHF}.

  \step \emph{Conclusion of the proof.}
  Steps~\ref{step:big-s-independence} and~\ref{step:big-s-stabilize} imply that
  $\cM(v;\betas;J)$ is transversally cut out, and consists
  of an odd number of points, for some relaxed almost-complex structure $J$, which we can take to be
  arbitrarily close to some fixed, admissible $J_0$ (sufficiently pinched).
  By a standard gluing argument, this moduli space can be glued to $u$ to give an odd number of holomorphic
  curves contributing to $D_{23}(\x_i)$. Thus, the number of points in
  the moduli space of $J$-holomorphic curves  $\cM(\x_i,\x_{i-1};(\rho_2,\rho_3);J)$
  is odd.  
  Let $J_0$ be a sufficiently pinched, 
  generic almost complex structure satisfying the more restrictive
  Condition~(\ref{item:J1}).
  Since the moduli space of $J_0$-holomorphic
  curves
  $\cM(\x_i,\x_{i-1};(\rho_2,\rho_3);J_0)$ is transversally cut out,
  it follows that
  $\#\cM(\x_i,\x_{i-1};(\rho_2,\rho_3);J_0)$ is odd, as well.
  Thus, $D_{23}(\x_i)=\x_{i-1}$, as desired.
  \end{steps}
\end{proof}

\begin{proof}[Proof of Theorem~\ref{thm:HFKtoHFD2}]
  We choose $n$ larger than $4c$, where $c$ is the constant from
  Lemma~\ref{lem:Exterior}, so that all exterior generators~$\x$ have
  $|S(\x)|< \frac{n-2}{4}$, and large enough that
  Lemmas~\ref{lem:OrderGens} and~\ref{lem:VanishingHomotopies} hold.
  
  Let $W=W^0\oplus W^1$ be the type $D$ module for
  $\CFDa(S^3\setminus\nbd(K))$.
  We see from Lemma~\ref{lem:NormalizeS} that $W^0_s=C(s)$, and its
  differential $D$ is clearly identified with the
  differential on knot Floer homology.
  
  Consider next generators $\x$ representing elements of $W^1$.  In the case
  where $|S(\x)|\geq \frac{n-2}{4}$, according to
  Lemmas~\ref{lem:NormalizeS} and~\ref{lem:Exterior}, all generators
  with the same $S$-grading are of the form $\x_{t}$ for
  $\x\in\Gen_K$, and
  $t$ uniquely determined by (the Alexander grading of) $\x$.
  This gives the identification of
  $W^1_s$ with $V^1_s$ for $\abs{s} > \frac{n}{4}$ as stated in the
  theorem.  For instance, for $s < -\frac{n}{4}$, $W^1_s$
  contains $\x_k$ for $S(\x) = s+\frac{n+1}{2}-k$ and can therefore be
  identified with $C\bigl(\leq s+\frac{n-1}{2}\bigr)$, as desired.
  Moreover,
  Lemma~\ref{lem:ProvincialDomains} then identifies the
  differential $D$ on $W^1_s$ with the stated differential $D$ on
  $V^1_s$. We return to studying $W^1$ with $|S(\x)|< \frac{n}{4}$ later in the proof. 
  We say what we can about the other coefficient maps first, though.

  Any domain representing $D_1$ must start at a generator of the
  form $\x_0$. By Lemma~\ref{lem:OrderGens}, the map $D_1$ is
  represented by the bigon from
  $\x_0$ to~$\x_1$. 
  Similarly, any domain representing $D_3$ is represented by the bigon
  from $\x_0$ to~$\x_{-1}$. We conclude that $D_1$ and $D_3$ are the maps from
  $\x_0$ to $\x_{1}$ and $\x_{-1}$, respectively, as desired. 

  Next, any domain representing $D_2$ must start at a generator of
  the form $\x_{s}$ with $s<0$ and end at a generator of the form
  $\y_0$.  By Lemma~\ref{lem:ProvincialDomains}, such a
  domain must be induced by some domain $\phi\in\piBig(\x,\y)$ with
  $n_z(\phi)=0$ and $n_w(\phi)=k\geq 1$.

  An analogous argument shows that $D_{123}$ has the stated form,
  except that we necessarily have $n_w(\phi) = 1$.
  
  There are no domains which could represent a map of the
  form $D_{12}$, for such a domain would need to have a cut going out
  to the intersection point $x_0$, which has the wrong idempotents.
  
  The fact that $D_{23}$ has the stated form on 
  $W^1_s$ for $s< -\frac{n}{4}$ was
  established in Lemma~\ref{lem:OrderGens}. The fact that 
  $D_{23}$ has the given form for $s>\frac{n}{4}$ follows from
  Lemma~\ref{lem:D23-big-s}. 

  We have now shown that all the coefficient maps have the right form
  for $W^1_s$ with $\abs{s}> \frac{n}{4}$ and for~$W^0_s$.
  We claim next that for all $|s|< \frac{n}{4}$,
  $H_*(W^1_s,D)\cong\Field$, and in fact 
  $$D_{23}\co
  W^1_{s}\to W^1_{s+1}$$ is a chain map (with respect to
  $D$) inducing an isomorphism on homology. We investigate the terms
  in Equation~\eqref{eq:Homotopy0123} and
  Equation~\eqref{eq:Homotopy2301}, bearing in mind
  Lemma~\ref{lem:VanishingHomotopies}. Specifically, for $\x$ with
  $|S(\x)|\leq \frac{n}{4}$, $D_{3}\circ D_{012}(\x)=0$,
  $D_{1}\circ D_{230}(\x)=0$, and  $D_{301}\circ D_2(\x) = 0$ for
  degree reasons (comparing
  Lemma~\ref{lem:DModDegrees} with the observation that $W^0_s=0$ for
  all $|s|\geq \frac{n}{4}$). Similarly, since
  $$D_0\co W^1_s\to W^0_{s-\frac{n+1}{2}},$$
  (by the
  same considerations as in Lemma~\ref{lem:DModDegrees}), it follows
  that $D_{123}\circ D_0(\x)=0$ for $|S(\x)|\leq \frac{n}{4}$. Thus
  $D_{01}$ and $D_{23}$ are homotopy inverses in this range, and so
  the map induced on homology by $D_{23}\co
  W^1_s\to W^1_{s+1}$ is an isomorphism.  
  Earlier in the proof we identified
  $W^1_{-\frac{n+2}{4}}$ with $C(\leqo n)=(C,\partial_w)$, 
  a chain complex with homology
  isomorphic to $\Field$. We can conclude that
  $$H_*(W^1_s,D)\cong\Field $$ for all $|s|< \frac{n}{4}$. Let 
  $\phi\co W^1_s\to \Field$ denote a (graded)
  quasi-isomorphism.

  It follows that for $|s|< \frac{n}{4}$, we have 
  a chain map $\phi\co W^1_s\to \Field$ inducing an isomorphism
  on homology. The map 
  $$\Phi\co W\to V$$
  given by 
  $$\Phi(\x)=
  \begin{cases}
    \phi(\x) & {\text{if $\x\in W^1_s$ with $|s|<\frac{n}{4}$}} \\
    \x & {\text{otherwise}}
  \end{cases}$$
  induces a homotopy equivalence between the type $D$ module 
  determined by $W$ and the one determined by $V$.

  The rest of the statement about gradings is clear from the
  construction.
\end{proof}

\begin{remark}
  \label{rmk:HFKtoHFD2Gradings}
  Most of the grading on $V$ can be read off immediately from the
  Heegaard diagram for $S^3\setminus\nbd(K)$: there is a one-to-one
  correspondence between $V^0_s$ and $W^0_s$ and between $V^1_s$ and
  $W^1_s$ as long as $|s|\geq \frac{n}{4}$.
\end{remark}

\begin{proposition}
  \label{prop:GradingHFKtoHFD}
  There is an element $s_0$ of the grading set with the property that for any generator $\y$
  for $\HFKa(K)$,
  the corresponding generator $\y_0$ for $\CFDa(S^3\setminus K)$ has grading given by
  \[ \gr(\y_0)=\lambda^{M(\y)-2A(\y)}\gr(\rho_{23})^{-A(\y)}\cdot s_0.\]
\end{proposition}

\begin{proof}
  Let $\x$ and $\y$ be any generators for $\HFKa(K)$.
  After switching the roles of $\x$ and $\y$, we can assume that
  $A(\x)\leq A(\y)$. We can find some ${B}\in
  \pi_2(\x,\y)$ (with $n_z({B})=0$). By
  Equation~\eqref{eq:AlexanderUpToTranslation}, $n_w({B})=A(\y)-A(\x)\geq 0$; also, by Equation~\eqref{eq:MaslovGrading},
  followed by the formula for the Maslov index~\cite[Corollary 4.3]{Lipshitz06:CylindricalHF},
  \begin{align*}
    M(\x)-M(\y) &= \Mas({B})-2 n_{w}({B}) \\
    &= e({B}) + n_{\x}({B})+ n_{\y}(B)
    -2 n_{\w}({B}).
  \end{align*}
  The domain ${B}$ induces a domain
  $\widetilde{B}\in\pi_2(\x_0,\y_0)$ with $n_z(\widetilde{B})=0$ and
  \begin{align*}
    n_{\x_0}(\widetilde{B}) &= n_{\x}(B) &
    \partial^{\partial} \widetilde{B}& =(0,n_{w}({B}),n_{w}({B})) \\
    n_{\y_0}(\widetilde{B}) &= n_{\y}({B}) &
    e(\widetilde{B}) &= e({B})-\frac{1}{2} n_{w}({B});
  \end{align*}
  so 
  \begin{align*}
    g(\widetilde{B})&=(-e(\widetilde{B})-n_{\x_0}(\widetilde{B})-n_{\y_0}(\widetilde{B});\partial^{\partial} \widetilde{B}) \\
    &= \left(-e({B})-n_{\x}({B})-n_{\y}({B})
    + \frac{1}{2} n_{w}({B}); \partial^{\partial}\widetilde{B}\right) \\
    &= \left(M(\y)-M(\x) -\frac{3}{2} n_{w}({B}); \partial^{\partial} \widetilde{B} \right) \\
    &= \left(M(\y)-\frac{3}{2}A(\y) -M(\x)+\frac{3}{2}A(\x); 0, n_{w}({B})\right).
  \end{align*}
  (In the expression for $\bdy^\bdy \widetilde{B}$ we list the multiplicities
  of $[\rho_1]$, $[\rho_2]$ and $[\rho_3]$. In the final equation, we are writing the group element in
  coordinates for~$G$
  as in Equation~\eqref{eq:grading-torus-alg}.)
  So,
  \begin{equation}
    \label{eq:TrefoilCFDGradings}
  \begin{aligned}
    \gr(\y) &= \gr(R(\widetilde{B}))\cdot \gr(\x) \\
    &= \left(M(\y)-\frac{3}{2}A(\y) -M(\x)+\frac{3}{2}A(\x); 0, -n_{w}({B})\right)\cdot\gr(\x) \\
    &= \lambda^{M(\y)-2A(\y)-M(\x)+2A(\x)} \gr(\rho_{23})^{A(\x)-A(\y)}\cdot \gr(\x).
  \end{aligned}
  \end{equation}
  Thus, 
  \[ \lambda^{-M(\y)+2A(\y)}\gr(\rho_{23})^{A(\y)}\cdot \gr(\y)
  = \lambda^{-M(\x)+2A(\x)}\gr(\rho_{23})^{A(\x)}\cdot \gr(\x).\]
  Since $\x$ and $\y$ were arbitrary, the result follows.
\end{proof}

\section{Proof of Theorem~\ref{thm:HFKtoHFD}}
\label{sec:HFKtoHFDproof}

We next deduce Theorem~\ref{thm:HFKtoHFD} from
Theorem \ref{thm:HFKtoHFD2}. But first, we show that $\CFKm(K)$ can be
horizontally or vertically reduced, as in
Definition~\ref{def:Simplified} (and as needed for the statement of Theorem~\ref{thm:HFKtoHFD}). 
To this end, we have the following preliminary lemma:

\begin{lemma}\label{lemma:map-filtered-bases}
   Let $V\!$ and $W\!$ be finite-dimensional filtered vector spaces and
   $f \co V \to W\!$ be a (possibly un-filtered) isomorphism between
   the underlying vector spaces. Then there is a filtered basis
   for~$V\!$ that maps under $f$ to a filtered basis for~$W\!$.
\end{lemma}

\begin{proof}
  The filtration $W_j\subset W$ of $W$ 
  induces a second filtration on $V$, by the subsets
  $f^{-1}W_j$, which, together with the original filtration on $V$
  gives a $(\ZZ\oplus\ZZ)$-filtration on $V$ by subcomplexes $V_{s,t}=V_s\cap f^{-1} W_t$.
  Correspondingly,  the {\em filtration level} of an element $v\in V$
  is defined to be the pair of integers $(s,t)$ where 
  $s=\inf\{i\in\ZZ\mid v\in V_i\}$, $t=\inf \{j\in\ZZ\mid f(v)\in W_j\}$.
  Let $\llbracket v \rrbracket$ denote the projection of $v$ to $V_{s,t}/(V_{s-1,t}+V_{s,t-1})$.
  We can form the associated $\ZZ\oplus\ZZ$-graded vector space $\gr_{\ZZ\oplus\ZZ}(V)=\bigoplus_{s,t} V_{s,t}/(V_{s-1,t}+V_{s,t-1})$.
  A {\em $(\ZZ\oplus\ZZ)$-filtered basis} is a basis for $V$ whose projection under $\llbracket\cdot\rrbracket$
  is a basis for $\gr_{\ZZ\oplus\ZZ}(V)$.
  
  Now, pick a $(\ZZ\oplus\ZZ)$-filtered basis for $V$ (which can be
  obtained by taking a basis for $\gr_{\ZZ\oplus\ZZ}V$ and lifting it
  to a basis for $V$). Such a basis is {\em a fortiori} a filtered basis
  for either the filtration by $V$ or the filtration induced by $W$; i.e., it
  can be thought of as a ($\ZZ$-)filtered basis for $V$ whose image under $f$ is
  a filtered basis of $W$, as needed.
\end{proof}

\begin{proposition}
  \label{prop:SimplifyComplex}
  Let $C$ be a $\ZZ$-filtered, $\ZZ$-graded, finitely-generated chain
  complex over $\Field[U]$ which is free as an $\Field[U]$-module.
  Then $C$ is $\ZZ$-filtered, $\ZZ$-graded
  homotopy equivalent to a chain complex $C'$ which is reduced. Further,
  one can choose a basis for $C'$ over $\Field[U]$ which is vertically simplified or,
  if one prefers, a basis which is horizontally simplified instead.
\end{proposition}

\begin{proof}
  We show that $C$ can be reduced by induction
  on its rank.  If $C$ is not reduced, by definition it admits 
  a non-trivial differential which does not change the Alexander
  filtration, so there is some $\xi\in C$ with $\partial\xi \ne 0$ and $A(\xi) =
  A(\partial\xi)$.
  Then $\xi$ and $\partial \xi$ generate a subcomplex of $C$ whose
  quotient complex $Q$ is filtered homotopy equivalent to $C$ with rank two
  less than the rank of~$C$.

  Next, we argue that $C$ can be vertically simplified.  Consider
  $\Cvert=C/(U\cdot C)$ with its differential $\partialvert$.
  By Lemma~\ref{lemma:map-filtered-bases} applied to the
  restriction of $\partialvert$ to a map $\Cvert/\Ker(\partialvert) \to \Image(\partialvert)$, we can find
  a filtered basis $\{\xi_2,\xi_4,\dots,\xi_{2m}\}$ for $\Image(\partialvert)$ and vectors
  $\{\xi_1,\xi_3,\dots,\xi_{2m-1}\}$ in $\Cvert$ so that
  $\partialvert(\xi_{2i-1}) = \xi_{2i}$ and
  $\{\xi_1,\xi_3,\dots,\xi_{2m-1}\}$ descends to a filtered basis of
  $\Cvert/\Ker(\partialvert)$.
  The set $\{\xi_1,\xi_2,\dots,\xi_{2m}\}$ is linearly independent,
  since $\partialvert$ is a differential.
  Extend $\{\xi_2,\xi_4,\dots,\xi_{2m}\}$ to a filtered basis
  $\{\xi_2,\xi_4,\dots,\xi_{2m}, \xi_{2m+1},\dots,\xi_n\}$ for
  $\Ker(\partialvert)$. Then $\{\xi_1,\dots,\xi_n\}$ is
  a filtered basis for $\Cvert$,
  and can be thought of as a basis for $C$ over $\Field[U]$, as
  well.  As a basis for~$C$, $\{\xi_1,\dots,\xi_n\}$ is vertically
  simplified.
  
  A similar argument applies for $\Chor=\gr_0(C\otimes_{\Field[U]}\Field[U,U^{-1}])$
  (where we are now taking the associated graded object with respect
  to the Alexander grading)
  to give a filtered basis
  $\{y_1,\dots,y_n\}$ for $\Chor$ so that $\partialhor(y_{2i-1}) = y_{2i}$ for $1
  \le i \le m$ and $\partialhor(y_i) = 0$ for $i > 2m$.
  Pick (arbitrary) lifts $\hat{y}_i$ of
  $y_i$ to~$C^\infty=C\otimes_{\Field[U]} \Field[U,U^{-1}]$; by
  assumption $A(\hat{y}_i) = 0$, and each $\hat{y}_i$ can be written
  uniquely as $\hat{y}_i = U^{k_i} \eta_i$ 
  for some $\eta_i\in C$ with non-zero projection to $C/U C$.
  Then the $\eta_i$ are the desired horizontally
  simplified basis.
\end{proof}

\begin{proof}[Proof of Theorem~\ref{thm:HFKtoHFD}.]
  Let $C$ be a reduced complex homotopy equivalent to $\CFKm(K)$, and
  let $\{\xi_i\}$ and $\{\eta_i\}$ be bases of $C$ which are
  vertically and horizontally simplified, respectively.  These bases
  are related by
  \begin{equation}
    \label{eq:xi-eta-reln}
    \xi_i = \sum_{j} \Lambda^j_i \eta_j
  \end{equation}
  for a suitable change of basis matrix~$\Lambda$.  
  Let  $\CFDa$ be $\CFDa(S^3\setminus \nbd(K))$ with a suitably
  large negative framing~$-n$, as computed in Theorem~\ref{thm:HFKtoHFD2}.
  In fact, we will think of $\CFDa$ as associated to $C$ (rather than
  $\CFKm(K)$), in view of Proposition~\ref{prop:HotopyEquivalentDStructures}.
  As in
  Theorem~\ref{thm:HFKtoHFD2}, $\{\xi_i\}$ and $\{\eta_i\}$ can also
  be thought of as two different bases for $V^0_s$ over~$\Field$ or,
  alternatively, two bases for $\iota_0\CFDa$ over $\Alg(\Torus)$. 
  We
  will modify and extend these to two bases for all of $\CFDa$ which
  \begin{itemize}
  \item  on $\iota_0\CFDa$ are still related by
    Formula~\eqref{eq:xi-eta-reln},
  \item agree on $\iota_1\CFDa$, and
  \item split as a direct sum of a contractible complex and a complex
    with the coefficient maps from Theorem~\ref{thm:HFKtoHFD}.
  \end{itemize}
  Here, to define the coefficient maps, we view $\CFDa$ as a type $D$
  structure as in Remark~\ref{rmk:UnderlyingTypeD}.

  To give an initial basis for $\iota_1\CFDa$, recall that all $\xi_i$
  except for $\xi_0$ come in pairs $\xi_{j}$ and $\xi_{j+1}$, where
  $\partial_z \xi_{j}=\xi_{j+1}$.  Fix $i$ so that there is a
  height~$\ell$ vertical arrow from $\xi_{i}$ to $\xi_{i+1}$.  Then
  the generators for $V^1_s$ for $s > 0$ coming from this pair have
  coefficient maps like
  \begin{equation}\label{eq:coeff-map-vert}
  \mathcenter{
  \xymatrix@-0.3pc{
    \xi_i\ar[r]^{D_{1}}&\kappa^i_1&\cdots\ar[l]_{D_{23}}&
      \kappa^i_k\ar[l]_{D_{23}}&\kappa^i_{k+1}\ar[l]_{D_{23}}&
      \cdots\ar[l]_{D_{23}}&
      \kappa^i_\ell\ar[l]_{D_{23}}&\alpha^i\ar[l]_{D_{23}}\ar[d]^{D_\emptyset}&
      \cdots\ar[l]_{D_{23}}\ar[d]^{D_\emptyset}\\
    &&&&&&\xi_{i+1}\ar[r]^{D_{1}}&\beta^i&\cdots\ar[l]_{D_{23}}
    }
    }
  \end{equation}
  That is, the only time we get a contribution from this pair to the
  $D_\emptyset$-homology of $V^1_s$ for $s > 0$ is for
  $A(\xi_{i+1})+\left(\frac{n-1}{2}\right)< s\le
  A(\xi_{i})+\left(\frac{n-1}{2}\right)$.  We call these generators
  $\kappa^{i}_k$ as in the statement of the theorem, where
  $k=A(\xi_{i})-s+\left(\frac{n-1}{2}\right)+1$.  For
  $s=A(\xi_{i+1})+\left(\frac{n-1}{2}\right)$, we have a canceling
  pair of generators $\alpha^i$ (from $\xi_i$) and $\beta^i$ (from
  $\xi_{i+1}$).  We also have canceling pairs for smaller~$s$.

  Similarly, for a horizontal arrow from $\eta_i$ to $\eta_{i+1}$ of
  width~$\ell$, we get generators and coefficient maps like
  \begin{equation}\label{eq:coeff-map-horiz}
  \mathcenter{
  \xymatrix@-0.3pc{
    \eta_i\ar[r]^{D_3}&\lambda^i_1\ar[r]^{D_{23}}&
      \cdots\ar[r]^{D_{23}}&
      \lambda^i_k\ar[r]^{D_{23}}&\lambda^i_{k+1}\ar[r]^{D_{23}}&
      \cdots\ar[r]^{D_{23}}&\lambda^i_\ell\ar[r]^{D_{23}}\ar[d]^{D_2}&
      \gamma^i\ar[r]^{D_{23}}\ar[d]^{D_\emptyset}&\cdots\ar[d]^{D_\emptyset}\\
    &&&&&&\eta_{i+1}\ar[r]^{D_3}&\delta^i\ar[r]^{D_{23}}&\cdots
  }
  }
  \end{equation}
  (The $\lambda_k^i$ are elements of $V^1_s$ with $s<0$.)
  In particular, there is an undesirable term in $D_3(\eta_{i+1})$.
     In addition, we usually have $D_{123}(\eta_i) = 0$; the
     exception is when the arrow has width equal to one, when
     \begin{equation}
       \label{eq:NontrivialD123}
       D_{123}(\eta_i)=D_{1}(\eta_{i+1}).
     \end{equation}
     
     Consider next $\xi_0$, which has no vertical differential either
     into or out of it. We see that its corresponding element comes
     from a chain

     $$\dots 
     \stackrel{D_{23}}{\longrightarrow}\varphi_3
     \stackrel{D_{23}}{\longrightarrow}\varphi_2
     \stackrel{D_{23}}{\longrightarrow}\varphi_1
     \stackrel{D_1}{\longleftarrow} \xi_0,
     $$
     where the chain of elements $\varphi_i$ extends back to
     $V^1_s$ where $s$ is  small in magnitude
     (and hence $V^1_s$ is
     one-dimensional).  Similarly, consider the element
     $\eta_0$ which has no horizontal differential either into or
     out of it. (Remember that $\xi_0$ and $\eta_0$ can coincide.)
     Associated to $\eta_0$ is a similar string
     $$\eta_0 
     \stackrel{D_3}{\longrightarrow}
     {\mu_1}
     \stackrel{D_{23}}{\longrightarrow}
     {\mu_2}
     \stackrel{D_{23}}{\longrightarrow}
     {\mu_3}
     \stackrel{D_{23}}{\longrightarrow}
     \cdots
     $$
     where the chain of elements $\mu_i$ extends forward
     to $V^1_s$ where $s$ is sufficiently small in magnitude 
     that $V^1_s$
     is one-dimensional).  Thus, for some sufficiently large $a$ and $b$,
     we must have that $\varphi_a=\mu_b$; i.e., we can connect the two
     chains to a single one of the form
     $${\eta_0}
     \stackrel{D_3}{\longrightarrow}
     {\mu_1}
     \stackrel{D_{23}}{\longrightarrow}
     {\mu_2}
     \stackrel{D_{23}}{\longrightarrow}
     {\mu_3}
     \stackrel{D_{23}}{\longrightarrow}
     {\cdots}
     \stackrel{D_{23}}{\longrightarrow}
     {\mu_m}
     \stackrel{D_{1}}{\longleftarrow}
     {\xi_0}
     $$

     We have now constructed bases for all the $V^1_s$:
     \begin{itemize}
     \item $\lambda^i_k$ and $\mu_k$ for $s \le -n/4$,
     \item $\phi_k = \mu_k$ for $-n/4 < s < n/4$, and
     \item $\kappa^i_k$ and $\phi_k$ for $s \ge n/4$.
     \end{itemize}
     However, we have extra generators for the modules and extra terms in
     the coefficient maps.
     To bring the module to the desired form, we must change basis
     so as to
     \begin{itemize}
     \item split off the chain of canceling generators
       in~\eqref{eq:coeff-map-horiz} by arranging for
       $D_3(\eta_{i+1})$ and $D_{23}(\lambda^i_\ell)$ to vanish;
     \item eliminate the generators with canceling differentials
       in~\eqref{eq:coeff-map-vert}, in particular $\alpha^j$ and~$\beta^j$; and
     \item bring $D_{123}$ to the desired form, by eliminating
       the non-trivial terms belonging to width one arrows
       (Equation~\eqref{eq:NontrivialD123}).
     \end{itemize}
     We do this by successively changing basis (and, at some point,
     dropping out $\beta^j$).

     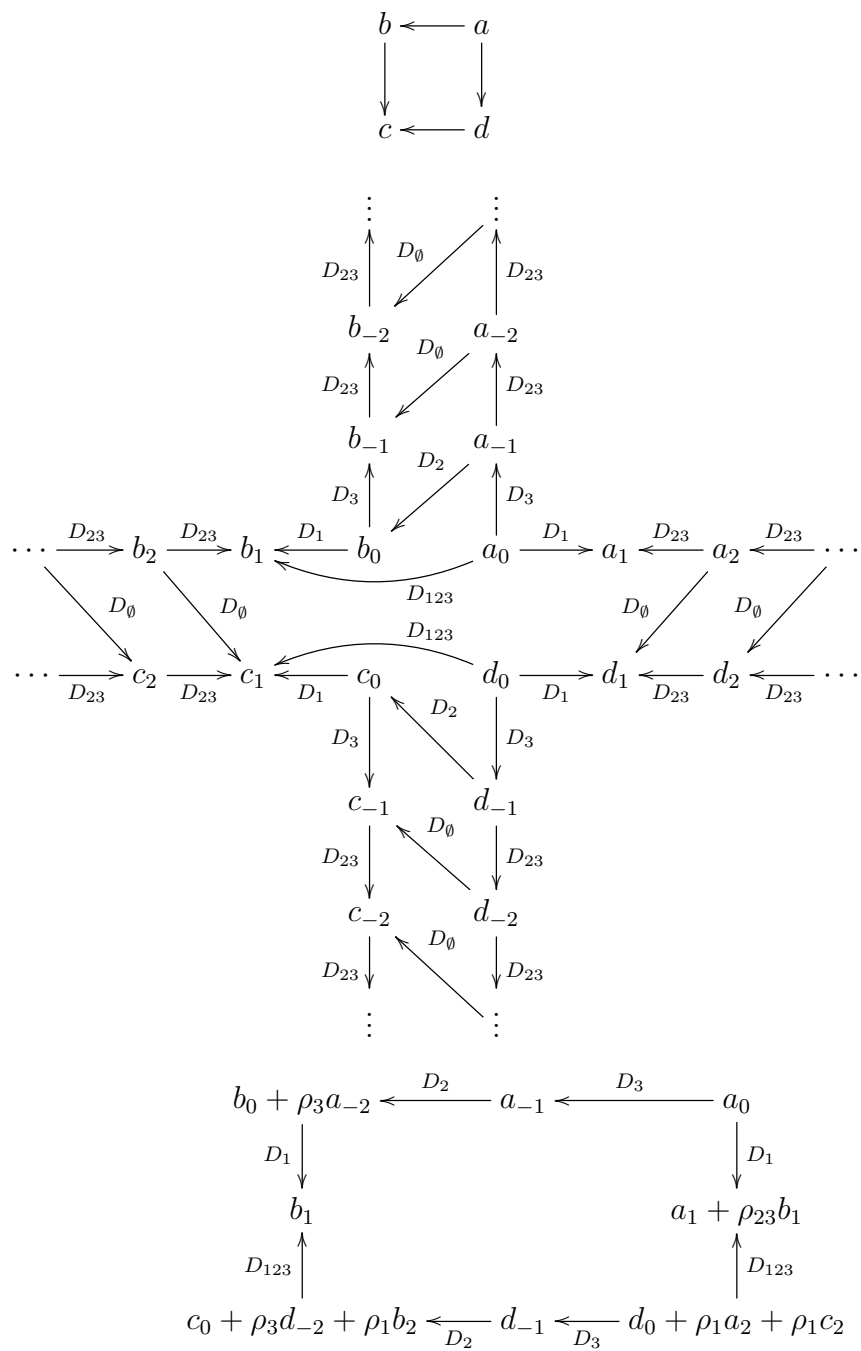
\begin{figure}
       \centering
       \[
       \xymatrix{
         b\ar[d] & a\ar[d]\ar[l]\\
         c & d\ar[l] }
       \]
       \smallskip
       \[
       \xymatrix{
         &&&\vdots&\vdots\ar[dl]_{D_{\emptyset}}\\
         &&&b_{-2}\ar[u]^{D_{23}}&a_{-2}\ar[dl]_(0.35){D_{\emptyset}}\ar[u]_{D_{23}}\\
         &&&b_{-1}\ar[u]^{D_{23}}&a_{-1}\ar[dl]_(0.35){D_2}\ar[u]_{D_{23}}\\
         \cdots\ar[r]^{D_{23}}\ar[dr]^(0.6){D_\emptyset}&
           b_2\ar[r]^{D_{23}}\ar[dr]^(0.6){D_\emptyset}&
            b_1&
             b_0\ar[l]_{D_1}\ar[u]^{D_3}&
              a_0\ar[r]^{D_1}\ar[u]_{D_3}\ar@/^1pc/[ll]^(0.3){D_{123}}&
               a_1&
                a_2\ar[l]_{D_{23}}\ar[dl]_(0.6){D_\emptyset}&
                 \cdots\ar[l]_{D_{23}}\ar[dl]_(0.6){D_\emptyset}\\
         \rule{0pt}{21pt}\cdots\ar[r]_{D_{23}}&
          c_2\ar[r]_{D_{23}}&
           c_1&
            c_0\ar[d]_{D_3}\ar[l]^{D_1}&
             d_0\ar[d]^{D_3}\ar[r]_{D_1}\ar@/_1pc/[ll]_(0.3){D_{123}}&
              d_1&
               d_2\ar[l]^{D_{23}}&
                \cdots\ar[l]^{D_{23}}\\
         &&&c_{-1}\ar[d]_{D_{23}}&d_{-1}\ar[d]^{D_{23}}\ar[ul]_(0.6){D_2}\\
         &&&c_{-2}\ar[d]_{D_{23}}&d_{-2}\ar[d]^{D_{23}}\ar[ul]_(0.6){D_\emptyset}\\
         &&&\vdots&\vdots\ar[ul]_(0.6){D_\emptyset}
       }
       \]
       \smallskip
       \[
       \xymatrix{ & & b_0+\rho_3a_{-2} \ar[d]_{D_1}&
         a_{-1}\ar[l]_<>(0.5){D_2} &
         a_0\ar[l]_{D_3}\ar[d]^{D_1}\\
         & & b_1 & & a_1+\rho_{23}b_1 \\
         & & c_0+\rho_3d_{-2}+\rho_1b_2\ar[u]^{D_{123}} &
         d_{-1}\ar[l]^<>(.5){D_2} &
         d_0+\rho_1a_2+\rho_1c_2\ar[l]^<>(.5){D_3}\ar[u]_{D_{123}} }
       \]
       \caption[Illustration of the proof of
       Theorem~\ref{thm:HFKtoHFD}]{\textbf{An illustration of the
           proof of Theorem \ref{thm:HFKtoHFD}.} The top illustrates a summand in the knot complex; the middle 
       illustrates a corresponding piece coming from
       Theorem~\ref{thm:HFKtoHFD2}; the bottom is obtained by making
       the substitutions coming from the proof of
       Theorem~\ref{thm:HFKtoHFD} and dropping acyclic summands.}
       \label{fig:torus-cancel}
     \end{figure}
          
     We start by modifying the base to eliminate the unwanted terms
     in~$D_3$ (from $\eta_{i+1}$ when there is a horizontal arrow from
     $\eta_i$). To this end, we change basis by
     $$\eta_{i+1}'\coloneqq\eta_{i+1}+\rho_3\cdot \gamma^i.$$
     With this substitution, we have that
     $D_3(\eta_{i+1}')=0$, $D_2(\lambda^i_\ell)=\eta_{i+1}'$, and
     $D_{23}(\lambda^i_\ell)=\penalty1000 0$.   Let $\{\xi_j'\}$ denote the
     corresponding change of basis applied to the $\{\xi_j\}$ so that
     Formula~\eqref{eq:xi-eta-reln} remains true.
     Note that the other coefficient maps are unchanged;
     in particular $D_1(\xi_j')=D_1(\xi_j)$, and there are no
     coefficient maps other than $D_2$ which enter $\eta_{i+1}'$.
     Also, the chain of canceling generators from this horizontal
     arrow now splits off as a direct sum, so we can drop them.

     We pass now to a submodule to eliminate $\alpha^i$ and $\beta^i$.
     To do this, for each vertical arrow from $\xi_i$ to $\xi_{i+1}$
     of height~$\ell$
     (and corresponding new elements $\xi_i'$ and $\xi_{i+1}'$),
     we replace $\xi_{i+1}'$ by
     \[\xi_{i+1}''\coloneqq\xi_{i+1}' + \rho_1 \cdot\alpha^i.\]
     Now we have $D_{123}(\xi_{i+1}'')=\kappa^{i}_\ell+D_{123}(\xi_{i+1})$. 
     Again, there is a corresponding modification to the other
     generators, replacing $\eta_i'$ by $\eta_i''$. Now the submodule
     without the chain of canceling generators from~\eqref{eq:coeff-map-vert}
     is easily seen to be homotopy equivalent to the original module.

     Finally, we wish to bring $D_{123}$ to the desired form,
     eliminating the terms from horizontal arrows of
     length one. This is done as follows. Suppose that $\xi_i$
     is the initial point of some vertical arrow. Replace
     $\kappa^{i}_1$ by
     $${\widetilde\kappa}^{i}_1\coloneqq\kappa^{i}_1+
     \rho_{23}\cdot D_{123}(\xi_i).$$
     With this new basis, then, $D_{123}(\xi_i'')=0$ (as its
     contribution is absorbed into~$D_1$). Moreover, 
     this change does not affect any of the 
     other coefficient maps.
     
     It remains to make one last change of basis corresponding
     to terminal points~$\xi_{i+1}$ of vertical arrows so that 
     $D_{123}(\xi_{i+1})$ is the element taking the place of
     $\kappa^{i}_\ell$.
     \emph{A priori} we have
     $D_{123}(\xi_{i+1}'')=\kappa^{i}_\ell+D_{123}(\xi_{i+1})$.
     In the case where the length of the arrow from $\xi_i$ to $\xi_{i+1}$
     is greater than one, we make the change by letting
     $${\widetilde \kappa}^{i}_\ell\coloneqq\kappa^{i}_\ell+D_{123}(\xi_{i+1}).$$
     Since $D_{123}(\xi_{i+1})$ is in the kernel of $D_{23}$
     (from the description in Theorem~\ref{thm:HFKtoHFD2}),
     all other coefficient maps remain unchanged. 

     In the case where the length of the arrow from $\xi_i$ to $\xi_{i+1}$
     is one, proceed as follows.
     First, we find an element $\zeta_{i}$ so that
     \[
     \bdy\zeta_i = \rho_{23}\cdot D_{123}(\xi_{i+1}).
     \]
     In particular, $D_\emptyset(\zeta_{i})=0$ and
     $D_{23}(\zeta_i)=D_{123}(\xi_{i+1})$.
     This element
     is found as follows. Recall that $\xi_i$ and $\xi_{i+1}$
     correspond to elements $\x^i$ and $\x^{i+1}$ in $C(s)$ and
     $C(s-1)$, respectively. Let $\y\coloneqq \partial^1(\x^i)\in
     C(s+1)$, i.e., the image of $\x^i$ under the map which
     counts holomorphic disks which cross $w$ exactly once.  (Recall
     from Theorem~\ref{thm:HFKtoHFD2}
     that $D_{123}(\xi_{i+1})$ is $\partial^1(\x^{i+1})$ considered as
     an element of $C(\geqo s)$.)
     We then let $\zeta_i=\partial_z(\y)$, where we consider $\y$
     as a chain in $C(\geq s-1)$.
     Then $\zeta_i$,
     thought of as a generator of the type $D$ module, has the desired
     properties.
     Making the substitution
     \[\xi_{i+1}'''\coloneqq\xi_{i+1}''+\rho_1\cdot \zeta_i,\]
     we see that 
     $D_1(\xi_{i+1}''')=D_1(\xi_{i+1}'')$ and
     $D_{123}(\xi_{i+1}''')$ has the desired form.  
 
     This module now has the
     form promised in Theorem~\ref{thm:HFKtoHFD}; except that we must relate 
     the length $m$ of the chain of
     $\mu_i$ to the framing parameter $n$ and the knot invariant
     $\tau(K)$, and we must compute the gradings.
     We do both of these, as follows.
%
     On the one hand, the stable chain of length 
     $m$ ensures that
     \[ \lambda^{-m+1}\gr(\eta_0)=\gr(\rho_{3})\cdot \gr(\rho_{23})^{m-1} \cdot \gr(\rho_{1})^{-1}\cdot \gr (\xi_0).\]
     On the other hand, by Equation~\eqref{eq:GradingsOnXiEta},
     Proposition~\ref{prop:GradingHFKtoHFD}
     guarantees that there is some element of the grading set ${\mathbf s}$ so that
     \begin{align*}
       \gr(\xi_0)&=\lambda^{-2\tau(K)}\gr(\rho_{23})^{-\tau(K)}{\mathbf s} \\
       \gr(\eta_0)&=\gr(\rho_{23})^{\tau(K)}{\mathbf s}.
     \end{align*}
     Substituting these into the relation from the stable chain (and making a straightforward computation in $G$), we see that
     \begin{align*}
       {\mathbf s} &= \lambda^{m-1-2\tau(K)}
       \gr(\rho_{23})^{-\tau(K)}\gr(\rho_{3})\cdot \gr(\rho_{23})^{m-1}\cdot \gr(\rho_{1})^{-1}\gr(\rho_{23})^{-\tau(K)} \cdot
       {\mathbf s} \\
       &= \lambda^{-1}\cdot \gr(\rho_{23})^{m-2\tau(K)} \gr(\rho_{12})^{-1}\cdot {\mathbf s}.
     \end{align*}
     This gives a periodic domain with local multiplicities $(0,-1,m-2\tau-1,m-2\tau)$ at the four regions
     $(0,1,2,3)$ at the boundary. We conclude that $n=m-2\tau$ (see
     the discussion leading to Formula~\eqref{eq:GradingMap}), and
     \[ G \cdot {\mathbf s} \cong G/\lambda^{-1} \gr(\rho_{23})^n \gr(\rho_{12})^{-1}. \]
     The rest of the grading statement follows from Proposition~\ref{prop:GradingHFKtoHFD}.
\end{proof}

In Figure~\ref{fig:torus-cancel} we have illustrated some of the changes
of basis from
Theorem~\ref{thm:HFKtoHFD} in the case where the knot complex has a summand
which is a
square, as illustrated on the top of the figure.
(Such a summand appears, e.g., in the knot Floer complex for the
figure eight knot,
on the bottom of Figure~\ref{fig:HFKtoHFD}.)
Initially, Theorem~\ref{thm:HFKtoHFD2} gives a corresponding summand
as illustrated in
the middle of Figure~\ref{fig:torus-cancel}. Going through the cancellations
prescribed in the proof of Theorem~\ref{thm:HFKtoHFD}, we end up with
the complex on the bottom of Figure~\ref{fig:torus-cancel}.

\section{Satellites revisited}
\label{sec:CablesAgain}
In this section, we illustrate 
Theorems~\ref{thm:PairingKnot} and~\ref{thm:HFKtoHFD} by
showing how to compute $\HFKm$ of a knot $K$ obtained as 
a
cable of the left-handed trefoil $T$ in $S^3$. (Specifically, we compute
$\HFKm$ of the $(2,-3)$-cable of $T$.) 
Although the knot Floer homology groups of this
particular knot were already known \cite{HeddenThesis}, the techniques
here extend to give a formula for $\HFKm$
of any particular satellite in terms of $\CFKm$ of the pattern.

\begin{figure}
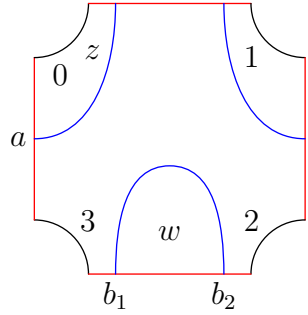

$\mfigb{torus-50}$
\caption[Diagram for the $(2,1)$ cable]{\label{fig:Cable}
  {\bf{Heegaard diagram for computing the type $A$ module for the $(2,1)$ cable.}}
  This is a doubly-pointed Heegaard diagram for the $(2,1)$ cable (of the unknot),
  thought of as a knot in the solid torus. 
  The basepoint $z$ lies in the region marked with a $0$.}\label{fig:21cable}
\end{figure}

Figure~\ref{fig:21cable} is a Heegaard diagram for the $(2,1)$ cabling
operation. Since this diagram has genus equal to one, the  holomorphic
curves are straightforward to count. 
The type $A$ module $\CFAm(C)$ associated
to this diagram is given as follows. It has three generators $a$, $b_1$ and
$b_2$, and the following relations:
\begin{equation}
\begin{aligned}
m_1(b_1) &= U \cdot b_2 \\
m_2(a, \rho_1) &= b_2 \\
m_{3+i}(a,\rho_3,\overbrace{\rho_{23},\dots,\rho_{23}}^{i},\rho_2)
&=U^{2i+2}\cdot a, & i&\ge 0\\
m_{4+i}(a, \rho_3, \overbrace{\rho_{23},\dots,\rho_{23}}^{i},
\rho_2,\rho_1) &= U^{2i+1} \cdot b_1, & i&\ge 0.
\end{aligned}\label{eq:CFAmC}
\end{equation}
Letting $u=gr(U)=(0;0,0;-1)$ in $\eG=\smallGroup(\PMC)\times \ZZ$, we find that 
the grading set is isomorphic to 
\[ 
\gr(a)\cdot \eG \cong  u^{-2} \lambda\gr(\rho_3)\gr(\rho_2)\backslash \eG=
u^{-2}\gr(\rho_{23})\backslash \eG;
\]
and
\[ \gr(a)=e, \qquad \gr(b_1)=\lambda u \gr(\rho_1) \qquad \gr(b_2)= \gr(\rho_1);\]
or, in coordinates, the grading set is
$(\OneHalf;0,1;-2)\backslash \eG$ and
\[ \gr(a)=(0;0,0;0), \qquad \gr(b_1)=(\OneHalf;\OneHalf,-\OneHalf;-1)
\qquad \gr(b_2)=(-\OneHalf;\OneHalf,-\OneHalf;0).\]
Here, the first three components parameterize the usual grading group
($\smallGroup$) and the last component represents the additional $\ZZ$ enhancement.

\begin{figure}
\centerline{\input{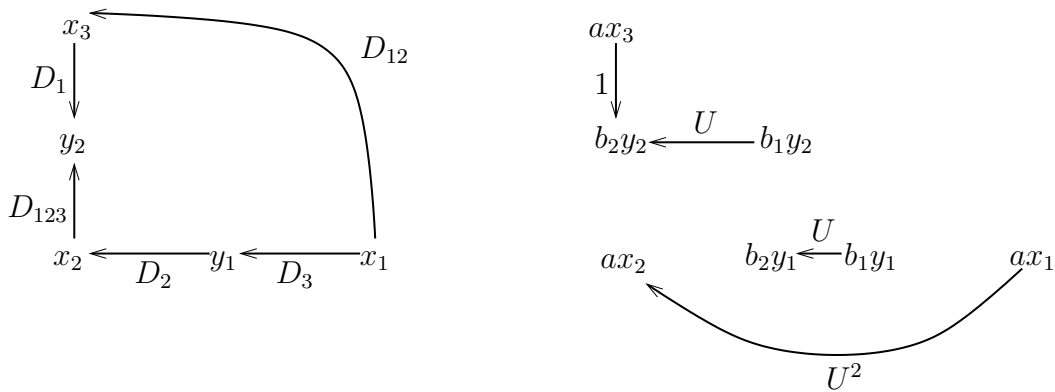}}
\caption[Cable of the trefoil]{\label{fig:CableTrefoil}
  {\bf{Cable of the trefoil.}}  On the left, the coefficient maps for
  the left-handed trefoil.  On the right, the result of cabling.
  Arrows here represent differentials; they are labeled
  by their corresponding coefficients (powers of $U$).
  }
\end{figure}

Consider next the type $D$ module for the left-handed trefoil knot
$T$ with framing $-2$, represented graphically in
Figure~\ref{fig:CableTrefoil}. This module has generators
$x_1,x_2,x_3$ in the $\iota_0$ idempotent and $y_1,y_2$ in the $\iota_1$
idempotent, and coefficient maps:
\begin{align*}
  D_3(x_1)&=y_1&
  D_{12}(x_1)&=x_3&
  D_2(y_1)&=x_2\\
  D_{123}(x_2)&=y_2&
  D_1(x_3)&=y_2.
\end{align*}
(Although this framing is not large enough for
Theorem~\ref{thm:HFKtoHFD} to apply, $\CFDa(T)$ may still be computed;
see~ Theorem~\ref{thm:HFKtoHFDframed}, whose proof uses results from
\cite{LOT2}.)

Gradings of the elements $x_i$ are computed in
Equation~\eqref{eq:gr-on-CFK-to-CFD}, with the exception of $y_1$ and
$y_2$. But the gradings of these are related to the gradings of the
$x_i$ according to the coefficient maps.  Explicitly, the generators
$x_1$, $x_2$, and $x_3$ correspond to knot Floer generators
with $(M,A)$-bigradings given by
$(0,-1)$, $(1,0)$, and $(2,1)$ respectively.
Since $\gr(x_2)=\lambda$ and $\partial y_1=\rho_2\cdot x_2$,
it follows that
\[ \gr(y_1)=\lambda^2 \gr(\rho_2).\]
Indeed, in coordinates:
\begin{align*}
  \gr(x_1)&=({\textstyle\frac{3}{2}};0,1) & \gr(x_2)&= (1;0,0) &
  \gr(x_3)&=(\OneHalf;0,-1) \\
  \gr(y_1)&=({\textstyle\frac{3}{2}};\OneHalf,\OneHalf) &
  \gr(y_2)&=(\OneHalf;-\OneHalf,-\OneHalf).&&
\end{align*}

Let $K$ denote the $(2,-3)$-cable of the (0-framed) left-handed
trefoil (which is also the $(2,1)$-cable of the $-2$-framed
left-handed trefoil).
By
Theorem~\ref{thm:PairingKnot}, $\gCFKm(K)$ is computed by the tensor
product $\CFAm(C)\DT\CFDa(T)$, which has generators $a\otimes x_i$ and
$b_k\otimes y_j$, which we abbreviate by dropping the tensor product
sign. The grading set is identified with the double-coset space
\[ u^{-2}\gr(\rho_{23})\backslash \eG / \lambda^{-1} \gr(\rho_{23})^2 \gr(\rho_{12})^{-1}.\]
Every element of the double-coset space is equivalent to a uniquely specified
element of the form $\lambda^a u^b$ for $a, b\in \ZZ$:
use right translation to eliminate $\gr(\rho_{12})$ from the $\SpinC$ component and
then left translation to eliminate $\gr(\rho_{23})$. For example, 
\begin{align*}
\gr(b_{1}y_{1})&=\gr(b_1)\gr(y_1)=\lambda^3 u \gr(\rho_{12})\\
&\sim (u^{4}\gr(\rho_{23})^{-2})
(\lambda^3 u \gr(\rho_{12}))
(\lambda^{-1}\gr(\rho_{23})^2 \gr(\rho_{12})^{-1})
=\lambda^6 u^{-5}.
\end{align*}
Proceeding in this manner, we find:
\begin{align}
  \gr(ax_1)&=(2;0,0;-2)&\gr(ax_2)&=(1;0,0;0)&\gr(ax_3)&=(0;0,0;2) \nonumber \\
  \gr(b_1 y_1)&=(6;0,0;-5)& \gr(b_2 y_1) &= (5;0,0;-4)& \gr(b_1 y_2) &= (0;0,0;1) \label{eq:GradedTensorProductCable} \\
  \gr(b_2 y_2) &= (-1;0,0;2) \nonumber
\end{align}

The first component here is the $z$-normalized Maslov grading $N$ of
Equation~\eqref{eq:RenormalizedMaslov}, {\em a priori} up to an
overall translation. Similarly, the last component is the Alexander
grading, again up to an overall translation. In fact, we graded both
the knot complement and the cabling piece so that the $z$-normalized
Maslov grading $N$ is computed absolutely; see Remark~\ref{rmk:Normalization}.

The differentials are readily computed; they are:
\begin{align*}
  \bdy(a\otimes x_1)&=U^2\cdot a\otimes x_2&
  \bdy(a\otimes x_2)&=0&
  \bdy(a\otimes x_3)&=b_2\otimes y_2\\
  \bdy(b_1\otimes y_1)&=U\cdot b_2\otimes y_1&
  \bdy(b_2\otimes y_1)&=0&
  \bdy(b_1\otimes y_2)&=U\cdot b_2\otimes y_2\\
  \bdy(b_2\otimes y_2)&=0.
\end{align*}

To pin down the indeterminacy in the Alexander grading, we proceed as follows.
Observe that the Poincar{\'e} polynomial for the
tensor product, after setting $U=0$,
is given by
$$
\sum_{d,r} q_N^d \cdot t^r \cdot \rank H_{\substack{A=r \\N=d}}(\CFAm(C)\DT\CFDa(T))=
t^{1} q_N^0 + t^{0} q_N + t^{-2}q_N^2 + t^{-4} q_N^5 + t^{-5} q_N^6.
$$
Here, the formal variable $t$ records the Alexander factor of the
grading, while the formal variable $q_N$ records the $z$-normalized
grading $N$ of Equation~\eqref{eq:RenormalizedMaslov}, which is the
first coordinate of Equation~\eqref{eq:GradedTensorProductCable}, up
to an additive constant. Further specializing to $q_N=(-1)$ gives a polynomial in $t$
which is not symmetric under exchanging $t$ and $t^{-1}$, but 
$t^{2}$ times it is. So, by Equation~\eqref{eq:AlexanderSymmetry}, the
absolute Alexander gradings are $2$ less than the last components in
Equation~\eqref{eq:GradedTensorProductCable}. With respect to this
translation, the Poincar\'e polynomial is
$$
\sum_{d,r} q_N^d \cdot t^r \cdot \rank H_{\substack{ A=r \\ N=d }}(\CFAm(C)\DT\CFDa(T))=
t^{-3} q_N^6 + t^{-2} q_N^5 + q_N^2 + t^2 q_N + t^3.
$$

Since $N=M-2A$, we can recover the Poincar{\'e} polynomial of $K$ in terms of the usual 
Alexander and Maslov gradings:
$$
\sum_{d,r} q^d \cdot t^r \cdot \rank H_{d}(\CFAm(C)\DT\CFDa(T),r)=
t^{-3} + t^{-2} q + q^2 + t^2 q^5 + t^3 q^6.
$$
Here, the formal variable $t$ records Alexander gradings as before,
while $q$ records the usual Maslov gradings. This answer is in
agreement with~\cite[Theorem 1.0.6]{HeddenThesis}.


\appendix
\chapter{Bimodules and change of framing}
\label{app:Bimodules}

The modules $\CFAa(Y)$ and $\CFDa(Y)$ depend (up to homotopy
equivalence) not only 
on the $3$-manifold $Y$ but also on the parametrization of $\bdy Y$ by
a reference surface $F$; this dependence can already be seen for the
modules associated to solid tori computed in
Section~\ref{sec:surg-exact-triangle}.  The result of
reparametrization (e.g., change of framing on a knot complement)
is captured by tensoring with certain
bimodules. Furthermore, the modules
$\CFAa(Y)$ and $\CFDa(Y)$ are related to each other by certain dualizing
bimodules. These topics are explored in detail in~\cite{LOT2}. 
In this chapter we present highlights from
that paper, and then exhibit the bimodules relevant to the case of
torus boundary. When combined with Theorem~\ref{thm:HFKtoHFD}, this
gives a complete computation of the bordered Floer invariants of knot
complements in $S^3$, with any framing.

We start by briefly stating the main results on bimodules in
Section~\ref{sec:bimod-results}. We then sketch the definitions of the
bimodules in Section~\ref{sec:bimod-construct}. We give explicitly the
bimodules for a set of generators of the mapping class group of the
torus in Section~\ref{sec:comp-torus-bimods} and use these bimodules
in Section~\ref{sec:hfk-to-CFD-general} to finish the computation of
$\CFDa$ of knot complements in terms of the knot Floer homology.

\section{Statement of results}\label{sec:bimod-results}
Fix a pointed matched circle $\PMC$.
Recall from Definition~\ref{def:bordered-3-mfld} that a bordered
  $3$-manifold with boundary $F(\PMC)$ is a triple 
$\gls*{borderedmfld}$
where $Y$ is a compact, oriented $3$-manifold with boundary, $\PMC$ is
a pointed matched circle and $\phi\co F(\PMC)\to\bdy Y$ is an
orientation-preserving diffeomorphism.  We have associated to
$(Y,\phi)$ homotopy equivalence classes of modules $\CFDa(Y,-\phi)$
(Chapter~\ref{chap:type-d-mod}) and $\CFAa(Y,\phi)$
(Chapter~\ref{chap:type-a-mod}) over $\Alg(-\PMC)$ and $\Alg(\PMC)$,
respectively.
Here, $-\phi$ is the map $\phi$ but viewed as a map $-F(\PMC)\to -\bdy
Y$; compare Theorem~\ref{intro:D-invariance}.

\begin{theorem}\label{thm:bimodule-pairing}
\index{pairing theorem!for change of framing}%
  Given a diffeomorphism $\psi\co F(\PMC)\to
  F(\PMC)$, there is an $\Ainf$
  bimodule 
  $\gls*{CFDApsi}$,
  well-defined up to ($\Ainf$) homotopy equivalence, such that for any 
  bordered three-manifold $(Y,\phi\co F(\PMC)\to\partial Y)$,
  \begin{align*}
    \CFAa(Y,\phi\circ\psi^{-1})&\simeq \CFAa(Y,\phi)\DTP_{\Alg(\PMC)}\CFDAa(\psi) \\
    \CFDa(-Y,\phi\circ\psi)&\simeq
    \CFDAa(\psi)\DTP_{\Alg(\PMC)}\CFDa(-Y,\phi).
  \end{align*}
  Moreover, given another map $\xi\co F(\PMC)\to F(\PMC)$,
  \[
  \CFDAa(\psi\circ\xi)\simeq \CFDAa(\xi)\DTP_{\Alg(\PMC)}\CFDAa(\psi).
  \]
\end{theorem}
As we will discuss in Section~\ref{sec:bimod-construct},
$\CFDAa(\psi)$ is defined in terms of a Heegaard
diagram with two boundary components, one of which is treated in type
$D$ fashion and the other of which is treated in type $A$ fashion.

Next we turn to the duality between $\CFDa$ and $\CFAa$. By
$\gls*{ModRS}$
we mean a bimodule $M$ with commuting right actions by rings
$R$ and $S$; similarly, by
$\gls*{RSMod}$
we mean a bimodule $M$ with
commuting left actions by $R$ and $S$.
\begin{theorem}\label{thm:dualizing-bimods}
\index{pairing theorem!relating $\CFAa$ and $\CFDa$}%
\index{duality theorem}%
  For any pointed matched circle~$\PtdMatchCirc$ there
  are differential graded bimodules
  $\gls*{CFAAId}$
  and
  $\gls*{CFDDId}$
  such that for any bordered $(Y,\phi)$,
  \begin{align*}
    \CFDa(Y,-\phi)&\simeq
    \CFAa(Y,\phi)\DTP_{\Alg(\PtdMatchCirc)}\CFDDa(\Id)\\
    \CFAa(Y,\phi)&\simeq\CFAAa(\Id)\DTP_{\Alg(\PtdMatchCirc)}\CFDa(Y,-\phi).
  \end{align*}
\end{theorem}
As the notation suggests, these dualizing bimodules correspond to the
identity map $\Id\co F\to F$. In fact, for any mapping class $\phi\co
F\to F$ there are associated bimodules $\CFDDa(\phi)$ and
$\CFAAa(\phi)$; see~\cite{LOT2}.

\section{Sketch of the construction}\label{sec:bimod-construct}
In this section we outline the construction of the bimodules discussed
above. Details can be found in~\cite{LOT2}, although the constructions
are close enough to the construction of $\CFDa$ and $\CFAa$ that the
reader might find it amusing to fill them in on his or her own. We
start by discussing bordered Heegaard diagrams representing (mapping
cylinders of) diffeomorphisms (Section~\ref{sec:hd-for-diffeo}). We
then give the definitions of the bimodules $\CFDAa(\phi)$
(Section~\ref{sec:CFDA-def}) and $\CFDDa(\Id)$ and $\CFAAa(\Id)$
(Section~\ref{sec:def-CFAA-CFDD-id}).

\subsection{Bordered Heegaard diagrams for diffeomorphisms}\label{sec:hd-for-diffeo}
The bimodule for a diffeomorphism $\phi$, $\CFDAa(\phi)$, is associated to a bordered Heegaard diagram
for $\phi$. Such a diagram is a slight generalization of the notion
from Chapter~\ref{chap:heegaard-diagrams-boundary} to allow more than
one boundary component. (Such diagrams already made a brief
appearance in Section~\ref{sec:CFD-to-HFK}.)

Specifically, an \emph{arced bordered Heegaard diagram with two boundary
components}\index{arced bordered Heegaard diagram}\index{Heegaard diagram!bordered!arced}
is a quadruple
$\gls*{abHD}$
where
\begin{itemize}
\item $\overline{\Sigma}$ is an oriented surface of genus $g$, with
  two boundary components $\bdy_L\overline{\Sigma}$ and $\bdy_R\overline{\Sigma}$;\glsadd{bdyLR}
\item
  $\overline\alphas=\{\overline{\alpha}_1^{a,L}\!\!,\dots,\overline{\alpha}_{2k}^{a,L}\!\!,\
  \overline{\alpha}_1^{a,R}\!\!,\dots,\overline{\alpha}_{2k}^{a,R}\!\!,\
  \alpha_1^c,\dots,\alpha_{g-2k}^c\}$
  where the $\overline{\alpha}_i^{a,L}$ are embedded arcs with
  endpoints on $\bdy_L\overline{\Sigma}$, the
  $\overline{\alpha}_i^{a,R}$ are embedded arcs with endpoints on
  $\bdy_R\overline{\Sigma}$, and the $\alpha_i^c$ are embedded circles
  in $\Sigma$; and all of the $\overline{\alpha}_i^{a,L}$\!\!,
  $\overline{\alpha}_i^{a,R}$ and $\alpha_i^c$ are pairwise disjoint,
  and $\Sigma\setminus\alphas$ is
  connected;
\item $\betas=\{\beta_1,\dots,\beta_g\}$ are circles in $\Sigma$ such that
  $\Sigma\setminus\betas$ is connected; and
\item $\arcz$ is an arc in $\overline{\Sigma}\setminus(\overline{\alphas}\cup\betas)$
  connecting $\bdy_L\overline{\Sigma}$ and $\bdy_R\overline{\Sigma}$.
\end{itemize}
As for ordinary bordered Heegaard diagrams, we let
$\Sigma=\overline{\Sigma}\setminus\bdy\overline{\Sigma}$; we will be
sloppy about the distinction between $\Sigma$ and $\overline{\Sigma}$.

The data $\HD$ specifies a $3$-manifold $Y$ with two boundary
components in a manner exactly analogous to
Section~\ref{sec:BorderedDiagrams}. Moreover, setting $\PMC_i$\glsadd{PMCLR}
to be
the pointed matched circle specified by $\bdy_i\Sigma$
($i\in\{L,R\}$), the diagram $\HD$
specifies a homeomorphism $\phi_i$\glsadd{phiLR}
from
$F_i=F(\PMC_i)$\glsadd{FLR}
to $\bdy_i Y$,
as well as a framed arc connecting
$\bdy_LY$ and $\bdy_RY$. We call $\HD$ a
bordered Heegaard diagram for $(Y,\phi_L,\phi_R)$.

Assume that $Y$ is homeomorphic to $[0,1]\times F_R$.
There
is a homeomorphism $\Phi\co [0,1]\times F_R \stackrel{\cong}{\to} Y$
extending $\phi_R\co F_R \to Y$.
Moreover, this homeomorphism is unique up to
homeomorphisms $\Xi\co [0,1]\times F_R\to [0,1]\times F_R$ such that
$\Xi|_{\{1\}\times F_R}$ is the identity. It follows that
$\Phi|_{\{0\}\times F_R}$ is well-defined up to homotopy, and hence
(because $F_R$ is a surface) up to isotopy. (Note that
$\Phi|_{\{0\}\times F_R}$ is orientation-reversing.)  Consequently, 
$(\Phi|_{\{0\}\times F_R})^{-1}\circ \phi_L$ is an orientation-reversing homeomorphism
$F_L\to F_R$, well-defined up to isotopy.
If we assume that $F_L=-F_R$ then $(\Phi|_{\{0\}\times F_R})^{-1}\circ
\phi_L$ can be regarded as an element 
$\Upsilon(\HD)$
of the mapping class group
$\MCG(F_R)$.  Paying attention to the arc $\arcz$, we obtain a map
\[
\Upsilon\co\{\text{bordered Heegaard diagrams for
  $([0,1]\times F,-F,F)$}\}\to \MCG_0(F),
\]
where
$\gls*{MCGzero}$
denotes the strongly based mapping class group
\index{mapping class group}%
\index{strongly based mapping class group}%
of any surface $F=F(\PMC)$; here, ``strongly based'' refers to 
the mapping class group of diffeomorphisms fixing a small disk in
$F$.

\begin{figure}
  \centering
  \includegraphics[scale=.83333]{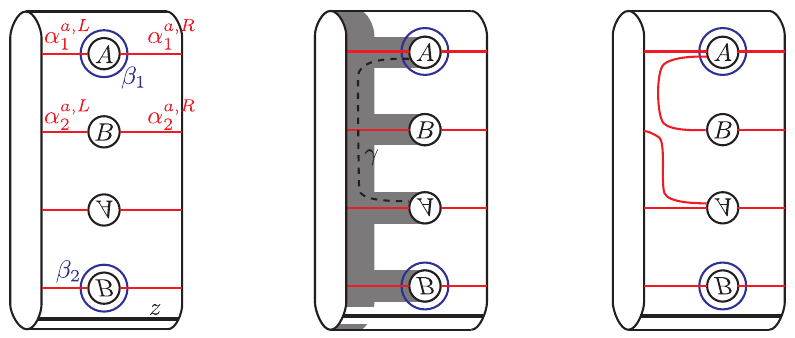}
  \caption[Bordered Heegaard diagrams for diffeomorphisms]{\textbf{Bordered Heegaard diagrams for diffeomorphisms.}
    Left: the identity bordered Heegaard diagram for the genus $1$
    surface. Center: the same diagram, with a regular neighborhood of
    $\alpha_1^{a,L}\cup\dots\cup\alpha_{2k}^{a,L}\cup\bdy_L\overline{\Sigma}\setminus\{\arcz\}$ in
    dark gray and a closed curve $\gamma$ as a dashed line. Right: a
    bordered Heegaard diagram for the Dehn twist around the curve
    $\gamma$.}
  \label{fig:hd-for-mcg}
\end{figure}

Every element of $\MCG_0(F)$ arises this
way. Indeed, given a diffeomorphism $\psi\co F\to F$, we can
explicitly construct a bordered Heegaard diagram realizing $\psi$ as
follows. Given a pointed matched circle $\PtdMatchCirc$, let
$\Id_{\PtdMatchCirc}=(\bSigma_{k},\balphas,\betas,\arcz)$ be the bordered
Heegaard diagram with $-\bdy_L\Id_\PMC=\bdy_R\Id_\PMC=\PMC$, with no
$\alpha$-circles and where $\balpha_i^{a,L}$
and $\balpha_i^{a,R}$ each
intersect $\beta_i$ in a single point and are disjoint from $\beta_j$
for $i\neq j$. (See Figure~\ref{fig:hd-for-mcg}.) Then
$\Id_{\PtdMatchCirc}$ is a bordered Heegaard diagram for the identity
map of $F(\PtdMatchCirc)$. Further, a regular neighborhood of
\[
\overline{\alpha}_1^{a,L}\cup\dots\cup{\overline{\alpha}}_{2k}^{a,L}\cup\bdy_L\overline{\Sigma}
\]
in $\overline{\Sigma}$ is canonically homeomorphic to~$F$ minus two
disks.  By removing the arc~$\arcz$ from this
neighborhood, we obtain $F$ minus one disc. Under this identification,
let $\overline{\psi}\co\bSigma\to\bSigma$ be the result of extending
$\psi\in\MCG_0(F)$ by the identity. Then,
$(\overline{\Sigma},\overline{\psi}(\balphas),\betas,\arcz)$ is a bordered Heegaard
diagram for $\psi$. (Equivalently, we could have applied
$\overline{\psi}^{-1}$ to $\betas$.)

The hypothesis that $F_L=-F_R$ in the above construction is not crucial.
In the more general case (where $F_L$ and $-F_R$ are different), one obtains Heegaard diagrams
corresponding to elements of a \emph{mapping class groupoid.}
See~\cite{LOT2} for further details.
\index{mapping class group!-oid}%

\subsection{The definition of \textalt{$\CFDAa$}{CFDA\textasciicircum}}\label{sec:CFDA-def}
Let $\HD=(\bSigma_g,\balphas,\betas,\arcz)$ be a bordered Heegaard diagram for
$\psi\in\MCG_0(F_R)$. The bimodule $\CFDAa(\psi)=\CFDAa(\HD)$ is defined
by counting holomorphic curves in $\Sigma\times[0,1]\times\RR$. As
usual, we choose an almost complex structure so that the two boundary
components of $\overline{\Sigma}$ become
cylindrical ends of $\Sigma$, and hence of
$\Sigma\times[0,1]\times\RR$.  We call the end corresponding to
$\bdy_L\overline{\Sigma}$ the \emph{left boundary}\index{left boundary}
and the end
corresponding to $\bdy_R\overline{\Sigma}$ \emph{right boundary}.
\index{right boundary}%

The bimodule $\CFDAa(\HD)$ is defined by treating the left boundary
as a type $D$ end (i.e., analogously to
Chapter~\ref{chap:type-d-mod}) and the right boundary a type $A$ end (i.e.,
analogously to Chapter~\ref{chap:type-a-mod}). 

More precisely, let
$\gls*{BimGens}$
denote the set of $g$-tuples of points $\x$ in $\Sigma$ so
that  exactly one $x_i$ lies on each $\alpha$- and each $\beta$-circle,
and no two $x_i$ lie on the same $\alpha$-arc.
Given $\x\in\S(\HD)$, let $o^L(\x)$ denote the set of $\alpha^L$-arcs
occupied by $\x$ and $o^R(\x)$ the set of $\alpha^R$-arcs occupied by
$\x$.\glsadd{oLxoRx} Let $X(\HD)$ be the $\Field$--vector space spanned by
$\S(\HD)$. We define left and right actions of
$\Idem(F_R)$ on $\S(\HD)$ by
\[
I(\SetS)\cdot\x\cdot I(\SetT)\coloneqq
\begin{cases}
  \x & \SetS=[2k]\setminus o^L(\x)\text{ and }\SetT=o^R(\x)\\
  0 &\text{otherwise},
\end{cases}
\]
where $\SetS$ and $\SetT$ are subsets of $[2k]$, as in
Formula~\eqref{eq:idempotent-def-2}.

As a left module, $\CFDAa(\HD)$ is just
$\Alg(F_L)\otimes_{\Idem}X(\HD)$. Note that, for the first time in
this book, we have a module on which $\Alg(\PtdMatchCirc,i)$ may act
non-trivially for $i\neq 0$.

Our next task is to define the $\Ainf$ bimodule structure on
$\CFDAa(\HD)$. (See, for instance, \cite[Section 3]{AinftyAlg}
for basic notions relating to $\Ainf$ bimodules.) The
left and right actions strongly commute, i.e.,
\[
\gls*{mCFDA}=0
\]
if $i,j > 0$, it only remains to describe the right action. The bimodule
$\CFDAa(\HD)$ is also strictly
unital, i.e.,
\begin{align*} 
  m_{0,1,1}(\x,\Unit) &= \x \\
  m_{0,1,n}(\x,\dots,\Unit,\dots) &= 0,\quad\textrm{$n>1$}.
\end{align*}

The definition of the operations
$m_{0,1,n}(\x,a(\rhos_1),\dots,a(\rhos_n))$ uses holomorphic curves. Let 
$
\Mod^B(\x,\y;\vec{\rho}^D;\vec\rhos^A)
$
denote the moduli space of embedded holomorphic curves in the homology
class $B$, asymptotic to
the sequence of Reeb chords in $\vec{\rho}^D$ at west $\infty$ (with no height
constraints), and to the ordered partition
$\vec\rhos^A = (\rhos_1^A,\dots,\rhos_n^A)$ of Reeb chords at east
$\infty$.
\index{infinity!east}\index{infinity!west}%
\index{east infinity!bimodule case}%
\index{west infinity!bimodule case}%
The rest of the bimodule structure is determined by
\[
m_{0,1,n}(\x,a(\rhos_1),\dots,a(\rhos_n))\coloneqq
  \!\!\!\!\sum_{\y\in\S(\HD)}
   \!\!\!\!\!\!\sum_{\substack{B\in\pi_2(\x,\y)\\
        \{\vec\rho^D\mid\ind(B;\vec{\rho}^D;\vec\rhos^A)= 1\}}}\!\!\!\!\!\!\!\!\!\!
    \big(\#\Mod^B(\x,\y;\vec{\rho}^D;\vec\rhos^A)\big)a(-\vec{\rho}^D)\y.
\]
\glsadd{mCFDA}%
(Compare Definitions~\ref{def:Dmod-boundary}
and~\ref{def:Amod-mult}. The case $n=0$
is essentially the differential on $\CFDa$.)

We prove in~\cite{LOT2} that this does, indeed, define an
$\Ainf$ bimodule. Invariance, the first part of
Theorem~\ref{thm:bimodule-pairing}, is proved similarly to invariance
of $\CFDa$ or $\CFAa$. Both of the proofs of the Pairing Theorem (Chapters
\ref{chap:nice-diagrams} and~\ref{chap:tensor-prod}) extend easily to
prove the second and third parts of Theorem~\ref{thm:bimodule-pairing}. Again,
this is explained in more detail in~\cite{LOT2}.

\begin{remark}\index{grading!on $\CFDAa$}
  Just as $\CFDa$ and $\CFAa$ have been graded by left or right
  $\smallGroup$-sets, bimodules are graded by sets with left and right
  actions by $\smallGroup$. The bimodule $\CFDAa(\psi)$ is graded by
  $\smallGroup(F_R)$, thought of as a set with commuting left and
  right actions by $\smallGroup(F_R)\cong\smallGroup(F_L)$. The
  construction of the grading
  is a simple modification of the discussions in
  Sections~\ref{sec:typeD-gradings} and~\ref{sec:typeA-gradings}.
\end{remark}

\subsection{The definition of \textalt{$\CFDDa(\Id)$}{CFDD\textasciicircum(Id)} and \textalt{$\CFAAa(\Id)$}{CFAA\textasciicircum(Id)}}\label{sec:def-CFAA-CFDD-id}
The bimodules $\CFDDa(\Id)$ and $\CFAAa(\Id)$ are both defined in terms
of holomorphic curves with respect to the identity Heegaard diagram
$\Id_{\PtdMatchCirc}=(\bSigma,\balphas,\betas,\arcz)$.
The bimodule $\CFDDa(\Id)$ is defined by treating both boundary
components as type $D$ boundary components, and the bimodule
$\CFAAa(\Id)$ is defined by treating both boundary components as
type $A$ boundary.

\glsadd{CFDDId}%
More precisely, $\CFDDa(\Id)$ is defined as follows.
Recall that $-\PtdMatchCirc$ is the pointed matched circle obtained by
reversing the orientation on $\PtdMatchCirc$.\index{orientation reversal}  The bimodule
$\CFDDa(\Id)$ is generated over
$\Alg(\PtdMatchCirc) \times \Alg(-\PtdMatchCirc)$
by the $2^{2k}=\sum_i\binom{2k}{i}$ ways of choosing $i$ points from
$\{\overline{\alpha}_1^{L,a}\cap\beta_1,\overline{\alpha}_2^{L,a}\cap\beta_2,\dots,\overline{\alpha}_{2k}^{L,a}\cap\beta_{2k}\}$
(and the $2k-i$ points on the complementary $\alpha_i^{R,a}$); we can
identify this with
the set $\S(\Id_{\PtdMatchCirc})$ from
Section~\ref{sec:CFDA-def} by looking at $o^L$, the set of strands
occupied on the left. There are left actions of
$\Idem(-\PtdMatchCirc)$ and $\Idem(\PtdMatchCirc)$ on
$X(\Id_{\PtdMatchCirc})$ by, respectively,
\begin{align*}
I(\SetS)\bullet\x&\coloneqq
\begin{cases}
  \x & \SetS=[2k]\setminus o^L(\x)\\
  0 &\text{otherwise},
\end{cases}\\
I(\SetS)\circ\x&\coloneqq
\begin{cases}
  \x & \SetS=[2k]\setminus o^R(\x)\\
  0 &\text{otherwise}.
\end{cases}
\end{align*}
(Note that $o^R(\x) = [2k]\setminus o^L(\x)$.)
These actions commute. Consequently, extending scalars from the
idempotents to the whole algebra
gives commuting left actions of $\Alg(\PtdMatchCirc)$ and $\Alg(-\PMC)$; this is the
bimodule structure on $\CFDDa(\Id)$.

The differential on $\CFDDa(\Id)$ is given as follows. Let
$
\Mod^B(\x,\y;\vec\rho^L,\vec\rho^R)
$
denote the moduli space of embedded holomorphic curves in the homology
class $B$, asymptotic to the sequence of Reeb chords
$\vec\rho^L$ on the left and $\vec\rho^R$ on the right.  (We do not
constrain the relative ordering of Reeb chords on the left and right.)
Then\glsadd{CFDDIdbdy}
\begin{equation}
\bdy\x\coloneqq
  \sum_{\y\in\S(\HD)}
  \!\sum_{\substack{B\in\pi_2(\x,\y)\\
      \{\vec\rho^L,\vec\rho^R\mid\ind(B;\vec{\rho}^L;\vec\rho^R)= 1\}}}\!\!\!\!\!
  \#\left(\Mod^B(\x,\y;\vec{\rho}^L;\vec\rho^R)\right)a(-\vec{\rho}^L)\bullet
  a(-\vec\rho^R)\circ\y.\label{eq:dd-id-diff}
\end{equation}

Next we turn to $\CFAAa(\Id)$.
\glsadd{CFAAId}%
As an $\FF_2$-vector space, the bimodule $\CFAAa(\Id)$ is just
$X(\Id_{\PtdMatchCirc})$.  There are commuting right actions of
$\Idem(-\PtdMatchCirc)$ and $\Idem(\PtdMatchCirc)$ on
$X(\Id_{\PtdMatchCirc})$ by
\begin{align*}
\x\bullet I(\SetS)&\coloneqq
\begin{cases}
  \x & \SetS=o^L(\x)\\
  0 &\text{otherwise},
\end{cases}\\
\x\circ I(\SetS)&\coloneqq
\begin{cases}
  \x & \SetS=o^R(\x)\\
  0 &\text{otherwise}.
\end{cases}\\
\end{align*}

The bimodule $\CFAAa(\Id)$ is strictly unital, i.e., 
\begin{align*}
  m_{1,1,0}(\x;\Unit;)&\coloneqq m_{1,0,1}(\x;;\Unit)\coloneqq \x \\
  m_{1,m,n}(\x;\dots,\Unit,\dots;\dots) &\coloneqq m_{1,m,n}(\x;\dots;\dots,\Unit,\dots)\coloneqq 0,\quad\textrm{$m+n>1$}.
\end{align*}\glsadd{mCFAA}The other higher multiplications are defined
by counting holomorphic curves.
Let
$
\Mod^B(\x,\y;\vec\rhos^L;\vec\etas^R)
$
denote the moduli space of embedded holomorphic curves in the homology
class $B$, asymptotic to the ordered partition
$\vec\rhos = (\rhos_1,\dots,\rhos_m)$ at
$(\bdy_L\overline{\Sigma})\times[0,1]\times\RR$ and asymptotic to the
ordered partition $\vec\etas = (\etas_1,\dots,\etas_n)$ at
$(\bdy_R\overline{\Sigma})\times[0,1]\times\RR$. Then define\glsadd{mCFAA}
\begin{multline*}
m_{1,m,n}(\x;a(\rhos_1),\dots,a(\rhos_m);a(\etas_1),\dots,a(\etas_n))
\coloneqq\\
  \sum_{\y\in\S(\HD)}
   \sum_{\substack{B\in\pi_2(\x,\y)\\
        \ind(B;\vec{\rhos};\vec\etas)= 1}}\!\!\!
    \#\left(\Mod^B(\x,\y;\vec{\rhos};\vec\etas)\right)\y.
\end{multline*}

Theorem~\ref{thm:dualizing-bimods} now follows from a version of the
pairing theorem, which can be proved using either the nice diagrams
technique of Chapter~\ref{chap:nice-diagrams} or the time dilation
argument of Chapter~\ref{chap:tensor-prod}.

\begin{remark}\index{grading!on $\CFAAa(\Id)$ and $\CFDDa(\Id)$}
  The gradings on $\CFAAa(\Id)$ and $\CFDDa(\Id)$ are defined
  analogously to Sections~\ref{sec:typeA-gradings}
  and~\ref{sec:typeD-gradings} respectively. In the case of
  $\CFAAa(\Id)$ (respectively $\CFDDa(\Id)$) the $\smallGroup$-grading
  set is $\smallGroup(\PMC)$, viewed as having two right (respectively
  left) actions, one by $\smallGroup(\PMC)$ and one by
  $\smallGroup(-\PMC)\cong \smallGroup(\PMC)^\op$. In fact, both of
  these actions are the standard action of the group on itself by
  translation.
\end{remark}

\begin{remark}
  For each of $\CFDAa(\HD)$, $\CFDDa(\Id)$ and $\CFAAa(\Id)$, the
  summands $\Alg(\PtdMatchCirc,i)$ inside $\Alg(\PtdMatchCirc)$ act
  non-trivially even for $i\neq 0$, and the bimodules decompose as direct
  sums
  \begin{align*}
    \CFDAa(\HD)&=\bigoplus_{i=-k}^{k}\CFDAa(\HD,i)\\
    \CFDDa(\Id)&=\bigoplus_{i=-k}^{k}\CFDDa(\Id,i)\\    
    \CFAAa(\Id)&=\bigoplus_{i=-k}^{k}\CFAAa(\Id,i).
  \end{align*}
  Of course, only the summands corresponding to $i=0$ contribute to
  tensor products with $\CFDa$ or $\CFAa$ as in
  Theorem~\ref{thm:bimodule-pairing}. The other summands do arise
  naturally if one studies self-gluing of Heegaard diagrams, which in
  turn relates to the Floer homology of open books.
\end{remark}

\section{Computations for \textalt{$3$}{3}-manifolds with torus boundary}\label{sec:comp-torus-bimods}

In this section we state the relevant bimodules for
three-manifolds with torus boundary. These are used in
Section~\ref{sec:hfk-to-CFD-general} to
deduce a version of Theorem~\ref{thm:HFKtoHFD} for
arbitrary integral framings of knots,
Theorem~\ref{thm:HFKtoHFDframed}. (In fact, the data from this
section can be used to calculate the result for arbitrary rational framings;
but that will be left to the interested reader.) The computations of the
bimodules are given in~\cite{LOT2}.

There is a unique pointed matched circle representing the surface of
genus one. Let $\Alg=\Alg(\PMC,0)$ denote the $i=0$ summand of the
corresponding algebra. Label the generators of $\Alg$ as in
Chapter~\ref{chap:TorusBoundary}.

\subsection{Type \DAm\ invariants of Dehn twists}
The mapping class group is generated by Dehn twists\index{Dehn twists}
$\gls*{taumu}$ and $\gls*{taulambda}$
along a meridian and a longitude,
respectively (i.e., $\tau_\mu$ takes an $n$-framed knot complement to
an $n+1$-framed knot complement).
Recall from Section~\ref{sec:CFD-to-HFK} that $\alpha_2^a$ is a
meridian. Thus, $\tau_\mu$ corresponds to a Dehn twist around a curve 
$\mu$ dual to the curve $\alpha_1^a$ from Figure~\ref{fig:torus-diagrams}. 

Most of this section is devoted to describing the $i=0$ part
$\CFDAa(\cdot,0)$ of the type \DAm\ bimodules for these Dehn twists
and their inverses.  First, we dispense with the $i=\pm 1$ parts. The
summand $\Alg(\Torus,-1)$ of $\Alg(\Torus)$ is isomorphic to $\Field$, and
$\CFDAa(\phi,-1)\cong \Field$ for any mapping class $\phi$. The summand
$\Alg(\Torus,1)$ of $\Alg(\Torus)$ is $7$-dimensional, but is
quasi-isomorphic to $\Field$. Thus, $\Alg(\Torus,1)$-bimodules are
determined (up to quasi-isomorphism) by corresponding
$\Field$-bimodules. For any mapping class $\phi$, the $\Field$-bimodule
corresponding to $\CFDAa(\phi,1)$ is, again, $\Field$.

\begin{figure}
  \input{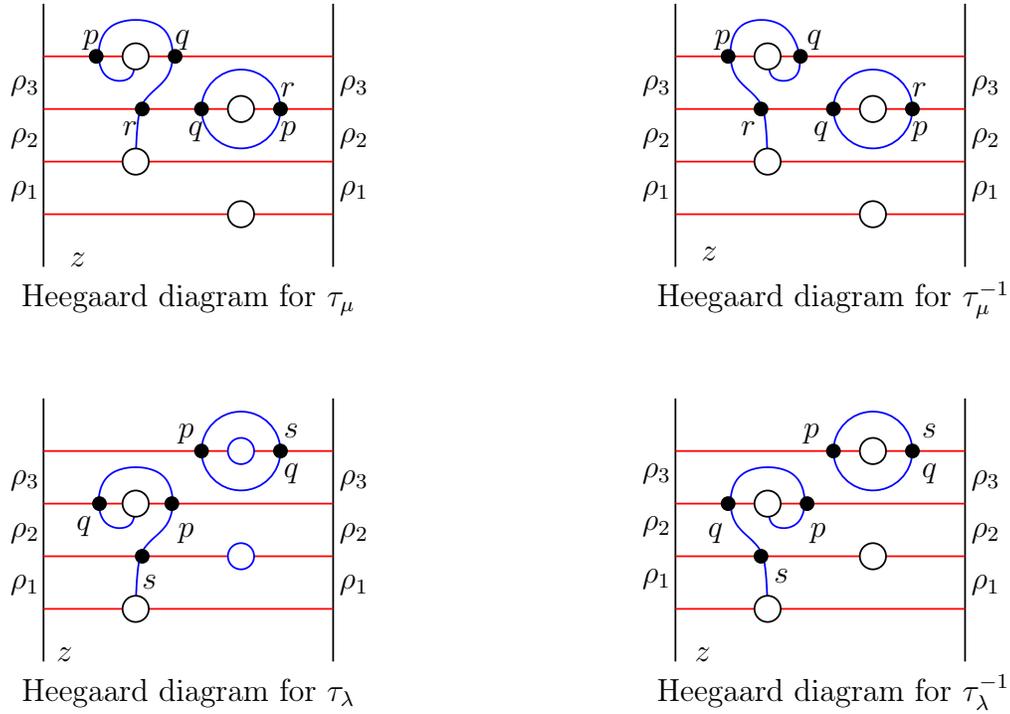}
  \caption[Heegaard diagrams for
  Dehn twists of the torus]{\label{fig:DehnTwistsGenusOne} {\bf{Heegaard diagrams for
        mapping class group elements.}}  Genus $2$ diagrams for
    $\tau_\mu$, $\tau_\mu^{-1}$, $\tau_\lambda$ and $\tau_\lambda^{-1}$
    are shown. In each of the four
    diagrams, there are three
    generators in the $i=0$ summand.}
\end{figure}

\begin{figure}
\begin{center}
  \begin{tikzpicture}[y=45pt,x=.8in]
    \node at (0,2) (p) {${\mathbf p}$} ;
    \node at (2,2) (q) {${\mathbf q}$} ;
    \node at (1,0) (r) {${\mathbf r}$} ;
    \node at (1,1.25) (label) {$\boldsymbol{\tau_\mu}$};
    \draw[->] (p) to node[above,sloped] {\lab{\rho_1\otimes\rho_1+\rho_{123}\otimes\rho_{123}+\rho_3\otimes (\rho_{3},\rho_{23})}}  (q)  ;
    \draw[->] (p) [bend left=15] to node[above,sloped] {\lab{\rho_{123}\otimes\rho_{12}+
\rho_3\otimes (\rho_{3},\rho_{2})}} (r) ;
    \draw[->] (q) [bend left=15] to node[below,sloped] {\lab{\rho_{23}\otimes \rho_{2}}} (r) ;
    \draw[->] (q) [loop] to node[above] {\lab{\rho_{23}\otimes\rho_{23}}}
                 (q) ;
    \draw[->] (r) [bend left=15] to node[below,sloped] {\lab{\rho_2\otimes 1}} (p) ;
    \draw[->] (r) [bend left=15] to node[below,sloped] {\lab{1\otimes \rho_3}} (q) ;
  \end{tikzpicture}
  \begin{tikzpicture}[y=45pt,x=.8in]
    \node at (0,2) (p) {${\mathbf p}$} ;
    \node at (2,2) (q) {${\mathbf q}$} ;
    \node at (1,0) (r) {${\mathbf r}$} ;
    \node at (1,1.25) (label) {$\boldsymbol{\tau_\mu^{-1}}$};
    \draw[->] (p) [bend left=15] to node[above,sloped] {\lab{\rho_1\otimes \rho_1+\rho_{123}\otimes\rho_{123}}}  (q)  ;
    \draw[->] (p) [loop] to node[above] {\lab{\rho_{12}\otimes(\rho_{123},\rho_2)}} (p) ;
    \draw[->] (p) [bend left=15] to node[above,sloped]  {\lab{\rho_3\otimes 1 + \rho_1\otimes \rho_{12}}} (r) ;
    \draw[->] (q) [bend left=15] to node[below,sloped] {\lab{\rho_2\otimes(\rho_{23},\rho_2)}} (p) ;
    \draw[->] (q) [bend left=15] to node[below,sloped] {\lab{1\otimes \rho_2}} 
                  (r) ;
    \draw[->] (q) [loop] to node[above] {\lab{\rho_{23}\otimes\rho_{23}}} (q) ;
    \draw[->] (r) [bend left=15] to node[below,sloped] {\lab{\rho_2\otimes (\rho_3,\rho_2)}} (p) ;
    \draw[->] (r) [bend left=15] to node[below,sloped] {\lab{\rho_{23}\otimes \rho_3}} (q) ;
  \end{tikzpicture} \\
  \begin{tikzpicture}[y=45pt,x=.8in]
    \node at (0,2) (p) {${\mathbf q}$} ;
    \node at (2,2) (q) {${\mathbf p}$} ;
    \node at (1,0) (r) {${\mathbf s}$} ;
    \node at (1,1.25) (label) {$\boldsymbol{\tau_\lambda}$};
    \draw[->] (p) [bend left=15] to node[above,sloped] {\lab{\rho_2\otimes(\rho_2,\rho_{12})}}  (q)  ;
    \draw[->] (p) [bend left=15] to node[above,sloped] {\lab{\rho_2\otimes (\rho_2,\rho_1)}} (r) ;
    \draw[->] (p) [loop] to node[above]  {\lab{\rho_{23}\otimes (\rho_2,\rho_{123})}}  (p) ;
    \draw[->] (q) [bend left=15] to node[below,sloped] {\lab{\rho_{12}\otimes\rho_1}} (r) ;
    \draw[->] (q) [loop] to node[above] {\lab{\rho_{12}\otimes\rho_{12}}} (q) ;
    \draw[->] (q) [bend left=15] to node[below,sloped] {\lab{\rho_3\otimes\rho_3+\rho_{123}\otimes\rho_{123}}} (p) ;
    \draw[->] (r) [bend left=15] to node[below,sloped] {\lab{\rho_1\otimes 1+\rho_3\otimes\rho_{23}}} (p) ;
    \draw[->] (r) [bend left=15] to node[below,sloped] {\lab{1\otimes \rho_2}} (q) ;
  \end{tikzpicture}
  \begin{tikzpicture}[y=45pt,x=.8in]
    \node at (0,2) (p) {${\mathbf q}$} ;
    \node at (2,2) (q) {${\mathbf p}$} ;
    \node at (1,0) (r) {${\mathbf s}$} ;
    \node at (1,1.25) (label) {$\boldsymbol{\tau_\lambda^{-1}}$};
    \draw[->] (p) [bend left=15] to node[above,sloped] {\lab{\rho_2\otimes 1}}  (r)  ;
    \draw[->] (q) [loop] to node[above] {\lab{\rho_{12}\otimes\rho_{12}}}
                (q) ;
    \draw[->] (q) [bend left=15] to node[below,sloped]  {\lab{1\otimes \rho_1}} (r) ;
    \draw[->] (q)  to node[above,sloped] {\lab{\rho_3\otimes\rho_3+\rho_{123}\otimes\rho_{123}+\rho_1\otimes(\rho_{12},\rho_1)}} (p) ;
    \draw[->] (r) [bend left=15] to node[below,sloped] {\lab{\rho_1\otimes (\rho_2,\rho_1)+\rho_{123}\otimes\rho_{23}}} (p) ;
    \draw[->] (r) [bend left=15] to node[above,sloped] {\lab{\rho_{12}\otimes\rho_2}} (q) ;
  \end{tikzpicture}
\end{center}
\caption[\DAm\ bimodules for generators of the genus-one mapping
class group] {
\label{fig:DAbimodules}
Type \DAm\ bimodules for torus mapping class group action.  These are
the bimodules associated to $\tau_\mu$, $\tau_\mu^{-1}$,
$\tau_{\lambda}$, $\tau_{\lambda}^{_1}$ respectively. The
notation is as follows. Consider the bimodule for $\tau_\mu$. The
label $\rho_1\otimes\rho_1$ on the horizontal arrow indicates that
$m_{0,1,1}(\mathbf{p},\rho_1)$ contains a term of the form
$\rho_1\otimes\mathbf{q}$. Similarly, the label
$\rho_3\otimes(\rho_3,\rho_{23})$ on that arrow indicates that
$m_{0,1,2}(\mathbf{p},\rho_3,\rho_{23})$ contains a term of the form
$\rho_3\otimes \mathbf{q}$.}
\end{figure}
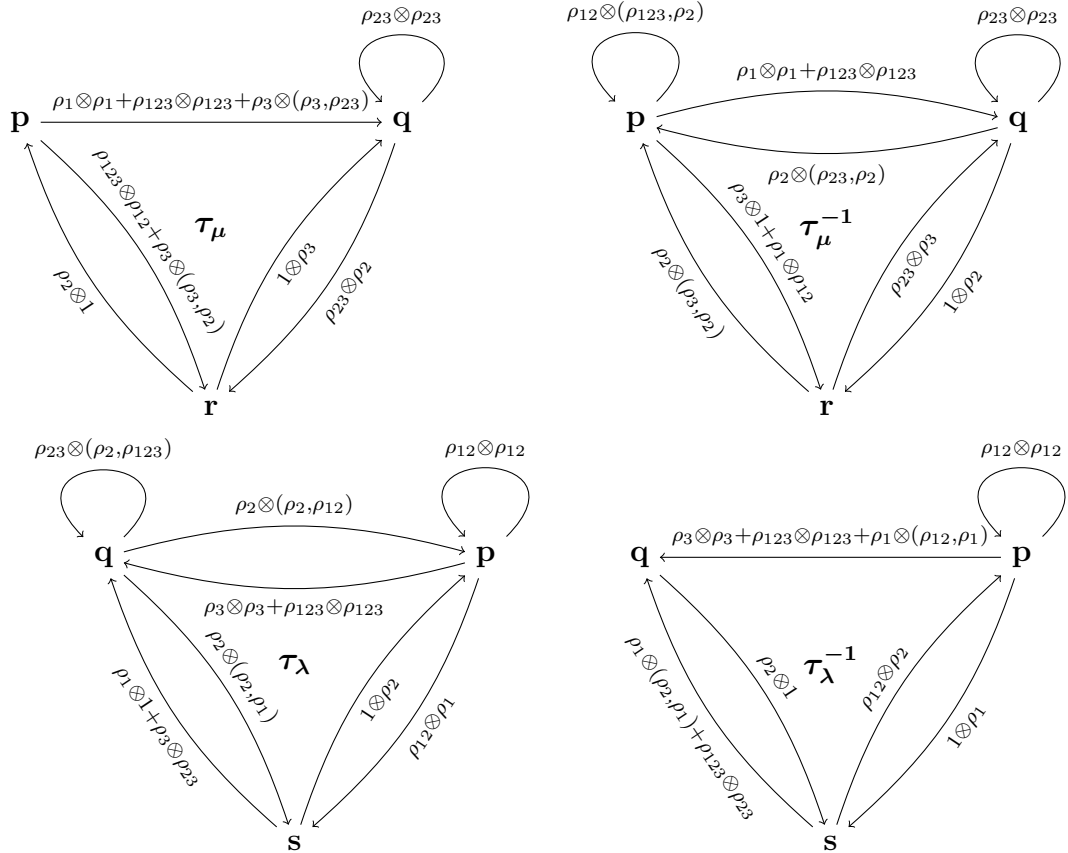

Heegaard diagrams for $\tau_\mu$, $\tau_\mu^{-1}$, $\tau_\lambda$ and
$\tau_\lambda^{-1}$ are illustrated in
Figure~\ref{fig:DehnTwistsGenusOne}.  Each of the four bimodules
$\CFDAa(\tau_\mu,0)$, $\CFDAa(\tau_\mu^{-1},0)$,
$\CFDAa(\tau_\lambda,0)$ and $\CFDAa(\tau_\lambda^{-1},0)$
has the form $\Alg\otimes_\Idem X$
where $X$ is $3$-dimensional over $\Field$. The vector space $X$ has
a basis consisting of $\mathbf{p}$, $\mathbf{q}$ and either
$\mathbf{r}$ or $\mathbf{s}$.
\glsadd{pqrsgens}%
The generators are compatible with the idempotents as follows:
\begin{equation}
\iota_0\cdot {\mathbf p}\cdot \iota_0 = {\mathbf p} \qquad
\iota_1\cdot {\mathbf q}\cdot \iota_1 = {\mathbf q} \qquad
\iota_1\cdot {\mathbf r}\cdot \iota_0 = {\mathbf r} \qquad
\iota_0\cdot {\mathbf s}\cdot \iota_1 = {\mathbf s}. \label{eq:pqrsIdems}
\end{equation}
As in Section~\ref{sec:CFDA-def}, the left module structure is
induced from the isomorphism $\CFDAa\cong \Alg\otimes_\Idem X$. Higher
products of the form $m_{i,1,j}(a_1,\dots,a_i,\x,b_1,\dots,b_j)$
vanish if $ij>0$. So, it remains to give the higher products of the
form $m_{0,1,j}$ for $j\geq 0$. We do this for each of the four
bimodules in turn; see also Figure~\ref{fig:DAbimodules}.

The bimodule $\CFDAa(\tau_\mu,0)$ is generated by $\mathbf{p}$,
$\mathbf{q}$ and $\mathbf{r}$. The idempotents act as above; the other
non-trivial algebra actions are given as follows:
\begin{align*}
  m_{0,1,1}({\mathbf p},\rho_1) &= \rho_1\otimes {\mathbf q} &
  m_{0,1,1}({\mathbf p}, \rho_{12})&= \rho_{123}\otimes {\mathbf r} \\
  m_{0,1,1}({\mathbf p},\rho_{123}) &= \rho_{123}\otimes {\mathbf q} &
  m_{0,1,2}({\mathbf p}, \rho_3,\rho_2) &= \rho_3\otimes {\mathbf r} \\
  m_{0,1,2}({\mathbf p}, \rho_3,\rho_{23}) &= \rho_3 \otimes {\mathbf q} &
  m_{0,1,1}({\mathbf q}, \rho_2) &= \rho_{23}\otimes {\mathbf r} \\
  m_{0,1,1}({\mathbf q}, \rho_{23}) &= \rho_{23}\otimes {\mathbf q} &
  m_{0,1,0}({\mathbf r})&=\rho_2\otimes {\mathbf p} \\
  m_{0,1,1}({\mathbf r},\rho_3)&={\mathbf q}.
\end{align*}

The bimodule $\CFDAa(\tau_\mu^{-1},0)$ is generated by $\mathbf{p}$,
$\mathbf{q}$ and $\mathbf{r}$. The idempotents act as above; the other
non-trivial algebra actions are given as follows:
\begin{align*}
  m_{0,1,0}({\mathbf p}) &= \rho_{3}\otimes {\mathbf r} &
  m_{0,1,1}({\mathbf p}, \rho_1) &= \rho_1\otimes {\mathbf q} \\
  m_{0,1,1}({\mathbf p},\rho_{12}) &= \rho_{1}\otimes {\mathbf r} &
  m_{0,1,1}({\mathbf p}, \rho_{123}) &= \rho_{123}\otimes {\mathbf q} \\
  m_{0,1,2}({\mathbf p},\rho_{123},\rho_2) &= \rho_{12}\otimes {\mathbf p} &
  m_{0,1,1}({\mathbf q},\rho_{2}) &= {\mathbf r} \\
  m_{0,1,1}({\mathbf q},\rho_{23}) &= \rho_{23}\otimes {\mathbf q} &
  m_{0,1,2}({\mathbf q},\rho_{23},\rho_2) &= \rho_{2}\otimes {\mathbf p} \\
  m_{0,1,1}({\mathbf r},\rho_{3}) &= \rho_{23}\otimes {\mathbf q} &
  m_{0,1,2}({\mathbf r},\rho_{3},\rho_2) &= \rho_{2}\otimes {\mathbf p}.
\end{align*}

The bimodule $\CFDAa(\tau_\lambda,0)$  is generated by $\mathbf{p}$,
$\mathbf{q}$ and $\mathbf{s}$. The idempotents act as above; the other
non-trivial algebra actions are given as follows:
\begin{align*}
  m_{0,1,2}({\mathbf q},\rho_2,\rho_{1})&= \rho_{2}\otimes {\mathbf s} &
  m_{0,1,2}({\mathbf q},\rho_2,\rho_{12})&= \rho_{2}\otimes {\mathbf p} \\
  m_{0,1,2}({\mathbf q},\rho_2,\rho_{123})&= \rho_{23}\otimes {\mathbf q} &
  m_{0,1,1}({\mathbf p}, \rho_1) &= \rho_{12}\otimes {\mathbf s} \\
  m_{0,1,1}({\mathbf p},\rho_{12}) &= \rho_{12}\otimes {\mathbf p} &
  m_{0,1,1}({\mathbf p},\rho_{123})&= \rho_{123}\otimes {\mathbf q} \\
  m_{0,1,1}({\mathbf p},\rho_3) &= \rho_3\otimes {\mathbf q} &
  m_{0,1,0}({\mathbf s}) &= \rho_1 \otimes {\mathbf q} \\
  m_{0,1,1}({\mathbf s},\rho_2) &= {\mathbf p} &
  m_{0,1,1}({\mathbf s},\rho_{23}) &= \rho_3 \otimes {\mathbf q}.
\end{align*}

The bimodule $\CFDAa(\tau_\lambda^{-1},0)$  is generated by $\mathbf{p}$,
$\mathbf{q}$ and $\mathbf{s}$. The idempotents act as above; the other
non-trivial algebra actions are given as follows:
\begin{align*}
  m_{0,1,0}({\mathbf q})&= \rho_2\otimes {\mathbf s} &
  m_{0,1,1}({\mathbf p},\rho_{1}) &= {\mathbf s} \\
  m_{0,1,1}({\mathbf p},\rho_{12}) &= \rho_{12}\otimes {\mathbf p} &
  m_{0,1,1}({\mathbf p},\rho_{123}) &= \rho_{123}\otimes {\mathbf q} \\
  m_{0,1,1}({\mathbf p},\rho_{3}) &= \rho_{3}\otimes {\mathbf q} &
  m_{0,1,2}({\mathbf p},\rho_{12},\rho_1) &= \rho_{1}\otimes {\mathbf q} \\
  m_{0,1,1}({\mathbf s},\rho_{2}) &= \rho_{12}\otimes {\mathbf p} &
  m_{0,1,1}({\mathbf s},\rho_{23}) &= \rho_{123}\otimes {\mathbf q} \\
  m_{0,1,2}({\mathbf s},\rho_{2},\rho_1) &= \rho_{1}\otimes {\mathbf q}.
\end{align*}

\begin{remark}
The computations leading to the results above are eased by noting that
\begin{itemize}
\item $\CFDAa(\tau_\mu^{-1}, 0)$ is inverse to $\CFDAa(\tau_\mu, 0)$,
  i.e., 
  \[
  \CFDAa(\tau_\mu^{-1}, 0) \DT \CFDAa(\tau_\mu, 0) \simeq
  \CFDAa(\Id, 0);
  \]
\item the Heegaard diagram for $\tau_\lambda$ is obtained from
  that for $\tau_\mu^{-1}$ by a horizontal reflection (across the $x$-axis); and
\item the Heegaard diagram for $\tau_\lambda^{-1}$ is obtained from
  that for $\tau_\mu$ by a horizontal reflection.
\end{itemize}
Thus there is essentially only one computation to do.
\end{remark}
\begin{remark}
  As the reader may well have expected, up to homotopy
  equivalence the type \DAm\  bimodule $\CFDAa(\Id,0)$
  associated to the identity map is isomorphic to $\Alg$ as an
  $\Alg$--$\Alg$ bimodule.
\end{remark}
\subsection{The type \AAm\ invariant of the identity}
\begin{figure}
\begin{center}
\input{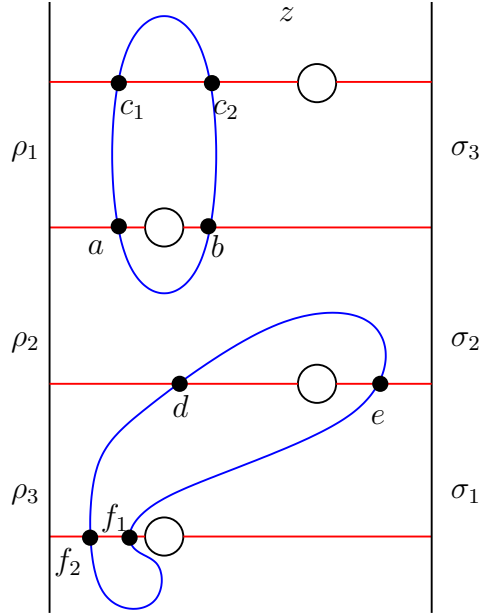}
\end{center}
\caption[Isotoped Heegaard diagram for identity map of torus]{\textbf{A Heegaard diagram for the genus $1$ identity map.}
  This diagram has four more intersection points in
  $\alphas\cap\betas$ than the standard Heegaard diagram for the
  identity map, introduced by two finger moves. As a result of these
  finger moves, the diagram is admissible, not merely provincially
  admissible, and is useful for computing $\CFAAa(\Id)$.}
\label{fig:AAinvpair}
\end{figure}

To compute the type \AAm\  and \DDm\  identity bimodules for the torus, it
is convenient to use the Heegaard diagram pictured in
Figure~\ref{fig:AAinvpair}.  Note that the boundary segments on both
sides have been labeled to be consistent with the type $A$ conventions.

The resulting type \AAm\ bimodule has six generators, denoted $bf_1$,
$bf_2$, $c_1e$, $c_2e$, $ae$ and $bd$. Its non-trivial
$\Ainf$-operations are given by:
\begin{align*}
  m_{1,1,0}(bf_1;\sigma_1;)&=ae&
  m_{1,0,0}(bf_1;;)&=bf_2\\
  m_{1,1,1}(bf_1;\sigma_{12};\rho_2)&=bd&
  m_{1,1,1}(bf_1;\sigma_{12};\rho_{23})&=bf_2\\
  m_{1,1,1}(bf_1;\sigma_{123};\rho_{2})&=c_2e&
  m_{1,3,1}(bf_1;\sigma_3,\sigma_2,\sigma_1;\rho_2)&=c_2e\\
  m_{1,0,0}(c_1e;;)&=c_2e&
  m_{1,0,1}(c_1e;;\rho_1)&=ae\\
  m_{1,1,1}(c_1e;\sigma_2;\rho_{12})&=bd&
  m_{1,1,1}(c_1e;\sigma_{23};\rho_{12})&=c_2e\\
  m_{1,1,1}(c_1e;\sigma_2;\rho_{123})&=bf_2&
  m_{1,1,1}(ae;\sigma_2;\rho_{23})&=bf_2\\
  m_{1,1,1}(ae;\sigma_2;\rho_2)&=bd&
  m_{1,1,1}(ae;\sigma_{23};\rho_2)&=c_2e\\
  m_{1,0,1}(bd;;\rho_3)&=bf_2&
  m_{1,1,0}(bd;\sigma_3;)&=c_2e.
\end{align*}
Here, to avoid confusion, we denote algebra elements acting on one
side (which we think of as $\Alg(-\PMC,0)$) by $\sigma_i$ while on the
other side (which we think of as $\Alg(\PMC,0)$) they are denoted by
$\rho_i$ (as indicated in Figure~\ref{fig:AAinvpair}).

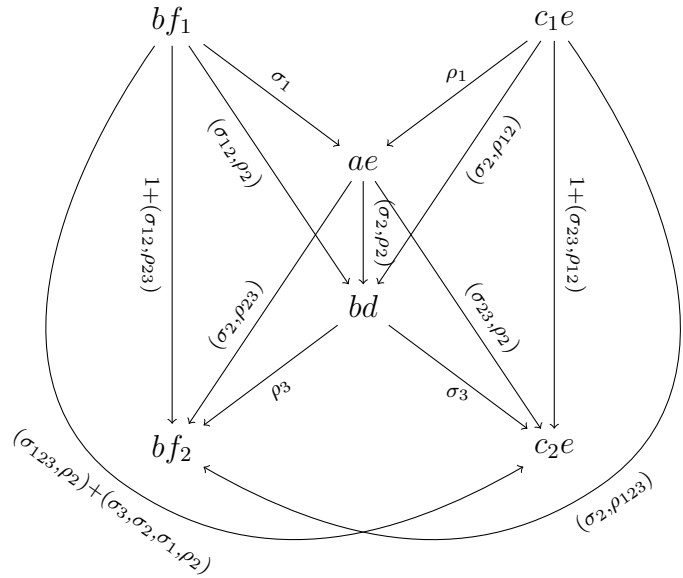
\begin{figure}
\begin{center}
  \begin{tikzpicture}[y=45pt,x=.8in]
    \node at (0,3) (w1) {$bf_1$} ;
    \node at (2,3) (z1) {$c_1 e$} ;
    \node at (1,2) (y)  {$ae$} ;
    \node at (1,1) (x)  {$bd$} ;
    \node at (0,0) (w2) {$bf_2$} ;
    \node at (2,0) (z2) {$c_2 e$} ;
    \draw[->] (w1) to node[above,sloped] {\lab{\sigma_1}} (y) ;
    \draw[->] (z1) to node[above] {\lab{\rho_1}} (y) ;
    \draw[->] (y)  to node[above,sloped] {\lab{(\sigma_2,\rho_2)}} (x) ;
    \draw[->] (x)  to node[below,sloped] {\lab{\rho_3}} (w2) ;
    \draw[->] (x)  to node[below] {\lab{\sigma_3}} (z2) ;
    \draw[->] (w1) to node[below,sloped] {\lab{1+(\sigma_{12},\rho_{23})}} (w2) ;
    \draw[->] (z1) to node[above,sloped] {\lab{1+(\sigma_{23},\rho_{12})}} (z2) ;
    \draw[->] (w1) to[pos=0.4] node[below,sloped] {\lab{(\sigma_{12},\rho_2)}} (x) ;
    \draw[->] (y)  to[pos=0.6] node[above,sloped] {\lab{(\sigma_2,\rho_{23})}} (w2) ;
    \draw[->] (z1) to[pos=0.4] node[below,sloped] {\lab{(\sigma_2,\rho_{12})}} (x) ;
    \draw[->] (y)  to[pos=0.6] node[above,sloped] {\lab{(\sigma_{23},\rho_2)}} (z2) ;
    \draw[->] (w1) to[out=-125,in=145] (-0.25,-0.25) to[out=-35,in=-150] node[pos=0,below,sloped]
       {\lab{(\sigma_{123},\rho_2)+(\sigma_3,\sigma_2,\sigma_1,\rho_2)}} (z2) ;
    \draw[->] (z1) to[out=-55,in=35] (2.25,-0.25) to[out=-145,in=-30] node[pos=0,below,sloped]
       {\lab{(\sigma_{2},\rho_{123})}} (w2) ;
  \end{tikzpicture}
\end{center}
\caption[The bimodule $\CFAAa(\Id,0)$ for the torus]{\textbf{The type \AAm\  bimodule $\CFAAa(\Id,0)$.}
  The labels on the arrows indicate the $\Ainf$
  operations; for instance, the label $(\sigma_{23},\rho_2)$ on the
  arrow from $ae$ to~$c_2 e$ means that $m_{1,1,1}(ae,\sigma_{23},\rho_2)$
  contains a term~$c_2 e$.}
\label{fig:AAmodAns}
\end{figure}

\begin{remark}
  The type \AAm\ identity bimodule of the torus has a model which is
  $2$-dimensional over $\Field$. However, this model has infinitely
  many non-trivial higher products. These can be described explicitly;
  see~\cite[Section 8.4]{LOT4}.
\end{remark}

\subsection{The type \DDm\ invariant of the identity}\label{sec:dd-of-id}
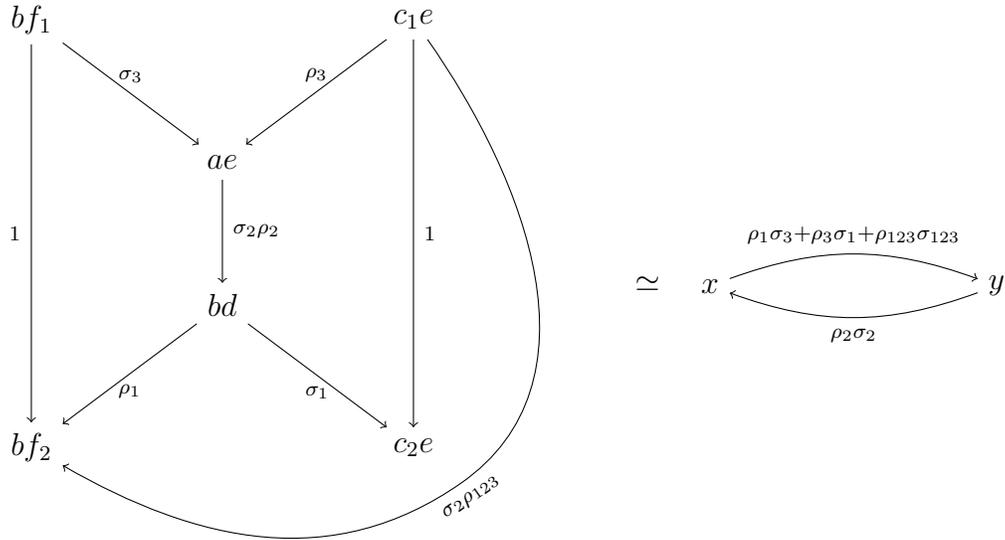
\begin{figure}
\begin{center}
  \[
  \mathcenter{
  \begin{tikzpicture}[y=45pt,x=.8in]
    \node at (0,3) (w1) {$bf_1$} ;
    \node at (2,3) (z1) {$c_1 e$} ;
    \node at (1,2) (y)  {$ae$} ;
    \node at (1,1) (x)  {$bd$} ;
    \node at (0,0) (w2) {$bf_2$} ;
    \node at (2,0) (z2) {$c_2 e$} ;
    \draw[->] (w1) to node[above] {\lab{\sigma_3}} (y) ;
    \draw[->] (z1) to node[above] {\lab{\rho_3}} (y) ;
    \draw[->] (y)  to node[right] {\lab{\sigma_2 \rho_2}} (x) ;
    \draw[->] (x)  to node[below] {\lab{\rho_1}} (w2) ;
    \draw[->] (x)  to node[below] {\lab{\sigma_1}} (z2) ;
    \draw[->] (w1) to node[left] {\lab{1}} (w2) ;
    \draw[->] (z1) to node[right] {\lab{1}} (z2) ;
    \draw[->] (z1) to[out=-55,in=35] (2.25,-0.25) to[out=-145,in=-30] node[pos=0,below,sloped]
       {\lab{\sigma_{2}\rho_{123}}} (w2) ;
  \end{tikzpicture}}
  \quad
  \mathcenter{\simeq}
  \quad
  \mathcenter{\begin{tikzpicture}[y=45pt,x=.8in]
    \node at (0,0) (x) {$x$} ;
    \node at (1.5,0) (y) {$y$} ;
    \draw[->, bend left=20] (x) to node[above] {\lab{\rho_1\sigma_3+\rho_3\sigma_1+\rho_{123}\sigma_{123}}} (y) ;
    \draw[->, bend left=20] (y) to node[below] {\lab{\rho_2\sigma_2}} (x) ;
  \end{tikzpicture}}
  \]
\end{center}
\caption[The bimodule $\CFDDa(\Id,0)$ for the torus]{\textbf{The type \DDm\  bimodule $\CFDDa(\Id,0)$.}
  The labels on the arrows indicate differentials; 
  for instance, the arrow from $ae$ to $bd$ signifies a term
  $(\rho_2\cdot \sigma_2)\otimes bd$ in $\partial ae$. The form on
  the left comes from the Heegaard diagram in
  Figure~\ref{fig:AAinvpair}, and is bounded in an analogous sense to
  Definition~\ref{def:BoundedTypeD}. The form on the right is homotopy 
  equivalent, but no longer bounded.}
\label{fig:DDmodAns}
\end{figure}

The type \DDm\  bimodule can also be calculated using the same
Heegaard diagram. The result, which is given graphically in
Figure~\ref{fig:DDmodAns}, is given as follows:
\begin{align*}
  \bdy bf_1 &= (\sigma_3\otimes 1)\otimes ae+bf_2&
  \bdy c_1e &=(1\otimes \rho_3)\otimes ae+c_2e+(\sigma_2\otimes \rho_{123})\otimes bf_2\\
  \bdy ae &=(\sigma_2\otimes \rho_2)\otimes bd&
  \bdy bd &=(1\otimes \rho_1)\otimes bf_2+(\sigma_1\otimes 1)\otimes c_2e\\
  \bdy bf_2 &=0 &
  \bdy c_2e &=0.
\end{align*}

(We have relabeled the
$\rho_i$ and $\sigma_i$ to conform to our conventions for
type $D$ boundary.)
This simplifies to a bimodule with two generators $x$ and $y$ with
\begin{align*}
  \partial x &= (\rho_1\sigma_3+ \rho_3\sigma_1 + \rho_{123}\sigma_{123}) \otimes y &
  \partial y &= (\rho_2\sigma_2)\otimes x.
\end{align*}

\section{From \textalt{$\HFK$}{HFK} to \textalt{$\CFDa$}{CFD\textasciicircum} for arbitrary integral framings}\label{sec:hfk-to-CFD-general}

We turn now to the generalization
of Theorem~\ref{thm:HFKtoHFD} for arbitrary (integral) framings.

\begin{theorem}
  \label{thm:HFKtoHFDframed}
  Let $K\subset S^3$ be a knot and let
  $\CFKm(K)$ be a reduced model for the knot Floer complex of $K$.
  Let $Y$ be the bordered three-manifold
  $S^3-\nbd{K}$, given any integral framing~$n$. The associated type
  $D$ module $\CFDa(Y)$ can be extracted from $\CFKm(K)$ using the
  procedure described in Theorem~\ref{thm:HFKtoHFD}, except for the
  unstable chain, whose precise form depends on the framing parameter
  $n$.  There are three cases for the form of the unstable chain. 
  When $n<2\tau(K)$, the unstable chain has the form
  $$\xi_0\stackrel{D_1}{\longrightarrow}\mu_1\stackrel{D_{23}}\longleftarrow
  \mu_2\stackrel{D_{23}}\longleftarrow
  \cdots
  \stackrel{D_{23}}\longleftarrow
  \mu_m
  \stackrel{D_3}\longleftarrow
  \eta_0,$$  where $m=2\tau(K)-n$.
  When $n=2\tau(K)$, the
  unstable chain has the form 
  $$\xi_0\stackrel{D_{12}}{\longrightarrow}
  \eta_0.$$ 
  Finally, when $n>2\tau(K)$ the unstable chain has the form
  $$\xi_0 \stackrel{D_{123}}{\longrightarrow} \mu_1 
  \stackrel{D_{23}}{\longrightarrow} \mu_2
  \dots
  \stackrel{D_{23}}{\longrightarrow} \mu_{m}
  \stackrel{D_{2}}{\longrightarrow} \eta_0,
  $$
  where $m=n-2\tau(K)$.
  
  The gradings are determined as follows:
  \begin{itemize}
  \item The grading set is $G/\lambda^{-1} \gr(\rho_{23})^{-n} \gr(\rho_{12})^{-1}$.
  \item Recall that any element $\x_0$ of $V^0$ is represented by a generator
    $\x$ for the knot Floer complex. The grading of $\x_0$ in the above grading
    set is determined by the Maslov grading $M$
    and Alexander grading $A$ of $\x$ by the formula
    \begin{equation}\label{eq:appendix:gr-on-CFK-to-CFD}
    \gr(\x_0)=\lambda^{M(\x)-2A(\x)}\left(\gr(\rho_{23})\right)^{-A(\x)}
    = (M-\frac{3}{2}A;0,-A).
    \end{equation}
  \end{itemize}
\end{theorem}
(The reader is warned that $n$ here has the opposite sign from
Theorem~\ref{thm:HFKtoHFD}, since the framing there is $-n$.)

\begin{proof}
  In view of Theorem~\ref{thm:bimodule-pairing}, we need only see what
  happens as we iteratively tensor the $D$ module for a knot complement with
  sufficiently small framing
  parameter, as calculated in Theorem~\ref{thm:HFKtoHFD}, with the
  bimodule $\CFDAa(\tau_{\mu},0)$.

  We denote the generators in the idempotent $\iota_0$ for $\CFDa(Y)$
  by $x_i$, and the generators in idempotent $\iota_1$ by $y_j$. In
  particular, we write $x_0$ and $x_1$ for the elements $\xi_0$ and
  $\eta_0$ in the statement of the above theorem.

  In the tensor product, each generator $x_i$ gives rise to a pair of
  generators ${\mathbf p}\DT x_i$ and ${\mathbf r}\DT x_i$, while each generator $y_j$ gives
  rise to a generator ${\mathbf q}\DT y_j$.  Note that the generators ${\mathbf p}\DT x_i$ are
  the generators in the $\iota_0$-idempotent for the tensor product.
  
  If we think of the $x_i$ as the generators in a
  horizontally-simplified basis in the sense of
  Definition~\ref{def:Simplified}, then if there are $2k+1$ basis
  elements, $k$ of them are initial points of horizontal arrows.
  For each initial point of a horizontal arrow $x_i$,
  we can cancel ${\mathbf r}\DT x_i$ with some corresponding element of the
  form ${\mathbf q}\DT y_1$; i.e., if we have a chain
  \[\begin{CD}
    x_i 
  @>{D_{3}}>> 
  y_1
  @>{D_{23}}>>
  \cdots
  @>{D_{23}}>>
  y_\ell
  @>{D_{2}}>>
  x_j;
  \end{CD}\]
  this is transformed to
\[  \begin{CD}
  {\mathbf p}\DT x_i
  @>{D_{3}}>> 
  {\mathbf q}\DT y_2
  @>{D_{23}}>>
  {\mathbf q}\DT y_3 
  @>{D_{23}}>>
  \cdots
  @>{D_{23}}>>
  {\mathbf q}\DT y_\ell
  @>{D_{23}}>>
  {\mathbf r}\DT x_j \\
  @A{D_2}AA @A{D_{23}}AA & & & & & & @VV{D_2}V \\
  {\mathbf r}\DT x_i @>{D_{\emptyset}}>> {\mathbf q}\DT y_1 & & & & & & & & {\mathbf p}\DT x_j
  \end{CD}.
  \]
  Setting $y_i'={\mathbf q}\DT y_{i+1}$ for $i=1\dots\ell-1$ and $y'_\ell={\mathbf r}\DT
  x_j$ and dropping the canceling pair ${\mathbf r}\DT x_i$ and ${\mathbf q}\DT y_1$,
  our new chain has the same form as the original chain.

  Vertical chains are taken 
  to chains of the same form (without any need for
  cancellation).

  Finally, we must observe the behavior of the ``unstable chain'' under tensor
  product. 
  To this end, 
  recall the generator $x_0$ which is
  neither the initial nor the final intersection point
  of a horizontal arrow. There is then a chain of 
  the form
  $${x_0}
  \stackrel{D_3}{\longrightarrow} {y_1}
  \stackrel{D_{23}}{\longrightarrow} {y_2}
  \stackrel{D_{23}}{\longrightarrow} {y_3}
  \stackrel{D_{23}}{\longrightarrow} {\cdots}
  \stackrel{D_{23}}{\longrightarrow} {y_m}
  \stackrel{D_{1}}{\longleftarrow} {x_1}
  $$
  with $m>1$. Under the tensor product, 
  canceling
  a differential from ${\mathbf r}\DT x_0$ to ${\mathbf q}\DT y_1$ as before,
  this chain is carried to a
  similar chain from 
  with length one less, from ${\mathbf p}\DT x_0$ to
  ${\mathbf p}\DT x_1$. 

  In the case where $m=1$, we have the chain
  $${x_0}
  \stackrel{D_3}{\longrightarrow} {y_1}
  \stackrel{D_{1}}{\longleftarrow} {x_1}.
  $$
  This chain is transformed to another configuration where we can once again 
  cancel the pair ${\mathbf r}\DT x_0$ and ${\mathbf q}\DT y_1$, to get the chain
  $${{\mathbf p}\DT x_0}  \stackrel{D_{12}}{\longleftarrow} {{\mathbf p}\DT x_1}.  $$

  Now if we start with a chain of the form
  $${x_0}  \stackrel{D_{12}}{\longleftarrow} {x_1},  $$
  this in turn is transformed into 
  $$
  {\mathbf p}\DT x_0  \stackrel{D_{2}}{\longleftarrow} {\mathbf r}\DT x_0 
  \stackrel{D_{123}}{\longleftarrow} {\mathbf p}\DT x_1
  \stackrel{D_{2}}{\longleftarrow} {\mathbf r}\DT x_1. 
  $$ 
  The arrow from 
  ${\mathbf r}\DT x_1$ should be viewed as part of another chain.

  The remaining elements lie in an unstable chain of the
  form
  $$x_0 \stackrel{D_{2}}{\longleftarrow} y_1 
  \stackrel{D_{23}}{\longleftarrow} y_2
  \dots
  \stackrel{D_{23}}{\longleftarrow} y_{m}
  \stackrel{D_{123}}{\longleftarrow} x_1
  $$
  with $m=1$. Now it is easy to see
  inductively that tensoring with our bimodule takes such an unstable
  chain with parameter $m$ to one with parameter $m+1$.

  We leave verification of the statement about gradings, which
  requires the gradings on the bimodules from~\cite{LOT2}, to the
  interested reader.  The precise relationship between the length (and
  type) of the unstable chain, the framing parameter, and $\tau$
  follows as in the proof of Theorem~\ref{thm:HFKtoHFD}.
\end{proof}


\backmatter

{\footnotesize \printglossary[title=Index of Notation,style=hangtree]}

\renewcommand{\indexname}{Index of Definitions}
\printindex

\bibliographystyle{hamsalpha}\bibliography{heegaardfloer}




\begin{flushright}
  \vspace{0.5cm}
  January 21, 2021
  \vspace{0.5cm}
\end{flushright}

\begin{center}
  \textsc{Errata}
\end{center}
\phantomsection
\addcontentsline{toc}{chapter}{Errata}
\renewcommand\thefigure{Err.\arabic{figure}}
\renewcommand\theHfigure{Err.\arabic{figure}}
\setcounter{figure}{0}
\noindent

\textbf{Generalized coefficient maps}

The authors are grateful to Wenzhao Chen for pointing out an error in the
proof of Proposition 11.30 and the construction of the maps $D_{0123}$
and $D_{2301}$ appearing there.

This error can be corrected by adding a term to the definition of
$D_{0123}$ at the beginning of Section~11.6, as follows. In addition
to the terms written,
the function $D_{0123}$ should also include a count of points in
a moduli space
${\mathcal M}^B({\mathbf x},{\mathbf y},\{\mathfrak o\})$, where
${\mathfrak o}$ denotes the simple Reeb orbit that covers the entire
boundary with local multiplicity~$1$. (The same term must be added
also to $D_{1230}$, $D_{2301}$, and $D_{3012}$.)

\begin{figure}
  \centering
  \includegraphics[scale=.666667]{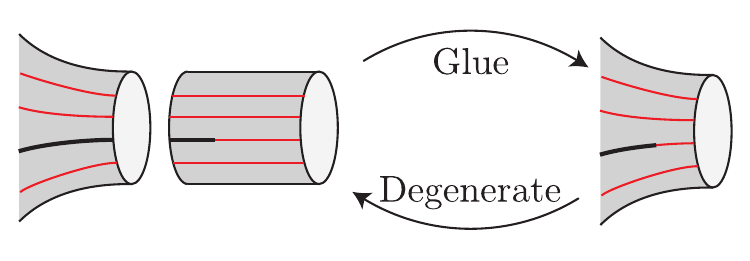}
  \caption{\textbf{An orbit curve.} Gluing the orbit curve (at east
    infinity) to a curve in $\Sigma\times[0,1]\times\mathbb{R}$ asymptotic to
    a length-4 Reeb chord gives a curve in
    $\Sigma\times[0,1]\times\mathbb{R}$ asymptotic to a Reeb orbit.}
  \label{fig:orbit-curve}
\end{figure}

This change also affects the proof of Proposition 11.30.  We had
neglected to mention that there are, in addition to
the two-story and boundary degeneration ends (which were accounted
for), split curve ends. The number of such join curve ends is identified with the count of points in moduli spaces
${\mathcal M}^B({\mathbf x},{\mathbf y},\{\rho\})$, taken over all
possible Reeb chords $\rho$ that cover the boundary with multiplicity one. These ends in turn cancel
against the ends of the moduli spaces ${\mathcal M}^B({\mathbf
  x},{\mathbf y},\{\mathfrak o\})$, where the orbit degenerates to the
boundary.  Counting this latter degeneration involves a new kind of rigid curve at
east infinity, called an \emph{orbit curve}.  (See
Figure~\ref{fig:orbit-curve}. These curves did not
appear in our analysis of curves at east infinity in Chapter 5, since
they do cover the basepoint $z$, with multiplicity $1$; compare
Figure~5.10.)

With this correction to the coefficient maps, the statement of
Proposition~11.30 and the rest of the discussion remain valid without
change.

This phenomenon, the cancellation of orbit curve ends against join
curve ends, plays an important role in upcoming work, wherein we
generalize the bordered Floer homology with torus boundary to $HF^-$.

\vspace{1em}
\noindent
\textbf{Other typos and minor corrections}

We have corrected minor typos on pages 68, 69, 70, 71, 114, 157, 158,
182, 216, and 230 of this version. Thanks to Vinicius Ambrosi,
Ina Petkova, Mike Wong,
and everyone else who
has pointed out typos.



\end{document}